\documentclass{amsart}
\usepackage{amssymb,amsmath,cite,geometry}

\usepackage{color}

\numberwithin{equation}{section}

\newtheorem{notation}{Notation}[section]

\newtheorem{theorem}{Theorem}[section]
\newtheorem{lemma}{Lemma}[section]
\newtheorem{corollary}{Corollary}[section]
\newtheorem{prop}{Proposition}[section]
\newtheorem{remark}{Remark}[section]
\newtheorem{assumption}{Assumption}[section]

\newcommand\y{{{\langle y \rangle}}}
\allowdisplaybreaks[4]

\begin{document}
\title[Well-posedness of Prandtl system with degenerate critical points]
{Well-posedness of the two-dimensional unsteady Prandtl system in Sobolev space with degenerate critical points}

\author{Shi-yong Zhu}
\address{Shi-yong Zhu}
\address{College of Mathematics and Systems Science, Shandong University of Science and Technology, Shandong, 266590, P. R. China.}
\email{shiyong\_zhu@163.com}

\author{Ya-Guang Wang}
\address{Ya-Guang Wang
	\newline\indent
	School of Mathematical Sciences, MOE-LSC and SHL-MAC, Shanghai Jiao Tong University,
	Shanghai, P. R. China}
\email{ygwang@sjtu.edu.cn  }

\date{}




\begin{abstract}
This paper is devoted to the well-posedness of classical Prandtl equations in a finite order Sobolev space. For a initial data with degenerate critical points and general outflow, we obtain the local-in-time existence and uniqueness of the solution to the Prandtl equations in a Sobolev space, by introducing a new iteration scheme and linear cancelation. This result shows that Oleinik's monotonicity condition is not a necessary condition for the Prandtl equations to be well-posed in Sobolev spaces and provides evidence to demonstrate that zero shear stress does not necessarily lead to boundary layer separation in two-dimensional unsteady boundary layers.
\end{abstract}

\maketitle

\tableofcontents

\section{Introduction}\setcounter{section}{1}\setcounter{equation}{0}

Boundary layer theory is one of the cornerstones of modern fluid mechanics theory, which is widely applied in various fields, such as aeronautics and astronautics, marine engineering, industrial design, and so on. It was first proposed by Prandtl \cite{Prandtl} in 1904 to address the disconnect between theoretical fluid mechanics and experimental hydraulics. According to Prandtl's boundary layer theory, the flow field around a body can be divided into two parts, in a thin layer close to the boundary whose thickness dependent on viscosity, both viscous force and inertial force play a dominant role in the fluid flow, while away from the boundary, the viscous force can be ignored compared to inertial force, and the flow can be considered as inviscid one. This thin layer is the boundary layer. By the balance of viscous force and inertial force in the boundary layer, Prandtl derived the thickness of the boundary layers and obtained the boundary layer equations, that is, the famous Prandtl equations, by simplifying the Navier-Stokes equations. One can consult \cite{SG} for more physical background about the boundary layer.

Due to its important practical significance, the boundary layer theory has received widespread attention since its proposal.
Although more than a century has passed, the boundary layer theory has yet to be fully and rigorously proven. For steady flow, the boundary layer theory is proved rigorously recently, cf.\cite{GM1,GI,GN-2,GZ,I}, while for unsteady flow, the boundary layer theory was proved only in some special cases, cf.\cite{Lop,LW,SC2,WWZ,Mae,GMM}, and there exist some results where the expansion fails under certain conditions, cf.\cite{Grenier,GN1,GN2,GGN,GN-1}. To prove boundary layer theory mathematically, a key issue is to establish the well-posedness of the Prandtl equations. In this paper, we focus on the two-dimensional unsteady Prandtl equations. Let $\Omega=\{(x,y)|x\in\mathbb{T},y\in[0,+\infty)\}$, where $\mathbb{T}$ is periodic. For a constant $T>0$, the two-dimensional unsteady Prandtl equations with non-slip boundary condition in $\Omega_{T}=[0,T]\times\Omega$ reads
\begin{equation}\label{Prandtl}
\begin{cases}
\partial_{t}u+u\partial_{x}u+v\partial_{y}u=\partial_{y}^{2}u-\partial_{x}P,\\
\partial_{x}u+\partial_{y}v=0,\\
u|_{y=0}=v|_{y=0}=0,~\lim\limits_{y\to+\infty} u=U(t,x),\\
u|_{t=0}=u_{in}(x,y)
\end{cases}
\end{equation}
where $(u(t,x,y),v(t,x,y))$ is the velocity field in the boundary layer, $U(t,x)$ and $P(t,x)$ are traces at the boundary $\{y=0\}$ of the tangential velocity and pressure of the Euler outer flow respectively, interrelated through Bernoulli's law
\begin{equation}\label{B}
\partial_{t}U+U\partial_{x}U=-\partial_{x}P.
\end{equation}
Due to the property of degenerate parabolicy and non-locality, there is a loss of tangential derivatives in the nonlocal term $v\partial_{y}u$, making it difficult to establish the well-posedness theory for the problem \eqref{Prandtl}. Under a special structural assumption, the well-posedness theory can be fulfilled. The first well-posedness result was obtained by Oleinik for the laminar boundary layer on a flat plat under the monotonicity assumption $\partial_{y}u>0$ in \cite{Oleinik}(see also \cite{Ole}). By applying Crocco transformation, the Prandtl equations is converted to a single degenerate parabolic equation, then by the maximum principle for the parabolic equation Oleinik  establishes the local well-posedness theory for the Prandtl equations. Later, the local well-posedness
result in the Oleinik's monotonic class was rebuilt in the Sobolev spaces. With the help of Nash-Moser-H\"{o}rmander iteration, the local well-posedness result was obtained in the Oleinik's monotonic class by using a direct energy method in \cite{AWXY}, and in \cite{MW}, the local well-posedness was also established in the Oleinik's monotonic class by taking advantage of the nonlinear cancelation property of the Prandtl equations. This two works utilized essentially identical good unknowns to overcome the loss of tangential derivative caused by degenerate parabolicy and non-locality. If the pressure is favorable, the global existence of  weak solutions to the problem \eqref{Prandtl} is obtained in the Oleinik's monotonic class in \cite{XZ}. Without any structural assumption, there are also some well-posedness results in analytic space, cf.\cite{LCS,KV,IV,SC1,ZZ}, since the loss of tangential derivative can be overcome by analyticity in  analytic spaces. An interesting work is \cite{KMVW}, in which the local existence and uniqueness of the problem \eqref{Prandtl} was obtained, for a class of initial data that are monotone on a number of intervals and analytic on the complement of these intervals.
The Oleinik's monotonicity assumption was once considered a necessary condition for the well-posedness of Prandtl equations in Sobolev spaces. Without Oleinik's monotonicity assumption, there are some ill-posedness and blow-up results for the Prandtl equations in Sobolev spaces, cf. \cite{EE,GD,G-N,GN-1,KVW,LY2}. Under the assumption that the initial velocity has a non-degenerate critical point over $\mathbb{R}_{+}$, the strong linear ill-posedness of the Prandtl equations in a Sobolev space is established in \cite{GD}. Later, this strong ill-posedness was extended to the  nonlinear case. These ill-posedness results show that without Oleinik's monotonicity, the well-posedness of Prandtl equations can not be obtained for general data in a finite order Sobolev space. Inspired by these results, many scholars endeavor to establish the well-posedness of Prandtl equations in Gevrey class, cf.\cite{DG,GM,LY,LWX}.

The result in \cite{GD} do not rule out the well-posedness of the Prandtl equation in the Sobolev space with the critical point being degenerate, so it is possible to get the well-posedness of the Prandtl equations under the assumption that the initial velocity has degenerate critical points.  This is the main purpose of our paper. For an initial data with degenerate critical points and general outflow, we obtain the local-in-time existence and uniqueness of the solution to the Prandtl equation in Sobolev space. Notice that the appearance of degenerate critical points indicates a lack of Oleinik's monotonicity, the methods introduced in \cite{AWXY} and \cite{MW} cannot work. To overcome this difficulty, we introduce a new iteration scheme and take advantage of the linear cancelation property of corrected increment problems. By a direct energy method we obtain the local-in-time existence and uniqueness of the solution to the problem \eqref{Prandtl}. This result can be extended directly to the case of finitely many degenerate critical points.

It is worth noting that our result does not conflict with \cite{Grenier2}. In \cite{Grenier2}, for a background flow with a degenerate critical point initially, a nonlinear instability result of the Prandtl boundary layers is obtained. This instability result shows that the Prandtl asymptotic expansion is not valid for some special case, but does not negate the well-posedness of Prandtl equations.

From a physical perspective, the existence of critical points implies the presence of zero shear stress within the flow field. Our result rigorously proves that although zero shear stress occurs within the flow field, the problem \eqref{Prandtl} can still be well-posed, at least for a short time. This means that zero shear stress does not necessarily lead to boundary layer separation in two-dimensional unsteady boundary layers, which is quite different from the steady boundary layer. This is consistent with experimental observations and numerical computations.

This paper is organized as follows. In Section 2, we introduce the main assumptions and the result of this paper, as well as some notation for later use. In Section 3, we introduce our methodology and iteration scheme  and construct the approximate solutions to  problem \eqref{Prandtl} by corrected increment solutions. In Section 4, we establish the well-posedness of corrected increment problems using uniform energy estimates. In Section 5, we give  a posteriori estimates of the approximate solutions to ensure that the iterative process can proceed. Finally, we use the well-posedness established in Section 3 for the corrected increment problems and the a posteriori estimates of the approximate solutions to get the local-in-time existence and uniqueness of the solution to the problem \eqref{Prandtl} in Section 6.

\section{Preliminary and main result}

First, let us introduce the solution space which we shall work in.
\begin{notation} (1)
For a given nonnegative function $\xi(y)\in C^{\infty}[0,+\infty)$, and an interger $s\ge 0$, define the weighted Sobolev space as
\begin{align*}
H_{\xi}^{s}(\Omega)=\{f:\Omega\to\mathbb{R}|~\|f\|_{H_{\xi}^{s}(\Omega)}<+\infty\},
\end{align*}
with $\Omega=\{(x,y)| x\in\mathbb{T}, y\ge 0\}$, and the norm being defined as
\begin{align*}
\|f\|_{H_{\xi}^{s}(\Omega)}
=\sum_{\substack{\alpha_{1}+\alpha_{2}\leq s}}\|\partial_{x}^{\alpha_{1}}\partial_{y}^{\alpha_{2}}f\|_{L_{\xi}^{2}(\Omega)}, \quad
\|f\|_{L_{\xi}^{2}(\Omega)}=\left(\int_{\Omega}(\xi(y)f(x,y))^{2}dxdy\right)^{\frac{1}{2}}.
\end{align*}
	
(2) Define the space $\hat{H}_{\psi}^{s}(\Omega)$ by
\begin{align*}
\hat{H}_{\psi}^{s}(\Omega)=\{f:\Omega\to\mathbb{R}|~\|f\|_{\hat{H}_{\psi}^{s}(\Omega)}<+\infty\},
\end{align*}
with the norm
\begin{align}\label{cn}
\|f\|_{\hat{H}_{\psi}^{s}(\Omega)}
= \sum_{\substack{\alpha_{1}+\alpha_{2}\leq s,\\ \alpha_{1}\leq s-1}}\|\partial_{x}^{\alpha_{1}}\partial_{y}^{\alpha_{2}}f\|_{L_{\psi}^{2}(\Omega)}
 +\|\partial_{x}^{s-1}\partial_{y}^{2}f\|_{L_{\hat{\psi}}^{2}(\Omega)},
\end{align}
where
\begin{align}\label{psi}
\psi(t,y)=e^{(\frac{3}{4}-\Lambda t)\delta\y},
\end{align}
for a positive constant $\Lambda$ to be determined later, and $\hat{\psi}(t,y)=\y^{-\frac{1}{2}}\psi(t,y)$ with $\y=1+y$.
	
\end{notation}

The aim of introducing weight $\psi(t,y)$ in the above space is to close energy estimates, in order to control the nonlocal term appeared in  \eqref{Prandtl}.

Next, we introduce two structural assumptions on the initial data  given in the problem \eqref{Prandtl}. The first one says that the monotonic initial velocity $u_{in}(x,y)$ has a degenerate critical point in $\{y>0\}$.

\begin{assumption}\label{inA'}
Assume that there exist $y_\ast>0$, and $0<\check{y}<y_{*}<\hat{y}<+\infty$ such that
the initial data $u_{in}(x,y)$ satisfies
\begin{align}\label{inA1}
\partial_{y}u_{in}(x,y_{*})=\partial_{y}^{2}u_{in}(x,y_{*})=0,~\forall x\in\mathbb{T},
\end{align}
\begin{align}\label{inA2}
\partial_{y}^{3}u_{in}\geq\mathcal{C}>0~\text{in}~[\check{y},\hat{y}],
\end{align}
\begin{align}\label{inA3}
\partial_{y}u_{in}\geq c_{in}\varpi (y) e^{-\delta\y}, \quad \forall x\in  \mathbb{T}, ~y\ge 0,
\end{align}
where $\varpi(y)\in C^\infty[0,+\infty)$ is a nonnegative function satisfying
\begin{align}\label{pi}
\varpi(y)=
\begin{cases}
\mathcal{C}(y-y_{*})^{2},~\text{in}~[\check{y},\hat{y}],\\
1,~\text{in}~[0,+\infty)\setminus[\check{y},\hat{y}].
\end{cases}
\end{align}
with $\delta, c_{in}$ and $ \mathcal{C}$ positive constants.

\end{assumption}

To study the  problem \eqref{Prandtl} in the space $\hat{H}_{\psi}^{s}(\Omega)$, we need to homogenize the limit state as $y\to +\infty$,  meanwhile also need to pay attention on the degeneracy of the initial velocity $u_{in}(x,y)$ at $y_{*}$. For this purpose, choose $u_{0}\in C^{\infty}(\Omega_{T})$ be a nonnegative function such that
\begin{align}\label{u0}
u_{0}(t,x,y)=
\begin{cases}
c_{0}\phi(y)+c_{0}ty,~&\text{in}~[0,T]\times\mathbb{T}\times[0,\hat{y}],\\
\phi(y)U(t,x),  ~&\text{in}~[0,T]\times\mathbb{T}\times[\hat{y}+1,+\infty),
\end{cases}
\end{align}
where $\phi(y)\in C^{\infty}[0,+\infty)$ is a nonnegative  function satisfying
\begin{align}\label{phi}
\begin{cases}
\phi(0)=0,\\
\phi(y)=1-e^{-\delta\y},\quad
\forall y\in[\hat{y}+1,+\infty),\\
\phi'(y)=\mathcal{C}_{0}(y-y_{*})^{2},\quad  \forall y\in[\check{y},\hat{y}],\\
\phi'(y)>0, \quad \forall y\in[0,\check{y}],
\end{cases}
\end{align}
with $c_{0}, \mathcal{C}_{0}$ being small positive constants such that $\partial_{y}u_{0}>0$ in $\Omega_{T}\setminus ([0,T]\times\mathbb{T}\times\{y_{*}\})$.

The second structural assumption concerns the detail behavior of the initial velocity $u_{in}(x,y)$ near the critical point $y_{*}$. To state this assumption, we first introduce the following weight functions.
For any $(x,y)\in\mathbb{T}\times\left([0,+\infty)\setminus\{y_{*}\}\right),$ let
$$\mathcal{P}_{in}(x,y)=\frac{\chi(y)}{\sqrt{\partial_{y}u_{in}(x,y)}}+1-\chi(y),
\qquad
\rho_{in}(x,y)=\frac{\partial_{y}^{2}u_{in}(x,y)}{\partial_{y}u_{in}(x,y)},$$
and
$$\omega_{a}^{in}(y)=\frac{\chi(y)}{(y-y_{*})^{2a}}+1-\chi(y),$$
with a positive index $a$, where $\chi(y)\in C^\infty(0,+\infty)$ is a nonnegative function satisfying
\begin{equation}\label{chi}
{\bf supp}~\chi\subset(\check{y},\hat{y}), \quad
\chi(y)\equiv 1, \quad {\rm on}~ [\check{y}_{*},\hat{y}_{*}]
\end{equation}
for fixed   $\check{y}_{*},\hat{y}_{*}\in(\check{y},\hat{y})$ such that $y_{*}\in(\check{y}_{*},\hat{y}_{*})$.

\begin{assumption}\label{inA}
For the initial data $u_{in}(x,y)$, let $\tilde{u}_{in}(x,y)=u_{in}(x,y)-u_{0}(0,x,y)$, with $u_0$ given in \eqref{u0}.
For some $s\in\mathbb{N}_{+}$, assume there exists positive constants $\varepsilon_{0}\in\left(0,\frac{1}{4}\right)$ and $\mathcal{C}_{in}$ such that
\begin{align}\label{inA4}
&\|\tilde{u}_{in}\|_{\hat{H}_{\psi_{in}}^{s}(\Omega)}
+\|\omega_{\varepsilon_{0}}^{in}\mathcal{P}_{in}\partial_{y}\partial_{x}^{s}\tilde{u}_{in}\|_{L_{\hat{\psi}_{in}}^{2}(\Omega)}
+\|\omega_{\varepsilon_{0}}^{in}\mathcal{P}_{in}^{2}\partial_{x}^{s}\tilde{u}_{in}\|_{L_{\hat{\psi}_{in}}^{2}(\Omega)}\\
&+\|\omega_{\varepsilon_{0}}^{in}\mathcal{P}_{in}\partial_{y}\partial_{x}^{s+1}\tilde{u}_{in}\|_{L_{\check{\psi}_{in}}^{2}(\Omega)}
+\|\omega_{\varepsilon_{0}}^{in}\mathcal{P}_{in}^{2}\partial_{x}^{s+1}\tilde{u}_{in}\|_{L_{\check{\psi}_{in}}^{2}(\Omega)}
\leq \mathcal{C}_{in},\nonumber
\end{align}
with $\psi_{in}=\psi(0,y)=e^{\frac{3}{4}\delta\y}$, $\hat{\psi}_{in}=\y^{-\frac{1}{2}} \psi_{in}$, $\check{\psi}_{in}=\y^{-1} \psi_{in}$,
and
\begin{align}\label{inA5}
|D^{\gamma}\tilde{u}_{in}|\leq \mathcal{C}_{in}\y^{-1}e^{-\delta\y}~\text{in}~\Omega,
\end{align}
for any $\gamma\in\Gamma_{9}$, with $\Gamma_{m}=\{(\gamma_{1},\gamma_{2})\in\mathbb{N}^{2}| \gamma_{1}+\gamma_{2}\leq m\}$  denoting a set of multiindex. Moreover, assume that
\begin{align}\label{inA6}
\begin{cases}
|\partial_{y}\rho_{in}|\leq \mathcal{C}_{in}\y^{-1}, |\partial_{x}\rho_{in}|\leq \mathcal{C}_{in},~\text{in}~\mathbb{T}\times[\hat{y},+\infty),\\
\rho_{in}\leq -\delta,~\text{in}~\mathbb{T}\times[Y,+\infty),
\end{cases}
\end{align}
for fixed constants $Y>\hat{y}$.

\end{assumption}

The main result of this paper is the following one.
\begin{theorem}\label{MR}
Assume that the initial data $u_{in}$ of the problem \eqref{Prandtl} satisfy  Assumptions \ref{inA'} and \ref{inA} for  a fixed odd integer $s\geq 11$, and the compatibility conditions of the problem \eqref{Prandtl} up to order $s-1$. Moreover, the outflow
$U\in C([0,T];H^{s+2}(\mathbb{T}))\bigcap C^{1}([0,T];H^{s+1}(\mathbb{T}))\bigcap\limits_{i=2}^{\frac{s-3}{2}}C^{i}([0,T];H^{s-2i}(\mathbb{T}))$ is positive.
Then, there exists $t_{*}\in(0,T]$ such that the problem \eqref{Prandtl} admits a unique classical solution $u$ within $[0,t_{*}]$ satisfying that
$$w:=u-u_{0}\in L^{\infty}((0, t_{*}), \hat{H}_{\psi}^{s}(\Omega))\cap L^{2}((0, t_{*}), H_{\hat{\psi}}^{s}\cap H_{\check{\psi}}^{s+1}(\Omega)),$$
where $ \check{\psi}(t,y)=\y^{-1}\psi(t,y)$, and the function $u_0(t,x,y)$ is given in \eqref{u0}.
\end{theorem}

\begin{notation}
For convenience, we introduce some notations that will be used throughout this paper. For any $y\in(0,+\infty)$, set $\hat{\Omega}^{y}=\mathbb{T}\times[0,y]$,  $\check{\Omega}^{y}=\mathbb{T}\times[y,+\infty)$, and for any $t\in[0,T]$, let
$\hat{\Omega}_{t}^{y}=[0,t]\times\hat{\Omega}^{y}$,
$\check{\Omega}_{t}^{y}=[0,t]\times\check{\Omega}^{y}$ respectively.
Set ${\Omega}^{*}=\mathbb{T}\times[\check{y}_{*},\hat{y}_{*}]$ and denote $\hat{\Omega}^{\hat{y}}, \check{\Omega}^{\check{y}}$ by $\hat{\Omega}, \check{\Omega}$, respectively. For any $t\in[0,T]$, let $\Omega_{t}=[0,t]\times\Omega$, and accordingly, $\hat{\Omega}_{t}^{y}=[0,t]\times\hat{\Omega}^{y}$, $\check{\Omega}_{t}^{y}=[0,t]\times\check{\Omega}^{y}$.
We use $C$ to represent a positive constant which may be different from line to line.
For any constant $\sigma$, we use $C_{\sigma}$ to represent a positive constant depending only on $\sigma$ which may be different from line to line. We use $\tau$ to stand for a constant to be determined that belongs to $(0,1)$.

\end{notation}

\section{Methodology and iteration scheme}

In this section, we explain the main idea and steps to construct a solution in the Sobolev space to the problem \eqref{Prandtl}, when the initial tangential velocity $u_{in}$ is monotonic with a degenerate critical point in the normal variable as given in Assumption \ref{inA'}, which broken the strict monotonicity assumption required in the classical work, cf. \cite{Oleinik, Ole, AWXY, MW}.

The Prandtl equations in \eqref{Prandtl} is degenerate parabolic for tangential velocity $u$, and there is a loss of tangential derivative due to the nonlocal term $v\partial_{y}u$ and $v$ depending on $\partial_{x}u$ from the divergence-free constrain, so the problem \eqref{Prandtl} might be ill-posed  in a finite order Sobolev space if without any structural assumption on $u$, e.g. it was shown in \cite{GD} that the linearized problem of \eqref{Prandtl} is ill-posed in Sobolev spaces if the monotonic background velocity has a non-degenerate critical point. In Oleinik's pioneering work, she obtained the local well-posed of the problem \eqref{Prandtl} in H\"{o}lder spaces by requiring a strictly monotonicity condition $\partial_{y}u>0$. Recently, under Oleinik's same monotonic assumption, a cencellation mechanism was observed by introducing a good unknown in \cite{AWXY} and \cite{MW} independently to overcome the difficulty caused by the loss of derivative,
\begin{align}\label{Good}
\mathcal{W}:=\partial_{y}u\partial_{y}\left(\frac{u-U}{\partial_{y}u}\right)=\partial_{y}u-\frac{\partial_{y}^{2}u}{\partial_{y}u}(u-U),
\end{align}
and the well-posedness of \eqref{Prandtl} in the finite order Sobolev space is obtained by a direct energy method.

In this paper, we shall study the well-posedness of the problem \eqref{Prandtl} when the initial tangential velocity is monotonic in the normal variable but has a degenerated critical point, that is, \eqref{inA1}-\eqref{inA3} holds. In particular, the assumption \eqref{inA1} violates Oleinik's strict monotonicity condition, thus the approach of introducing the good unknown \eqref{Good} cannot work anymore. One idea for overcoming this difficulty is to add a "patch" to the left-hand sides of \eqref{Prandtl}$_{1}$ to gain Oleinik's monotonicity, and then to add an artificial viscosity in addition to help control the added term. In more detail, we substitute the following equation for \eqref{Prandtl}$_{1}$,
\begin{align}\label{PwP}
\partial_{t}u+u\partial_{x}u-\partial_{y}^{-1}[\partial_{x}u](\partial_{y}u+\varsigma)
=\partial^2_{y}u+\eta\partial_{x}^{2}(u-u_{0})-\partial_{x}P,
\end{align}
where $\varsigma=\varsigma(y)$ is a nonnegative function supported near the critical points $y=y_{*}$ of $u$, with upper bound $\bar{\varsigma}$, and $\eta$ is a positive constant. In following proof of the main result, to obtain the decay estimates of the solution, $\eta$ is set to be a function of $(t,y)$ exponential decaying in the variable $y$, and the decay rate increases over time. We shall need $\bar{\varsigma}$ to be much smaller than $\eta^{2}$ in the support of $\varsigma$. With the "patch" $\varsigma\partial_{y}^{-1}[\partial_{x}u]$ given in \eqref{PwP}, Oleinik's monotonicity condition is fulfilled, and thus one can introduce an unknown
\begin{align*}
\mathbf{w}:=(\partial_{y}u+\varsigma)\partial_{y}\left(\frac{\partial_{x}^{s}(u-u_{0})}{\partial_{y}u+\varsigma}\right)
           =\partial_{y}\partial_{x}^{s}(u-u_{0})-\frac{\partial_{y}^{2}u+\varsigma'}{\partial_{y}u+\varsigma}\partial_{x}^{s}(u-u_{0}).
\end{align*}
to gain the highest-order tangential derivative estimate of $u$, while the low order derivative estimates shall be deduced by a direct energy estimate on $u$.  However, the unknown $\mathbf{w}$ cannot give a good control over $\partial_{x}^{s}u$ when $\bar{\varsigma}$ tends to zero. Moreover, although from \eqref{PwP} the estimates can be obtained  for $\mathbf{w}$  by using the standard energy method, this estimate is not uniform in $\eta$ and $\bar{\varsigma}$, thus one cannot deduce the well-posedness of \eqref{Prandtl} by  taking the limits $\bar{\varsigma}\to0$ and $\eta\to0$.

One of our important observations is that the unknown $\mathfrak{P}^{2}\mathbf{w}$ can  control  over $\mathfrak{P}^{2}\partial_{x}^{s}(u-u_{0})$ uniformly in $\bar{\varsigma}$ and $\eta$, where $\mathfrak{P}=\frac{\chi}{\sqrt{\partial_{y}u+\varsigma}}+1-\chi$ with $\chi$ being given in \eqref{chi} of Section 2, see Lemma \ref{nc}. Thus, if one has the uniform estimates of $\mathfrak{P}^{2}\mathbf{w}$, the uniform estimates of $u$ follow. However, to get uniform estimates of $\mathfrak{P}^{2}\mathbf{w}$ is almost impossible, due to the strong singularity of $\mathfrak{P}^{2}$ when $\bar{\varsigma}\to 0$ at the initial moment. The key important observation is that under the assumption \eqref{inA1}--\eqref{inA3}, Oleinik's monotonicity grows out immediately after the initial moment. In fact, by \eqref{Prandtl}$_{1}$ one has
$$
\partial_{t}\partial_{y}u=\partial^{3}_{y}u-u\partial_{x}\partial_{y}u+\partial_{y}^{-1}[\partial_{x}u]\partial_{y}^{2}u,
$$
which implies that $\partial_{t}\partial_{y}u>0$ near $y=y_{*}$ at least within a short time interval after the initial moment if the problem \eqref{Prandtl} admits a solution in a high order Sobolev space. Thus, $\partial_{y}u$ behaves as $C((y-y_{*})^{2}+t)$ and $\mathfrak{P}\leq \frac{C}{\sqrt{(y-y_{*})^{2}+t}}$ near $y=y_{*}$. This observation enables us to use a suitable weight to control $\mathfrak{P}$. Therefore, by introducing a new unknown $\mathcal{U}=\mathfrak{P}\mathbf{w}$, one could establish a uniform weighted energy estimate for $\mathcal{U}$, which gives the uniform weighted energy estimates of $\partial_{x}^{s}(u-u_{0})$ by the Lemma \ref{nc}.

To construct the solution to the problem \eqref{Prandtl} and derive some uniform estimates, instead of establishing the energy estimates for \eqref{PwP} directly, we use a kind of incremental iterative method to decompose the solution $u$ of problem \eqref{Prandtl} into a sum of a series of corrected increment:
$u=\sum_{n=0}^{\infty}u_{n}$, where $u_{0}$ is given by \eqref{u0} in Section 2, and $u_{n}$ for $n\ge 1$ satisfy a series of semilinear problems. Now, we describe the iterative procedure in detail.

We start with the iteration from $u_{0}$. Take $u_{0}$ as the approximate solution of the zero-th order to the problem \eqref{Prandtl}, and denote it by $\mathfrak{U}_{0}$. If $u$ is a classical solution of the problem \eqref{Prandtl}, then $w_{1}=u-u_{0}$ satisfies the following problem
\begin{equation}\label{E0}
\begin{cases}
\partial_{t}w_{1}
+w_{1}\partial_{x}(\mathfrak{U}_{0}+w_{1})
+\mathfrak{U}_{0}\partial_{x}w_{1}
-\partial_{y}^{-1}[\partial_{x}w_{1}]\partial_{y}(\mathfrak{U}_{0}+w_{1})
-\partial_{y}^{-1}[\partial_{x}\mathfrak{U}_{0}]\partial_{y}w_{1}
=\partial_{y}^{2}w_{1}
+F_{0},\\
w_{1}|_{y=0}=0,~\lim\limits_{y\to+\infty} w_{1}=0,\\
w_{1}|_{t=0}=\tilde{u}_{in}(x,y),
\end{cases}\tag{E$_{1}$}
\end{equation}
where
\begin{align}\label{F0}
F_{0}
=-\partial_{t}u_{0}
 -u_{0}\partial_{x}u_{0}
 +\partial_{y}^{-1}[\partial_{x}u_{0}]\partial_{y}u_{0}
 +\partial_{y}^{2}u_{0}
 -\partial_{x}P.
\end{align}
The problem \eqref{E0} is still nonlinear, nonlocal and degenerate parabolic as \eqref{Prandtl}. Now, we try to find a problem close to \eqref{E0}, but is linear without nonlocality and degeneracy so that its well-posedness can be easily obtained, and then study the difference between the two problems.
To get a liner problem from \eqref{E0}, a natural way is to remove the nonlinear terms  $w_{1}\partial_{x}w_{1}$ and $-\partial_{y}^{-1}[\partial_{x}w_{1}]\partial_{y}w_{1}$ as source terms. Removing the term $w_{1}\partial_{x}w_{1}$ will not cause any trouble when $w_1$ is small, but removing the nonlocal term $-\partial_{y}^{-1}[\partial_{x}w_{1}]\partial_{y}w_{1}$ will bring certain difficulty, as it shall yield the loss of decay, preventing the iteration from proceeding. For this reason, we abandon the pursuit of nonlocality and introduce a function $h=h(y)\in C^{\infty}(\mathbb{R}_{+})$, which is a nondecreasing function such that
\begin{align}\label{h}
{\bf supp}~h(y)\in(\hat{y},+\infty), \quad\text{and}~h(y)\equiv1 \quad {\rm for~all}~ y\in[\hat{y}+1,+\infty).
\end{align}
We retain the nonlocal term in the equation by multiplying it by $h(y)$, then we can get a liner problem in the domain $\hat{\Omega}=\mathbb{T}\times[0,\hat{y}]$,   and avoid the problem of lossing decay. Based on these considerations, by adding a "patch" to the linear nonlocal term and an artificial viscosity, we obtain the following modified increment problem instead of \eqref{E0},
\begin{equation}\label{AE1}
\begin{cases}
\partial_{t}u_{1}
+\mathfrak{U}_{0}\partial_{x}u_{1}
+u_{1}\partial_{x}\mathfrak{U}_{0}
-\partial_{y}^{-1}[\partial_{x}u_{1}](\partial_{y}\mathfrak{U}_{0}+\varsigma_{1})
-\partial_{y}^{-1}[\partial_{x}\mathfrak{U}_{0}]\partial_{y}u_{1}
-h\partial_{y}^{-1}[\partial_{x}u_{1}]\partial_{y}u_{1}\\
\hspace{3.6in}
=\partial_{y}^{2}u_{1}
+\epsilon_{0}^{2}\vartheta\partial_{x}^{2}u_{1}
+F_{0},\\
u_{1}|_{y=0}=0,~\lim\limits_{y\to+\infty} u_{1}=0,\\
u_{1}|_{t=0}=\tilde{u}_{in}(x,y).
\end{cases}\tag{AE$_{1}$}
\end{equation}
Here $\varsigma_{1}=\bar{\varsigma}_{1}\chi_{c}(y)$, where $\bar{\varsigma}_{1}$ and $\epsilon_{0}\in(0,1)$ are constants,
with the function $\chi_{c}(y)\in C^\infty(0,+\infty)$  nonnegative and satisfying
\begin{equation}\label{chi-c}
 {\bf supp}~\chi_{c} \subset(\check{y}_{*},\hat{y}_{*}), \quad \chi_{c}(y)\equiv 1\quad  {\rm on} \quad [\check{y}^{*},\hat{y}^{*}]
\end{equation}
for some $\check{y}^{*},\hat{y}^{*}$ satisfying $0<\check{y}_{*}<\check{y}^{*}<y_*<\hat{y}^{*}<\hat{y}_{*}$,
and
$\vartheta(t,y)=\vartheta_{c}^{2}(t,y)$, with $\vartheta_{c}(t,y)=\vartheta_{o}(\theta ty)$, where $\theta$ is a positive constant to be determined later, and $\vartheta_{o}(y)\in C^\infty[0,+\infty)$ is a positive function satisfying
\begin{align}\label{theta}
\begin{cases}
\vartheta_{o}(y)=1,~ \forall y\in[0,Y],\\
\vartheta_{o}'\leq 0~\text{and}~|\vartheta_{o}^{(m)}|\leq c_{0}\vartheta_{c},~ \forall y\in [0,+\infty),\\
\vartheta_{o}\leq c_{1}e^{-\delta\y},~ \forall y\in [0,+\infty),
\end{cases}
\end{align}
for any integer $1\leq m\leq s+1$, with $c_{0},c_{1}$ being two positive constants.
The special form of the function $\vartheta(t,y)$ is designed for the decay estimate.
By the solution to the problem \eqref{AE1}, let  $\mathfrak{U}_{1}=u_{0}+u_{1}$ be the first order approximate solution to the problem \eqref{Prandtl}.

To get the second order approximate solution, let $w_{2}=u-u_0-u_{1}$, then $w_{2}$ obeys the following problem in $\Omega_{T}$,
\begin{equation}\label{E1}
\begin{cases}
\partial_{t}w_{2}
+w_{2}\partial_{x}(\mathfrak{U}_{1}+w_{2})
+\mathfrak{U}_{1}\partial_{x}w_{2}
-\partial_{y}^{-1}[\partial_{x}w_{2}]\partial_{y}(\mathfrak{U}_{1}+w_{2})
-\partial_{y}^{-1}[\partial_{x}\mathfrak{U}_{1}]\partial_{y}w_{2}\\
\hspace{3in}
=\partial_{y}^{2}w_{2}
-\epsilon_{0}^{2}\vartheta\partial_{x}^{2}\mathbb{U}_{1}
+F_{1},\\
w_{2}|_{y=0}=0,~\lim\limits_{y\to+\infty} w_{2}=0,\\
w_{2}|_{t=0}=0,
\end{cases}\tag{E$_{2}$}
\end{equation}
where  $\mathbb{U}_{1}=\mathfrak{U}_{1}-u_{0}=u_{1}$, and
\begin{align*}
F_{1}
=-\partial_{y}^{-1}[\partial_{x}u_{1}]\varsigma_{1}
 -u_{1}\partial_{x}u_{1}
 +(1-h)\partial_{y}^{-1}[\partial_{x}u_{1}]\partial_{y}u_{1}.
\end{align*}

In an argument similar from \eqref{E0} to \eqref{AE1}, instead of the problem \eqref{E1}, we study the following corrected increment problem in $\Omega_{T}$,
\begin{equation}\label{AE2}
\begin{cases}
\partial_{t}u_{2}
+\mathfrak{U}_{1}\partial_{x}u_{2}
+u_{2}\partial_{x}\mathfrak{U}_{1}
-\partial_{y}^{-1}[\partial_{x}u_{2}](\partial_{y}\mathfrak{U}_{1}+\varsigma_{2})
-\partial_{y}^{-1}[\partial_{x}\mathfrak{U}_{1}]\partial_{y}u_{2}
-h\partial_{y}^{-1}[\partial_{x}u_{2}]\partial_{y}u_{2}\\
\hspace{2.5in}
=\partial_{y}^{2}u_{2}+\epsilon_{0}^{4}\vartheta\partial_{x}^{2}u_{2}
+(\epsilon_{0}^{4}-\epsilon_{0}^{2})\vartheta\partial_{x}^{2}\mathbb{U}_{1}
+F_{1},\\
u_{2}|_{y=0}=0,~\lim\limits_{y\to+\infty} u_{2}=0,\\
u_{2}|_{t=0}=0,
\end{cases}\tag{AE$_{2}$}
\end{equation}
where $\varsigma_{2}=\bar{\varsigma}_{2}\chi_{c}(y)$ with $\bar{\varsigma}_{2}$ a constant.
Here we add the term $\epsilon_{0}^{4}\vartheta\partial_{x}^{2}\mathbb{U}_{1}$ to ensure that the parabolic problems obeyed by the approximate solutions is non-degenerate, see \eqref{An}.
The  second order approximate solution of the problem \eqref{Prandtl}  is constructed by $\mathfrak{U}_{2}=\sum\limits_{j=0}^{2}u_{j}$, and set $\mathbb{U}_{2}=\mathfrak{U}_{2}-u_{0}$.

The higher-order approximate solutions can be constructed successively. For any integer $n\geq 2$, assume that we have constructed the $n-$th order approximate solution
$\mathfrak{U}_{n}=\sum\limits_{j=0}^{n}u_{j}$
of \eqref{Prandtl},
let $w_{n+1}=u-\mathfrak{U}_{n}$, then $w_{n+1}$ obeys the following problem in $\Omega_{T}$,
\begin{equation}\label{En}
\begin{cases}
\partial_{t}w_{n+1}
+w_{n+1}\partial_{x}(\mathfrak{U}_{n}+w_{n+1})
+\mathfrak{U}_{n}\partial_{x}w_{n+1}
-\partial_{y}^{-1}[\partial_{x}w_{n+1}]\partial_{y}(\mathfrak{U}_{n}+w_{n+1})
-\partial_{y}^{-1}[\partial_{x}\mathfrak{U}_{n}]\partial_{y}w_{n+1}\\
\hspace{3.8in}
=\partial_{y}^{2}w_{n+1}
-\epsilon_{0}^{2n}\vartheta\partial_{x}^{2}\mathbb{U}_{n}
+F_{n},\\
w_{n+1}|_{y=0}=0,~\lim\limits_{y\to+\infty} w_{n+1}=0,\\
w_{n+1}|_{t=0}=0,
\end{cases}\tag{E$_{n+1}$}
\end{equation}
where  $\mathbb{U}_{n}=\mathfrak{U}_{n}-u_{0}$, and
\begin{align}\label{Fn}
F_{n}
=-\partial_{y}^{-1}[\partial_{x}u_{n}]\varsigma_{n}
 -u_{n}\partial_{x}u_{n}
 +(1-h)\partial_{y}^{-1}[\partial_{x}u_{n}]\partial_{y}u_{n}.
\end{align}

As from \eqref{E1} to \eqref{AE2},  instead of the problem \eqref{En}, we study the following corrected increment problem in $\Omega_{T}$,
\begin{equation}\label{AEn+1}
\begin{cases}
\partial_{t}u_{n+1}
+\mathfrak{U}_{n}\partial_{x}u_{n+1}+u_{n+1}\partial_{x}\mathfrak{U}_{n}
-\partial_{y}^{-1}[\partial_{x}u_{n+1}](\partial_{y}\mathfrak{U}_{n}+\varsigma_{n+1})
-\partial_{y}^{-1}[\partial_{x}\mathfrak{U}_{n}]\partial_{y}u_{n+1}
-h\partial_{y}^{-1}[\partial_{x}u_{n+1}]\partial_{y}u_{n+1}\\
\hspace{3in}
=\partial_{y}^{2}u_{n+1}+\epsilon_{0}^{2n+2}\vartheta\partial_{x}^{2}u_{n+1}
+(\epsilon_{0}^{2n+2}-\epsilon_{0}^{2n})\vartheta\partial_{x}^{2}\mathbb{U}_{n}
+F_{n},\\
u_{n+1}|_{y=0}=0,~\lim\limits_{y\to+\infty} u_{n+1}=0,\\
u_{n+1}|_{t=0}=0.
\end{cases}\tag{AE$_{n+1}$}
\end{equation}
After solving $u_{n+1}$ from the above problem, we construct the $(n+1)-$th order approximate solution
of \eqref{Prandtl} as
$\mathfrak{U}_{n+1}=\sum\limits_{j=0}^{n+1}u_{j}$.

It is worth noting that for any $n\in\mathbb{N}_{+}$, the $n-$th order approximate solution
 $\mathfrak{U}_{n}$ satisfies the following problem in $\Omega_{T}$,.
\begin{equation}\label{An}
\begin{cases}
\partial_{t}\mathfrak{U}_{n}
+\mathfrak{U}_{n}\partial_{x}\mathfrak{U}_{n}
-\partial_{y}^{-1}[\partial_{x}\mathfrak{U}_{n}]\partial_{y}\mathfrak{U}_{n}
=\partial_{y}^{2}\mathfrak{U}_{n}+\epsilon_{0}^{2n}\vartheta\partial_{x}^{2}\mathbb{U}_{n}-\partial_{x}P-F_{n},\\
\mathfrak{U}_{n}|_{y=0}=0,~\lim\limits_{y\to+\infty}\mathfrak{U}_{n}=U,\\
\mathfrak{U}_{n}|_{t=0}=u_{in}(x,y),
\end{cases}\tag{A$_n$}
\end{equation}
with $\mathbb{U}_{n}=\mathfrak{U}_{n}-u_{0}$.
Thus, if the sequence $\{\mathfrak{U}_{n}\}_{n\in\mathbb{N}}$ converges strongly in a specific energy space as $n\to\infty$, and
$$\lim_{n\to\infty}\epsilon_{0}^{2n}\vartheta\partial_{x}^{2}\mathbb{U}_{n}\to0, \qquad
\lim_{n\to\infty}F_{n}\to0,$$
then the limit $\lim\limits_{n\to\infty}\mathfrak{U}_{n}$ is a solution to the original nonlinear problem \eqref{Prandtl}.

For any $n\in\mathbb{N}_{+}$, under a certain iterative assumption on $\mathfrak{U}_{n}$ and $F_{n}$, we can establish the uniform (in $n$) a priori estimates of  the solution $u_{n+1}$ to the problem \eqref{AEn+1}, then get the existence and uniqueness of $u_{n+1}$. With a posterior estimates we show by induction that the iterative assumption also holds for $\mathfrak{U}_{n+1}$ and $F_{n+1}$, and thus the iteration can proceed step by step.

To get convergence of the sequence $\{\mathfrak{U}_{n}\}_{n\in\mathbb{N}}$ as $n\to\infty$, we show that there exist a positive constant $\mathcal{Z}$ and $t_{*}>0$ uniformly in $n$, such that for any $n\in\mathbb{N}_{+}$,
\begin{align}\label{tar}
\|u_{n}\|_{L_{t}^{\infty}\hat{H}_{\psi}^{s}(\Omega_{t_{*}})}\leq \frac{1-\epsilon_{0}}{\epsilon_{0}}\epsilon_{0}^{n}\mathcal{Z},
\end{align}
thus $\{\mathfrak{U}_{n}\}_{n\in\mathbb{N}}$ converge strongly in the space $L_{t}^{\infty}\hat{H}_{\psi}^{s}(\Omega_{t_{*}})$ as $n\to\infty$.
The estimate \eqref{tar} is achievable, since $u_{n}|_{t=0}=0$ for $n\geq 2$ and all terms in $F_{n-1}$ are nonlinear in $u_{n-1}$ except
$\partial_{y}^{-1}[\partial_{x}u_{n-1}]\varsigma_{n-1}$, $\|F_{n-1}\|_{L_{t}^{\infty}\hat{H}_{\psi}^{s}(\Omega_{t_{*}})}$ can be of order $\epsilon_{0}^{n}$ from the induction estimate of $u_{n-1}$ and choosing $\varsigma_{n-1}$ small. Thus, by letting $t_{*}$ small, the estimate \eqref{tar} can be fulfilled from that of $u_{n-1}$.

\section{Well-posedness of the corrected increment problems }

In this section, we study the well-posedness of the corrected increment  problems (AE$_{n}$) for any $n\in\mathbb{N}_{+}$.
Note that for different $n$, the structures of corrected increment problems (AE$_{n}$) and the approximate problems \eqref{An} remain the same, so it suffices to study the following problem in $\Omega_{T}$,
\begin{equation}\label{app}
\begin{cases}
\partial_{t}u
+\mathfrak{U}\partial_{x}u+u\partial_{x}\mathfrak{U}
-\partial_{y}^{-1}[\partial_{x}u](\partial_{y}\mathfrak{U}+\varsigma)
-\partial_{y}^{-1}[\partial_{x}\mathfrak{U}]\partial_{y}u
-h\partial_{y}^{-1}[\partial_{x}u]\partial_{y}u\\
\quad
=\partial_{y}^{2}u+\eta\partial_{x}^{2}u+\kappa\partial_{x}^{2}\mathbb{U}+F,\\
u|_{y=0}=0,~\lim\limits_{y\to+\infty} u=0,\\
u|_{t=0}=\check{u}(x,y).
\end{cases}
\end{equation}
where $\check{u}(x,y)$ and $F(t,x,y)$ are two known functions, $h(y)$ is the cutoff function given in \eqref{h},  $\varsigma(y)=\bar{\varsigma}\chi_{c}(y)$ with $\chi_c$ given in \eqref{chi-c}, $\eta(t,y)=\bar{\eta}\vartheta_{c}^{2}(t,y)$ with $\vartheta_{c}(t,y)$ given above \eqref{theta}, and $\kappa(t,y)=\eta(t,y)-\epsilon_{0}^{-2}\eta(t,y)$,
with $\bar{\varsigma}, \bar{\eta}>0$ two constants,  $\mathbb{U}=\mathfrak{U}-u_{0}$. It should be noted that when $n=1$, $\check{u}=\tilde{u}_{in}, \mathfrak{U}(t,x,y)=u_{0}$, while for $n\geq 2$, $\check{u}=0$ and $\mathfrak{U}(t,x,y)$ is the solution to the following problem,
\begin{equation}\label{App}
\begin{cases}
\partial_{t}\mathfrak{U}
+\mathfrak{U}\partial_{x}\mathfrak{U}
-\partial_{y}^{-1}[\partial_{x}\mathfrak{U}]\partial_{y}\mathfrak{U}
=\partial_{y}^{2}\mathfrak{U}+\epsilon_{0}^{-2}\eta\partial_{x}^{2}\mathbb{U}-\partial_{x}P-F,\\
\mathfrak{U}|_{y=0}=0,~\lim\limits_{y\to+\infty}\mathfrak{U}=U(t,x),\\
\mathfrak{U}|_{t=0}=u_{in}(x,y).
\end{cases}
\end{equation}

For convenience,  let $\bar{\kappa}=\epsilon_{0}^{-2}\bar{\eta}-\bar{\eta}$,
and $\bar{\varsigma}=\iota\bar{\eta}^{2}$, that is $\iota=\bar{\varsigma}\bar{\eta}^{-2}$.
To use the cancellation property of \eqref{app} to study the well-posedness in a Sobolev space, for all $(t,x,y)\in\Omega_{T}$ we let
\begin{align}\label{rhoPd}
\rho(t,x,y)=\frac{\partial_{y}\dot{\mathfrak{U}}+\varsigma'}{\dot{\mathfrak{U}}+\varsigma}, \quad
\mathcal{P}(t,x,y)=\frac{\chi(y)}{\sqrt{\dot{\mathfrak{U}}+\varsigma}}+1-\chi(y),
\end{align}
where $\dot{\mathfrak{U}}=\partial_{y}\mathfrak{U}$, and $\chi(y)$ is given in \eqref{chi}.

For $r=s$ or $s+1$, define
$$\mathcal{U}_{r}(t,x,y)=\mathcal{P}(\partial_{y}\partial_{x}^{r}u-\rho\partial_{x}^{r}u),$$
and
$$\mathbf{U}_{r}(t,x,y)=\mathcal{P}\left(\partial_{y}\partial_{x}^{r}\mathbb{U}-\rho\partial_{x}^{r}\mathbb{U}\right). $$

\subsection{Construction of a weight function}

To get a priori energy estimates, we need an appropriate weight, which is bounded and provides a good control of $\mathcal{P}$. As mentioned in Section 3, if $\mathfrak{U}$ is a classical solution to the problem \eqref{App} in $\Omega_{T}$ under the assumption in Theorem \ref{MR}, then $\dot{\mathfrak{U}}(t,x,y)\geq C((y-y_{*})^{2}+t)$ in $\Omega_{T}^{*}$ for some constant $C$. Motivated by this observation, we are going to find a bounded weight function $\varphi(t,y)$ satisfying
\begin{align}\label{tage}
\partial_{t}\varphi^{2}+\partial_{y}^{2}\varphi^{2}\leq-\frac{C_{L}}{(y-y_{*})^{2}+(\bar{\varsigma}+t)}\varphi^{2},\quad
\forall (t,y)\in [0,T]\times[\check{y}_{*},\hat{y}_{*}],
\end{align}
where $C_{L}$ is a constant independent of $\bar{\varsigma}$.
It needs careful analysis as $(t,y)$ close to $(0,y_{*})$.
Let $\zeta_{\lambda}^{oc}(t,y)$ be given by
\begin{align*}
\zeta_{\lambda}^{oc}(t,y)=\exp\left({-\frac{16\lambda^{4}(\bar{\varsigma}+t)}{(y-y_{*})^{2}+16\lambda^{2}(\bar{\varsigma}+t)}}\right)
\end{align*}
where $\lambda$ is a large even number to be determined later. By a direct calculation, one can easily verify  that the function $\varphi=(\zeta_{\lambda}^{oc})^{\frac 12}$ satisfies the inequality \eqref{tage} when $\lambda$ is large and $|y-y_{*}|\geq \lambda\sqrt{\bar{\varsigma}+t}$.

On the other hand, when $|y-y_{*}|< \lambda\sqrt{\bar{\varsigma}+t}$, define  $\zeta_{\lambda}^{c}(t,y)$ as
\begin{align*}
\zeta_{\lambda}^{c}(t,y)=\zeta(z),
\end{align*}
with
\begin{align*}
z=\frac{y-y_{*}}{\sqrt{\bar{\varsigma}+t}}+\frac{1}{(\bar{\varsigma}+t)^{\varepsilon_{1}}}, \quad
\varepsilon_{1}=\frac{\varepsilon_{0}}{\lambda}.
\end{align*}

By a direct calculation, it gives that
\begin{align}\label{zetae}
\partial_{t}\zeta_\lambda^{c}+\partial_{y}^{2}\zeta_\lambda^{c}
=\frac{1}{(\bar{\varsigma}+t)}\left[\zeta''-\frac{1}{2}\left(z-\frac{1-2\varepsilon_{1}}{(\bar{\varsigma}+t)^{\varepsilon_{1}}}\right)\zeta'\right].
\end{align}
 In order for the function $\sqrt{\zeta^{c}_\lambda}(t,y)$ to satisfy \eqref{tage}, one considers the following equation,
\begin{align}\label{hermite}
\zeta''-\frac{1}{2}(z-z_{0})\zeta'+\frac{1}{2}\lambda\zeta=0,
\end{align}
with $z_{0}=\frac{1-2\varepsilon_{1}}{(\bar{\varsigma}+t)^{\varepsilon_{1}}}$. The  equation \eqref{hermite} admits a solution in the form of $\zeta(z)=H_{\lambda}(\frac{z-z_{0}}{2})$, where $H_{\lambda}(\cdot)$ is the following $\lambda$-th Hermite polynomial(e.g. see \cite{Sz}),
$$
H_{\lambda}(x)=\lambda!\sum_{i=0}^{\frac{\lambda}{2}}\frac{(-1)^{i}}{i!(\lambda-2i)!}(2x)^{\lambda-2i}.
$$
It is known that all zeros of the polynomial $H_{\lambda}(\cdot)$ lie in the interval $(-\sqrt{2\lambda+1},\sqrt{2\lambda+1})$, thus $\zeta(z)$ is positive as long as \begin{align*}
z-z_{0}=\frac{y-y_{*}}{\sqrt{\bar{\varsigma}+t}}+\frac{2\varepsilon_{1}}{(\bar{\varsigma}+t)^{\varepsilon_{1}}}>2\sqrt{2\lambda+1}.
\end{align*}
For any fixed $\lambda$, this can be achieved for any $y\in(y_{*}-\lambda\sqrt{\bar{\varsigma}+t},y_{*}+\lambda\sqrt{\bar{\varsigma}+t})$ by controlling $\bar{\varsigma}$ and $t$ to be sufficiently small. Thus, we find a suitable function
\begin{align}\label{zetac}
\zeta_{\lambda}^{c}(t,y)
=H_{\lambda}\left(\frac{y-y_{*}}{2\sqrt{(\bar{\varsigma}+t)}}+\frac{\varepsilon_{1}}{(\bar{\varsigma}+t)^{\varepsilon_{1}}}\right).
\end{align}
such that $\sqrt{\zeta_{\lambda}^{c}}(t,y)$ satisfies \eqref{tage}
when $|y-y_{*}|< \lambda\sqrt{\bar{\varsigma}+t}$.

To balance the degree of the weight function with respect to the factor $\bar{\varsigma}+t$ inside and outside the interval $(y_{*}-\lambda\sqrt{\bar{\varsigma}+t},y_{*}+\lambda\sqrt{\bar{\varsigma}+t})$, we let
\begin{align}\label{zetao}
\zeta_{\lambda}^{o}(t,y)=\frac{\zeta_{\lambda}^{oc}(t,y)}{[\frac{\lambda^{-2}}{16}(y-y_{*})^{2}+\bar{\varsigma}+t]^{\varepsilon_{0}}}.
\end{align}
Based on the above analysis, we take the weight function as
\begin{equation}\label{varphi}
\varphi(t,y)=\y^{-\frac{1}{2}}\psi(t,y)\sqrt{\zeta_{\lambda}(t,y)},
\end{equation}
with $\psi(t,y)$ given in Notation 2.1, and
\begin{align*}
\zeta_{\lambda}(t,y)
=\chi(y)\left[\zeta_{\lambda}^{o}(t,y)+\left(\frac{\lambda}{\varepsilon_{0}}\right)^{\lambda}e^{-\lambda^{2}}\zeta_{\lambda}^{c}(t,y)\chi_{\lambda}\left(\frac{y-y_{*}}{\sqrt{\bar{\varsigma}+t}}\right)\right]
 +1-\chi(y).
\end{align*}
Here $\chi(y)$ is given in \eqref{chi} of Section 2, and $\chi_{\lambda}(y)\in C^{\infty}(\mathbb{R})$ is a nonnegative function such that
\begin{align*}
\chi_{\lambda}(y)=
\begin{cases}
1,~\forall y\in[-\lambda+1,\lambda-1],\\
0,~\forall y\in\mathbb{R}\setminus[-\lambda,\lambda].\\
\end{cases}
\end{align*}
The factor $\y^{-\frac{1}{2}}$ in $\varphi$ is introduced to handle the estimation of nonlocal terms.

By choosing $\lambda$ large enough and $\bar{\varsigma}, t$ appropriately small, it can be proved that the weight function $\varphi$ defined in \eqref{varphi} satisfies the inequality \eqref{tage}. In this regard, we have the following lemma.

\begin{lemma}\label{weighte}
Let $\varphi$ be defined as in \eqref{varphi}, with $\lambda$ a even positive integer. If $\lambda$ is large enough, there exists positive constants $\Lambda_{\lambda}$ and $\iota_{\lambda}$ depending only on $\lambda$ such that for any
$\Lambda>\Lambda_{\lambda}, \iota\in(0,\iota_{\lambda})$ one can find a constant $T_{\lambda,\Lambda}\in(0,T)$ depending only on $\lambda$ and $\Lambda$ to obtain the following inequalities in $[0,T_{\lambda,\Lambda}]\times[0,+\infty)$,
\begin{align}\label{We}
\partial_{t}\varphi^{2}+\partial_{y}^{2}\varphi^{2}
\leq -\frac{\lambda}{32}\omega_{\lambda}\varphi^{2}
     -\frac{3}{2}\Lambda\delta\y\varphi^{2},
\end{align}
moreover,
\begin{align}\label{varphiy}
|\partial_{y}\varphi^{2}|
\leq C_{l}\varphi^{2}+C_{l}\lambda\sqrt{\omega_{\lambda}}\varphi^{2},
\end{align}
where $C_{l}$ is a positive constant independent of $\lambda,\iota$ and $\Lambda$, and
$$\omega_{\lambda}(t,y)=\frac{1}{\frac{\lambda^{-2}}{16}(y-y_{*})^{2}+\bar{\varsigma}+t}.$$
\end{lemma}

\begin{proof}
For any  $(t,y)\in[0,T]\times\mathbb{R}$, to simplify the notation, we let $z=\frac{1}{4\lambda}(y-y_{*}), \tau=\bar{\varsigma}+t$, and let
\begin{align*}
\zeta_{\lambda}^{*}(t,y)
= \zeta_{\lambda}^{o}(t,y)
 +\bar{\zeta}_{\lambda}^{c}(t,y)\bar{\chi}_{\lambda}(t,y),
\end{align*}
with
\begin{align*}
\bar{\chi}_{\lambda}(t,y)=\chi_{\lambda}\left(\frac{y-y_{*}}{\sqrt{\bar{\varsigma}+t}}\right),~
\bar{\zeta}_{\lambda}^{c}(t,y)=\left(\frac{\lambda}{\varepsilon_{0}}\right)^{\lambda}e^{-\lambda^{2}}\zeta_{\lambda}^{c}(t,y).
\end{align*}
By a direct calculation, it follows that
\begin{align*}
\partial_{t}\zeta_{\lambda}^{o}+\partial_{y}^{2}\zeta_{\lambda}^{o}
=&~\left[-\frac{\lambda^{2} z^{2}}{(z^{2}+\tau)^{2}}-\frac{z^{2}\tau}{2(z^{2}+\tau)^{3}}
         +\frac{\tau}{8(z^{2}+\tau)^{2}}+\frac{\lambda^{2}\tau^{2}z^{2}}{4(z^{2}+\tau)^{4}}\right]\zeta_{\lambda}^{o}\\
 &+\left[-\frac{\varepsilon_{0}}{(z^{2}+\tau)}-\frac{\varepsilon_{0}}{8\lambda^{2}(z^{2}+\tau)}
         +\frac{\varepsilon_{0}(1+\varepsilon_{0})z^{2}}{4\lambda^{2}(z^{2}+\tau)^{2}}\right]\zeta_{\lambda}^{o}\\
\leq&~\left[-\frac{3\lambda^{2} z^{2}}{4(z^{2}+\tau)^{2}}+\frac{\tau^{3}}{8(z^{2}+\tau)^{4}}\right]\zeta_{\lambda}^{o},
\end{align*}
Thus for any $y\in\mathbb{R}\setminus(y_{*}-4\sqrt{\bar{\varsigma}+t},y_{*}+4\sqrt{\bar{\varsigma}+t})$, one has
\begin{align}\label{zetaoe-o-0}
\partial_{t}\zeta_{\lambda}^{o}+\partial_{y}^{2}\zeta_{\lambda}^{o}
\leq -\frac{5(y-y_{*})^{2}}{128(\frac{\lambda^{-2}}{16}(y-y_{*})^{2}+\bar{\varsigma}+t)^{2}}\zeta_{\lambda}^{o},
\end{align}
and for any $y\in\left(y_{*}-4\sqrt{\bar{\varsigma}+t},y_{*}+4\sqrt{\bar{\varsigma}+t}\right)$, it follows that
\begin{align}\label{zetaoe-i}
\partial_{t}\zeta_{\lambda}^{o}+\partial_{y}^{2}\zeta_{\lambda}^{o}
\leq \frac{1}{8(\bar{\varsigma}+t)}\zeta_{\lambda}^{o}.
\end{align}

From \eqref{zetae} and \eqref{hermite} one knows that for any $y\in(y_{*}-\lambda\sqrt{\bar{\varsigma}+t},y_{*}+\lambda\sqrt{\bar{\varsigma}+t})$,
\begin{align}\label{zetace}
\partial_{t}\bar{\zeta}_{\lambda}^{c}+\partial_{y}^{2}\bar{\zeta}_{\lambda}^{c}
=-\frac{\lambda}{2(\bar{\varsigma}+t)}\bar{\zeta}_{\lambda}^{c}.
\end{align}
moreover,
\begin{align}\label{zetacb}
\left(\frac{3}{2}\right)^{\lambda}\frac{e^{-\lambda^{2}}}{{(\bar{\varsigma}+t)^{\varepsilon_{0}}}}
\leq \bar{\zeta}_{\lambda}^{c}(t,y)\leq
\left(\frac{5}{2}\right)^{\lambda}\frac{e^{-\lambda^{2}}}{{(\bar{\varsigma}+t)^{\varepsilon_{0}}}},
\end{align}
provided $\iota+t\leq \mathcal{C}_{\lambda}$, for some sufficiently small positive constant $\mathcal{C}_{\lambda}$ which depends only on $\lambda$.

One can check that if $\lambda$ is large enough, then for any $y\in\left(y_{*}-\sqrt{\lambda}\sqrt{\bar{\varsigma}+t},y_{*}+\sqrt{\lambda}\sqrt{\bar{\varsigma}+t}\right)$,
\begin{align}\label{olc}
\zeta_{\lambda}^{o}(t,y)\leq\bar{\zeta}_{\lambda}^{c}(t,y).
\end{align}
Thus for $y\in\left(y_{*}-\sqrt{\lambda}\sqrt{\bar{\varsigma}+t},y_{*}+\sqrt{\lambda}\sqrt{\bar{\varsigma}+t}\right)$,
combining \eqref{zetaoe-o-0}, \eqref{zetaoe-i} and \eqref{zetace} with \eqref{olc} one has
\begin{align}\label{zeta*1}
\partial_{t}\zeta_{\lambda}^{*}+\partial_{y}^{2}\zeta_{\lambda}^{*}
=&~\left(\partial_{t}+\partial_{y}^{2}\right)(\zeta_{\lambda}^{o}+\bar{\zeta}_{\lambda}^{c})\\
\leq& \frac{1}{8(\bar{\varsigma}+t)}\zeta_{\lambda}^{o}-\frac{\lambda}{2(\bar{\varsigma}+t)}\bar{\zeta}_{\lambda}^{c}\nonumber\\
\leq& -\frac{\lambda}{4(\frac{\lambda^{-2}}{16}(y-y_{*})^{2}+\bar{\varsigma}+t)}\zeta_{\lambda}^{*}.\nonumber
\end{align}

For $y\in\left(y_{*}-(\lambda-1)\sqrt{\bar{\varsigma}+t},y_{*}+(\lambda-1)\sqrt{\bar{\varsigma}+t}\right)\setminus(y_{*}-\sqrt{\lambda}\sqrt{\bar{\varsigma}+t},y_{*}+\sqrt{\lambda}\sqrt{\bar{\varsigma}+t})$,
from \eqref{zetaoe-o-0} one can deduce that
\begin{align}\label{zetaoe-o-1}
\partial_{t}\zeta_{\lambda}^{o}+\partial_{y}^{2}\zeta_{\lambda}^{o}
\leq -\frac{9\lambda}{256(\frac{\lambda^{-2}}{16}(y-y_{*})^{2}+\bar{\varsigma}+t)}\zeta_{\lambda}^{o},
\end{align}
by choosing $\lambda$ large enough. Combining \eqref{zetace} wih \eqref{zetaoe-o-1} one gets
\begin{align}\label{zeta*2}
\partial_{t}\zeta_{\lambda}^{*}+\partial_{y}^{2}\zeta_{\lambda}^{*}
=&~\left(\partial_{t}+\partial_{y}^{2}\right)(\zeta_{\lambda}^{o}+\bar{\zeta}_{\lambda}^{c})\\
\leq& -\frac{9\lambda}{256(\frac{\lambda^{-2}}{16}(y-y_{*})^{2}+\bar{\varsigma}+t)}\zeta_{\lambda}^{o}
      -\frac{\lambda}{2(\bar{\varsigma}+t)}\bar{\zeta}_{\lambda}^{c}\nonumber\\
\leq& -\frac{9\lambda}{256(\frac{\lambda^{-2}}{16}(y-y_{*})^{2}+\bar{\varsigma}+t)}\zeta_{\lambda}^{*}.\nonumber
\end{align}

For $y\in[y_{*}-\lambda\sqrt{\bar{\varsigma}+t},y_{*}-(\lambda-1)\sqrt{\bar{\varsigma}+t}]\cup[y_{*}+(\lambda-1)\sqrt{\bar{\varsigma}+t},y_{*}+\lambda\sqrt{\bar{\varsigma}+t}]$,
by the property of the Hermite polynomial $H_{\lambda}'(x)=2\lambda H_{\lambda-1}(x)$, one has
\begin{align}\label{chie}
\bar{\zeta}_{\lambda}^{c}\left(\partial_{t}+\partial_{y}^{2}\right)\bar{\chi}_{\lambda}+2\partial_{y}\bar{\zeta}_{\lambda}^{c}\partial_{y}\bar{\chi}_{\lambda}
=&\left[-\frac{y-y_{*}}{2\sqrt{\bar{\varsigma}+t}}\frac{\chi_{\lambda}'\left(\frac{y-y_{*}}{\sqrt{\bar{\varsigma}+t}}\right)}{\bar{\varsigma}+t}
  +\frac{1}{\bar{\varsigma}+t}\frac{\chi_{\lambda}''\left(\frac{y-y_{*}}{\sqrt{\bar{\varsigma}+t}}\right)}{\bar{\varsigma}+t}\right]\bar{\zeta}_{\lambda}^{c}\\
 &+2\lambda\left(\frac{\lambda}{\varepsilon_{0}}\right)^{\lambda}e^{-\lambda^{2}}\frac{\chi_{\lambda}'\left(\frac{y-y_{*}}{\sqrt{\bar{\varsigma}+t}}\right)}{(\bar{\varsigma}+t)} H_{\lambda-1}\left(\frac{y-y_{*}}{2\sqrt{(\bar{\varsigma}+t)}}+\frac{\varepsilon_{1}}{(\bar{\varsigma}+t)^{\varepsilon_{1}}}\right)\nonumber\\
\leq&~\frac{C\lambda}{\bar{\varsigma}+t}\bar{\zeta}_{\lambda}^{c}.\nonumber
\end{align}
As a consequence of \eqref{zetaoe-o-0}, \eqref{zetace}, \eqref{chie} and \eqref{zetacb}, we have
\begin{align}\label{zeta*3}
\partial_{t}\zeta_{\lambda}^{*}+\partial_{y}^{2}\zeta_{\lambda}^{*}
=&~\left(\partial_{t}+\partial_{y}^{2}\right)\zeta_{\lambda}^{o}
  +\bar{\chi}_{\lambda}\left(\partial_{t}+\partial_{y}^{2}\right)\bar{\zeta}_{\lambda}^{c}
  +\bar{\zeta}_{\lambda}^{c}\left(\partial_{t}+\partial_{y}^{2}\right)\bar{\chi}_{\lambda}
  +2\partial_{y}\bar{\zeta}_{\lambda}^{c}\partial_{y}\bar{\chi}_{\lambda}\\
\leq&-\frac{5(y-y_{*})^{2}}{128(\frac{\lambda^{-2}}{16}(y-y_{*})^{2}+\bar{\varsigma}+t)^{2}}\zeta_{\lambda}^{o}
     +\frac{C\lambda}{\bar{\varsigma}+t}\bar{\zeta}_{\lambda}^{c}\nonumber\\
\leq&-\frac{9(\lambda-1)^{2}}{256(\frac{\lambda^{-2}}{16}(y-y_{*})^{2}+\bar{\varsigma}+t)^{2}}\zeta_{\lambda}^{o}\nonumber\\
\leq&-\frac{9\lambda}{256(\frac{\lambda^{-2}}{16}(y-y_{*})^{2}+\bar{\varsigma}+t)}\zeta_{\lambda}^{*}.\nonumber
\end{align}
by choosing $\lambda$ large enough.

Combining \eqref{zeta*1}, \eqref{zeta*2}, \eqref{zeta*3} and \eqref{zetaoe-o-0} yields
\begin{align}\label{zeta*ty}
\partial_{t}\zeta_{\lambda}^{*}+\partial_{y}^{2}\zeta_{\lambda}^{*}
\leq -\frac{9\lambda}{256(\frac{\lambda^{-2}}{16}(y-y_{*})^{2}+\bar{\varsigma}+t)}\zeta_{\lambda}^{*},~
\forall (t,y)\in [0,T]\times\mathbb{R},
\end{align}
provided $\lambda$ large enough and $\mathcal{C}_{\lambda}$  sufficiently small.

Recall the definition of $\varphi$, it is easy to obtain that
\begin{align}\label{varphity}
\partial_{t}\varphi^{2}+\partial_{y}^{2}\varphi^{2}
=&~\left(\partial_{t}\zeta_{\lambda}^{*}+\partial_{y}^{2}\zeta_{\lambda}^{*}\right)\y^{-1}\psi^{2}\chi
  -2\Lambda\delta\y\varphi^{2}
  +a_{1}\y^{-1}\psi^{2}\partial_{y}\zeta_{\lambda}^{*}\\
 &+a_{2}\y^{-1}\psi^{2}\zeta_{\lambda}^{*}
  +a_{3}\y^{-1}\psi^{2},\nonumber
\end{align}
where
$$a_{1}=-2\y^{-1}\chi+\left(3-4\Lambda t\right)\delta\chi+2\chi',$$
$$a_{2}=2\y^{-2}\chi+\delta^{2}\left(\frac{3}{2}-2\Lambda t\right)^{2}\chi
         +\chi''-\y^{-1}\delta\left(3-4\Lambda t\right)\chi
         -\y^{-1}\chi'
         +\delta\left(3-4\Lambda t\right)\chi',$$
and
\begin{align*}
a_{3}=&~2\y^{-2}(1-\chi)-\chi''+\delta^{2}\left(\frac{3}{2}-2\Lambda t\right)^{2}(1-\chi)+2\y^{-1}\chi'\\
      &-\delta\y^{-1}\left(3-4\Lambda t\right)(1-\chi)
       -\delta\left(3-4\Lambda t\right)\chi'.
\end{align*}
A direct calculation gives
\begin{align*}
\partial_{y}\zeta_{\lambda}^{*}
=&~\partial_{y}\zeta_{\lambda}^{o}
 +\partial_{y}\bar{\zeta}_{\lambda}^{c}(t,y)\bar{\chi}_{\lambda}(t,y)
 +\bar{\zeta}_{\lambda}^{c}(t,y)\partial_{y}\bar{\chi}_{\lambda}(t,y)\\
=&~\left[\frac{(y-y_{*})(\bar{\varsigma}+t)}{8(\frac{\lambda^{-2}}{16}(y-y_{*})^{2}+\bar{\varsigma}+t)^{2}}
         -\frac{\varepsilon_{0}(y-y_{*})}{8\lambda^{2}(\frac{\lambda^{-2}}{16}(y-y_{*})^{2}+\bar{\varsigma}+t)}\right]\zeta_{\lambda}^{o}\\
 &+\lambda\left(\frac{\lambda}{\varepsilon_{0}}\right)^{\lambda}e^{-\lambda^{2}}\frac{\chi_{\lambda}\left(\frac{y-y_{*}}{\sqrt{\bar{\varsigma}+t}}\right)}{\sqrt{\bar{\varsigma}+t}}
   H_{\lambda-1}\left(\frac{y-y_{*}}{2\sqrt{(\bar{\varsigma}+t)}}+\frac{\varepsilon_{1}}{(\bar{\varsigma}+t)^{\varepsilon_{1}}}\right)\nonumber\\
 &+\frac{\chi_{\lambda}'\left(\frac{y-y_{*}}{\sqrt{\bar{\varsigma}+t}}\right)}{\sqrt{\bar{\varsigma}+t}}\bar{\zeta}_{\lambda}^{c}(t,y).\nonumber
\end{align*}
thus, by taking $\lambda$ large enough and $\mathcal{C}_{\lambda}$ sufficiently small, one has
\begin{align}\label{zeta*y}
|\partial_{y}\zeta_{\lambda}^{*}|
\leq&~\frac{\lambda}{\sqrt{\frac{\lambda^{-2}}{16}(y-y_{*})^{2}+\bar{\varsigma}+t}}\zeta_{\lambda}^{o}
     +C\frac{\chi_{\lambda}\left(\frac{y-y_{*}}{\sqrt{\bar{\varsigma}+t}}\right)+\chi_{\lambda}'\left(\frac{y-y_{*}}{\sqrt{\bar{\varsigma}+t}}\right)}{\sqrt{\frac{\lambda^{-2}}{16}(y-y_{*})^{2}+\bar{\varsigma}+t}}\bar{\zeta}_{\lambda}^{c}(t,y)\\
\leq&\frac{C\lambda}{\sqrt{\frac{\lambda^{-2}}{16}(y-y_{*})^{2}+\bar{\varsigma}+t}}\zeta_{\lambda}^{*}.\nonumber
\end{align}
Substituting \eqref{zeta*ty} and \eqref{zeta*y} into \eqref{varphity}, it follows that
\begin{align}
\partial_{t}\varphi^{2}+\partial_{y}^{2}\varphi^{2}
\leq&-\frac{9\lambda\y^{-1}\psi^{2}\chi\zeta_{\lambda}^{*}}{256(\frac{\lambda^{-2}}{16}(y-y_{*})^{2}+\bar{\varsigma}+t)}
  -2\Lambda\delta\y\varphi^{2}
  +\frac{C\lambda\y^{-1}\psi^{2}\zeta_{\lambda}^{*}}{\sqrt{\frac{\lambda^{-2}}{16}(y-y_{*})^{2}+\bar{\varsigma}+t}}\\
 &+C\y^{-1}\psi^{2}\zeta_{\lambda}^{*}
  +C(1-\chi+\chi')\y^{-1}\psi^{2}\nonumber\\
\leq&-\frac{\lambda\y^{-1}\psi^{2}\chi\zeta_{\lambda}^{*}}{32(\frac{\lambda^{-2}}{16}(y-y_{*})^{2}+\bar{\varsigma}+t)}
     -2\Lambda\delta\y\varphi^{2}
     +C\lambda\y^{-1}\psi^{2}\zeta_{\lambda}^{*}\nonumber\\
    &+C(1-\chi+\chi')\y^{-1}\psi^{2}\nonumber\\
\leq&-\frac{\lambda}{32(\frac{\lambda^{-2}}{16}(y-y_{*})^{2}+\bar{\varsigma}+t)}\varphi^{2}
     -2\Lambda\delta\y\varphi^{2}
     +C\lambda\varphi^{2}\nonumber\\
\leq&-\frac{\lambda}{32(\frac{\lambda^{-2}}{16}(y-y_{*})^{2}+\bar{\varsigma}+t)}\varphi^{2}
     -\frac{3}{2}\Lambda\delta\y\varphi^{2},\nonumber
\end{align}
provided $\Lambda$ large enough and $t\leq\frac{3}{4\Lambda}$.

Since
\begin{align*}
\partial_{y}\varphi^{2}
=\left[\delta\left(\frac{3}{2}-2\Lambda t\right)-\y^{-1}\right]\varphi^{2}
 +\y^{-1}\psi^{2}(\chi'\zeta_{\lambda}^{*}-\chi')+\y^{-1}\psi^{2}\chi\partial_{y}\zeta_{\lambda}^{*},
\end{align*}
one can use \eqref{zeta*y} to obtain
\begin{align*}
|\partial_{y}\varphi^{2}|
\leq C\varphi^{2}
    +\frac{C\lambda\y^{-1}\psi^{2}\chi\zeta_{\lambda}^{*}}{\sqrt{\frac{\lambda^{-2}}{16}(y-y_{*})^{2}+\bar{\varsigma}+t}}
\leq C\varphi^{2}+\frac{C\lambda}{\sqrt{\frac{\lambda^{-2}}{16}(y-y_{*})^{2}+\bar{\varsigma}+t}}\varphi^{2}.
\end{align*}
Thus \eqref{varphiy} holds.
\end{proof}

\subsection{Well-posed result of the problem \eqref{app}}

With the weights above, for any function $f(x,y)$ defined in $\Omega$ and $s\in\mathbb{N}_{+}$, we introduce the norm
\begin{align*}
\|f\|_{\mathcal{H}_{\psi,\varphi}^{s}(\Omega)}
= \|f\|_{\hat{H}_{\psi}^{s}(\Omega)}
 +\|\mathcal{P}(\partial_{y}\partial_{x}^{s}f-\rho\partial_{x}^{s}f)\|_{L_{\varphi}^{2}(\Omega)},
\end{align*}
with the norm $\|\cdot\|_{\hat{H}_{\psi}^{s}(\Omega)}$ given in \eqref{cn}, and correspondingly, the function space
\begin{align*}
\mathcal{H}_{\psi,\varphi}^{s}(\Omega)=\{f:\Omega\to\mathbb{R}|\|f\|_{\mathcal{H}_{\psi,\varphi}^{s}(\Omega)}<+\infty\}.
\end{align*}
For convenience of notation, for any nonnegative function $\xi(y)\in C^{\infty}(\mathbb{R}_{+})$, we introduce
\begin{align*}
\|f\|_{\hat{H}_{\psi,\xi}^{s}(\Omega)}
= \sum_{\substack{|\alpha|\leq s,\\ \alpha_{1}<s}}\|\sqrt{\xi}\partial_{x}^{\alpha_{1}}\partial_{y}^{\alpha_{2}}f\|_{L_{\psi}^{2}(\Omega)}
 +\|\sqrt{\xi}\partial_{x}^{s-1}\partial_{y}^{2}f\|_{L_{\hat{\psi}}^{2}(\Omega)},
\end{align*}
For any nonnegative integer $p$, define
\begin{align*}
\|f\|_{\mathring{H}^{p}(\mathbb{T})}=\sum_{0\leq i+j\leq p}\|\partial_{x}^{i}\partial_{y}^{j}f|_{y=0}\|_{L_{x}^{2}(\mathbb{T})}.
\end{align*}
and for any positive integer $q$, define
\begin{align*}
\|f\|_{\ddot{H}^{q}(\Omega)}=\|f\|_{H^{q-1}(\Omega)}+\|\partial_{y}f\|_{H^{q-1}(\Omega)}.
\end{align*}

Let $\hat{\varphi}=\y^{-\frac{1}{2}}\varphi$. To get the existence of the solution to the problem \eqref{app}, the main assumption we need is as follows.
\begin{assumption}\label{MA}
Let $\rho$ and $\mathcal{P}$ be given in \eqref{rhoPd} and $\mathbb{U}=\mathfrak{U}-u_{0}$. Assume that there exist $\hat{T}_{*}\in(0,T]$  such that the following inequalities hold,
$$\dot{\mathfrak{U}}\geq C\varpi(U-\mathfrak{U})\geq C\varpi e^{-(1+\theta t)\delta\y}~\text{in}~\Omega_{\hat{T}_{*}},$$
$$\dot{\mathfrak{U}}\geq C(\varpi+t)~\text{in}~\Omega_{\hat{T}_{*}}^{*},$$
$$|\partial_{x}\dot{\mathfrak{U}}|\leq C(1+t\y)\dot{\mathfrak{U}},
  |\partial_{x}^{2}\dot{\mathfrak{U}}|\leq C(1+t\y)\dot{\mathfrak{U}},
  |\partial_{y}\partial_{x}^{2}\dot{\mathfrak{U}}|\leq C(1+t\y)\dot{\mathfrak{U}}^{\frac{1}{2}}~\text{in}~\Omega_{\hat{T}_{*}},$$
$$|\vartheta_{c}^{2}\partial_{x}^{2}\mathbb{U}|\leq C_{*}\y^{-1}\dot{\mathfrak{U}}~\text{in}~\check{\Omega}_{\hat{T}_{*}}^{\hat{y}},$$
\begin{align*}
& \|\mathbb{U}\|_{L_{t}^{\infty}\mathcal{H}_{\psi,\varphi}^{s}(\Omega_{\hat{T}_{*}})}
 +\|\partial_{y}\mathbb{U}\|_{L_{t}^{2}\hat{H}_{\psi}^{s}(\Omega_{\hat{T}_{*}})}
 +\|\partial_{y}\mathbf{U}_{s}\|_{L_{t}^{2}L_{\varphi}^{2}(\Omega_{\hat{T}_{*}})}
 +\|\mathcal{P}\mathbf{U}_{s}\|_{L_{t}^{2}L_{\varphi}^{2}(\Omega_{\hat{T}_{*}})}\\
&+\|\mathbf{U}_{s+1}\|_{L_{t}^{2}L_{\hat{\varphi}}^{2}(\Omega_{\hat{T}_{*}})}
 +\|\mathcal{P}\mathbf{U}_{s+1}\|_{L_{t}^{2}L_{\hat{\varphi}}^{2}(\Omega_{\hat{T}_{*}})}
 +\epsilon^{-2}\|\sqrt{\eta}\partial_{x}\mathbf{U}_{s+1}\|_{L_{t}^{2}L_{\hat{\varphi}}^{2}(\Omega_{\hat{T}_{*}})}
\leq \mathcal{Z},
\end{align*}
and
$$\rho\leq -\frac{15}{16}\delta,~\text{in}~\check{\Omega}_{T_{*}}^{Y},$$
\begin{align*}
|\rho|\leq C\mathcal{P},~|\partial_{x}\rho|\leq C\mathcal{P},~
|\partial_{y}\rho|\leq C\mathcal{P}^{2},~
\left|\partial_{x}(\mathcal{P}\partial_{x}\rho)\right|\leq C(1+t\y)^{2}\mathcal{P}^{2},~\text{in}~\Omega_{T_{*}},
\end{align*}
where $\varpi$ and $\vartheta_{c}$ are given in \eqref{pi} and Section 3 respectively.
Moreover, the following constraints on $F$ hold,
\begin{align*}
|\partial_{y}F|\leq C_{*}\iota^{-1}(\dot{\mathfrak{U}}+\bar{\varsigma})~\text{in}~\Omega_{\hat{T}_{*}}~\text{and}~
\|\mathcal{P}(\partial_{y}^{2}F-\rho\partial_{y}F)\|_{L^{\infty}(\Omega_{\hat{T}_{*}})}\leq C_{*}\iota^{-1},
\end{align*}
$$
\|\partial_{t}F\|_{L_{t}^{\infty}H_{\psi}^{s-4}(\Omega_{\hat{T}_{*}})}+\|F\|_{L_{t}^{\infty}H_{\psi}^{s-2}(\Omega_{\hat{T}_{*}})}
\leq C_{*}\bar{\eta},
$$
$$
|D^{\gamma}F|\leq C_{*}e^{-\frac{4}{3}\delta\y}~\text{in}~\Omega_{\hat{T}_{*}}, \forall \gamma\in\Gamma_{2},
$$
and
\begin{align*}
\|\mathbf{F}_{s}\|_{L_{t}^{2}L_{\hat{\psi}}^{2}(\Omega_{\hat{T}_{*}})}^{2}
+\|F\|_{L_{t}^{2}\hat{H}_{\psi}^{s}(\Omega_{\hat{T}_{*}})}^{2}
+\sum_{i=0}^{\frac{s-1}{2}}\|\partial_{t}^{i}F\|_{L_{t}^{\infty}\mathring{H}^{s-1-2i}([0,\hat{T}_{*}]\times\mathbb{T})}^{2}
+\|\partial_{x}^{s}F\|_{L_{t}^{\infty}\mathring{H}^{0}([0,\hat{T}_{*}]\times\mathbb{T})}^{2}
\leq C_{*}\bar{\eta},
\end{align*}
where $\mathbf{F}_{s}=\mathcal{P}\left(\partial_{y}\partial_{x}^{s}F-\rho \partial_{x}^{s}F\right)$.
\end{assumption}
\begin{remark}
Assumption \ref{MA} is the iterative assumption. When $n=1$, the conditions in Assumption \ref{MA} is automatically satisfied by the construction of $u_{0}$, since $\mathfrak{U}=u_{0}$ and $\mathbb{U}=0$ in this case. The detailed verification will be given in Section 6.
\end{remark}

Now we can state the main result of this section.
\begin{theorem}\label{wpae}
There exists a positive constant $\lambda_{*}$ such that for any $\lambda>\lambda_{*}$, there are positive constants $\Lambda_{*}, \iota^{*}$ and $\mathcal{Z}_{*}$ depending only on $\lambda$ for which, when $\Lambda\geq\Lambda_{*}$, $\theta>0$, $\iota\in(0,\iota^{*})$ and $\mathcal{Z}\geq\mathcal{Z}_{*}$, there is $t^{*}\in(0,\min\{T_{\lambda,\Lambda},\frac{1}{12\Lambda},T\}]$ independent of $\bar{\eta}$ such that if Assumption \ref{MA} holds for any $t_{v}\in(0,t^{*}]$, moreover, the initial date $\check{u}$ satisfies
$$\|\check{u}\|_{\hat{H}_{\psi_{in}}^{s}(\Omega)}
 +\|\mathcal{P}(0,x,y)(\partial_{y}\partial_{x}^{s}\check{u}-\rho(0,x,y)\partial_{x}^{s}\check{u})\|_{L_{\varphi_{in}}^{2}(\Omega)}
\leq \frac{(1-\epsilon_{0})\sqrt{\bar{\eta}}\mathcal{Z}}{8\epsilon_{0}},~\text{in}~\Omega,$$
and $\check{u}\mathbb{U}=0$ in $\Omega_{t_{v}}$, where $\varphi_{in}=\varphi(0,y)$, then the problem \eqref{app} admits a unique solution
$u\in L_{t}^{\infty}\mathcal{H}_{\psi,\varphi}^{s}(\Omega_{t_{v}})$ satisfying
\begin{align}\label{priore}
& \|u\|_{L_{t}^{\infty}\mathcal{H}_{\psi,\varphi}^{s}(\Omega_{t_{v}})}
 +\|u\|_{L_{t}^{2}\hat{H}_{\psi,\y}^{s}(\Omega_{t_{v}})}
 +\|\partial_{y}u\|_{L_{t}^{2}\hat{H}_{\psi}^{s}(\Omega_{t_{v}})}
 +\|\partial_{x}u\|_{L_{t}^{2}\hat{H}_{\psi,\eta}^{s}(\Omega_{t_{v}})}\\
&+\|\mathcal{P}\mathcal{U}_{s}\|_{L_{t}^{2}L_{\varphi}^{2}(\Omega_{t_{v}})}
 +\|\sqrt{\y}\mathcal{U}_{s}\|_{L_{t}^{2}L_{\varphi}^{2}(\Omega_{t_{v}})}
 +\|\partial_{y}\mathcal{U}_{s}\|_{L_{t}^{2}L_{\varphi}^{2}(\Omega_{t_{v}})}
 +\|\sqrt{\eta}\partial_{x}\mathcal{U}_{s}\|_{L_{t}^{2}L_{\varphi}^{2}(\Omega_{t_{v}})}
\leq \frac{1-\epsilon_{0}}{\epsilon_{0}}\sqrt{\bar{\eta}}\mathcal{Z}.\nonumber
\end{align}
Here $T_{\lambda,\Lambda}$ is given in Lemma \ref{weighte}.
\end{theorem}

\begin{remark}
Here we should give an explanation for the condition $\check{u}\mathbb{U}=0$ in $\Omega_{t_{v}}$. In the estimate we established to gain \eqref{priore}, the coefficient of the initial term depends on $\mathbb{U}$, see \eqref{t*}. To get \eqref{priore}, we need $\check{u}$ to be zero when $\mathbb{U}$ is non-zero, and $\mathbb{U}$ to be zero when $\check{u}$ is non-zero, that is, $\check{u}\mathbb{U}=0$.
This condition is achievable during the iteration process, since $\mathbb{U}=\mathbb{U}_{0}=0$ when $n=1$ and the initial data of n-th order corrected increment solution $u_{n}(0,x,y)=0$ for any $n\geq 2$.
\end{remark}

In the remaining part of this paper, we use $C_{*}$ to stand for a positive constant depending only on $\lambda$ and $\mathcal{Z}$ which may be different from line to line.
~~~~~~~~~~

\subsection{Several lemmas}

Before the proof of Theorem \ref{wpae}, we shall introduce some essential preliminary lemmas.

The following lemma gives a Sobolev's type inequality, which is quoted from the reference \cite{MW}. As mentioned in \cite{MW}, the lemma can be proved by the standard extension argument and the Fourier's inversion formula. We quoted it without proof here.
\begin{lemma}\label{Linfty}
For any $f: \Omega\to\mathbb{R}$, there exists a universal positive constant $C$ such that
\begin{align*}
\|f\|_{L^{\infty}(\Omega)}\leq C\left(\|f\|_{L^{2}(\Omega)}+\|\partial_{x}f\|_{L^{2}(\Omega)}+\|\partial_{y}^{2}f\|_{L^{2}(\Omega)}\right).
\end{align*}
\end{lemma}

In deriving the energy estimates for the problem \eqref{app} by integration by parts, some boundary integral terms involving higher-order derivatives can not be controlled by a simple trace estimate. To handle these terms, we introduce the following two lemmas.
\begin{lemma}\label{Bc}
Assume that $T_{c}\in (0,T)$, and $\mathfrak{U}(t,x,y)$ is a solution to the problem \eqref{App} such that
$\mathbb{U}=(\mathfrak{U}-u_{0})\in L_{t}^{\infty}\hat{H}_{\psi}^{s}(\Omega_{T_{c}})$. For any $k\in\mathbb{N}_{+}$ and $n\in\mathbb{N}$ such that $n+2k\leq s+1$, there exists a positive constant $C$ independent of $\mathfrak{U}$, $F$ and $\partial_{x}P$ such that for any $t\in [0,T_{c}]$,
\begin{align}\label{Bc1}
\|\partial_{t}^{k}\mathfrak{U}\|_{\mathring{H}^{n}(\mathbb{T})}(t)
\leq&~C\sum_{\substack{p_{ij}>0,\\ i+j\leq k-1,\\0<i+j+p_{ij}+q_{j}\leq k}}\prod_{i=0}^{k-1}
       \left(\|\partial_{t}^{i}F\|_{\mathring{H}^{n+2j}(\mathbb{T})}+\bar{\tilde{\eta}}\|\partial_{t}^{i}u_{0}\|_{\mathring{H}^{n+2+2j}(\mathbb{T})}\right)^{p_{ij}}
       \|\mathfrak{U}\|_{\ddot{H}^{n+2j+1}(\hat{\Omega})}^{q_{j}}\\
    &+C\sum_{i=0}^{k}\|\mathfrak{U}\|_{\ddot{H}^{n+2i+1}(\hat{\Omega})}^{k-i+1}.\nonumber
\end{align}
Moreover, for any $k\in \mathbb{N}_{+}$ and $m,n\in \mathbb{N}$ such that $2k+m+2n\leq s+1$, there exists a positive constant $C$ and a binary polynomial function $\mathcal{Q}(\cdot,\cdot)$ independent of $\mathfrak{U}$, $F$ and $\partial_{x}P$ such that for any $t\in [0,T_{c}]$,
\begin{align}\label{Bc2}
\|\partial_{t}^{k}\partial_{y}^{2n}\partial_{x}^{m}\mathfrak{U}\|_{\mathring{H}^{0}(\mathbb{T})}
\leq&\mathcal{Q}\left(\|\mathfrak{U}\|_{\ddot{H}^{2n+2k+m-1}(\hat{\Omega})},
                      \sum_{i=0}^{n-2}\|\partial_{t}^{i+k}F\|_{\mathring{H}^{2n+m-2-2i}(\mathbb{T})}\right)\\
    &+C\sum_{i=0}^{n-1}\|\partial_{t}^{i+k}\partial_{x}P\|_{H^{2n+m-2-2i}(\mathbb{T})}.\nonumber
\end{align}
\end{lemma}

\begin{lemma}\label{bc}
Assume that $T_{c}\in (0,T)$ and $u\in L_{t}^{\infty}\hat{H}_{\psi}^{s}(\Omega_{T_{c}})$ is a solution to the problem \eqref{app},
and $\mathfrak{U}$ is a solution to the problem \eqref{App} such that $\mathbb{U}\in L_{t}^{\infty}\hat{H}_{\psi}^{s}(\Omega_{T_{c}})$.
For any $n\in\mathbb{N}_{+}$ and $m\in\mathbb{N}$ such that $2n+m\leq s+1$, there exist a constant $c_{*}$ depending only on $\|\mathfrak{U}\|_{L_{t}^{\infty}\ddot{H}^{2n+m-1}(\hat{\Omega}_{T_{c}})}$ and
$\sum\limits_{i=0}^{n-2}\|\partial_{t}^{i}F\|_{L_{t}^{\infty}\mathring{H}^{2n+m-2i-4}([0,T_{c}]\times\mathbb{T})}$,
and a constant  $c^{*}$ depending only on $\|\mathfrak{U}\|_{L_{t}^{\infty}\ddot{H}^{2n+m-1}(\hat{\Omega}_{T_{c}})}$,
$\sum\limits_{i=0}^{n-3}\|\partial_{t}^{i}F\|_{L_{t}^{\infty}\mathring{H}^{2n+m-2i-2}([0,T_{c}]\times\mathbb{T})}$,
and $\sum\limits_{i=0}^{n-2}\|\partial_{t}^{i}\partial_{x}P\|_{L_{t}^{\infty}H^{2n+m-2-2i}([0,T_{c}]\times\mathbb{T})}$,
such that for any $t\in [0,T_{c}]$,
\begin{align*}
\|\partial_{y}^{2n}\partial_{x}^{m}u\|_{\mathring{H}^{0}(\mathbb{T})}
\leq c_{*}\|u\|_{\ddot{H}^{2n+m-1}(\hat{\Omega})}
    +c_{*}\sum_{i=0}^{n-1}\|\partial_{t}^{i}F\|_{\mathring{H}^{2n+m-2i-2}(\mathbb{T})}
    +c^{*}\bar{\kappa}.
\end{align*}
\end{lemma}
The proof of Lemma \ref{Bc} and Lemma \ref{bc} can be found in the Appendix A.

The following lemma is critical to the proof of Theorem \ref{wpae}, which can help us to get the estimate of $\partial_{x}^{r}u$ and $\partial_{x}^{r}\mathbb{U}$ form $\mathcal{U}_{r}(t,x,y)$ and $\mathbf{U}_{r}(t,x,y)$.
\begin{lemma}\label{nc}
Assume $\psi_{c}(y)\in C^{\infty}(\mathbb{R}_{+})$ to be a positive nondecreasing function,
and $\mathfrak{g}(y)\in C^{1}(\mathbb{R}_{+})\cap L_{\psi_{c}}^{2}(\mathbb{R}_{+})$ to be a positive function.
Let $\Phi(y)=\psi_{c}(y)\omega_{c}(y)\left[\frac{\chi(y)}{\sqrt{\mathfrak{g}(y)}}+1-\chi(y)\right]^{2}$, with
$$\omega_{c}(y)=\frac{\tilde{\chi}(y)}{((y-y_{*})^{2}+\mathfrak{c})^{\frac{\varepsilon_{0}}{2}}}+1-\tilde{\chi}(y),$$
where $\epsilon_{0}$ is given as in Assumption \ref{inA}, $\mathfrak{c}$ is a positive constant, $\chi(y)$ is given in Section 2, and $\tilde{\chi}(y)\in C^{\infty}(\mathbb{R}_{+})$ is a function such that
$0\leq\tilde{\chi}(y)\leq 1$.
For any $f(y)\in L_{\psi_{c}}^{2}(\mathbb{R}_{+})$ satisfying $f(0)=0$, define $$\dot{f}(y)=f'(y)-\varrho(y)f(y),$$
with $\varrho(y)=\frac{\mathfrak{g}'(y)}{\mathfrak{g}(y)}$.
Assume there exist a constant $\ell\in(0,1)$ and negative constants $\check{\mathfrak{c}}$ and $\hat{\mathfrak{c}}$ such that
 $\check{\mathfrak{c}}\leq\varrho\leq -\ell^{-1}\frac{\psi_{c}'(y)}{\psi_{c}(y)}\leq\hat{\mathfrak{c}}$ for any $y\in[Y,+\infty)$,
and $\dot{f}(y)\in L_{\Phi}^{2}(\mathbb{R}_{+})$,
then for any nonnegative function $\xi(y)\in C^{1}[0,+\infty)$ such that $\xi'(y)\leq0$,
there exist a positive constant $\mathcal{N}$ depending on the values of $\mathfrak{g}$ on $[0,Y]\setminus[\check{y}_{*},\hat{y}_{*}]$, $\psi_{c}$ and $\ell$ but independent of $f$ and $\mathfrak{c}$ such that
\begin{align}\label{nce0}
\|\xi f\|_{L_{\Phi}^{2}(\mathbb{R}_{+})}\leq \mathcal{N}\|\xi\dot{f}\|_{L_{\Phi}^{2}(\mathbb{R}_{+})}.
\end{align}
\end{lemma}

\begin{proof}
For convenience, in the proof below we denote by $C^{*}$ a positive constant depending only on $\psi_{c}$ and the values of $\mathfrak{g}$ on $[0,Y]\setminus[\check{y}_{*},\hat{y}_{*}]$, which may be different from line to line.

Since $f(0)=0$, by noticing that $\dot{f}(y)=\mathfrak{g}(y)\left(\frac{f(y)}{\mathfrak{g}(y)}\right)'$ one has
\begin{align*}
\int_{0}^{Y}\Phi^{2}(y)\xi^{2}(y)f^{2}(y)dy
\leq&~C^{*}\int_{0}^{Y}\omega_{c}^{2}\frac{\xi^{2}(y)f^{2}(y)}{\mathfrak{g}^{2}(y)}dy\\
\leq&~C^{*}\int_{0}^{Y}\omega_{c}^{2}\int_{0}^{y}\left[\frac{\xi^{2}(s)f^{2}(s)}{\mathfrak{g}^{2}(s)}\right]'dsdy\nonumber\\
=&~2C^{*}\int_{0}^{Y}\omega_{c}^{2}\int_{0}^{y}\xi^{2}\frac{\dot{f}(s)f(s)}{\mathfrak{g}^{2}(s)}dsdy
  +2C^{*}\int_{0}^{Y}\omega_{c}^{2}\int_{0}^{y}\frac{\xi(s)\xi'(s)f^{2}(s)}{\mathfrak{g}^{2}(s)}dsdy\nonumber\\
\leq&~C^{*}\|\omega_{c}\|_{L^{2}[0,Y]}^{2}\|\xi\mathfrak{g}^{-1}f\|_{L^{2}[0,Y]}\|\xi\mathfrak{g}^{-1}\dot{f}\|_{L^{2}[0,Y]}\nonumber\\
\leq&~C^{*}\|\xi f\|_{L_{\Phi}^{2}[0,Y]}\|\xi\dot{f}\|_{L_{\Phi}^{2}[0,Y]},\nonumber
\end{align*}
which gives
\begin{align}\label{nc1}
\|\xi f\|_{L_{\Phi}^{2}[0,Y]}^{2}
\leq C^{*}\|\xi\dot{f}\|_{L_{\Phi}^{2}[0,Y]}^{2}.
\end{align}
Similarly, one has
\begin{align}\label{nc2}
\frac{\xi^{2}(Y)f^{2}(Y)}{\mathfrak{g}^{2}(Y)}
=&~\int_{0}^{Y}\left[\frac{\xi^{2}(y)f^{2}(y)}{\mathfrak{g}^{2}(y)}\right]'dy\\
=&~2\int_{0}^{Y}\frac{\xi^{2}(y)f(y)\dot{f}(y)}{\mathfrak{g}^{2}(y)}dy
  +2\int_{0}^{Y}\frac{\xi(y)\xi'(y)f(y)^{2}}{\mathfrak{g}^{2}(y)}dy\nonumber\\
\leq&~C^{*}\|\xi f\|_{L_{\Phi}^{2}[0,Y]}\|\xi\dot{f}\|_{L_{\Phi}^{2}[0,Y]}\nonumber\\
\leq&~C^{*}\|\xi\dot{f}\|_{L_{\Phi}^{2}[0,Y]}^{2}.\nonumber
\end{align}
Since $f(y)\in L_{\psi_{c}}^{2}(\mathbb{R}_{+})$, one has $\lim\limits_{y\to+\infty}\psi_{c}(y)f(y)=0$, thus
\begin{align*}
-\psi_{c}^{2}(Y)\xi(Y)^{2}f^{2}(Y)
=&\int_{Y}^{+\infty}\left[\psi_{c}^{2}(y)\xi^{2}(y)f^{2}(y)\right]'dy\\
=& 2\int_{Y}^{+\infty}\psi_{c}(y)\psi_{c}'(y)\xi^{2}(y)f^{2}(y)dy
  +2\int_{Y}^{+\infty}\psi_{c}^{2}(y)\xi(y)\xi'(y)f^{2}(y)dy\nonumber\\
 &+2\int_{Y}^{+\infty}\psi_{c}^{2}(y)\varrho(y)\xi^{2}(y)f^{2}(y)dy
  +2\int_{Y}^{+\infty}\psi_{c}^{2}(y)\xi^{2}(y)f(y)\dot{f}(y)dy.\nonumber
\end{align*}
As a result, one has
\begin{align}\label{nc3}
-\int_{Y}^{+\infty}\psi_{c}^{2}(y)\xi^{2}(y)f(y)\dot{f}(y)dy
=& \frac{1}{2}\psi_{c}^{2}(Y)\xi(Y)^{2}f^{2}(Y)
  +\int_{Y}^{+\infty}\psi_{c}(y)\psi_{c}'(y)\xi^{2}(y)f^{2}(y)dy\\
 &+\int_{Y}^{+\infty}\psi_{c}^{2}(y)\xi(y)\xi'(y)f^{2}(y)dy
  +\int_{Y}^{+\infty}\psi_{c}^{2}(y)\varrho(y)\xi^{2}(y)f^{2}(y)dy.\nonumber
\end{align}
Since $\frac{\psi_{c}'(y)}{\psi_{c}(y)}\leq \ell|\varrho(y)|$ for any $y\in[Y,+\infty)$, from \eqref{nc3} we know
\begin{align*}
&(1-\ell)\int_{Y}^{+\infty}\psi_{c}^{2}(y)|\varrho(y)|\xi^{2}(y)f^2(y)dy
-\int_{Y}^{+\infty}\psi_{c}^{2}(y)\xi(y)\xi'(y)f^{2}(y)dy\\
\leq& \int_{Y}^{+\infty}\psi_{c}^{2}(y)\xi^{2}(y)|f(y)\dot{f}(y)|dy
    +\frac{1}{2}\xi(Y)^{2}\psi_{c}^{2}(Y)f^{2}(Y) ,
\end{align*}
hence  from $\xi'(y)\leq0$ we can deduce
\begin{align}\label{nc4}
\int_{Y}^{+\infty}\psi_{c}^{2}(y)|\varrho(y)|\xi^{2}(y)f^2(y)dy
\leq& \frac{1}{(1-\ell)^{2}}\int_{Y}^{+\infty}\psi_{c}^{2}(y)|\varrho^{-1}(y)|\xi^{2}(y)\dot{f}^{2}(y)dy\\
    &+\frac{\xi(Y)^{2}\psi_{c}^{2}(Y)f^{2}(Y)}{1-\ell}.\nonumber
\end{align}
Substituting \eqref{nc2} into \eqref{nc4}, we have
\begin{align}\label{nc5}
\|\xi f\|_{L_{\Phi}^{2}(Y,+\infty)}^{2}
\leq C^{*}\|\xi\dot{f}\|_{L_{\Phi}^{2}(\mathbb{R}_{+})}^{2}.
\end{align}
Combining \eqref{nc1} with \eqref{nc5}, we get \eqref{nce0}.
\end{proof}

\begin{corollary}\label{ncC1}
Assume $\tilde{\xi}(y)\in C^{1}[0,+\infty)$ to be a nonnegative bounded function. Let $\mathcal{N}$ be given in Lemma \ref{nc}.
Under the assumption of Lemma \ref{nc}, if $|\tilde{\xi}'(y)|\leq \frac{1}{2\mathcal{N}}\tilde{\xi}(y)$ for any $y\in[0,+\infty)$, then
\begin{align*}
\|\tilde{\xi} f\|_{L_{\Phi}^{2}(\mathbb{R}_{+})}\leq 2\mathcal{N}\|\tilde{\xi}\dot{f}\|_{L_{\Phi}^{2}(\mathbb{R}_{+})}.
\end{align*}

\end{corollary}
\begin{proof}
Taking $\xi\equiv1$ in \eqref{nce0}, one has
\begin{align}\label{nce0'}
\|f\|_{L_{\Phi}^{2}(\mathbb{R}_{+})}\leq C\|\dot{f}\|_{L_{\Phi}^{2}(\mathbb{R}_{+})}.
\end{align}
By applying \eqref{nce0'}, one knows that
\begin{align}\label{ncc1}
\|\tilde{\xi}f\|_{L_{\Phi}^{2}(\mathbb{R}_{+})}
\leq \mathcal{N}\|\tilde{\xi}'f+\tilde{\xi}\dot{f}\|_{L_{\Phi}^{2}(\mathbb{R}_{+})}
\leq \mathcal{N}\|\tilde{\xi}'f\|_{L_{\Phi}^{2}(\mathbb{R}_{+})}
    +\mathcal{N}\|\tilde{\xi}\dot{f}\|_{L_{\Phi}^{2}(\mathbb{R}_{+})}
\end{align}
If $|\tilde{\xi}'(y)|\leq \frac{1}{2\mathcal{N}}\tilde{\xi}$ for any $y\in[0,+\infty)$, then
$\mathcal{N}\|\tilde{\xi}'f\|_{L_{\Phi}^{2}(\mathbb{R}_{+})}\leq \frac{1}{2}\|\tilde{\xi}f\|_{L_{\Phi}^{2}(\mathbb{R}_{+})}$.
From \eqref{ncc1} one can deduce
\begin{align*}
\|\tilde{\xi} f\|_{L_{\Phi}^{2}(\mathbb{R}_{+})}
\leq 2\mathcal{N}\|\tilde{\xi}\dot{f}\|_{L_{\Phi}^{2}(\mathbb{R}_{+})}.
\end{align*}
which gives the conclusion.
\end{proof}

Let $\rho$ and $\mathcal{P}$ be given in \eqref{rhoPd}. For any positive integer $r$ such that $r\leq s+1$, and function $f:\Omega_{T}\to\mathbb{R}$, let $$\mathcal{F}_{r}=\mathcal{P}(\partial_{y}\partial_{x}^{r}f-\rho\partial_{x}^{r}f).$$
Define $\Psi(t,x,y)=\mathcal{P}^{2}(t,x,y)\varphi_{o}(t,y)$ and $\Psi_{c}(t,x,)=\mathcal{P}^{2}(t,x,y)\psi(t,y)$, with
$$\omega_{o}(t,y)=\frac{\chi(y)}{((y-y_{*})^{2}+\bar{\varsigma}+t)^{\frac{\varepsilon_{0}}{2}}}+1-\chi(y),~\text{and,}~\varphi_{o}(t,y)=\psi(t,y)\omega_{o}(t,y).$$
Recall the definition of $\varphi$ given in \eqref{varphi}, one can check that there exists positive constants $C_{\lambda,1}, C_{\lambda,2}$ depending only on $\lambda$ such that $C_{\lambda,1}\varphi\leq\varphi_{o}\y^{-\frac{1}{2}}\leq C_{\lambda,2}\varphi$.
By Lemma \ref{nc}, it follows the following corollaries, which will be used frequently in the sequel.
\begin{corollary}\label{ncC2}
Assume there exist $T_{c}\in(0,T)$ such that
$$\dot{\mathfrak{U}}(t,x,\cdot)\in C^{1}(\mathbb{R}_{+})\cap L_{\psi}^{2}(\mathbb{R}_{+}),~\forall (t,x)\in[0,T_{c}]\times\mathbb{T},$$
and
$$\mathcal{P}^{-1}\leq C,~|\rho|\leq C\mathcal{P},~|\partial_{x}\rho|\leq C\mathcal{P},~|\partial_{y}\rho|\leq C\mathcal{P}^{2},~\text{in}~ \Omega_{T_{c}},$$
then for any positive integer $r$ such that $r\leq s+1$ and any nonnegative function $\xi(y)\in C^{1}[0,+\infty)$ such that $\xi'(y)\leq0$, the following norm inequalities hold for any $t\in[0,T_{c}]$,
\begin{align}\label{ncc21}
\|\xi\partial_{x}^{r}f\|_{L_{\Psi}^{2}(\Omega)}
\leq C\|\xi\mathcal{P}\mathcal{F}_{r}\|_{L_{\varphi_{o}}^{2}(\Omega)},~
\|\xi\mathcal{P}\partial_{y}\partial_{x}^{r}f\|_{L_{\varphi_{o}}^{2}(\Omega)}
\leq C\|\xi\mathcal{F}_{r}\|_{L_{\varphi_{o}}^{2}(\Omega)}
    +C\|\xi\mathcal{P}\mathcal{F}_{r}\|_{L_{\varphi_{o}}^{2}(\Omega)},
\end{align}
and
\begin{align}\label{ncc22}
\|\xi\partial_{y}^{2}\partial_{x}^{r}f\|_{L_{\varphi_{o}}^{2}(\Omega)}
\leq C\|\xi\mathcal{F}_{r}\|_{L_{\varphi_{o}}^{2}(\Omega)}
    +C\|\xi\partial_{y}\mathcal{F}_{r}\|_{L_{\varphi_{o}}^{2}(\Omega)}
    +C\|\xi\mathcal{P}\mathcal{F}_{r}\|_{L_{\varphi_{o}}^{2}(\Omega)},
\end{align}
moreover,
\begin{align}\label{ncc23}
\|\xi\partial_{x}^{r}f\|_{L_{\Psi_{c}}^{2}(\Omega)}
\leq C\|\xi\mathcal{P}\mathcal{F}_{r}\|_{L_{\psi}^{2}(\Omega)},~
\|\xi\mathcal{P}\partial_{y}\partial_{x}^{r}f\|_{L_{\Psi_{c}}^{2}(\Omega)}
\leq C\|\xi\mathcal{F}_{r}\|_{L_{\psi}^{2}(\Omega)}
    +C\|\xi\mathcal{P}\mathcal{F}_{r}\|_{L_{\psi}^{2}(\Omega)},
\end{align}
\begin{align}\label{ncc24}
\|\xi\partial_{y}^{2}\partial_{x}^{r}f\|_{L_{\Psi_{c}}^{2}(\Omega)}
\leq C\|\xi\mathcal{F}_{r}\|_{L_{\psi}^{2}(\Omega)}
    +C\|\xi\partial_{y}\mathcal{F}_{r}\|_{L_{\psi}^{2}(\Omega)}
    +C\|\xi\mathcal{P}\mathcal{F}_{r}\|_{L_{\psi}^{2}(\Omega)}.
\end{align}
\end{corollary}
\begin{proof}
By using Lemma \ref{nc} with $\tilde{\chi}=\chi$, one can obtain
\begin{align*}
\|\xi\partial_{x}^{r}f\|_{L_{\Psi}^{2}(\Omega)}
\leq C\|\xi\mathcal{P}^{-1}\mathcal{F}_{r}\|_{L_{\Psi}^{2}(\Omega)}
\leq C\|\xi\mathcal{P}\mathcal{F}_{r}\|_{L_{\varphi_{o}}^{2}(\Omega)},
\end{align*}
and recalling the definition of $\mathcal{F}_{r}$ one has,
\begin{align*}
\|\xi\mathcal{P}\partial_{y}\partial_{x}^{r}f\|_{L_{\varphi_{o}}^{2}(\Omega)}
\leq& \|\xi\mathcal{F}_{r}\|_{L_{\varphi_{o}}^{2}(\Omega)}
     +\|\xi\mathcal{P}\rho\partial_{x}^{r}f\|_{L_{\varphi_{o}}^{2}(\Omega)}\\
\leq& C\|\xi\mathcal{F}_{r}\|_{L_{\varphi_{o}}^{2}(\Omega)}
     +C\|\xi\partial_{x}^{r}f\|_{L_{\Psi}^{2}(\Omega)}\nonumber\\
\leq& C\|\xi\mathcal{F}_{r}\|_{L_{\varphi_{o}}^{2}(\Omega)}
     +C\|\xi\mathcal{P}\mathcal{F}_{r}\|_{L_{\varphi_{o}}^{2}(\Omega)}.
\end{align*}
Since
\begin{align*}
\partial_{y}^{2}\partial_{x}^{r}f
= \mathcal{P}^{-1}\partial_{y}\mathcal{F}_{r}
 -\mathcal{P}^{-2}\partial_{y}\mathcal{P}\mathcal{F}_{r}
 +\partial_{y}\rho\partial_{x}^{r}f
 +\rho\partial_{y}\partial_{x}^{r}f
\end{align*}
we have
\begin{align}\label{ncC2-1}
\|\xi\partial_{y}^{2}\partial_{x}^{r}f\|_{L_{\varphi_{o}}^{2}(\Omega)}
\leq&~\|\xi\mathcal{P}^{-1}\partial_{y}\mathcal{F}_{r}\|_{L_{\varphi_{o}}^{2}(\Omega)}
     +\|\xi\mathcal{P}^{-2}\partial_{y}\mathcal{P}\mathcal{F}_{r}\|_{L_{\varphi_{o}}^{2}(\Omega)}\\
    &+\|\xi\partial_{y}\rho\partial_{x}^{r}f\|_{L_{\varphi_{o}}^{2}(\Omega)}
     +\|\xi\rho\partial_{y}\partial_{x}^{r}f\|_{L_{\varphi_{o}}^{2}(\Omega)}.\nonumber
\end{align}
Since $|\rho|\leq C\mathcal{P}$ one has $|\partial_{y}\mathcal{P}|\leq C\mathcal{P}^{2}$, which together with $|\partial_{y}\rho|\leq C\mathcal{P}^{2}$ and Lemma \ref{nc} gives
\begin{align}\label{ncC2-11}
\|\xi\mathcal{P}^{-2}\partial_{y}\mathcal{P}\mathcal{F}_{r}\|_{L_{\varphi_{o}}^{2}(\Omega)}
+\|\xi\partial_{y}\rho\partial_{x}^{r}f\|_{L_{\varphi_{o}}^{2}(\Omega)}
\leq& C\|\xi\mathcal{F}_{r}\|_{L_{\varphi_{o}}^{2}(\Omega)}
     +C\|\xi\mathcal{P}^{2}\partial_{x}^{r}f\|_{L_{\varphi_{o}}^{2}(\Omega)}\\
\leq& C\|\xi\mathcal{F}_{r}\|_{L_{\varphi_{o}}^{2}(\Omega)}
     +C\|\xi\partial_{x}^{r}f\|_{L_{\Psi}^{2}(\Omega)}\nonumber\\
\leq& C\|\xi\mathcal{F}_{r}\|_{L_{\varphi_{o}}^{2}(\Omega)}
     +C\|\xi\mathcal{P}\mathcal{F}_{r}\|_{L_{\varphi_{o}}^{2}(\Omega)}.\nonumber
\end{align}
By using $|\rho|\leq C\mathcal{P}$ again, we have
\begin{align}\label{ncC2-12}
\|\xi\rho\partial_{y}\partial_{x}^{r}f\|_{L_{\varphi_{o}}^{2}(\Omega)}
\leq& \|\xi\rho(\partial_{y}\partial_{x}^{r}f-\rho\partial_{x}^{r}f)\|_{L_{\varphi_{o}}^{2}(\Omega)}
     +\|\xi\rho^{2}\partial_{x}^{r}f\|_{L_{\varphi_{o}}^{2}(\Omega)}\\
\leq& C\|\xi\mathcal{F}_{r}\|_{L_{\varphi_{o}}^{2}(\Omega)}
     +C\|\xi\partial_{x}^{r}f\|_{L_{\Psi}^{2}(\Omega)}\nonumber\\
\leq& C\|\xi\mathcal{F}_{r}\|_{L_{\varphi_{o}}^{2}(\Omega)}
     +C\|\xi\mathcal{P}\mathcal{F}_{r}\|_{L_{\varphi_{o}}^{2}(\Omega)}.\nonumber
\end{align}
Combining \eqref{ncC2-11} and \eqref{ncC2-12} with \eqref{ncC2-1}, we have \eqref{ncc22}.

By the same argument as above, one can obtain the inequality \eqref{ncc23} and \eqref{ncc24} by using Lemma \ref{nc} with $\tilde{\chi}=0$. Here we omit the detail.
\end{proof}

\begin{corollary}\label{ncC3}
Assume there exist $T_{c}\in(0,T)$ such that
$$\dot{\mathfrak{U}}(t,x,\cdot)\in C^{1}(\mathbb{R}_{+})\cap L_{\psi}^{2}(\mathbb{R}_{+}),~\forall (t,x)\in[0,T_{c}]\times\mathbb{T},$$
and
$$\mathcal{P}^{-1}\leq C,~\partial_{x}\mathcal{P}^{-1}\leq C,~|\rho|\leq C\mathcal{P},~|\partial_{x}\rho|\leq C\mathcal{P},~\text{in}~ \Omega_{T_{c}}.$$
Assume that $\xi(y)\in C^{1}[0,+\infty)$ is a nonnegative function such that $\xi'(y)\leq0$, then one has the following norm inequalities for any $t\in[0,T_{c}]$,
\begin{align}\label{ncc3}
\|\xi\partial_{x}^{r+1}f\|_{L_{\Psi}^{2}(\Omega)}
\leq C\|\mathcal{P}\|_{L^{\infty}(\Omega)}\|\xi\partial_{x}\mathcal{F}_{r}\|_{L_{\varphi_{o}}^{2}(\Omega)}
    +C\|\mathcal{P}\|_{L^{\infty}(\Omega)}\|\xi\mathcal{P}\mathcal{F}_{r}\|_{L_{\varphi_{o}}^{2}(\Omega)},
\end{align}
and
\begin{align}\label{ncc3'}
\|\xi\partial_{x}^{r+1}f\|_{L_{\Psi_{c}}^{2}(\Omega)}
\leq C\|\mathcal{P}\omega_{o}^{-1}\|_{L^{\infty}(\Omega)}\|\xi\partial_{x}\mathcal{F}_{r}\|_{L_{\varphi_{o}}^{2}(\Omega)}
    +C\|\mathcal{P}\omega_{o}^{-1}\|_{L^{\infty}(\Omega)}\|\xi\mathcal{P}\mathcal{F}_{r}\|_{L_{\varphi_{o}}^{2}(\Omega)},
\end{align}

\end{corollary}
\begin{proof}
Since
\begin{align}\label{ncc30}
\mathcal{P}^{-1}\mathcal{F}_{r+1}
=\partial_{y}\partial_{x}^{r+1}f-\rho\partial_{x}^{r+1}f
=\partial_{x}(\partial_{y}\partial_{x}^{r}f-\rho\partial_{x}^{r}f)+\partial_{x}\rho\partial_{x}^{r}f
=\partial_{x}(\mathcal{P}^{-1}\mathcal{F}_{r})+\partial_{x}\rho\partial_{x}^{r}f,
\end{align}
by applying \eqref{ncc21}, we can deduce from \eqref{ncc30} that
\begin{align}\label{ncc31}
\|\xi\partial_{x}^{r+1}f\|_{L_{\Psi}^{2}(\Omega)}
\leq&~C\|\xi\mathcal{P}^{-1}\mathcal{F}_{r+1}\|_{L_{\Psi}^{2}(\Omega)}\\
\leq&~C|\xi\partial_{x}(\mathcal{P}^{-1}\mathcal{F}_{r})\|_{L_{\Psi}^{2}(\Omega)}
     +C\|\xi\partial_{x}\rho\partial_{x}^{r}f\|_{L_{\Psi}^{2}(\Omega)}\nonumber\\
\leq&~C\|\xi\mathcal{P}\partial_{x}\mathcal{F}_{r}\|_{L_{\varphi_{o}}^{2}(\Omega)}
     +C\|\xi\partial_{x}\mathcal{P}^{-1}\mathcal{P}^{2}\mathcal{F}_{r}\|_{L_{\varphi_{o}}^{2}(\Omega)}
     +C\|\xi\mathcal{P}\partial_{x}^{r}f\|_{L_{\Psi}^{2}(\Omega)}\nonumber\\
\leq&~C\|\xi\mathcal{P}\partial_{x}\mathcal{F}_{r}\|_{L_{\varphi_{o}}^{2}(\Omega)}
     +C\|\xi\mathcal{P}^{2}\mathcal{F}_{r}\|_{L_{\varphi_{o}}^{2}(\Omega)}
     +C\|\xi\mathcal{P}\partial_{x}^{r}f\|_{L_{\Psi}^{2}(\Omega)}\nonumber\\
\leq&~C\|\mathcal{P}\|_{L^{\infty}(\Omega)}\|\xi\partial_{x}\mathcal{F}_{r}\|_{L_{\varphi_{o}}^{2}(\Omega)}
     +C\|\mathcal{P}\|_{L^{\infty}(\Omega)}\|\xi\mathcal{P}\mathcal{F}_{r}\|_{L_{\varphi_{o}}^{2}(\Omega)}\nonumber\\
    &+C\|\mathcal{P}\|_{L^{\infty}(\Omega)}\|\xi\partial_{x}^{r}f\|_{L_{\Psi}^{2}(\Omega)}\nonumber\\
\leq&~C\|\mathcal{P}\|_{L^{\infty}(\Omega)}\|\xi\partial_{x}\mathcal{F}_{r}\|_{L_{\varphi_{o}}^{2}(\Omega)}
     +C\|\mathcal{P}\|_{L^{\infty}(\Omega)}\|\xi\mathcal{P}\mathcal{F}_{r}\|_{L_{\varphi_{o}}^{2}(\Omega)}.\nonumber
\end{align}
Similarly, by using \eqref{ncc23} one has
\begin{align}\label{ncc32}
\|\xi\partial_{x}^{r+1}f\|_{L_{\Psi_{c}}^{2}(\Omega)}
\leq&~C\|\xi\mathcal{P}^{-1}\mathcal{F}_{r+1}\|_{L_{\Psi_{c}}^{2}(\Omega)}\\
\leq&~C\|\xi\mathcal{P}\partial_{x}\mathcal{F}_{r}\|_{L_{\psi}^{2}(\Omega)}
     +C\|\xi\mathcal{P}^{2}\mathcal{F}_{r}\|_{L_{\psi}^{2}(\Omega)}
     +C\|\xi\mathcal{P}\partial_{x}^{r}f\|_{L_{\Psi_{c}}^{2}(\Omega)}\nonumber\\
\leq&~C\|\mathcal{P}\omega_{o}^{-1}\|_{L^{\infty}(\Omega)}\|\xi\partial_{x}\mathcal{F}_{r}\|_{L_{\varphi_{o}}^{2}(\Omega)}
     +C\|\mathcal{P}\omega_{o}^{-1}\|_{L^{\infty}(\Omega)}\|\xi\mathcal{P}\mathcal{F}_{r}\|_{L_{\varphi_{o}}^{2}(\Omega)}\nonumber\\
    &+C\|\mathcal{P}\omega_{o}^{-1}\|_{L^{\infty}(\Omega)}\|\xi\partial_{x}^{r}f\|_{L_{\Psi}^{2}(\Omega)}\nonumber\\
\leq&~C\|\mathcal{P}\omega_{o}^{-1}\|_{L^{\infty}(\Omega)}\|\xi\partial_{x}\mathcal{F}_{r}\|_{L_{\varphi_{o}}^{2}(\Omega)}
     +C\|\mathcal{P}\omega_{o}^{-1}\|_{L^{\infty}(\Omega)}\|\xi\mathcal{P}\mathcal{F}_{r}\|_{L_{\varphi_{o}}^{2}(\Omega)}.\nonumber
\end{align}
\end{proof}

\begin{corollary}\label{ncC4}
Assume there exist $T_{c}\in(0,T)$ such that
$$\dot{\mathfrak{U}}(t,x,\cdot)\in C^{1}(\mathbb{R}_{+})\cap L_{\psi}^{2}(\mathbb{R}_{+}),~\forall (t,x)\in[0,T_{c}]\times\mathbb{T},$$
and
$$\mathcal{P}^{-1}\leq C,~\partial_{x}\mathcal{P}^{-1}\leq C,~|\rho|\leq C\mathcal{P},~|\partial_{x}\rho|\leq C\mathcal{P},~\text{in}~ \Omega_{T_{c}}.$$
Assume $\tilde{\xi}(y)\in C^{1}[0,+\infty)$ is a nonnegative bounded function. Under the assumption of Corollary \ref{ncC1}, there exists a positive constant $\tilde{l}$ such that if $|\tilde{\xi}'(y)|\leq \tilde{l}\tilde{\xi}(y)$ for any $y\in[0,+\infty)$, then the following norm inequalities hold for any $t\in[0,T_{c}]$,
\begin{align}\label{ncc41}
\|\tilde{\xi}\partial_{x}^{r+1}f\|_{L_{\Psi}^{2}(\Omega)}
\leq C\|\mathcal{P}\|_{L^{\infty}(\Omega)}\|\tilde{\xi}\partial_{x}\mathcal{F}_{r}\|_{L_{\varphi_{o}}^{2}(\Omega)}
    +C\|\mathcal{P}\|_{L^{\infty}(\Omega)}\|\tilde{\xi}\mathcal{P}\mathcal{F}_{r}\|_{L_{\varphi_{o}}^{2}(\Omega)}.
\end{align}
and
\begin{align}\label{ncc42}
\|\tilde{\xi}\partial_{x}^{r+1}f\|_{L_{\Psi_{c}}^{2}(\Omega)}
\leq C\|\mathcal{P}\omega_{o}^{-1}\|_{L^{\infty}(\Omega)}\|\tilde{\xi}\partial_{x}\mathcal{F}_{r}\|_{L_{\varphi_{o}}^{2}(\Omega)}
    +C\|\mathcal{P}\omega_{o}^{-1}\|_{L^{\infty}(\Omega)}\|\tilde{\xi}\mathcal{P}\mathcal{F}_{r}\|_{L_{\varphi_{o}}^{2}(\Omega)}.
\end{align}
\end{corollary}
\begin{proof}
By using Corollary \ref{ncC1}, one knows that there exists a positive constant $\tilde{l}$ such that if $|\tilde{\xi}'(y)|\leq \tilde{l}\tilde{\xi}(y)$ for any $y\in[0,+\infty)$, then
\begin{align}\label{ncc41'}
\|\tilde{\xi}\partial_{x}^{r+1}f\|_{L_{\Psi}^{2}(\Omega)}
\leq C\|\tilde{\xi}\mathcal{P}^{-1}\mathcal{F}_{r+1}\|_{L_{\Psi}^{2}(\Omega)},
\end{align}
and
\begin{align}\label{ncc42'}
\|\xi\partial_{x}^{r+1}f\|_{L_{\Psi_{c}}^{2}(\Omega)}
\leq C\|\xi\mathcal{P}^{-1}\mathcal{F}_{r+1}\|_{L_{\Psi_{c}}^{2}(\Omega)}.
\end{align}
With the same approach as in \eqref{ncc31} and \eqref{ncc32}, one can obtain \eqref{ncc41} and \eqref{ncc42} from \eqref{ncc41'} and \eqref{ncc42'} respectively.
\end{proof}

Now we are in a position to prove Theorem \ref{wpae}. Note that the problem \eqref{app} is nondegenerate parabolic, once the uniform (in $\bar{\eta}$) a priori estimate is obtained, the existence of solution in a time interval independent of $\bar{\eta}$ follows immediately by the well-posedness theorem for classical nonlinear parabolic problem, and the uniqueness of the solution follows by the $L^{2}$ comparison. So to prove Theorem \ref{wpae}, we just need to establish the uniform (in $\bar{\eta}$) a priori estimate, that is what we shall do in the remainder of this section.

To determine $\lambda_{*}$, and thereby $\Lambda_{*},\iota^{*},\mathcal{Z}_{*},t^{*}$ such that Theorem \ref{wpae} holds, we first choose $\lambda,\Lambda$ large enough and $\iota$ small enough, all to be determined later, such that Lemma \ref{weighte} holds. For such $\lambda$ and $\Lambda$, let $T_{\lambda,\Lambda}$ be given as in Lemma \ref{weighte}. Assume that Assumption \ref{MA} hold for
$T_{*}=\min\{T_{\lambda,\Lambda},\frac{1}{12\Lambda},T\}$, then for a sufficient large $\lambda$, we can find some $t^{*}\in(0,T_{*}]$, and establish a uniform energy estimate \eqref{priore} on $[0,t^{*}]$ for the problem \eqref{app}, by taking suitable $\Lambda,\iota$ and $\mathcal{Z}$. Finally, note that for any $t_{v}\in[0,t^{*}]$, to establish the uniform energy estimate on $[0,t_{v}]$, it is sufficient to assume that Assumption \ref{MA} hold for $t_{v}$ and $\check{u}\mathbb{U}=0$ in $\Omega_{t_{v}}$, we get Theorem \ref{wpae}.

~~~~~~~~~~

\subsection{Uniform estimate with low-order $x$-derivatives}~

This section devotes to the uniform estimate of $\|u\|_{\hat{H}_{\psi}^{s}(\Omega)}$, which can be done by using the standard energy method.
For any multiindex $\gamma=(\gamma_{1},\gamma_{2})$ such that
$$\gamma\in\Gamma=\{(\gamma_{1},\gamma_{2})|(\gamma_{1},\gamma_{2})\in\mathbb{N}^{2},\gamma_{1}+\gamma_{2}\leq s,\gamma<s\}\cup\{(s-1,2)\},$$
let $D^{\gamma}=\partial_{x}^{\gamma_{1}}\partial_{y}^{\gamma_{2}}$. By acting $D^{\gamma}$ on \eqref{app}$_{1}$, one can obtain
\begin{align}\label{lde}
\partial_{t}D^{\gamma}u
-\partial_{y}^{2}D^{\gamma}u
-\eta\partial_{x}^{2}D^{\gamma}u
=&-D^{\gamma}[\mathfrak{U}\partial_{x}u]
  -D^{\gamma}[u\partial_{x}\mathfrak{U}]
  +D^{\gamma}[\partial_{y}^{-1}[\partial_{x}u](\partial_{y}\mathfrak{U}+\varsigma)]
  +D^{\gamma}[\partial_{y}^{-1}[\partial_{x}\mathfrak{U}]\partial_{y}u]\\
 &+D^{\gamma}[h\partial_{y}^{-1}[\partial_{x}u]\partial_{y}u]
  +D^{\gamma}[\kappa\partial_{x}^{2}\mathbb{U}]
  +D^{\gamma}F
  +\left(D^{\gamma}[\eta\partial_{x}^{2}u]-\eta\partial_{x}^{2}D^{\gamma}u\right).\nonumber
\end{align}
Let $\psi_{\gamma}=\psi$ for $|\gamma|\leq s$ and $\psi_{\gamma}=\hat{\psi}$ for $|\gamma|=s+1$, moreover, $\hat{\psi}_{\gamma}=\y^{-\frac{1}{2}}\psi_{\gamma}$.
For any $t\in[0,T_{*}]$, multiplying \eqref{lde} by $\psi_{\gamma}^{2}D^{\gamma}u$ and integrating over $\Omega$, using the integration by parts we get
\begin{align}\label{le}
&\frac{1}{2}\frac{d}{dt}\|D^{\gamma}u \|_{L_{\psi_{\gamma}}^{2}(\Omega)}^{2}
+\Lambda\delta\|\sqrt{\y}D^{\gamma}u \|_{L_{\psi_{\gamma}}^{2}(\Omega)}^{2}
+\|\partial_{y}D^{\gamma}u\|_{L_{\psi_{\gamma}}^{2}(\Omega)}^{2}
+\|\sqrt{\eta}\partial_{x}D^{\gamma}u\|_{L_{\psi_{\gamma}}^{2}(\Omega)}^{2}\\
=&-(D^{\gamma}[\mathfrak{U}\partial_{x}u],\psi_{\gamma}^{2} D^{\gamma}u)
  -(D^{\gamma}[u\partial_{x}\mathfrak{U}],\psi_{\gamma}^{2} D^{\gamma}u)
  +(D^{\gamma}[\partial_{y}^{-1}[\partial_{x}u](\partial_{y}\mathfrak{U}+\varsigma)],\psi_{\gamma}^{2}D^{\gamma}u)\nonumber\\
 &+(D^{\gamma}[\partial_{y}^{-1}[\partial_{x}\mathfrak{U}]\partial_{y}u],\psi_{\gamma}^{2}D^{\gamma}u)
  +(D^{\gamma}[h\partial_{y}^{-1}[\partial_{x}u]\partial_{y}u],\psi_{\gamma}^{2}D^{\gamma}u)
  +(D^{\gamma}[\kappa\partial_{x}^{2}\mathbb{U}],\psi_{\gamma}^{2}D^{\gamma}u)\nonumber\\
 &+(D^{\gamma}F,\psi_{\gamma}^{2}D^{\gamma}u)
  +\left(D^{\gamma}[\eta\partial_{x}^{2}u]-\eta\partial_{x}^{2}D^{\gamma}u,\psi_{\gamma}^{2}D^{\gamma}u\right)
  -2(\psi_{\gamma}\partial_{y}\psi_{\gamma}\partial_{y}D^{\gamma}u,D^{\gamma}u)\nonumber\\
 &-(\partial_{y}D^{\gamma}u,\psi_{\gamma}^{2}D^{\gamma}u)_{L_{x}^{2}}|_{y=0}
:=\sum_{i=1}^{10}\mathcal{I}_{i}.\nonumber
\end{align}
We shall estimate $\mathcal{I}_{i}(i=1,...,10)$ term by term.

\noindent{\bf \underline{Estimate of~$\mathcal{I}_{1}$.}}
Notice that
\begin{align*}
\mathcal{I}_{1}
=-\sum_{0\leq\beta\leq\gamma}\binom{\gamma}{\beta}(D^{\beta}\mathfrak{U}D^{\gamma-\beta}[\partial_{x}u],\psi_{\gamma}^{2} D^{\gamma}u).
\end{align*}

\noindent{\bf Case 1:} $|\beta|=0$. By using integration by parts and Lemma \ref{Linfty} one can obtain
\begin{align*}
-(D^{\beta}\mathfrak{U}D^{\gamma-\beta}[\partial_{x}u],\psi_{\gamma}^{2}D^{\gamma}u)
=&\frac{1}{2}(\partial_{x}\mathfrak{U}D^{\gamma}u,\psi_{\gamma}^{2}D^{\gamma}u)\\
\leq& \frac{1}{2}\|\partial_{x}\mathfrak{U}\|_{L^{\infty}(\Omega)}\|D^{\gamma}u\|_{L_{\psi_{\gamma}}^{2}(\Omega)}^{2}\\
\leq& C_{*}\|u\|_{\mathcal{H}_{\psi,\varphi}^{s}(\Omega)}^{2}.
\end{align*}

\noindent{\bf Case 2:} $1\leq|\beta|\leq s-3$. By noticing $\beta\leq\gamma$, with the help of Corollary \ref{ncC2} one can get
\begin{align*}
&-(D^{\beta}\mathfrak{U}D^{\gamma-\beta}[\partial_{x}u],\psi_{\gamma}^{2}D^{\gamma}u)\\
\leq& \|D^{\beta}\mathfrak{U}\|_{L^{\infty}(\Omega)}
      \|D^{\gamma-\beta}[\partial_{x}u]\|_{L_{\psi_{\gamma}}^{2}(\Omega)}
      \|D^{\gamma}u\|_{L_{\psi_{\gamma}}^{2}(\Omega)}\\
\leq& C\|D^{\beta}\mathfrak{U}\|_{L^{\infty}(\Omega)}
       \left(\|u\|_{\hat{H}_{\psi}^{s}(\Omega)}
             +\|\partial_{x}^{s}u\|_{L_{\hat{\psi}}^{2}(\Omega)}
             +\|\partial_{y}\partial_{x}^{s}u\|_{L_{\hat{\psi}}^{2}(\Omega)}\right)
       \|\sqrt{\y}D^{\gamma}u\|_{L_{\psi_{\gamma}}^{2}(\Omega)}\\
\leq& C_{*}\|D^{\beta}\mathfrak{U}\|_{L^{\infty}(\Omega)}
       (\|u\|_{\mathcal{H}_{\psi,\varphi}^{s}(\Omega)}+\|\mathcal{P}\mathcal{U}_{s}\|_{L_{\varphi}^{2}(\Omega)})
       \|\sqrt{\y}D^{\gamma}u\|_{L_{\psi_{\gamma}}^{2}(\Omega)}\\
\leq& C_{*}\|u\|_{\mathcal{H}_{\psi,\varphi}^{s}(\Omega)}^{2}
     +C_{*}\|u\|_{\hat{H}_{\psi,\y}^{s}(\Omega)}^{2}
     +\frac{2}{3(s+1)(s+2)}\|\mathcal{P}\mathcal{U}_{s}\|_{L_{\varphi}^{2}(\Omega)}^{2}.
\end{align*}

\noindent{\bf Case 3:} $s-2\leq|\beta|\leq s+1$. By the Sobolev imbedding, we have
\begin{align*}
-(D^{\beta}\mathfrak{U}D^{\gamma-\beta}[\partial_{x}u],\psi_{\gamma}^{2} D^{\gamma}u)
=&-(D^{\beta}\mathbb{U}D^{\gamma-\beta}[\partial_{x}u],\psi_{\gamma}^{2} D^{\gamma}u)
  -(D^{\beta}u_{0}D^{\gamma-\beta}[\partial_{x}u],\psi_{\gamma}^{2} D^{\gamma}u)\\
\leq& \|D^{\beta}\mathbb{U}\|_{L_{\psi_{\gamma}}^{2}(\Omega)}
      \|D^{\gamma-\beta}[\partial_{x}u]\|_{L^{\infty}(\Omega)}
      \|D^{\gamma}u\|_{L_{\psi_{\gamma}}^{2}(\Omega)}\\
    &+\|D^{\beta}u_{0}\|_{L^{\infty}(\Omega)}
      \|D^{\gamma-\beta}[\partial_{x}u]\|_{L_{\psi_{\gamma}}^{2}(\Omega)}
      \|D^{\gamma}u\|_{L_{\psi_{\gamma}}^{2}(\Omega)}\\
\leq& C_{*}\|u\|_{\mathcal{H}_{\psi,\varphi}^{s}(\Omega)}^{2}.
\end{align*}

In conclusion, one has
\begin{align*}
\mathcal{I}_{1}
\leq C_{*}\|u\|_{\mathcal{H}_{\psi,\varphi}^{s}(\Omega)}^{2}
    +C_{*}\|u\|_{\hat{H}_{\psi,\y}^{s}(\Omega)}^{2}
    +\frac{2}{3(s+1)(s+2)}\|\mathcal{P}\mathcal{U}_{s}\|_{L_{\varphi}^{2}(\Omega)}^{2}.
\end{align*}

\noindent{\bf \underline{Estimate of~$\mathcal{I}_{2}$.}}
Notice that
\begin{align*}
\mathcal{I}_{2}
=-\sum_{0\leq\beta\leq\gamma}\binom{\gamma}{\beta}(D^{\beta}[\partial_{x}\mathfrak{U}]D^{\gamma-\beta}u,\psi_{\gamma}^{2} D^{\gamma}u).
\end{align*}

\noindent{\bf Case 1:} $0\leq|\beta|\leq s-3$. By using the Lemma \ref{Linfty} one can obtain
\begin{align*}
-(D^{\beta}[\partial_{x}\mathfrak{U}]D^{\gamma-\beta}u,\psi_{\gamma}^{2}D^{\gamma}u)
\leq& \|D^{\beta}[\partial_{x}\mathfrak{U}]\|_{L^{\infty}(\Omega)}
      \|D^{\gamma-\beta}u\|_{L_{\psi_{\gamma}}^{2}(\Omega)}
      \|D^{\gamma}u\|_{L_{\psi_{\gamma}}^{2}(\Omega)}\\
\leq& C_{*}\|u\|_{\mathcal{H}_{\psi,\varphi}^{s}(\Omega)}^{2}.
\end{align*}

\noindent{\bf Case 2:} $s-2\leq|\beta|\leq s+1$. Noticing $\beta\leq \gamma$, by the Lemma \ref{Linfty} and Corollary \ref{ncC2} one has
\begin{align*}
&-(D^{\beta}[\partial_{x}\mathfrak{U}]D^{\gamma-\beta}u,\psi_{\gamma}^{2}D^{\gamma}u)\\
=&-(D^{\beta}[\partial_{x}\mathbb{U}]D^{\gamma-\beta}u,\psi_{\gamma}^{2}D^{\gamma}u)
  -(D^{\beta}[\partial_{x}u_{0}]D^{\gamma-\beta}u,\psi_{\gamma}^{2}D^{\gamma}u)\\
\leq& \|D^{\beta}[\partial_{x}\mathbb{U}]\|_{L_{\hat{\psi}}^{2}(\Omega)}\|\sqrt{\y}D^{\gamma-\beta}u\|_{L^{\infty}(\Omega)}
      \|D^{\gamma}u\|_{L_{\psi_{\gamma}}^{2}(\Omega)}
     +\|D^{\beta}[\partial_{x}u_{0}]\|_{L^{\infty}(\Omega)}\|D^{\gamma-\beta}u\|_{L_{\psi_{\gamma}}^{2}(\Omega)}
      \|D^{\gamma}u\|_{L_{\psi_{\gamma}}^{2}(\Omega)}\\
\leq& C\left(\|\mathbb{U}\|_{\hat{H}_{\psi}^{s}(\Omega)}
        +\|\partial_{y}\mathbb{U}\|_{\hat{H}_{\psi}^{s}(\Omega)}
        +\|\partial_{x}^{s}\mathbb{U}\|_{L_{\hat{\psi}}^{2}(\Omega)}
        +\|\partial_{y}\partial_{x}^{s}\mathbb{U}\|_{L_{\hat{\psi}}^{2}(\Omega)}
        +\|\partial_{y}^{2}\partial_{x}^{s}\mathbb{U}\|_{L_{\hat{\psi}}^{2}(\Omega)}+1\right)
       \|u\|_{\mathcal{H}_{\psi,\varphi}^{s}(\Omega)}^{2}\\
\leq& C_{*}\left(\|\mathbb{U}\|_{\mathcal{H}_{\psi,\varphi}^{s}(\Omega)}
                +\|\partial_{y}\mathbb{U}\|_{\hat{H}_{\psi}^{s}(\Omega)}
                +\|\partial_{y}\mathbf{U}_{s}\|_{L_{\varphi}^{2}(\Omega)}
                +\|\mathcal{P}\mathbf{U}_{s}\|_{L_{\varphi}^{2}(\Omega)}+1\right)
     \|u\|_{\mathcal{H}_{\psi,\varphi}^{s}(\Omega)}^{2}.
\end{align*}

To sum up, one has
\begin{align*}
\mathcal{I}_{2}
\leq C_{*}\left(\|\mathbb{U}\|_{\mathcal{H}_{\psi,\varphi}^{s}(\Omega)}
                +\|\partial_{y}\mathbb{U}\|_{\hat{H}_{\psi}^{s}(\Omega)}
                +\|\partial_{y}\mathbf{U}_{s}\|_{L_{\varphi}^{2}(\Omega)}
                +\|\mathcal{P}\mathbf{U}_{s}\|_{L_{\varphi}^{2}(\Omega)}+1\right)
          \|u\|_{\mathcal{H}_{\psi,\varphi}^{s}(\Omega)}^{2}.
\end{align*}

\noindent{\bf \underline{Estimate of~$\mathcal{I}_{3}$.}}
It is easy to see that
\begin{align*}
\mathcal{I}_{3}
=\sum_{0\leq\beta\leq\gamma}\binom{\gamma}{\beta}(D^{\beta}[\partial_{y}\mathfrak{U}+\varsigma]D^{\gamma-\beta}\partial_{y}^{-1}[\partial_{x}u],\psi_{\gamma}^{2}D^{\gamma}u).
\end{align*}

\noindent{\bf Case 1:} $|\beta|\geq s$.
Noticing $|\partial_{y}\psi_{\gamma}|\leq C\psi_{\gamma}$ and $\gamma-\beta=(0,0)$ or $(0,1)$, by using integration by parts, Lemma \ref{Linfty} and Sobolev imbedding, one has
\begin{align*}
&(D^{\beta}[\partial_{y}\mathfrak{U}+\varsigma]D^{\gamma-\beta}\partial_{y}^{-1}[\partial_{x}u],\psi_{\gamma}^{2}D^{\gamma}u)\\
=& (D^{\beta}[\partial_{y}\mathbb{U}]D^{\gamma-\beta}\partial_{y}^{-1}[\partial_{x}u],\psi_{\gamma}^{2} D^{\gamma}u)
  +(D^{\beta}[\partial_{y}u_{0}]D^{\gamma-\beta}\partial_{y}^{-1}[\partial_{x}u],\psi_{\gamma}^{2} D^{\gamma}u)
  +(D^{\beta}\varsigma D^{\gamma-\beta}\partial_{y}^{-1}[\partial_{x}u],\psi_{\gamma}^{2} D^{\gamma}u)\\
=&-(D^{\beta}\mathbb{U}D^{\gamma-\beta}[\partial_{x}u],\psi_{\gamma}^{2} D^{\gamma}u)
  -(D^{\beta}\mathbb{U}D^{\gamma-\beta}\partial_{y}^{-1}[\partial_{x}u],\psi_{\gamma}^{2}\partial_{y}D^{\gamma}u)
  -2(D^{\beta}\mathbb{U}D^{\gamma-\beta}\partial_{y}^{-1}[\partial_{x}u],\psi_{\gamma}\partial_{y}\psi_{\gamma}D^{\gamma}u)\\
 &+(D^{\beta}[\partial_{y}u_{0}]D^{\gamma-\beta}\partial_{y}^{-1}[\partial_{x}u],\psi_{\gamma}^{2} D^{\gamma}u)
  +(D^{\beta}\varsigma D^{\gamma-\beta}\partial_{y}^{-1}[\partial_{x}u],\psi_{\gamma}^{2}D^{\gamma}u)\\
\leq& C\|D^{\beta}\mathbb{U}\|_{L_{\psi_{\gamma}}^{2}(\Omega)}
       \|D^{\gamma-\beta}[\partial_{x}u]\|_{L^{\infty}(\Omega)}
       \|D^{\gamma}u\|_{L_{\psi_{\gamma}}^{2}(\Omega)}
     +C\|D^{\beta}\mathbb{U}\|_{L_{\psi_{\gamma}}^{2}(\Omega)}
       \|\psi_{\gamma}D^{\gamma-\beta}[\partial_{x}u]\|_{L_{x}^{\infty}L_{y}^{2}(\Omega)}
       \|\partial_{y}D^{\gamma}u\|_{L_{\psi_{\gamma}}^{2}(\Omega)}\\
    &+\|D^{\beta}\mathbb{U}\|_{L_{\psi_{\gamma}}^{2}(\Omega)}
      \|\psi_{\gamma}D^{\gamma-\beta}[\partial_{x}u]\|_{L_{x}^{\infty}L_{y}^{2}(\Omega)}
      \|D^{\gamma}u\|_{L_{\psi_{\gamma}}^{2}(\Omega)}
     +\|D^{\beta}\partial_{y}u_{0}\|_{L_{\psi_{\gamma}}^{2}(\Omega)}
      \|\psi_{\gamma}D^{\gamma-\beta}[\partial_{x}u]\|_{L_{x}^{\infty}L_{y}^{2}(\Omega)}
      \|D^{\gamma}u\|_{L_{\psi_{\gamma}}^{2}(\Omega)}\\
    &+\|D^{\beta}\varsigma\|_{L_{\psi_{\gamma}}^{2}(\mathbb{R}_{+})}
      \|D^{\gamma-\beta}[\partial_{x}u]\|_{L_{\psi_{\gamma}}^{2}(\Omega)}
      \|D^{\gamma}u\|_{L_{\psi_{\gamma}}^{2}(\Omega)}\\
\leq& C_{*}\|u\|_{\mathcal{H}_{\psi,\varphi}^{s}(\Omega)}^{2}
     +\frac{1}{16}\|\partial_{y}D^{\gamma}u\|_{L_{\psi_{\gamma}}^{2}(\Omega)}^{2}.
\end{align*}

\noindent{\bf Case 2:} $s-1\geq|\beta|\geq s-2$.
By Sobolev imbedding one has
\begin{align*}
&(D^{\beta}[\partial_{y}\mathfrak{U}+\varsigma]D^{\gamma-\beta}\partial_{y}^{-1}[\partial_{x}u],\psi_{\gamma}^{2}D^{\gamma}u)\\
=& (D^{\beta}[\partial_{y}\mathfrak{U}]D^{\gamma-\beta}\partial_{y}^{-1}[\partial_{x}u],\psi_{\gamma}^{2}D^{\gamma}u)
  +(D^{\beta}\varsigma D^{\gamma-\beta}\partial_{y}^{-1}[\partial_{x}u],\psi_{\gamma}^{2}D^{\gamma}u)\\
\leq& \|D^{\gamma-\beta}\partial_{y}^{-1}[\partial_{x}u]\|_{L^{\infty}(\Omega)}
      \|D^{\beta}[\partial_{y}\mathfrak{U}]\|_{L_{\psi}^{2}(\Omega)}
      \|D^{\gamma}u\|_{L_{\psi_{\gamma}}^{2}(\Omega)}\\
    &+\|D^{\gamma-\beta}\partial_{y}^{-1}[\partial_{x}u]\|_{L_{x}^{2}L_{y}^{\infty}(\Omega)}
      \|D^{\beta}\varsigma\|_{L_{\psi}^{2}(\mathbb{R}_{+})}\|D^{\gamma}u\|_{L_{\psi_{\gamma}}^{2}(\Omega)}\\
\leq& C_{*}\|u\|_{\mathcal{H}_{\psi,\varphi}^{s}(\Omega)}^{2}.
\end{align*}

\noindent{\bf Case 3: $s-3\geq|\beta|\geq 0$}.
If $\gamma_{2}-\beta_{2}=0$, by noticing $\gamma_{1}<s$ one can deduce from Corollary \ref{ncC2} that
\begin{align*}
&(D^{\beta}[\partial_{y}\mathfrak{U}+\varsigma]D^{\gamma-\beta}\partial_{y}^{-1}[\partial_{x}u],\psi_{\gamma}^{2}D^{\gamma}u)\\
=& (D^{\beta}[\partial_{y}\mathfrak{U}]D^{\gamma-\beta}\partial_{y}^{-1}[\partial_{x}u],\psi_{\gamma}^{2}D^{\gamma}u)
  +(D^{\beta}\varsigma D^{\gamma-\beta}\partial_{y}^{-1}[\partial_{x}u],\psi_{\gamma}^{2}D^{\gamma}u)\\
\leq& \|\psi_{\gamma}D^{\beta}[\partial_{y}\mathfrak{U}]\|_{L_{x}^{\infty}L_{y}^{2}(\Omega)}
      \|D^{\gamma-\beta}\partial_{y}^{-1}[\partial_{x}u]\|_{L_{x}^{2}L_{y}^{\infty}(\Omega)}
      \|D^{\gamma}u\|_{L_{\psi_{\gamma}}^{2}(\Omega)}\\
    &+\|D^{\beta}\varsigma\|_{L_{\psi}^{2}(\mathbb{R}_{+})}
      \|D^{\gamma-\beta}\partial_{y}^{-1}[\partial_{x}u]\|_{L_{x}^{2}L_{y}^{\infty}(\Omega)}
      \|D^{\gamma}u\|_{L_{\psi_{\gamma}}^{2}(\Omega)}\\
\leq& C_{*}\left(\|\partial_{x}^{s}u\|_{L_{\varphi}^{2}(\Omega)}+\|u\|_{\hat{H}_{\psi}^{s}(\Omega)}\right)
       \|D^{\gamma}u\|_{L_{\psi_{\gamma}}^{2}(\Omega)}\\
\leq& C_{*}\left(\|\mathcal{P}\mathcal{U}_{s}\|_{L_{\varphi}^{2}(\Omega)}+\|u\|_{\hat{H}_{\psi}^{s}(\Omega)}\right)
       \|D^{\gamma}u\|_{L_{\psi_{\gamma}}^{2}(\Omega)}\\
\leq& \frac{2}{3(s+1)(s+2)}\|\mathcal{P}\mathcal{U}_{s}\|_{L_{\varphi}^{2}(\Omega)}^{2}
     +C_{*}\|u\|_{\mathcal{H}_{\psi,\varphi}^{s}(\Omega)}^{2},
\end{align*}
else, by using Corollary \ref{ncC2} one can obtain
\begin{align*}
&(D^{\beta}[\partial_{y}\mathfrak{U}+\varsigma]D^{\gamma-\beta}\partial_{y}^{-1}[\partial_{x}u],\psi_{\gamma}^{2}D^{\gamma}u)\\
=& (D^{\beta}[\partial_{y}\mathfrak{U}]D^{\gamma-\beta}\partial_{y}^{-1}[\partial_{x}u],\psi_{\gamma}^{2}D^{\gamma}u)
  +(D^{\beta}\varsigma D^{\gamma-\beta}\partial_{y}^{-1}[\partial_{x}u],\psi_{\gamma}^{2}D^{\gamma}u)\\
\leq& \|\sqrt{\y}D^{\beta}[\partial_{y}\mathfrak{U}]\|_{L^{\infty}(\Omega)}
      \|D^{\gamma-\beta}\partial_{y}^{-1}[\partial_{x}u]\|_{L_{\hat{\psi}_{\gamma}}^{2}(\Omega)}
      \|D^{\gamma}u\|_{L_{\psi_{\gamma}}^{2}(\Omega)}\\
    &+\|\sqrt{\y}D^{\beta}\varsigma\|_{L^{\infty}(\mathbb{R}_{+})}
      \|D^{\gamma-\beta}\partial_{y}^{-1}[\partial_{x}u]\|_{L_{\hat{\psi}_{\gamma}}^{2}(\Omega)}
      \|D^{\gamma}u\|_{L_{\psi_{\gamma}}^{2}(\Omega)}\\
\leq& C_{*}\left(\|\partial_{y}\partial_{x}^{s}u\|_{L_{\hat{\psi}}^{2}(\Omega)}
            +\|\partial_{x}^{s}u\|_{L_{\hat{\psi}}^{2}(\Omega)}
            +\|u\|_{\hat{H}_{\psi}^{s}(\Omega)}\right)
       \|D^{\gamma}u\|_{L_{\psi_{\gamma}}^{2}(\Omega)}\\
\leq& C_{*}\left(\|\mathcal{P}\mathcal{U}_{s}\|_{L_{\varphi}^{2}(\Omega)}+\|u\|_{\mathcal{H}_{\psi,\varphi}^{s}(\Omega)}\right)
       \|D^{\gamma}u\|_{L_{\psi_{\gamma}}^{2}(\Omega)}\\
\leq& \frac{2}{3(s+1)(s+2)}\|\mathcal{P}\mathcal{U}_{s}\|_{L_{\varphi}^{2}(\Omega)}^{2}
     +C_{*}\|u\|_{\mathcal{H}_{\psi,\varphi}^{s}(\Omega)}^{2}.
\end{align*}

In summary, one can obtain
\begin{align*}
\mathcal{I}_{3}
\leq& C_{*}\|u\|_{\mathcal{H}_{\psi,\varphi}^{s}(\Omega)}^{2}
     +\frac{1}{16}\|\partial_{y}D^{\gamma}u\|_{L_{\psi_{\gamma}}^{2}(\Omega)}^{2}
     +\frac{2}{3(s+1)(s+2)}\|\mathcal{P}\mathcal{U}_{s}\|_{L_{\varphi}^{2}(\Omega)}^{2}.
\end{align*}

\noindent{\bf\underline{Estimate of~$\mathcal{I}_{4}$.}}
\begin{align*}
\mathcal{I}_{4}
=& (D^{\gamma}[\partial_{y}^{-1}[\partial_{x}\mathfrak{U}]\partial_{y}u],\psi_{\gamma}^{2} D^{\gamma}u)\\
=&\sum_{0\leq\beta\leq\gamma}\binom{\gamma}{\beta}(D^{\beta}\partial_{y}^{-1}[\partial_{x}\mathfrak{U}]D^{\gamma-\beta}[\partial_{y}u],\psi_{\gamma}^{2}D^{\gamma}u).
\end{align*}

\noindent{\bf Case 1:} $|\beta|\geq s$. Since $\gamma_{1}<s$, we know $\beta_{2}\geq1$. With the help of Corollary \ref{ncC2}, we have
\begin{align*}
&(D^{\beta}\partial_{y}^{-1}[\partial_{x}\mathfrak{U}]D^{\gamma-\beta}[\partial_{y}u],\psi_{\gamma}^{2}D^{\gamma}u)\\
=& (\partial_{x}^{\beta_{1}+1}\partial_{y}^{\beta_{2}-1}\mathbb{U}D^{\gamma-\beta}[\partial_{y}u],\psi_{\gamma}^{2}D^{\gamma}u)
  +(\partial_{x}^{\beta_{1}+1}\partial_{y}^{\beta_{2}-1}u_{0}D^{\gamma-\beta}[\partial_{y}u],\psi_{\gamma}^{2}D^{\gamma}u)\\
\leq& \|\partial_{x}^{\beta_{1}+1}\partial_{y}^{\beta_{2}-1}\mathbb{U}\|_{L_{\hat{\psi}_{\gamma}}^{2}(\Omega)}
      \|\sqrt{\y}D^{\gamma-\beta}[\partial_{y}u]\|_{L^{\infty}(\Omega)}
      \|D^{\gamma}u\|_{L_{\psi_{\gamma}}^{2}(\Omega)}\\
    &+\|\partial_{x}^{\beta_{1}+1}\partial_{y}^{\beta_{2}-1}u_{0}\|_{L^{\infty}(\Omega)}
      \|D^{\gamma-\beta}[\partial_{y}u]\|_{L_{\psi_{\gamma}}^{2}(\Omega)}
      \|D^{\gamma}u\|_{L_{\psi_{\gamma}}^{2}(\Omega)}\\
\leq& C\left(\|\mathbb{U}\|_{\hat{H}_{\psi}^{s}(\Omega)}
        +\|\partial_{x}^{s}\mathbb{U}\|_{L_{\hat{\psi}}^{2}(\Omega)}
        +\|\partial_{y}\partial_{x}^{s}\mathbb{U}\|_{L_{\hat{\psi}}^{2}(\Omega)}+1\right)
       \|u\|_{\mathcal{H}_{\psi,\varphi}^{s}(\Omega)}^{2}\\
\leq& C_{*}\left(\|\mathcal{P}\mathbf{U}_{s}\|_{L_{\varphi}^{2}(\Omega)}+\|\mathbb{U}\|_{\mathcal{H}_{\psi,\varphi}^{s}(\Omega)}+1\right)
       \|u\|_{\mathcal{H}_{\psi,\varphi}^{s}(\Omega)}^{2}.
\end{align*}

\noindent{\bf Case 2:} $s-2\leq|\beta|\leq s-1$. By using Sobolev imbedding and Corollary \ref{ncC2}, we have
\begin{align*}
&(D^{\beta}\partial_{y}^{-1}[\partial_{x}\mathfrak{U}]D^{\gamma-\beta}[\partial_{y}u],\psi_{\gamma}^{2}D^{\gamma}u)\\
=& (D^{\beta}\partial_{y}^{-1}[\partial_{x}\mathbb{U}]D^{\gamma-\beta}[\partial_{y}u],\psi_{\gamma}^{2}D^{\gamma}u)
  +(D^{\beta}\partial_{y}^{-1}[\partial_{x}u_{0}]D^{\gamma-\beta}[\partial_{y}u],\psi_{\gamma}^{2}D^{\gamma}u)\\
\leq& \|D^{\beta}\partial_{y}^{-1}[\partial_{x}\mathbb{U}]\|_{L_{x}^{2}L_{y}^{\infty}(\Omega)}
      \|\psi_{\gamma}D^{\gamma-\beta}[\partial_{y}u]\|_{L_{x}^{\infty}L_{y}^{2}(\Omega)}
      \|D^{\gamma}u\|_{L_{\psi_{\gamma}}^{2}(\Omega)}\\
    &+\|\y^{-1}D^{\beta}\partial_{y}^{-1}[\partial_{x}u_{0}]\|_{L_{x}^{2}L_{y}^{\infty}(\Omega)}
      \|\psi_{\gamma}\sqrt{\y}D^{\gamma-\beta}[\partial_{y}u]\|_{L_{x}^{\infty}L_{y}^{2}(\Omega)}
      \|\sqrt{\y}D^{\gamma}u\|_{L_{\psi_{\gamma}}^{2}(\Omega)}\\
\leq& C\left(\|\mathbb{U}\|_{\hat{H}_{\psi}^{s}(\Omega)}
        +\|\partial_{x}^{s}\mathbb{U}\|_{L_{\hat{\psi}}^{2}(\Omega)}\right)
       \|u\|_{\mathcal{H}_{\psi,\varphi}^{s}(\Omega)}^{2}
     +C\|u\|_{\hat{H}_{\psi,\y}^{s}(\Omega)}^{2}\\
\leq& C_{*}\left(\|\mathbb{U}\|_{\hat{H}_{\psi}^{s}(\Omega)}+\|\mathcal{P}\mathbf{U}_{s}\|_{L_{\varphi}^{2}(\Omega)}\right)
       \|u\|_{\mathcal{H}_{\psi,\varphi}^{s}(\Omega)}^{2}
     +C\|u\|_{\hat{H}_{\psi,\y}^{s}(\Omega)}^{2}.
\end{align*}

\noindent{\bf Case 3:} $0<|\beta|\leq s-3$. Noticing that $\psi_{\gamma}=\psi$ when $|\gamma|\leq s$ and $\psi_{\gamma}=\hat{\psi}$ when $|\gamma|=s+1$, we have
\begin{align*}
\|\sqrt{\y}D^{\gamma-\beta}[\partial_{y}u]\|_{L_{\psi_{\gamma}}^{2}(\Omega)}
\leq \|\sqrt{\y}u\|_{\hat{H}_{\psi}^{s}(\Omega)}+\|\partial_{y}u\|_{\hat{H}_{\psi}^{s}(\Omega)},
\end{align*}
thus,
\begin{align*}
&(D^{\beta}\partial_{y}^{-1}[\partial_{x}\mathfrak{U}]D^{\gamma-\beta}[\partial_{y}u],\psi_{\gamma}^{2}D^{\gamma}u)\\
=& (D^{\beta}\partial_{y}^{-1}[\partial_{x}\mathbb{U}]D^{\gamma-\beta}[\partial_{y}u],\psi_{\gamma}^{2}D^{\gamma}u)
  +(D^{\beta}\partial_{y}^{-1}[\partial_{x}u_{0}]D^{\gamma-\beta}[\partial_{y}u],\psi_{\gamma}^{2}D^{\gamma}u)\\
\leq& \|D^{\beta}\partial_{y}^{-1}[\partial_{x}\mathbb{U}]\|_{L^{\infty}(\Omega)}
      \|D^{\gamma-\beta}[\partial_{y}u]\|_{L_{\psi_{\gamma}}^{2}(\Omega)}
      \|D^{\gamma}u\|_{L_{\psi_{\gamma}}^{2}(\Omega)}\\
    &+\|\y^{-1}D^{\beta}\partial_{y}^{-1}[\partial_{x}u_{0}]\|_{L^{\infty}(\Omega)}
      \|\sqrt{\y}D^{\gamma-\beta}[\partial_{y}u]\|_{L_{\psi_{\gamma}}^{2}(\Omega)}
      \|\sqrt{\y}D^{\gamma}u\|_{L_{\psi_{\gamma}}^{2}(\Omega)}\\
\leq& C\|\mathbb{U}\|_{\hat{H}_{\psi}^{s}(\Omega)}
       \left(\|u\|_{\hat{H}_{\psi}^{s}(\Omega)}+\|\partial_{y}u\|_{\hat{H}_{\psi}^{s}(\Omega)}\right)
       \|D^{\gamma}u\|_{L_{\psi_{\gamma}}^{2}(\Omega)}\\
    &+C\left(\|\sqrt{\y}u\|_{\hat{H}_{\psi}^{s}(\Omega)}+\|\partial_{y}u\|_{\hat{H}_{\psi}^{s}(\Omega)}\right)
       \|\sqrt{\y}D^{\gamma}u\|_{L_{\psi_{\gamma}}^{2}(\Omega)}\\
\leq& C_{*}\|u\|_{\mathcal{H}_{\psi,\varphi}^{s}(\Omega)}^{2}
     +C\|u\|_{\hat{H}_{\psi,\y}^{s}(\Omega)}^{2}
     +\frac{1}{8(s+1)(s+2)}\|\partial_{y}u\|_{\hat{H}_{\psi}^{s}(\Omega)}^{2}.
\end{align*}

\noindent{\bf Case 4:} $|\beta|=0$. By using integration by parts, one has
\begin{align*}
&(\partial_{y}^{-1}[\partial_{x}\mathfrak{U}]D^{\gamma}[\partial_{y}u],\psi_{\gamma}^{2}D^{\gamma}u)\\
=&-\frac{1}{2}(\partial_{x}\mathfrak{U}D^{\gamma}u,\psi_{\gamma}^{2}D^{\gamma}u)
  -(\partial_{y}^{-1}[\partial_{x}\mathbb{U}]D^{\gamma}u,\psi_{\gamma}\partial_{y}\psi_{\gamma}D^{\gamma}u)
  -(\partial_{y}^{-1}[\partial_{x}u_{0}]D^{\gamma}u,\psi_{\gamma}\partial_{y}\psi_{\gamma}D^{\gamma}u)\\
\leq&~C\left(\|\partial_{x}\mathfrak{U}\|_{L^{\infty}(\Omega)}+\|\partial_{y}^{-1}[\partial_{x}\mathbb{U}]\|_{L^{\infty}(\Omega)}\right)
       \|D^{\gamma}u\|_{L_{\psi_{\gamma}}^{2}(\Omega)}^{2}
     +\|\y^{-1}\partial_{y}^{-1}[\partial_{x}u_{0}]\|_{L^{\infty}(\Omega)}
      \|\sqrt{\y}D^{\gamma}u\|_{L_{\psi_{\gamma}}^{2}(\Omega)}^{2}\\
\leq&~C_{*}\|D^{\gamma}u\|_{L_{\psi_{\gamma}}^{2}(\Omega)}^{2}
     +C\|\sqrt{\y}D^{\gamma}u\|_{L_{\psi_{\gamma}}^{2}(\Omega)}^{2}.
\end{align*}

To sum up, we have
\begin{align*}
\mathcal{I}_{4}
\leq C_{*}(1+\|\mathbb{U}\|_{\mathcal{H}_{\psi,\varphi}^{s}(\Omega)}+\|\mathcal{P}\mathbf{U}_{s}\|_{L_{\varphi}^{2}(\Omega)})
      \|u\|_{\mathcal{H}_{\psi,\varphi}^{s}(\Omega)}^{2}
    +C\|u\|_{\hat{H}_{\psi,\y}^{s}(\Omega)}^{2}
    +\frac{1}{8(s+1)(s+2)}\|\partial_{y}u\|_{\hat{H}_{\psi}^{s}(\Omega)}^{2}.
\end{align*}

\noindent{\bf\underline{Estimate of~$\mathcal{I}_{5}$.}}
\begin{align*}
\mathcal{I}_{5}
=& (D^{\gamma}[h\partial_{y}^{-1}[\partial_{x}u]\partial_{y}u],\psi_{\gamma}^{2} D^{\gamma}u)\\
=&\sum_{0\leq\beta\leq\gamma}\sum_{0\leq\alpha\leq\beta}\binom{\gamma}{\beta}\binom{\beta}{\alpha}
  (D^{\beta-\alpha}hD^{\alpha}\partial_{y}^{-1}[\partial_{x}u]D^{\gamma-\beta}[\partial_{y}u],\psi_{\gamma}^{2}D^{\gamma}u).
\end{align*}

\noindent{\bf Case 1:} $|\alpha|\geq s$. Since $\gamma_{1}<s$, we know $\alpha_{2}\geq1$. With the help of Corollary \ref{ncC2}, we have
\begin{align*}
(D^{\beta-\alpha}hD^{\alpha}\partial_{y}^{-1}[\partial_{x}u]D^{\gamma-\beta}[\partial_{y}u],\psi_{\gamma}^{2}D^{\gamma}u)
=& (D^{\beta-\alpha}h\partial_{x}^{\alpha_{1}+1}\partial_{y}^{\alpha_{2}-1}uD^{\gamma-\beta}[\partial_{y}u],\psi_{\gamma}^{2}D^{\gamma}u)\\
\leq& C\|\partial_{x}^{\alpha_{1}+1}\partial_{y}^{\alpha_{2}-1}u\|_{L_{\hat{\psi}_{\gamma}}^{2}(\Omega)}
       \|\sqrt{\y}D^{\gamma-\beta}[\partial_{y}u]\|_{L^{\infty}(\Omega)}
       \|D^{\gamma}u\|_{L_{\psi_{\gamma}}^{2}(\Omega)}\\
\leq& C\left(\|u\|_{\hat{H}_{\psi}^{s}(\Omega)}
        +\|\partial_{x}^{s}u\|_{L_{\hat{\psi}}^{2}(\Omega)}
        +\|\partial_{y}\partial_{x}^{s}u\|_{L_{\hat{\psi}}^{2}(\Omega)}\right)
       \|u\|_{\mathcal{H}_{\psi,\varphi}^{s}(\Omega)}^{2}\\
\leq& C_{*}\left(\|u\|_{\mathcal{H}_{\psi,\varphi}^{s}(\Omega)}+\|\mathcal{P}\mathcal{U}_{s}\|_{L_{\varphi}^{2}(\Omega)}\right)
       \|u\|_{\mathcal{H}_{\psi,\varphi}^{s}(\Omega)}^{2}\\
\leq& C_{*}\|u\|_{\mathcal{H}_{\psi,\varphi}^{s}(\Omega)}^{4}
     +C_{*}\|u\|_{\mathcal{H}_{\psi,\varphi}^{s}(\Omega)}^{2}
     +\frac{2}{3(s+1)(s+2)}\|\mathcal{P}\mathcal{U}_{s}\|_{L_{\varphi}^{2}(\Omega)}^{2}.
\end{align*}

\noindent{\bf Case 2:} $s-2\leq|\alpha|\leq s-1$. By using Lemma \ref{Linfty}, Sobolev imbedding and Corollary \ref{ncC2}, we have
\begin{align*}
(D^{\beta-\alpha}hD^{\alpha}\partial_{y}^{-1}[\partial_{x}u]D^{\gamma-\beta}[\partial_{y}u],\psi_{\gamma}^{2}D^{\gamma}u)
\leq& C\|D^{\alpha}\partial_{y}^{-1}[\partial_{x}u]\|_{L_{x}^{2}L_{y}^{\infty}(\Omega)}
       \|\psi_{\gamma}D^{\gamma-\beta}[\partial_{y}u]\|_{L_{x}^{\infty}L_{y}^{2}(\Omega)}
       \|D^{\gamma}u\|_{L_{\psi_{\gamma}}^{2}(\Omega)}\\
\leq& C\left(\|u\|_{\hat{H}_{\psi}^{s}(\Omega)}
        +\|\partial_{x}^{s}u\|_{L_{\hat{\psi}}^{2}(\Omega)}\right)
       \|u\|_{\mathcal{H}_{\psi,\varphi}^{s}(\Omega)}^{2}\\
\leq& C_{*}\left(\|u\|_{\hat{H}_{\psi}^{s}(\Omega)}+\|\mathcal{P}\mathcal{U}_{s}\|_{L_{\varphi}^{2}(\Omega)}\right)
       \|u\|_{\mathcal{H}_{\psi,\varphi}^{s}(\Omega)}^{2}\\
\leq& C_{*}\|u\|_{\mathcal{H}_{\psi,\varphi}^{s}(\Omega)}^{4}
     +C_{*}\|u\|_{\mathcal{H}_{\psi,\varphi}^{s}(\Omega)}^{2}
     +\frac{1}{(s+1)(s+2)}\|\mathcal{P}\mathcal{U}_{s}\|_{L_{\varphi}^{2}(\Omega)}^{2}.
\end{align*}

\noindent{\bf Case 3:} $0\leq|\alpha|\leq s-3$. We have
\begin{align*}
(D^{\beta-\alpha}hD^{\alpha}\partial_{y}^{-1}[\partial_{x}u]D^{\gamma-\beta}[\partial_{y}u],\psi_{\gamma}^{2}D^{\gamma}u)
\leq& \|D^{\alpha}\partial_{y}^{-1}[\partial_{x}u]\|_{L^{\infty}(\Omega)}
      \|D^{\gamma-\beta}[\partial_{y}u]\|_{L_{\psi_{\gamma}}^{2}(\Omega)}
      \|D^{\gamma}u\|_{L_{\psi_{\gamma}}^{2}(\Omega)}\\
\leq& C\|u\|_{\hat{H}_{\psi}^{s}(\Omega)}
       \left(\|u\|_{\hat{H}_{\psi}^{s}(\Omega)}+\|\partial_{y}u\|_{\hat{H}_{\psi}^{s}(\Omega)}\right)
       \|D^{\gamma}u\|_{L_{\psi_{\gamma}}^{2}(\Omega)}\\
\leq& C_{*}\|u\|_{\mathcal{H}_{\psi,\varphi}^{s}(\Omega)}^{4}
     +C_{*}\|u\|_{\mathcal{H}_{\psi,\varphi}^{s}(\Omega)}^{2}
     +\frac{1}{8(s+1)(s+2)}\|\partial_{y}u\|_{\hat{H}_{\psi}^{s}(\Omega)}^{2}.
\end{align*}

To sum up, we have
\begin{align*}
\mathcal{I}_{5}
\leq C_{*}\|u\|_{\mathcal{H}_{\psi,\varphi}^{s}(\Omega)}^{4}
     +C_{*}\|u\|_{\mathcal{H}_{\psi,\varphi}^{s}(\Omega)}^{2}
     +\frac{2}{3(s+1)(s+2)}\|\mathcal{P}\mathcal{U}_{s}\|_{L_{\varphi}^{2}(\Omega)}^{2}
     +\frac{1}{8(s+1)(s+2)}\|\partial_{y}u\|_{\hat{H}_{\psi}^{s}(\Omega)}^{2}.
\end{align*}

\noindent{\bf \underline{Estimate of~$\mathcal{I}_{6}$.}}

\noindent{\bf Case 1:} $|\gamma|=s+1$.
In this case, by using Corollary \ref{ncC2}, one has
\begin{align}\label{I6c1}
\mathcal{I}_{6}
=&~(\kappa\partial_{y}^{2}\partial_{x}^{s+1}\mathbb{U},\hat{\psi}^{2}D^{\gamma}u)
  +2(\partial_{y}\kappa\partial_{y}\partial_{x}^{s+1}\mathbb{U},\hat{\psi}^{2}D^{\gamma}u)
  +(\partial_{y}^{2}\kappa\partial_{x}^{s+1}\mathbb{U},\hat{\psi}^{2}D^{\gamma}u)\\
\leq&~C_{\theta}\bar{\kappa}\left(\|\y^{-\frac{1}{2}}\partial_{y}^{2}\partial_{x}^{s+1}\mathbb{U}\|_{L_{\hat{\psi}}^{2}(\Omega)}
                        +\|\y^{-\frac{1}{2}}\partial_{y}\partial_{x}^{s+1}\mathbb{U}\|_{L_{\hat{\psi}}^{2}(\Omega)}
                        +\|\y^{-\frac{1}{2}}\partial_{x}^{s+1}\mathbb{U}\|_{L_{\hat{\psi}}^{2}(\Omega)}\right)
                  \|\sqrt{\y}D^{\gamma}u\|_{L_{\hat{\psi}}^{2}(\Omega)}\nonumber\\
\leq&~\bar{\kappa}^{2}\left(\|\mathbf{U}_{s+1}\|_{L_{\hat{\varphi}}^{2}(\Omega)}^{2}
                     +\tau\|\mathcal{P}\mathbf{U}_{s+1}\|_{L_{\hat{\varphi}}^{2}(\Omega)}^{2}
                     +\tau\|\partial_{y}\mathbf{U}_{s+1}\|_{L_{\hat{\varphi}}^{2}(\Omega)}^{2}\right)
     +C_{*}C_{\tau}C_{\theta}\|\sqrt{\y}D^{\gamma}u\|_{L_{\hat{\psi}}^{2}(\Omega)}^{2}.\nonumber
\end{align}

\noindent{\bf Case 2:} $|\gamma|=s$. If $\gamma_{1}=s-1$, by Corollary \ref{ncC2} one has
\begin{align}\label{I6c21}
\mathcal{I}_{6}
=&~(\kappa\partial_{y}\partial_{x}^{s+1}\mathbb{U},\psi^{2}D^{\gamma}u)
  +(\partial_{y}\kappa\partial_{x}^{s+1}\mathbb{U},\psi^{2}D^{\gamma}u)\\
\leq&~C_{\theta}\bar{\kappa}\left(\|\partial_{y}\partial_{x}^{s+1}\mathbb{U}\|_{L_{\hat{\psi}}^{2}(\Omega)}
                                  +\|\partial_{x}^{s+1}\mathbb{U}\|_{L_{\hat{\psi}}^{2}(\Omega)}\right)
                            \|\sqrt{\y}D^{\gamma}u\|_{L_{\psi}^{2}(\Omega)}\nonumber\\
\leq&~C_{\theta}\bar{\kappa}
      \left(\|\mathbf{U}_{s+1}\|_{L_{\varphi}^{2}(\Omega)}+\|\mathcal{P}\mathbf{U}_{s+1}\|_{L_{\varphi}^{2}(\Omega)}\right)
      \|\sqrt{\y}D^{\gamma}u\|_{L_{\psi}^{2}(\Omega)}\nonumber\\
\leq&~\bar{\kappa}^{2}\tau\left(\|\sqrt{\y}\mathbf{U}_{s+1}\|_{L_{\hat{\varphi}}^{2}(\Omega)}^{2}
                            +\|\mathcal{P}\mathbf{U}_{s+1}\|_{L_{\hat{\varphi}}^{2}(\Omega)}^{2}\right)
     +C_{*}C_{\tau}C_{\theta}\|\sqrt{\y}D^{\gamma}u\|_{L_{\psi}^{2}(\Omega)}^{2}.\nonumber
\end{align}
Here we use the fact
\begin{align*}
\|\mathcal{P}\mathbf{U}_{s+1}\|_{L_{\varphi}^{2}(\Omega)}
\leq C\left(\|\sqrt{\y}\mathbf{U}_{s+1}\|_{L_{\hat{\varphi}}^{2}(\Omega)}^{2}+\|\mathcal{P}\mathbf{U}_{s+1}\|_{L_{\hat{\varphi}}^{2}(\Omega)}^{2}\right).
\end{align*}
If $\gamma_{1}=s-2$, the estimate of $\mathcal{I}_{6}$ is similar to that in Case 1, one has
\begin{align}\label{I6c22}
\mathcal{I}_{6}
=&~(\kappa\partial_{y}^{2}\partial_{x}^{s}\mathbb{U},\psi^{2}D^{\gamma}u)
  +2(\partial_{y}\kappa\partial_{y}\partial_{x}^{s}\mathbb{U},\psi^{2}D^{\gamma}u)
  +(\partial_{y}^{2}\kappa\partial_{x}^{s}\mathbb{U},\psi^{2}D^{\gamma}u)\\
\leq&~C_{\theta}\bar{\kappa}\left(\|\partial_{y}^{2}\partial_{x}^{s}\mathbb{U}\|_{L_{\hat{\psi}}^{2}(\Omega)}
                        +\|\partial_{y}\partial_{x}^{s}\mathbb{U}\|_{L_{\hat{\psi}}^{2}(\Omega)}
                        +\|\partial_{x}^{s}\mathbb{U}\|_{L_{\hat{\psi}}^{2}(\Omega)}\right)
                  \|\sqrt{\y}D^{\gamma}u\|_{L_{\psi}^{2}(\Omega)}\nonumber\\
\leq&~\bar{\kappa}^{2}\left(\|\mathbf{U}_{s}\|_{L_{\varphi}^{2}(\Omega)}^{2}
                +\tau\|\mathcal{P}\mathbf{U}_{s}\|_{L_{\varphi}^{2}(\Omega)}^{2}
                +\tau\|\partial_{y}\mathbf{U}_{s}\|_{L_{\varphi}^{2}(\Omega)}^{2}\right)
     +C_{*}C_{\tau}C_{\theta}\|\sqrt{\y}D^{\gamma}u\|_{L_{\psi}^{2}(\Omega)}^{2}.\nonumber
\end{align}
If $\gamma_{1}<s-2$, by applying integration by parts it follows
\begin{align}\label{I6c23}
\mathcal{I}_{6}
=&(\kappa\partial_{y}^{\gamma_{2}}\partial_{x}^{\gamma_{1}+2}\mathbb{U}, \psi^{2}D^{\gamma}u)
  +\sum_{0<\beta_{2}\leq\gamma_{2}}\binom{\gamma_{2}}{\beta_{2}}
   (\partial_{y}^{\beta_{2}}\kappa\partial_{y}^{\gamma_{2}-\beta_{2}}\partial_{x}^{\gamma_{1}+2}\mathbb{U},\psi^{2}D^{\gamma}u)\\
=&-(\kappa\partial_{y}^{\gamma_{2}-1}\partial_{x}^{\gamma_{1}+2}\mathbb{U}, \psi^{2} \partial_{y}D^{\gamma}u)
  -(\partial_{y}\kappa\partial_{y}^{\gamma_{2}-1}\partial_{x}^{\gamma_{1}+2}\mathbb{U}, \psi^{2}D^{\gamma}u)
  -(\kappa\partial_{y}^{\gamma_{2}-1}\partial_{x}^{\gamma_{1}+2}\mathbb{U}, \psi\partial_{y}\psi D^{\gamma}u)\nonumber\\
 &-\psi_{0}^{2}(\kappa\partial_{y}^{\gamma_{2}-1}\partial_{x}^{\gamma_{1}+2}\mathbb{U},\partial_{y}^{\gamma_{2}}\partial_{x}^{\gamma_{1}}u)_{L_{x}^{2}}|_{y=0}
  -\sum_{0<\beta_{2}\leq\gamma_{2}}\binom{\gamma_{2}}{\beta_{2}}
   (\partial_{y}^{\beta_{2}}\kappa\partial_{y}^{\gamma_{2}-\beta_{2}}\partial_{x}^{\gamma_{1}+1}\mathbb{U},\psi^{2}\partial_{x}D^{\gamma}u),\nonumber
\end{align}
where $\psi_{0}=\psi(t,0)$. Noticing $\gamma_{2}\geq3$ when $\gamma_{1}<s-2$, we have
\begin{align}\label{I6c231}
&-(\kappa\partial_{y}^{\gamma_{2}-1}\partial_{x}^{\gamma_{1}+2}\mathbb{U}, \psi^{2} \partial_{y}D^{\gamma}u)
 -(\partial_{y}\kappa\partial_{y}^{\gamma_{2}-1}\partial_{x}^{\gamma_{1}+2}\mathbb{U}, \psi^{2}D^{\gamma}u)
 -(\kappa\partial_{y}^{\gamma_{2}-1}\partial_{x}^{\gamma_{1}+2}\mathbb{U}, \psi\partial_{y}\psi D^{\gamma}u)\\
\leq& \bar{\kappa}\|\partial_{y}\mathbb{U}\|_{\hat{H}_{\psi}^{s}(\Omega)}\|\partial_{y}D^{\gamma}u\|_{L_{\psi}^{2}(\Omega)}
     +C\bar{\kappa}\|\partial_{y}\mathbb{U}\|_{\hat{H}_{\psi}^{s}(\Omega)}\|D^{\gamma}u\|_{L_{\psi}^{2}(\Omega)}\nonumber\\
\leq&  \frac{15\tau\bar{\kappa}^{2}}{16}\|\partial_{y}\mathbb{U}\|_{\hat{H}_{\psi}^{s}(\Omega)}
     +\frac{5}{16}\|\partial_{y}D^{\gamma}u\|_{L_{\psi}^{2}(\Omega)}^{2}
     +C_{\tau}\|D^{\gamma}u\|_{L_{\psi}^{2}(\Omega)}^{2}.\nonumber
\end{align}
The boundary term in \eqref{I6c23} will be estimate in different ways when $\gamma_{2}$ is odd or even, by Lemma \ref{Bc} and \ref{bc} respectively. If $\gamma_{2}$ is odd, thanks to Lemma \ref{Bc} and Sobolev imbedding, we have
\begin{align}\label{I6c232}
&-\psi_{0}^{2}(\kappa\partial_{y}^{\gamma_{2}-1}\partial_{x}^{\gamma_{1}+2}\mathbb{U},\partial_{y}^{\gamma_{2}}\partial_{x}^{\gamma_{1}}u)_{L_{x}^{2}}|_{y=0}\\
\leq&~\bar{\kappa}\psi_{0}^{2}\|\partial_{y}^{\gamma_{2}-1}\partial_{x}^{\gamma_{1}+2}\mathbb{U}\|_{\mathring{H}^{0}(\mathbb{T})}
                      \|\partial_{y}^{\gamma_{2}}\partial_{x}^{\gamma_{1}}u\|_{\mathring{H}^{0}(\mathbb{T})}\nonumber\\
\leq&~\bar{\kappa}^{2}\|\partial_{y}^{\gamma_{2}-1}\partial_{x}^{\gamma_{1}+2}\mathbb{U}\|_{\mathring{H}^{0}(\mathbb{T})}^{2}
     +C\|D^{\gamma}u\|_{L_{\psi}^{2}(\Omega)}^{2}
     +\frac{1}{16}\|\partial_{y}D^{\gamma}u\|_{L_{\psi}^{2}(\Omega)}^{2}\nonumber\\
\leq&~\bar{\kappa}^{2}
      \mathcal{Q}\left(\|\mathbb{U}\|_{\ddot{H}^{s}(\hat{\Omega})},
                       \sum_{i=0}^{\frac{\gamma_{2}-5}{2}}\|\partial_{t}^{i}F\|_{\mathring{H}^{s-1-2i}(\mathbb{T})}\right)
     +C\|D^{\gamma}u\|_{L_{\psi}^{2}(\Omega)}^{2}
     +\frac{1}{16}\|\partial_{y}D^{\gamma}u\|_{L_{\psi}^{2}(\Omega)}^{2}.\nonumber
\end{align}
If $\gamma_{2}$ is even, by using integration by parts, Lemma \ref{bc} and Sobolev imbedding, we have
\begin{align}\label{I6c233}
&-\psi_{0}^{2}(\kappa\partial_{y}^{\gamma_{2}-1}\partial_{x}^{\gamma_{1}+2}\mathbb{U},\partial_{y}^{\gamma_{2}}\partial_{x}^{\gamma_{1}}u)_{L_{x}^{2}}|_{y=0}\\
=&\bar{\kappa}\psi_{0}^{2}(\partial_{y}^{\gamma_{2}-1}\partial_{x}^{\gamma_{1}+1}\mathbb{U},\partial_{y}^{\gamma_{2}}\partial_{x}^{\gamma_{1}+1}u)_{L_{x}^{2}}|_{y=0}\nonumber\\
\leq& \bar{\kappa}\psi_{0}^{2}\|\partial_{y}^{\gamma_{2}-1}\partial_{x}^{\gamma_{1}+1}\mathbb{U}\|_{\mathring{H}^{0}(\mathbb{T})}
                      \|\partial_{y}^{\gamma_{2}}\partial_{x}^{\gamma_{1}+1}u\|_{\mathring{H}^{0}(\mathbb{T})}\nonumber\\
\leq& \bar{\kappa}^{2}\left(\|\mathbb{U}\|_{\hat{H}_{\psi}^{s}(\hat{\Omega})}^{2}
                +\frac{\tau}{16}\|\partial_{y}\mathbb{U}\|_{\hat{H}_{\psi}^{s}(\hat{\Omega})}^{2}\right)
     +C_{\tau}\|\partial_{y}^{\gamma_{2}}\partial_{x}^{\gamma_{1}+1}u\|_{\mathring{H}^{0}(\mathbb{T})}^{2}\nonumber\\
\leq& \bar{\kappa}^{2}\left(\|\mathbb{U}\|_{\hat{H}_{\psi}^{s}(\Omega)}^{2}
                +\frac{\tau}{16}\|\partial_{y}\mathbb{U}\|_{\hat{H}_{\psi}^{s}(\Omega)}^{2}+C_{\tau}c^{*}\right)
     +C_{\tau}c_{*}\|u\|_{\ddot{H}^{s}(\Omega)}^{2}
     +C_{\tau}c_{*}\sum_{i=0}^{\frac{\gamma_{2}}{2}-1}\|\partial_{t}^{i}F\|_{\mathring{H}^{s-1-2i}(\mathbb{T})}^{2}\nonumber\\
\leq& \bar{\kappa}^{2}\left(\|\mathbb{U}\|_{\hat{H}_{\psi}^{s}(\Omega)}^{2}
                +\frac{\tau}{16}\|\partial_{y}\mathbb{U}\|_{\hat{H}_{\psi}^{s}(\Omega)}^{2}+C_{\tau}c^{*}\right)
     +C_{\tau}c_{*}\|u\|_{\mathcal{H}_{\psi,\varphi}^{s}(\Omega)}^{2}
     +C_{\tau}c_{*}\sum_{i=0}^{\frac{\gamma_{2}}{2}-1}\|\partial_{t}^{i}F\|_{\mathring{H}^{s-1-2i}(\mathbb{T})}^{2}\nonumber.
\end{align}
where $c_{*}$ is a constant depending on $\|\mathbb{U}\|_{\ddot{H}^{s}(\hat{\Omega})}$,
$\sum\limits_{i=0}^{\frac{\gamma_{2}}{2}-2}\|\partial_{t}^{i}F\|_{\mathring{H}^{s-3-2i}(\mathbb{T})}$, and
$c^{*}$ is a constant depending on $\|\mathbb{U}\|_{\ddot{H}^{s}(\hat{\Omega})}$,
$\sum\limits_{i=0}^{\frac{\gamma_{2}}{2}-3}\|\partial_{t}^{i}F\|_{\mathring{H}^{s-1-2i}(\mathbb{T})}$,
and $\sum\limits_{i=0}^{\frac{\gamma_{2}}{2}-2}\|\partial_{t}^{i}P\|_{H^{s-2i}(\mathbb{T})}$.
Since $\bar{\kappa}\leq \epsilon^{-2}\bar{\eta}$, one has
\begin{align}\label{I6c234}
&-\sum_{0<\beta_{2}\leq\gamma_{2}}\binom{\gamma_{2}}{\beta_{2}}
   (\partial_{y}^{\beta_{2}}\kappa\partial_{y}^{\gamma_{2}-\beta_{2}}\partial_{x}^{\gamma_{1}+1}\mathbb{U},\psi^{2}\partial_{x}D^{\gamma}u)\\
\leq& C\epsilon^{-2}\sqrt{\bar{\eta}}\sum_{0<\beta_{2}\leq\gamma_{2}}\binom{\gamma_{2}}{\beta_{2}}
       \|\partial_{y}^{\gamma_{2}-\beta_{2}}\partial_{x}^{\gamma_{1}+1}\mathbb{U}\|_{L_{\psi}^{2}(\Omega)}
       \|\sqrt{\eta}\partial_{x}D^{\gamma}u\|_{L_{\psi_{\gamma}}^{2}(\Omega)}\nonumber\\
\leq& C\bar{\eta}\|\mathbb{U}\|_{\hat{H}_{\psi}^{s}(\Omega)}^{2}
     +\frac{1}{16}\|\sqrt{\eta}\partial_{x}D^{\gamma}u\|_{L_{\psi_{\gamma}}^{2}(\Omega)}^{2}.\nonumber
\end{align}
A combination of \eqref{I6c231}, \eqref{I6c232}, \eqref{I6c233} and \eqref{I6c234} with \eqref{I6c23} yields that
\begin{align}\label{I6c23'}
\mathcal{I}_{6}
\leq& C_{\tau}\|D^{\gamma}u\|_{L_{\psi}^{2}(\Omega)}^{2}
     +\frac{1}{16}\|\partial_{y}D^{\gamma}u\|_{L_{\psi}^{2}(\Omega)}^{2}
     +\frac{1}{16}\|\sqrt{\eta}\partial_{x}D^{\gamma}u\|_{L_{\psi_{\gamma}}^{2}(\Omega)}^{2}
     +C_{\tau}\tilde{c}_{*}\|u\|_{\mathcal{H}_{\psi,\varphi}^{s}(\Omega)}^{2}\\
    &+\bar{\kappa}^{2}\left(\|\mathbb{U}\|_{\hat{H}_{\psi}^{s}(\Omega)}^{2}
                +\tau\|\partial_{y}\mathbb{U}\|_{\hat{H}_{\psi}^{s}(\Omega)}^{2}+C_{\tau}\tilde{c}^{*}\right)
     +C\bar{\eta}\|\mathbb{U}\|_{\hat{H}_{\psi}^{s}(\Omega)}^{2}
     +C_{\tau}\tilde{c}_{*}\sum_{i=0}^{\frac{s-3}{2}}\|\partial_{t}^{i}F\|_{\mathring{H}^{s-1-2i}(\mathbb{T})}^{2},\nonumber
\end{align}
where $\tilde{c}_{*}$ is a constant depending on $\|\mathbb{U}\|_{\ddot{H}^{s}(\hat{\Omega})}$,
$\sum\limits_{i=0}^{\frac{s-5}{2}}\|\partial_{t}^{i}F\|_{\mathring{H}^{s-3-2i}(\mathbb{T})}$,
and $\tilde{c}^{*}$ is a constant depending on $\|\mathbb{U}\|_{\ddot{H}^{s}(\hat{\Omega})}$,
$\sum\limits_{i=0}^{\frac{s-7}{2}}\|\partial_{t}^{i}F\|_{\mathring{H}^{s-1-2i}(\mathbb{T})}$,
and $\sum\limits_{i=0}^{\frac{s-5}{2}}\|\partial_{t}^{i}P\|_{H^{s-2i}(\mathbb{T})}$.

According to \eqref{I6c21}, \eqref{I6c22} and \eqref{I6c23'}, we have estimate for $\mathcal{I}_{6}$ when $|\gamma|=s$ as follows,
\begin{align*}
\mathcal{I}_{6}
\leq&~C_{*}C_{\tau}C_{\theta}\|D^{\gamma}u\|_{L_{\psi}^{2}(\Omega)}^{2}
     +\frac{3}{8}\|\partial_{y}D^{\gamma}u\|_{L_{\psi}^{2}(\Omega)}^{2}
     +\frac{1}{16}\|\sqrt{\eta}\partial_{x}D^{\gamma}u\|_{L_{\psi_{\gamma}}^{2}(\Omega)}^{2}
     +C_{\tau}\tilde{c}_{*}\|u\|_{\mathcal{H}_{\psi,\varphi}^{s}(\Omega)}^{2}\\
    &+\bar{\kappa}^{2}\left(\|\mathbf{U}_{s+1}\|_{L_{\hat{\varphi}}^{2}(\Omega)}^{2}
                +\tau\|\mathcal{P}\mathbf{U}_{s+1}\|_{L_{\hat{\varphi}}^{2}(\Omega)}^{2}
                +\tau\|\partial_{y}\mathbf{U}_{s}\|_{L_{\varphi}^{2}(\Omega)}^{2}
                +\tau\|\mathcal{P}\mathbf{U}_{s}\|_{L_{\varphi}^{2}(\Omega)}^{2}\right)
     +C\bar{\eta}\|\mathbb{U}\|_{\hat{H}_{\psi}^{s}(\Omega)}^{2}\\
    &+\bar{\kappa}^{2}\left(\|\mathbb{U}\|_{\mathcal{H}_{\psi,\varphi}^{s}(\Omega)}^{2}
                +\tau\|\partial_{y}\mathbb{U}\|_{\hat{H}_{\psi}^{s}(\Omega)}^{2}
                +C_{\tau}\tilde{c}^{*}\right)
     +C_{\tau}\tilde{c}_{*}\sum_{i=0}^{\frac{s-3}{2}}\|\partial_{t}^{i}F\|_{\mathring{H}^{s-1-2i}(\mathbb{T})}^{2}.
\end{align*}

\noindent{\bf Case 3:} $|\gamma|=s-1$. If $\gamma_{1}=s-1$, we have
\begin{align*}
\mathcal{I}_{6}
=&(\kappa\partial_{x}^{s+1}\mathbb{U}, \psi^{2}D^{\gamma}u)\\
\leq& C_{*}\bar{\kappa}\|\partial_{x}^{s+1}\mathbb{U}\|_{L_{\hat{\psi}}^{2}(\Omega)}\|\sqrt{\y}D^{\gamma}u\|_{L_{\psi}^{2}(\Omega)}\\
\leq& C_{*}\bar{\kappa}\|\mathcal{P}\mathbf{U}_{s+1}\|_{L_{\varphi}^{2}(\Omega)}\|\sqrt{\y}D^{\gamma}u\|_{L_{\psi}^{2}(\Omega)}\\
\leq& \bar{\kappa}^{2}\tau
      \left(\|\sqrt{\y}\mathbf{U}_{s+1}\|_{L_{\hat{\varphi}}^{2}(\Omega)}^{2}+\|\mathcal{P}\mathbf{U}_{s+1}\|_{L_{\hat{\varphi}}^{2}(\Omega)}^{2}\right)
     +C_{*}C_{\tau}\|\sqrt{\y}D^{\gamma}u\|_{L_{\psi}^{2}(\Omega)}^{2}.
\end{align*}
If $\gamma_{1}=s-2$, if follows from Corollary \ref{ncC2} that
\begin{align*}
\mathcal{I}_{6}
=& (\kappa\partial_{y}\partial_{x}^{s}\mathbb{U}, \psi^{2}D^{\gamma}u)
  +(\partial_{y}\kappa\partial_{x}^{s}\mathbb{U}, \psi^{2}D^{\gamma}u)\\
\leq& C_{\theta}\bar{\kappa}\|\partial_{y}\partial_{x}^{s}\mathbb{U}\|_{L_{\hat{\psi}}^{2}(\Omega)}\|\sqrt{\y}D^{\gamma}u\|_{L_{\psi}^{2}(\Omega)}
     +C_{\theta}\bar{\kappa}\|\partial_{x}^{s}\mathbb{U}\|_{L_{\hat{\psi}}^{2}(\Omega)}\|\sqrt{\y}D^{\gamma}u\|_{L_{\psi}^{2}(\Omega)}\\
\leq& C_{*}\bar{\kappa}\left(\|\mathbf{U}_{s}\|_{L_{\varphi}^{2}(\Omega)}+\|\mathcal{P}\mathbf{U}_{s}\|_{L_{\varphi}^{2}(\Omega)}\right)
           \|\sqrt{\y}D^{\gamma}u\|_{L_{\psi}^{2}(\Omega)}\\
\leq& \bar{\kappa}^{2}\left(\|\mathbf{U}_{s}\|_{L_{\varphi}^{2}(\Omega)}^{2}
                            +\tau\|\mathcal{P}\mathbf{U}_{s}\|_{L_{\varphi}^{2}(\Omega)}^{2}\right)
     +C_{*}C_{\tau}C_{\theta}\|\sqrt{\y}D^{\gamma}u\|_{L_{\psi}^{2}(\Omega)}^{2}.
\end{align*}
If $\gamma_{1}<s-2$, then $\gamma_{1}+2\leq s-1$, it is clear that
\begin{align*}
\mathcal{I}_{6}
=& (\kappa\partial_{y}^{\gamma_{2}}\partial_{x}^{\gamma_{1}+2}\mathbb{U},\psi^{2}D^{\gamma}u)
  +\sum_{0<\beta_{2}\leq\gamma_{2}}\binom{\gamma_{2}}{\beta_{2}}
   (\partial_{y}^{\beta_{2}}\kappa\partial_{y}^{\gamma_{2}-\beta_{2}}\partial_{x}^{\gamma_{1}+2}\mathbb{U},\psi^{2}D^{\gamma}u)\\
\leq& \bar{\kappa}^{2}\left(\|\mathbb{U}\|_{\hat{H}_{\psi}^{s}(\Omega)}^{2}+\tau\|\partial_{y}\mathbb{U}\|_{\hat{H}_{\psi}^{s}(\Omega)}^{2}\right)
     +C_{*}C_{\tau}C_{\theta}\|D^{\gamma}u\|_{L_{\psi}^{2}(\Omega)}^{2}.
\end{align*}
Thus for $|\gamma|=s-1$ we have
\begin{align*}
\mathcal{I}_{6}
\leq& \bar{\kappa}^{2}\left(\tau\|\mathcal{P}\mathbf{U}_{s+1}\|_{L_{\hat{\varphi}}^{2}(\Omega)}^{2}
                +\tau\|\sqrt{\y}\mathbf{U}_{s+1}\|_{L_{\hat{\varphi}}^{2}(\Omega)}^{2}
                +\tau\|\partial_{y}\mathbb{U}\|_{\hat{H}_{\psi}^{s}(\Omega)}^{2}
                +\tau\|\mathcal{P}\mathbf{U}_{s}\|_{L_{\hat{\varphi}}^{2}(\Omega)}^{2}
                +\|\mathbb{U}\|_{\mathcal{H}_{\psi,\varphi}^{s}(\Omega)}^{2}\right)\\
    &+C_{\tau}C_{*}C_{\theta}\|\sqrt{\y}D^{\gamma}u\|_{L_{\psi}^{2}(\Omega)}^{2}.
\end{align*}

\noindent{\bf Case 4:} $|\gamma|\leq s-2$.
If $\gamma_{1}=s-2$, with the help of Corollary \ref{ncC2}, one has
\begin{align*}
\mathcal{I}_{6}
=&(\kappa\partial_{x}^{s}\mathbb{U}, \psi^{2}D^{\gamma}u)\\
\leq&~C\bar{\kappa}\|\partial_{x}^{s}\mathbb{U}\|_{L_{\hat{\psi}}^{2}(\Omega)}\|\sqrt{\y}D^{\gamma}u\|_{L_{\psi}^{2}(\Omega)}\\
\leq&~\bar{\kappa}^{2}\tau\|\mathcal{P}\mathbf{U}_{s}\|_{L_{\varphi}^{2}(\Omega)}^{2}
     +C_{\tau}C_{*}\|\sqrt{\y}D^{\gamma}u\|_{L_{\psi}^{2}(\Omega)}^{2},
\end{align*}
else, by the Cauchy inequality,
\begin{align*}
\mathcal{I}_{6}
=&(\kappa\partial_{y}^{\gamma_{2}}\partial_{x}^{\gamma_{1}+2}\mathbb{U},\psi^{2}D^{\gamma}u)
  +\sum_{0<\beta_{2}\leq\gamma_{2}}\binom{\gamma_{2}}{\beta_{2}}
   (\partial_{y}^{\beta_{2}}\kappa\partial_{y}^{\gamma_{2}-\beta_{2}}\partial_{x}^{\gamma_{1}+2}\mathbb{U},\psi^{2}D^{\gamma}u)\\
\leq& \bar{\kappa}^{2}\|\partial_{y}^{\gamma_{2}}\partial_{x}^{\gamma_{1}+2}\mathbb{U}\|_{L_{\psi}^{2}(\Omega)}^{2}
     +C\bar{\kappa}^{2}\sum_{0<\beta_{2}\leq\gamma_{2}}\|\partial_{y}^{\gamma_{2}-\beta_{2}}\partial_{x}^{\gamma_{1}+2}\mathbb{U}\|_{L_{\psi}^{2}(\Omega)}^{2}
     +C\|D^{\gamma}u\|_{L_{\psi}^{2}(\Omega)}^{2}\\
\leq& \bar{\kappa}^{2}\|\mathbb{U}\|_{\hat{H}_{\psi}^{s}(\Omega)}^{2}
     +C\|D^{\gamma}u\|_{L_{\psi}^{2}(\Omega)}^{2}.
\end{align*}
Accordingly, for $|\gamma|\leq s-2$ we have
\begin{align*}
\mathcal{I}_{6}
\leq\bar{\kappa}^{2}\left(\|\mathbb{U}\|_{\hat{H}_{\psi}^{s}(\Omega)}^{2}+\tau\|\mathcal{P}\mathbf{U}_{s}\|_{L_{\varphi}^{2}(\Omega)}^{2}\right)
    +C_{\tau}C_{*}\|D^{\gamma}u\|_{L_{\psi}^{2}(\Omega)}^{2}.
\end{align*}

Collecting all the estimates of $\mathcal{I}_{6}$ above, we can deduce
\begin{align*}
\mathcal{I}_{6}
\leq&~C_{\tau}C_{*}C_{\theta}\|\sqrt{\y}D^{\gamma}u\|_{L_{\psi_{\gamma}}^{2}(\Omega)}^{2}
     +\frac{3}{8}\|\partial_{y}D^{\gamma}u\|_{L_{\psi_{\gamma}}^{2}(\Omega)}^{2}
     +\frac{1}{16}\|\sqrt{\eta}\partial_{x}D^{\gamma}u\|_{L_{\psi_{\gamma}}^{2}(\Omega)}^{2}
     +C_{\tau}\tilde{c}_{*}\|u\|_{\mathcal{H}_{\psi,\varphi}^{s}(\Omega)}^{2}\\
    &+\bar{\kappa}^{2}\left(\|\mathbf{U}_{s+1}\|_{L_{\hat{\varphi}}^{2}(\Omega)}^{2}
                     +\tau\|\mathcal{P}\mathbf{U}_{s+1}\|_{L_{\hat{\varphi}}^{2}(\Omega)}^{2}
                     +\tau\|\partial_{y}\mathbf{U}_{s+1}\|_{L_{\hat{\varphi}}^{2}(\Omega)}^{2}
                     +\tau\|\partial_{y}\mathbf{U}_{s}\|_{L_{\varphi}^{2}(\Omega)}^{2}
                     +\tau\|\mathcal{P}\mathbf{U}_{s}\|_{L_{\varphi}^{2}(\Omega)}^{2}\right)\\
    &+\bar{\kappa}^{2}\left(\|\mathbb{U}\|_{\mathcal{H}_{\psi,\varphi}^{s}(\Omega)}^{2}
                +\tau\|\partial_{y}\mathbb{U}\|_{\hat{H}_{\psi}^{s}(\Omega)}^{2}
                +C_{\tau}\tilde{c}^{*}\right)
     +C\bar{\eta}\|\mathbb{U}\|_{\hat{H}_{\psi}^{s}(\Omega)}^{2}
     +C_{\tau}\tilde{c}_{*}\sum_{i=0}^{\frac{s-3}{2}}\|\partial_{t}^{i}F\|_{\mathring{H}^{s-1-2i}(\mathbb{T})}^{2}.
\end{align*}

\noindent{\bf \underline{Estimate of~$\mathcal{I}_{7}$.}} Using the Cauchy inequality, one has
\begin{align*}
\mathcal{I}_{7}\leq \tau\|D^{\gamma}F\|_{L_{\psi_{\gamma}}^{2}(\Omega)}^{2}+C_{\tau}\|D^{\gamma}u\|_{L_{\psi_{\gamma}}^{2}(\Omega)}^{2}.
\end{align*}

\noindent{\bf \underline{Estimate of~$\mathcal{I}_{8}$.}} Note that, under Assumption \ref{MA}, $\sqrt{t\theta\y}\vartheta_{c}\leq C$ and $\mathcal{P}\leq\frac{C_{*}}{\sqrt{t+\bar{\varsigma}}}$ hold, by using integration by parts and Corollary \ref{ncC2}, one can get
\begin{align*}
\mathcal{I}_{8}
=&-\sum_{\substack{\gamma\in\Gamma,\\0<\beta_{2}\leq\gamma_{2}}}\binom{\gamma_{2}}{\beta_{2}}(\partial_{y}^{\beta_{2}}\eta\partial_{x}D^{\gamma-\beta}u,\psi_{\gamma}^{2}\partial_{x}D^{\gamma}u)\\
\leq&~C\bar{\eta}(t\theta)^{\beta_{2}-\frac{1}{2}}\sum_{\substack{\gamma\in\Gamma,\\0<\beta_{2}\leq\gamma_{2}}}
       \|\partial_{x}D^{\gamma-\beta}u\|_{L_{\hat{\psi}_{\gamma}}^{2}(\Omega)}
       \|\vartheta_{c}\partial_{x}D^{\gamma}u\|_{L_{\psi_{\gamma}}^{2}(\Omega)}\\
\leq&~C\sqrt{\bar{\eta}}(t\theta)^{\beta_{2}-\frac{1}{2}}\sum_{\substack{\gamma\in\Gamma,\\0<\beta_{2}\leq\gamma_{2}}}
       \left(\|u\|_{\hat{H}_{\psi}^{s}(\Omega)}
        +\|\partial_{x}^{s}u\|_{L_{\hat{\psi}}^{2}(\Omega)}
        +\|\partial_{y}\partial_{x}^{s}u\|_{L_{\hat{\psi}}^{2}(\Omega)}\right)
         \|\sqrt{\eta}\partial_{x}D^{\gamma}u\|_{L_{\psi_{\gamma}}^{2}(\Omega)}\\
\leq&~C\sqrt{\bar{\eta}}(t\theta)^{\beta_{2}-\frac{1}{2}}\sum_{\substack{\gamma\in\Gamma,\\0<\beta_{2}\leq\gamma_{2}}}
       \left(\|u\|_{\hat{H}_{\psi}^{s}(\Omega)}
        +\|\mathcal{U}_{s}\|_{L_{\varphi}^{2}(\Omega)}
        +\|\mathcal{P}\mathcal{U}_{s}\|_{L_{\varphi}^{2}(\Omega)}\right)
       \|\sqrt{\eta}\partial_{x}D^{\gamma}u\|_{L_{\psi_{\gamma}}^{2}(\Omega)}\\
\leq&~C_{*}C_{\theta}\left(\|\mathcal{U}_{s}\|_{L_{\varphi}^{2}(\Omega)}+\|u\|_{\hat{H}_{\psi}^{s}(\Omega)}\right)
       \|\sqrt{\eta}\partial_{x}D^{\gamma}u\|_{L_{\psi_{\gamma}}^{2}(\Omega)}\\
\leq&~\frac{1}{16}\|\sqrt{\eta}\partial_{x}D^{\gamma}u\|_{L_{\psi_{\gamma}}^{2}(\Omega)}^{2}
     +C_{*}C_{\theta}\|u\|_{\mathcal{H}_{\psi,\varphi}^{s}(\Omega)}^{2}.
\end{align*}

\noindent{\bf \underline{Estimate of~$\mathcal{I}_{9}$.}} Noticing that $|\partial_{y}\psi_{\gamma}|\leq C\psi_{\gamma}$, one has
\begin{align*}
\mathcal{I}_{9}
=&-2(\partial_{y}D^{\gamma}u,\psi_{\gamma}\partial_{y}\psi_{\gamma}D^{\gamma}u)\\
\leq& \frac{1}{16}\|\partial_{y}D^{\gamma}u\|_{L_{\psi_{\gamma}}^{2}(\Omega)}^{2}
     +C\|D^{\gamma}u\|_{L_{\psi_{\gamma}}^{2}(\Omega)}^{2}.
\end{align*}

\noindent{\bf \underline{Estimate of~$\mathcal{I}_{10}$.}}

\noindent{\bf Case 1:} $|\gamma|\leq s-1$. By using the trace theorem one can get
\begin{align*}
\mathcal{I}_{10}
\leq& C(\|u\|_{\hat{H}_{\psi}^{s}(\Omega)}+\|\partial_{y}u\|_{\hat{H}_{\psi}^{s}(\Omega)})
       \|u\|_{\hat{H}_{\psi}^{s}(\Omega)}\\
\leq& C\|u\|_{\mathcal{H}_{\psi,\varphi}^{s}(\Omega)}^{2}
     +\frac{1}{8(s+1)(s+2)}\|\partial_{y}u\|_{\hat{H}_{\psi}^{s}(\Omega)}^{2}.
\end{align*}

\noindent{\bf Case 2:} $|\gamma|=s$.  In this case, if $\gamma_{2}$ is odd, one can use Lemma \ref{bc} and the trace theorem to obtain
\begin{align*}
\mathcal{I}_{10}
=&-(\psi^{2}\partial_{y}^{\gamma_{2}+1}\partial_{x}^{s-\gamma_{2}}u, \partial_{y}^{\gamma_{2}}\partial_{x}^{s-\gamma_{2}}u)_{L_{x}^{2}}|_{y=0}\\
\leq& C\psi_{0}^{2}\|\partial_{y}^{\gamma_{2}+1}\partial_{x}^{s-\gamma_{2}}u\|_{\mathring{H}^{0}(\mathbb{T})}
       \|\partial_{y}^{\gamma_{2}}\partial_{x}^{s-\gamma_{2}}u\|_{\mathring{H}^{0}(\mathbb{T})}\\
\leq& \psi_{0}^{2}\left(c_{*}\|u\|_{\ddot{H}^{s}(\hat{\Omega})}
                        +c_{*}\sum_{i=0}^{\frac{\gamma_{2}-1}{2}}\|\partial_{t}^{i}F\|_{\mathring{H}^{s-1-2i}(\mathbb{T})}
                        +c^{*}\bar{\kappa}\right)
                  \left(\|u\|_{\hat{H}_{\psi}^{s}(\hat{\Omega})}+\|\partial_{y}u\|_{\hat{H}_{\psi}^{s}(\hat{\Omega})}\right)\\
\leq& c_{*}\|u\|_{\mathcal{H}_{\psi,\varphi}^{s}(\Omega)}^{2}
     +\frac{1}{8(s+2)(s+1)}\|\partial_{y}u\|_{\hat{H}_{\psi}^{s}(\Omega)}^{2}
     +c_{*}\sum_{i=0}^{\frac{\gamma_{2}-1}{2}}\|\partial_{t}^{i}F\|_{\mathring{H}^{s-2i-1}(\mathbb{T})}^{2}
     +c^{*}\bar{\kappa}^{2},
\end{align*}
where $c_{*}$ depending only on $\|\mathbb{U}\|_{\ddot{H}^{s}(\hat{\Omega})}$ and
$\sum\limits_{i=0}^{\frac{\gamma_{2}-3}{2}}\|\partial_{t}^{i}F\|_{\mathring{H}^{s-3-2i}(\mathbb{T})}$,
and a constant  $c^{*}$ depending only on $\|\mathbb{U}\|_{\ddot{H}^{s}(\hat{\Omega})}$,
$\sum\limits_{i=0}^{\frac{\gamma_{2}-5}{2}}\|\partial_{t}^{i}F\|_{\mathring{H}^{s-1-2i}(\mathbb{T})}$,
and $\sum\limits_{i=0}^{\frac{\gamma_{2}-3}{2}}\|\partial_{t}^{i}P\|_{H^{s-2i}(\mathbb{T})}$.

If $\gamma_{2}$ is even, then $s-\gamma_{2}>0$, one can use the integration by parts and Lemma \ref{bc} to obtain
\begin{align*}
\mathcal{I}_{10}
=&(\psi^{2}\partial_{y}^{\gamma_{2}+1}\partial_{x}^{s-\gamma_{2}-1}u, \partial_{y}^{\gamma_{2}}\partial_{x}^{s-\gamma_{2}+1}u)_{L_{x}^{2}}|_{y=0}\\
\leq& C\psi_{0}^{2}\|\partial_{y}^{\gamma_{2}+1}\partial_{x}^{s-\gamma_{2}-1}u\|_{\mathring{H}^{0}(\mathbb{T})}
       \|\partial_{y}^{\gamma_{2}}\partial_{x}^{s-\gamma_{2}+1}u\|_{\mathring{H}^{0}(\mathbb{T})}\\
\leq& \psi_{0}^{2}\left(\|u\|_{\hat{H}_{\psi}^{s}(\Omega)}+\|\partial_{y}u\|_{\hat{H}_{\psi}^{s}(\hat{\Omega})}\right)
       \left(c_{*}\|u\|_{\ddot{H}^{s}(\hat{\Omega})}
     +c_{*}\sum_{i=0}^{\frac{\gamma_{2}}{2}}\|\partial_{t}^{i}F\|_{\mathring{H}^{s-1-2i}(\mathbb{T})}
     +c^{*}\bar{\kappa}\right)\\
\leq& c_{*}\|u\|_{\mathcal{H}_{\psi,\varphi}^{s}(\Omega)}^{2}
     +\frac{1}{8(s+2)(s+1)}\|\partial_{y}u\|_{\hat{H}_{\psi}^{s}(\Omega)}^{2}
     +c_{*}\sum_{i=0}^{\frac{\gamma_{2}}{2}}\|\partial_{t}^{i}F\|_{\mathring{H}^{s-1-2i}(\mathbb{T})}^{2}
     +c^{*}\bar{\kappa}^{2},
\end{align*}
where $c_{*}$ depending only on $\|\mathbb{U}\|_{\ddot{H}^{s}(\hat{\Omega})}$ and
$\sum\limits_{i=0}^{\frac{\gamma_{2}-4}{2}}\|\partial_{t}^{i}F\|_{\mathring{H}^{s-3-2i}(\mathbb{T})}$,
and a constant  $c^{*}$ depending only on $\|\mathbb{U}\|_{\ddot{H}^{s}(\hat{\Omega})}$,
$\sum\limits_{i=0}^{\frac{\gamma_{2}-6}{2}}\|\partial_{t}^{i}F\|_{\mathring{H}^{s-1-2i}(\mathbb{T})}$,
and $\sum\limits_{i=0}^{\frac{\gamma_{2}-4}{2}}\|\partial_{t}^{i}\partial_{x}P\|_{H^{s-2i-1}(\mathbb{T})}$,

\noindent{\bf Case 3:} $|\gamma|=s+1$. Applying integration by parts twice, noticing $\partial_{y}^{2}u(t,x,0)=F(t,x,0)$, by the trace theorem we have
\begin{align*}
\mathcal{I}_{10}
=(\hat{\psi}^{2}\partial_{y}^{3}\partial_{x}^{s-1}u, \partial_{y}^{2}\partial_{x}^{s-1}u)_{L_{x}^{2}}|_{y=0}
=& (\psi^{2}\partial_{y}^{3}\partial_{x}^{s-3}u, \partial_{y}^{2}\partial_{x}^{s+1}u)_{L_{x}^{2}}|_{y=0}\\
=&-(\psi^{2}\partial_{y}^{3}\partial_{x}^{s-3}u, \partial_{x}^{s+1}F)_{L_{x}^{2}}|_{y=0}\\
\leq& C\psi_{0}^{2}\|\partial_{y}^{3}\partial_{x}^{s-3}u\|_{\mathring{H}^{0}(\mathbb{T})}
       \|\partial_{x}^{s+1}F\|_{\mathring{H}^{0}(\mathbb{T})}\\
\leq& C\psi_{0}^{2}\left(\|u\|_{\hat{H}_{\psi}^{s}(\hat{\Omega})}+\|\partial_{y}u\|_{\hat{H}_{\psi}^{s}(\hat{\Omega})}\right)
       \|\partial_{x}^{s+1}F\|_{\mathring{H}^{0}(\mathbb{T})}\\
\leq& C\|u\|_{\mathcal{H}_{\psi,\varphi}^{s}(\Omega)}^{2}
     +\frac{1}{8(s+2)(s+1)}\|\partial_{y}u\|_{\hat{H}_{\psi}^{s}(\Omega)}^{2}
     +C\|\partial_{x}^{s+1}F\|_{\mathring{H}^{0}(\mathbb{T})}^{2}.
\end{align*}

In summary, one has
\begin{align*}
\mathcal{I}_{10}
\leq \bar{c}_{*}\|u\|_{\mathcal{H}_{\psi,\varphi}^{s}(\Omega)}^{2}
    +\frac{1}{8(s+1)(s+2)}\|\partial_{y}u\|_{\hat{H}_{\psi}^{s}(\Omega)}^{2}
    +\bar{c}^{*}\bar{\kappa}^{2}
    +\sum_{i=0}^{\frac{s-1}{2}}\|\partial_{t}^{i}F\|_{\mathring{H}^{s-1-2i}(\mathbb{T})}^{2}
    +C\|\partial_{x}^{s+1}F\|_{\mathring{H}^{0}(\mathbb{T})}^{2},
\end{align*}
where $\bar{c}_{*}$ is a constant depending on $\|\mathbb{U}\|_{\ddot{H}^{s}(\hat{\Omega})}$,
$\sum\limits_{i=0}^{\frac{s-3}{2}}\|\partial_{t}^{i}F\|_{\mathring{H}^{s-3-2i}(\mathbb{T})}$,
and $\bar{c}^{*}$ is a constant depending on $\|\mathbb{U}\|_{\ddot{H}^{s}(\hat{\Omega})}$,
$\sum\limits_{i=0}^{\frac{s-5}{2}}\|\partial_{t}^{i}F\|_{\mathring{H}^{s-1-2i}(\mathbb{T})}$,
and $\sum\limits_{i=0}^{\frac{s-3}{2}}\|\partial_{t}^{i}\partial_{x}P\|_{H^{s-2i-1}(\mathbb{T})}$.

Now substituting all the estimates of $\mathcal{I}_{1},\mathcal{I}_{2},...,\mathcal{I}_{10}$ above into \eqref{le}, under the assumption of Theorem \ref{wpae}, for any $t\in[0,T_{*}]$ and $\gamma\in\Gamma$ one has
\begin{align*}
&\frac{1}{2}\frac{d}{dt}\|D^{\gamma}u \|_{L_{\psi_{\gamma}}^{2}(\Omega)}^{2}
+\Lambda\delta\|\sqrt{\y}D^{\gamma}u \|_{L_{\psi_{\gamma}}^{2}(\Omega)}^{2}
+\|\partial_{y}D^{\gamma}u\|_{L_{\psi_{\gamma}}^{2}(\Omega)}^{2}
+\|\sqrt{\eta}\partial_{x}D^{\gamma}u\|_{L_{\psi_{\gamma}}^{2}(\Omega)}^{2}\\
\leq& C_{*}\left(\|\mathbb{U}\|_{\mathcal{H}_{\psi,\varphi}^{s}(\Omega)}
                +\|\partial_{y}\mathbb{U}\|_{\hat{H}_{\psi}^{s}(\Omega)}
                +\|\partial_{y}\mathbf{U}_{s}\|_{L_{\varphi}^{2}(\Omega)}
                +\|\mathcal{P}\mathbf{U}_{s}\|_{L_{\varphi}^{2}(\Omega)}
                +C_{\tau}C_{\theta}\right)
          \|u\|_{\mathcal{H}_{\psi,\varphi}^{s}(\Omega)}^{2}\\
    &+C_{*}\|u\|_{\mathcal{H}_{\psi,\varphi}^{s}(\Omega)}^{4}
     +\frac{2}{(s+1)(s+2)}\|\mathcal{P}\mathcal{U}_{s}\|_{L_{\varphi}^{2}(\Omega)}^{2}
     +C_{*}C_{\tau}C_{\theta}\|u\|_{\hat{H}_{\psi,\y}^{s}(\Omega)}^{2}
     +\frac{3}{8(s+2)(s+1)}\|\partial_{y}u\|_{\hat{H}_{\psi}^{s}(\Omega)}^{2}\\
    &+\frac{1}{16}\|\sqrt{\eta}\partial_{x}D^{\gamma}u\|_{L_{\psi_{\gamma}}^{2}(\Omega)}^{2}
     +\frac{1}{2}\|\partial_{y}D^{\gamma}u\|_{L_{\psi_{\gamma}}^{2}(\Omega)}^{2}
     +\bar{\kappa}^{2}\left(\|\mathbb{U}\|_{\mathcal{H}_{\psi,\varphi}^{s}(\Omega)}^{2}
                +\tau\|\partial_{y}\mathbb{U}\|_{\hat{H}_{\psi}^{s}(\Omega)}^{2}
                +C_{\tau}C_{*}\right)\\
    &+\bar{\kappa}^{2}\left(\|\mathbf{U}_{s+1}\|_{L_{\hat{\varphi}}^{2}(\Omega)}^{2}
                     +\tau\|\mathcal{P}\mathbf{U}_{s+1}\|_{L_{\hat{\varphi}}^{2}(\Omega)}^{2}
                     +\tau\|\partial_{y}\mathbf{U}_{s+1}\|_{L_{\hat{\varphi}}^{2}(\Omega)}^{2}
                     +\tau\|\partial_{y}\mathbf{U}_{s}\|_{L_{\varphi}^{2}(\Omega)}^{2}
                     +\tau\|\mathcal{P}\mathbf{U}_{s}\|_{L_{\varphi}^{2}(\Omega)}^{2}\right)\\
    &+\tau\bar{\kappa}^{2}\|\sqrt{\y}\mathbf{U}_{s+1}\|_{L_{\hat{\varphi}}^{2}(\Omega)}^{2}
     +C\epsilon^{-4}\bar{\eta}\|\mathbb{U}\|_{\hat{H}_{\psi}^{s}(\Omega)}^{2}
     +\tau\|D^{\gamma}F\|_{L_{\psi_{\gamma}}^{2}(\Omega)}^{2}
     +C_{*}\sum_{i=0}^{\frac{s-1}{2}}\|\partial_{t}^{i}F\|_{\mathring{H}^{s-1-2i}(\mathbb{T})}^{2}
     +C\|\partial_{x}^{s}F\|_{\mathring{H}^{0}(\mathbb{T})}^{2},
\end{align*}
taking summation over all $\gamma\in\Gamma$, it follows that
\begin{prop}\label{lowe}
Let $\lambda$ and $\Lambda$ be large enough and $\iota$ small enough such that Lemma \ref{weighte} holds on the time interval $[0,T_{*}]$. Assume that Assumption \ref{MA} holds for $T_{*}$, then for any $\theta>0, \tau>0$ and $t\in[0,T_{*}]$, we have the following a priori estimate for the solution to the problem \eqref{app},
\begin{align}\label{le'}
&\frac{1}{2}\frac{d}{dt}\|u\|_{\hat{H}_{\psi}^{s}(\Omega)}^{2}
+\Lambda\delta\|u\|_{\hat{H}_{\psi,\y}^{s}(\Omega)}^{2}
+\frac{5}{16}\|\partial_{y}u\|_{\hat{H}_{\psi}^{s}(\Omega)}^{2}
+\frac{15}{16}\|\partial_{x}u\|_{\hat{H}_{\psi,\eta}^{s}(\Omega)}^{2}\\
\leq& C_{*}\left(\|\mathbb{U}\|_{\mathcal{H}_{\psi,\varphi}^{s}(\Omega)}
                +\|\partial_{y}\mathbb{U}\|_{\hat{H}_{\psi}^{s}(\Omega)}
                +\|\partial_{y}\mathbf{U}_{s}\|_{L_{\varphi}^{2}(\Omega)}
                +\|\mathcal{P}\mathbf{U}_{s}\|_{L_{\varphi}^{2}(\Omega)}
                +C_{\tau}C_{\theta}\right)
          \|u\|_{\mathcal{H}_{\psi,\varphi}^{s}(\Omega)}^{2}\nonumber\\
    &+C_{*}\|u\|_{\mathcal{H}_{\psi,\varphi}^{s}(\Omega)}^{4}
     +\|\mathcal{P}\mathcal{U}_{s}\|_{L_{\varphi}^{2}(\Omega)}^{2}
     +C_{*}C_{\tau}C_{\theta}\|u\|_{\hat{H}_{\psi,\y}^{s}(\Omega)}^{2}
     +C\tau\bar{\kappa}^{2}\|\sqrt{\y}\mathbf{U}_{s+1}\|_{L_{\hat{\varphi}}^{2}(\Omega)}^{2}\nonumber\\
    &+C\bar{\kappa}^{2}\left(\|\mathbf{U}_{s+1}\|_{L_{\hat{\varphi}}^{2}(\Omega)}^{2}
                     +\tau\|\mathcal{P}\mathbf{U}_{s+1}\|_{L_{\hat{\varphi}}^{2}(\Omega)}^{2}
                     +\tau\|\partial_{y}\mathbf{U}_{s+1}\|_{L_{\hat{\varphi}}^{2}(\Omega)}^{2}
                     +\tau\|\partial_{y}\mathbf{U}_{s}\|_{L_{\varphi}^{2}(\Omega)}^{2}
                     +\tau\|\mathcal{P}\mathbf{U}_{s}\|_{L_{\varphi}^{2}(\Omega)}^{2}\right)\nonumber\\
    &+\bar{\kappa}^{2}\left(\|\mathbb{U}\|_{\mathcal{H}_{\psi,\varphi}^{s}(\Omega)}^{2}
                +C\tau\|\partial_{y}\mathbb{U}\|_{\hat{H}_{\psi}^{s}(\Omega)}^{2}
                +C_{\tau}C_{*}\right)
     +C\epsilon^{-4}\bar{\eta}\|\mathbb{U}\|_{\hat{H}_{\psi}^{s}(\Omega)}^{2}
     +\tau\|F\|_{\hat{H}_{\psi}^{s}(\Omega)}^{2}\nonumber\\
    &+C_{*}\sum_{i=0}^{\frac{s-1}{2}}\|\partial_{t}^{i}F\|_{\mathring{H}^{s-1-2i}(\mathbb{T})}^{2}
     +C\|\partial_{x}^{s}F\|_{\mathring{H}^{0}(\mathbb{T})}^{2}.\nonumber
\end{align}
\end{prop}

Since the right hand side of \eqref{le'} involves the high-order derivative term $\mathcal{U}_{s}$, to close the estimate we need to establish the uniform estimate of it. This will be done in the following two subsections.

~~~~~~~~~~

\subsection{Uniform $L_{t}^{\infty}$ estimate of the corrected increment solution}~

In establishing the uniform in $\bar{\eta}$ estimate for $\mathcal{U}_{s}$, some terms involve $\mathcal{P}$ are difficult to control, since $\mathcal{P}$ tends to $\frac{C}{\sqrt{\bar{\varsigma}}}$ as $t\to0$ and $y\to y_{*}$. Under Assumption \ref{MA}, by observing that $\mathcal{P}\leq\frac{C}{\sqrt{t+\bar{\varsigma}}}$ , we shall obtain a multiplier $\sqrt{t}$ through the following uniform $L_{t}^{\infty}$ estimate of $u$ to help control these terms.
\begin{lemma}\label{Lt}
Assume
$$\|\partial_{t}F\|_{L_{t}^{\infty}H_{\psi}^{s-4}(\Omega_{T_{*}})}+\|F\|_{L_{t}^{\infty}H_{\psi}^{s-2}(\Omega_{T_{*}})}
\leq C_{*}\bar{\eta},$$
then for any $t\in[0,T_{*}]$ we have
\begin{align}\label{ltp}
\|u\|_{H_{\psi}^{s-4}(\hat{\Omega})}(t)
\leq C_{*}t^{\frac{3}{2}}\|u\|_{L_{t}^{2}\mathcal{H}_{\psi,\varphi}^{s}(\Omega_{t})}
     +C_{*}\bar{\eta}t
     +C_{*}t\check{c}+\check{c}.
\end{align}
and
\begin{align}\label{lt}
\|u\|_{H_{\psi}^{s-3}(\Omega)}(t)
\leq& C\sqrt{t}\|u\|_{L_{t}^{4}\mathcal{H}_{\psi,\varphi}^{s}(\Omega_{t})}^{2}
     +C_{*}\sqrt{t}\|u\|_{L_{t}^{2}\mathcal{H}_{\psi,\varphi}^{s}(\Omega_{t})}
     +C_{*}\bar{\eta}t
     +C\check{c},
\end{align}
where $\check{c}=\|\check{u}\|_{\hat{H}_{\psi}^{s}(\Omega)}$.

\end{lemma}
\begin{proof}
Let $\gamma=(\gamma_{1},\gamma_{2})$ be any multiindex in $\Gamma_{s-4}$.
Noticing that ${\bf supp}~h(y)\in(\hat{y},+\infty)$, one knows that $u$ obey the following equation in $\hat{\Omega}_{T_{*}}$,
\begin{equation}\label{partu}
\partial_{t}u
+\mathfrak{U}\partial_{x}u+u\partial_{x}\mathfrak{U}
-\partial_{y}^{-1}[\partial_{x}u](\partial_{y}\mathfrak{U}+\varsigma)
-\partial_{y}^{-1}[\partial_{x}\mathfrak{U}]\partial_{y}u
=\partial_{y}^{2}u+\eta\partial_{x}^{2}u+\kappa\partial_{x}^{2}\mathbb{U}+F,\\
\end{equation}
with the initial condition $u|_{t=0}=\check{u}(x,y)$.
For $k\in\{1,2\}$, by acting $\partial_{t}^{k-1}D^{\gamma}$ on \eqref{partu} one can obtain
\begin{align}\label{partde}
\partial_{t}^{k}D^{\gamma}u
=& \partial_{y}^{2}\partial_{t}^{k-1}D^{\gamma}u
  +\partial_{t}^{k-1}D^{\gamma}[\eta\partial_{x}^{2}u]
  -\partial_{t}^{k-1}D^{\gamma}[\mathfrak{U}\partial_{x}u]
  -\partial_{t}^{k-1}D^{\gamma}[u\partial_{x}\mathfrak{U}]\\
 &+\partial_{t}^{k-1}D^{\gamma}[\partial_{y}^{-1}[\partial_{x}u](\partial_{y}\mathfrak{U}+\varsigma)]
  +\partial_{t}^{k-1}D^{\gamma}[\partial_{y}^{-1}[\partial_{x}\mathfrak{U}]\partial_{y}u]
  +\partial_{t}^{k-1}D^{\gamma}[\kappa\partial_{x}^{2}\mathbb{U}]
  +\partial_{t}^{k-1}D^{\gamma}F.\nonumber
\end{align}
which gives
\begin{align}\label{partdte}
\|\partial_{t}^{k}D^{\gamma}u\|_{L_{\psi}^{2}(\hat{\Omega})}
\leq C_{*}\|u\|_{\hat{H}_{\psi}^{s}(\Omega)}
    +C_{*}\bar{\kappa}
    +C_{*}\sum_{i=0}^{k-1}\|\partial_{t}^{i}F\|_{H_{\psi}^{s-2-2i}(\Omega)},~\forall t\in[0,T_{*}], k=1,2.
\end{align}
By using Taylor's formula, for any $t\in[0,T_{*}]$ we have
\begin{align*}
D^{\gamma}u(t,x,y)
= \int_{0}^{t}\partial_{t}^{2}D^{\gamma}u(s,x,y)(t-s)ds
 +\partial_{t}D^{\gamma}u(0,x,y)t
 +D^{\gamma}\check{u}(x,y),
\end{align*}
utilizing \eqref{partdte}, it follows from the assumption of Lemma \ref{Lt} that for any $t\in[0,T_{*}]$
\begin{align}\label{ltp'}
\|D^{\gamma}u\|_{L_{\psi}^{2}(\hat{\Omega})}(t)
\leq& \int_{0}^{t}\|\partial_{t}^{2}D^{\gamma}u(\sigma,\cdot)\|_{L_{\psi}^{2}(\hat{\Omega})}(t-\sigma)d\sigma
     +\|\partial_{t}D^{\gamma}u(0,x,y)\|_{L_{\psi}^{2}(\hat{\Omega})}t
     +\|D^{\gamma}\check{u}\|_{L_{\psi}^{2}(\hat{\Omega})}\\
\leq& C_{*}\left(t\|u\|_{L_{t}^{1}\hat{H}_{\psi}^{s}(\Omega)}+\bar{\kappa}t^{2}
                 +t^{2}\|\partial_{t}F\|_{L_{t}^{\infty}H_{\psi}^{s-4}(\Omega_{T^{*}})}
                 +t^{2}\|F\|_{L_{t}^{\infty}H_{\psi}^{s-2}(\Omega_{T^{*}})}\right)\nonumber\\
    &+C_{*}t(\|\check{u}\|_{\hat{H}_{\psi}^{s}(\Omega)}+\bar{\kappa}+\|{F}\|_{L_{t}^{\infty}H_{\psi}^{s-2}(\Omega_{T^{*}})})
     +\|\check{u}\|_{\hat{H}_{\psi}^{s}(\Omega)}\nonumber\\
\leq& C_{*}t^{\frac{3}{2}}\|u\|_{L_{t}^{2}\mathcal{H}_{\psi,\varphi}^{s}(\Omega)}
     +C_{*}\bar{\eta}t
     +C_{*}t\check{c}
     +\check{c}.\nonumber
\end{align}
By the arbitrariness of $\gamma$ in the set $\Gamma_{s-4}$, \eqref{ltp} follows from \eqref{ltp'}.

Similarly, for any $\gamma$ in $\Gamma_{s-3}$, one can  act $D^{\gamma}$ on \eqref{app}$_{1}$ to obtain the following equation in $\Omega_{T_{*}}$,
\begin{align*}
\partial_{t}D^{\gamma}u
=& \partial_{y}^{2}D^{\gamma}u
  +D^{\gamma}[\eta\partial_{x}^{2}u]
  -D^{\gamma}[\mathfrak{U}\partial_{x}u]
  -D^{\gamma}[u\partial_{x}\mathfrak{U}]
  +D^{\gamma}[\partial_{y}^{-1}[\partial_{x}u](\partial_{y}\mathfrak{U}+\varsigma)]\\
 &+D^{\gamma}[\partial_{y}^{-1}[\partial_{x}\mathfrak{U}]\partial_{y}u]
  +D^{\gamma}[h\partial_{y}^{-1}[\partial_{x}u]\partial_{y}u]
  +D^{\gamma}[\kappa\partial_{x}^{2}\mathbb{U}]
  +D^{\gamma}F,\nonumber
\end{align*}
which gives
\begin{align}\label{dte}
\|\partial_{t}D^{\gamma}u\|_{L_{\psi}^{2}(\Omega)}
\leq C\|u\|_{\hat{H}_{\psi}^{s}(\Omega)}^{2}
    +C_{*}\|u\|_{\hat{H}_{\psi}^{s}(\Omega)}
    +C_{*}\bar{\kappa}
    +\|F\|_{H_{\psi}^{s-3}(\Omega)},
\end{align}
for any $t\in[0,T_{*}]$. Thus for any $t\in[0,T_{*}]$,
\begin{align}\label{lt'}
\|D^{\gamma}u\|_{L_{\psi}^{2}(\Omega)}(t)
\leq&~\int_{0}^{t}\|\partial_{t}D^{\gamma}u\|_{L_{\psi}^{2}}(\sigma)d\sigma
      +\|D^{\gamma}\check{u}\|_{L_{\psi}^{2}}\\
\leq&~C\|u\|_{L_{t}^{2}\mathcal{H}_{\psi,\varphi}^{s}(\Omega_{t})}^{2}
    +C_{*}\|u\|_{L_{t}^{1}\mathcal{H}_{\psi,\varphi}^{s}(\Omega_{t})}
    +C_{*}\bar{\kappa}t
    +\|F\|_{L_{t}^{1}H_{\psi}^{s-3}(\Omega_{t})}
    +\|D^{\gamma}\check{u}\|_{L_{\psi}^{2}}\nonumber\\
\leq&~C\|u\|_{L_{t}^{2}\mathcal{H}_{\psi,\varphi}^{s}(\Omega_{t})}^{2}
     +C_{*}\|u\|_{L_{t}^{1}\mathcal{H}_{\psi,\varphi}^{s}(\Omega_{t})}
     +C_{*}\bar{\eta}t
     +C\check{c}\nonumber\\
\leq&~C\sqrt{t}\|u\|_{L_{t}^{4}\mathcal{H}_{\psi,\varphi}^{s}(\Omega_{t})}^{2}
     +C_{*}\sqrt{t}\|u\|_{L_{t}^{2}\mathcal{H}_{\psi,\varphi}^{s}(\Omega_{t})}
     +C_{*}\bar{\eta}t
     +C\check{c}.\nonumber
\end{align}
The arbitrariness of $\gamma$ in the set $\Gamma_{s-3}$ and \eqref{lt'} gives \eqref{lt}.
\end{proof}

With Lemma \ref{Lt}, we can obtain the following norm estimate, which will play an important role in the establishment of the uniform estimate for $\mathcal{U}_{s}$.
\begin{corollary}\label{Ltc}
Under the assumption of Lemma \ref{Lt}, for any $t\in[0,T_{*}]$ and $\gamma$ in $\Gamma_{s-6}$, we have
\begin{align}
&\|\varphi\mathcal{P}D^{\gamma}u\|_{L_{x}^{\infty}L_{y}^{2}(\Omega)}(t)+\|\varphi \mathcal{P}^{2}D^{\gamma}u\|_{L_{x}^{\infty}L_{y}^{2}(\Omega)}(t)\\
\leq&~C_{*}\sqrt{t}\|u\|_{L_{t}^{2}\mathcal{H}_{\psi,\varphi}^{s}(\Omega_{t})}
     +C\|u\|_{\hat{H}_{\psi}^{s}(\Omega)}
     +C_{*}\bar{\eta}
     +C_{*}\check{c}\iota^{-1}.\nonumber
\end{align}
\end{corollary}
\begin{proof}
Since $\bar{\varsigma}=\epsilon_{0}^{4}\iota$ when $\check{c}\neq0$, by using Sobolev embedding estimate and \eqref{ltp}, for any $t\in[0,T_{*}]$ and $\gamma$ in $\Gamma_{s-6}$ we have
\begin{align*}
&\|\varphi\mathcal{P}D^{\gamma}u\|_{L_{x}^{\infty}L_{y}^{2}(\Omega)}(t)
+\|\varphi\mathcal{P}^{2}D^{\gamma}u\|_{L_{x}^{\infty}L_{y}^{2}(\Omega)}(t)\\
\leq&~C_{*}\left(\|\varphi\mathcal{P}\|_{L^{2}(0,\hat{y})}+\|\varphi\mathcal{P}^{2}\|_{L^{2}(0,\hat{y})}\right)
          \|D^{\gamma}u(t,\cdot)\|_{L^{\infty}(\hat{\Omega})}
     +C\|\psi D^{\gamma}u(t,\cdot)\|_{L^{\infty}(\Omega\setminus\hat{\Omega})}\\
\leq&~\frac{C_{*}}{t+\bar{\varsigma}}\|u(t,\cdot)\|_{H_{\psi}^{s-4}(\hat{\Omega})}
     +C\|D^{\gamma}u(t,\cdot)\|_{L^{\infty}(\Omega\setminus\hat{\Omega})}\\
\leq&~C_{*}\sqrt{t}\|u\|_{L_{t}^{2}\mathcal{H}_{\psi,\varphi}^{s}(\Omega_{t})}
     +C\|u\|_{\hat{H}_{\psi}^{s}(\Omega)}
     +C_{*}\bar{\eta}
     +C_{*}\check{c}\iota^{-1}.
\end{align*}
\end{proof}
~~~~~~~~~~

\subsection{Uniform estimate with high-order $x$-derivatives}~

In this subsection, we shall establish the uniform (in $\bar{\eta}$) estimate of $\mathcal{U}_{s}$ with the help of Lemma \ref{Lt} and Corollary \ref{Ltc} obtained in subsection 4.3. The problem obeyed by $\mathcal{U}_{s}$ can be derived from the problem \eqref{app}. Since the function $\mathfrak{U}$ in the problem \eqref{app} varies with the iteration number $n$, that is, $\mathfrak{U}=u_{0}$ when $n=1$, and  for $n\geq 2$,  $\mathfrak{U}$ is the solution to the problem \eqref{App}, the treatment of the estimate for $\mathcal{U}_{s}$ will be quite different. In what follows, we discuss these two case separately.

\noindent{\bf \underline{Case 1: $n=1$.}}
According to our setting, $\mathfrak{U}=u_{0}$ and $\bar{\varsigma}=\iota\epsilon_{0}^{4}$.
Consequently, $\mathcal{P}\leq \frac{C}{\sqrt{t+\bar{\varsigma}}}\leq C\iota^{-\frac{1}{2}}$ and $\varphi\leq C_{\lambda}\iota^{-\frac{\varepsilon_{0}}{2}}$ in $\Omega_{T_{*}}$. The estimate for $\mathcal{U}_{s}$ is relatively straightforward.

From \eqref{app} one knows that $\mathcal{U}_{s}$ obeys the following problem in $\Omega_{T_{*}}$,
\begin{equation}\label{us0}
\begin{cases}
\partial_{t}\mathcal{U}_{s}
+\mathcal{L}_{1}^{o}\partial_{x}\mathcal{U}_{s}
+\mathcal{L}_{2}^{o}\partial_{y}\mathcal{U}_{s}
+\mathcal{L}_{0}^{o}\mathcal{U}_{s}
-\partial_{y}^{2}\mathcal{U}_{s}-\eta\partial_{x}^{2}\mathcal{U}_{s}
=\mathbf{F}_{s}+\sum\limits_{i=1}^{6}\mathcal{R}_{i}^{o},\\
(\partial_{y}\mathcal{U}_{s}+\rho\mathcal{U}_{s}+\partial_{x}^{s}F)|_{y=0}=0,~
\lim\limits_{y\to+\infty} \mathcal{U}_{s}=0,\\
\mathcal{U}_{s}|_{t=0}=\mathcal{P}(0,x,y)(\partial_{y}\partial_{x}^{s}\tilde{u}_{in}-\rho(0,x,y)\partial_{x}^{s}\tilde{u}_{in}),
\end{cases}
\end{equation}
where $\mathbf{F}_{s}$ is defined as in Assumption \ref{MA},
\begin{align*}
\mathcal{L}_{0}^{o}
= ~s\partial_{x}\mathfrak{U}-2\partial_{y}\rho
  -2\left(\frac{\partial_{y}\mathcal{P}}{\mathcal{P}}\right)^{2}-2\eta\left(\frac{\partial_{x}\mathcal{P}}{\mathcal{P}}\right)^{2}
  -\frac{\partial_{t}\mathcal{P}+\mathfrak{U}\partial_{x}\mathcal{P}-\partial_{y}^{-1}[\partial_{x}\mathfrak{U}]\partial_{y}\mathcal{P}-\partial_{y}^{2}\mathcal{P}-\eta\partial_{x}^{2}\mathcal{P}}{\mathcal{P}},
\end{align*}
and $\mathcal{L}_{1}^{o}=\mathfrak{U}+2\eta\frac{\partial_{x}\mathcal{P}}{\mathcal{P}}$,
$\mathcal{L}_{2}^{o}=2\frac{\partial_{y}\mathcal{P}}{\mathcal{P}}-\partial_{y}^{-1}\left[\partial_{x}\mathfrak{U}\right]$,
and $\mathcal{R}_{1}^{o}=\mathcal{P}(-\partial_{y}\bar{\mathcal{R}}_{1}^{o}+\rho \bar{\mathcal{R}}_{1}^{o})$, with
\begin{align*}
\bar{\mathcal{R}}_{1}^{o}(t,x,y)
=& \sum_{i=0}^{s-2}\binom{s}{i}\partial_{x}^{s-i}\mathfrak{U}\partial_{x}^{i+1}u
  +\sum_{i=0}^{s-1}\binom{s}{i}\partial_{x}^{s-i+1}\mathfrak{U}\partial_{x}^{i}u\\
 &-\sum_{i=0}^{s-2}\binom{s}{i}\partial_{y}^{-1}[\partial_{x}^{i+1}u]\partial_{x}^{s-i}\partial_{y}\mathfrak{U}
  -\sum_{i=0}^{s-1}\binom{s}{i}\partial_{y}^{-1}[\partial_{x}^{s-i+1}\mathfrak{U}]\partial_{x}^{i}\partial_{y}u,
\end{align*}
$\mathcal{R}_{2}^{o}=-\mathcal{P}\bar{\mathcal{R}}_{2}^{o}\partial_{x}^{s}u$, with
\begin{align*}
\bar{\mathcal{R}}_{2}^{o}(t,x,y)
= \partial_{t}\rho+\mathfrak{U}\partial_{x}\rho-\partial_{y}^{-1}[\partial_{x}\mathfrak{U}]\partial_{y}\rho-2\rho\partial_{y}\rho
 -\partial_{y}^{2}\rho-\eta\partial_{x}^{2}\rho+\partial_{x}\partial_{y}\mathfrak{U}-\rho\partial_{x}\mathfrak{U},
\end{align*}
and $\mathcal{R}_{i}(i=3,...,7)$ are given below,
\begin{align*}
\mathcal{R}_{3}^{o}(t,x,y)
=s\mathcal{P}\left(\partial_{x}\partial_{y}^{2}\mathfrak{U}-\rho\partial_{x}\partial_{y}\mathfrak{U}\right)\partial_{y}^{-1}[\partial_{x}^{s}u],
\end{align*}
\begin{align*}
\mathcal{R}_{4}^{o}(t,x,y)
=2\eta\partial_{x}\rho\mathcal{P}\partial_{x}^{s+1}u,
\end{align*}
\begin{align*}
\mathcal{R}_{5}^{o}(t,x,y)
= \mathcal{P}\partial_{y}\partial_{x}^{s}(h\partial_{y}^{-1}[\partial_{x}u]\partial_{y}u)
 -\mathcal{P}\rho\partial_{x}^{s}(h\partial_{y}^{-1}[\partial_{x}u]\partial_{y}u).
\end{align*}
\begin{align*}
\mathcal{R}_{6}^{o}(t,x,y)
= \mathcal{P}\partial_{y}\eta\partial_{x}^{s+2}u.
\end{align*}

Let $\varphi$ be defined as in \eqref{varphi}. For any $t\in[0,T_{*}]$, multiplying \eqref{us0}$_{1}$ by $\varphi^{2} \mathcal{U}_{s}$, then integrating over $\Omega$ and applying integration by parts one has
\begin{align}\label{he0}
&\frac{1}{2}\frac{d}{dt}\|\mathcal{U}_{s}\|_{L_{\varphi}^{2}(\Omega)}^{2}
+\|\partial_{y}\mathcal{U}_{s}\|_{L_{\varphi}^{2}(\Omega)}^{2}
+\|\sqrt{\eta}\partial_{x}\mathcal{U}_{s}\|_{L_{\varphi}^{2}(\Omega)}^{2}\\
=&~\frac{1}{2}(\mathcal{U}_{s},(\partial_{t}\varphi^{2}+\partial_{y}^{2}\varphi^{2})\mathcal{U}_{s})
  -(\mathcal{L}_{1}^{o}\partial_{x}\mathcal{U}_{s}+\mathcal{L}_{2}^{o}\partial_{y}\mathcal{U}_{s}+\mathcal{L}_{0}^{o}\mathcal{U}_{s},\varphi^{2}\mathcal{U}_{s})\nonumber\\
 &+(\mathcal{U}_{s}\partial_{y}\varphi-\varphi\partial_{y}\mathcal{U}_{s},\varphi\mathcal{U}_{s})_{L_{x}^{2}}|_{y=0}
  +\sum_{i=1}^{6}(\mathcal{R}_{i}^{o},\varphi^{2}\mathcal{U}_{s})
  +(\mathbf{F}_{s},\varphi^{2}\mathcal{U}_{s})\nonumber\\
:=&\sum_{i=0}^{9}I_{i}^{o},\nonumber
\end{align}
To get the estimate of $\mathcal{U}_{s}$, we shall estimate $I_{0}^{o},I_{1}^{o},...,I_{9}^{o}$ step by step.

\noindent{\bf \underline{Estimate of~$I_{0}^{o}$.}}
With the help of Lemma \ref{weighte}, we can obtain
\begin{align*}
I_{0}^{o}\leq -\frac{\lambda}{64}\|\sqrt{\omega_{\lambda}}\mathcal{U}_{s}\|_{L_{\varphi}^{2}(\Omega)}^{2}
          -\frac{3}{4}\Lambda\delta\|\sqrt{\y}\mathcal{U}_{s}\|_{L_{\varphi}^{2}(\Omega)}^{2}.
\end{align*}

\noindent{\bf \underline{Estimate of~$I_{1}^{o}$.}}
Note that $\mathcal{L}_{2}^{o}|_{y=0}=0$, by using integration by parts, one can get
\begin{align*}
I_{1}^{o}
=&\frac{1}{2}(\partial_{x}\mathcal{L}_{1}^{o}\mathcal{U}_{s}+\partial_{y}\mathcal{L}_{2}^{o}\mathcal{U}_{s}-2\mathcal{L}_{0}^{o}\mathcal{U}_{s},\varphi^{2}\mathcal{U}_{s})\\
\leq&~(\|\partial_{x}\mathcal{L}_{1}^{o}\|_{L^{\infty}(\Omega)}+\|\partial_{y}\mathcal{L}_{2}^{o}\|_{L^{\infty}(\Omega)}+2\|\mathcal{L}_{0}^{o}\|_{L^{\infty}(\Omega)})
      \|\mathcal{U}_{s}\|_{L_{\varphi}^{2}(\Omega)}^{2}\\
\leq&~ C_{\iota}\|\mathcal{U}_{s}\|_{L_{\varphi}^{2}(\Omega)}^{2}.
\end{align*}

\noindent{\bf \underline{Estimate of~$I_{2}^{o}$.}}
Using the boundary condition in \eqref{us0} and trace estimate one has
\begin{align*}
I_{2}^{o}
=&~((\partial_{y}\varphi+\rho\varphi)\mathcal{U}_{s},\varphi\mathcal{U}_{s})_{L_{x}^{2}}|_{y=0}
  +(\partial_{x}^{s}F,\varphi^{2}\mathcal{U}_{s})_{L_{x}^{2}}|_{y=0}\\
\leq&~C\|\mathcal{U}_{s}\|_{\mathring{H}^{0}(\mathbb{T})}^{2}+\|\partial_{x}^{s}F\|_{\mathring{H}^{0}(\mathbb{T})}^{2}\\
\leq&~C\|\mathcal{U}_{s}\|_{L_{\varphi}^{2}(\Omega)}^{2}
     +\frac{1}{4}\|\partial_{y}\mathcal{U}_{s}\|_{L_{\varphi}^{2}(\Omega)}^{2}
     +\|\partial_{x}^{s}F\|_{\mathring{H}^{0}(\mathbb{T})}^{2}.
\end{align*}

\noindent{\bf \underline{Estimate of~$I_{3}^{o}$.}}
By using Sobolev imbedding inequality, with \eqref{u0} one can deduce that
\begin{align*}
I_{3}^{o}
\leq C_{\lambda}C_{\iota}\|U\|_{H^{s+2}(\mathbb{T})}\|u\|_{\hat{H}_{\psi}^{s}(\Omega)}\|\mathcal{U}_{s}\|_{L_{\varphi}^{2}(\Omega)}
\leq C_{*}\|u\|_{\mathcal{H}_{\psi,\varphi}^{s}(\Omega)}^{2}+C_{\iota}\|\mathcal{U}_{s}\|_{L_{\varphi}^{2}(\Omega)}^{2}.
\end{align*}

\noindent{\bf \underline{Estimate of~$I_{4}^{o}$.}} Utilizing Corollary \ref{ncC2},
\begin{align*}
I_{4}^{o}
\leq&~CC_{\iota}\|\partial_{x}^{s}u\|_{L_{\varphi}^{2}(\Omega)}\|\mathcal{U}_{s}\|_{L_{\varphi}^{2}(\Omega)}\\
\leq&~C_{\lambda}C_{\iota}\|\partial_{x}^{s}u\|_{L_{\hat{\Psi}}^{2}(\Omega)}\|\mathcal{U}_{s}\|_{L_{\varphi}^{2}(\Omega)}\\
\leq&~C_{\lambda}C_{\iota}\|\mathcal{P}\mathcal{U}_{s}\|_{L_{\varphi}^{2}(\Omega)}\|\mathcal{U}_{s}\|_{L_{\varphi}^{2}(\Omega)}\\
\leq&~C_{*}C_{\iota}\|\mathcal{U}_{s}\|_{L_{\varphi}^{2}(\Omega)}^{2}.
\end{align*}
Here and henceforth, we denote $\y^{-\frac{1}{2}}\Psi$ by $\hat{\Psi}$, where $\Psi$ is given below Corollary \ref{ncC1}.

\noindent{\bf \underline{Estimate of~$I_{5}^{o}$.}} Note that $\partial_{x}\partial_{y}^{2}\mathfrak{U}-\rho\partial_{x}\partial_{y}\mathfrak{U}=\partial_{x}\rho(\dot{\mathfrak{U}}+\varsigma)$, also by Corollary \ref{ncC2}, one has
\begin{align*}
I_{5}^{o}
\leq&~C_{\lambda}C_{\iota}\|\varphi(\dot{\mathfrak{U}}+\varsigma)\|_{L_{x}^{\infty}L_{y}^{2}(\Omega)}
       \|\partial_{x}^{s}u\|_{L_{\varphi}^{2}(\Omega)}\|\mathcal{U}_{s}\|_{L_{\varphi}^{2}(\Omega)}\\
\leq&~C_{\lambda}C_{\iota}\|\partial_{x}^{s}u\|_{L_{\hat{\Psi}}^{2}(\Omega)}\|\mathcal{U}_{s}\|_{L_{\varphi}^{2}(\Omega)}\\
\leq&~C_{\lambda}C_{\iota}\|\mathcal{P}\mathcal{U}_{s}\|_{L_{\varphi}^{2}(\Omega)}\|\mathcal{U}_{s}\|_{L_{\varphi}^{2}(\Omega)}\\
\leq&~C_{*}C_{\iota}\|\mathcal{U}_{s}\|_{L_{\varphi}^{2}(\Omega)}^{2}.
\end{align*}

\noindent{\bf \underline{Estimate of~$I_{6}^{o}$.}} With Corollary \ref{ncC3} we can deduce that
\begin{align*}
I_{6}^{o}
\leq&~C\|\vartheta_{c}\partial_{x}^{s+1}u\|_{L_{\varphi}^{2}(\Omega)}\|\mathcal{U}_{s}\|_{L_{\varphi}^{2}(\Omega)}\\
\leq&~C_{\lambda}C_{\iota}\|\vartheta_{c}\partial_{x}^{s+1}u\|_{L_{\hat{\Psi}}^{2}(\Omega)}\|\mathcal{U}_{s}\|_{L_{\varphi}^{2}(\Omega)}\\
\leq&~C_{\lambda}C_{\iota}\left(\|\vartheta_{c}\partial_{x}\mathcal{U}_{s}\|_{L_{\varphi}^{2}(\Omega)}+\|\mathcal{P}\mathcal{U}_{s}\|_{L_{\varphi}^{2}(\Omega)}\right)
                 \|\mathcal{U}_{s}\|_{L_{\varphi}^{2}(\Omega)}\\
\leq&~C_{*}C_{\iota}\|\mathcal{U}_{s}\|_{L_{\varphi}^{2}(\Omega)}^{2}+\frac{1}{4}\|\sqrt{\eta}\partial_{x}\mathcal{U}_{s}\|_{L_{\varphi}^{2}(\Omega)}^{2}.
\end{align*}

\noindent{\bf \underline{Estimate of~$I_{7}^{o}$.}}
In the following estimate for $I_{7}^{o}$, we use $\frac{C}{\sqrt{t+\bar{\varsigma}}}$ to control $\mathcal{P}$ instead of $C$, so that the estimate remains valid for $n\geq 2$.

Since $\mathcal{P}=1$ for $y\in{\bf supp}~h(y)$, by a direct calculation one has
\begin{align*}
I_{7}^{o}
=& (\partial_{y}\partial_{x}^{s}(h\partial_{y}^{-1}[\partial_{x}u]\partial_{y}u),\varphi^{2}\mathcal{U}_{s})
  -(\rho\partial_{x}^{s}(h\partial_{y}^{-1}[\partial_{x}u]\partial_{y}u),\varphi^{2}\mathcal{U}_{s})\\
=& (h'\partial_{x}^{s}(\partial_{y}^{-1}[\partial_{x}u]\partial_{y}u),\varphi^{2}\mathcal{U}_{s})
  +(h\partial_{x}^{s}(\partial_{x}u\partial_{y}u),\varphi^{2}\mathcal{U}_{s})\\
 &+(h\partial_{x}^{s}(\partial_{y}^{-1}[\partial_{x}u]\partial_{y}^{2}u),\varphi^{2}\mathcal{U}_{s})
  -(h\rho\partial_{x}^{s}(\partial_{y}^{-1}[\partial_{x}u]\partial_{y}u),\varphi^{2}\mathcal{U}_{s})\\
=& (\partial_{y}^{-1}[\partial_{x}^{s+1}u](h'\partial_{y}u+h\partial_{y}^{2}u-h\rho\partial_{y}u),\varphi^{2}\mathcal{U}_{s})
  +(h\partial_{y}^{-1}[\partial_{x}u](\partial_{x}^{s}\partial_{y}^{2}u-\rho\partial_{x}^{s}\partial_{y}u),\varphi^{2}\mathcal{U}_{s})\\
 &+s(\partial_{y}^{-1}[\partial_{x}^{s}u](h'\partial_{x}\partial_{y}u+h\partial_{x}\partial_{y}^{2}u-h\partial_{x}\partial_{y}u),\varphi^{2}\mathcal{U}_{s})
  +((h'\partial_{y}^{-1}[\partial_{x}u]+h\partial_{x}u)\partial_{x}^{s}\partial_{y}u,\varphi^{2}\mathcal{U}_{s})\\
 &+(h\partial_{x}^{s+1}u\partial_{y}u,\varphi^{2}\mathcal{U}_{s})
  +s(h\partial_{x}^{s}u\partial_{x}\partial_{y}u,\varphi^{2}\mathcal{U}_{s})
  +\mathfrak{R}\\
=:& \sum_{i=1}^{6}\mathfrak{I}_{i}+\mathfrak{R},
\end{align*}
where
\begin{align*}
\mathfrak{R}
=& \sum_{i=1}^{s-2}\binom{s}{i}(h'\partial_{y}^{-1}[\partial_{x}^{i+1}u]\partial_{x}^{s-i}\partial_{y}u,\varphi^{2}\mathcal{U}_{s})
  +\sum_{i=1}^{s-2}\binom{s}{i}(h\partial_{x}^{i+1}u\partial_{x}^{s-i}\partial_{y}u,\varphi^{2}\mathcal{U}_{s})\\
 &+\sum_{i=1}^{s-2}\binom{s}{i}(h\partial_{y}^{-1}[\partial_{x}^{i+1}u]\partial_{x}^{s-i}\partial_{y}^{2}u,\varphi^{2}\mathcal{U}_{s})
  +\sum_{i=1}^{s-2}\binom{s}{i}(h\rho\partial_{y}^{-1}[\partial_{x}^{i+1}u]\partial_{x}^{s-i}\partial_{y}u,\varphi^{2}\mathcal{U}_{s}).
\end{align*}
Noticing that $\bar{\varsigma}=\epsilon_{0}^{4}\iota$ when $\check{c}\neq0$, ${\bf supp}~h(y)\in(\hat{y},+\infty)$ and $|\partial_{y}\rho|\leq C\mathcal{P}^{2}$  in $\Omega_{T_{*}}$, thanks to Corollary \ref{ncC2}, Corollary \ref{ncC3} and Lemma \ref{Lt}, we can deduce that
\begin{align*}
\mathfrak{I}_{1}
=&~(\partial_{y}^{-1}[\partial_{x}^{s+1}u](h'\partial_{y}u+h\partial_{y}^{2}u-h\rho\partial_{y}u),\varphi^{2}\mathcal{U}_{s})\\
\leq&~C_{*}\|\vartheta_{c}\partial_{x}^{s+1}u\|_{L_{\hat{\Psi}}^{2}(\Omega)}
           (\|\hat{\psi}\partial_{y}u\|_{L_{x}^{\infty}L_{y}^{2}(\Omega)}+\|\hat{\psi}\partial_{y}^{2}u\|_{L_{x}^{\infty}L_{y}^{2}(\Omega)})
           \|\mathcal{U}_{s}\|_{L_{\varphi}^{2}(\Omega)}\\
\leq&~C_{*}\|\vartheta_{c}\partial_{x}^{s+1}u\|_{L_{\hat{\Psi}}^{2}(\Omega)}
           \left[\sqrt{t}\left(\|u\|_{L_{t}^{4}\mathcal{H}_{\psi,\varphi}^{s}(\Omega_{t})}^{2}+\|u\|_{L_{t}^{2}\mathcal{H}_{\psi,\varphi}^{s}(\Omega_{t})}\right)
                 +\bar{\eta}t+\check{c}\iota^{-1}\right]
           \|\mathcal{U}_{s}\|_{L_{\varphi}^{2}(\Omega)}\\
\leq&~C_{*}\left(\bar{\eta}^{-\frac{1}{2}}\|\sqrt{\eta}\partial_{x}\mathcal{U}_{s}\|_{L_{\varphi}^{2}(\Omega)}+\|\mathcal{P}\mathcal{U}_{s}\|_{L_{\varphi}^{2}(\Omega)}\right)
           \left[\|u\|_{L_{t}^{4}\mathcal{H}_{\psi,\varphi}^{s}(\Omega_{t})}^{2}+\|u\|_{L_{t}^{2}\mathcal{H}_{\psi,\varphi}^{s}(\Omega_{t})}
                 +\bar{\eta}\sqrt{t}+\check{c}\iota^{-\frac{3}{2}}\right]
           \|\mathcal{U}_{s}\|_{L_{\varphi}^{2}(\Omega)}\\
\leq&~C_{*}\iota^{-3}\|\mathcal{U}_{s}\|_{L_{\varphi}^{2}(\Omega)}^{2}
     +\frac{1}{16}\|\sqrt{\eta}\partial_{x}\mathcal{U}_{s}\|_{L_{\varphi}^{2}(\Omega)}^{2}
     +\frac{1}{2}\|\mathcal{P}\mathcal{U}_{s}\|_{L_{\varphi}^{2}(\Omega)}^{2}
     +C_{*}\bar{\eta}^{-1}\|u\|_{L_{t}^{4}\mathcal{H}_{\psi,\varphi}^{s}(\Omega_{t})}^{4}\|\mathcal{U}_{s}\|_{L_{\varphi}^{2}(\Omega)}^{2}\\
    &+C_{*}\bar{\eta}^{-1}\|u\|_{L_{t}^{2}\mathcal{H}_{\psi,\varphi}^{s}(\Omega_{t})}^{2}\|\mathcal{U}_{s}\|_{L_{\varphi}^{2}(\Omega)}^{2}.
\end{align*}
and
\begin{align*}
\mathfrak{I}_{2}
=&~(h\partial_{y}^{-1}[\partial_{x}u](\partial_{y}\mathcal{U}_{s}+\partial_{y}\rho\partial_{x}^{s}u),\varphi^{2}\mathcal{U}_{s})\\
\leq&~C_{*}\|\hat{\psi}\partial_{x}u\|_{L_{x}^{\infty}L_{y}^{2}(\Omega)}
      \left(\|\partial_{y}\mathcal{U}_{s}\|_{L_{\varphi}^{2}(\Omega)}+\|\partial_{x}^{s}u\|_{L_{\hat{\Psi}}^{2}(\Omega)}\right)
      \|\mathcal{U}_{s}\|_{L_{\varphi}^{2}(\Omega)}\\
\leq&~C_{*}\left[\sqrt{t}\left(\|u\|_{L_{t}^{4}\mathcal{H}_{\psi,\varphi}^{s}(\Omega_{t})}^{2}+\|u\|_{L_{t}^{2}\mathcal{H}_{\psi,\varphi}^{s}(\Omega_{t})}\right)
                 +\bar{\eta}t+\check{c}\iota^{-1}\right]
      \left(\|\partial_{y}\mathcal{U}_{s}\|_{L_{\varphi}^{2}(\Omega)}+\|\mathcal{P}\mathcal{U}_{s}\|_{L_{\varphi}^{2}(\Omega)}\right)
      \|\mathcal{U}_{s}\|_{L_{\varphi}^{2}(\Omega)}\\
\leq&~C_{*}\iota^{-3}\|\mathcal{U}_{s}\|_{L_{\varphi}^{2}(\Omega)}^{2}
     +\frac{1}{16}\|\partial_{y}\mathcal{U}_{s}\|_{L_{\varphi}^{2}(\Omega)}^{2}
     +C_{*}\|u\|_{L_{t}^{4}\mathcal{H}_{\psi,\varphi}^{s}(\Omega_{t})}^{4}\|\mathcal{U}_{s}\|_{L_{\varphi}^{2}(\Omega)}^{2}
     +C_{*}\|u\|_{L_{t}^{2}\mathcal{H}_{\psi,\varphi}^{s}(\Omega_{t})}^{2}\|\mathcal{U}_{s}\|_{L_{\varphi}^{2}(\Omega)}^{2}.
\end{align*}
Also by using Corollary \ref{ncC2} and Lemma \ref{Lt}, we get
\begin{align*}
\mathfrak{I}_{3}
\leq& C_{*}\|\partial_{x}^{s}u\|_{L_{\hat{\Psi}}^{2}(\Omega)}
           (\|\hat{\psi}\partial_{x}\partial_{y}u\|_{L_{x}^{\infty}L_{y}^{2}(\Omega)}+\|\hat{\psi}\partial_{x}\partial_{y}^{2}u\|_{L_{x}^{\infty}L_{y}^{2}(\Omega)})
           \|\mathcal{U}_{s}\|_{L_{\varphi}^{2}(\Omega)}\\
\leq& C_{*}\|\mathcal{P}\mathcal{U}_{s}\|_{L_{\varphi}^{2}(\Omega)}
           \left[\sqrt{t}\left(\|u\|_{L_{t}^{4}\mathcal{H}_{\psi,\varphi}^{s}(\Omega_{t})}^{2}+\|u\|_{L_{t}^{2}\mathcal{H}_{\psi,\varphi}^{s}(\Omega_{t})}\right)
                 +\bar{\eta}t+\check{c}\iota^{-1}\right]
           \|\mathcal{U}_{s}\|_{L_{\varphi}^{2}(\Omega)}\\
\leq& C_{*}\iota^{-\frac{3}{2}}\|\mathcal{U}_{s}\|_{L_{\varphi}^{2}(\Omega)}^{2}
     +C_{*}\|u\|_{L_{t}^{4}\mathcal{H}_{\psi,\varphi}^{s}(\Omega_{t})}^{4}\|\mathcal{U}_{s}\|_{L_{\varphi}^{2}(\Omega)}^{2}
     +C_{*}\|u\|_{L_{t}^{2}\mathcal{H}_{\psi,\varphi}^{s}(\Omega_{t})}^{2}\|\mathcal{U}_{s}\|_{L_{\varphi}^{2}(\Omega)}^{2}.
\end{align*}
Similarly, since $|\rho|\leq C\mathcal{P}$ in $\Omega_{T_{*}}$, we have
\begin{align*}
\mathfrak{I}_{4}
=&~((h'\partial_{y}^{-1}[\partial_{x}u]+h\partial_{x}u)\partial_{x}^{s}\partial_{y}u,\varphi^{2}\mathcal{U}_{s})\\
=&~((h'\partial_{y}^{-1}[\partial_{x}u]+h\partial_{x}u)(\mathcal{U}_{s}+\rho\partial_{x}^{s}u),\varphi^{2}\mathcal{U}_{s})\\
\leq&~C_{*}\left(\|\hat{\psi}\partial_{x}u\|_{L_{x}^{\infty}L_{y}^{2}(\Omega)}+\|\partial_{x}u\|_{L^{\infty}(\Omega)}\right)
           \left(\|\mathcal{U}_{s}\|_{L_{\varphi}^{2}(\Omega)}+\|\partial_{x}^{s}u\|_{L_{\hat{\Psi}}^{2}(\Omega)}\right)
           \|\mathcal{U}_{s}\|_{L_{\varphi}^{2}(\Omega)}\\
\leq&~C_{*}\left[\sqrt{t}\left(\|u\|_{L_{t}^{4}\mathcal{H}_{\psi,\varphi}^{s}(\Omega_{t})}^{2}+\|u\|_{L_{t}^{2}\mathcal{H}_{\psi,\varphi}^{s}(\Omega_{t})}\right)
                 +\bar{\eta}t+\check{c}\iota^{-1}\right]
           \left(\|\mathcal{U}_{s}\|_{L_{\varphi}^{2}(\Omega)}+\|\mathcal{P}\mathcal{U}_{s}\|_{L_{\varphi}^{2}(\Omega)}\right)
           \|\mathcal{U}_{s}\|_{L_{\varphi}^{2}(\Omega)}\\
\leq&~C_{*}\iota^{-\frac{3}{2}}\|\mathcal{U}_{s}\|_{L_{\varphi}^{2}(\Omega)}^{2}
     +C_{*}\|u\|_{L_{t}^{4}\mathcal{H}_{\psi,\varphi}^{s}(\Omega_{t})}^{4}\|\mathcal{U}_{s}\|_{L_{\varphi}^{2}(\Omega)}^{2}
     +C_{*}\|u\|_{L_{t}^{2}\mathcal{H}_{\psi,\varphi}^{s}(\Omega_{t})}^{2}\|\mathcal{U}_{s}\|_{L_{\varphi}^{2}(\Omega)}^{2}.
\end{align*}
By virtue of Corollary \ref{ncC3} and Lemma \ref{Lt}, we get
\begin{align*}
\mathfrak{I}_{5}
\leq&~C_{*}\|\partial_{y}u\|_{L^{\infty}(\Omega)}\|\partial_{x}^{s+1}u\|_{L_{\hat{\Psi}}^{2}(\Omega)}\|\mathcal{U}_{s}\|_{L_{\varphi}^{2}(\Omega)}\\
\leq&~C_{*}\left[\sqrt{t}\left(\|u\|_{L_{t}^{4}\mathcal{H}_{\psi,\varphi}^{s}(\Omega_{t})}^{2}+\|u\|_{L_{t}^{2}\mathcal{H}_{\psi,\varphi}^{s}(\Omega_{t})}\right)
                 +\bar{\eta}t+\check{c}\iota^{-1}\right]
           \|\partial_{x}^{s+1}u\|_{L_{\hat{\Psi}}^{2}(\Omega)}
           \|\mathcal{U}_{s}\|_{L_{\varphi}^{2}(\Omega)}\\
\leq&~C_{*}\left[\|u\|_{L_{t}^{4}\mathcal{H}_{\psi,\varphi}^{s}(\Omega_{t})}^{2}+\|u\|_{L_{t}^{2}\mathcal{H}_{\psi,\varphi}^{s}(\Omega_{t})}
                 +\bar{\eta}\sqrt{t}+\check{c}\iota^{-\frac{3}{2}}\right]
          \left(\|\partial_{x}\mathcal{U}_{s}\|_{L_{\varphi}^{2}(\Omega)}+\|\mathcal{P}\mathcal{U}_{s}\|_{L_{\varphi}^{2}(\Omega)}\right)
           \|\mathcal{U}_{s}\|_{L_{\varphi}^{2}(\Omega)}\\
\leq&~C_{*}\iota^{-3}\|\mathcal{U}_{s}\|_{L_{\varphi}^{2}(\Omega)}^{2}
     +\frac{1}{16}\|\sqrt{\eta}\partial_{x}\mathcal{U}_{s}\|_{L_{\varphi}^{2}(\Omega)}^{2}
     +\frac{1}{2}\|\mathcal{P}\mathcal{U}_{s}\|_{L_{\varphi}^{2}(\Omega)}^{2}
     +C_{*}\bar{\eta}^{-1}\|u\|_{L_{t}^{4}\mathcal{H}_{\psi,\varphi}^{s}(\Omega_{t})}^{4}\|\mathcal{U}_{s}\|_{L_{\varphi}^{2}(\Omega)}^{2}\\
    &+C_{*}\bar{\eta}^{-1}\|u\|_{L_{t}^{2}\mathcal{H}_{\psi,\varphi}^{s}(\Omega_{t})}^{2}\|\mathcal{U}_{s}\|_{L_{\varphi}^{2}(\Omega)}^{2}.
\end{align*}
Again by Corollary \ref{ncC2} and Lemma \ref{Lt},
\begin{align*}
\mathfrak{I}_{6}
\leq&~C_{*}\|\partial_{x}^{s}u\|_{L_{\hat{\Psi}}^{2}(\Omega)}
           \|\partial_{x}\partial_{y}u\|_{L_{x}^{\infty}L_{y}^{2}(\Omega)}
           \|\mathcal{U}_{s}\|_{L_{\varphi}^{2}(\Omega)}\\
\leq&~C_{*}\|\mathcal{P}\mathcal{U}_{s}\|_{L_{\varphi}^{2}(\Omega)}
           \left[\sqrt{t}\left(\|u\|_{L_{t}^{4}\mathcal{H}_{\psi,\varphi}^{s}(\Omega_{t})}^{2}+\|u\|_{L_{t}^{2}\mathcal{H}_{\psi,\varphi}^{s}(\Omega_{t})}\right)
                 +\bar{\eta}t+\check{c}\iota^{-1}\right]
           \|\mathcal{U}_{s}\|_{L_{\varphi}^{2}(\Omega)}\\
\leq&~C_{*}\iota^{-\frac{3}{2}}\|\mathcal{U}_{s}\|_{L_{\varphi}^{2}(\Omega)}^{2}
     +C_{*}\|u\|_{L_{t}^{4}\mathcal{H}_{\psi,\varphi}^{s}(\Omega_{t})}^{4}\|\mathcal{U}_{s}\|_{L_{\varphi}^{2}(\Omega)}^{2}
     +C_{*}\|u\|_{L_{t}^{2}\mathcal{H}_{\psi,\varphi}^{s}(\Omega_{t})}^{2}\|\mathcal{U}_{s}\|_{L_{\varphi}^{2}(\Omega)}^{2}.
\end{align*}
Finally, it follows from Lemma \ref{Lt} that
\begin{align*}
\mathfrak{R}
\leq&~C_{*}\|u\|_{\hat{H}_{\psi}^{s}(\Omega)}
           \left[\sqrt{t}\left(\|u\|_{L_{t}^{4}\mathcal{H}_{\psi,\varphi}^{s}(\Omega_{t})}^{2}+\|u\|_{L_{t}^{2}\mathcal{H}_{\psi,\varphi}^{s}(\Omega_{t})}\right)
                 +\bar{\eta}t+\check{c}\iota^{-1}\right]
           \|\mathcal{U}_{s}\|_{L_{\varphi}^{2}(\Omega)}\\
\leq&~C_{*}\iota^{-2}\|\mathcal{U}_{s}\|_{L_{\varphi}^{2}(\Omega)}^{2}
     +C\|u\|_{\hat{H}_{\psi}^{s}(\Omega)}^{2}
     +C_{*}\|u\|_{L_{t}^{4}\mathcal{H}_{\psi,\varphi}^{s}(\Omega_{t})}^{4}\|u\|_{\hat{H}_{\psi}^{s}(\Omega)}^{2}
     +C_{*}\|u\|_{L_{t}^{2}\mathcal{H}_{\psi,\varphi}^{s}(\Omega_{t})}^{2}\|u\|_{\hat{H}_{\psi}^{s}(\Omega)}^{2}.
\end{align*}

To sum up, we have
\begin{align*}
I_{7}^{o}
\leq&~C_{*}\iota^{-3}\|\mathcal{U}_{s}\|_{L_{\varphi}^{2}(\Omega)}^{2}
     +C\|u\|_{\hat{H}_{\psi}^{s}(\Omega)}^{2}
     +\frac{1}{8}\|\sqrt{\eta}\partial_{x}\mathcal{U}_{s}\|_{L_{\varphi}^{2}(\Omega)}^{2}
     +\frac{1}{16}\|\partial_{y}\mathcal{U}_{s}\|_{L_{\varphi}^{2}(\Omega)}^{2}
     +\|\mathcal{P}\mathcal{U}_{s}\|_{L_{\varphi}^{2}(\Omega)}^{2}\\
    &+C_{*}\bar{\eta}^{-1}\|u\|_{L_{t}^{4}\mathcal{H}_{\psi,\varphi}^{s}(\Omega_{t})}^{4}\|\mathcal{U}_{s}\|_{L_{\varphi}^{2}(\Omega)}^{2}
     +C_{*}\bar{\eta}^{-1}\|u\|_{L_{t}^{2}\mathcal{H}_{\psi,\varphi}^{s}(\Omega_{t})}^{2}\|\mathcal{U}_{s}\|_{L_{\varphi}^{2}(\Omega)}^{2}.
\end{align*}

\noindent{\bf \underline{Estimate of~$I_{8}^{o}$.}} Since $\mathcal{P}=1$ when $\partial_{y}\eta\neq 0$, by using integration by parts we have
\begin{align*}
I_{8}^{o}
=(\partial_{y}\eta\partial_{x}^{s+2}u,\varphi^{2}\mathcal{U}_{s})
=&-(\partial_{y}\eta\partial_{x}^{s+1}u,\varphi^{2}\partial_{x}\mathcal{U}_{s})\\
\leq&~C\bar{\eta}\theta t\|\vartheta_{c}\partial_{x}^{s+1}u\|_{L_{\varphi}^{2}(\Omega)}
                         \|\vartheta_{c}\partial_{x}\mathcal{U}_{s}\|_{L_{\varphi}^{2}(\Omega)}\\
\leq&~C_{\lambda}\bar{\eta}\theta t\|\vartheta_{c}\partial_{x}^{s+1}u\|_{L_{\hat{\Psi}}^{2}(\Omega)}
                                   \|\vartheta_{c}\partial_{x}\mathcal{U}_{s}\|_{L_{\varphi}^{2}(\Omega)}\\
\leq&~C_{\lambda}\bar{\eta}\theta t
     \left(\|\vartheta_{c}\partial_{x}\mathcal{U}_{s}\|_{L_{\varphi}^{2}(\Omega)}+\|\mathcal{P}\mathcal{U}_{s}\|_{L_{\varphi}^{2}(\Omega)}\right)
     \|\vartheta_{c}\partial_{x}\mathcal{U}_{s}\|_{L_{\varphi}^{2}(\Omega)}\\
\leq&~C_{\lambda}\theta t\|\sqrt{\eta}\partial_{x}\mathcal{U}_{s}\|_{L_{\varphi}^{2}(\Omega)}^{2}
     +C_{*}\|\mathcal{U}_{s}\|_{L_{\varphi}^{2}(\Omega)}^{2}.
\end{align*}

\noindent{\bf \underline{Estimate of~$I_{9}^{o}$.}} Applying the Cauchy inequality, one has
\begin{align*}
I_{9}^{o}\leq \tau\|\mathbf{F}_{s}\|_{L_{\hat{\psi}}^{2}}+C_{*}C_{\tau}\|\mathcal{U}_{s}\|_{L_{\varphi}^{2}(\Omega)}^{2}.
\end{align*}

Combining the above estimates for $I_{0}^{o},I_{1}^{o},...,I_{9}^{o}$, from \eqref{he0} we conclude that when $n=1$, for any $t\in[0,T_{*}]$ the following estimate for $\mathcal{U}_{s}$ holds,
\begin{align}\label{He0}
&\frac{1}{2}\frac{d}{dt}\|\mathcal{U}_{s}\|_{L_{\varphi}^{2}(\Omega)}^{2}
+\frac{\lambda}{64}\|\sqrt{\omega_{\lambda}}\mathcal{U}_{s}\|_{L_{\varphi}^{2}(\Omega)}^{2}
+\frac{3}{4}\Lambda\delta\|\sqrt{\y}\mathcal{U}_{s}\|_{L_{\varphi}^{2}(\Omega)}^{2}
+\frac{1}{2}\|\partial_{y}\mathcal{U}_{s}\|_{L_{\varphi}^{2}(\Omega)}^{2}
+\frac{5}{8}\|\sqrt{\eta}\partial_{x}\mathcal{U}_{s}\|_{L_{\varphi}^{2}(\Omega)}^{2}\\
\leq&~C_{*}\|u\|_{\mathcal{H}_{\psi,\varphi}^{s}(\Omega)}^{2}
     +C_{*}C_{\tau}C_{\iota}\|\mathcal{U}_{s}\|_{L_{\varphi}^{2}(\Omega)}^{2}
     +C\|u\|_{\hat{H}_{\psi}^{s}(\Omega)}^{2}
     +\|\mathcal{P}\mathcal{U}_{s}\|_{L_{\varphi}^{2}(\Omega)}^{2}
     +C_{\lambda}\theta t\|\sqrt{\eta}\partial_{x}\mathcal{U}_{s}\|_{L_{\varphi}^{2}(\Omega)}^{2}\nonumber\\
    &+C_{*}\bar{\eta}^{-1}\|u\|_{L_{t}^{4}\mathcal{H}_{\psi,\varphi}^{s}(\Omega_{t})}^{4}\|\mathcal{U}_{s}\|_{L_{\varphi}^{2}(\Omega)}^{2}
     +C_{*}\bar{\eta}^{-1}\|u\|_{L_{t}^{2}\mathcal{H}_{\psi,\varphi}^{s}(\Omega_{t})}^{2}\|\mathcal{U}_{s}\|_{L_{\varphi}^{2}(\Omega)}^{2}
     +\|\partial_{x}^{s}F\|_{\mathring{H}^{0}(\mathbb{T})}^{2}
     +\tau\|\mathbf{F}_{s}\|_{L_{\hat{\psi}}^{2}}.\nonumber
\end{align}

\noindent{\bf \underline{Case 2: $n\geq 2$.}}
Note that $\check{u}=0$ when $n\ge 2$, from \eqref{app} and \eqref{App} one knows that $\mathcal{U}_{s}$ obeys the following problem in $\Omega_{T_{*}}$,
\begin{equation}\label{us}
\begin{cases}
\partial_{t}\mathcal{U}_{s}
+\mathcal{L}_{1}\partial_{x}\mathcal{U}_{s}
+\mathcal{L}_{2}\partial_{y}\mathcal{U}_{s}
+\mathcal{L}_{0}\mathcal{U}_{s}
-\partial_{y}^{2}\mathcal{U}_{s}-\eta\partial_{x}^{2}\mathcal{U}_{s}
=\mathbf{F}_{s}+\sum\limits_{i=1}^{7}\mathcal{R}_{i},\\
(\partial_{y}\mathcal{U}_{s}+\rho\mathcal{U}_{s}+\partial_{x}^{s}F)|_{y=0}=0,~
\lim\limits_{y\to+\infty} \mathcal{U}_{s}=0,\\
\mathcal{U}_{s}|_{t=0}=0,
\end{cases}
\end{equation}
where $\mathbf{F}_{s}$ is defined as in Assumption \ref{MA},
\begin{align*}
\mathcal{L}_{0}
=&~s\partial_{x}\mathfrak{U}-2\partial_{y}\rho
  -2\left(\frac{\partial_{y}\mathcal{P}}{\mathcal{P}}\right)^{2}-2\eta\left(\frac{\partial_{x}\mathcal{P}}{\mathcal{P}}\right)^{2}
  +\frac{3\chi}{4\mathcal{P}\sqrt{\dot{\mathfrak{U}}+\varsigma}}\left[\rho^{2}+\eta\left(\frac{\partial_{x}\dot{\mathfrak{U}}}{\dot{\mathfrak{U}}+\varsigma}\right)^{2}\right]\\
 &-\frac{\partial_{y}^{-1}\left[\partial_{x}\mathfrak{U}\right]\varsigma'+\varsigma''-\partial_{y}\eta\partial_{x}^{2}\mathbb{U}+\eta\partial_{x}^{2}\partial_{y}u_{0}+\partial_{y}F+\partial_{y}[\kappa\partial_{x}^{2}\mathbb{U}]}{2(\dot{\mathfrak{U}}+\varsigma)}\\
 &+(1-\chi)\frac{\partial_{y}^{-1}\left[\partial_{x}\mathfrak{U}\right]\varsigma'+\varsigma''-\partial_{y}\eta\partial_{x}^{2}\mathbb{U}+\eta\partial_{x}^{2}\partial_{y}u_{0}+\partial_{y}F+\partial_{y}[\kappa\partial_{x}^{2}\mathbb{U}]}{2\mathcal{P}(\dot{\mathfrak{U}}+\varsigma)}\\
 &-\frac{-\partial_{y}^{-1}[\partial_{x}\mathfrak{U}]\chi'+\chi'\rho-\chi''}{\mathcal{P}\sqrt{\dot{\mathfrak{U}}+\varsigma}}
  -\frac{\partial_{y}^{-1}\left[\partial_{x}\mathfrak{U}\right]\chi'+\chi''}{\mathcal{P}},
\end{align*}
and $\mathcal{L}_{1}=\mathfrak{U}+2\eta\frac{\partial_{x}\mathcal{P}}{\mathcal{P}}$,
$\mathcal{L}_{2}=2\frac{\partial_{y}\mathcal{P}}{\mathcal{P}}-\partial_{y}^{-1}\left[\partial_{x}\mathfrak{U}\right]$,
and $\mathcal{R}_{1}=\mathcal{P}(-\partial_{y}\bar{\mathcal{R}}_{1}+\rho \bar{\mathcal{R}}_{1})$, with
\begin{align*}
\bar{\mathcal{R}}_{1}(t,x,y)
=& \sum_{i=0}^{s-2}\binom{s}{i}\partial_{x}^{s-i}\mathfrak{U}\partial_{x}^{i+1}u
  +\sum_{i=0}^{s-1}\binom{s}{i}\partial_{x}^{s-i+1}\mathfrak{U}\partial_{x}^{i}u\\
 &-\sum_{i=0}^{s-2}\binom{s}{i}\partial_{y}^{-1}[\partial_{x}^{i+1}u]\partial_{x}^{s-i}\partial_{y}\mathfrak{U}
  -\sum_{i=0}^{s-1}\binom{s}{i}\partial_{y}^{-1}[\partial_{x}^{s-i+1}\mathfrak{U}]\partial_{x}^{i}\partial_{y}u,
\end{align*}
$\mathcal{R}_{2}=-\mathcal{P}\bar{\mathcal{R}}_{2}\partial_{x}^{s}u$, with
\begin{align*}
\bar{\mathcal{R}}_{2}(t,x,y)
=& 2\eta\frac{\partial_{x}\dot{\mathfrak{U}}}{\dot{\mathfrak{U}}+\varsigma}\partial_{x}\rho
  +\rho\cdot\frac{\partial_{y}F+\partial_{y}[\kappa\partial_{x}^{2}\mathbb{U}]+\partial_{y}^{-1}[\partial_{x}\mathfrak{U}]\varsigma'+\varsigma''}{\dot{\mathfrak{U}}+\varsigma}\\
 &-\frac{\partial_{y}^{2}F+\partial_{y}^{2}[\kappa\partial_{x}^{2}\mathbb{U}]-\varsigma\partial_{x}\dot{\mathfrak{U}}+\varsigma'\partial_{x}\mathfrak{U}+\varsigma'''+\partial_{y}^{-1}[\partial_{x}\mathfrak{U}]\varsigma''}{\dot{\mathfrak{U}}+\varsigma},
\end{align*}
and $\mathcal{R}_{i}(i=3,...,7)$ are given below,
\begin{align*}
\mathcal{R}_{3}(t,x,y)
=s\mathcal{P}\left(\partial_{x}\partial_{y}^{2}\mathfrak{U}-\rho\partial_{x}\partial_{y}\mathfrak{U}\right)\partial_{y}^{-1}[\partial_{x}^{s}u],
\end{align*}
\begin{align*}
\mathcal{R}_{4}(t,x,y)
=(\varsigma+2\eta\partial_{x}\rho)\mathcal{P}\partial_{x}^{s+1}u,
\end{align*}
\begin{align*}
\mathcal{R}_{5}(t,x,y)
=\kappa\mathcal{P}\partial_{y}\partial_{x}^{s+2}\mathbb{U}
 -\kappa\mathcal{P}\rho\partial_{x}^{s+2}\mathbb{U},
\end{align*}
\begin{align*}
\mathcal{R}_{6}(t,x,y)
= \mathcal{P}\partial_{y}\partial_{x}^{s}(h\partial_{y}^{-1}[\partial_{x}u]\partial_{y}u)
 -\mathcal{P}\rho\partial_{x}^{s}(h\partial_{y}^{-1}[\partial_{x}u]\partial_{y}u).
\end{align*}
\begin{align*}
\mathcal{R}_{7}(t,x,y)
= \mathcal{P}\partial_{y}\eta\partial_{x}^{s+2}u
 +\mathcal{P}\partial_{y}\kappa\partial_{x}^{s+2}\mathbb{U}.
\end{align*}

Let $\varphi$ be defined as in \eqref{varphi}.
For any $t\in[0,T_{*}]$, multiplying \eqref{us}$_{1}$ by $\varphi^{2} \mathcal{U}_{s}$, then integrating over $\Omega$ and applying integration by parts one has
\begin{align}\label{he}
&\frac{1}{2}\frac{d}{dt}\|\mathcal{U}_{s}\|_{L_{\varphi}^{2}(\Omega)}^{2}
+\|\partial_{y}\mathcal{U}_{s}\|_{L_{\varphi}^{2}(\Omega)}^{2}
+\|\sqrt{\eta}\partial_{x}\mathcal{U}_{s}\|_{L_{\varphi}^{2}(\Omega)}^{2}\\
=&~\frac{1}{2}(\mathcal{U}_{s},(\partial_{t}\varphi^{2}+\partial_{y}^{2}\varphi^{2})\mathcal{U}_{s})
  -(\mathcal{L}_{1}\partial_{x}\mathcal{U}_{s},\varphi^{2}\mathcal{U}_{s})
  -(\mathcal{L}_{2}\partial_{y}\mathcal{U}_{s},\varphi^{2}\mathcal{U}_{s})
  -(\mathcal{L}_{0}\mathcal{U}_{s},\varphi^{2}\mathcal{U}_{s})\nonumber\\
 &+(\mathcal{U}_{s},\varphi\partial_{y}\varphi\mathcal{U}_{s})_{L_{x}^{2}}|_{y=0}
  -(\partial_{y}\mathcal{U}_{s},\varphi^{2}\mathcal{U}_{s})_{L_{x}^{2}}|_{y=0}
  +\sum_{i=1}^{7}(\mathcal{R}_{i},\varphi^{2}\mathcal{U}_{s})
  +(\mathbf{F}_{s},\varphi^{2}\mathcal{U}_{s})\nonumber\\
:=&\sum_{i=0}^{13}I_{i},\nonumber
\end{align}

To get the estimate of $\mathcal{U}_{s}$, we shall conduct an item-by-item analysis hereinafter.

\noindent{\bf \underline{Estimate of~$I_{0}$.}}
With the help of Lemma \ref{weighte}, we can obtain
\begin{align*}
I_{0}\leq -\frac{\lambda}{64}\|\sqrt{\omega_{\lambda}}\mathcal{U}_{s}\|_{L_{\varphi}^{2}(\Omega)}^{2}
          -\frac{3}{4}\Lambda\delta\|\sqrt{\y}\mathcal{U}_{s}\|_{L_{\varphi}^{2}(\Omega)}^{2}.
\end{align*}

\noindent{\bf \underline{Estimate of~$I_{1}$.}}
Applying integration by parts, one has
\begin{align*}
I_{1}
=&~\frac{1}{2}(\partial_{x}\mathcal{L}_{1}\mathcal{U}_{s},\varphi^{2}\mathcal{U}_{s})\\
=&~\frac{1}{2}(\partial_{x}\mathfrak{U}\mathcal{U}_{s},\varphi^{2}\mathcal{U}_{s})
  +\left(\eta\partial_{x}\left(\frac{\partial_{x}\mathcal{P}}{\mathcal{P}}\right)\mathcal{U}_{s},\varphi^{2}\mathcal{U}_{s}\right)\\
\leq&~C_{*}\|\mathcal{U}_{s}\|_{L_{\varphi}^{2}(\Omega)}^{2}
    -\frac{1}{2}\int_{\mathbb{T}}\int_{\check{y}_{*}}^{\hat{y}_{*}}\eta\varphi^{2}\left[\frac{\partial_{x}^{2}\dot{\mathfrak{U}}}{\dot{\mathfrak{U}}+\varsigma}-\left(\frac{\partial_{x}\dot{\mathfrak{U}}}{\dot{\mathfrak{U}}+\varsigma}\right)^{2}\right]\mathcal{U}_{s}^{2}dydx\nonumber\\
   &+\int_{\mathbb{T}}\int_{[\check{y},\hat{y}]\setminus[\check{y}_{*},\hat{y}_{*}]}\eta\varphi^{2}\partial_{x}\left(\frac{\partial_{x}\mathcal{P}}{\mathcal{P}}\right)\mathcal{U}_{s}^{2}dydx\nonumber\\
\leq&~C_{*}\|\mathcal{U}_{s}\|_{L_{\varphi}^{2}(\Omega)}^{2},\nonumber
\end{align*}
since $|\partial_{x}\dot{\mathfrak{U}}|\leq C\y\dot{\mathfrak{U}}$, $|\partial_{x}^{2}\dot{\mathfrak{U}}|\leq C\y\dot{\mathfrak{U}}$ on $[\check{y}_{*},\hat{y}_{*}]$, and $\left|\partial_{x}\left(\frac{\partial_{x}\mathcal{P}}{\mathcal{P}}\right)\right|\leq C_{*}$ on $[\check{y},\hat{y}]\setminus[\check{y}_{*},\hat{y}_{*}]$.

\noindent{\bf \underline{Estimate of~$I_{2}$.}}
Recall the definition of $\mathcal{L}_{2}$, one has
\begin{align}\label{I20}
I_{2}= \left(\partial_{y}^{-1}[\partial_{x}\mathfrak{U}]\partial_{y}\mathcal{U}_{s},\varphi^{2}\mathcal{U}_{s}\right)
      -2\left(\frac{\partial_{y}\mathcal{P}}{\mathcal{P}}\partial_{y}\mathcal{U}_{s},\varphi^{2}\mathcal{U}_{s}\right)
    :=I_{21}+I_{22}.
\end{align}
Applying integration by parts and \eqref{varphiy}, one can deduce that
\begin{align}\label{I21}
I_{21}
=&~-\frac{1}{2}(\partial_{x}\mathfrak{U}\mathcal{U}_{s},\varphi^{2}\mathcal{U}_{s})
  -\frac{1}{2}(\partial_{y}^{-1}\left[\partial_{x}\mathfrak{U}\right]\mathcal{U}_{s},\partial_{y}\varphi^{2}\mathcal{U}_{s})\\
=&~-\frac{1}{2}(\partial_{x}\mathfrak{U}\mathcal{U}_{s},\varphi^{2}\mathcal{U}_{s})
  -\frac{1}{2}(\chi\partial_{y}^{-1}\left[\partial_{x}\mathfrak{U}\right]\mathcal{U}_{s},\partial_{y}\varphi^{2}\mathcal{U}_{s})
  -\frac{1}{2}((1-\chi)\partial_{y}^{-1}\left[\partial_{x}\mathfrak{U}\right]\mathcal{U}_{s},\partial_{y}\varphi^{2}\mathcal{U}_{s})\nonumber\\
\leq&~C_{*}\|\mathcal{U}_{s}\|_{L_{\varphi}^{2}(\Omega)}^{2}
     +C_{*}(\mathcal{U}_{s},\sqrt{\omega_{\lambda}}\varphi^{2}\mathcal{U}_{s})
     +C_{*}(\mathcal{U}_{s},\varphi^{2}\mathcal{U}_{s})
     +C_{*}\|\y^{-1}\partial_{y}^{-1}\left[\partial_{x}\mathfrak{U}\right]\|_{L^{\infty}(\Omega)}
       \|\sqrt{\y}\mathcal{U}_{s}\|_{L_{\varphi}^{2}(\Omega)}^{2}\nonumber\\
\leq&~C_{*}\|\mathcal{U}_{s}\|_{L_{\varphi}^{2}(\Omega)}^{2}
     +\frac{1}{16}\|\sqrt{\omega_{\lambda}}\mathcal{U}_{s}\|_{L_{\varphi}^{2}(\Omega)}^{2}
     +C_{*}\|\sqrt{\y}\mathcal{U}_{s}\|_{L_{\varphi}^{2}(\Omega)}^{2}.\nonumber
\end{align}

In view of the fact ${\bf supp}~\partial_{y}\mathcal{P}\subset(\check{y},\hat{y})$, one knows
\begin{align}\label{I22}
I_{22}
=-2\int_{\mathbb{T}}\int_{\check{y}_{*}}^{\hat{y}_{*}}\varphi^{2}\frac{\partial_{y}\mathcal{P}}{\mathcal{P}}\mathcal{U}_{s}\partial_{y}\mathcal{U}_{s}dydx
 -2\int_{\mathbb{T}}\int_{[\check{y},\hat{y}]\setminus[\check{y}_{*},\hat{y}_{*}]}\varphi^{2}\frac{\partial_{y}\mathcal{P}}{\mathcal{P}}\mathcal{U}_{s}\partial_{y}\mathcal{U}_{s}dydx.
\end{align}
Since $\frac{\partial_{y}\mathcal{P}}{\mathcal{P}}=-\frac{1}{2}\rho$
and $\left|\frac{\partial_{y}\mathcal{P}}{\mathcal{P}}\right|\leq C_{*}$ on $[\check{y},\hat{y}]\setminus[\check{y}_{*},\hat{y}_{*}]$, one has
\begin{align}\label{I221}
-2\int_{\mathbb{T}}\int_{\check{y}_{*}}^{\hat{y}_{*}}\varphi^{2}\frac{\partial_{y}\mathcal{P}}{\mathcal{P}}\mathcal{U}_{s}\partial_{y}\mathcal{U}_{s}dydx
\leq&~\frac{1}{2}\|\partial_{y}\mathcal{U}_{s}\|_{L_{\varphi}^{2}(\Omega)}^{2}
     +2\|\rho\mathcal{U}_{s}\|_{L_{\varphi}^{2}(\Omega)}^{2}\\
\leq&~\frac{1}{2}\|\partial_{y}\mathcal{U}_{s}\|_{L_{\varphi}^{2}(\Omega)}^{2}
     +C\|\mathcal{P}\mathcal{U}_{s}\|_{L_{\varphi}^{2}(\Omega)}^{2},\nonumber
\end{align}
and
\begin{align}\label{I222}
-2\int_{\mathbb{T}}\int_{[\check{y},\hat{y}]\setminus[\check{y}_{*},\hat{y}_{*}]}\varphi^{2}\frac{\partial_{y}\mathcal{P}}{\mathcal{P}}\mathcal{U}_{s}\partial_{y}\mathcal{U}_{s}dydx
\leq \frac{1}{16}\|\partial_{y}\mathcal{U}_{s}\|_{L_{\varphi}^{2}(\Omega)}^{2}
    +C_{*}\|\mathcal{U}_{s}\|_{L_{\varphi}^{2}(\Omega)}^{2}.
\end{align}
Substituting \eqref{I221} and \eqref{I222} into \eqref{I22}, one has
\begin{align}\label{I22'}
I_{22}
\leq \frac{9}{16}\|\partial_{y}\mathcal{U}_{s}\|_{L_{\varphi}^{2}(\Omega)}^{2}
    +C_{*}\|\mathcal{U}_{s}\|_{L_{\varphi}^{2}(\Omega)}^{2}
    +C\|\mathcal{P}\mathcal{U}_{s}\|_{L_{\varphi}^{2}(\Omega)}^{2}.
\end{align}

Combining \eqref{I21} and \eqref{I22'} with \eqref{I20} one has
\begin{align*}
I_{2}
\leq~C_{*}\|\mathcal{U}_{s}\|_{L_{\varphi}^{2}(\Omega)}^{2}
    +\frac{1}{16}\|\sqrt{\omega_{\lambda}}\mathcal{U}_{s}\|_{L_{\varphi}^{2}(\Omega)}^{2}
    +C_{*}\|\sqrt{\y}\mathcal{U}_{s}\|_{L_{\varphi}^{2}(\Omega)}^{2}
    +\frac{9}{16}\|\partial_{y}\mathcal{U}_{s}\|_{L_{\varphi}^{2}(\Omega)}^{2}
    +C\|\mathcal{P}\mathcal{U}_{s}\|_{L_{\varphi}^{2}(\Omega)}^{2}
\end{align*}

\noindent{\bf \underline{Estimate of~$I_{3}$.}}
Since $|\partial_{y}\rho|\leq C\mathcal{P}^{2}$, one has
\begin{align}\label{I31}
2\partial_{y}\rho-s\partial_{x}\mathfrak{U}\leq C\mathcal{P}^{2}+C_{*}.
\end{align}
By a direct calculation, one has
\begin{align*}
\frac{\partial_{y}\mathcal{P}}{\mathcal{P}}
= \frac{\chi'}{\mathcal{P}\sqrt{\dot{\mathfrak{U}}+\varsigma}}
 +\frac{\chi\rho}{2\mathcal{P}\sqrt{\dot{\mathfrak{U}}+\varsigma}}
 -\frac{\chi'}{\mathcal{P}},~\text{and,}~
\frac{\partial_{x}\mathcal{P}}{\mathcal{P}}
=\frac{\chi}{2\mathcal{P}\sqrt{\dot{\mathfrak{U}}+\varsigma}}\frac{\partial_{x}\dot{\mathfrak{U}}}{\dot{\mathfrak{U}}+\varsigma}.
\end{align*}
Thus, by noticing $\mathcal{P}\sqrt{\dot{\mathfrak{U}}+\varsigma}=1$ for $[\check{y}_{*},\hat{y}_{*}]$, one can get
\begin{align}\label{I32}
&2\left(\frac{\partial_{y}\mathcal{P}}{\mathcal{P}}\right)^{2}
+2\eta\left(\frac{\partial_{x}\mathcal{P}}{\mathcal{P}}\right)^{2}
-\frac{3\chi}{4\mathcal{P}\sqrt{\dot{\mathfrak{U}}+\varsigma}}\left[\rho^{2}+\eta\left(\frac{\partial_{x}\dot{\mathfrak{U}}}{\dot{\mathfrak{U}}+\varsigma}\right)^{2}\right]\\
\leq&-\frac{\chi}{4\mathcal{P}\sqrt{\dot{\mathfrak{U}}+\varsigma}}\left[\rho^{2}+\eta\left(\frac{\partial_{x}\dot{\mathfrak{U}}}{\dot{\mathfrak{U}}+\varsigma}\right)^{2}\right]
     +C_{*}
\leq C_{*}.\nonumber
\end{align}

Observe that ${\bf supp}\chi' \subset(\check{y},\hat{y})\setminus[\check{y}_{*},\hat{y}_{*}]$ and
${\bf supp}(1-\chi) \subset\mathbb{R}\setminus[\check{y}_{*},\hat{y}_{*}]$,
by noticing the definition of $\eta,\kappa,u_{0}$  and using the conditions $|\partial_{y}F|\leq C_{*}\iota^{-1}(\dot{\mathfrak{U}}+\bar{\varsigma})$, $|\rho|\leq C\mathcal{P}$ and $|\vartheta_{c}^{2}\partial_{x}^{2}\mathbb{U}|\leq C_{*}\dot{\mathfrak{U}}~\text{in}~\check{\Omega}_{T_{*}}^{\hat{y}}$, we can deduce that
\begin{align}\label{I33}
&\frac{\partial_{y}^{-1}\left[\partial_{x}\dot{\mathfrak{U}}\right]\varsigma'+\varsigma''-\partial_{y}\eta\partial_{x}^{2}\mathbb{U}+\eta\partial_{x}^{2}\partial_{y}u_{0}+\partial_{y}F+\partial_{y}[\kappa\partial_{x}^{2}\mathbb{U}]}{2(\dot{\mathfrak{U}}+\varsigma)}\\
&-(1-\chi)\frac{\partial_{y}^{-1}\left[\partial_{x}\mathfrak{U}\right]\varsigma'+\varsigma''-\partial_{y}\eta\partial_{x}^{2}\mathbb{U}+\eta\partial_{x}^{2}\partial_{y}u_{0}+\partial_{y}F+\partial_{y}[\kappa\partial_{x}^{2}\mathbb{U}]}{2\mathcal{P}(\dot{\mathfrak{U}}+\varsigma)}\nonumber\\
&\frac{-\partial_{y}^{-1}[\partial_{x}\mathfrak{U}]\chi'+\chi'\rho-\chi''}{\mathcal{P}\sqrt{\dot{\mathfrak{U}}+\varsigma}}
+\frac{\partial_{y}^{-1}\left[\partial_{x}\mathfrak{U}\right]\chi'+\chi''}{\mathcal{P}}
\leq C_{*}\iota^{-1}.\nonumber
\end{align}
A combination of \eqref{I31}--\eqref{I33} leads to
\begin{align*}
-\mathcal{L}_{0}\leq C\mathcal{P}^{2}+C_{*}\iota^{-1},
\end{align*}
which gives
\begin{align*}
I_{3}
\leq C_{*}\iota^{-1}\|\mathcal{U}_{s}\|_{L_{\varphi}^{2}(\Omega)}^{2}
    +C\|\mathcal{P}\mathcal{U}_{s}\|_{L_{\varphi}^{2}(\Omega)}^{2}.
\end{align*}

\noindent{\bf \underline{Estimate of~$I_{4}$.}}
Applying the trace theorem, one has
\begin{align*}
I_{4}
=(\mathcal{U}_{s},\varphi\partial_{y}\varphi\mathcal{U}_{s})_{L_{x}^{2}}|_{y=0}
\leq C\|\mathcal{U}_{s}\|_{L^{2}(\Omega)}^{2}
     +\frac{1}{16}\|\partial_{y}\mathcal{U}_{s}\|_{L_{\varphi}^{2}(\Omega)}^{2}.
\end{align*}

\noindent{\bf \underline{Estimate of~$I_{5}$.}}
Recall the boundary condition in \eqref{us}, by using trace estimate one has
\begin{align*}
I_{5}
=&-(\varphi^{2}\partial_{y}\mathcal{U}_{s},\mathcal{U}_{s})_{L_{x}^{2}}|_{y=0}\\
=& \varphi_{0}^{2}(\rho \mathcal{U}_{s},\mathcal{U}_{s})_{L_{x}^{2}}|_{y=0}
  +\left.\varphi_{0}^{2}(\partial_{x}^{s}F,\mathcal{U}_{s})_{L_{x}^{2}}\right|_{y=0}\\
\leq& C_{*}\|\mathcal{U}_{s}\|_{L_{\varphi}^{2}(\Omega)}^{2}
     +\frac{1}{16}\|\partial_{y}\mathcal{U}_{s}\|_{L_{\varphi}^{2}(\Omega)}^{2}
     +\|\partial_{x}^{s}F\|_{\mathring{H}^{0}(\mathbb{T})}^{2},
\end{align*}
where $\varphi_{0}=\varphi(t,0)$.

\noindent{\bf \underline{Estimate of~$I_{6}$.}}
Notice that
\begin{align*}
\bar{\mathcal{R}}_{1}(t,x,y)
=& (s+1)\partial_{x}^{s}\mathbb{U}\partial_{x}u
  +u\partial_{x}^{s+1}\mathbb{U}
  -\partial_{y}^{-1}[\partial_{x}u]\partial_{x}^{s}\partial_{y}\mathbb{U}
  -\partial_{y}^{-1}[\partial_{x}^{s+1}\mathbb{U}]\partial_{y}u
  -s\partial_{y}^{-1}[\partial_{x}^{s}\mathbb{U}]\partial_{x}\partial_{y}u\\
 &+\bar{\mathcal{R}}_{1R1}
  +\bar{\mathcal{R}}_{1R2},
\end{align*}
where
\begin{align*}
\bar{\mathcal{R}}_{1R1}(t,x,y)
=& \sum_{i=1}^{s-2}\binom{s}{i}\partial_{x}^{s-i}\mathfrak{U}\partial_{x}^{i+1}u
  +\sum_{i=2}^{s-1}\binom{s}{i}\partial_{x}^{s-i+1}\mathfrak{U}\partial_{x}^{i}u\\
 &+\sum_{i=1}^{s-2}\binom{s}{i}\partial_{y}^{-1}[\partial_{x}^{i+1}u]\partial_{x}^{s-i}\partial_{y}\mathfrak{U}
  +\sum_{i=2}^{s-1}\binom{s}{i}\partial_{y}^{-1}[\partial_{x}^{s-i+1}\mathbb{U}]\partial_{x}^{i}\partial_{y}u\\
 &+(s+1)\partial_{x}^{s}u_{0}\partial_{x}u
  +u\partial_{x}^{s+1}u_{0}
  -\partial_{y}^{-1}[\partial_{x}u]\partial_{x}^{s}\partial_{y}u_{0}
\end{align*}
and
\begin{align*}
\bar{\mathcal{R}}_{1R2}(t,x,y)
=&-\partial_{y}^{-1}[\partial_{x}^{s+1}u_{0}]\partial_{y}u
  -s\partial_{y}^{-1}[\partial_{x}^{s}u_{0}]\partial_{x}\partial_{y}u
  +\sum_{i=2}^{s-1}\binom{s}{i}\partial_{y}^{-1}[\partial_{x}^{s-i+1}u_{0}]\partial_{x}^{i}\partial_{y}u.
\end{align*}

A direct calculation gives
\begin{align*}
\mathcal{P}(\rho\bar{\mathcal{R}}_{1}-\partial_{y}\bar{\mathcal{R}}_{1})
=&-s\partial_{x}u\mathbf{U}_{s}
  -u\mathbf{U}_{s+1}
  +\mathcal{P}\partial_{x}^{s}\mathbb{U}(\rho\partial_{x}u-\partial_{y}\partial_{x}u)\\
 &+\mathcal{P}\partial_{y}^{-1}[\partial_{x}u](\partial_{y}^{2}\partial_{x}^{s}\mathbb{U}-\rho\partial_{y}\partial_{x}^{s}\mathbb{U})
  +\mathcal{P}\partial_{y}^{-1}[\partial_{x}^{s+1}\mathbb{U}](\partial_{y}^{2}u-\rho\partial_{y}u)\\
 &+s\mathcal{P}\partial_{y}^{-1}[\partial_{x}^{s}\mathbb{U}](\partial_{y}^{2}\partial_{x}u-\rho\partial_{y}\partial_{x}u)
  +\mathcal{P}(\rho\bar{\mathcal{R}}_{1R1}-\partial_{y}\bar{\mathcal{R}}_{1R1})\\
 &+\mathcal{P}(\rho\bar{\mathcal{R}}_{1R2}-\partial_{y}\bar{\mathcal{R}}_{1R2}),
\end{align*}
accordingly,
\begin{align}\label{I6}
I_{6}
=(\mathcal{P}(\rho\bar{\mathcal{R}}_{1}-\partial_{y}\bar{\mathcal{R}}_{1}),\varphi^{2}\mathcal{U}_{s})
=&-s(\partial_{x}u\mathbf{U}_{s},\varphi^{2}\mathcal{U}_{s})
  -(u\mathbf{U}_{s+1},\varphi^{2}\mathcal{U}_{s})\\
 &-(\mathcal{P}\partial_{x}^{s}\mathbb{U}(\rho\partial_{x}u-\partial_{y}\partial_{x}u),\varphi^{2}\mathcal{U}_{s})\nonumber\\
 &+(\mathcal{P}\partial_{y}^{-1}[\partial_{x}u](\partial_{y}^{2}\partial_{x}^{s}\mathbb{U}-\rho\partial_{y}\partial_{x}^{s}\mathbb{U}),\varphi^{2}\mathcal{U}_{s})\nonumber\\
 &+(\mathcal{P}\partial_{y}^{-1}[\partial_{x}^{s+1}\mathbb{U}](\partial_{y}^{2}u-\rho\partial_{y}u),\varphi^{2}\mathcal{U}_{s})\nonumber\\
 &+s(\mathcal{P}\partial_{y}^{-1}[\partial_{x}^{s}\mathbb{U}](\partial_{y}^{2}\partial_{x}u-\rho\partial_{y}\partial_{x}u),\varphi^{2}\mathcal{U}_{s})\nonumber\\
 &+(\mathcal{P}(\rho\bar{\mathcal{R}}_{1R1}-\partial_{y}\bar{\mathcal{R}}_{1R1}),\varphi^{2}\mathcal{U}_{s})\nonumber\\
 &+(\mathcal{P}(\rho\bar{\mathcal{R}}_{1R2}-\partial_{y}\bar{\mathcal{R}}_{1R2}),\varphi^{2}\mathcal{U}_{s}).\nonumber
\end{align}
Utilizing Corollary \ref{ncC2}, one can obtain
\begin{align}\label{I61}
-s(\partial_{x}u\mathbf{U}_{s},\varphi^{2}\mathcal{U}_{s})
\leq& C\|\partial_{x}u\|_{L^{\infty}(\Omega)}
      \|\mathbf{U}_{s}\|_{L_{\varphi}^{2}(\Omega)}
      \|\mathcal{U}_{s}\|_{L_{\varphi}^{2}(\Omega)}\\
\leq& C\|\mathbf{U}_{s}\|_{L_{\varphi}^{2}(\Omega)}
       \|u\|_{\mathcal{H}_{\psi,\varphi}^{s}(\Omega)}^{2}.\nonumber
\end{align}
and
\begin{align}\label{I62}
-(u\mathbf{U}_{s+1},\varphi^{2}\mathcal{U}_{s})
\leq& C\|\sqrt{\y}u\|_{L^{\infty}(\Omega)}
      \|\mathbf{U}_{s+1}\|_{L_{\hat{\varphi}}^{2}(\Omega)}
      \|\mathcal{U}_{s}\|_{L_{\varphi}^{2}(\Omega)}\\
\leq& C\|\mathbf{U}_{s+1}\|_{L_{\hat{\varphi}}^{2}(\Omega)}^{2}
       \|u\|_{\mathcal{H}_{\psi,\varphi}^{s}(\Omega)}^{2}.\nonumber
\end{align}
With the help of Corollary \ref{ncC2} and Sobolev imbedding inequality, one has
\begin{align}\label{I63}
&-(\mathcal{P}\partial_{x}^{s}\mathbb{U}(\rho\partial_{x}u-\partial_{y}\partial_{x}u),\varphi^{2}\mathcal{U}_{s})\\
\leq&~C_{*}\left(\|\partial_{x}u\|_{L^{\infty}(\Omega)}+\|\partial_{x}\partial_{y}u\|_{L^{\infty}(\Omega)}\right)
      \|\mathcal{P}^{2}\partial_{x}^{s}\mathbb{U}\|_{L_{\varphi}^{2}(\Omega)}
      \|\mathcal{U}_{s}\|_{L_{\varphi}^{2}(\Omega)}\nonumber\\
\leq&~C_{*}\left(\|\partial_{x}u\|_{L^{\infty}(\Omega)}+\|\partial_{x}\partial_{y}u\|_{L^{\infty}(\Omega)}\right)
       \|\mathcal{P}\mathbf{U}_{s}\|_{L_{\varphi}^{2}(\Omega)}
       \|\mathcal{U}_{s}\|_{L_{\varphi}^{2}(\Omega)}\nonumber\\
\leq&~C_{*}\|\mathcal{P}\mathbf{U}_{s}\|_{L_{\varphi}^{2}(\Omega)}\|u\|_{\mathcal{H}_{\psi,\varphi}^{s}(\Omega)}^{2}.\nonumber
\end{align}
Again by using Corollary \ref{ncC2}, one can obtain,
\begin{align}\label{I64'}
&(\mathcal{P}\partial_{y}^{-1}[\partial_{x}u](\partial_{y}^{2}\partial_{x}^{s}\mathbb{U}-\rho\partial_{y}\partial_{x}^{s}\mathbb{U}),\varphi^{2}\mathcal{U}_{s})\\
\leq&~C\left(\|\partial_{y}^{2}\partial_{x}^{s}\mathbb{U}\|_{L_{\varphi}^{2}(\Omega)}+\|\mathcal{P}\partial_{y}\partial_{x}^{s}\mathbb{U}\|_{L_{\varphi}^{2}(\Omega)}\right)
       \|u\|_{\hat{H}_{\psi}^{s}(\Omega)}\|\mathcal{P}\mathcal{U}_{s}\|_{L_{\varphi}^{2}(\Omega)}\nonumber\\
\leq&~C_{*}(\|\partial_{y}\mathbf{U}_{s}\|_{L_{\varphi}^{2}(\Omega)}^{2}
        +\|\mathcal{P}\mathbf{U}_{s}\|_{L_{\varphi}^{2}(\Omega)}^{2}
        +\|\mathbf{U}_{s}\|_{L_{\varphi}^{2}(\Omega)}^{2})
       \|u\|_{\mathcal{H}_{\psi,\varphi}^{s}(\Omega)}^{2}\nonumber\\
    &+\frac{1}{4}\|\mathcal{P}\mathcal{U}_{s}\|_{L_{\varphi}^{2}(\Omega)}^{2}.\nonumber
\end{align}
It follows from Corollary \ref{ncC2} and Corollary \ref{Ltc} that
\begin{align}\label{I65}
&(\mathcal{P}\partial_{y}^{-1}[\partial_{x}^{s+1}\mathbb{U}](\partial_{y}^{2}u-\rho\partial_{y}u),\varphi^{2}\mathcal{U}_{s})\\
\leq& C_{*}(\|\varphi\mathcal{P}^{2}\partial_{y}u\|_{L_{x}^{\infty}L_{y}^{2}(\Omega)}
      +\|\varphi\mathcal{P}\partial_{y}^{2}u\|_{L_{x}^{\infty}L_{y}^{2}(\Omega)})
      \|\partial_{x}^{s+1}\mathbb{U}\|_{L_{\hat{\varphi}}^{2}(\Omega)}
      \|\mathcal{U}_{s}\|_{L_{\varphi}^{2}(\Omega)}\nonumber\\
\leq& C_{*}\|\mathcal{P}\mathbf{U}_{s+1}\|_{L_{\hat{\varphi}}^{2}(\Omega)}^{2}
      \|\mathcal{U}_{s}\|_{L_{\varphi}^{2}(\Omega)}^{2}
     +\|\varphi\mathcal{P}^{2}\partial_{y}u\|_{L_{x}^{\infty}L_{y}^{2}(\Omega)}^{2}
     +\|\varphi\mathcal{P}\partial_{y}^{2}u\|_{L_{x}^{\infty}L_{y}^{2}(\Omega)}^{2}\nonumber\\
\leq& C_{*}\|\mathcal{P}\mathbf{U}_{s+1}\|_{L_{\hat{\varphi}}^{2}(\Omega)}^{2}\|\mathcal{U}_{s}\|_{L_{\varphi}^{2}(\Omega)}^{2}
     +C_{*}t\|u\|_{L_{t}^{2}\mathcal{H}_{\psi,\varphi}^{s}(\Omega)}^{2}
     +C_{*}\|u\|_{\hat{H}_{\psi}^{s}(\Omega)}^{2}
     +C_{*}\bar{\eta}^{2}
     +C_{*}\check{c}^{2}\iota^{-2}.\nonumber
\end{align}
and similarly,
\begin{align}\label{I66}
&s(\mathcal{P}\partial_{y}^{-1}[\partial_{x}^{s}\mathbb{U}](\partial_{y}^{2}\partial_{x}u-\rho\partial_{y}\partial_{x}u),\varphi^{2}\mathcal{U}_{s})\\
\leq& C_{*}(\|\varphi\mathcal{P}^{2}\partial_{y}\partial_{x}u\|_{L_{x}^{\infty}L_{y}^{2}(\Omega)}
       +\|\varphi\mathcal{P}\partial_{y}^{2}\partial_{x}u\|_{L_{x}^{\infty}L_{y}^{2}(\Omega)})
      \|\partial_{x}^{s}\mathbb{U}\|_{L_{\varphi}^{2}(\Omega)}
      \|\mathcal{U}_{s}\|_{L_{\varphi}^{2}(\Omega)}\nonumber\\
\leq& C_{*}\|\mathcal{P}\mathbf{U}_{s}\|_{L_{\varphi}^{2}(\Omega)}^{2}
      \|\mathcal{U}_{s}\|_{L_{\varphi}^{2}(\Omega)}^{2}
     +\|\varphi\mathcal{P}^{2}\partial_{y}\partial_{x}u\|_{L_{x}^{\infty}L_{y}^{2}(\Omega)}^{2}
     +\|\varphi\mathcal{P}\partial_{y}^{2}\partial_{x}u\|_{L_{x}^{\infty}L_{y}^{2}(\Omega)}^{2}\nonumber\\
\leq& C_{*}\|\mathcal{P}\mathbf{U}_{s}\|_{L_{\varphi}^{2}(\Omega)}^{2}\|\mathcal{U}_{s}\|_{L_{\varphi}^{2}(\Omega)}^{2}
     +C_{*}t\|u\|_{L_{t}^{2}\mathcal{H}_{\psi,\varphi}^{s}(\Omega)}^{2}
     +C_{*}\|u\|_{\hat{H}_{\psi}^{s}(\Omega)}^{2}
     +C_{*}\bar{\eta}^{2}
     +C_{*}\check{c}^{2}\iota^{-2}.\nonumber
\end{align}
Note that $\dot{\mathfrak{U}}\geq C(\varpi+t)$ on $[\check{y},\hat{y}]$, we have
\begin{align}\label{I671}
(\mathcal{P}\rho\bar{\mathcal{R}}_{1R1}\varphi^{2},\mathcal{U}_{s})
=& \int_{\mathbb{T}}\int_{\check{y}}^{\hat{y}}\varphi^{2}\mathcal{P}\rho\bar{\mathcal{R}}_{1R1}\mathcal{U}_{s}dydx
  +\int_{\mathbb{T}}\int_{\mathbb{R}_{+}\setminus[\check{y},\hat{y}]}\varphi^{2}\mathcal{P}\rho\bar{\mathcal{R}}_{1R1}\mathcal{U}_{s}dydx\\
\leq& C_{*}\|\bar{\mathcal{R}}_{1R1}\|_{L^{2}_{x}(\mathbb{T})L^{\infty}_{y}[\check{y},\hat{y}]}\|\rho\mathcal{U}_{s}\|_{L_{\varphi}^{2}(\Omega)}
       \left(\int_{\check{y}_{*}}^{\hat{y}_{*}}\frac{1}{((y-y_{*})^2+\bar{\varsigma}+t)^{1+\varepsilon_{0}}}dy\right)^{\frac{1}{2}}\nonumber\\
    &+C_{*}\|\bar{\mathcal{R}}_{1R1}\|_{L_{\psi}^{2}(\Omega)}\|\mathcal{U}_{s}\|_{L_{\varphi}^{2}(\Omega)}\nonumber\\
\leq& \frac{1}{4}\|\mathcal{P}\mathcal{U}_{s}\|_{L_{\varphi}^{2}(\Omega)}^{2}
     +\frac{C_{*}}{(\bar{\varsigma}+t)^{\frac{1}{2}+\varepsilon_{0}}}\|\bar{\mathcal{R}}_{1R1}\|_{L_{x}^{2}(\mathbb{T})L_{y}^{\infty}[\check{y}_{*},\hat{y}_{*}]}^{2}
     +C_{*}\|\bar{\mathcal{R}}_{1R1}\|_{L_{\varphi}^{2}(\Omega)}\|\mathcal{U}_{s}\|_{L_{\varphi}^{2}(\Omega)}\nonumber\\
\leq& \frac{1}{4}\|\mathcal{P}\mathcal{U}_{s}\|_{L_{\varphi}^{2}(\Omega)}^{2}
     +\frac{C_{*}}{(\bar{\varsigma}+t)^{\frac{1}{2}+\varepsilon_{0}}}\left(\|\mathbb{U}\|_{\hat{H}_{\psi}^{s}(\Omega)}^{2}+1\right)\|u\|_{\mathcal{H}_{\psi,\varphi}^{s}(\Omega)}^{2}\nonumber\\
    &+C_{*}\left(\|\mathbb{U}\|_{\hat{H}_{\psi}^{s}(\Omega)}^{2}+1\right)\|u\|_{\mathcal{H}_{\psi,\varphi}^{s}(\Omega)}^{2}\nonumber\\
\leq& \frac{1}{4}\|\mathcal{P}\mathcal{U}_{s}\|_{L_{\varphi}^{2}(\Omega)}^{2}
     +\frac{C_{*}}{(\bar{\varsigma}+t)^{\frac{1}{2}+\varepsilon_{0}}}\|u\|_{\mathcal{H}_{\psi,\varphi}^{s}(\Omega)}^{2}
     +C_{*}\|u\|_{\mathcal{H}_{\psi,\varphi}^{s}(\Omega)}^{2}.\nonumber
\end{align}
and
\begin{align}\label{I672}
-(\mathcal{P}\partial_{y}\bar{\mathcal{R}}_{1R1},\varphi^{2}\mathcal{U}_{s})
\leq&~\frac{C_{*}}{(\bar{\varsigma}+t)^{\varepsilon_{0}}}\|\partial_{y}\bar{\mathcal{R}}_{1R1}\|_{L_{\hat{\psi}}^{2}(\Omega)}^{2}
     +\frac{1}{4}\left\|\mathcal{P}\mathcal{U}_{s}\right\|_{L_{\varphi}^{2}(\Omega)}^{2}\\
\leq&~\frac{C_{*}}{(\bar{\varsigma}+t)^{\varepsilon_{0}}}\left(\|\mathbb{U}\|_{\hat{H}_{\psi}^{s}(\Omega)}^{2}+1\right)\|u\|_{\mathcal{H}_{\psi,\varphi}^{s}(\Omega)}^{2}
     +\frac{1}{4}\|\mathcal{P}\mathcal{U}_{s}\|_{L_{\varphi}^{2}(\Omega)}^{2}\nonumber\\
\leq&~\frac{C_{*}}{(\bar{\varsigma}+t)^{\varepsilon_{0}}}\|u\|_{\mathcal{H}_{\psi,\varphi}^{s}(\Omega)}^{2}
     +\frac{1}{4}\|\mathcal{P}\mathcal{U}_{s}\|_{L_{\varphi}^{2}(\Omega)}^{2}.\nonumber
\end{align}
Similarly, one has
\begin{align}\label{I681}
(\mathcal{P}\rho\bar{\mathcal{R}}_{1R2}\varphi^{2},\mathcal{U}_{s})
=&~\int_{\mathbb{T}}\int_{\check{y}}^{\hat{y}}\varphi^{2}\mathcal{P}\rho\bar{\mathcal{R}}_{1R2}\mathcal{U}_{s}dydx
  +\int_{\mathbb{T}}\int_{\mathbb{R}_{+}\setminus[\check{y},\hat{y}]}\varphi^{2}\mathcal{P}\rho\bar{\mathcal{R}}_{1R2}\mathcal{U}_{s}dydx\\
\leq&~C_{*}\|\bar{\mathcal{R}}_{1R2}\|_{L^{2}_{x}(\mathbb{T})L^{\infty}_{y}[\check{y},\hat{y}]}\|\rho\mathcal{U}_{s}\|_{L_{\varphi}^{2}(\Omega)}
       \left(\int_{\check{y}_{*}}^{\hat{y}_{*}}\frac{1}{((y-y_{*})^2+\bar{\varsigma}+t)^{1+\varepsilon_{0}}}dy\right)^{\frac{1}{2}}\nonumber\\
    &+C_{*}\|\bar{\mathcal{R}}_{1R2}\|_{L_{\hat{\psi}}^{2}(\Omega)}\|\mathcal{U}_{s}\|_{L_{\varphi}^{2}(\Omega)}\nonumber\\
\leq&~\frac{1}{4}\|\mathcal{P}\mathcal{U}_{s}\|_{L_{\varphi}^{2}(\Omega)}^{2}
     +\frac{C_{*}}{(\bar{\varsigma}+t)^{\frac{1}{2}+\varepsilon_{0}}}\|\bar{\mathcal{R}}_{1R2}\|_{L_{x}^{2}(\mathbb{T})L_{y}^{\infty}[\check{y}_{*},\hat{y}_{*}]}^{2}
     +C_{*}\|u\|_{\hat{H}_{\psi,\y}^{s}(\Omega)}\|\mathcal{U}_{s}\|_{L_{\varphi}^{2}(\Omega)}\nonumber\\
\leq&~\frac{1}{4}\|\mathcal{P}\mathcal{U}_{s}\|_{L_{\varphi}^{2}(\Omega)}^{2}
     +\frac{C_{*}}{(\bar{\varsigma}+t)^{\frac{1}{2}+\varepsilon_{0}}}\|u\|_{\mathcal{H}_{\psi,\varphi}^{s}(\Omega)}^{2}
     +C_{*}\|u\|_{\hat{H}_{\psi,\y}^{s}(\Omega)}^{2}
     +C_{*}\|\mathcal{U}_{s}\|_{L_{\varphi}^{2}(\Omega)}^{2}\nonumber,
\end{align}
and
\begin{align}\label{I682}
(\mathcal{P}\partial_{y}\bar{\mathcal{R}}_{1R2}\varphi^{2},\mathcal{U}_{s})
=& \int_{\mathbb{T}}\int_{\check{y}}^{\hat{y}}\varphi^{2}\mathcal{P}\partial_{y}\bar{\mathcal{R}}_{1R2}\mathcal{U}_{s}dydx
  +\int_{\mathbb{T}}\int_{\mathbb{R}_{+}\setminus[\check{y},\hat{y}]}\varphi^{2}\mathcal{P}\partial_{y}\bar{\mathcal{R}}_{1R2}\mathcal{U}_{s}dydx\\
\leq& C_{*}\|\partial_{y}\bar{\mathcal{R}}_{1R2}\|_{L^{2}_{x}(\mathbb{T})L^{\infty}_{y}[\check{y},\hat{y}]}
           \|\mathcal{P}\mathcal{U}_{s}\|_{L_{\varphi}^{2}(\Omega)}\nonumber\\
    &+C_{*}\|\partial_{y}\bar{\mathcal{R}}_{1R2}\|_{L_{\hat{\psi}}^{2}(\Omega)}\|\mathcal{U}_{s}\|_{L_{\varphi}^{2}(\Omega)}\nonumber\\
\leq& \frac{1}{4}\|\mathcal{P}\mathcal{U}_{s}\|_{L_{\varphi}^{2}(\Omega)}^{2}
     +C_{*}\|u\|_{\hat{H}_{\psi}^{s}(\Omega)}^{2}
     +C_{*}\|u\|_{\hat{H}_{\psi,\y}^{s}(\Omega)}\|\mathcal{U}_{s}\|_{L_{\varphi}^{2}(\Omega)}\nonumber\\
\leq& \frac{1}{4}\|\mathcal{P}\mathcal{U}_{s}\|_{L_{\varphi}^{2}(\Omega)}^{2}
     +C_{*}\|u\|_{\mathcal{H}_{\psi,\varphi}^{s}(\Omega)}^{2}
     +C_{*}\|u\|_{\hat{H}_{\psi,\y}^{s}(\Omega)}^{2}
     +C_{*}\|\mathcal{U}_{s}\|_{L_{\varphi}^{2}(\Omega)}^{2}\nonumber,
\end{align}
since $\partial_{x}u_{0}=0$ for $y\in[0,\hat{y}]$.

Substituting \eqref{I61}-\eqref{I63},\eqref{I64'} and \eqref{I65}-\eqref{I682} into \eqref{I6}, it follows
\begin{align*}
I_{6}
\leq& C_{*}\left(1
        +\|\mathbf{U}_{s}\|_{L_{\varphi}^{2}(\Omega)}^{2}
        +\|\mathcal{P}\mathbf{U}_{s}\|_{L_{\varphi}^{2}(\Omega)}^{2}
        +\|\partial_{y}\mathbf{U}_{s}\|_{L_{\varphi}^{2}(\Omega)}^{2}\right)
       \|u\|_{\mathcal{H}_{\psi,\varphi}^{s}(\Omega)}^{2}\\
    &+C_{*}\left(\|\mathbf{U}_{s+1}\|_{L_{\hat{\varphi}}^{2}(\Omega)}^{2}
        +\|\mathcal{P}\mathbf{U}_{s+1}\|_{L_{\hat{\varphi}}^{2}(\Omega)}^{2}\right)\|u\|_{\mathcal{H}_{\psi,\varphi}^{s}(\Omega)}^{2}
        +\frac{C_{*}}{(\bar{\varsigma}+t)^{\frac{1}{2}+\varepsilon_{0}}}\|u\|_{\mathcal{H}_{\psi,\varphi}^{s}(\Omega)}^{2}
     +\|\mathcal{P}\mathcal{U}_{s}\|_{L_{\varphi}^{2}(\Omega)}^{2}\\
    &+C_{*}\|u\|_{\hat{H}_{\psi,\y}^{s}(\Omega)}^{2}
     +C_{*}t\|u\|_{L_{t}^{2}\mathcal{H}_{\psi,\varphi}^{s}(\Omega)}^{2}
     +C_{*}\bar{\eta}^{2}
     +C_{*}\check{c}^{2}\iota^{-2}.
\end{align*}

\noindent{\bf \underline{Estimate of~$I_{7}$.}}
Recall the definition of $\mathcal{R}_{2}$, one has
\begin{align}\label{I7}
I_{7}
=&\left((\partial_{y}^{2}F-\rho\partial_{y}F)\mathcal{P}\frac{\partial_{x}^{s}u}{\dot{\mathfrak{U}}+\varsigma},\varphi^{2}\mathcal{U}_{s}\right)
  -2\left(\eta\mathcal{P}\partial_{x}\dot{\mathfrak{U}}\partial_{x}\rho\frac{\partial_{x}^{s}u}{\dot{\mathfrak{U}}+\varsigma},\varphi^{2}\mathcal{U}_{s}\right)\\
 &+\left(\mathcal{P}\left(\varsigma'\partial_{x}\mathfrak{U}+\varsigma'''+\partial_{y}^{-1}[\partial_{x}\mathfrak{U}]\varsigma''-\varsigma\partial_{x}\dot{\mathfrak{U}}
  -\rho\left(\partial_{y}^{-1}[\partial_{x}\mathfrak{U}]\varsigma'+\varsigma''\right)\right)\frac{\partial_{x}^{s}u}{\dot{\mathfrak{U}}+\varsigma},\varphi^{2}\mathcal{U}_{s}\right)\nonumber\\
 &+\left((\partial_{y}^{2}[\kappa\partial_{x}^{2}\mathbb{U}]-\rho\partial_{y}[\kappa\partial_{x}^{2}\mathbb{U}])\mathcal{P}\frac{\partial_{x}^{s}u}{\dot{\mathfrak{U}}+\varsigma},\varphi^{2}\mathcal{U}_{s}\right).\nonumber
\end{align}
In view of $\|\mathcal{P}(\partial_{y}^{2}F-\rho\partial_{y}F)\|_{L^{\infty}(\Omega_{T_{*}})}\leq C_{*}\iota^{-1}$, and for any $\gamma\in\Gamma_{2}, |D^{\gamma}F|\leq C_{*}e^{-\frac{4}{3}\delta\y}~\text{in}~\Omega_{T_{*}},$
by using Corollary \ref{ncC2} we have
\begin{align}\label{I71}
\left((\partial_{y}^{2}F-\rho\partial_{y}F)\mathcal{P}\frac{\partial_{x}^{s}u}{\dot{\mathfrak{U}}+\varsigma},\varphi^{2}\mathcal{U}_{s}\right)
\leq&~C_{*}\iota^{-1}\|\partial_{x}^{s}u\|_{L_{\hat{\Psi}}^{2}(\Omega)}\|\mathcal{U}_{s}\|_{L_{\varphi}^{2}(\Omega)}\\
\leq&~C_{*}\iota^{-1}\|\mathcal{P}\mathcal{U}_{s}\|_{L_{\varphi}^{2}(\Omega)}\|\mathcal{U}_{s}\|_{L_{\varphi}^{2}(\Omega)}\nonumber\\
\leq&~\frac{1}{3}\|\mathcal{P}\mathcal{U}_{s}\|_{L_{\varphi}^{2}(\Omega)}^{2}
     +C_{*}\iota^{-2}\|\mathcal{U}_{s}\|_{L_{\varphi}^{2}(\Omega)}^{2}.\nonumber
\end{align}
Owing to $|\partial_{x}\rho|\leq C\mathcal{P}$ and $|\partial_{x}\dot{\mathfrak{U}}|\leq C_{*}(1+t\y)\dot{\mathfrak{U}}$ in $\Omega_{T_{*}}$, we get
\begin{align}\label{I72}
-2\left(\eta\mathcal{P}\partial_{x}\dot{\mathfrak{U}}\partial_{x}\rho\frac{\partial_{x}^{s}u}{\dot{\mathfrak{U}}+\varsigma},\varphi^{2}\mathcal{U}_{s}\right)
\leq& C_{*}C_{\theta}\bar{\eta}\|\partial_{x}^{s}u\|_{L_{\hat{\Psi}}^{2}(\Omega)}\|\mathcal{U}_{s}\|_{L_{\varphi}^{2}(\Omega)}\\
\leq& C_{*}C_{\theta}\bar{\eta}\|\mathcal{P}\mathcal{U}_{s}\|_{L_{\varphi}^{2}(\Omega)}\|\mathcal{U}_{s}\|_{L_{\varphi}^{2}(\Omega)}\nonumber\\
\leq& C_{*}C_{\theta}\iota^{-\frac{1}{2}}\|\mathcal{U}_{s}\|_{L_{\varphi}^{2}(\Omega)}^{2}.\nonumber
\end{align}
Due to
$$\varsigma'\partial_{x}\mathfrak{U}+\varsigma'''+\partial_{y}^{-1}[\partial_{x}\mathfrak{U}]\varsigma''-\varsigma\partial_{x}\dot{\mathfrak{U}}
  -\rho\left(\partial_{y}^{-1}[\partial_{x}\mathfrak{U}]\varsigma'+\varsigma''\right)\leq C_{*}\sqrt{\bar{\varsigma}},$$
and $\sqrt{\bar{\varsigma}}\mathcal{P}\leq C_{*}$, using Corollary \ref{ncC2} one has
\begin{align}\label{I73}
&\left(\mathcal{P}\left(\varsigma'\partial_{x}\mathfrak{U}+\varsigma'''+\partial_{y}^{-1}[\partial_{x}\mathfrak{U}]\varsigma''-\varsigma\partial_{x}\dot{\mathfrak{U}}
  -\rho\left(\partial_{y}^{-1}[\partial_{x}\mathfrak{U}]\varsigma'+\varsigma''\right)\right)\frac{\partial_{x}^{s}u}{\dot{\mathfrak{U}}+\varsigma},\varphi^{2}\mathcal{U}_{s}\right)\\
\leq&~C_{*}\|\partial_{x}^{s}u\|_{L_{\hat{\Psi}}^{2}(\Omega)}\|\mathcal{U}_{s}\|_{L_{\varphi}^{2}(\Omega)}\nonumber\\
\leq&~\frac{1}{3}\|\mathcal{P}\mathcal{U}_{s}\|_{L_{\varphi}^{2}(\Omega)}^{2}
     +C_{*}\|\mathcal{U}_{s}\|_{L_{\varphi}^{2}(\Omega)}^{2}.\nonumber
\end{align}
Since $|\partial_{x}^{2}\dot{\mathfrak{U}}|\leq C(1+t\y)\dot{\mathfrak{U}}, |\partial_{y}\partial_{x}^{2}\dot{\mathfrak{U}}|\leq C(1+t\y)\dot{\mathfrak{U}}^{\frac{1}{2}}$ in $\Omega_{T_{*}}$, $\partial_{x}u_{0}=0$ for any $y\in[0,\hat{y}]$,
 and
$|\vartheta_{c}^{2}\partial_{x}^{2}\mathbb{U}|\leq C_{*}\dot{\mathfrak{U}}~\text{in}~\check{\Omega}_{T_{*}}^{\hat{y}}$,
recall the definition of $\kappa$, also by using Corollary \ref{ncC2} one gets
\begin{align}\label{I74}
&\left((\partial_{y}^{2}[\kappa\partial_{x}^{2}\mathbb{U}]-\rho\partial_{y}[\kappa\partial_{x}^{2}\mathbb{U}])\mathcal{P}\frac{\partial_{x}^{s}u}{\dot{\mathfrak{U}}+\varsigma},\varphi^{2}\mathcal{U}_{s}\right)\\
\leq&~C_{*}C_{\theta}\|\partial_{x}^{s}u\|_{L_{\hat{\Psi}}^{2}(\Omega)}\|\mathcal{U}_{s}\|_{L_{\varphi}^{2}(\Omega)}\nonumber\\
\leq&~\frac{1}{3}\|\mathcal{P}\mathcal{U}_{s}\|_{L_{\varphi}^{2}(\Omega)}^{2}
     +C_{*}C_{\theta}\|\mathcal{U}_{s}\|_{L_{\varphi}^{2}(\Omega)}^{2}.\nonumber
\end{align}

In virtue of \eqref{I71}-\eqref{I74}, from \eqref{I7} one can obtain
\begin{align*}
I_{7}\leq \|\mathcal{P}\mathcal{U}_{s}\|_{L_{\varphi}^{2}(\Omega)}^{2}
         +C_{*}C_{\theta}\iota^{-2}\|\mathcal{U}_{s}\|_{L_{\varphi}^{2}(\Omega)}^{2}.
\end{align*}

\noindent{\bf \underline{Estimate of~$I_{8}$.}}
Recall that
\begin{align*}
I_{8}
=s\left(\mathcal{P}\left(\partial_{x}\partial_{y}^{2}\mathfrak{U}-\rho\partial_{x}\partial_{y}\mathfrak{U}\right)\partial_{y}^{-1}[\partial_{x}^{s}u],\varphi^{2}\mathcal{U}_{s}\right)
\end{align*}
Noticing $\partial_{x}\partial_{y}^{2}\mathfrak{U}-\rho\partial_{x}\partial_{y}\mathfrak{U}=\partial_{x}\rho(\dot{\mathfrak{U}}+\varsigma)$,
$|\partial_{x}\rho|\leq C\mathcal{P}$,
by using Corollary \ref{ncC2} we have
\begin{align*}
I_{8}
=&s\left((\dot{\mathfrak{U}}+\varsigma)\mathcal{P}\partial_{x}\rho\partial_{y}^{-1}[\partial_{x}^{s}u],\varphi^{2}\mathcal{U}_{s}\right)\\
\leq& C_{*}\|\varphi(\dot{\mathfrak{U}}+\varsigma)\mathcal{P}^{2}\|_{L_{x}^{\infty}(\mathbb{T})L_{y}^{2}(\mathbb{R_{+}})}
       \|\partial_{x}^{s}u\|_{L_{\varphi}^{2}(\Omega)}\|\mathcal{U}_{s}\|_{L_{\varphi}^{2}(\Omega)}\\
\leq& C_{*}\|\partial_{x}^{s}u\|_{L_{\hat{\Psi}}^{2}(\Omega)}\|\mathcal{U}_{s}\|_{L_{\varphi}^{2}(\Omega)}\\
\leq& C_{*}\|\mathcal{P}\mathcal{U}_{s}\|_{L_{\varphi}^{2}(\Omega)}\|\mathcal{U}_{s}\|_{L_{\varphi}^{2}(\Omega)}\\
\leq& \|\mathcal{P}\mathcal{U}_{s}\|_{L_{\varphi}^{2}(\Omega)}^{2}
     +C_{*}\|\mathcal{U}_{s}\|_{L_{\varphi}^{2}(\Omega)}^{2}.
\end{align*}

\noindent{\bf \underline{Estimate of~$I_{9}$.}}
By using integration by parts, one has
\begin{align}\label{I9}
I_{9}
=&((\varsigma+2\eta\partial_{x}\rho)\mathcal{P}\partial_{x}^{s+1}u,\varphi^{2}\mathcal{U}_{s})\\
=&-2(\eta\partial_{x}(\mathcal{P}\partial_{x}\rho)\partial_{x}^{s}u,\varphi^{2}\mathcal{U}_{s})
  -2(\eta\partial_{x}\rho\mathcal{P}\partial_{x}^{s}u,\varphi^{2}\partial_{x}\mathcal{U}_{s})\nonumber\\
 &-(\varsigma\mathcal{P}\partial_{x}^{s}u,\varphi^{2}\partial_{x}\mathcal{U}_{s})
  -(\varsigma\partial_{x}\mathcal{P}\partial_{x}^{s}u,\varphi^{2}\mathcal{U}_{s}).\nonumber
\end{align}
In view of $\left|\partial_{x}\rho\right|\leq C\mathcal{P}$
and $\left|\partial_{x}(\mathcal{P}\partial_{x}\rho)\right|\leq C(1+t\y)^{2}\mathcal{P}^{2}$,
by using Corollary \ref{ncC2} we have
\begin{align}\label{I91}
&-2(\eta\partial_{x}(\mathcal{P}\partial_{x}\rho)\partial_{x}^{s}u,\varphi^{2}\mathcal{U}_{s})\\
\leq& C_{*}\|\eta\partial_{x}(\mathcal{P}\partial_{x}\rho)\partial_{x}^{s}u\|_{L_{\varphi}^{2}(\Omega)}^{2}
     +C\|\mathcal{U}_{s}\|_{L_{\varphi}^{2}(\Omega)}^{2}\nonumber\\
\leq& C_{*}C_{\theta}\bar{\eta}^{2}\|\partial_{x}^{s}u\|_{L_{\hat{\Psi}}^{2}(\Omega)}^{2}
     +C\|\mathcal{U}_{s}\|_{L_{\varphi}^{2}(\Omega)}^{2}\nonumber\\
\leq& C_{*}C_{\theta}\bar{\eta}^{2}\|\mathcal{P}\mathcal{U}_{s}\|_{L_{\varphi}^{2}(\Omega)}^{2}
     +C\|\mathcal{U}_{s}\|_{L_{\varphi}^{2}(\Omega)}^{2}\nonumber\\
\leq& C_{*}C_{\theta}\iota^{-1}\|\mathcal{U}_{s}\|_{L_{\varphi}^{2}(\Omega)}^{2},\nonumber
\end{align}
and
\begin{align}\label{I92}
-2(\eta\partial_{x}\rho\mathcal{P}\partial_{x}^{s}u,\varphi^{2}\partial_{x}\mathcal{U}_{s})
\leq& C\|\sqrt{\eta}\partial_{x}\rho\mathcal{P}\partial_{x}^{s}u\|_{L_{\varphi}^{2}(\Omega)}^{2}
    +\frac{1}{16}\|\sqrt{\eta}\partial_{x}\mathcal{U}_{s}\|_{L_{\varphi}^{2}(\Omega)}^{2}\\
\leq& C_{*}\bar{\eta}\|\partial_{x}^{s}u\|_{L_{\hat{\Psi}}^{2}(\Omega)}^{2}
     +\frac{1}{16}\|\sqrt{\eta}\partial_{x}\mathcal{U}_{s}\|_{L_{\varphi}^{2}(\Omega)}^{2}\nonumber\\
\leq& C_{*}\bar{\eta}\|\mathcal{P}\mathcal{U}_{s}\|_{L_{\varphi}^{2}(\Omega)}^{2}
     +\frac{1}{16}\|\sqrt{\eta}\partial_{x}\mathcal{U}_{s}\|_{L_{\varphi}^{2}(\Omega)}^{2}\nonumber\\
\leq& C_{*}\iota^{-\frac{1}{2}}\|\mathcal{U}_{s}\|_{L_{\varphi}^{2}(\Omega)}\|\mathcal{P}\mathcal{U}_{s}\|_{L_{\varphi}^{2}(\Omega)}
     +\frac{1}{16}\|\sqrt{\eta}\partial_{x}\mathcal{U}_{s}\|_{L_{\varphi}^{2}(\Omega)}^{2}\nonumber\\
\leq& C_{*}\iota^{-1}\|\mathcal{U}_{s}\|_{L_{\varphi}^{2}(\Omega)}^{2}
     +\|\mathcal{P}\mathcal{U}_{s}\|_{L_{\varphi}^{2}(\Omega)}^{2}
     +\frac{1}{16}\|\sqrt{\eta}\partial_{x}\mathcal{U}_{s}\|_{L_{\varphi}^{2}(\Omega)}^{2}.\nonumber
\end{align}
Similarly, noticing that ${\bf supp}~\varsigma\subset(\check{y},\hat{y})$ one has
\begin{align}\label{I93}
-(\varsigma\mathcal{P}\partial_{x}^{s}u,\varphi^{2}\partial_{x}\mathcal{U}_{s})
\leq& C\bar{\eta}^{-1}\bar{\varsigma}^{2}\|\mathcal{P}\partial_{x}^{s}u\|_{L_{\varphi}^{2}(\Omega)}^{2}
     +\frac{1}{16}\|\sqrt{\eta}\partial_{x}\mathcal{U}_{s}\|_{L_{\varphi}^{2}(\Omega)}^{2}\\
\leq& C\bar{\eta}^{-1}\bar{\varsigma}^{2}\|\partial_{x}^{s}u\|_{L_{\hat{\Psi}}^{2}(\Omega)}^{2}
     +\frac{1}{16}\|\sqrt{\eta}\partial_{x}\mathcal{U}_{s}\|_{L_{\varphi}^{2}(\Omega)}^{2}\nonumber\\
\leq& C_{*}\bar{\eta}^{-1}\bar{\varsigma}^{2}\|\mathcal{P}\mathcal{U}_{s}\|_{L_{\varphi}^{2}(\Omega)}^{2}
     +\frac{1}{16}\|\sqrt{\eta}\partial_{x}\mathcal{U}_{s}\|_{L_{\varphi}^{2}(\Omega)}^{2}\nonumber\\
\leq& C_{*}\|\mathcal{U}_{s}\|_{L_{\varphi}^{2}(\Omega)}^{2}
     +\frac{1}{16}\|\sqrt{\eta}\partial_{x}\mathcal{U}_{s}\|_{L_{\varphi}^{2}(\Omega)}^{2}.\nonumber
\end{align}
Note that $|\partial_{x}\mathcal{P}|\leq C\mathcal{P}$, we have
\begin{align}\label{I94}
-(\varsigma\partial_{x}\mathcal{P}\partial_{x}^{s}u,\varphi^{2}\mathcal{U}_{s})
\leq& \|\varsigma\partial_{x}\mathcal{P}\partial_{x}^{s}u\|_{L_{\varphi}^{2}(\Omega)}^{2}
     +\|\mathcal{U}_{s}\|_{L_{\varphi}^{2}(\Omega)}^{2}\\
\leq& C\bar{\varsigma}^{2}\|\partial_{x}^{s}u\|_{L_{\hat{\Psi}}^{2}(\Omega)}^{2}
     +\|\mathcal{U}_{s}\|_{L_{\varphi}^{2}(\Omega)}^{2}\nonumber\\
\leq& C_{*}\bar{\varsigma}^{2}\|\mathcal{P}\mathcal{U}_{s}\|_{L_{\varphi}^{2}(\Omega)}^{2}
     +\|\mathcal{U}_{s}\|_{L_{\varphi}^{2}(\Omega)}^{2}\nonumber\\
\leq& C_{*}\|\mathcal{U}_{s}\|_{L_{\varphi}^{2}(\Omega)}^{2}.\nonumber
\end{align}

Combining \eqref{I91}-\eqref{I94} with \eqref{I9}, we can conclude
\begin{align*}
I_{9}
\leq \frac{1}{8}\|\sqrt{\eta}\partial_{x}\mathcal{U}_{s}\|_{L_{\varphi}^{2}(\Omega)}^{2}
    +C_{*}C_{\theta}\iota^{-1}\|\mathcal{U}_{s}\|_{L_{\varphi}^{2}(\Omega)}^{2}.
\end{align*}

\noindent{\bf \underline{Estimate of~$I_{10}$.}}
A direct calculation gives
\begin{align*}
I_{10}
=& (\kappa\mathcal{P}(\partial_{y}\partial_{x}^{s+2}\mathbb{U}-\rho\partial_{x}^{s+2}\mathbb{U}),\varphi^{2}\mathcal{U}_{s})\\
=& (\kappa\mathcal{P}\partial_{x}(\partial_{y}\partial_{x}^{s+1}\mathbb{U}-\rho\partial_{x}^{s+1}\mathbb{U}),\varphi^{2}\mathcal{U}_{s})
  +(\kappa\mathcal{P}\partial_{x}\rho\partial_{x}^{s+1}\mathbb{U},\varphi^{2}\mathcal{U}_{s})\\
=& (\kappa\mathcal{P}\partial_{x}(\mathcal{P}^{-1}\mathbf{U}_{s+1}),\varphi^{2}\mathcal{U}_{s})
  +(\kappa\mathcal{P}\partial_{x}\rho\partial_{x}^{s+1}\mathbb{U},\varphi^{2}\mathcal{U}_{s}).
\end{align*}
Noticing that $|\mathcal{P}\partial_{x}\mathcal{P}^{-1}|\leq C$, one has
\begin{align*}
(\kappa\mathcal{P}\partial_{x}(\mathcal{P}^{-1}\mathbf{U}_{s+1}),\varphi^{2}\mathcal{U}_{s})
=&~(\kappa\partial_{x}\mathbf{U}_{s+1},\varphi^{2}\mathcal{U}_{s})
  +(\kappa\mathcal{P}\partial_{x}\mathcal{P}^{-1}\mathbf{U}_{s+1},\varphi^{2}\mathcal{U}_{s})\\
\leq&~\tau\bar{\kappa}^{2}\|\vartheta_{c}\partial_{x}\mathbf{U}_{s+1}\|_{L_{\hat{\varphi}}^{2}(\Omega)}^{2}
     +\bar{\kappa}^{2}\|\mathbf{U}_{s+1}\|_{L_{\hat{\varphi}}^{2}(\Omega)}^{2}
     +C_{*}C_{\tau}\|\sqrt{\y}\mathcal{U}_{s}\|_{L_{\varphi}^{2}(\Omega)}^{2},
\end{align*}
Since $|\partial_{x}\rho|\leq C\mathcal{P}$, one can obtain
\begin{align*}
(\kappa\mathcal{P}\partial_{x}\rho\partial_{x}^{s+1}\mathbb{U},\varphi^{2}\mathcal{U}_{s})
\leq&~C_{*}\bar{\kappa}\|\mathcal{P}^{2}\partial_{x}^{s+1}\mathbb{U}\|_{L_{\hat{\varphi}}^{2}(\Omega)}
                 \|\sqrt{\y}\mathcal{U}_{s}\|_{L_{\varphi}^{2}(\Omega)}\\
\leq&~C_{*}\bar{\kappa}\|\mathcal{P}\mathbf{U}_{s+1}\|_{L_{\hat{\varphi}}^{2}(\Omega)}\|\sqrt{\y}\mathcal{U}_{s}\|_{L_{\varphi}^{2}(\Omega)}\\
\leq&~\bar{\kappa}^{2}\tau\|\mathcal{P}\mathbf{U}_{s+1}\|_{L_{\hat{\varphi}}^{2}(\Omega)}^{2}
     +C_{*}C_{\tau}\|\sqrt{\y}\mathcal{U}_{s}\|_{L_{\varphi}^{2}(\Omega)}^{2}.
\end{align*}
To sum up, we have
\begin{align*}
I_{10}
\leq C_{*}C_{\tau}\|\sqrt{\y}\mathcal{U}_{s}\|_{L_{\varphi}^{2}(\Omega)}^{2}
    +\bar{\kappa}^{2}\left(\tau\|\vartheta_{c}\partial_{x}\mathbf{U}_{s+1}\|_{L_{\hat{\varphi}}^{2}(\Omega)}^{2}
               +\|\mathbf{U}_{s+1}\|_{L_{\hat{\varphi}}^{2}(\Omega)}^{2}
               +\tau\|\mathcal{P}\mathbf{U}_{s+1}\|_{L_{\hat{\varphi}}^{2}(\Omega)}^{2}\right).
\end{align*}

\noindent{\bf \underline{Estimate of~$I_{11}$.}} The estimate for $I_{11}$ is exactly the same as that for $I_{7}^{o}$ obtained in the case $n=1$ above, and we have
\begin{align*}
I_{11}
\leq&~C_{*}\iota^{-3}\|\mathcal{U}_{s}\|_{L_{\varphi}^{2}(\Omega)}^{2}
     +C\|u\|_{\hat{H}_{\psi}^{s}(\Omega)}^{2}
     +\frac{1}{8}\|\sqrt{\eta}\partial_{x}\mathcal{U}_{s}\|_{L_{\varphi}^{2}(\Omega)}^{2}
     +\frac{1}{16}\|\partial_{y}\mathcal{U}_{s}\|_{L_{\varphi}^{2}(\Omega)}^{2}
     +\|\mathcal{P}\mathcal{U}_{s}\|_{L_{\varphi}^{2}(\Omega)}^{2}\\
    &+C_{*}\bar{\eta}^{-1}\|u\|_{L_{t}^{4}\mathcal{H}_{\psi,\varphi}^{s}(\Omega_{t})}^{4}\|\mathcal{U}_{s}\|_{L_{\varphi}^{2}(\Omega)}^{2}
     +C_{*}\bar{\eta}^{-1}\|u\|_{L_{t}^{2}\mathcal{H}_{\psi,\varphi}^{s}(\Omega_{t})}^{2}\|\mathcal{U}_{s}\|_{L_{\varphi}^{2}(\Omega)}^{2}.
\end{align*}

\noindent{\bf \underline{Estimate of~$I_{12}$.}}
Observe that $\mathcal{P}=1$ when $\partial_{y}\vartheta_{c}\neq0$, utilizing integration by parts one has
\begin{align*}
I_{12}
= (\partial_{y}\eta\partial_{x}^{s+2}u,\varphi^{2}\mathcal{U}_{s})
 +(\partial_{y}\kappa\partial_{x}^{s+2}\mathbb{U},\varphi^{2}\mathcal{U}_{s})
=-(\partial_{y}\eta\partial_{x}^{s+1}u,\varphi^{2}\partial_{x}\mathcal{U}_{s})
 -(\partial_{y}\kappa\partial_{x}^{s+1}\mathbb{U},\varphi^{2}\partial_{x}\mathcal{U}_{s})
\end{align*}
By using \eqref{theta} and Corollary \ref{ncC3}, we have
\begin{align*}
\|\partial_{y}\vartheta_{c}\partial_{x}^{s+1}u\|_{L_{\varphi}^{2}(\Omega)}
\leq&~C_{\lambda}t\theta\|\y^{-\frac{1}{2}}\vartheta_{c}\partial_{x}^{s+1}u\|_{L_{\varphi_{o}}^{2}(\Omega)}\\
\leq&~C_{\lambda}\sqrt{t}\theta\|\vartheta_{c}\partial_{x}\mathcal{U}_{s}\|_{L_{\varphi}^{2}(\Omega)}
     +C_{\lambda}\theta\|\mathcal{U}_{s}\|_{L_{\varphi}^{2}(\Omega)},
\end{align*}
which gives
\begin{align*}
-(\partial_{y}\eta\partial_{x}^{s+1}u,\varphi^{2}\partial_{x}\mathcal{U}_{s})
\leq& 2\bar{\eta}\|\partial_{y}\vartheta_{c}\partial_{x}^{s+1}u\|_{L_{\varphi}^{2}(\Omega)}
      \|\vartheta_{c}\partial_{x}\mathcal{U}_{s}\|_{L_{\varphi}^{2}(\Omega)}\\
\leq& C_{\lambda}\bar{\eta}\theta
      \left(\sqrt{t}\|\vartheta_{c}\partial_{x}\mathcal{U}_{s}\|_{L_{\varphi}^{2}(\Omega)}+\|\mathcal{U}_{s}\|_{L_{\varphi}^{2}(\Omega)}\right)
      \|\vartheta_{c}\partial_{x}\mathcal{U}_{s}\|_{L_{\varphi}^{2}(\Omega)}\\
\leq& \left(\frac{1}{16}+C_{\lambda}\theta\sqrt{t}\right)\|\sqrt{\eta}\partial_{x}\mathcal{U}_{s}\|_{L_{\varphi}^{2}(\Omega)}^{2}
     +C_{*}C_{\theta}\|\mathcal{U}_{s}\|_{L_{\varphi}^{2}(\Omega)}^{2}.
\end{align*}
Using Corollary \ref{ncC2} and $\sqrt{\theta t\y}\vartheta_{c}\leq C$, one can deduce that
\begin{align*}
-(\partial_{y}\kappa\partial_{x}^{s+1}\mathbb{U},\varphi^{2}\partial_{x}\mathcal{U}_{s})
\leq& C_{\lambda}\sqrt{\bar{\eta}}\theta t\|\vartheta_{c}\partial_{x}^{s+1}\mathbb{U}\|_{L_{\varphi}^{2}(\Omega)}
                 \|\sqrt{\eta}\partial_{x}\mathcal{U}_{s}\|_{L_{\varphi}^{2}(\Omega)}\\
\leq& C_{\lambda}\sqrt{\bar{\eta}\theta t}\|\partial_{x}^{s+1}\mathbb{U}\|_{L_{\hat{\varphi}}^{2}(\Omega)}
                 \|\sqrt{\eta}\partial_{x}\mathcal{U}_{s}\|_{L_{\varphi}^{2}(\Omega)}\\
\leq& C_{\lambda}\sqrt{\bar{\eta}\theta t}\|\mathcal{P}\mathbf{U}_{s+1}\mathbb{U}\|_{L_{\hat{\varphi}}^{2}(\Omega)}
                 \|\sqrt{\eta}\partial_{x}\mathcal{U}_{s}\|_{L_{\varphi}^{2}(\Omega)}\\
\leq& C_{\lambda}C_{\theta}\bar{\eta}\|\mathbf{U}_{s+1}\|_{L_{\hat{\varphi}}^{2}(\Omega)}^{2}
      +\frac{1}{16}\|\sqrt{\eta}\partial_{x}\mathcal{U}_{s}\|_{L_{\varphi}^{2}(\Omega)}.
\end{align*}

In summary,
\begin{align*}
I_{12}
\leq \left(\frac{1}{8}+C_{\lambda}\theta\sqrt{t}\right)\|\sqrt{\eta}\partial_{x}\mathcal{U}_{s}\|_{L_{\varphi}^{2}(\Omega)}^{2}
    +C_{*}C_{\theta}\|\mathcal{U}_{s}\|_{L_{\varphi}^{2}(\Omega)}^{2}
    +C_{\lambda}C_{\theta}\bar{\eta}\|\mathbf{U}_{s+1}\|_{L_{\hat{\varphi}}^{2}(\Omega)}^{2}.
\end{align*}

\noindent{\bf \underline{Estimate of~$I_{13}$.}}
Since $\varphi=\hat{\psi}\sqrt{\zeta_{\lambda}}\leq \frac{C_{\lambda}}{(\bar{\varsigma}+t)^{\frac{\varepsilon_{0}}{2}}}\hat{\psi}$, by using the Cauchy inequality one knows that,
\begin{align*}
I_{13}=(\mathbf{F}_{s},\varphi^{2}\mathcal{U}_{s})
\leq&~\tau\|\mathbf{F}_{s}\|_{L_{\hat{\psi}}^{2}(\Omega)}^{2}
     +\frac{C_{*}C_{\tau}}{(\bar{\varsigma}+t)^{\varepsilon_{0}}}\|\mathcal{U}_{s}\|_{L_{\varphi}^{2}(\Omega)}^{2}\\
\leq&~\tau\|\mathbf{F}_{s}\|_{L_{\hat{\psi}}^{2}(\Omega)}^{2}
     +\frac{C_{*}C_{\tau}}{(\bar{\varsigma}+t)^{\frac{1}{2}+\varepsilon_{0}}}\|u\|_{\mathcal{H}_{\psi,\varphi}^{s}(\Omega)}^{2}.
\end{align*}

Collecting all the estimates of $I_{0}, I_{2},..., I_{13}$ above, from \eqref{he} and \eqref{He0} we can obtain the following proposition.
\begin{prop}\label{highe}
Let $\lambda$ and $\Lambda$ be large enough and $\iota$ small enough such that Lemma \ref{weighte} holds on the time interval $[0,T_{*}]$. Assume that Assumption \ref{MA} holds for $T_{*}$, then for any $\theta>0, \tau>0$ and $t\in[0,T_{*}]$, we have the following a priori estimate for the solution to the problem \eqref{us},
\begin{align}\label{he'}
&\frac{1}{2}\frac{d}{dt}\|\mathcal{U}_{s}\|_{L_{\varphi}^{2}(\Omega)}^{2}
+\frac{\lambda}{128}\|\sqrt{\omega_{\lambda}}\mathcal{U}_{s}\|_{L_{\varphi}^{2}(\Omega)}^{2}
+\frac{3}{4}\Lambda\delta\|\sqrt{\y}\mathcal{U}_{s}\|_{L_{\varphi}^{2}(\Omega)}^{2}
+\frac{1}{4}\|\partial_{y}\mathcal{U}_{s}\|_{L_{\varphi}^{2}(\Omega)}^{2}
+\frac{5}{8}\|\sqrt{\eta}\partial_{x}\mathcal{U}_{s}\|_{L_{\varphi}^{2}(\Omega)}^{2}\\
\leq&~C\|\mathcal{P}\mathcal{U}_{s}\|_{L_{\varphi}^{2}(\Omega)}^{2}
     +C_{*}C_{\tau}C_{\theta}C_{\iota}\|\mathcal{U}_{s}\|_{L_{\varphi}^{2}(\Omega)}^{2}
     +\frac{C_{*}C_{\tau}}{(\bar{\varsigma}+t)^{\frac{1}{2}+\varepsilon_{0}}}\|u\|_{\mathcal{H}_{\psi,\varphi}^{s}(\Omega)}^{2}
     +C_{\lambda}\theta\sqrt{t}\|\sqrt{\eta}\partial_{x}\mathcal{U}_{s}\|_{L_{\varphi}^{2}(\Omega)}^{2}\nonumber\\
    &+C_{*}\left(1+\|\mathbf{U}_{s+1}\|_{L_{\hat{\varphi}}^{2}(\Omega)}^{2}
        +\|\mathcal{P}\mathbf{U}_{s+1}\|_{L_{\hat{\varphi}}^{2}(\Omega)}^{2}
        +\|\mathcal{P}\mathbf{U}_{s}\|_{L_{\varphi}^{2}(\Omega)}^{2}
        +\|\mathbf{U}_{s}\|_{L_{\varphi}^{2}(\Omega)}^{2}
        +\|\partial_{y}\mathbf{U}_{s}\|_{L_{\varphi}^{2}(\Omega)}^{2}\right)
       \|u\|_{\mathcal{H}_{\psi,\varphi}^{s}(\Omega)}^{2}\nonumber\\
    &+C_{*}\|u\|_{\hat{H}_{\psi,\y}^{s}(\Omega)}^{2}
     +C_{*}C_{\tau}\|\sqrt{\y}\mathcal{U}_{s}\|_{L_{\varphi}^{2}(\Omega)}^{2}
     +C_{*}t\|u\|_{L_{t}^{2}\mathcal{H}_{\psi,\varphi}^{s}(\Omega_{t})}^{2}
     +C_{*}\bar{\eta}^{2}
     +C_{*}\check{c}^{2}\iota^{-2}\nonumber\\
    &+\bar{\kappa}^{2}\left(\tau\|\vartheta_{c}\partial_{x}\mathbf{U}_{s+1}\|_{L_{\hat{\varphi}}^{2}(\Omega)}^{2}
               +\|\mathbf{U}_{s+1}\|_{L_{\hat{\varphi}}^{2}(\Omega)}^{2}
               +\tau\|\mathcal{P}\mathbf{U}_{s+1}\|_{L_{\hat{\varphi}}^{2}(\Omega)}^{2}\right)
     +C_{\lambda}C_{\theta}\bar{\eta}\|\mathbf{U}_{s+1}\|_{L_{\hat{\varphi}}^{2}(\Omega)}^{2}\nonumber\\
    &+C_{*}\bar{\eta}^{-1}\|u\|_{L_{t}^{4}\mathcal{H}_{\psi,\varphi}^{s}(\Omega_{t})}^{4}\|\mathcal{U}_{s}\|_{L_{\varphi}^{2}(\Omega)}^{2}
     +C_{*}\bar{\eta}^{-1}\|u\|_{L_{t}^{2}\mathcal{H}_{\psi,\varphi}^{s}(\Omega_{t})}^{2}\|\mathcal{U}_{s}\|_{L_{\varphi}^{2}(\Omega)}^{2}\nonumber\\
    &+\tau\|\mathbf{F}_{s}\|_{L_{\hat{\psi}}^{2}(\Omega)}^{2}
     +C\|\partial_{x}^{s}F\|_{\mathring{H}^{0}(\mathbb{T})}^{2}.\nonumber
\end{align}
\end{prop}
~~~~~~~~~~~

\subsection{Uniform a priori estimates of the corrected increment solutions}~

In this subsection, we shall establish the uniform a priori estimates of the corrected increment solutions by using the estimates established in subsection 4.2 and 4.4.

Combining Proposition \ref{lowe} with Proposition \ref{highe}, one can obtain that under the assumption of Theorem \ref{wpae}, for any $\theta>0, \tau>0$ and $t\in[0,T_{*}]$ it follows that,
\begin{align}\label{upe-1}
& \frac{1}{2}\frac{d}{dt}\|u\|_{\mathcal{H}_{\psi,\varphi}^{s}(\Omega)}^{2}
 +\Lambda\delta\|u\|_{\hat{H}_{\psi,\y}^{s}(\Omega)}^{2}
 +\frac{3}{4}\Lambda\delta\|\sqrt{\y}\mathcal{U}_{s}\|_{L_{\varphi}^{2}(\Omega)}^{2}
 +\frac{\lambda}{128}\|\sqrt{\omega_{\lambda}}\mathcal{U}_{s}\|_{L_{\varphi}^{2}(\Omega)}^{2}\\
&+\frac{5}{16}\|\partial_{y}u\|_{\hat{H}_{\psi}^{s}(\Omega)}^{2}
 +\frac{1}{4}\|\partial_{y}\mathcal{U}_{s}\|_{L_{\varphi}^{2}(\Omega)}^{2}
 +\frac{15}{16}\|\partial_{x}u\|_{\hat{H}_{\psi,\eta}^{s}(\Omega)}^{2}
 +\frac{5}{8}\|\sqrt{\eta}\partial_{x}\mathcal{U}_{s}\|_{L_{\varphi}^{2}(\Omega)}^{2}\nonumber\\
\leq&~\hat{C}\|\mathcal{P}\mathcal{U}_{s}\|_{L_{\varphi}^{2}(\Omega)}^{2}
     +C_{*}\|u\|_{\hat{H}_{\psi,\y}^{s}(\Omega)}^{2}
     +C_{*}C_{\tau}\|\sqrt{\y}\mathcal{U}_{s}\|_{L_{\varphi}^{2}(\Omega)}^{2}
     +C_{\lambda}\theta\sqrt{t}\|\sqrt{\eta}\partial_{x}\mathcal{U}_{s}\|_{L_{\varphi}^{2}(\Omega)}^{2}\nonumber\\
    &+C_{*}\left(\|\mathbb{U}\|_{\mathcal{H}_{\psi,\varphi}^{s}(\Omega)}^{2}
                +\|\partial_{y}\mathbb{U}\|_{\hat{H}_{\psi}^{s}(\Omega)}^{2}
                +\|\partial_{y}\mathbf{U}_{s}\|_{L_{\varphi}^{2}(\Omega)}^{2}
                +\|\mathcal{P}\mathbf{U}_{s}\|_{L_{\varphi}^{2}(\Omega)}^{2}
                +C_{\tau}C_{\theta}C_{\iota}\right)
          \|u\|_{\mathcal{H}_{\psi,\varphi}^{s}(\Omega)}^{2}\nonumber\\
    &+C_{*}\left(\|\mathbf{U}_{s+1}\|_{L_{\hat{\varphi}}^{2}(\Omega)}^{2}+\|\mathcal{P}\mathbf{U}_{s+1}\|_{L_{\hat{\varphi}}^{2}(\Omega)}^{2}\right)
       \|u\|_{\mathcal{H}_{\psi,\varphi}^{s}(\Omega)}^{2}
     +\frac{C_{*}C_{\tau}}{(\bar{\varsigma}+t)^{\frac{1}{2}+\varepsilon_{0}}}\|u\|_{\mathcal{H}_{\psi,\varphi}^{s}(\Omega)}^{2}
     +C_{*}\|u\|_{\mathcal{H}_{\psi,\varphi}^{s}(\Omega)}^{4}\nonumber\\
    &+C_{*}\bar{\eta}^{-1}\|u\|_{L_{t}^{4}\mathcal{H}_{\psi,\varphi}^{s}(\Omega_{t})}^{4}\|\mathcal{U}_{s}\|_{L_{\varphi}^{2}(\Omega)}^{2}
     +C_{*}\bar{\eta}^{-1}\|u\|_{L_{t}^{2}\mathcal{H}_{\psi,\varphi}^{s}(\Omega_{t})}^{2}\|\mathcal{U}_{s}\|_{L_{\varphi}^{2}(\Omega)}^{2}
     +C_{*}t\|u\|_{L_{t}^{2}\mathcal{H}_{\psi,\varphi}^{s}(\Omega_{t})}^{2}\nonumber\\
    &+\tau\|\mathbf{F}_{s}\|_{L_{\hat{\psi}}^{2}(\Omega)}^{2}
     +\tau\|F\|_{\hat{H}_{\psi}^{s}(\Omega)}^{2}
     +C_{*}\sum_{i=0}^{\frac{s-1}{2}}\|\partial_{t}^{i}F\|_{\mathring{H}^{s-1-2i}(\mathbb{T})}^{2}
     +C\|\partial_{x}^{s}F\|_{\mathring{H}^{0}(\mathbb{T})}^{2}\nonumber\\
    &+C\bar{\kappa}^{2}\left(\|\mathbb{U}\|_{\mathcal{H}_{\psi,\varphi}^{s}(\Omega)}^{2}
                +\tau\|\partial_{y}\mathbb{U}\|_{\hat{H}_{\psi}^{s}(\Omega)}^{2}
                +C_{\tau}C_{*}
                +\tau\|\partial_{y}\mathbf{U}_{s}\|_{L_{\varphi}^{2}(\Omega)}^{2}
                +\tau\|\mathcal{P}\mathbf{U}_{s}\|_{L_{\varphi}^{2}(\Omega)}^{2}\right)\nonumber\\
    &+C\bar{\kappa}^{2}\left(\|\mathbf{U}_{s+1}\|_{L_{\hat{\varphi}}^{2}(\Omega)}^{2}
                     +\tau\|\mathcal{P}\mathbf{U}_{s+1}\|_{L_{\hat{\varphi}}^{2}(\Omega)}^{2}
                     +\tau\|\vartheta_{c}\partial_{x}\mathbf{U}_{s+1}\|_{L_{\hat{\varphi}}^{2}(\Omega)}^{2}
                     +\tau\|\partial_{y}\mathbf{U}_{s+1}\|_{L_{\hat{\varphi}}^{2}(\Omega)}^{2}
                     \right)\nonumber\\
    &+C_{\lambda}C_{\theta}\bar{\eta}\|\mathbf{U}_{s+1}\|_{L_{\hat{\varphi}}^{2}(\Omega)}^{2}
     +C\bar{\eta}\|\mathbb{U}\|_{\hat{H}_{\psi}^{s}(\Omega)}^{2}
     +C_{*}\bar{\eta}^{2}
     +C_{*}\check{c}\iota^{-2}\nonumber,
\end{align}
where $\hat{C}$ is a positive constant independent of $\lambda$ and $\mathcal{Z}$.

Since $\mathcal{P}\leq C\sqrt{\omega_{\lambda}}+1$, one can take $\lambda$ large enough such that
\begin{align}\label{pusc}
\|\mathcal{P}\mathcal{U}_{s}\|_{L_{\varphi}^{2}(\Omega)}^{2}
\leq \frac{\lambda}{256}\|\sqrt{\omega_{\lambda}}\mathcal{U}_{s}\|_{L_{\varphi}^{2}(\Omega)}^{2}
    +\|\mathcal{U}_{s}\|_{L_{\varphi}^{2}(\Omega)}^{2}.
\end{align}
moreover,
\begin{align*}
\hat{C}\|\mathcal{P}\mathcal{U}_{s}\|_{L_{\varphi}^{2}(\Omega)}^{2}
\leq \frac{\lambda}{256}\|\sqrt{\omega_{\lambda}}\mathcal{U}_{s}\|_{L_{\varphi}^{2}(\Omega)}^{2}
    +\hat{C}\|\mathcal{U}_{s}\|_{L_{\varphi}^{2}(\Omega)}^{2}.
\end{align*}
With $\lambda$ fixed, one can take $\hat{t}^{*}\in(0,T_{*}]$ small enough such that $C_{\lambda}\theta\sqrt{t}\leq\frac{1}{8}$ holds for any $t\in[0,\hat{t}^{*}]$.
Then for any $t\in[0,\hat{t}^{*}]$, from \eqref{upe-1} one can deduce that
\begin{align}\label{upe-2}
& \frac{d}{dt}\|u\|_{\mathcal{H}_{\psi,\varphi}^{s}(\Omega)}^{2}
 +2\Lambda\delta\|u\|_{\hat{H}_{\psi,\y}^{s}(\Omega)}^{2}
 +\frac{3}{2}\Lambda\delta\|\sqrt{\y}\mathcal{U}_{s}\|_{L_{\varphi}^{2}(\Omega)}^{2}
 +\frac{\lambda}{128}\|\sqrt{\omega_{\lambda}}\mathcal{U}_{s}\|_{L_{\varphi}^{2}(\Omega)}^{2}\\
&+\frac{5}{8}\|\partial_{y}u\|_{\hat{H}_{\psi}^{s}(\Omega)}^{2}
 +\frac{1}{2}\|\partial_{y}\mathcal{U}_{s}\|_{L_{\varphi}^{2}(\Omega)}^{2}
 +\frac{15}{8}\|\partial_{x}u\|_{\hat{H}_{\psi,\eta}^{s}(\Omega)}^{2}
 +\|\sqrt{\eta}\partial_{x}\mathcal{U}_{s}\|_{L_{\varphi}^{2}(\Omega)}^{2}\nonumber\\
\leq&~C_{*}\|u\|_{\hat{H}_{\psi,\y}^{s}(\Omega)}^{2}
     +C_{*}C_{\tau}\|\sqrt{\y}\mathcal{U}_{s}\|_{L_{\varphi}^{2}(\Omega)}^{2}
     +G_{1}(t)\|u\|_{\mathcal{H}_{\psi,\varphi}^{s}(\Omega)}^{2}
     +C_{*}\|u\|_{\mathcal{H}_{\psi,\varphi}^{s}(\Omega)}^{4}\nonumber\\
    &+C_{*}\bar{\eta}^{-1}\|u\|_{L_{t}^{4}\mathcal{H}_{\psi,\varphi}^{s}(\Omega_{t})}^{4}\|\mathcal{U}_{s}\|_{L_{\varphi}^{2}(\Omega)}^{2}
     +C_{*}\bar{\eta}^{-1}\|u\|_{L_{t}^{2}\mathcal{H}_{\psi,\varphi}^{s}(\Omega_{t})}^{2}\|\mathcal{U}_{s}\|_{L_{\varphi}^{2}(\Omega)}^{2}
     +C_{*}t\|u\|_{L_{t}^{2}\mathcal{H}_{\psi,\varphi}^{s}(\Omega)}^{2}
     +G_{2}(t)\nonumber,
\end{align}
where
\begin{align*}
G_{1}(t)
=& C_{*}\left(\|\mathbb{U}\|_{\mathcal{H}_{\psi,\varphi}^{s}(\Omega)}^{2}
                +\|\partial_{y}\mathbb{U}\|_{\hat{H}_{\psi}^{s}(\Omega)}^{2}
                +\|\partial_{y}\mathbf{U}_{s}\|_{L_{\varphi}^{2}(\Omega)}^{2}
                +\|\mathcal{P}\mathbf{U}_{s}\|_{L_{\varphi}^{2}(\Omega)}^{2}
                +C_{\tau}C_{\iota}\right)\\
 &+C_{*}\left(\|\mathbf{U}_{s+1}\|_{L_{\hat{\varphi}}^{2}(\Omega)}^{2}+\|\mathcal{P}\mathbf{U}_{s+1}\|_{L_{\hat{\varphi}}^{2}(\Omega)}^{2}\right)
  +\frac{C_{*}C_{\tau}}{(\bar{\varsigma}+t)^{\frac{1}{2}+\varepsilon_{0}}},
\end{align*}
and
\begin{align*}
G_{2}(t)
=&2\tau\|\mathbf{F}_{s}\|_{L_{\hat{\psi}}^{2}(\Omega)}^{2}
     +2\tau\|F\|_{\hat{H}_{\psi}^{s}(\Omega)}^{2}
     +C_{*}\sum_{i=0}^{\frac{s-1}{2}}\|\partial_{t}^{i}F\|_{\mathring{H}^{s-1-2i}(\mathbb{T})}^{2}
     +C\|\partial_{x}^{s}F\|_{\mathring{H}^{0}(\mathbb{T})}^{2},\\
    &+C\bar{\kappa}^{2}\left(\|\mathbb{U}\|_{\mathcal{H}_{\psi,\varphi}^{s}(\Omega)}^{2}
                +\tau\|\partial_{y}\mathbb{U}\|_{\hat{H}_{\psi}^{s}(\Omega)}^{2}
                +C_{*}C_{\tau}
                +\tau\|\partial_{y}\mathbf{U}_{s}\|_{L_{\varphi}^{2}(\Omega)}^{2}
                +\tau\|\mathcal{P}\mathbf{U}_{s}\|_{L_{\varphi}^{2}(\Omega)}^{2}\right)\\
    &+C\bar{\kappa}^{2}\left(\|\mathbf{U}_{s+1}\|_{L_{\hat{\varphi}}^{2}(\Omega)}^{2}
                     +\tau\|\mathcal{P}\mathbf{U}_{s+1}\|_{L_{\hat{\varphi}}^{2}(\Omega)}^{2}
                     +\tau\|\vartheta_{c}\partial_{x}\mathbf{U}_{s+1}\|_{L_{\hat{\varphi}}^{2}(\Omega)}^{2}
                     +\tau\|\partial_{y}\mathbf{U}_{s+1}\|_{L_{\hat{\varphi}}^{2}(\Omega)}^{2}
                     \right)\\
    &+C_{\lambda}C_{\theta}\bar{\eta}\|\mathbf{U}_{s+1}\|_{L_{\hat{\varphi}}^{2}(\Omega)}^{2}
     +C\bar{\eta}\|\mathbb{U}\|_{\hat{H}_{\psi}^{s}(\Omega)}^{2}
     +C_{*}\bar{\eta}^{2}
     +C_{*}\check{c}\iota^{-2}\nonumber.
\end{align*}
Let $g(t)=e^{\int_{0}^{t}G_{1}(\sigma)d\sigma}$. By taking $\tau$ and $\check{t}^{*}\in(0,\hat{t}^{*}]$ small enough, one can obtain
\begin{align}\label{G2}
g(t)\int_{0}^{t}G_{2}(\sigma)g^{-1}(\sigma)d\sigma
\leq \frac{(1-\epsilon_{0})^{2}\bar{\eta}\mathcal{Z}^{2}}{32\epsilon_{0}^{2}},~\forall t\in(0,\check{t}^{*}].
\end{align}
With $\tau$ fixed, one can choose $\Lambda>\delta^{-1}$ large enough in \eqref{upe-2} to obtain the following inequality in $(0,\check{t}^{*}]$,
\begin{align}\label{upe-3}
& \frac{d}{dt}\|u\|_{\mathcal{H}_{\psi,\varphi}^{s}(\Omega)}^{2}
 +\Lambda\delta\|u\|_{\hat{H}_{\psi,\y}^{s}(\Omega)}^{2}
 +\Lambda\delta\|\sqrt{\y}\mathcal{U}_{s}\|_{L_{\varphi}^{2}(\Omega)}^{2}
 +\frac{\lambda}{128}\|\sqrt{\omega_{\lambda}}\mathcal{U}_{s}\|_{L_{\varphi}^{2}(\Omega)}^{2}\\
&+\frac{5}{8}\|\partial_{y}u\|_{\hat{H}_{\psi}^{s}(\Omega)}^{2}
 +\frac{1}{2}\|\partial_{y}\mathcal{U}_{s}\|_{L_{\varphi}^{2}(\Omega)}^{2}
 +\frac{15}{8}\|\partial_{x}u\|_{\hat{H}_{\psi,\eta}^{s}(\Omega)}^{2}
 +\|\sqrt{\eta}\partial_{x}\mathcal{U}_{s}\|_{L_{\varphi}^{2}(\Omega)}^{2}\nonumber\\
\leq&~G_{1}(t)\|u\|_{\mathcal{H}_{\psi,\varphi}^{s}(\Omega)}^{2}
     +C_{*}\|u\|_{\mathcal{H}_{\psi,\varphi}^{s}(\Omega)}^{4}
     +C_{*}\bar{\eta}^{-1}\|u\|_{L_{t}^{4}\mathcal{H}_{\psi,\varphi}^{s}(\Omega_{t})}^{4}\|\mathcal{U}_{s}\|_{L_{\varphi}^{2}(\Omega)}^{2}\nonumber\\
    &+C_{*}\bar{\eta}^{-1}\|u\|_{L_{t}^{2}\mathcal{H}_{\psi,\varphi}^{s}(\Omega_{t})}^{2}\|\mathcal{U}_{s}\|_{L_{\varphi}^{2}(\Omega)}^{2}
     +C_{*}t\|u\|_{L_{t}^{2}\mathcal{H}_{\psi,\varphi}^{s}(\Omega)}^{2}
     +G_{2}(t)\nonumber.
\end{align}

For any $t\in [0,\check{t}^{*}]$, let $S(t)=\|u\|_{\mathcal{H}_{\psi,\varphi}^{s}(\Omega)}^{2}$, and
\begin{align*}
W(t)
=& \Lambda\delta\|u\|_{\hat{H}_{\psi,\y}^{s}(\Omega)}^{2}
  +\|\partial_{y}u\|_{\hat{H}_{\psi}^{s}(\Omega)}^{2}
  +\|\partial_{x}u\|_{\hat{H}_{\psi,\eta}^{s}(\Omega)}^{2}
  +\frac{\lambda}{64}\|\sqrt{\omega_{\lambda}}\mathcal{U}_{s}\|_{L_{\varphi}^{2}(\Omega)}^{2}\\
 &+\Lambda\delta\|\sqrt{\y}\mathcal{U}_{s}\|_{L_{\varphi}^{2}(\Omega)}^{2}
  +\|\partial_{y}\mathcal{U}_{s}\|_{L_{\varphi}^{2}(\Omega)}^{2}
  +\|\sqrt{\eta}\partial_{x}\mathcal{U}_{s}\|_{L_{\varphi}^{2}(\Omega)}^{2},
\end{align*}
then \eqref{upe-3} gives that for any $t\in [0,\check{t}^{*}]$, there holds
\begin{align}\label{upe-4}
\frac{d}{dt}S(t)+\frac{1}{2}W(t)
\leq& G_{1}(t)S(t)
     +C_{*}S^{2}(t)
     +C_{*}\bar{\eta}^{-1}S(t)\int_{0}^{t}S^{2}(\sigma)d\sigma\\
    &+C_{*}\bar{\eta}^{-1}S(t)\int_{0}^{t}S(\sigma)d\sigma
     +C_{*}t\int_{0}^{t}S(\sigma)d\sigma
     +G_{2}(t).\nonumber
\end{align}
Let $\mathcal{S}(t)=\bar{\eta}^{-1}g^{-1}(t)S(t)$, $\mathcal{W}(t)=\bar{\eta}^{-1}g^{-1}(t)W(t)$. From \eqref{upe-4} one knows,
\begin{align}\label{upe-5}
\frac{d}{dt}\mathcal{S}(t)+\frac{1}{2}\mathcal{W}(t)
\leq&~C_{*}\bar{\eta}g(t)\mathcal{S}^{2}(t)
     +C_{*}\bar{\eta}\int_{0}^{t}g^{2}(\sigma)\mathcal{S}^{2}(\sigma)d\sigma\mathcal{S}(t)
     +C_{*}\int_{0}^{t}g(\sigma)\mathcal{S}(\sigma)d\sigma\mathcal{S}(t)\\
    &+C_{*}tg^{-1}(t)\int_{0}^{t}g(\sigma)\mathcal{S}(\sigma)d\sigma
     +\bar{\eta}^{-1}g^{-1}(t)G_{2}(t),~\forall t\in [0,\check{t}^{*}].\nonumber
\end{align}
For any $t\in[0,\check{t}^{*}]$, integration \eqref{upe-5} over $[0,t]$, one gets
\begin{align}\label{upe-6}
\mathcal{S}(t)+\frac{1}{2}\int_{0}^{t}\mathcal{W}(\sigma)d\sigma
\leq&~C_{*}\bar{\eta}\int_{0}^{t}g(\sigma)\mathcal{S}^{2}(\sigma)d\sigma
     +C_{*}\bar{\eta}\int_{0}^{t}g^{2}(\sigma)\mathcal{S}^{2}(\sigma)d\sigma\int_{0}^{t}\mathcal{S}(\sigma)d\sigma\\
    &+C_{*}\int_{0}^{t}g(\sigma)\mathcal{S}(\sigma)d\sigma\int_{0}^{t}\mathcal{S}(\sigma)d\sigma
     +C_{*}\int_{0}^{t}\sigma g^{-1}(\sigma)d\sigma\int_{0}^{t}g(\sigma)\mathcal{S}(\sigma)d\sigma\nonumber\\
    &+\int_{0}^{t}\bar{\eta}^{-1}G_{2}(\sigma)g^{-1}(\sigma)d\sigma
     +\mathcal{S}(0).\nonumber
\end{align}

Under the assumption of Theorem \ref{wpae}, one has $\mathcal{S}(0)\leq C_{*}$ and
\begin{align*}
1\leq g(t)\leq C_{*},~\forall t\in [0,\check{t}^{*}],
\end{align*}
with \eqref{G2} in addition, one can deduce from \eqref{upe-6} that for any $t\in[0,\check{t}^{*}]$,
\begin{align}\label{upe-6'}
\mathcal{S}(t)
\leq&~C_{*}\int_{0}^{t}\mathcal{S}^{2}(\sigma)d\sigma
     +C_{*}\int_{0}^{t}\mathcal{S}^{2}(\sigma)d\sigma\int_{0}^{t}\mathcal{S}(\sigma)d\sigma
     +C_{*}\left(\int_{0}^{t}\mathcal{S}(\sigma)d\sigma\right)^{2}\\
    &+C_{*}\int_{0}^{t}\mathcal{S}(\sigma)d\sigma
     +C_{*}\nonumber\\
\leq&~C_{*}\left(\int_{0}^{t}\mathcal{S}^{2}(\sigma)d\sigma\right)^{\frac{3}{2}}
     +C_{*}\int_{0}^{t}\mathcal{S}^{2}(\sigma)d\sigma
     +C_{*}\left(\int_{0}^{t}\mathcal{S}^{2}(\sigma)d\sigma\right)^{\frac{1}{2}}
     +C_{*}\nonumber
\end{align}

Let $\mathfrak{S}(t)=\int_{0}^{t}\mathcal{S}^{2}(\sigma)d\sigma$. From \eqref{upe-6'} one knows that there exist a positive constant $\bar{C}$ such that for any $t\in[0,\bar{t}]$ with  $\bar{t}^{*}=\min\{\check{t}^{*},\frac{1}{4\bar{C}}\}$, there holds
\begin{align*}
\frac{d}{dt}\mathfrak{S}(t)
\leq&~\bar{C}\mathfrak{S}^{3}(t)
     +\bar{C}\mathfrak{S}^{2}(t)
     +\bar{C}\mathfrak{S}(t)
     +\bar{C}\\
\leq&~\bar{C}(\mathfrak{S}(t)+1)^{3},
\end{align*}
which gives,
\begin{align}\label{upe-7}
\mathfrak{S}(t)\leq \frac{\bar{C}t}{1-2\bar{C}t}, ~\forall t\in [0,\bar{t}^{*}].
\end{align}
Substitute \eqref{upe-7} into \eqref{upe-6}, one knows that there is a positive constant $\mathfrak{C}$ such that
\begin{align*}
\mathcal{S}(t)+\frac{1}{2}\int_{0}^{t}\mathcal{W}(\sigma)d\sigma
\leq& \frac{\mathfrak{C}\bar{\eta}t}{1-2\bar{C}t}
     +\frac{\mathfrak{C}t^{2}}{1-2\bar{C}t}
     +\int_{0}^{t}\bar{\eta}^{-1}G_{2}(\mathfrak{\sigma})g^{-1}(\sigma)d\sigma
     +\mathcal{S}(0),~\forall t\in[0,\bar{t}^{*}],
\end{align*}
which gives by the definition of $\mathcal{S}(t)$ and $\mathcal{W}(t)$ that
\begin{align*}
S(t)+\frac{1}{2}g(t)\int_{0}^{t}g^{-1}(\sigma)W(\sigma)d\sigma
\leq \frac{\mathfrak{C}\bar{\eta}g(t)(\bar{\eta}t+t^{2})}{1-2\bar{C}t}
     +g(t)\int_{0}^{t}G_{2}(\sigma)g^{-1}(\sigma)d\sigma
     +S(0)g(t),~\forall t\in[0,\bar{t}^{*}],
\end{align*}
Noticing that $g(t)$ is nondecreasing, one has
\begin{align}\label{upe-8}
S(t)+\frac{1}{2}\int_{0}^{t}W(\sigma)d\sigma
\leq \frac{\mathfrak{C}\bar{\eta}g(t)(\bar{\eta}t+t^{2})}{1-2\bar{C}t}
     +g(t)\int_{0}^{t}G_{2}(\sigma)g^{-1}(\sigma)d\sigma
     +S(0)g(t),~\forall t\in[0,\bar{t}^{*}].
\end{align}
Since $\check{u}\mathbb{U}=0$ in $\Omega_{T^{*}}$, recall the definition of $g(t)$ and $S(t)$, one can take $t^{*}\in(0,\bar{t}^{*}]$ small enough such that
\begin{align}\label{t*}
\frac{\mathfrak{C}\bar{\eta}g(t)(\bar{\eta}t+t^{2})}{1-2\bar{C}t}+S(0)g(t)
\leq \frac{(1-\epsilon_{0})^{2}\bar{\eta}\mathcal{Z}^{2}}{32\epsilon_{0}^{2}},~\forall t\in[0,t^{*}].
\end{align}
Substituting \eqref{G2} and \eqref{t*} into \eqref{upe-8}, one can obtain
\begin{align*}
S(t)+\frac{1}{2}\int_{0}^{t}W(\sigma)d\sigma
\leq \frac{(1-\epsilon_{0})^{2}\bar{\eta}\mathcal{Z}^{2}}{16\epsilon_{0}^{2}},~\forall t\in[0,t^{*}],
\end{align*}
which together with \eqref{pusc} gives \eqref{priore} for $t_{v}=t^{*}$.

To establish \eqref{priore} for any $t_{v}\in[0,t^{*}]$, as can be seen from our proof above, it is sufficient to assume that Assumption \ref{MA} hold for $t_{v}$ and $\check{u}\mathbb{U}=0$ in $\Omega_{t_{v}}$, so we get Theorem \ref{wpae}.

~~~~~~~~~~

\section{A posteriori estimates of the approximate solutions}

In Section 4, we establish the uniform estimate (in $\bar{\eta}$) for the corrected increment solution $u$ under some conditions on the approximate solution $\mathfrak{U}$. The approximate solution $\mathfrak{U}$ to the problem \eqref{Prandtl} is obtained by a iterative method introduced in Section 3. To ensure that the iterative process can proceed, we need to prove that these conditions are satisfied by the approximate solution of various order throughout the iteration. For this purpose, we shall establish the uniform a posteriori estimates of $\mathfrak{U}$ in this section.

Assume that $u$ is the solution to the following problem in $\Omega_{t^{*}}$ obtained by using Theorem \ref{wpae},
\begin{equation}\label{app'}
\begin{cases}
\partial_{t}u
+\tilde{\mathfrak{U}}\partial_{x}u+u\partial_{x}\tilde{\mathfrak{U}}
-\partial_{y}^{-1}[\partial_{x}u](\partial_{y}\tilde{\mathfrak{U}}+\tilde{\varsigma})
-\partial_{y}^{-1}[\partial_{x}\tilde{\mathfrak{U}}]\partial_{y}u
-h\partial_{y}^{-1}[\partial_{x}u]\partial_{y}u\\
\quad
=\partial_{y}^{2}u+\eta\partial_{x}^{2}u+\kappa\partial_{x}^{2}\tilde{\mathbb{U}}+F,\\
u|_{y=0}=0,~\lim\limits_{y\to+\infty} u=0,\\
u|_{t=0}=\check{u}(x,y).
\end{cases}
\end{equation}
where $\tilde{\mathfrak{U}}$ is a approximate solution to the problem \eqref{Prandtl} which is obtained by the iterative method introduced in Section 3, $\tilde{\varsigma}(y)=\tilde{\bar{\varsigma}}\chi_{c}(y)$, $\eta(y)=\bar{\eta}\vartheta_{c}^{2}(y)$ and $\kappa(y)=\kappa_{*}\vartheta_{c}^{2}(y)$ with $\tilde{\bar{\varsigma}}>0, \bar{\eta}>0, \kappa_{*}$ constants.

Utilizing $u$ and $\tilde{\mathfrak{U}}$, one can obtain the next order approximate solution to the problem \eqref{Prandtl}, that is, $\mathfrak{U}=u+\tilde{\mathfrak{U}}$, which obeys the following problem in $\Omega_{t^{*}}$.
\begin{equation}\label{App+1}
\begin{cases}
\partial_{t}\mathfrak{U}
+\mathfrak{U}\partial_{x}\mathfrak{U}
-\partial_{y}^{-1}[\partial_{x}\mathfrak{U}]\partial_{y}\mathfrak{U}
=\partial_{y}^{2}\mathfrak{U}+\eta\partial_{x}^{2}\mathbb{U}-\partial_{x}P-\hat{F},\\
\mathfrak{U}|_{y=0}=0,~\lim\limits_{y\to+\infty}\mathfrak{U}=U(t,x),\\
\mathfrak{U}|_{t=0}=u_{in}(x,y).
\end{cases}
\end{equation}
where
\begin{align}\label{hatF}
\hat{F}=-\partial_{y}^{-1}[\partial_{x}u]\tilde{\varsigma}-u\partial_{x}u+(1-h)\partial_{y}^{-1}[\partial_{x}u]\partial_{y}u,
\end{align}
and $\mathbb{U}=\mathfrak{U}-u_{0}$.

Let $\tilde{\mathbb{U}}=\mathbb{U}-u, \tilde{\eta}=\epsilon_{0}^{-2}\eta$ and $\tilde{\varsigma}=\epsilon_{0}^{-4}\varsigma$. For convenience, we write $\tilde{\bar{\eta}}=\epsilon_{0}^{-2}\bar{\eta}$ and $\tilde{\bar{\varsigma}}=\epsilon_{0}^{-4}\bar{\varsigma}$.
Let $\rho$ and $\mathcal{P}$ be defined as in \eqref{rhoPd}.
For $r=s$ or $s+1$, let $\mathbf{U}_{r}$ be defined as at the beginning of Section 4, and correspondingly, let
$\tilde{\mathbf{U}}_{r}=\tilde{\mathcal{P}}\left(\partial_{y}\partial_{x}^{r}\tilde{\mathbb{U}}-\tilde{\rho}\partial_{x}^{r}\tilde{\mathbb{U}}\right)$,
with
$$\tilde\rho(t,x,y)=\frac{\partial_{y}^{2}\tilde{\mathfrak{U}}+\tilde{\varsigma}'}{\partial_{y}\tilde{\mathfrak{U}}+\tilde{\varsigma}},
~\text{and,}~
\tilde{\mathcal{P}}(t,x,y)=\frac{\chi(y)}{\sqrt{\partial_{y}\tilde{\mathfrak{U}}+\tilde{\varsigma}}}+1-\chi(y).$$
Regarding the weight function, we set $\tilde{\varphi}(t,y)=\y^{-\frac{1}{2}}\psi(t,y)\sqrt{\tilde{\zeta}_{\lambda}(t,y)}$, where
\begin{align*}
\tilde{\zeta}_{\lambda}(t,y)
=&\chi(y)\left[\frac{e^{-\frac{16\lambda^{4}(\tilde{\bar{\varsigma}}+t)}{(y-y_{*})^{2}+16\lambda^{2}(\tilde{\bar{\varsigma}}+t)}}}{[\frac{\lambda^{-2}}{16}(y-y_{*})^{2}+\tilde{\bar{\varsigma}}+t]^{\varepsilon_{0}}}
  +\left(\frac{\lambda}{\varepsilon_{0}}\right)^{\lambda}e^{-\lambda^{2}}H_{\lambda}\left(\frac{y-y_{*}}{\sqrt{2(\tilde{\bar{\varsigma}}+t)}}+\frac{\sqrt{2}\varepsilon_{1}}{(\tilde{\bar{\varsigma}}+t)^{\varepsilon_{1}}}\right)\chi_{\lambda}\left(\frac{y-y_{*}}{\sqrt{\tilde{\bar{\varsigma}}+t}}\right)\right]\\
 &+1-\chi(y).
\end{align*}
In fact, $\tilde{\varphi}$ can be obtained from $\varphi$ by replacing $\bar{\varsigma}$ with $\tilde{\bar{\varsigma}}$. One can check that there exists positive constants $c_{\lambda}$ and $C_{\lambda}$ depending only $\lambda$ such that
$c_{\lambda}\varphi\leq\tilde{\varphi}\leq C_{\lambda}\varphi$.

Let
$\mathcal{U}_{s}(t,x,y)=\tilde{\mathcal{P}}\left(\partial_{y}\partial_{x}^{s}u-\tilde{\rho}\partial_{x}^{s}u\right)$,
from Theorem \ref{wpae} one has
\begin{align}\label{ue}
\|u\|_{L_{t}^{\infty}\mathcal{H}_{\psi,\tilde{\varphi}}^{s}(\Omega_{t^{*}})}
+\|\partial_{y}\mathcal{U}_{s}\|_{L_{t}^{2}L_{\tilde{\varphi}}^{2}(\Omega_{t^{*}})}
+\|\sqrt{\eta}\partial_{x}\mathcal{U}_{s}\|_{L_{t}^{2}L_{\tilde{\varphi}}^{2}(\Omega_{t^{*}})}
+\|\tilde{\mathcal{P}}\mathcal{U}_{s}\|_{L_{t}^{2}L_{\tilde{\varphi}}^{2}(\Omega_{t^{*}})}
\leq C_{*}\sqrt{\bar{\eta}}.
\end{align}
As a consequence of \eqref{lt} and \eqref{ue}, one has
\begin{align}\label{us-3}
\|u\|_{H_{\psi}^{s-3}(\Omega)}\leq C\check{c}+C_{*}\sqrt{\bar{\eta}}t,~\forall t\in[0,t^{*}].
\end{align}
Moreover, note that the approximate solution is the superposition of the corrected increment solutions, from \eqref{ue} we can obtain $\|\mathbb{U}\|_{L_{t}^{\infty}\hat{H}_{\psi}^{s}(\Omega_{t^{*}})}\leq C_{*}$, where $C_{*}$ is a positive constant depending only on $\lambda$ and $\mathcal{Z}$ as before.

~~~~~~~~~~

\subsection{Low-order estimates of the approximate solutions}~

The main purpose of this section is to establish low-order estimates of $\mathfrak{U}$ by using the maximum principle. These estimates give us the low-order estimates of $\rho$, which are vital to the uniform estimate of corrected increment solutions and approximate solutions.

The main assumption we need in this section is as follows.
\begin{assumption}\label{FA}
Let $\hat{F}$ be given  in \eqref{hatF}. For some $\hat{t}_{*}\in(0,T]$ and any multiindex $\gamma=(\gamma_{1},\gamma_{2})\in\Gamma_{9}$, assume $\hat{F}$ satisfies
\begin{align}\label{Fd}
|D^{\gamma}\hat{F}|\leq C_{*}e^{-\frac{4}{3}\delta\y},~\text{in}~\Omega_{\hat{t}^{*}}.
\end{align}
\end{assumption}

By using Assumption \ref{FA},  Assumption \ref{inA'} and \ref{inA}, we can establish the low-order estimates of $\mathfrak{U}$ and $\rho$ below.
\begin{lemma}\label{DI}
Assume the initial date $\tilde{u}_{in}$ in the problem \eqref{App+1} satisfying Assumption \ref{inA}.
There exists a positive constant $\theta_{*}$ such that for any $\theta\geq\theta_{*}$, there corresponds a constant $\bar{t}_{0}\in(0,t^{*}]$ such that for any multiindex $\gamma\in\Gamma_{6}$ and any $t_{v}\in(0,\bar{t}_{0}]$,
if $\mathfrak{U}$ is a solution to the problem \eqref{App+1} on the time interval $[0,t_{v}]$ such that
$\|\mathbb{U}\|_{L_{t}^{\infty}\hat{H}_{\psi}^{s}(\Omega_{t_{v}})}\leq C_{*}$,
and Assumption \ref{FA} holds for $t_{v}$, then
\begin{align}\label{Di}
|D^{\gamma}\mathbb{U}|\leq C_{*}\y^{-1}e^{-(1-\theta t)\delta\y}, ~\text{in}~\Omega_{t_{v}}.
\end{align}
\end{lemma}
\begin{proof}
This lemma will be proved in two step. First, we obtain a rough decay estimate of $D^{\gamma}\mathbb{U}$ for any $\gamma\in\Gamma_{7}$ by using energy method. Then, with the help of the rough decay estimate, we prove \eqref{Di} by using the maximal principle for parabolic equation.

\textbf{Step 1.} Let $\psi_{\theta}(t,y)=\y^{-1}e^{(1-\theta t)\delta\y}$ for $\theta$ to be determined and $\bar{t}_{0}=\min\{t^{*},\frac{1}{2\theta}\}$. Assume $\mathfrak{U}$ is a solution to the problem \eqref{App+1} on the time interval $[0,\bar{t}_{0}]$ such that $\|\mathbb{U}\|_{L_{t}^{\infty}\hat{H}_{\psi}^{s}(\Omega_{\bar{t}_{0}})}\leq C_{*}$,
and the source term $\hat{F}$ satisfying Assumption \ref{FA} for $\bar{t}_{0}$.
From \eqref{App+1} one knows that $\mathbb{U}$ obeys the following problem in $\Omega_{\bar{t}_{0}}$,
\begin{equation}\label{mathbbU}
\begin{cases}
\partial_{t}\mathbb{U}
+\mathfrak{U}\partial_{x}\mathbb{U}
-\partial_{y}^{-1}[\partial_{x}\mathfrak{U}]\partial_{y}\mathbb{U}
+\partial_{x}u_{0}\mathbb{U}
=\partial_{y}^{2}\mathbb{U}+\eta\partial_{x}^{2}\mathbb{U}+\mathfrak{R}_{0},\\
\mathbb{U}|_{y=0}=0,~\lim\limits_{y\to+\infty}\mathbb{U}=0,\\
\mathbb{U}|_{t=0}=\tilde{u}_{in}(x,y),
\end{cases}
\end{equation}
where
\begin{align*}
\mathfrak{R}_{0}
=\partial_{y}^{2}u_{0}-\partial_{t}u_{0}-u_{0}\partial_{x}u_{0}+\partial_{y}^{-1}[\partial_{x}\mathfrak{U}]\partial_{y}u_{0}-\partial_{x}P-\hat{F}.
\end{align*}
For any multiindex $\gamma\in\Gamma_{9}$, by acting $D^{\gamma}$ on \eqref{mathbbU}$_{1}$ and multiplying the resulting equation by $\psi_{\theta}^{2}D^{\gamma}\mathbb{U}$ then integrating over $\Omega$, one gets the following equality for any  $t\in[0,\bar{t}_{0}]$,
\begin{align}\label{DmathbbUE}
&\frac{1}{2}\frac{d}{dt}\|D^{\gamma}\mathbb{U}\|_{L_{\psi_{\theta}}^{2}(\Omega)}^{2}
+\theta\delta\|\sqrt{\y}D^{\gamma}\mathbb{U}\|_{L_{\psi_{\theta}}^{2}(\Omega)}^{2}
+\|\partial_{y}D^{\gamma}\mathbb{U}\|_{L_{\psi_{\theta}}^{2}(\Omega)}^{2}
+\|\sqrt{\eta}\partial_{x}D^{\gamma}\mathbb{U}\|_{L_{\psi_{\theta}}^{2}(\Omega)}^{2}\\
=&-(D^{\gamma}[\mathfrak{U}\partial_{x}\mathbb{U}],\psi_{\theta}^{2}D^{\gamma}\mathbb{U})
  +(D^{\gamma}[\partial_{y}^{-1}[\partial_{x}\mathfrak{U}]\partial_{y}\mathbb{U}],\psi_{\theta}^{2}D^{\gamma}\mathbb{U})
  +(D^{\gamma}[\partial_{x}u_{0}\mathbb{U}],\psi_{\theta}^{2}D^{\gamma}\mathbb{U})\nonumber\\
 &+(D^{\gamma}[\eta\partial_{x}^{2}\mathbb{U}]-\eta\partial_{x}^{2}D^{\gamma}\mathbb{U},\psi_{\theta}^{2}D^{\gamma}\mathbb{U})
  +(D^{\gamma}\mathfrak{R}_{0},\psi_{\theta}^{2}D^{\gamma}\mathbb{U})
  -2(\partial_{y}D^{\gamma}U,\psi_{\theta}\partial_{y}\psi_{\theta}D^{\gamma}\mathbb{U})\nonumber\\
 &-(\partial_{y}D^{\gamma}U,\psi_{\theta}^{2}D^{\gamma}\mathbb{U})_{L_{x}^{2}}|_{y=0}:=\sum_{i=1}^{7}\mathfrak{B}_{i}.\nonumber
\end{align}
Now we shall estimate $\mathfrak{B}_{1},...,\mathfrak{B}_{7}$ one by one.

\noindent{\bf \underline{Estimate of~$\mathfrak{B}_{1}$.}}
It is clear to see that
\begin{align*}
\mathfrak{B}_{1}
=\sum_{0\leq\beta\leq\gamma}\binom{\gamma}{\beta}(D^{\beta}\mathfrak{U}D^{\gamma-\beta}[\partial_{x}\mathbb{U}],\psi_{\theta}^{2}D^{\gamma}\mathbb{U}).
\end{align*}
If $\beta=0$, one can use integration by parts to obtain
\begin{align*}
(D^{\beta}\mathfrak{U}D^{\gamma-\beta}[\partial_{x}\mathbb{U}],\psi_{\theta}^{2}D^{\gamma}\mathbb{U})
=-\frac{1}{2}(\partial_{x}\mathfrak{U}D^{\gamma}\mathbb{U},\psi_{\theta}^{2}D^{\gamma}\mathbb{U})
\leq C_{*}\|\mathbb{U}\|_{H_{\psi_{\theta}}^{9}(\Omega)}^{2}.
\end{align*}
else, one has
\begin{align*}
(D^{\beta}\mathfrak{U}D^{\gamma-\beta}[\partial_{x}\mathbb{U}],\psi_{\theta}^{2}D^{\gamma}\mathbb{U})
\leq C_{*}\|D^{\gamma-\beta}[\partial_{x}\mathbb{U}]\|_{L_{\psi_{\theta}}^{2}(\Omega)}\|D^{\gamma}\mathbb{U}\|_{L_{\psi_{\theta}}^{2}(\Omega)}
\leq C_{*}\|\mathbb{U}\|_{H_{\psi_{\theta}}^{9}(\Omega)}^{2}.
\end{align*}
Thus,
\begin{align*}
\mathfrak{B}_{1}
\leq C_{*}\|\mathbb{U}\|_{H_{\psi_{\theta}}^{9}(\Omega)}^{2}.
\end{align*}

\noindent{\bf \underline{Estimate of~$\mathfrak{B}_{2}$.}}
It is clear to see that
\begin{align*}
\mathfrak{B}_{2}
=\sum_{0\leq\beta\leq\gamma}\binom{\gamma}{\beta}(D^{\beta}\partial_{y}^{-1}[\partial_{x}\mathfrak{U}]D^{\gamma-\beta}[\partial_{y}\mathbb{U}],\psi_{\theta}^{2}D^{\gamma}\mathbb{U}).
\end{align*}
If $\beta=0$, one can use integration by parts to obtain
\begin{align*}
(D^{\beta}\partial_{y}^{-1}[\partial_{x}\mathfrak{U}]D^{\gamma-\beta}[\partial_{y}\mathbb{U}],\psi_{\theta}^{2}D^{\gamma}\mathbb{U})
=&-\frac{1}{2}(\partial_{x}\mathfrak{U}D^{\gamma}\mathbb{U},\psi_{\theta}^{2}D^{\gamma}\mathbb{U})
 +(\partial_{y}^{-1}[\partial_{x}\mathfrak{U}]D^{\gamma}\mathbb{U},\psi_{\theta}\partial_{y}\psi_{\theta}D^{\gamma}\mathbb{U})\\
\leq&~C_{*}\|\mathbb{U}\|_{H_{\psi_{\theta}}^{9}(\Omega)}^{2}
     +C_{*}\|\sqrt{\y}D^{\gamma}\mathbb{U}\|_{L_{\psi_{\theta}}^{2}(\Omega)}^{2}.
\end{align*}
else, one has
\begin{align*}
(D^{\beta}\partial_{y}^{-1}[\partial_{x}\mathfrak{U}]D^{\gamma-\beta}[\partial_{y}\mathbb{U}],\psi_{\theta}^{2}D^{\gamma}\mathbb{U})
\leq C_{*}\|\sqrt{\y}D^{\gamma-\beta}[\partial_{y}\mathbb{U}]\|_{L_{\psi_{\theta}}^{2}(\Omega)}
          \|\sqrt{\y}D^{\gamma}\mathbb{U}\|_{L_{\psi_{\theta}}^{2}(\Omega)}^{2}
\leq C_{*}\|\mathbb{U}\|_{H_{\psi_{\theta},\y}^{9}(\Omega)}^{2}.
\end{align*}
To sum up,
\begin{align*}
\mathfrak{B}_{2}
\leq C_{*}\|\mathbb{U}\|_{H_{\psi_{\theta},\y}^{9}(\Omega)}^{2}.
\end{align*}

\noindent{\bf \underline{Estimate of~$\mathfrak{B}_{3}$.}}
\begin{align*}
\mathfrak{B}_{3}
=\sum_{0\leq\beta\leq\gamma}\binom{\gamma}{\beta}(D^{\beta}[\partial_{x}u_{0}]D^{\gamma-\beta}\mathbb{U},\psi_{\theta}^{2}D^{\gamma}\mathbb{U})
\leq C_{*}\|\mathbb{U}\|_{H_{\psi_{\theta}}^{9}(\Omega)}^{2}
\end{align*}

\noindent{\bf \underline{Estimate of~$\mathfrak{B}_{4}$.}}
By using integration by parts, one has
\begin{align*}
\mathfrak{B}_{4}
=&-\sum_{0<\beta\leq\gamma}\binom{\gamma}{\beta}
  (D^{\beta}\eta D^{\gamma-\beta}[\partial_{x}\mathbb{U}],\psi_{\theta}^{2}\partial_{x}D^{\gamma}\mathbb{U})\\
\leq&~C\sum_{0<\beta\leq\gamma}\|D^{\gamma-\beta}[\partial_{x}\mathbb{U}]\|_{L_{\psi_{\theta}}^{2}(\Omega)}
       \|\sqrt{\eta}\partial_{x}D^{\gamma}\mathbb{U}\|_{L_{\psi_{\theta}}^{2}(\Omega)}\\
\leq&~\frac{1}{2}\|\sqrt{\eta}\partial_{x}D^{\gamma}\mathbb{U}\|_{L_{\psi_{\theta}}^{2}(\Omega)}^{2}
     +C\|\mathbb{U}\|_{H_{\psi_{\theta}}^{9}(\Omega)}^{2}.
\end{align*}
\noindent{\bf \underline{Estimate of~$\mathfrak{B}_{5}$.}}
It follows from \eqref{B} that
$$\partial_{t}u_{0}+u_{0}\partial_{x}u_{0}+\partial_{x}P=\partial_{t}(u_{0}-U)+(u_{0}-U)\partial_{x}u_{0}+U\partial_{x}(u_{0}-U),$$
thus by using \eqref{Fd} and \eqref{phi} one has $\mathfrak{R}_{0}\in L_{t}^{\infty}H_{\psi_{\theta}}^{9}(\Omega_{\bar{t}_{0}})$, and
\begin{align*}
\mathfrak{B}_{5}
\leq \|D^{\gamma}\mathfrak{R}_{0}\|_{L_{\psi_{\theta}}^{2}(\Omega)}\|D^{\gamma}\mathbb{U}\|_{L_{\psi_{\theta}}^{2}(\Omega)}
\leq \|D^{\gamma}\mathfrak{R}_{0}\|_{L_{\psi_{\theta}}^{2}(\Omega)}^{2}+\|D^{\gamma}\mathbb{U}\|_{L_{\psi_{\theta}}^{2}(\Omega)}^{2}.
\end{align*}
\noindent{\bf \underline{Estimate of~$\mathfrak{B}_{6}$.}}
By using the Young inequality, it is easy to obtain
\begin{align*}
\mathfrak{B}_{6}
\leq \frac{1}{2}\|\partial_{y}D^{\gamma}\mathbb{U}\|_{L_{\psi_{\theta}}^{2}(\Omega)}^{2}
    +C\|D^{\gamma}\mathbb{U}\|_{L_{\psi_{\theta}}^{2}(\Omega)}^{2}.
\end{align*}
\noindent{\bf \underline{Estimate of~$\mathfrak{B}_{7}$.}}
Since $\|\mathbb{U}\|_{L_{t}^{\infty}\hat{H}_{\psi}^{s}(\Omega_{t^{*}})}\leq C_{*}$, by using Sobolev imbedding theorems, one has
\begin{align*}
\mathfrak{B}_{7}
\leq C\|\partial_{y}D^{\gamma}\mathbb{U}\|_{L^{\infty}(\Omega_{t^{*}})}\|D^{\gamma}\mathbb{U}\|_{L^{\infty}(\Omega_{t^{*}})}
\leq C_{*}.
\end{align*}

Collect all the estimate of $\mathfrak{B}_{1},...,\mathfrak{B}_{7}$ above, from \eqref{DmathbbUE} it follows that for any $t\in(0,\bar{t}_{0}]$,
\begin{align}\label{mathbbUe}
&\frac{1}{2}\frac{d}{dt}\|D^{\gamma}\mathbb{U}\|_{L_{\psi_{\theta}}^{2}(\Omega)}^{2}
+\theta\delta\|\sqrt{\y}D^{\gamma}\mathbb{U}\|_{L_{\psi_{\theta}}^{2}(\Omega)}^{2}
+\frac{1}{2}\|\partial_{y}D^{\gamma}\mathbb{U}\|_{L_{\psi_{\theta}}^{2}(\Omega)}^{2}
+\frac{1}{2}\|\sqrt{\eta}\partial_{x}D^{\gamma}\mathbb{U}\|_{L_{\psi_{\theta}}^{2}(\Omega)}^{2}\\
\leq&~C_{*}\|\mathbb{U}\|_{H_{\psi_{\theta}}^{9}(\Omega)}^{2}
     +C_{*}\|\mathbb{U}\|_{H_{\psi_{\theta},\y}^{9}(\Omega)}^{2}
     +\|D^{\gamma}\mathfrak{R}_{0}\|_{L_{\psi_{\theta}}^{2}(\Omega)}^{2}
     +C_{*}.\nonumber
\end{align}
Taking summation over all $\gamma\in\Gamma_{8}$, \eqref{mathbbUe} gives that for any $t\in(0,\bar{t}_{0}]$
\begin{align}\label{mathbbUE}
\frac{1}{2}\frac{d}{dt}\|\mathbb{U}\|_{H_{\psi_{\theta}}^{9}(\Omega)}^{2}
+\theta\|\mathbb{U}\|_{H_{\psi_{\theta},\y}^{9}(\Omega)}^{2}
\leq C_{*}\|\mathbb{U}\|_{H_{\psi_{\theta}}^{9}(\Omega)}^{2}
    +C_{*}\|\mathbb{U}\|_{H_{\psi_{\theta},\y}^{9}(\Omega)}^{2}
    +\|\mathfrak{R}_{0}\|_{H_{\psi_{\theta}}^{9}(\Omega)}^{2}
    +C_{*}.
\end{align}
By choosing $\theta$ large enough in \eqref{mathbbUE}, one can deduce that for any $t\in(0,\bar{t}_{0}]$,
\begin{align}\label{mathbbUE'}
\frac{d}{dt}\|\mathbb{U}\|_{H_{\psi_{\theta}}^{9}(\Omega)}^{2}
\leq C_{*}\|\mathbb{U}\|_{H_{\psi_{\theta}}^{9}(\Omega)}^{2}
    +C\|\mathfrak{R}_{0}\|_{H_{\psi_{\theta}}^{9}(\Omega)}^{2}
    +C_{*}.
\end{align}
Due to \eqref{inA3} and \eqref{phi}, one has
\begin{align*}
\|\tilde{u}_{in}\|_{H_{\psi_{\theta}^{0}}^{9}(\Omega)}^{2}\leq C,
\end{align*}
where $\psi_{\theta}^{0}=\y^{-1}e^{\delta\y}$. By using Grownwall inequality, from \eqref{mathbbUE'} one can deduce
\begin{align*}
\|\mathbb{U}\|_{L_{t}^{\infty}H_{\psi_{\theta}}^{9}(\Omega_{\bar{t}_{0}})}^{2}
\leq C_{*}\|\tilde{u}_{in}\|_{H_{\psi_{\theta}^{0}}^{9}(\Omega)}^{2}
    +C_{*}\|\mathfrak{R}_{0}\|_{L_{t}^{\infty}H_{\psi_{\theta}}^{9}(\Omega_{\bar{t}_{0}})}^{2}
    +C_{*}
\leq C_{*},
\end{align*}
which together with Sobolev imbedding theorems gives
\begin{align}\label{Di0}
|D^{\gamma}\mathbb{U}|\leq C_{*}\y e^{-(1-\theta t)\delta\y}~\text{in}~\Omega_{\bar{t}_{0}},
\end{align}
for any multiindex $\gamma\in\Gamma_{7}$.

\textbf{Step 2. }Now with the help of \eqref{Di0}, we are in a position to prove \eqref{Di}. The proof will be carried out by induction. First, we prove that \eqref{Di} for the case $|\gamma|=0$ by using maximum principle for parabolic equations, then under the induction hypothesis that \eqref{Di} holds for any $\gamma\in\Gamma_{k-1}$, we prove that \eqref{Di} holds for any $\gamma\in\Gamma_{k}, k\in\{1,2,...,6\}$ by induction.

To prove \eqref{Di} for the case $|\gamma|=0$, we define
$$B_{c}^{\pm}(t,x,y)=\pm e^{-M_{c}t}\mathbb{U}-N_{c}\y^{-1}e^{-(1-\theta t)\delta\y}-\nu,$$
where $M_{c}=\|\partial_{x}u_{0}\|_{L^{\infty}(\Omega_{t^{*}})}$, and $N_{c},\nu$ are positive constants which will be determined later.
From \eqref{mathbbU} one can see that $B_{c}^{\pm}$ obeys the following problem in $\Omega_{\bar{t}_{0}}$,
\begin{equation}\label{BE}
\begin{cases}
\partial_{t}B_{c}^{\pm}
+\mathfrak{U}\partial_{x}B_{c}^{\pm}
-\partial_{y}^{-1}[\partial_{x}\mathfrak{U}]\partial_{y}B_{c}^{\pm}
+(M_{c}+\partial_{x}u_{0})B_{c}^{\pm}
=\partial_{y}^{2}B_{c}^{\pm}+\eta\partial_{x}^{2}B_{c}^{\pm}\pm e^{-M_{c}t}\mathfrak{R}_{0}+\mathcal{R}_{Bc},\\
B_{c}^{\pm}|_{y=0}=-N_{c}e^{-(1-\theta t)\delta}-\nu,~\lim\limits_{y\to+\infty}B_{c}^{\pm}=-\nu,\\
B_{c}^{\pm}|_{t=0}=\pm\tilde{u}_{in}(x,y)-N_{c}\y^{-1}e^{-\delta\y}-\nu,
\end{cases}
\end{equation}
where
\begin{align*}
\mathcal{R}_{Bc}
=&-N_{c}\theta\delta e^{-(1-\theta t)\delta\y}
  -N_{c}(1-\theta t)\delta\y^{-1}\partial_{y}^{-1}[\partial_{x}\mathfrak{U}]e^{-(1-\theta t)\delta\y}
  -N_{c}\y^{-2}\partial_{y}^{-1}[\partial_{x}\mathfrak{U}]e^{-(1-\theta t)\delta\y}\\
 &-(M_{c}+\partial_{x}u_{0})(N_{c}\y^{-1}e^{-(1-\theta t)\delta\y}+\nu)
  +N_{c}\y^{-1}(1-\theta t)^{2}\delta^{2}e^{-(1-\theta t)\delta\y}\\
 &+2N_{c}\y^{-3}e^{-(1-\theta t)\delta\y}
  +N_{c}\y^{-2}(1-\theta t)\delta e^{-(1-\theta t)\delta\y}.
\end{align*}
By choosing $N_{c}$ large enough, one has $B_{c}^{\pm}|_{t=0}<0$ in $\Omega$. Noticing that $|\partial_{y}^{-1}[\partial_{x}\mathfrak{U}]|\leq C_{*}\y$,
one can take $\theta$ large enough to ensure $\mathcal{R}_{Bc}\pm e^{-M_{c}t}\mathfrak{R}_{0}<0$ in $\Omega_{\bar{t}_{0}}$.
Due to \eqref{Di0}, one has
$$|\pm e^{-M_{c}t}\mathbb{U}-N_{c}\y^{-1}e^{-(1-\theta t)\delta\y}|\leq C_{*}\y e^{-(1-\theta t)\delta\y}~\text{in}~\Omega_{\bar{t}_{0}},$$
thus for any fixed positive constant $\nu$, one knows $B_{c}^{\pm}<0$ for any $y>r_{c}=\frac{4}{\delta}\ln{\frac{C_{*}}{\nu\delta}}$ in $\Omega_{\bar{t}_{0}}$.
Applying the classical maximum principle for parabolic equations to the problem \eqref{BE}, it gives that
$B_{c}^{\pm}<0$ in $\hat{\Omega}_{\bar{t}_{0}}^{r_{c}}$.
Taking the limit $\nu\to 0^{+}$, one has $\lim\limits_{\nu\to 0^{+}}B_{c}^{\pm}\leq0 ~\text{in}~\Omega_{\bar{t}_{0}}$, which gives
\begin{align*}
|\mathbb{U}|\leq C_{*}\y^{-1}e^{-(1-\theta t)\delta\y} ~\text{in}~\Omega_{\bar{t}_{0}}.
\end{align*}

Assume that for $k\in\{1,2,...,6\}$ the inequality \eqref{Di} holds for any $\gamma\in\Gamma_{k-1}$. We shall prove that \eqref{Di} holds for any $\gamma\in\Gamma_{k}$. For this purpose, we define
$$B_{i}^{\pm}=\pm e^{-M_{i}t}\partial_{x}^{i}\partial_{y}^{k-i}\mathbb{U}-N_{i}\y^{-1}e^{-(1-\theta t)\delta\y}-\nu_{i},$$
for any $i\in\{0,1,2,...,k\}$, where $M_{i},N_{i},\nu_{i}$ are positive constants to be determined later.

We start from the case $i=0$. From \eqref{mathbbU} one can check that $B_{0}^{\pm}$ obeys the following problem in $\Omega_{\bar{t}_{0}}$,
\begin{equation}\label{B0}
\begin{cases}
\partial_{t}B_{0}^{\pm}
+\mathfrak{U}\partial_{x}B_{0}^{\pm}
-\partial_{y}^{-1}[\partial_{x}\mathfrak{U}]\partial_{y}B_{0}^{\pm}
+(M_{0}+\partial_{x}u_{0}-k\partial_{x}\mathfrak{U})B_{0}^{\pm}\\
\quad=\partial_{y}^{2}B_{0}^{\pm}+\eta\partial_{x}^{2}B_{0}^{\pm}\pm e^{-M_{0}t}\partial_{y}^{k}\mathfrak{R}_{0}+\mathcal{R}_{B0},\\
B_{0}^{\pm}|_{y=0}=\pm e^{-M_{0}t}\partial_{y}^{k}\mathbb{U}|_{y=0}-N_{0}e^{-(1-\theta t)\delta}-\nu_{0},
~\lim\limits_{y\to+\infty}B_{0}^{\pm}=-\nu_{0},\\
B_{0}^{\pm}|_{t=0}=\pm\partial_{y}^{k}\tilde{u}_{in}(x,y)-N_{0}\y^{-1}e^{-\delta\y}-\nu_{0},
\end{cases}
\end{equation}
where
\begin{align*}
\mathcal{R}_{B0}
=&-e^{-M_{0}t}\sum_{i=0}^{k-1}\binom{k}{i}\partial_{y}^{k-i}\mathfrak{U}\partial_{x}\partial_{y}^{i}\mathbb{U}
  +e^{-M_{0}t}\sum_{i=0}^{k-2}\binom{k}{i}\partial_{y}^{k-i-1}\partial_{x}\mathfrak{U}\partial_{y}^{i+1}\mathbb{U}\\
 &-e^{-M_{0}t}\sum_{i=0}^{k-1}\binom{k}{i}\partial_{y}^{k-i}\partial_{x}u_{0}\partial_{y}^{i}\mathbb{U}
  +e^{-M_{0}t}\sum_{i=0}^{k-1}\binom{k}{i}\partial_{y}^{k-i}\eta \partial_{y}^{i}\partial_{x}^{2}\mathbb{U}\\
 &-N_{0}\theta\delta e^{-(1-\theta t)\delta\y}
  -N_{0}(1-\theta t)\delta\y^{-1}\partial_{y}^{-1}[\partial_{x}\mathfrak{U}]e^{-(1-\theta t)\delta\y}
  -N_{0}\y^{-2}\partial_{y}^{-1}[\partial_{x}\mathfrak{U}]e^{-(1-\theta t)\delta\y}\\
 &-(M_{0}+\partial_{x}u_{0}-k\partial_{x}\mathfrak{U})(N_{0}\y^{-1}e^{-(1-\theta t)\delta\y}+\nu)
  +N_{0}\y^{-1}(1-\theta t)^{2}\delta^{2}e^{-(1-\theta t)\delta\y}\\
 &+2N_{0}\y^{-3}e^{-(1-\theta t)\delta\y}
  +N_{0}\y^{-2}(1-\theta t)\delta e^{-(1-\theta t)\delta\y}.
\end{align*}
Note that when the subscript exceeds the superscript in the second summation notation of $\mathcal{R}_{B0}$, the sum evaluates to zero by convention.

Take $M_{0}=\|\partial_{x}u_{0}\|_{L^{\infty}(\Omega_{\bar{t}_{0}})}+k\|\partial_{x}\mathfrak{U}\|_{L^{\infty}(\Omega_{\bar{t}_{0}})}$.
Since $\|\mathbb{U}\|_{L_{t}^{\infty}\hat{H}_{\psi}^{s}(\Omega_{t^{*}})}\leq C_{*}$ and \eqref{inA3}, one can choose $N_{0}$ large enough to obtain $B_{0}^{\pm}|_{y=0}<0$ in $\Omega_{\bar{t}_{0}}$ and $B_{0}^{\pm}|_{t=0}<0$ in $\Omega$. Recalling that $\mathfrak{U}=\mathbb{U}+u_{0}$, by using \eqref{Di0} one has
\begin{align}\label{B0e1}
|\partial_{y}^{k-i}\mathfrak{U}\partial_{x}\partial_{y}^{i}\mathbb{U}|
\leq C_{*}\y^{2}e^{-2(1-\theta t)\delta\y}
\leq C_{*}e^{-(1-\theta t)\delta\y}~\text{in}~\Omega_{\bar{t}_{0}},~\forall i\in[0,k-1]\cap\mathbb{Z}.
\end{align}
Also, \eqref{Di0} together with \eqref{theta} leads to
\begin{align}\label{B0e2}
\partial_{y}^{k-i}\eta \partial_{y}^{i}\partial_{x}^{2}\mathbb{U}
\leq C_{*}\theta t\y e^{-(1+\theta t)\delta\y}
\leq C_{*}e^{-(1-\theta t)\delta\y}~\text{in}~\Omega_{\bar{t}_{0}}, ~\forall i\in[0,k-1]\cap\mathbb{Z}.
\end{align}
By the induction hypothesis that the inequality \eqref{Di} holds for any $\gamma\in\Gamma_{k-1}$, one gets
\begin{align}\label{B0e3}
|\partial_{y}^{i}\mathbb{U}|\leq C_{*}\y^{-1}e^{-(1-\theta t)\delta\y}~\text{in}~\Omega_{\bar{t}_{0}},~\forall i\in[0,k-1]\cap\mathbb{Z}.
\end{align}
With the help of \eqref{B0e1}--\eqref{B0e3}, one can take $\theta$ large enough to obtain
 $\mathcal{R}_{B0}\pm e^{-M_{0}t}\partial_{y}^{k}\mathfrak{R}_{0}<0$ in $\Omega_{\bar{t}_{0}}$.
Again by using \eqref{Di0}, one has
$$|\pm e^{-M_{0}t}\partial_{y}^{k}\mathbb{U}-N_{0}\y^{-1}e^{-(1-\theta t)\delta\y}|
\leq C_{*}\y e^{-(1-\theta t)\delta\y}~\text{in}~\Omega_{\bar{t}_{0}},$$
thus for any fixed positive constant $\nu_{0}$, one knows $B_{0}^{\pm}<0$ for any $y>r_{c}=\frac{4}{\delta}\ln{\frac{C_{*}}{\nu_{0}\delta}}$ in $\Omega_{\bar{t}_{0}}$.
Applying the classical maximum principle for parabolic equations to the problem \eqref{B0}, it follows that
$B_{0}^{\pm}<0$ in $\hat{\Omega}_{\bar{t}_{0}}^{r_{c}}$.
Taking the limit $\nu_{0}\to 0^{+}$, one has $\lim\limits_{\nu_{0}\to 0^{+}}B_{0}^{\pm}\leq0 ~\text{in}~\Omega_{\bar{t}_{0}}$, which gives
\begin{align*}
|\partial_{y}^{k}\mathbb{U}|\leq C_{*}\y^{-1}e^{-(1-\theta t)\delta\y} ~\text{in}~\Omega_{\bar{t}_{0}}.
\end{align*}

Now assume that for $j\in[1,k]\cap\mathbb{Z}$, the following inequality holds for any $i\in[0,j-1]\cap\mathbb{Z}$,
\begin{align}\label{Dik}
|\partial_{x}^{i}\partial_{y}^{k-i}\mathbb{U}|\leq C_{*}\y^{-1}e^{-(1-\theta t)\delta\y}~\text{in}~\Omega_{\bar{t}_{0}}.
\end{align}
We shall prove that \eqref{Dik} holds for $i=j$.
To this end, define
$$B_{j}^{\pm}=\pm e^{-M_{j}t}\partial_{x}^{j}\partial_{y}^{k-j}\mathbb{U}-N_{j}\y^{-1}e^{-(1-\theta t)\delta\y}-\nu_{j},$$
where $M_{j},N_{j},\nu_{j}$ are positive constants to be determined later.
From \eqref{mathbbU} one can check that $B_{j}^{\pm}$ obeys the following problem in $\Omega_{\bar{t}_{0}}$,
\begin{equation}\label{Bj}
\begin{cases}
\partial_{t}B_{j}^{\pm}
+\mathfrak{U}\partial_{x}B_{j}^{\pm}
-\partial_{y}^{-1}[\partial_{x}\mathfrak{U}]\partial_{y}B_{j}^{\pm}
+(M_{j}+\partial_{x}u_{0}+(2j-k)\partial_{x}\mathfrak{U})B_{j}^{\pm}\\
\quad=\partial_{y}^{2}B_{j}^{\pm}+\eta\partial_{x}^{2}B_{j}^{\pm}+e^{-M_{j}t}\partial_{x}^{j}\partial_{y}^{k-j}\mathfrak{R}_{0}+\mathcal{R}_{Bj},\\
B_{j}^{\pm}|_{y=0}=\pm e^{-M_{j}t}\partial_{x}^{j}\partial_{y}^{k-j}\mathbb{U}|_{y=0}-N_{j}e^{-(1-\theta t)\delta}-\nu_{j},
~\lim\limits_{y\to+\infty}B_{j}^{\pm}=-\nu_{j},\\
B_{j}^{\pm}|_{t=0}=\pm\partial_{x}^{j}\partial_{y}^{k-j}\tilde{u}_{in}(x,y)-N_{j}\y^{-1}e^{-\delta\y}-\nu_{j},
\end{cases}
\end{equation}
where
\begin{align*}
\mathcal{R}_{Bj}
=&-e^{-M_{j}t}\left(\partial_{x}^{j}\partial_{y}^{k-j}[\mathfrak{U}\partial_{x}\mathbb{U}]
                    -\mathfrak{U}\partial_{x}^{j+1}\partial_{y}^{k-j}\mathbb{U}
                    -j\partial_{x}\mathfrak{U}\partial_{x}^{j}\partial_{y}^{k-j}\mathbb{U}\right)\\
 &+e^{-M_{j}t}\left(\partial_{x}^{j}\partial_{y}^{k-j}[\partial_{y}^{-1}[\partial_{x}\mathfrak{U}]\partial_{y}\mathbb{U}]
                    -\partial_{y}^{-1}[\partial_{x}\mathfrak{U}]\partial_{x}^{j}\partial_{y}^{k-j+1}\mathbb{U}
                    -(k-j)\partial_{x}\mathfrak{U}\partial_{x}^{j}\partial_{y}^{k-j}\mathbb{U}\right)\\
 &-e^{-M_{0}t}\left(\partial_{x}^{j}\partial_{y}^{k-j}[\partial_{x}u_{0}\mathbb{U}]-\partial_{x}u_{0}\partial_{x}^{j}\partial_{y}^{k-j}\mathbb{U}\right)
  +e^{-M_{0}t}\left(\partial_{x}^{j}\partial_{y}^{k-j}[\eta\partial_{x}^{2}\mathbb{U}]-\eta\partial_{x}^{j+2}\partial_{y}^{k-j}\mathbb{U}\right)\\
 &-N_{j}\theta\delta e^{-(1-\theta t)\delta\y}
  -N_{j}(1-\theta t)\delta\y^{-1}\partial_{y}^{-1}[\partial_{x}\mathfrak{U}]e^{-(1-\theta t)\delta\y}
  -N_{j}\y^{-2}\partial_{y}^{-1}[\partial_{x}\mathfrak{U}]e^{-(1-\theta t)\delta\y}\\
 &-(M_{j}+\partial_{x}u_{0}+(2j-k)\partial_{x}\mathfrak{U})(N_{j}e^{-(1-\theta t)\delta\y}+\nu)
  +N_{j}\y^{-1}(1-\theta t)^{2}\delta^{2}e^{-(1-\theta t)\delta\y}\nonumber\\
 &+2N_{j}\y^{-3}e^{-(1-\theta t)\delta\y}
  +N_{j}\y^{-2}(1-\theta t)\delta e^{-(1-\theta t)\delta\y}.
\end{align*}

Take $M_{j}=\|\partial_{x}u_{0}\|_{L^{\infty}(\Omega_{\bar{t}_{0}})}+|2j-k|\|\partial_{x}\mathfrak{U}\|_{L^{\infty}(\Omega_{\bar{t}_{0}})}$.
Under the conditions that $\|\mathbb{U}\|_{L_{t}^{\infty}\hat{H}_{\psi}^{s}(\Omega_{t^{*}})}\leq C_{*}$ and \eqref{inA3}, in addition, the induction hypothesis that the inequality \eqref{Di} holds for any $\gamma\in\Gamma_{k-1}$ and the inequality \eqref{Dik} holds for any $i\in[0,j-1]\cap\mathbb{Z}$, one can take $N_{j}$ and $\theta$ large enough to ensure $B_{j}^{\pm}|_{y=0}<0$ in $\Omega_{\bar{t}_{0}}$ and $B_{0}^{\pm}|_{t=0}<0$ in $\Omega$, and  $\mathcal{R}_{Bj}\pm e^{-M_{j}t}\partial_{x}^{j}\partial_{y}^{k-j}\mathfrak{R}_{0}<0$ in $\Omega_{\bar{t}_{0}}$.
By using \eqref{Di0}, one has
$$|\pm e^{-M_{j}t}\partial_{x}^{j}\partial_{y}^{k-j}\mathbb{U}-N_{j}\y^{-1}e^{-(1-\theta t)\delta\y}|
\leq C_{*}\y e^{-(1-\theta t)\delta\y}~\text{in}~\Omega_{\bar{t}_{0}},$$
thus for any fixed positive constant $\nu_{j}$, one knows $B_{j}^{\pm}<0$ for any $y>r_{c}=\frac{4}{\delta}\ln{\frac{C_{*}}{\nu_{j}\delta}}$ in $\Omega_{\bar{t}_{0}}$.
Applying the classical maximum principle for parabolic equations to the problem \eqref{B0}, it follows that
$B_{j}^{\pm}<0$ in $\hat{\Omega}_{\bar{t}_{0}}^{r_{c}}$.
Taking the limit $\nu_{j}\to 0^{+}$, one has $\lim\limits_{\nu_{j}\to 0^{+}}B_{0}^{\pm}\leq0 ~\text{in}~\Omega_{\bar{t}_{0}}$, which gives
\begin{align}\label{Dikj}
|\partial_{x}^{j}\partial_{y}^{k-j}\mathbb{U}|\leq C_{*}\y^{-1}e^{-(1-\theta t)\delta\y} ~\text{in}~\Omega_{\bar{t}_{0}}.
\end{align}
Combining \eqref{Dik} with \eqref{Dikj}, one knows that \eqref{Di} holds for any $\gamma\in\Gamma_{k}(k\in\{1,2,...,6\})$ and $t_{v}=\bar{t}_{0}$. To get \eqref{Di} for any $t_{v}\in[0,\bar{t}_{0}]$, one can check that it is sufficient to assume that $\|\mathbb{U}\|_{L_{t}^{\infty}\hat{H}_{\psi}^{s}(\Omega_{t_{v}})}\leq C_{*}$,
and  Assumption \ref{FA} holds for $t_{v}$. Therefore, the conclusion of Lemma \ref{DI} holds.

\end{proof}

\begin{lemma}\label{wdotu}
Let $\theta_{*}$ and $\bar{t}_{0}$ be obtained in Lemma \ref{DI}.
Assume the initial date $\tilde{u}_{in}$ in the problem \eqref{App+1} satisfies Assumption \ref{inA'} and \ref{inA}.
There exists a positive constant $\theta^{*}\geq \theta_{*}$ and positive constants $K_{1},K_{2}$ independent of $\mathcal{Z}$ such that for any $\theta\geq\theta^{*}$, there is a corresponding $t_{0}\in(0,\bar{t}_{0}]$ such that for any $t_{v}\in[0,t_{0}]$, if $\mathfrak{U}$ is a solution to the problem \eqref{App+1} on the time interval $[0,t_{v}]$ such that
$\|\mathbb{U}\|_{L_{t}^{\infty}\hat{H}_{\psi}^{s}(\Omega_{t_{v}})}\leq C_{*}$,
and Assumption \ref{FA} holds for $t_{v}$, then the following inequality holds in $\Omega_{t_{v}}$,
\begin{align}\label{dotu}
K_{1}\varpi e^{-(1+\theta t)\delta\y}\leq \dot{\mathfrak{U}}\leq K_{2}e^{-(1-\theta t)\delta\y}.
\end{align}
\end{lemma}
\begin{proof}
Assume $\mathfrak{U}$ is a solution to the problem \eqref{App+1} on the time interval $[0,\bar{t}_{0}]$ such that
$\|\mathbb{U}\|_{L_{t}^{\infty}\hat{H}_{\psi}^{s}(\Omega_{\bar{t}_{0}})}\leq C_{*}$,
and the source term $\hat{F}$ satisfying Assumption \ref{FA} for $\bar{t}_{0}$.
Restricting the equation \eqref{App+1}$_{1}$ to the boundary $y=0$, one has
\begin{align*}
\partial_{y}\mathfrak{U}(t,x,0)
=& \partial_{y}u_{in}(x,0)+\int_{0}^{t}\partial_{t}\partial_{y}\mathfrak{U}(\tau,x,0)d\tau\\
=& \partial_{y}u_{in}(x,0)+\int_{0}^{t}[\partial_{y}^{3}\mathfrak{U}(\tau,x,0)+\eta\partial_{x}^{2}\partial_{y}\mathbb{U}(\tau,x,0)]d\tau,
\end{align*}
Since $\|\mathbb{U}\|_{L_{t}^{\infty}\hat{H}_{\psi}^{s}(\Omega_{\bar{t}_{0}})}\leq C_{*}$, one can choose $\tilde{t}$ small enough such that
\begin{align}\label{dotu0}
2\max_{(t,x)\in[0,\tilde{t}]\times\mathbb{T}}\partial_{y}u_{in}(x,0)
\geq \partial_{y}\mathfrak{U}(t,x,0)
\geq \frac{1}{2}\min_{(t,x)\in[0,\tilde{t}]\times\mathbb{T}}\partial_{y}u_{in}(x,0), ~\text{in}~[0,\tilde{t}]\times\mathbb{T}.
\end{align}

To get \eqref{dotu}, we define
$$\mathcal{G}_{1}(t,x,y)=\dot{\mathfrak{U}}-\mu_{1}\varpi e^{-(1+\theta t)\delta\y}+\omega,$$
and
$$\mathcal{G}_{2}(t,x,y)=\dot{\mathfrak{U}}-\mu_{2}e^{-(1-\theta t)\delta\y}-\nu,$$
where $\mu_{1},\mu_{2},\theta,\nu$ and $\omega$ are positive constants which will be determined later, with $\theta\geq\theta_{*}$.
Let $t_{0}=\min\left\{\bar{t}_{0},\tilde{t},\frac{1}{4\theta}\right\}$.
One can check that $\mathcal{G}_{1}$ obey the following problem in $\Omega_{t_{0}}$.
\begin{equation}\label{pdotuge}
\begin{cases}
\partial_{t}\mathcal{G}_{1}
+\mathfrak{U}\partial_{x}\mathcal{G}_{1}
-\partial_{y}^{-1}[\partial_{x}\mathfrak{U}]\partial_{y}\mathcal{G}_{1}
-\partial_{y}^{2}\mathcal{G}_{1}
-\eta\partial_{x}^{2}\mathcal{G}_{1}
=\mathfrak{R}_{1},\\
\mathcal{G}_{1}|_{y=0}=\partial_{y}\mathfrak{U}|_{y=0}-\mu_{1}\varpi(0)e^{-(1+\theta t)\delta}+\omega,~
\lim\limits_{y\to+\infty}\mathcal{G}_{1}=\omega,\\
\mathcal{G}_{1}|_{t=0}=\partial_{y}u_{in}(x,y)-\mu_{1}\varpi e^{-\delta\y}+\omega,
\end{cases}
\end{equation}
where
\begin{align*}
\mathfrak{R}_{1}
=& \mu_{1}\varpi \theta\delta\y e^{-(1+\theta t)\delta\y}
  +\mu_{1}\partial_{y}^{-1}[\partial_{x}\mathfrak{U}]\left(\varpi'-\delta\varpi-\theta t\delta\varpi\right)e^{-(1+\theta t)\delta\y}
  +\mu_{1}\varpi''e^{-(1+\theta t)\delta\y}\\
 &-2\mu_{1}(1+\theta t)\delta\varpi'e^{-(1+\theta t)\delta\y}
  +\mu_{1}(1+\theta t)^{2}\delta^{2}\varpi e^{-(1+\theta t)\delta\y}
  -\partial_{y}\hat{F}
  -\eta\partial_{x}^{2}\partial_{y}u_{0}
  +\eta'\partial_{x}^{2}\mathbb{U},
\end{align*}
and $\mathcal{G}_{2}$ obey the following problem in $\Omega_{t_{0}}$.
\begin{equation}\label{pdotule}
\begin{cases}
\partial_{t}\mathcal{G}_{2}
+\mathfrak{U}\partial_{x}\mathcal{G}_{2}
-\partial_{y}^{-1}[\partial_{x}\mathfrak{U}]\partial_{y}\mathcal{G}_{2}
-\partial_{y}^{2}\mathcal{G}_{2}
-\eta\partial_{x}^{2}\mathcal{G}_{2}
=\mathfrak{R}_{2},\\
\mathcal{G}_{2}|_{y=0}=\partial_{y}\mathfrak{U}|_{y=0}-\mu_{2}e^{-(1-\theta t)\delta}-\nu,~
\lim\limits_{y\to+\infty}\mathcal{G}_{2}=-\nu,\\
\mathcal{G}_{2}|_{t=0}=\partial_{y}u_{in}(x,y)-\mu_{2}e^{-\delta\y}-\nu,
\end{cases}
\end{equation}
where
\begin{align*}
\mathfrak{R}_{2}
=&-\mu_{2}\theta\delta\y e^{-(\delta-N_{2}t)\y}
  -\mu_{2}\delta\partial_{y}^{-1}[\partial_{x}\mathfrak{U}]\left(1-\theta t\right) e^{-(1-\theta t)\delta\y}\\
 &+\mu_{2}(1-\theta t)^{2}\delta^{2}e^{-(1-\theta t)\delta\y}
  -\partial_{y}\hat{F}
  -\eta\partial_{x}^{2}\partial_{y}u_{0}
  +\eta'\partial_{x}^{2}\mathbb{U}.
\end{align*}

Since $\partial_{y}u_{in}(x,y)\geq c_{in}\varpi(y)e^{-\delta\y}$ in $\Omega$, from \eqref{dotu0} one has $\mathcal{G}_{1}|_{y=0}>0$ and $\mathcal{G}_{1}|_{t=0}>0$ in $\Omega_{t_{0}}$, by choosing $\mu_{1}$ small enough. Noticing the definition of $\varpi$, and
$$|\eta\partial_{x}^{2}\partial_{y}u_{0}|\leq Ce^{-(1+\theta t)\delta\y},
|\eta'\partial_{x}^{2}\mathbb{U}|\leq C_{*}e^{-(1+\theta t)\delta\y},~\text{in}~\Omega_{t_{0}},$$
one can take $\theta$ large enough to ensure $\mathfrak{R}_{1}>0$ in $\Omega_{t_{0}}$.
For any fixed positive constant $\omega$, since
$$|\dot{\mathfrak{U}}-\mu_{1}\varpi e^{-(1+\theta t)\delta\y}|\leq C_{*}e^{-\frac{2}{3}\delta\y}~\text{in}~\Omega_{t_{0}},$$
we know $\mathcal{G}_{1}>0$ for any $y>R_{0}=\frac{3}{2\delta}\ln{\frac{C_{*}}{\omega}}$ in $\Omega_{t_{0}}$.
Applying the classical maximum principle for parabolic equations to the problem \eqref{pdotuge}, it gives that
$\mathcal{G}_{1}>0$ in $\hat{\Omega}_{t_{0}}^{R_{0}}$.
Taking the limit $\omega\to 0^{+}$, one has $\lim\limits_{\omega\to 0^{+}}\mathcal{G}_{1}\geq 0 ~\text{in}~\Omega_{t_{0}}$, which gives the left inequality in \eqref{dotu} for $t_{v}=t_{0}$.

The proof of the right inequality in \eqref{dotu} is similar.
By choosing $\mu_{2}$ large enough, we can obtain that
$\mathcal{G}_{2}|_{y=0}<0$ and $\mathcal{G}_{2}|_{t=0}<0 ~\text{in}~\Omega_{t_{0}}.$
For any fixed positive constant $\nu$, since
$$|\dot{\mathfrak{U}}-\mu_{2}e^{-(1-\theta t)\delta\y}|\leq C_{*}e^{-\frac{2}{3}\delta\y},$$
we know that $\mathcal{G}_{2}<0$ for any $y>R_{1}=\frac{3}{2\delta}\ln{\frac{C_{*}}{\nu}}$ in $\Omega_{t_{0}}$.
By taking $\mu_{2}$ and $\theta$ large enough, one can has $\mathfrak{R}_{2}<0$.
By the classical maximum principle for parabolic equations one knows that $\mathcal{G}_{2}<0$ in $\hat{\Omega}_{t_{0}}^{R_{1}}.$
Taking the limit $\nu\to 0^{+}$, one has $\lim\limits_{\nu\to0^{+}}\mathcal{G}_{2}\leq0 ~\text{in}~\Omega_{t_{0}}$, which gives the right inequality in \eqref{dotu}  for $t_{v}=t_{0}$.

As can be seen from our proof of \eqref{dotu} for $t_{v}=t_{0}$ above,  to obtain \eqref{dotu} for any $t_{v}\in[0,t_{0}]$, it suffices to assume that $\|\mathbb{U}\|_{L_{t}^{\infty}\hat{H}_{\psi}^{s}(\Omega_{t_{v})}}\leq C_{*}$,
and Assumption \ref{FA} holds for $t_{v}$. Therefore, Lemma \ref{wdotu} holds.
\end{proof}

As a direct corollary of Lemma \ref{DI} and Lemma \ref{wdotu}, one has
\begin{corollary}\label{wdotuC'}
Under the assumption of Lemma \ref{wdotu}, for any multiindex $\gamma\in\Gamma_{6}$ and $t_{v}\in[0,t_{0}]$, one has
\begin{align*}
|\vartheta_{c}^{2}D^{\gamma}\mathbb{U}|
\leq C_{*}\y^{-1}e^{-(1+\theta t)\delta\y}
\leq C_{*}\y^{-1}\dot{\mathfrak{U}},~\text{in}~\check{\Omega}_{t_{v}}^{\hat{y}}.
\end{align*}
\end{corollary}

\begin{corollary}\label{wdotuC}
Under the assumption of Lemma \ref{wdotu}, for any $t_{v}\in[0,t_{0}]$ one has
\begin{align}\label{Umu}
\frac{K_{1}}{(1+\theta t)\delta}e^{-(1+\theta t)\delta\hat{y}}e^{-(1+\theta t)\delta\y}
\leq U-\mathfrak{U}
\leq \frac{K_{2}}{(1-\theta t)\delta}e^{-(1-\theta t)\delta\y}, ~\text{in}~\Omega_{t_{v}}.
\end{align}
\end{corollary}
\begin{proof}
By \eqref{dotu}, for any $t_{v}\in[0,t_{0}]$ one has
\begin{align*}
-K_{1}\varpi e^{-(1+\theta t)\delta\y}\geq \partial_{y}(U-\mathfrak{U})\geq -K_{2}e^{-(1-\theta t)\delta\y}, ~\text{in}~\Omega_{t_{v}}.
\end{align*}
which gives
\begin{align}\label{Umu1}
K_{1}\int_{y}^{+\infty}\varpi e^{-(1+\theta t)\delta(1+\mathfrak{y})}d\mathfrak{y}
\leq U-\mathfrak{U}
\leq K_{2}\int_{y}^{+\infty}e^{-(1-\theta t)\delta(1+\mathfrak{y})}d\mathfrak{y}, ~\text{in}~\Omega_{t_{v}}.
\end{align}
Noticing the definition of $\varpi$, for any $y>\hat{y}$, it follows from \eqref{Umu1} that
\begin{align}\label{Umu2}
\frac{K_{1}}{(1+\theta t)\delta}e^{-(1+\theta t)\delta\y}\leq U-\mathfrak{U}\leq \frac{K_{2}}{(1-\theta t)\delta}e^{-(1-\theta t)\delta\y}, ~\text{in}~\Omega_{t_{v}},
\end{align}
while for $y\in[0,\hat{y}]$ one has,
\begin{align}\label{Umu3}
\frac{K_{1}}{(1+\theta t)\delta}e^{-(1+\theta t)\delta\langle\hat{y}\rangle}
\leq U-\mathfrak{U}
\leq \frac{K_{2}}{(1-\theta t)\delta}e^{-(1-\theta t)\delta\y}, ~\text{in}~\Omega_{t_{v}}.
\end{align}
Since for any $y\in[0,\hat{y}]$,
\begin{align*}
e^{-(1+\theta t)\delta\langle\hat{y}\rangle}
= e^{(1+\theta t)\delta(y-\hat{y})}e^{-(1+\theta t)\delta\y}
\geq e^{-(1+\theta t)\delta\hat{y}}e^{-(1+\theta t)\delta\y},
\end{align*}
from \eqref{Umu3} one has
\begin{align}\label{Umu3'}
\frac{K_{1}}{(1+\theta t)\delta}e^{-(1+\theta t)\delta\hat{y}}e^{-(1+\theta t)\delta\y}
\leq U-\mathfrak{U}
\leq \frac{K_{2}}{(1-\theta t)\delta}e^{-(1-\theta t)\delta\y}, ~\text{in}~\Omega_{t_{v}}.
\end{align}
Combining \eqref{Umu2} with \eqref{Umu3'}, we get the conclusion.
\end{proof}

\begin{lemma}\label{sdotu}
Let $\Gamma_{*}=\{(m,k)|m,k\in\mathbb{N}, 0\leq m\leq 2, 0\leq m+k\leq3\}$.
Let $\theta^{*}$ and $t_{0}$ be obtained in Lemma \ref{wdotu}.
Assume the initial date $\tilde{u}_{in}$ in the problem \eqref{App+1} satisfies Assumption \ref{inA'} and \ref{inA}.
There exists positive constants $K_{mk}$ independent of $\mathcal{Z}$ such that for any $\theta\geq\theta^{*}$, there is a corresponding $t_{1}\in(0,t_{0}]$ such that for any $t_{v}\in[0,t_{1}]$, if $\mathfrak{U}$ is a solution to the problem \eqref{App+1} on the time interval $[0,t_{v}]$ such that
$\|\mathbb{U}\|_{L_{t}^{\infty}\hat{H}_{\psi}^{s}(\Omega_{t_{v}})}\leq C_{*}$,
and Assumption \ref{FA} holds for $t_{v}$, then the following inequalities hold in $\Omega_{t_{v}}$,
\begin{align}\label{dotul}
|\dot{\mathfrak{U}}|\leq K_{00}(U-\mathfrak{U}),~
|\partial_{y}\dot{\mathfrak{U}}|\leq K_{01}(U-\mathfrak{U}),~
|\partial_{y}^{2}\dot{\mathfrak{U}}|\leq K_{02}(U-\mathfrak{U}),~
|\partial_{y}^{3}\dot{\mathfrak{U}}|\leq K_{03}(U-\mathfrak{U}),
\end{align}
\begin{align}\label{dotux1}
|\partial_{x}\dot{\mathfrak{U}}|\leq K_{10}(1+\y t)(U-\mathfrak{U}),~
|\partial_{x}\partial_{y}\dot{\mathfrak{U}}|\leq K_{11}(1+\y t)(U-\mathfrak{U}),~
|\partial_{x}\partial_{y}^{2}\dot{\mathfrak{U}}|\leq K_{12}(1+\y t)(U-\mathfrak{U})
\end{align}
\begin{align}\label{dotux2}
|\partial_{x}^{2}\dot{\mathfrak{U}}|\leq K_{20}(1+\y t)(U-\mathfrak{U}),~
|\partial_{x}^{2}\partial_{y}\dot{\mathfrak{U}}|\leq K_{21}(1+\y t)(U-\mathfrak{U}),
\end{align}
in addition, there exists a positive constant $\bar{\iota}$ such that if
\begin{align}\label{Fe}
|\partial_{y}\hat{F}|\leq \bar{\iota}+C|y-y_{*}|+C_{*}t~\text{in}~\Omega_{t_{v}}^{*},
\end{align}
then there are positive constants $k_{0},k_{1},k_{2},k_{3}$ and $k_{4}$ independent of $\mathcal{Z}$ such that
\begin{align}\label{dotug}
\dot{\mathfrak{U}}\geq k_{0}\varpi(U-\mathfrak{U}),~\text{in}~\Omega_{t_{v}},
\end{align}
and
\begin{align}\label{dotub}
k_{1}\varpi+k_{2}t\leq \dot{\mathfrak{U}}\leq k_{3}\varpi+k_{4}t,~\text{in}~\hat{\Omega}_{t_{v}}\cap\check{\Omega}_{t_{v}}.
\end{align}
\end{lemma}

\begin{proof}
Here we only prove that for any $\theta\geq \theta^{*}$, there exist $t_{1}\in(0,t_{0}]$ such that \eqref{dotul}--\eqref{dotux2} and \eqref{dotug}--\eqref{dotub} hold for $t_{v}=t_{1}$, under the conditions that $\|\mathbb{U}\|_{L_{t}^{\infty}\hat{H}_{\psi}^{s}(\Omega_{t_{0}})}\leq C_{*}$, and Assumption \ref{FA} holds for $t_{0}$, in addition, \eqref{Fe} holds for $t_{v}=t_{0}$. One can use the same argument as in Lemma \ref{wdotu} to obtain the conclusion of Lemma \ref{sdotu}.

For $(m,k)\in\Gamma_{*}$, let
    $$G_{0k}(t,x,y)=e^{-\lambda_{0k}t}\left[\partial_{y}^{k}\dot{\mathfrak{U}}-\mu(1+Mt)(U-\mathfrak{U})\right]-\nu,$$
and
$$G_{mk}(t,x,y)=e^{-\lambda_{mk}t}\left[\partial_{x}^{m}\partial_{y}^{k}\dot{\mathfrak{U}}-\mu(1+M\y t)(U-\mathfrak{U})\right]-\nu,~m=1~\text{or}~2,$$
parallel,
    $$G_{0k}^{-}(t,x,y)=e^{-\lambda_{0k}t}\left[-\partial_{y}^{k}\dot{\mathfrak{U}}-\mu(1+Mt)(U-\mathfrak{U})\right]-\nu,$$

$$G_{mk}^{-}(t,x,y)=e^{-\lambda_{mk}t}\left[-\partial_{x}^{m}\partial_{y}^{k}\dot{\mathfrak{U}}-\mu(1+M\y t)(U-\mathfrak{U})\right]-\nu,~m=1~\text{or}~2,$$
where $$\lambda_{mk}=|m-k|\|\partial_{x}\mathfrak{U}\|_{L^{\infty}(\Omega_{t_{0}})},$$
and $M,\mu,\nu$ are positive constants which will be determined later.
Since
\begin{align*}
\partial_{x}^{m}\partial_{y}^{k}\dot{\mathfrak{U}}(t,x,0)
=& \partial_{x}^{m}\partial_{y}^{k+1}u_{in}(x,0)
  +\int_{0}^{t}\partial_{x}^{m}\partial_{y}^{k}\partial_{t}\dot{\mathfrak{U}}(\tau,x,0)d\tau\\
=& \partial_{x}^{m}\partial_{y}^{k+1}u_{in}(x,0)
  +\int_{0}^{t}\partial_{x}^{m}\partial_{y}^{k}[\partial_{y}^{2}\dot{\mathfrak{U}}(\tau,x,0)+\eta\partial_{x}^{2}\partial_{y}\mathbb{U}(\tau,x,0)-\partial_{y}\mathbb{F}]d\tau\\
 &-\int_{0}^{t}\partial_{x}^{m}\partial_{y}^{k}[\mathfrak{U}\partial_{x}\dot{\mathfrak{U}}-\partial_{y}^{-1}[\partial_{x}\mathfrak{U}]\partial_{y}\dot{\mathfrak{U}}](\tau,x,0)d\tau,
\end{align*}
one can choose $\tilde{t}$ small enough such that for any $(m,k)\in\Gamma_{*}$,
\begin{align}\label{boundary}
\max_{(t,x)\in[0,\tilde{t}]\times\mathbb{T}}|\partial_{x}^{m}\partial_{y}^{k+1}\mathfrak{U}(t,x,0)|
\leq \max_{(t,x)\in[0,\tilde{t}]\times\mathbb{T}}|\partial_{x}^{m}\partial_{y}^{k+1}u_{in}(x,0)|+1.
\end{align}
Let $\hat{t}=\min\{t_{0},\tilde{t},\frac{1}{M}\}$. For any $(m,k)\in\Gamma_{*}$, by Corollary \ref{wdotuC} and \eqref{boundary} one can choose $\mu$ large enough but independent of $\mathcal{Z}$ to gain
$$G_{mk}|_{y=0}<0, G_{mk}|_{t=0}<0 ~\text{in}~\Omega_{\hat{t}}.$$
For any fixed positive constant $\nu$, since
$$|\partial_{x}^{m}\partial_{y}^{k}\dot{\mathfrak{U}}-\mu(1+Mt)(U-\mathfrak{U})|\leq C_{*}e^{-\frac{2}{3}\delta\y},$$
one gets that $G_{mk}<0$ for any $y>R=\frac{3}{2\delta}\ln{\frac{C_{*}}{\nu}}$ in $\Omega_{\hat{t}}$.

\noindent{\bf \underline{Proof of \eqref{dotul}.}}
From \eqref{App+1} one can check that $G_{0k}(k=0,1,2,3)$ obeys the follows problem in $\Omega_{\hat{t}}$,
\begin{equation*}
\begin{cases}
\partial_{t}G_{0k}
+\mathfrak{U}\partial_{x}G_{0k}
-\partial_{y}^{-1}[\partial_{x}\mathfrak{U}]\partial_{y}G_{0k}
+(\lambda_{0k}-k\partial_{x}\mathfrak{U})G_{0k}
-\partial_{y}^{2}G_{0k}
-\eta\partial_{x}^{2}G_{0k}
=\mathfrak{R}_{0k},\\
G_{0k}|_{y=0}=e^{-\lambda_{0k}t}\left[\partial_{y}^{k+1}\mathfrak{U}|_{y=0}-\mu(1+Mt)U\right]-\nu,~
\lim\limits_{y\to+\infty}G_{0k}=-\nu,\\
G_{0k}|_{t=0}=\partial_{x}^{m}\partial_{y}u_{in}(x,y)-\mu(U(0,x)-u_{in}(x,y))-\nu.
\end{cases}
\end{equation*}
where
\begin{align*}
\mathfrak{R}_{00}
=\mathfrak{R}_{c0}
  -\partial_{y}\hat{F}
  -\eta\partial_{x}^{2}\partial_{y}u_{0}
  +\partial_{y}\eta\partial_{x}^{2}\mathbb{U},
\end{align*}
\begin{align*}
\mathfrak{R}_{01}
=& e^{-\lambda_{01}t}\mathfrak{R}_{c0}
  +\mu(1+Mt)e^{-\lambda_{01}t}\partial_{x}\mathfrak{U}(U-\mathfrak{U})
  -e^{-\lambda_{01}t}\dot{\mathfrak{U}}\partial_{x}\dot{\mathfrak{U}}
  -\nu(\lambda_{01}-\partial_{x}\mathfrak{U})\\
 &-e^{-\lambda_{01}t}\partial_{y}^{2}\hat{F}
  -\eta e^{-\lambda_{01}t}\partial_{x}^{2}\partial_{y}^{2}u_{0}
  +2\partial_{y}\eta e^{-\lambda_{01}t}\partial_{x}^{2}\partial_{y}\mathbb{U}
  +\partial_{y}^{2}\eta e^{-\lambda_{01}t}\partial_{x}^{2}\mathbb{U},
\end{align*}
\begin{align*}
\mathfrak{R}_{02}
=& e^{-\lambda_{02}t}\mathfrak{R}_{c0}
  +2\mu(1+Mt)e^{-\lambda_{02}t}\partial_{x}\mathfrak{U}(U-\mathfrak{U})
  -2e^{-\lambda_{02}t}\dot{\mathfrak{U}}\partial_{x}\partial_{y}\dot{\mathfrak{U}}\\
 &-\nu(\lambda_{02}-2\partial_{x}\mathfrak{U})
  -e^{-\lambda_{02}t}\partial_{y}^{3}\hat{F}
  -\eta e^{-\lambda_{02}t}\partial_{x}^{2}\partial_{y}^{3}u_{0}
  +3\partial_{y}\eta e^{-\lambda_{02}t}\partial_{x}^{2}\partial_{y}^{2}\mathbb{U}\\
 &+3\partial_{y}^{2}\eta e^{-\lambda_{02}t}\partial_{x}^{2}\partial_{y}\mathbb{U}
  +\partial_{y}^{3}\eta e^{-\lambda_{02}t}\partial_{x}^{2}\mathbb{U},
\end{align*}
and
\begin{align*}
\mathfrak{R}_{03}
=& e^{-\lambda_{03}t}\mathfrak{R}_{c0}
  +3\mu(1+Mt)e^{-\lambda_{03}t}\partial_{x}\mathfrak{U}(U-\mathfrak{U})
  +2e^{-\lambda_{03}t}\partial_{x}\dot{\mathfrak{U}}\partial_{y}^{2}\dot{\mathfrak{U}}
  -2e^{-\lambda_{03}t}\partial_{y}\dot{\mathfrak{U}}\partial_{x}\partial_{y}\dot{\mathfrak{U}}\\
 &-3e^{-\lambda_{03}t}\dot{\mathfrak{U}}\partial_{x}\partial_{y}^{2}\dot{\mathfrak{U}}
  -\nu(\lambda_{03}-3\partial_{x}\mathfrak{U})
  -e^{-\lambda_{03}t}\partial_{y}^{4}\hat{F}
  -e^{-\lambda_{03}t}\eta\partial_{x}^{2}\partial_{y}^{4}u_{0}
  +4\partial_{y}\eta e^{-\lambda_{03}t}\partial_{x}^{2}\partial_{y}^{3}\mathbb{U}\\
 &+6\partial_{y}^{2}\eta e^{-\lambda_{03}t}\partial_{x}^{2}\partial_{y}^{2}\mathbb{U}
  +6\partial_{y}^{3}\eta e^{-\lambda_{03}t}\partial_{x}^{2}\partial_{y}\mathbb{U}
  +\partial_{y}^{4}\eta e^{-\lambda_{03}t}\partial_{x}^{2}\mathbb{U},
\end{align*}
with
\begin{align*}
\mathfrak{R}_{c0}
= -\mu M(U-\mathfrak{U})
  +\mu\eta(1+Mt)\partial_{x}^{2}(U-u_{0})
  +\mu(1+Mt)\partial_{x}U(U-\mathfrak{U})
  -\mu(1+Mt)\hat{F}.
\end{align*}
Notice that for any $k\in\{1,2,3,4\}$, one has
\begin{align}\label{u-0}
|\eta\partial_{x}^{2}\partial_{y}^{k}u_{0}|\leq C_{*}e^{-(1+\theta t)\delta\y}\leq C_{*}(U-\mathfrak{U}),~\text{in}~\Omega_{\hat{t}},
\end{align}
and for any $(i,j)\in\{(i,j)|i,j\in\mathbb{N},i+j\leq4\}$,
\begin{align}\label{U}
|\partial_{y}^{i}\eta\partial_{x}^{2}\partial_{y}^{j}\mathbb{U}|\leq C_{*}e^{-(1+\theta t)\delta\y}
\leq C_{*}(U-\mathfrak{U}),~\text{in}~\Omega_{\hat{t}},
\end{align}
moreover, by \eqref{Fd}, for any $(m,k)\in\Gamma_{*}$ one has
\begin{align}\label{F}
|\partial_{x}^{m}\partial_{y}^{k+1}\hat{F}|
\leq C_{*}e^{-\frac{4}{3}\delta\y}
\leq C_{*}(U-\mathfrak{U}),~\text{in}~\Omega_{\hat{t}}.
\end{align}
In virtue of Corollary \ref{wdotuC}, using \eqref{u-0},\eqref{U} and \eqref{F} one can take $M$ large enough and $\mu$ large enough but independent of $\mathcal{Z}$ to get $\mathfrak{R}_{00}<0$. Also by  using \eqref{u-0},\eqref{U} and \eqref{F}, noticing taht
$$|\dot{\mathfrak{U}}\partial_{x}\dot{\mathfrak{U}}|\leq C_{*}e^{-\frac{4}{3}\delta\y}\leq C_{*}(U-\mathfrak{U}),$$
one can take $M$ large enough and $\mu$ large enough independent of $\mathcal{Z}$ to obtain $\mathfrak{R}_{01}<0$.

By the classical maximum principle for parabolic equations one knows that
$$G_{00}<0,G_{01}<0 ~\text{in}~\hat{\Omega}_{\hat{t}}^{R}.$$
Taking the limit $\nu\to 0^{+}$, one has
$$\lim\limits_{\nu\to 0^{+}}G_{00}\leq 0, \lim\limits_{\nu\to 0^{+}}G_{01}\leq 0 ~\text{in}~\Omega_{\hat{t}},$$
which gives that
\begin{align}\label{dotu01l}
\dot{\mathfrak{U}}\leq \mu_{00}(U-\mathfrak{U}), \partial_{y}\dot{\mathfrak{U}}\leq \mu_{01}(U-\mathfrak{U})~\text{in}~\Omega_{\hat{t}},
\end{align}
for some positives constants $\mu_{00},\mu_{01}$ independent of $\mathcal{Z}$.

Since
\begin{align*}
&|\dot{\mathfrak{U}}\partial_{x}\partial_{y}\dot{\mathfrak{U}}|\leq C_{*}e^{-\frac{4}{3}\delta\y}\leq C_{*}(U-\mathfrak{U}),~
 |\partial_{x}\dot{\mathfrak{U}}\partial_{y}^{2}\dot{\mathfrak{U}}|\leq C_{*}e^{-\frac{4}{3}\delta\y}\leq C_{*}(U-\mathfrak{U}), ~\text{in}~\Omega_{\hat{t}},\\
&|\partial_{y}\dot{\mathfrak{U}}\partial_{x}\partial_{y}\dot{\mathfrak{U}}|\leq C_{*}e^{-\frac{4}{3}\delta\y}\leq C_{*}(U-\mathfrak{U}),~
 |\dot{\mathfrak{U}}\partial_{x}\partial_{y}^{2}\dot{\mathfrak{U}}|\leq C_{*}e^{-\frac{4}{3}\delta\y}
 \leq C_{*}(U-\mathfrak{U}),~\text{in}~\Omega_{\hat{t}},
\end{align*}
with the help of \eqref{u-0} and \eqref{U}, one can use the same argument as above to obtain
\begin{align}\label{dotu23l}
\partial_{y}^{2}\dot{\mathfrak{U}}\leq \mu_{02}(U-\mathfrak{U}),~
\partial_{y}^{3}\dot{\mathfrak{U}}\leq \mu_{03}(U-\mathfrak{U})~\text{in}~\Omega_{\hat{t}},
\end{align}
for some positive constants $\mu_{02},\mu_{03}$ independent of $\mathcal{Z}$, by taking $M$ large enough depending on $\mathcal{Z}$ but $\mu$ large independent of $\mathcal{Z}$.

To get \eqref{dotul}, it remain to prove
\begin{align}\label{-dotu}
&-\dot{\mathfrak{U}}\leq \bar{\mu}_{00}(U-\mathfrak{U}),
 -\partial_{y}\dot{\mathfrak{U}}\leq \bar{\mu}_{01}(U-\mathfrak{U}),~\text{in}~\Omega_{\hat{t}},\\
&-\partial_{y}^{2}\dot{\mathfrak{U}}\leq \bar{\mu_{02}}(U-\mathfrak{U}),
 -\partial_{y}^{3}\dot{\mathfrak{U}}\leq \bar{\mu_{03}}(U-\mathfrak{U}),~\text{in}~\Omega_{\hat{t}},\nonumber
\end{align}
for some positive constants $\bar{\mu}_{00},\bar{\mu}_{01},\bar{\mu}_{02},\bar{\mu}_{03}$ independent of $\mathcal{Z}$. Parallel to \eqref{dotu01l} and \eqref{dotu23l}, \eqref{-dotu} can be obtained by applying the classical maximum principle for the parabolic equations to the problems which are obeyed by $G_{0k}^{-}, k=0,1,2,3$. The proof is quite similar to above, so we omit it here.

Combining \eqref{dotu01l} and \eqref{dotu23l} with \eqref{-dotu}, one can get \eqref{dotul}.

\noindent{\bf \underline{Proof of \eqref{dotux1}.}}
From \eqref{App} one knows that $G_{1k}(k=0,1,2)$ obeys the follows problem in $\Omega_{\hat{t}}$,
\begin{equation*}
\begin{cases}
\partial_{t}G_{1k}
+\mathfrak{U}\partial_{x}G_{1k}
-\partial_{y}^{-1}[\partial_{x}\mathfrak{U}]\partial_{y}G_{1k}
+[\lambda_{1k}+(1-k)\partial_{x}\mathfrak{U}]G_{1k}
-\partial_{y}^{2}G_{1k}
-\eta\partial_{x}^{2}G_{1k}
=\mathfrak{R}_{1k},\\
G_{1k}|_{y=0}=e^{-\lambda_{1k}t}\left[\partial_{x}\partial_{y}^{k+1}\mathfrak{U}|_{y=0}-\mu(1+Mt)U\right]-\nu,~
\lim\limits_{y\to+\infty}G_{1k}=-\nu,\\
G_{1k}|_{t=0}=\partial_{x}\partial_{y}^{k+1}u_{in}(x,y)-\mu(U(0,x)-u_{in}(x,y))-\nu.
\end{cases}
\end{equation*}
where
\begin{align*}
\mathfrak{R}_{10}
=& e^{-\lambda_{10}t}\mathfrak{R}_{c1}
  -\mu(1+M\y t)e^{-\lambda_{10}t}\partial_{x}\mathfrak{U}(U-\mathfrak{U})
  +e^{-\lambda_{10}t}\partial_{y}^{-1}[\partial_{x}^{2}\mathfrak{U}]\partial_{y}\dot{\mathfrak{U}}\\
 &-\nu(\lambda_{10}+\partial_{x}\mathfrak{U})
  -e^{-\lambda_{10}t}\partial_{y}\partial_{x}\hat{F}
  -\eta e^{-\lambda_{10}t}\partial_{x}^{3}\partial_{y}u_{0}
  +\partial_{y}\eta e^{-\lambda_{10}t}\partial_{x}^{3}\mathbb{U},
\end{align*}
\begin{align*}
\mathfrak{R}_{11}
=& \mathfrak{R}_{c1}
  -\partial_{x}\dot{\mathfrak{U}}\partial_{x}\dot{\mathfrak{U}}
  -\dot{\mathfrak{U}}\partial_{x}^{2}\dot{\mathfrak{U}}
  +\partial_{x}^{2}\mathfrak{U}\partial_{y}\dot{\mathfrak{U}}
  +\partial_{y}^{-1}[\partial_{x}^{2}\mathfrak{U}]\partial_{y}^{2}\dot{\mathfrak{U}}
  -\partial_{y}^{2}\partial_{x}\hat{F}
  -\eta\partial_{x}^{3}\partial_{y}^{2}u_{0}
  +2\partial_{y}\eta\partial_{x}^{3}\partial_{y}\mathbb{U}
  +\partial_{y}^{2}\eta\partial_{x}^{3}\mathbb{U},
\end{align*}
and
\begin{align*}
\mathfrak{R}_{12}
=& e^{-\lambda_{12}t}\mathfrak{R}_{c1}
  +\mu(1+M\y t)e^{-\lambda_{12}t}\partial_{x}\mathfrak{U}(U-\mathfrak{U})
  -2e^{-\lambda_{12}t}\partial_{x}\dot{\mathfrak{U}}\partial_{y}\partial_{x}\dot{\mathfrak{U}}
  -2e^{-\lambda_{12}t}\dot{\mathfrak{U}}\partial_{x}^{2}\partial_{y}\dot{\mathfrak{U}}\\
 &+2e^{-\lambda_{12}t}\partial_{x}^{2}\mathfrak{U}\partial_{y}^{2}\dot{\mathfrak{U}}
  +e^{-\lambda_{12}t}\partial_{y}^{-1}[\partial_{x}^{2}\mathfrak{U}]\partial_{y}^{3}\dot{\mathfrak{U}}
  -\nu(\lambda_{12}-\partial_{x}\mathfrak{U})
  -e^{-\lambda_{12}t}\partial_{y}^{3}\partial_{x}\hat{F}
  -\eta e^{-\lambda_{12}t}\partial_{x}^{3}\partial_{y}^{3}u_{0}\\
 &+3\partial_{y}\eta e^{-\lambda_{12}t}\partial_{x}^{3}\partial_{y}^{2}\mathbb{U}
  +3\partial_{y}^{2}\eta e^{-\lambda_{12}t}\partial_{x}^{3}\partial_{y}\mathbb{U}
  +\partial_{y}^{3}\eta e^{-\lambda_{12}t}\partial_{x}^{3}\mathbb{U},
\end{align*}
with
\begin{align*}
\mathfrak{R}_{c1}
=&-\mu M\y(U-\mathfrak{U})
  +\mu\eta(1+M\y t)\partial_{x}^{2}(U-u_{0})
  +\mu(1+M\y t)\partial_{x}U(U-\mathfrak{U})\\
 &-2\mu Mt\dot{\mathfrak{U}}
  -\mu Mt\partial_{y}^{-1}[\partial_{x}\mathfrak{U}](U-\mathfrak{U})
  -\mu(1+M\y t)\hat{F}.
\end{align*}

By using \eqref{dotul}, one has
\begin{align*}
&|\partial_{y}^{-1}[\partial_{x}^{2}\mathfrak{U}]\partial_{y}\dot{\mathfrak{U}}|\leq C_{*}\y(U-\mathfrak{U}),~
 |\partial_{x}^{2}\mathfrak{U}\partial_{y}\dot{\mathfrak{U}}|\leq C_{*}(U-\mathfrak{U}),~\text{in}~\Omega_{\hat{t}}\\
&|\partial_{y}^{-1}[\partial_{x}^{2}\mathfrak{U}]\partial_{y}^{2}\dot{\mathfrak{U}}|\leq C_{*}\y(U-\mathfrak{U}),~
 |\partial_{y}^{-1}[\partial_{x}^{2}\mathfrak{U}]\partial_{y}^{3}\dot{\mathfrak{U}}|\leq C_{*}\y(U-\mathfrak{U}),~\text{in}~\Omega_{\hat{t}}.
\end{align*}
In addition, in virtue of  \eqref{u0},\eqref{U},\eqref{F} and
\begin{align*}
&|\partial_{x}\dot{\mathfrak{U}}\partial_{x}\dot{\mathfrak{U}}|\leq C_{*}e^{-\frac{4}{3}\delta\y}\leq C_{*}(U-\mathfrak{U}),~
 |\dot{\mathfrak{U}}\partial_{x}^{2}\dot{\mathfrak{U}}|\leq C_{*}e^{-\frac{4}{3}\delta\y}\leq C_{*}(U-\mathfrak{U}), ~\text{in}~\Omega_{\hat{t}}\\
&|\partial_{x}\dot{\mathfrak{U}}\partial_{x}\partial_{y}\dot{\mathfrak{U}}|\leq C_{*}e^{-\frac{4}{3}\delta\y}\leq C_{*}(U-\mathfrak{U}),~
 |\dot{\mathfrak{U}}\partial_{x}^{2}\partial_{y}\dot{\mathfrak{U}}|
 \leq C_{*}e^{-\frac{4}{3}\delta\y}\leq C_{*}(U-\mathfrak{U}),~\text{in}~\Omega_{\hat{t}}.
\end{align*}
one can take $M$ large enough depending on $\mathcal{Z}$, but $\mu$ large enough independent of $\mathcal{Z}$ to obtain
$$\mathfrak{R}_{1k}<0,~G_{1k}|_{y=0}<0,~G_{1k}|_{t=0}<0,~\text{in}~\Omega_{\hat{t}},$$
for $k=0,1,2$. Using the same argument as in the proof of \eqref{dotul}, one can obtain
\begin{align*}
\partial_{x}\dot{\mathfrak{U}}\leq \mu_{10}(U-\mathfrak{U}),~
\partial_{x}\partial_{y}\dot{\mathfrak{U}}\leq \mu_{11}(U-\mathfrak{U}),~
\partial_{x}\partial_{y}^{2}\dot{\mathfrak{U}}\leq \mu_{12}(U-\mathfrak{U}),
\end{align*}
for some positive constants $\mu_{10},\mu_{11}$ and $\mu_{12}$ independent of $\mathcal{Z}$.

To get \eqref{dotux1}, one need to apply the classical maximum principle for the parabolic equations to the problems obeyed by $G_{1k}^{-}$ for $k=0,1,2$. The argument is parallel to above, so we omit the detail.

\noindent{\bf \underline{Proof of \eqref{dotux2}.}}
It follows from \eqref{App} that $G_{2k}(k=0,1)$ obeys the follows problem in $\Omega_{\hat{t}}$,
\begin{equation*}
\begin{cases}
\partial_{t}G_{2k}
+\mathfrak{U}\partial_{x}G_{2k}
-\partial_{y}^{-1}[\partial_{x}\mathfrak{U}]\partial_{y}G_{2k}
+[\lambda_{2k}+(2-k)\partial_{x}\mathfrak{U}]G_{2k}
-\partial_{y}^{2}G_{2k}
-\eta\partial_{x}^{2}G_{2k}
=\mathfrak{R}_{2k},\\
G_{2k}|_{y=0}=e^{-\lambda_{2k}t}\left[\partial_{x}^{2}\partial_{y}^{k+1}\mathfrak{U}|_{y=0}-\mu(1+Mt)U\right]-\nu,~
\lim\limits_{y\to+\infty}G_{2k}=-\nu,\\
G_{2k}|_{t=0}=\partial_{x}^{2}\partial_{y}^{k+1}u_{in}(x,y)-\mu(U(0,x)-u_{in}(x,y))-\nu.
\end{cases}
\end{equation*}
where
\begin{align*}
\mathfrak{R}_{20}
=& e^{-\lambda_{20}t}\mathfrak{R}_{c1}
  -2\mu(1+M\y t)e^{-\lambda_{20}t}\partial_{x}\mathfrak{U}(U-\mathfrak{U})
  -e^{-\lambda_{20}t}\partial_{x}^{2}\mathfrak{U}\partial_{x}\dot{\mathfrak{U}}
  +e^{-\lambda_{20}t}\partial_{y}^{-1}[\partial_{x}^{3}\mathfrak{U}]\partial_{y}\dot{\mathfrak{U}}\\
 &+2e^{-\lambda_{20}t}\partial_{y}^{-1}[\partial_{x}^{2}\mathfrak{U}]\partial_{y}\partial_{x}\dot{\mathfrak{U}}
  -\nu(\lambda_{20}+2\partial_{x}\mathfrak{U})
  -e^{-\lambda_{20}t}\partial_{x}^{2}\partial_{y}\hat{F}
  -\eta e^{-\lambda_{20}t}\partial_{x}^{4}\partial_{y}u_{0}
  +\partial_{y}\eta e^{-\lambda_{12}t}\partial_{x}^{4}\mathbb{U},
\end{align*}
and
\begin{align*}
\mathfrak{R}_{21}
=& e^{-\lambda_{21}t}\mathfrak{R}_{c1}
  -\mu(1+M\y t)e^{-\lambda_{21}t}\partial_{x}\mathfrak{U}(U-\mathfrak{U})
  -e^{-\lambda_{21}t}\sum_{i=0}^{2}\partial_{x}^{2-i}\dot{\mathfrak{U}}\partial_{x}^{i+1}\dot{\mathfrak{U}}
  +e^{-\lambda_{21}t}\partial_{x}^{3}\mathfrak{U}\partial_{y}\dot{\mathfrak{U}}\\
 &+e^{-\lambda_{21}t}\partial_{x}^{2}\mathfrak{U}\partial_{x}\partial_{y}\dot{\mathfrak{U}}
  +e^{-\lambda_{21}t}\partial_{y}^{-1}[\partial_{x}^{3}\mathfrak{U}]\partial_{y}^{2}\dot{\mathfrak{U}}
  +2e^{-\lambda_{21}t}\partial_{y}^{-1}[\partial_{x}^{2}\mathfrak{U}]\partial_{x}\partial_{y}^{2}\dot{\mathfrak{U}}\\
 &-\nu(\lambda_{21}+\partial_{x}\mathfrak{U})
  -e^{-\lambda_{21}t}\partial_{x}^{2}\partial_{y}^{2}\hat{F}
  -\eta e^{-\lambda_{21}t}\partial_{x}^{4}\partial_{y}^{2}u_{0}
  +2\partial_{y}\eta e^{-\lambda_{21}t}\partial_{x}^{4}\partial_{y}\mathbb{U}
  +\partial_{y}^{2}\eta e^{-\lambda_{21}t}\partial_{x}^{4}\mathbb{U}.
\end{align*}

Taking advantage of \eqref{dotul} and \eqref{dotux1}, one has
$$|\partial_{x}^{2}\mathfrak{U}\partial_{x}\dot{\mathfrak{U}}|\leq C_{*}(U-\mathfrak{U}),~
  |\partial_{y}^{-1}[\partial_{x}^{3}\mathfrak{U}]\partial_{y}\dot{\mathfrak{U}}|\leq C_{*}\y(U-\mathfrak{U}),~
  |\partial_{y}^{-1}[\partial_{x}^{2}\mathfrak{U}]\partial_{x}\partial_{y}\dot{\mathfrak{U}}|\leq C_{*}\y(U-\mathfrak{U}),
  ~\text{in}~\Omega_{\hat{t}},$$
$$\left|\sum_{i=0}^{2}\partial_{x}^{2-i}\dot{\mathfrak{U}}\partial_{x}^{i+1}\dot{\mathfrak{U}}\right|\leq C_{*}(U-\mathfrak{U}),~
  |\partial_{x}^{3}\mathfrak{U}\partial_{y}\dot{\mathfrak{U}}|\leq C_{*}(U-\mathfrak{U}),~
  |\partial_{x}^{2}\mathfrak{U}\partial_{x}\partial_{y}\dot{\mathfrak{U}}|\leq C_{*}(U-\mathfrak{U}),~\text{in}~\Omega_{\hat{t}},$$
$$|\partial_{y}^{-1}[\partial_{x}^{3}\mathfrak{U}]\partial_{y}^{2}\dot{\mathfrak{U}}|\leq C_{*}\y(U-\mathfrak{U}),~
  |\partial_{y}^{-1}[\partial_{x}^{2}\mathfrak{U}]\partial_{x}\partial_{y}^{2}\dot{\mathfrak{U}}|\leq C_{*}\y(U-\mathfrak{U}),
  ~\text{in}~\Omega_{\hat{t}}.$$
Similar to the proof of \eqref{dotul} and \eqref{dotux1}, one can take $M$ large enough depending on $\mathcal{Z}$, but $\mu$ large enough independent of $\mathcal{Z}$ to obtain \eqref{dotux2} by applying the classical maximum principle for the parabolic equations to the problems obeyed by $G_{20},G_{21}$ and $G_{20}^{-},G_{21}^{-}$. The detail is also omitted.

\noindent{\bf \underline{Proof of \eqref{dotug}.}}
To prove \eqref{dotug}, let
$$G_{0}(t,x,y)=\dot{\mathfrak{U}}-\frac{\mu\varpi}{1+Mt}(U-\mathfrak{U})+\omega,$$
where $M,\mu,\nu$ are positive constants which will be determined later.
From \eqref{App} one knows that $G_{0}$ obeys the follows problem in $\Omega_{\hat{t}}$,
\begin{equation}\label{G0}
\begin{cases}
\partial_{t}G_{0}
+\mathfrak{U}\partial_{x}G_{0}
-\partial_{y}^{-1}[\partial_{x}\mathfrak{U}]\partial_{y}G_{0}
-\partial_{y}^{2}G_{0}
-\eta\partial_{x}^{2}G_{0}
=\mathfrak{R}_{g},\\
G_{0}|_{y=0}=\partial_{y}\mathfrak{U}|_{y=0}-\frac{\mu\varpi(0)U}{1+Mt}+\omega,~
\lim\limits_{y\to+\infty}G_{0}=\omega,\\
G_{0}|_{t=0}=\partial_{y}u_{in}(x,y)-\mu\varpi(U(0,x)-u_{in}(x,y))+\omega.
\end{cases}
\end{equation}
where
\begin{align*}
\mathfrak{R}_{g}
=& \frac{\mu M\varpi}{(1+Mt)^{2}}(U-\mathfrak{U})
  +\frac{\mu\eta\varpi}{1+Mt}\partial_{x}^{2}(U-u_{0})
  +\frac{\mu\varpi}{1+Mt}\partial_{x}U(U-\mathfrak{U})
  +\frac{\mu\varpi'}{1+Mt}\partial_{y}^{-1}[\partial_{x}\mathfrak{U}]\\
 &-\frac{2\mu\varpi'}{1+Mt}\dot{\mathfrak{U}}
  +\frac{\mu\varpi''}{1+Mt}(U-\mathfrak{U})
  -\frac{\mu\varpi}{1+Mt}\hat{F}
  -\partial_{y}\hat{F}
  -\eta\partial_{x}^{2}\partial_{y}u_{0}
  +\eta'\partial_{x}^{2}\mathbb{U}.
\end{align*}

Taking $\mu$ small enough, one can obtain
$$G_{0}|_{t=0}>0,~G_{0}|_{y=0}>0~\text{in}~\Omega_{\hat{t}}.$$
By the definition of $\varpi$, one can get $\mathfrak{R}_{g}>0$ in $\Omega_{\hat{t}}$ by taking $M$ large enough, if $\bar{\iota}$ is small enough.
For any fixed positive constant $\omega$, since
$$\left|\dot{\mathfrak{U}}-\frac{\mu\varpi}{1+Mt}(U-\mathfrak{U})\right|\leq C_{*}e^{-\frac{2}{3}\delta\y},$$
one knows that $G_{0}>0$ for any $y>R_{0}=\frac{3}{2\delta}\ln{\frac{C_{*}}{\omega}}$ in $\Omega_{\hat{t}}$.
By the classical maximum principle for parabolic equations one knows that
$$G_{0}>0,~\text{in}~\hat{\Omega}_{\hat{t}}^{R_{0}}.$$
Taking the limit $\omega\to 0^{+}$, one has
 $$\lim\limits_{\omega\to 0^{+}}G_{0}\geq 0,~\text{in}~\Omega_{\hat{t}},$$
which gives that
\begin{align}\label{dotu01g}
\dot{\mathfrak{U}}\geq k_{0}\varpi(U-\mathfrak{U}),~\text{in}~\Omega_{\hat{t}},
\end{align}
for some positives constants $k_{0}$ independent of $\mathcal{Z}$.

\noindent{\bf \underline{Proof of \eqref{dotub}.}}
Let $\mathcal{Q}(t,x,y)=\dot{\mathfrak{U}}-\nu\varpi-\omega t$, where $\nu$ and $\omega$ are positive constant which will be determined later.
Let $\Omega_{\hat{t},\varepsilon}=[0,\hat{t}]\times\mathbb{T}\times[y_{*}-\varepsilon,y_{*}+\varepsilon]$, where $\varepsilon$ is a positive constant satisfying $0<\varepsilon<\min\{y_{*}-\check{y}_{*},\hat{y}_{*}-y_{*}\}$, which will be determined later.
From \eqref{App} one knows that $\mathcal{Q}$ obeys the following problem in $\Omega_{\hat{t},\varepsilon}$,
\begin{equation}\label{Q}
\begin{cases}
\partial_{t}\mathcal{Q}
+\mathfrak{U}\partial_{x}\mathcal{Q}
-\partial_{y}^{-1}[\partial_{x}\mathfrak{U}]\partial_{y}\mathcal{Q}
=\partial_{y}^{2}\mathcal{Q}+\eta\partial_{x}^{2}\mathcal{Q}+\mathfrak{R}_{\varepsilon},\\
\mathcal{Q}|_{y=y_{*}-\varepsilon}=\partial_{y}\mathfrak{U}|_{y=y_{*}-\varepsilon}-\nu\mathcal{C}\varepsilon^{2}-\omega t,\\
\mathcal{Q}|_{y=y_{*}+\varepsilon}=\partial_{y}\mathfrak{U}|_{y=y_{*}+\varepsilon}-\nu\mathcal{C}\varepsilon^{2}-\omega t,\\
\mathcal{Q}|_{t=0}=\partial_{y}u_{in}(x,y)-\nu(y-y_{*})^{2}.
\end{cases}
\end{equation}
where
\begin{align*}
\mathfrak{R}_{\varepsilon}
= 2\nu\mathcal{C}
 -\omega+\nu\partial_{y}^{-1}[\partial_{x}\mathfrak{U}](y-y_{*})
 -\partial_{y}\hat{F}.
\end{align*}
Here we use the fact that $\eta$ is a constant and $\partial_{x}^{2}u_{0}=0$ for $y\in[y_{*}-\varepsilon,y_{*}+\varepsilon]$.
Using \eqref{dotux1}, one has
\begin{align*}
|\partial_{y}^{-1}[\partial_{x}\mathfrak{U}]||y-y_{*}|
=|\partial_{y}^{-1}[\partial_{y}^{-1}[\partial_{x}\partial_{y}\mathfrak{U}]]||y-y_{*}|
\leq K_{10}\hat{y}_{*}^{3}U\varepsilon,~\text{in}~\Omega_{\hat{t},\varepsilon}.
\end{align*}
Recall that $\|\partial_{y}\hat{F}\|_{L^{\infty}(\Omega)}\leq \bar{\iota}+C|y-y_{*}|+C_{*}t$ in $\Omega_{t^{*}}^{*}$,
thus one can take $\varepsilon,\omega$ and $\bar{\iota}$ small enough depending on $\nu$ but independent of $\mathcal{Z}$ and $\bar{t}\leq\hat{t}$ small enough such that $\mathfrak{R}_{\varepsilon}>0$ in $\Omega_{\bar{t},\varepsilon}$. Then, in virtue of \eqref{dotug}, one can take $\nu$ small enough independent of $\mathcal{Z}$ to get
\begin{align*}
\mathcal{Q}|_{y=y_{*}-\varepsilon}>0,~\mathcal{Q}|_{y=y_{*}+\varepsilon}>0,~\mathcal{Q}|_{t=0}>0,~\text{in}~\Omega_{\bar{t},\varepsilon}.
\end{align*}
Applying the classical maximum principle for the parabolic equations to the problem \eqref{Q}, one can obtain
\begin{align*}
\mathcal{Q}\geq0,~\text{in}~\Omega_{\bar{t},\varepsilon},
\end{align*}
which gives
\begin{align}\label{pdotug'}
\dot{\mathfrak{U}}\geq \nu\varpi+\omega t,~\text{in}~\Omega_{\bar{t},\varepsilon}.
\end{align}
Combining \eqref{dotug} and \eqref{pdotug'} one can get left inequality in \eqref{dotub}.

To get the right inequality in \eqref{dotub}, we restrict \eqref{App} to the region $\hat{\Omega}_{\bar{t}}\cap\check{\Omega}_{\bar{t}}$ to get
\begin{align}
\partial_{t}\dot{\mathfrak{U}}
+\mathfrak{U}\partial_{x}\dot{\mathfrak{U}}
-\partial_{y}^{-1}[\partial_{x}\mathfrak{U}]\partial_{y}\dot{\mathfrak{U}}
=\partial_{y}^{2}\dot{\mathfrak{U}}+\eta\partial_{x}^{2}\dot{\mathfrak{U}}-\partial_{y}\hat{F},~\text{in}~\hat{\Omega}_{\bar{t}}\cap\check{\Omega}_{\bar{t}},
\end{align}
with the initial condition $\dot{\mathfrak{U}}|_{t=0}=\partial_{y}u_{in}$. Thus for any $t\in[0,\bar{t}]$, by using \eqref{dotul}, \eqref{dotux1} and \eqref{Fe}, we have
\begin{align*}
\dot{\mathfrak{U}}
\leq \partial_{y}u_{in}+\int_{0}^{t}\partial_{t}\dot{\mathfrak{U}}(\sigma,x,y)d\sigma
\leq C\varpi+Ct+C_{*}t^{2},~\text{in}~\hat{\Omega}_{\bar{t}}\cap\check{\Omega}_{\bar{t}},
\end{align*}
which gives the right inequality in \eqref{dotub} by choosing $t_{1}\leq\bar{t}$ small enough.

\end{proof}

Now we shall use Lemma \ref{DI}--\ref{sdotu} to obtain the estimate of $\rho$ and $\mathcal{P}$.
\begin{lemma}\label{rhoP}
Let $\theta^{*}$ be obtained in Lemma \ref{wdotu} and $t_{1}$ be obtained in Lemma \ref{sdotu}.
Assume the initial date $\tilde{u}_{in}$ in the problem \eqref{App+1} satisfies Assumption \ref{inA'} and \ref{inA}.
There exists positive constant $K$ independent of $\mathcal{Z}$ such that for any $\theta\geq\theta^{*}$, there is a corresponding $t_{2}\in(0,t_{1}]$ such that for any $t_{v}\in[0,t_{2}]$, if $\mathfrak{U}$ is a solution to the problem \eqref{App+1} on the time interval $[0,t_{v}]$ such that
$\|\mathbb{U}\|_{L_{t}^{\infty}\hat{H}_{\psi}^{s}(\Omega_{t_{v}})}\leq C_{*}$,
and Assumption \ref{FA} holds for $t_{v}$, then the following inequalities hold in $\Omega_{t_{v}}$,
\begin{align}\label{rhop1}
|\rho|\leq K\mathcal{P},~
|\partial_{y}\rho|\leq K\mathcal{P}^{2},~
|\partial_{x}\rho|\leq K(1+\y t)\mathcal{P},~
|\partial_{x}^{2}\rho|\leq K(1+\y t)^{2}\mathcal{P},~
~\text{in}~\Omega_{t_{v}},
\end{align}
\begin{align}\label{rhop2}
|\partial_{x}\dot{\mathfrak{U}}|\leq K(1+\y t)\dot{\mathfrak{U}},~
|\partial_{x}^{2}\dot{\mathfrak{U}}|\leq C(1+\y t)\dot{\mathfrak{U}},~
|\partial_{x}^{2}\partial_{y}\dot{\mathfrak{U}}|\leq C(1+\y t)\dot{\mathfrak{U}}^{\frac{1}{2}},~\text{in}~\Omega_{t_{v}}.
\end{align}
\begin{align}\label{rhop3}
|\partial_{x}\mathcal{P}|\leq K\mathcal{P},~
\left|\partial_{x}(\mathcal{P}\partial_{x}\rho)\right|\leq K(1+\y t)^{2}\mathcal{P}^{2}~\text{in}~\Omega_{t_{v}}.
\end{align}
\end{lemma}
\begin{proof}

As in the proof of Lemma \ref{sdotu}, to prove Lemma \ref{rhoP} we only need to prove that for any $\theta\geq \theta^{*}$, there exist $t_{2}\in(0,t_{1}]$ such that \eqref{rhop1}--\eqref{rhop3} hold for $t_{v}=t_{2}$, under the conditions that $\|\mathbb{U}\|_{L_{t}^{\infty}\hat{H}_{\psi}^{s}(\Omega_{t_{1}})}\leq C_{*}$, and Assumption \ref{FA} holds for $t_{1}$.

In the region $\Omega_{t_{1}}\setminus\Omega_{t_{1}}^{*}$, \eqref{rhop1} and \eqref{rhop2} are corollaries of Lemma \ref{sdotu}.
In fact, by using Lemma \ref{sdotu}, one has
\begin{align}\label{Partrho}
|\rho|
=\left|\frac{\partial_{y}\dot{\mathfrak{U}}+\varsigma'}{\dot{\mathfrak{U}}+\varsigma}\right|
\leq \frac{K_{01}(U-\mathfrak{U})+C\bar{\varsigma}}{k_{0}\varpi(U-\mathfrak{U})+\varsigma}
\leq C
\leq C\mathcal{P},~\text{in}~\Omega_{t_{1}}\setminus\Omega_{t_{1}}^{*}.
\end{align}
Using Lemma \ref{sdotu} and \eqref{Partrho}, it follows that
\begin{align}\label{Partrhoy}
|\partial_{y}\rho|
\leq \left|\frac{\partial_{y}^{2}\dot{\mathfrak{U}}+\varsigma''}{\dot{\mathfrak{U}}+\varsigma}\right|+\rho^{2}
\leq \frac{K_{02}(U-\mathfrak{U})+C\bar{\varsigma}}{k_{0}\varpi(U-\mathfrak{U})+\varsigma}+C\mathcal{P}^{2}
\leq C\mathcal{P}^{2},~\text{in}~\Omega_{t_{1}}\setminus\Omega_{t_{1}}^{*}.
\end{align}
A direct application of Lemma \ref{sdotu} gives
\begin{align*}
\left|\frac{\partial_{x}\dot{\mathfrak{U}}}{\dot{\mathfrak{U}}}\right|
\leq \frac{K_{10}(1+\y t)(U-\mathfrak{U})}{k_{0}\varpi(U-\mathfrak{U})}
\leq C(1+\y t), ~\text{in}~\Omega_{t_{1}}\setminus\Omega_{t_{1}}^{*},
\end{align*}
\begin{align*}
\left|\frac{\partial_{x}^{2}\dot{\mathfrak{U}}}{\dot{\mathfrak{U}}}\right|
\leq \frac{K_{20}(1+\y t)(U-\mathfrak{U})}{k_{0}\varpi(U-\mathfrak{U})}
\leq C(1+\y t), ~\text{in}~\Omega_{t_{1}}\setminus\Omega_{t_{1}}^{*},
\end{align*}
and
\begin{align*}
\left|\frac{\partial_{x}^{2}\partial_{y}\dot{\mathfrak{U}}}{\dot{\mathfrak{U}}}\right|
\leq \frac{K_{21}(1+\y t)(U-\mathfrak{U})}{k_{0}\varpi(U-\mathfrak{U})}
\leq C(1+\y t), ~\text{in}~\Omega_{t_{1}}\setminus\Omega_{t_{1}}^{*},
\end{align*}
which lead to
\begin{align}\label{Partux/u}
|\partial_{x}\dot{\mathfrak{U}}|\leq C(1+\y t)\dot{\mathfrak{U}},~
|\partial_{x}^{2}\dot{\mathfrak{U}}|\leq C(1+\y t)\dot{\mathfrak{U}},~\text{and}~
|\partial_{x}^{2}\partial_{y}\dot{\mathfrak{U}}|\leq C(1+\y t)\dot{\mathfrak{U}},~\text{in}~\Omega_{t_{1}}\setminus\Omega_{t_{1}}^{*}.
\end{align}
Moreover, by a direct calculation, from \eqref{Partux/u} one can deduce
\begin{align}\label{Partpx}
|\partial_{x}\mathcal{P}|
=\frac{1}{2}\frac{\chi}{\sqrt{\dot{\mathfrak{U}}+\varsigma}}\left|\frac{\partial_{x}\dot{\mathfrak{U}}}{\dot{\mathfrak{U}}+\varsigma}\right|
\leq C\frac{\chi(1+\y t)}{\sqrt{\dot{\mathfrak{U}}+\varsigma}}
\leq C\mathcal{P},~\text{in}~\Omega_{t_{1}}\setminus\Omega_{t_{1}}^{*}.
\end{align}
Utilizing \eqref{Partrho},\eqref{Partux/u} and Lemma \ref{sdotu}, one has
\begin{align}\label{Partrhox}
|\partial_{x}\rho|
\leq& \left|\frac{\partial_{x}\partial_{y}\dot{\mathfrak{U}}}{\dot{\mathfrak{U}}+\varsigma}\right|
    +|\rho|\left|\frac{\partial_{x}\dot{\mathfrak{U}}}{\dot{\mathfrak{U}}+\varsigma}\right|
\leq \frac{K_{11}(1+\y t)(U-\mathfrak{U})}{k_{0}\varpi(U-\mathfrak{U})+\varsigma}
    +C(1+\y t)\mathcal{P}\\
\leq& C(1+\y t)\mathcal{P},~\text{in}~\Omega_{t_{1}}\setminus\Omega_{t_{1}}^{*}.\nonumber
\end{align}
A direct calculation gives that
\begin{align}\label{Rhox2}
\partial_{x}^{2}\rho
= \frac{\partial_{x}^{2}\partial_{y}\dot{\mathfrak{U}}}{\dot{\mathfrak{U}}+\varsigma}
 -2\frac{\partial_{x}\partial_{y}\dot{\mathfrak{U}}\partial_{x}\dot{\mathfrak{U}}}{(\dot{\mathfrak{U}}+\varsigma)^{2}}
 -\rho\frac{\partial_{x}^{2}\dot{\mathfrak{U}}}{\dot{\mathfrak{U}}+\varsigma}
 +2\rho\left(\frac{\partial_{x}\dot{\mathfrak{U}}}{\dot{\mathfrak{U}}+\varsigma}\right)^{2}.
\end{align}
With \eqref{Partrho},\eqref{Partux/u} and Lemma \ref{sdotu}, one can deduce from \eqref{Rhox2} that
\begin{align}\label{Partrhox2}
|\partial_{x}^{2}\rho|
\leq& \left|\frac{\partial_{x}^{2}\partial_{y}\dot{\mathfrak{U}}}{\dot{\mathfrak{U}}+\varsigma}\right|
     +2\left|\frac{\partial_{x}\partial_{y}\dot{\mathfrak{U}}}{\dot{\mathfrak{U}}+\varsigma}\right|
      \left|\frac{\partial_{x}\dot{\mathfrak{U}}}{\dot{\mathfrak{U}}+\varsigma}\right|
     +|\rho|\left|\frac{\partial_{x}\dot{\mathfrak{U}}}{\dot{\mathfrak{U}}+\varsigma}\right|
     +2\rho\left|\frac{\partial_{x}\dot{\mathfrak{U}}}{\dot{\mathfrak{U}}+\varsigma}\right|^{2}\\
\leq& \frac{K_{21}(1+\y t)(U-\mathfrak{U})}{k_{0}\varpi(U-\mathfrak{U})+\varsigma}
     +C(1+\y t)^{2}\mathcal{P}\nonumber\\
\leq& C(1+\y t)^{2}\mathcal{P},~\text{in}~\Omega_{t_{1}}\setminus\Omega_{t_{1}}^{*}.\nonumber
\end{align}
By virtue of \eqref{Partpx}, \eqref{Partrhox} and \eqref{Partrhox2}, one can get that
\begin{align}\label{Partrhopx}
|\partial_{x}(\mathcal{P}\partial_{x}\rho)|
\leq |\partial_{x}\mathcal{P}||\partial_{x}\rho|+\mathcal{P}|\partial_{x}^{2}\rho|
\leq C(1+\y t)\mathcal{P}^{2}+\mathcal{P}|\partial_{x}^{2}\rho|
\leq C(1+\y t)^{2}\mathcal{P}^{2},~\text{in}~\Omega_{t_{1}}\setminus\Omega_{t_{1}}^{*}.
\end{align}

Now we focus on the proof of \eqref{rhop1}-\eqref{rhop3} in the domain $\Omega_{t_{1}}^{*}$. Recall the definition of $\vartheta_{c}$ and $u_{0}$, from \eqref{App+1} one knows that $\dot{\mathfrak{U}}$ obeys the follows equation in $\Omega_{t_{1}}^{*}$,
\begin{equation}\label{partdotu}
\partial_{t}\dot{\mathfrak{U}}
+\mathfrak{U}\partial_{x}\dot{\mathfrak{U}}
-\partial_{y}^{-1}[\partial_{x}\mathfrak{U}]\partial_{y}\dot{\mathfrak{U}}
=\partial_{y}^{2}\dot{\mathfrak{U}}+\eta\partial_{x}^{2}\dot{\mathfrak{U}}-\partial_{y}\hat{F},
\end{equation}
with the initial condition $\dot{\mathfrak{U}}|_{t=0}=\partial_{y}u_{in}(x,y)$.

Let $\gamma=(\gamma_{1},\gamma_{2})$ be any multiindex in $\Gamma_{3}$.
For $k\in\{1,2\}$, by acting $\partial_{t}^{k-1}D^{\gamma}$ on \eqref{partdotu} one can obtain
\begin{align*}
\partial_{t}^{k}D^{\gamma}\dot{\mathfrak{U}}
=& \partial_{y}^{2}\partial_{t}^{k-1}D^{\gamma}\dot{\mathfrak{U}}
  +\eta\partial_{x}^{2}\partial_{t}^{k-1}D^{\gamma}\dot{\mathfrak{U}}
  -\partial_{t}^{k-1}D^{\gamma}[\mathfrak{U}\partial_{x}\dot{\mathfrak{U}}]\\
 &+\partial_{t}^{k-1}D^{\gamma}[\partial_{y}^{-1}[\partial_{x}\mathfrak{U}]\partial_{y}\dot{\mathfrak{U}}]
  -\partial_{t}^{k-1}D^{\gamma}\partial_{y}\hat{F},\nonumber
\end{align*}
which gives
\begin{align*}
\|\partial_{t}^{k}D^{\gamma}\dot{\mathfrak{U}}\|_{L^{\infty}(\Omega^{*})}(t)\leq C_{*}~\text{and}~
\|\partial_{t}D^{\gamma}\dot{\mathfrak{U}}\|_{L^{\infty}(\Omega^{*})}(0)\leq C,~\forall t\in[0,t_{1}].
\end{align*}
Utilizing Taylor's formula, we have
\begin{align*}
D^{\gamma}\dot{\mathfrak{U}}(t,x,y)
= \int_{0}^{t}\partial_{t}^{2}D^{\gamma}\dot{\mathfrak{U}}(s,x,y)(t-s)ds
 +\partial_{t}D^{\gamma}\dot{\mathfrak{U}}(0,x,y)t
 +D^{\gamma}\partial_{y}u_{in}(x,y),~\text{in}~\Omega_{t_{1}}^{*}.
\end{align*}
It follows that
\begin{align*}
|D^{\gamma}\dot{\mathfrak{U}}(t,x,y)|
\leq& \int_{0}^{t}|\partial_{t}^{2}D^{\gamma}\dot{\mathfrak{U}}(s,x,y)|(t-s)ds
 +|\partial_{t}D^{\gamma}\dot{\mathfrak{U}}(0,x,y)|t
 +|D^{\gamma}\partial_{y}u_{in}(x,y)|\\
\leq& |D^{\gamma}\partial_{y}u_{in}(x,y)|+Ct+C_{*}t^{2},~\text{in}~\Omega_{t_{1}}^{*},
\end{align*}
from which we have
\begin{align}\label{ddotu}
|D^{\gamma}\dot{\mathfrak{U}}(t,x,y)|\leq |D^{\gamma}\partial_{y}u_{in}(x,y)|+Ct,~\text{in}~\Omega_{t_{2}}^{*},
\end{align}
for some $t_{2}\leq t_{1}$. By using \eqref{inA1}--\eqref{inA3} and \eqref{ddotu}, one has
\begin{align}\label{partrho}
|\rho|=\left|\frac{\partial_{y}\dot{\mathfrak{U}}+\varsigma'}{\dot{\mathfrak{U}}+\varsigma}\right|
\leq C\frac{|y-y_{*}|+t+\bar{\varsigma}}{\sqrt{\dot{\mathfrak{U}}+\varsigma}\sqrt{k_{1}(y-y_{*})^{2}+k_{2}t+\bar{\varsigma}}}
\leq C\frac{1}{\sqrt{\dot{\mathfrak{U}}+\varsigma}}
=C\mathcal{P},~\text{in}~\Omega_{t_{2}}^{*},
\end{align}
and with Lemma \ref{sdotu} in addition, one can obtain
\begin{align}\label{partux}
|\partial_{x}\dot{\mathfrak{U}}|
\leq |\partial_{x}\partial_{y}u_{in}(x,y)|+Ct
\leq C(y-y_{*})^{2}+Ct
\leq C\dot{\mathfrak{U}},~\text{in}~\Omega_{t_{2}}^{*}.
\end{align}
Similarly,
\begin{align}\label{partux2}
|\partial_{x}^{2}\dot{\mathfrak{U}}|
\leq |\partial_{x}^{2}\partial_{y}u_{in}(x,y)|+Ct
\leq C(y-y_{*})^{2}+Ct
\leq C\dot{\mathfrak{U}},~\text{in}~\Omega_{t_{2}}^{*}.
\end{align}
Moreover,
\begin{align}\label{partux2y}
|\partial_{x}^{2}\partial_{y}\dot{\mathfrak{U}}|
\leq |\partial_{x}^{2}\partial_{y}^{2}u_{in}(x,y)|+Ct
\leq C|y-y_{*}|+Ct
\leq C\dot{\mathfrak{U}}^{\frac{1}{2}},~\text{in}~\Omega_{t_{2}}^{*}.
\end{align}
The inequality \eqref{ddotu} gives  $|\partial_{y}^{2}\dot{\mathfrak{U}}|\leq|\partial_{y}^{3}u_{in}(x,y)|+Ct\leq C$ in $\Omega_{t_{2}}^{*}$,
as a consequence of \eqref{partrho}, one has
\begin{align}\label{partrhoy}
|\partial_{y}\rho|
\leq \left|\frac{\partial_{y}^{2}\dot{\mathfrak{U}}+\varsigma''}{\dot{\mathfrak{U}}+\varsigma}\right|+\rho^{2}
\leq C\mathcal{P}^{2},~\text{in}~\Omega_{t_{2}}^{*}.
\end{align}
With the help of \eqref{partux}, one can obtain
\begin{align}\label{partpx}
|\partial_{x}\mathcal{P}|
=\frac{1}{2}\mathcal{P}\left|\frac{\partial_{x}\dot{\mathfrak{U}}}{\dot{\mathfrak{U}}+\varsigma}\right|
\leq \frac{1}{2}\mathcal{P}\frac{C\y\dot{\mathfrak{U}}}{\dot{\mathfrak{U}}+\varsigma}
\leq C\mathcal{P},~\text{in}~\Omega_{t_{2}}^{*}.
\end{align}
From \eqref{ddotu} one has $|\partial_{x}\partial_{y}\dot{\mathfrak{U}}(t,x,y)|\leq |\partial_{x}\partial_{y}^{2}u_{in}(x,y)|+Ct$ in $\Omega_{t_{2}}^{*}$, which together with \eqref{inA1}--\eqref{inA3},\eqref{partrho} and \eqref{partux} gives
\begin{align}\label{partrhox}
|\partial_{x}\rho|
\leq \left|\frac{\partial_{x}\partial_{y}\dot{\mathfrak{U}}}{\dot{\mathfrak{U}}+\varsigma}\right|
    +|\rho|\left|\frac{\partial_{x}\dot{\mathfrak{U}}}{\dot{\mathfrak{U}}+\varsigma}\right|
\leq C\frac{|y-y_{*}|+t}{\dot{\mathfrak{U}}+\varsigma}+C|\rho|
\leq C\mathcal{P}+C|\rho|
\leq C\mathcal{P},~\text{in}~\Omega_{t_{2}}^{*}.
\end{align}
Again by \eqref{inA1}--\eqref{inA3}, \eqref{ddotu} and Lemma \ref{sdotu}, one has
\begin{align*}
\left|\frac{\partial_{x}^{2}\partial_{y}\dot{\mathfrak{U}}}{\dot{\mathfrak{U}}+\varsigma}\right|
\leq C\frac{|y-y_{*}|+t}{\dot{\mathfrak{U}}+\varsigma}
\leq C\mathcal{P},~\text{in}~\Omega_{t_{2}}^{*}.
\end{align*}
thus from \eqref{Rhox2}, \eqref{partrho}, \eqref{partux} and \eqref{partux2} one knows
\begin{align}\label{partrhox2}
|\partial_{x}^{2}\rho|
\leq \left|\frac{\partial_{x}^{2}\partial_{y}\dot{\mathfrak{U}}}{\dot{\mathfrak{U}}+\varsigma}\right|
    +2\left|\frac{\partial_{x}\partial_{y}\dot{\mathfrak{U}}}{\dot{\mathfrak{U}}+\varsigma}\right|
      \left|\frac{\partial_{x}\dot{\mathfrak{U}}}{\dot{\mathfrak{U}}+\varsigma}\right|
    +|\rho|\left|\frac{\partial_{x}\dot{\mathfrak{U}}}{\dot{\mathfrak{U}}+\varsigma}\right|
    +2\rho\left|\frac{\partial_{x}\dot{\mathfrak{U}}}{\dot{\mathfrak{U}}+\varsigma}\right|^{2}
\leq C\mathcal{P},~\text{in}~\Omega_{t_{2}}^{*}.
\end{align}
As a consequence of \eqref{partpx},\eqref{partrhox},\eqref{partux} and \eqref{partrhox2}, one has
\begin{align}\label{partrhopx}
|\partial_{x}(\mathcal{P}\partial_{x}\rho)|
\leq |\partial_{x}\mathcal{P}||\partial_{x}\rho|+|\mathcal{P}\partial_{x}^{2}\rho|
\leq C\mathcal{P}^{2}+\mathcal{P}|\partial_{x}^{2}\rho|
\leq C\mathcal{P}^{2},~\text{in}~\Omega_{t_{2}}^{*}.
\end{align}

Combining \eqref{Partrho}--\eqref{Partrhopx} with \eqref{partrho}--\eqref{partrhopx}, with the help of Lemma \ref{wdotu} we get the conclusion.
\end{proof}

\begin{lemma}\label{rhoP'}
Let $\theta^{*}$ be obtained in Lemma \ref{wdotu} and $t_{2}$ be obtained in Lemma \ref{rhoP}.
Assume the initial date $\tilde{u}_{in}$ in the problem \eqref{App+1} satisfies \eqref{inA1}--\eqref{inA3}, and Assumption \ref{inA}.
There exists positive constant $K$ independent of $C$ such that for any $\theta\geq\theta^{*}$, there is a corresponding $t_{3}\in(0,t_{2}]$ such that for any $t_{v}\in[0,t_{3}]$, if $\mathfrak{U}$ is a solution to the problem \eqref{App+1} on the time interval $[0,t_{v}]$ such that
$\|\mathbb{U}\|_{L_{t}^{\infty}\hat{H}_{\psi}^{s}(\Omega_{t_{v}})}\leq C_{*}$,
and Assumption \ref{FA} holds for $t_{v}$, then the following inequalities hold in $\Omega_{t_{v}}$,
\begin{align}\label{rhox}
|\partial_{x}\rho|\leq C\mathcal{P},~\text{in}~\Omega_{t_{v}},
\end{align}
and
\begin{align}\label{rhol}
\rho\leq -\frac{15}{16}\delta,~\text{in}~\check{\Omega}_{t_{v}}^{Y}.
\end{align}

\end{lemma}
\begin{proof}
By the same argument as in the proof of Lemma \ref{wdotu}, to obtain Lemma \ref{rhoP'} we only need to prove that for any $\theta\geq \theta^{*}$, there exist $t_{3}\in(0,t_{2}]$ such that \eqref{rhox}--\eqref{rhol} hold for $t_{v}=t_{3}$, under the conditions that $\|\mathbb{U}\|_{L_{t}^{\infty}\hat{H}_{\psi}^{s}(\Omega_{t_{2}})}\leq C_{*}$, and Assumption \ref{FA} holds for $t_{2}$.

By \eqref{App+1} one can check that $\rho$ obeys the follows equation in $\check{\Omega}_{t_{2}}^{\hat{y}}$,
\begin{align}\label{rhoe}
\partial_{t}\rho
+\mathfrak{U}\partial_{x}\rho
+\mathfrak{L}_{1}\partial_{y}\rho
+\mathfrak{L}_{0}\rho
-\partial_{y}^{2}\rho
=\mathfrak{R}_{\rho},
\end{align}
with the initial date $\rho|_{t=0}=\frac{\partial_{y}^{2}u_{in}}{\partial_{y}u_{in}}$, where
\begin{align*}
\mathfrak{L}_{1}=-\partial_{y}^{-1}[\partial_{x}\mathfrak{U}]-2\rho,
\end{align*}
\begin{align*}
\mathfrak{L}_{0}
=-\partial_{x}\mathfrak{U}
 -\frac{\partial_{y}\hat{F}-\partial_{y}\eta\partial_{x}^{2}\mathbb{U}-\eta\partial_{x}^{2}\partial_{y}\mathbb{U}}{\dot{\mathfrak{U}}},
\end{align*}
and
\begin{align*}
\mathfrak{R}_{\rho}
=&-\partial_{x}\partial_{y}\mathfrak{U}
  -\frac{\partial_{y}^{2}\hat{F}}{\dot{\mathfrak{U}}}
  +\frac{2\partial_{y}\eta\partial_{x}^{2}\partial_{y}\mathbb{U}+\partial_{y}^{2}\eta\partial_{x}^{2}\mathbb{U}+\eta\partial_{x}^{2}\partial_{y}^{2}\mathbb{U}}{\dot{\mathfrak{U}}}.
\end{align*}

\noindent{\bf \underline{Proof of \eqref{rhox}.}}
From Lemma \ref{rhoP} one can see that $|\partial_{x}\rho|\leq C\mathcal{P}$ holds in $\hat{\Omega}_{t_{2}}$. To get \eqref{rhox} we let $\mathfrak{m}=\|\y^{-1}\mathfrak{L}_{1}\|_{L^{\infty}(\Omega)}+2$, define
\begin{align*}
\mathfrak{D}^{y}(t,x,y)=e^{-Mt}[\partial_{y}\rho-\mu(1+Nt)\y^{-1}]-\nu\ln\y-\nu \mathfrak{m}t,
\end{align*}
and
\begin{align*}
\mathfrak{D}^{x}(t,x,y)=e^{-\bar{M}t}\y^{-1}[\partial_{x}\rho-\bar{\mu}(1+\bar{N}t)]-\bar{\nu}\ln\y-\bar{\nu}\mathfrak{m}t,
\end{align*}
where $M,N,\bar{M},\bar{N},\nu,\bar{\nu}$ are positive constant to be determined. From \eqref{rhoe} one knows that $\mathfrak{D}^{y}$ obeys the following equation in $\check{\Omega}_{t_{2}}^{\hat{y}}$,
\begin{align*}
\partial_{t}\mathfrak{D}^{y}
+\mathfrak{U}\partial_{x}\mathfrak{D}^{y}
+\mathfrak{L}_{1}\partial_{y}\mathfrak{D}^{y}
+(M+\mathfrak{L}_{0}+\partial_{y}\mathfrak{L}_{1})\mathfrak{D}^{y}
-\partial_{y}^{2}\mathfrak{D}^{y}
=\mathfrak{R}^{y},
\end{align*}
with the initial date
$\mathfrak{D}^{y}|_{t=0}=\frac{\partial_{y}^{3}u_{in}}{\partial_{y}u_{in}}-\left(\frac{\partial_{y}^{2}u_{in}}{\partial_{y}u_{in}}\right)^{2}
-\mu\y^{-1}-\nu\ln\y$, and
\begin{align*}
\mathfrak{R}^{y}=& e^{-Mt}\partial_{y}\mathfrak{R}_{\rho}
  -e^{-Mt}\partial_{y}\mathfrak{U}\partial_{x}\rho
  -e^{-Mt}\partial_{y}\mathfrak{L}_{0}\rho
  -\mu Ne^{-Mt}\y^{-1}\\
 &+\mu(1+Nt)e^{-Mt}\y^{-2}\mathfrak{L}_{1}
  -\mu(1+Nt)e^{-Mt}\y^{-1}(\mathfrak{L}_{0}+\partial_{y}\mathfrak{L}_{1})\\
 &+2e^{-Mt}\mu(1+Nt)\y^{-3}
  +\nu(2\y^{-2}-\y^{-1}\mathfrak{L}_{1}-\mathfrak{m})
  -\nu(M+\mathfrak{L}_{0}+\partial_{y}\mathfrak{L}_{1})(\ln\y+\mathfrak{m}t),
\end{align*}

Let $\hat{t}=\min\{\frac{1}{N},t_{2}\}$. For any $\nu>0$, using Lemma \ref{rhoP} one can take $\mu$ large enough but independent of $\mathcal{Z}$ to get
\begin{align*}
\mathfrak{D}^{y}|_{t=0}<0,~
\mathfrak{D}^{y}|_{y=\hat{y}}<0,~\text{in}~\check{\Omega}_{\hat{t}}^{\hat{y}}.
\end{align*}
Under the condition \eqref{Fd}, thanks to Lemma \ref{DI}--Lemma \ref{rhoP} and Corollary \eqref{wdotuC}, one can obtain
\begin{align*}
|\partial_{y}\mathfrak{R}_{\rho}|\leq C\y^{-1},~
|\partial_{y}\mathfrak{U}\partial_{x}\rho|\leq C\y^{-1},~
|\partial_{y}\mathfrak{L}_{0}\rho|\leq C\y^{-1},~
\mathfrak{L}_{0}\leq C_{*},~
\mathfrak{L}_{1}\leq C_{*}\y,~
\partial_{y}\mathfrak{L}_{1}\leq C_{*}~\text{in}~\check{\Omega}_{\hat{t}}^{\hat{y}}.
\end{align*}
Thus by taking $M$ and $N$ large enough, one can get $\mathfrak{R}^{y}<0$ in $\check{\Omega}_{\hat{t}}^{\hat{y}}$.
Lemma \ref{rhoP} tells that $|\partial_{y}\rho|\leq K$ in $\check{\Omega}_{\hat{t}}^{\hat{y}}$, from which one knows that for any $y>R_{\nu}=e^{\frac{K+2\mu}{\nu}}$, $\mathfrak{D}^{y}<0$ in $\check{\Omega}_{\hat{t}}^{\hat{y}}$.
By employing the classical maximum principle for the parabolic equations, one has
\begin{align*}
\mathfrak{D}^{y}<0,~\text{in}~\check{\Omega}_{\hat{t}}^{\hat{y}}\cap\hat{\Omega}_{\hat{t}}^{R_{\nu}}.
\end{align*}
Taking the limit $\nu\to 0^{+}$, one has $e^{-Mt}[\partial_{y}\rho-\mu(1+Nt)\y^{-1}]\leq0$ in $\check{\Omega}_{\hat{t}}^{\hat{y}}$, which gives
\begin{align}\label{rhoy+}
\partial_{y}\rho\leq\mu(1+Nt)\y^{-1}, ~\text{in}~\check{\Omega}_{\hat{t}}^{\hat{y}}.
\end{align}
A similar argument applied to the function
\begin{align*}
\mathfrak{D}^{y}_{-}=e^{-Mt}[-\partial_{y}\rho-\mu(1+Nt)\y^{-1}]-\nu\ln\y-\nu \mathfrak{m}t.
\end{align*}
gives
\begin{align}\label{rhoy-}
-\partial_{y}\rho\leq\mu(1+Nt)\y^{-1}, ~\text{in}~\check{\Omega}_{\hat{t}}^{\hat{y}},
\end{align}
provided that $N$ is large enough. Here we omit the detail. Combining \eqref{rhoy+} with \eqref{rhoy-}, one has
\begin{align}\label{rhoy}
|\partial_{y}\rho|\leq\mu(1+Nt)\y^{-1}, ~\text{in}~\check{\Omega}_{\hat{t}}^{\hat{y}}.
\end{align}

Now with the help of \eqref{rhoy}, we are able to prove \eqref{rhox}.
One can check that $\mathfrak{D}^{x}$ obeys the following equation in $\check{\Omega}_{t_{2}}^{\hat{y}}$,
\begin{align*}
&\partial_{t}\mathfrak{D}^{x}
+\mathfrak{L}_{1}\partial_{x}\mathfrak{D}^{x}
+(\mathfrak{L}_{2}-2\y^{-1})\partial_{y}\mathfrak{D}^{x}
+(\bar{M}+\mathfrak{L}_{0}+\partial_{x}\mathfrak{U}-2\partial_{y}\rho+\y^{-1}\mathfrak{L}_{1})\mathfrak{D}^{x}
-\partial_{y}^{2}\mathfrak{D}^{x}
=\mathfrak{R}^{x},
\end{align*}
with the initial date
$\mathfrak{D}^{x}|_{t=0}
=\frac{\partial_{x}\partial_{y}^{2}u_{in}}{\partial_{y}u_{in}}\y^{-1}
-\frac{\partial_{y}^{2}u_{in}\partial_{x}\partial_{y}u_{in}}{(\partial_{y}u_{in})^{2}}\y^{-1}
-\bar{\mu}\y^{-1}-\bar{\nu}\ln\y$, where
\begin{align*}
\mathfrak{R}^{x}
=& e^{-\bar{M}t}\y^{-1}\partial_{x}\mathfrak{R}_{\rho}
  -e^{-\bar{M}t}\y^{-1}\partial_{x}\mathfrak{L}_{0}\rho
  +e^{-\bar{M}t}\y^{-1}\partial_{y}^{-1}[\partial_{x}^{2}\mathfrak{U}]\partial_{y}\rho
  -\bar{\mu}\bar{N}e^{-\bar{M}t}\y^{-1}\\
 &-\bar{\mu}(1+\bar{N}t)e^{-\bar{M}t}(\mathfrak{L}_{1}-2\y^{-1})\y^{-2}
  -\bar{\mu}(1+\bar{N}t)e^{-\bar{M}t}\y^{-1}(\mathfrak{L}_{0}+\partial_{x}\mathfrak{U}-2\partial_{y}\rho+\y^{-1}\mathfrak{L}_{1})\\
 &+2e^{-\bar{M}t}\bar{\mu}(1+\bar{N}t)\y^{-3}
  -\bar{\nu}(\bar{M}+\mathfrak{L}_{0}+\partial_{x}\mathfrak{U}-2\partial_{y}\rho+\y^{-1}\mathfrak{L}_{1})(\ln\y+\mathfrak{m}t)\\
 &+\bar{\nu}(2\y^{-2}-\y^{-1}\mathfrak{L}_{1}-\mathfrak{m}).
\end{align*}
Let $\hat{t}_{3}=\min\{\frac{1}{\bar{N}},\hat{t}\}$.
For any $\bar{\nu}>0$, by taking $\bar{\mu}$ large enough but independent of $\mathcal{Z}$, one can has
\begin{align*}
\mathfrak{D}^{x}|_{t=0}<0,~
\mathfrak{D}^{x}|_{y=\hat{y}}<0,~\text{in}~\check{\Omega}_{\hat{t}_{3}}^{\hat{y}}.
\end{align*}
Also by the condition \eqref{Fd}, Lemma \ref{DI}--Lemma \ref{rhoP} and Corollary \eqref{wdotuC}, with the help of \eqref{rhoy} one can obtain
\begin{align*}
|\partial_{x}\mathfrak{R}_{\rho}|\leq C,~
|\partial_{x}\mathfrak{L}_{0}\rho|\leq C_{*},~
\partial_{y}^{-1}[\partial_{x}^{2}\mathfrak{U}]\partial_{y}\rho\leq C_{*}.
\end{align*}
Thus by taking $\bar{M}$ and $\bar{N}$ large enough one can get $\mathfrak{R}^{x}<0$ in $\check{\Omega}_{\hat{t}_{3}}^{\hat{y}}$. Noticing that for any $y>R_{\bar{\nu}}=e^{\frac{K+2\bar{\mu}}{\bar{\nu}}}$, $\mathfrak{D}^{x}<0$ in $\check{\Omega}_{\hat{t}_{3}}^{\hat{y}}$, by employing the classical maximum principle for the parabolic equations, one has
\begin{align*}
\mathfrak{D}^{x}<0,~\text{in}~\check{\Omega}_{\hat{t}_{3}}^{\hat{y}}\cap\hat{\Omega}_{\hat{t}_{3}}^{R_{\bar{\nu}}}.
\end{align*}
Taking the limit $\bar{\nu}\to 0^{+}$, one has $e^{-\bar{M}t}\y^{-1}[\partial_{x}\rho-\bar{\mu}(1+\bar{N}t)]<0$ in $\check{\Omega}_{\hat{t}_{3}}^{\hat{y}}$,
which gives
\begin{align}\label{rhox+}
\partial_{x}\rho\leq\bar{\mu}(1+\bar{N}t), ~\text{in}~\check{\Omega}_{\hat{t}_{3}}^{\hat{y}}.
\end{align}
A similar argument applied to the function
\begin{align*}
\mathfrak{D}_{-}^{x}=e^{-\bar{M}t}\y^{-1}[-\partial_{x}\rho-\bar{\mu}(1+\bar{N}t)]-\bar{\nu}\ln\y-\bar{\nu}\mathfrak{m}t.
\end{align*}
gives that
\begin{align}\label{rhox-}
-\partial_{x}\rho\leq\bar{\mu}(1+\bar{N}t), ~\text{in}~\check{\Omega}_{\hat{t}_{3}}^{\hat{y}}.
\end{align}
if $\bar{N}$ is large enough.
We omit the detail here. Combining \eqref{rhox+} with \eqref{rhox-}, one has
\begin{align}\label{rhox'}
|\partial_{x}\rho|\leq\bar{\mu}(1+\bar{N}t)\leq C\mathcal{P}, ~\text{in}~\check{\Omega}_{\hat{t}_{3}}^{\hat{y}},
\end{align}
Take $t_{3}\in(0,\hat{t}_{3}]$, \eqref{rhox} follows immediately.

\noindent{\bf \underline{Proof of \eqref{rhol}.}}
From \eqref{rhoe} one can see that
\begin{align}\label{rhoe'}
\partial_{t}\rho
= \mathfrak{R}_{\rho}
 -\mathfrak{U}\partial_{x}\rho
 -\mathfrak{L}_{1}\partial_{y}\rho
 -\mathfrak{L}_{0}\rho
 +\partial_{y}^{2}\rho,~\text{in}~\check{\Omega}_{\hat{t}_{3}}^{Y},
\end{align}
where $\mathfrak{L}_{0}, \mathfrak{L}_{1}$ and $\mathfrak{R}_{\rho}$ are given below \eqref{rhoe}.

Since
\begin{align*}
|\partial_{y}\hat{F}|
\leq C_{*}e^{-\frac{4}{3}\delta\y}
\leq C_{*}(U-\mathfrak{U}),~\text{in}~\check{\Omega}_{\hat{t}_{3}}^{Y},
\end{align*}
and
$$
|\eta'\partial_{x}^{2}\partial_{y}\mathbb{U}|+|\eta''\partial_{x}^{2}\mathbb{U}|+|\eta'\partial_{x}^{2}\mathbb{U}|
+|\eta\partial_{x}^{2}\partial_{y}^{2}\mathbb{U}|+|\eta\partial_{x}^{2}\partial_{y}\mathbb{U}|
\leq C_{*}e^{-(1+\theta)\delta\y}
\leq C_{*}(U-\mathfrak{U}),~\text{in}~\check{\Omega}_{\hat{t}_{3}}^{Y},
$$
by using Lemma \ref{sdotu} one has $|\mathfrak{R}_{\rho}|\leq C_{*}$ and $|\mathfrak{L}_{0}\rho|\leq C_{*}$ in $\check{\Omega}_{\hat{t}_{3}}^{Y}$. From Lemma \ref{rhoP} and \eqref{rhox'} one can deduce $|\mathfrak{U}\partial_{x}\rho|\leq C_{*}$ in $\check{\Omega}_{\hat{t}_{3}}^{Y}$.
From Lemma \eqref{rhoP} and \eqref{rhoy} it follows that $|\mathfrak{L}_{1}\partial_{y}\rho|\leq C_{*}$ in $\check{\Omega}_{\hat{t}_{3}}^{Y}$.
Finally, noticing that
$$\partial_{y}^{2}\rho=\frac{\partial_{y}^{3}\dot{\mathfrak{U}}}{\dot{\mathfrak{U}}}
                       -\frac{\rho\partial_{y}^{2}\dot{\mathfrak{U}}}{\dot{\mathfrak{U}}}
                       -2\rho\partial_{y}\rho,~\text{in}~\check{\Omega}_{\hat{t}_{3}}^{Y},$$
by using Lemma \ref{sdotu} and Lemma \ref{rhoP} one has $|\partial_{y}^{2}\rho|\leq C$ in $\check{\Omega}_{\hat{t}_{3}}^{Y}$. From \eqref{rhoe'} one knows that $|\partial_{t}\rho|\leq C_{*}$ in $\check{\Omega}_{\hat{t}_{3}}^{Y}$, thus there exists $t_{3}\in(0,\hat{t}_{3}]$ such that
\begin{align*}
\rho(t,x,y)
\leq&~\rho_{in}(x,y)+\int_{0}^{t}|\partial_{t}\rho(\sigma,x,y)|d\sigma\\
\leq&-\delta+C_{*}t
\leq -\frac{15}{16}\delta,~\text{in}~\check{\Omega}_{t_{3}}^{Y}.
\end{align*}

\end{proof}
~~~~~~~~~~

\subsection{High-order estimates of the approximate solutions}~

In the previous section, we established a uniform estimate for the corrected increment solution $u$ under certain conditions on $\mathbf{U}_{s}$ and $\mathbf{U}_{s+1}$. In order for the iteration to proceed, we need to establish a uniform estimate of $\mathbf{U}_{s}$ and $\mathbf{U}_{s+1}$ to guarantee that these conditions are satisfied throughout the iteration. This is the main purpose of this subsection.

The main assumption in this subsection is as follows.
\begin{assumption}\label{peA}
For some $\hat{t}_{*}\in(0,T]$, assume that
Assumption \ref{FA} holds for $\hat{t}_{*}$ such that  Lemma \ref{DI}--Lemma \ref{rhoP'} are valid on $[0,\hat{t}_{*}]$. Assume $u$ is a solution to the problem \eqref{app'} on the time interval $[0,\hat{t}_{*}]$ such that
\begin{align*}
 \|u\|_{L_{t}^{\infty}\mathcal{H}_{\psi,\tilde{\varphi}}^{s}(\Omega_{\hat{t}_{*}})}
+\|\partial_{y}\mathcal{U}_{s}\|_{L_{t}^{2}L_{\tilde{\varphi}}^{2}(\Omega_{\hat{t}_{*}})}
+\|\sqrt{\eta}\partial_{x}\mathcal{U}_{s}\|_{L_{t}^{2}L_{\tilde{\varphi}}^{2}(\Omega_{\hat{t}_{*}})}
+\|\tilde{\mathcal{P}}\mathcal{U}_{s}\|_{L_{t}^{2}L_{\tilde{\varphi}}^{2}(\Omega_{\hat{t}_{*}})}
\leq C_{*}\sqrt{\bar{\eta}},
\end{align*}
and $\mathfrak{U}$ is a solution to the problem \eqref{App+1} on the time interval $[0,\hat{t}_{*}]$ with $\hat{F}$ satisfying
$$|\partial_{y}\hat{F}|\leq C_{*}\iota^{-1}(\dot{\mathfrak{U}}+\bar{\varsigma})~\text{in}~\Omega_{\hat{t}_{*}},$$
and
$$\|\mathcal{P}(\partial_{y}^{2}\hat{F}-\rho\partial_{y}\hat{F})\|_{L^{\infty}(\Omega_{\hat{t}_{*}})}\leq C_{*},$$
such that $\mathbb{U}=\mathfrak{U}-u_{0}$ satisfies
$$\|\mathbb{U}\|_{L_{t}^{\infty}\hat{H}_{\psi}^{s}(\Omega_{\hat{t}_{*}})}\leq \mathcal{Z},
~\|\mathbb{U}\|_{L_{t}^{2}\hat{H}_{\psi,\y}^{s}(\Omega_{\hat{t}_{*}})}\leq \mathcal{Z}.$$
Moreover, assume that
\begin{align*}
 \|\sqrt{\tilde{\eta}}\partial_{x}\tilde{\mathbf{U}}_{s}\|_{L_{t}^{2}L_{\tilde{\varphi}}^{2}(\Omega_{\hat{t}_{*}})}
+\|\tilde{\mathcal{P}}\tilde{\mathbf{U}}_{s}\|_{L_{t}^{2}L_{\tilde{\varphi}}^{2}(\Omega_{\hat{t}_{*}})}
+\|\tilde{\mathbf{U}}_{s}\|_{L_{t}^{\infty}L_{\tilde{\varphi}}^{2}(\Omega_{\hat{t}_{*}})}
\leq \mathcal{Z},
\end{align*}
and
\begin{align*}
 \|\sqrt{\tilde{\eta}}\partial_{x}\tilde{\mathbf{U}}_{s+1}\|_{L_{t}^{2}L_{\hat{\tilde{\varphi}}}^{2}(\Omega_{\hat{t}_{*}})}
+\|\tilde{\mathcal{P}}\tilde{\mathbf{U}}_{s+1}\|_{L_{t}^{2}L_{\hat{\tilde{\varphi}}}^{2}(\Omega_{\hat{t}_{*}})}
+\|\tilde{\mathbf{U}}_{s+1}\|_{L_{t}^{\infty}L_{\hat{\tilde{\varphi}}}^{2}(\Omega_{\hat{t}_{*}})}
\leq \mathcal{Z},
\end{align*}
where $\hat{\tilde{\varphi}}=\y^{-\frac{1}{2}}\tilde{\varphi}$.
\end{assumption}

The main conclusion of this subsection can be summarized as follows.
\begin{theorem}\label{Phe}
Let $\theta^{*}$ be given in Lemma \ref{wdotu} and $t_{3}$ be given in Lemma \ref{rhoP'}. Assume that the assumption of Theorem \ref{MR} holds.
There exists a positive constant $\lambda^{*}$ such that for any $\lambda>\lambda^{*}$, there are positive constants $\Lambda^{*}, \iota_{*}$ and $\mathcal{Z}^{*}$ depending only on $\lambda$ for which, when $\Lambda\geq\Lambda^{*}$, $\theta\geq\theta^{*}$, $\iota\in(0,\iota_{*})$ and $\mathcal{Z}\geq\mathcal{Z}^{*}$, there exists $t_{*}\in(0,t_{3}]$ such that if Assumption \ref{peA} holds for $t_{*}$, then we have
\begin{align}\label{Hoe1}
&\|\mathbf{U}_{s}\|_{L_{t}^{\infty}L_{\varphi}^{2}(\Omega_{t_{*}})}
 +\Lambda\delta\|\sqrt{\y}\mathbf{U}_{s}\|_{L_{t}^{2}L_{\varphi}^{2}(\Omega_{t_{*}})}
 +\|\mathcal{P}\mathbf{U}_{s}\|_{L_{t}^{2}L_{\varphi}^{2}(\Omega_{t_{*}})}\\
&+\|\partial_{y}\mathbf{U}_{s}\|_{L_{t}^{2}L_{\varphi}^{2}(\Omega_{t_{*}})}
 +\|\sqrt{\eta}\partial_{x}\mathbf{U}_{s}\|_{L_{t}^{2}L_{\varphi}^{2}(\Omega_{t_{*}})}
\leq \mathcal{Z},\nonumber
\end{align}
and
\begin{align}\label{Hoe2}
&\|\mathbf{U}_{s+1}\|_{L_{t}^{\infty}L_{\hat{\varphi}}^{2}(\Omega_{t_{*}})}
 +\Lambda\delta\|\sqrt{\y}\mathbf{U}_{s+1}\|_{L_{t}^{2}L_{\hat{\varphi}}^{2}(\Omega_{t_{*}})}
 +\|\mathcal{P}\mathbf{U}_{s+1}\|_{L_{t}^{2}L_{\hat{\varphi}}^{2}(\Omega_{t_{*}})}\\
&+\|\partial_{y}\mathbf{U}_{s+1}\|_{L_{t}^{2}L_{\hat{\varphi}}^{2}(\Omega_{t_{*}})}
 +\|\sqrt{\eta}\partial_{x}\mathbf{U}_{s+1}\|_{L_{t}^{2}L_{\hat{\varphi}}^{2}(\Omega_{t_{*}})}
\leq \mathcal{Z}.\nonumber
\end{align}
\end{theorem}
\begin{proof}[Proof of Theorem \ref{Phe}]

We shall prove that the inequalities \eqref{Hoe1} and \eqref{Hoe2} holds for some $t_{*}\in (0,t_{3}]$, under the assumption that Assumption \ref{peA} holds for $t_{3}$, with $\lambda,\Lambda, \theta, \mathcal{Z}$ sufficiently large and $\iota$ sufficiently small. As can be seen from our proof, the inequalities \eqref{Hoe1} and \eqref{Hoe2} still hold for $t_{*}$ when Assumption \ref{peA} holds for $t_{*}$.

Let $r=s$ or $s+1$. From \eqref{App+1} one knows that $\mathbf{U}_{r}(t,x,y)$ obeys the following problem
\begin{equation}\label{Un}
\begin{cases}
\partial_{t}\mathbf{U}_{r}
+\mathbf{L}_{1}\partial_{x}\mathbf{U}_{r}
+\mathbf{L}_{2}\partial_{y}\mathbf{U}_{r}
+\mathbf{L}_{0}\mathbf{U}_{r}
-\partial_{y}^{2}\mathbf{U}_{r}
-\eta\partial_{x}^{2}\mathbf{U}_{r}
=\sum\limits_{i=1}^{6}\mathcal{T}_{i},\\
\left(\partial_{y}\mathbf{U}_{r}+\rho\mathbf{U}_{r}\right)|_{y=0}=\partial_{x}^{r+1}P,~
\lim\limits_{y\to+\infty} \mathbf{U}_{r}=0,\\
\mathbf{U}_{r}|_{t=0}=\mathcal{P}(0,x,y)(\partial_{y}\partial_{x}^{s}\tilde{u}_{in}-\rho(0,x,y)\partial_{x}^{s}\tilde{u}_{in}),
\end{cases}
\end{equation}
where
\begin{align*}
\mathbf{L}_{0}
=&~r\partial_{x}\mathfrak{U}-2\partial_{y}\rho
  -2\left(\frac{\partial_{y}\mathcal{P}}{\mathcal{P}}\right)^{2}-2\eta\left(\frac{\partial_{x}\mathcal{P}}{\mathcal{P}}\right)^{2}
  +\frac{3\chi}{4\mathcal{P}\sqrt{\dot{\mathfrak{U}}+\varsigma}}\left[\rho^{2}+\eta\left(\frac{\partial_{x}\dot{\mathfrak{U}}}{\dot{\mathfrak{U}}+\varsigma}\right)^{2}\right]\\
 &-\frac{\partial_{y}^{-1}\left[\partial_{x}\mathfrak{U}\right]\varsigma'+\varsigma''-\partial_{y}\eta\partial_{x}^{2}\mathbb{U}+\eta\partial_{x}^{2}\partial_{y}u_{0}+\partial_{y}\hat{F}}{2(\dot{\mathfrak{U}}+\varsigma)}\\
 &+(1-\chi)\frac{\partial_{y}^{-1}\left[\partial_{x}\mathfrak{U}\right]\varsigma'+\varsigma''-\partial_{y}\eta\partial_{x}^{2}\mathbb{U}+\eta\partial_{x}^{2}\partial_{y}u_{0}+\partial_{y}\hat{F}}{2\mathcal{P}(\dot{\mathfrak{U}}+\varsigma)}\\
 &-\frac{-\partial_{y}^{-1}[\partial_{x}\mathfrak{U}]\chi'+\chi'\rho-\chi''}{\mathcal{P}\sqrt{\dot{\mathfrak{U}}+\varsigma}}
  -\frac{\partial_{y}^{-1}\left[\partial_{x}\mathfrak{U}\right]\chi'+\chi''}{\mathcal{P}},
\end{align*}
$$\mathbf{L}_{1}=\mathfrak{U}+2\eta\frac{\partial_{x}\mathcal{P}}{\mathcal{P}},
~\mathbf{L}_{2}=2\frac{\partial_{y}\mathcal{P}}{\mathcal{P}}-\partial_{y}^{-1}\left[\partial_{x}\mathfrak{U}\right],$$
and $\mathcal{T}_{1}=\mathcal{P}(-\partial_{y}\bar{\mathcal{T}}_{1}+\rho\bar{\mathcal{T}}_{1})$, with
$$
\bar{\mathcal{T}}_{1}(t,x,y)=
 \sum_{i=2}^{r-1}\binom{r}{i}\partial_{x}^{r-i+1}\mathbb{U}\partial_{x}^{i}\mathbb{U}
-\sum_{i=1}^{r-2}\binom{r}{i}\partial_{y}^{-1}[\partial_{x}^{i+1}\mathbb{U}]\partial_{x}^{r-i}\partial_{y}\mathbb{U},
$$
and
$$
\mathcal{T}_{2}(t,x,y)
= \mathcal{P}\left(-\partial_{y}^{-1}[\partial_{x}^{r+1}\tilde{\mathfrak{U}}](\varsigma'-\rho\varsigma)
 +2\eta\partial_{x}\rho\partial_{x}^{r+1}\mathbb{U}+\eta'\partial_{x}^{r+2}\mathbb{U}\right),
$$
$$
\mathcal{T}_{3}(t,x,y)
=r\mathcal{P}\left(\partial_{x}\partial_{y}^{2}\mathfrak{U}-\rho\partial_{x}\partial_{y}\mathfrak{U}\right)\partial_{y}^{-1}[\partial_{x}^{r}\mathfrak{U}],
$$
and $\mathcal{T}_{4}(t,x,y)=\mathcal{P}\bar{\mathcal{T}}_{4}\frac{\partial_{x}^{r}\mathbb{U}}{\dot{\mathfrak{U}}+\varsigma}$, with
\begin{align*}
\bar{\mathcal{T}}_{4}(t,x,y)
=&\partial_{y}^{2}\hat{F}+\varsigma'\partial_{x}\mathfrak{U}+\varsigma'''+\partial_{y}^{-1}[\partial_{x}\mathfrak{U}]\varsigma''+\eta\partial_{x}^{2}\partial_{y}^{2}u_{0}-\varsigma\partial_{x}\partial_{y}\mathfrak{U}\\
 &-\rho\left(\partial_{y}\hat{F}+\partial_{y}^{-1}[\partial_{x}\mathfrak{U}]\varsigma'+\varsigma''+\eta\partial_{x}^{2}\partial_{y}u_{0}\right)
  -2\eta\partial_{x}\partial_{y}\mathfrak{U}\partial_{x}\rho,
\end{align*}
and
$$
\mathcal{T}_{5}=\mathcal{P}(\varsigma\partial_{x}^{r+1}u-\partial_{y}\partial_{x}^{r}\mathbb{F}+\rho\partial_{x}^{r}\mathbb{F}),
$$
with $\mathbb{F}=-u\partial_{x}u+(1-h)\partial_{y}^{-1}[\partial_{x}u]\partial_{y}u$, and finally, $\mathcal{T}_{6}=\mathcal{P}\bar{\mathcal{T}}_{6}$, with
\begin{align*}
\bar{\mathcal{T}}_{6}
=& \rho\sum_{i=2}^{r}\binom{r}{i}\partial_{x}^{r-i+1}\mathbb{U}\partial_{x}^{i}u_{0}
  +\rho\sum_{i=0}^{r-1}\binom{r}{i}\partial_{x}^{r-i+1}u_{0}\partial_{x}^{i}\mathbb{U}
  -\sum_{i=0}^{r}\binom{r}{i}\partial_{y}\partial_{x}^{r-i+1}u_{0}\partial_{x}^{i}u_{0}\\
 &-\sum_{i=0}^{r-1}\binom{r}{i}\partial_{y}\partial_{x}^{r-i+1}u_{0}\partial_{x}^{i}\mathbb{U}
  -\sum_{i=2}^{r}\binom{r}{i}\partial_{y}\partial_{x}^{r-i+1}\mathbb{U}\partial_{x}^{i}u_{0}\\
 &+\sum_{i=0}^{r-2}\binom{r}{i}\partial_{y}^{-1}[\partial_{x}^{i+1}\mathbb{U}]\partial_{x}^{r-i}\partial_{y}^{2}u_{0}
  +\sum_{i=1}^{r-2}\binom{r}{i}\partial_{y}^{-1}[\partial_{x}^{i+1}u_{0}]\partial_{x}^{r-i}\partial_{y}^{2}\mathbb{U}
  +\sum_{i=0}^{r-2}\binom{r}{i}\partial_{y}^{-1}[\partial_{x}^{i+1}u_{0}]\partial_{x}^{r-i}\partial_{y}^{2}u_{0}\\
 &-\rho\sum_{i=0}^{r-2}\binom{r}{i}\partial_{y}^{-1}[\partial_{x}^{i+1}\mathbb{U}]\partial_{x}^{r-i}\partial_{y}u_{0}
  -\rho\sum_{i=1}^{r-2}\binom{r}{i}\partial_{y}^{-1}[\partial_{x}^{i+1}u_{0}]\partial_{x}^{r-i}\partial_{y}\mathbb{U}
  -\rho\sum_{i=0}^{r-2}\binom{r}{i}\partial_{y}^{-1}[\partial_{x}^{i+1}u_{0}]\partial_{x}^{r-i}\partial_{y}u_{0}\\
 &+\rho\partial_{x}^{r}(\partial_{t}u_{0}+u_{0}\partial_{x}u_{0}+\partial_{x}P)
  +\partial_{y}^{3}\partial_{x}^{r}u_{0}
  -\rho\partial_{x}^{r}\partial_{y}^{2}u_{0}
  -\partial_{t}\partial_{y}\partial_{x}^{r}u_{0}.
\end{align*}

Let $\varphi_{r}=\varphi, \hat{\Psi}_{r}=\hat{\Psi}=\y^{-\frac{1}{2}}\Psi$ when $r=s$ and  $\varphi_{r}=\hat{\varphi},\hat{\Psi}_{r}=\y^{-1}\Psi$ when $r=s+1$, where $\Psi$ is given below Corollary \eqref{ncC1}.
Correspondingly, let $\tilde{\varphi}_{r}=\tilde{\varphi}, \tilde{\hat{\Psi}}_{r}=\y^{-\frac{1}{2}}\tilde{\mathcal{P}}^{2}(t,x,y)\varphi_{o}(t,y)$ when $r=s$ and  $\tilde{\varphi}_{r}=\y^{-\frac{1}{2}}\tilde{\varphi},\tilde{\hat{\Psi}}_{r}=\y^{-1}\tilde{\mathcal{P}}^{2}(t,x,y)\varphi_{o}(t,y)$ when $r=s+1$.
For any $t\in[0,t_{3}]$, multiplying $\varphi_{r}^{2}\mathbf{U}_{r}$ and then integrating over $\Omega$, by using the integration by parts we get
\begin{align}\label{UnE}
&\frac{1}{2}\frac{d}{dt}\|\mathbf{U}_{r}\|_{L_{\varphi_{r}}^{2}(\Omega)}^{2}
+\|\partial_{y}\mathbf{U}_{r}\|_{L_{\varphi_{r}}^{2}(\Omega)}^{2}
+\|\sqrt{\eta}\partial_{x}\mathbf{U}_{r}\|_{L_{\varphi_{r}}^{2}(\Omega)}^{2}\\
=&~\frac{1}{2}(\mathbf{U}_{r},(\partial_{t}\varphi_{r}^{2}+\partial_{y}^{2}\varphi_{r}^{2})\mathbf{U}_{r})
  -(\mathbf{L}_{1}\partial_{x}\mathbf{U}_{r},\varphi_{r}^{2}\mathbf{U}_{r})
  -(\mathbf{L}_{2}\partial_{y}\mathbf{U}_{r},\varphi_{r}^{2}\mathbf{U}_{r})
  -(\mathbf{L}_{0}\mathbf{U}_{r},\varphi_{r}^{2}\mathbf{U}_{r})\nonumber\\
 &+(\mathbf{U}_{r},\varphi_{r}\partial_{y}\varphi_{r}\mathbf{U}_{r})_{L_{x}^{2}}|_{y=0}
  -(\partial_{y}\mathbf{U}_{r},\varphi_{r}^{2}\mathbf{U}_{r})_{L_{x}^{2}}|_{y=0}
  +\sum_{i=1}^{6}(\mathcal{T}_{i},\varphi_{r}^{2}\mathbf{U}_{r})\nonumber\\
:=&\sum_{i=0}^{11}J_{i}.\nonumber
\end{align}

To get the a posterior estimates of the approximate solutions, we shall estimate $J_{1}$, $J_{2}$,...,$J_{13}$ term by term.
In what follows, we shall always take $\lambda, \Lambda, \theta, \iota, \mathcal{Z}$ to be positive constants to be determined later, satisfying $\lambda>\lambda_{*}, \Lambda>\Lambda_{*}, \theta>\theta^{*}, \iota\in(0,\iota^{*}), \mathcal{Z}>\mathcal{Z}_{*}$, where $\lambda_{*}, \Lambda_{*},\iota^{*}, \mathcal{Z}_{*}$ are given in Theorem \ref{wpae}, and $\theta^{*}$ is given in Lemma \ref{wdotu}.

\noindent{\bf \underline{Estimate of~$J_{0}$.}}
When $r=s$, parallel to the estimate of $I_{0}$ in Section 4.4, one has
\begin{align}\label{J00}
J_{0}\leq -\frac{\lambda}{64}\|\sqrt{\omega_{\lambda}}\mathbf{U}_{s}\|_{L_{\varphi}^{2}(\Omega)}^{2}
          -\frac{3}{4}\Lambda\delta\|\sqrt{\y}\mathbf{U}_{s}\|_{L_{\varphi}^{2}(\Omega)}^{2}.
\end{align}
When $r=s+1$, it follows from Lemma \ref{weighte} that
\begin{align}\label{J01}
J_{0}
=&~\frac{1}{2}(\mathbf{U}_{s+1},(\partial_{t}\varphi^{2}+\partial_{y}^{2}\varphi^{2})\y^{-1}\mathbf{U}_{s+1})
  +(\mathbf{U}_{s+1},\partial_{y}\varphi^{2}\partial_{y}\y^{-1}\mathbf{U}_{s+1})
  +\frac{1}{2}(\mathbf{U}_{s+1},\varphi^{2}\partial_{y}^{2}\y^{-1}\mathbf{U}_{s+1})\\
\leq&-\frac{\lambda}{64}\|\sqrt{\omega_{\lambda}}\mathbf{U}_{s+1}\|_{L_{\hat{\varphi}}^{2}(\Omega)}^{2}
     -\frac{3}{4}\Lambda\delta\|\sqrt{\y}\mathbf{U}_{s+1}\|_{L_{\hat{\varphi}}^{2}(\Omega)}^{2}
     +C\lambda(\mathbf{U}_{s+1},\sqrt{\omega_{\lambda}}\hat{\varphi}^{2}\mathbf{U}_{s+1})
     +C_{*}\|\mathbf{U}_{s+1}\|_{L_{\hat{\varphi}}^{2}(\Omega)}^{2}\nonumber\\
\leq&-\left(\frac{\lambda}{64}-\frac{1}{16}\right)\|\sqrt{\omega_{\lambda}}\mathbf{U}_{s+1}\|_{L_{\hat{\varphi}}^{2}(\Omega)}^{2}
     -\frac{3}{4}\Lambda\delta\|\sqrt{\y}\mathbf{U}_{s+1}\|_{L_{\hat{\varphi}}^{2}(\Omega)}^{2}
     +C_{*}\|\mathbf{U}_{s+1}\|_{L_{\hat{\varphi}}^{2}(\Omega)}^{2}.\nonumber
\end{align}
As a consequence of \eqref{J00} and \eqref{J01}, one has
\begin{align*}
J_{0}
\leq -\left(\frac{\lambda}{64}-\frac{1}{16}\right)\|\sqrt{\omega_{\lambda}}\mathbf{U}_{r}\|_{L_{\varphi_{r}}^{2}(\Omega)}^{2}
     -\frac{3}{4}\Lambda\delta\|\sqrt{\y}\mathbf{U}_{r}\|_{L_{\varphi_{r}}^{2}(\Omega)}^{2}
     +C_{*}\|\mathbf{U}_{r}\|_{L_{\varphi_{r}}^{2}(\Omega)}^{2}.
\end{align*}

\noindent{\bf \underline{Estimate of~$J_{1}$.}}
By using the integration by parts and Lemma \ref{rhoP}, one has
\begin{align*}
J_{1}
=\frac{1}{2}(\partial_{x}\mathbf{L}_{1}\mathbf{U}_{r},\varphi_{r}^{2}\mathbf{U}_{r})
=&~\frac{1}{2}(\partial_{x}\mathfrak{U}\mathbf{U}_{r},\varphi_{r}^{2}\mathbf{U}_{r})
  +\left(\eta\partial_{x}\left(\frac{\partial_{x}\mathcal{P}}{\mathcal{P}}\right)\mathbf{U}_{r},\varphi_{r}^{2}\mathbf{U}_{r}\right)\\
=&~\frac{1}{2}(\partial_{x}\mathfrak{U}\mathbf{U}_{r},\varphi_{r}^{2}\mathbf{U}_{r})
  -\frac{1}{2}\int_{\mathbb{T}}\int_{\check{y}_{*}}^{\hat{y}_{*}}\eta\varphi_{r}^{2}\left[\frac{\partial_{x}^{2}\dot{\mathfrak{U}}}{\dot{\mathfrak{U}}+\varsigma}-\left(\frac{\partial_{x}\dot{\mathfrak{U}}}{\dot{\mathfrak{U}}+\varsigma}\right)^{2}\right]\mathbf{U}_{r}^{2}dydx\\
 &+\int_{\mathbb{T}}\int_{[\check{y},\hat{y}]\setminus[\check{y}_{*},\hat{y}_{*}]}\eta\varphi_{r}^{2}\partial_{x}\left(\frac{\partial_{x}\mathcal{P}}{\mathcal{P}}\right)\mathbf{U}_{r}^{2}dydx\\
\leq&~C_{*}\|\mathbf{U}_{r}\|_{L_{\varphi_{r}}^{2}(\Omega)}^{2}.
\end{align*}

\noindent{\bf \underline{Estimate of~$J_{2}$.}}
By the definition of $\mathbf{L}_{2}$, one can decompose $J_{2}$ into two parts,
\begin{align*}
J_{2}
= (\partial_{y}^{-1}[\partial_{x}]\partial_{y}\mathfrak{U}\mathbf{U}_{r},\varphi_{r}^{2}\mathbf{U}_{r})
 -2\left(\frac{\partial_{y}\mathcal{P}}{\mathcal{P}} \partial_{y}\mathbf{U}_{r},\varphi_{r}^{2}\mathbf{U}_{r}\right)
:=J_{21}+J_{22}.
\end{align*}
The estimate of $J_{21}$ is parallel to that of $I_{21}$ in Section 4.6, and one has
\begin{align*}
J_{21}
=&~-\frac{1}{2}(\partial_{x}\mathfrak{U}\mathbf{U}_{r},\varphi_{r}^{2}\mathbf{U}_{r})
  -\frac{1}{2}(\partial_{y}^{-1}\left[\partial_{x}\mathfrak{U}\right]\mathbf{U}_{r},\partial_{y}\varphi_{r}^{2}\mathbf{U}_{r})\\
=&~-\frac{1}{2}(\partial_{x}\mathfrak{U}\mathbf{U}_{r},\varphi_{r}^{2}\mathbf{U}_{r})
  -\frac{1}{2}(\chi\partial_{y}^{-1}\left[\partial_{x}\mathfrak{U}\right]\mathbf{U}_{r},\partial_{y}\varphi_{r}^{2}\mathbf{U}_{r})
  -\frac{1}{2}((1-\chi)\partial_{y}^{-1}\left[\partial_{x}\mathfrak{U}\right]\mathbf{U}_{r},\partial_{y}\varphi_{r}^{2}\mathbf{U}_{r})\\
\leq&~C_{*}\|\mathbf{U}_{r}\|_{L_{\varphi_{r}}^{2}(\Omega)}^{2}
     +C_{*}\lambda(\mathbf{U}_{r},\sqrt{\omega_{\lambda}}\varphi_{r}^{2}\mathbf{U}_{r})
     +C_{*}(\mathbf{U}_{r},\varphi_{r}^{2}\mathbf{U}_{r})
     +C\|\y^{-1}\partial_{y}^{-1}\left[\partial_{x}\mathfrak{U}\right]\|_{L^{\infty}(\Omega)}
       \|\sqrt{\y}\mathbf{U}_{r}\|_{L_{\varphi_{r}}^{2}(\Omega)}^{2}\\
\leq&~C_{*}\|\mathbf{U}_{r}\|_{L_{\varphi_{r}}^{2}(\Omega)}^{2}
     +\frac{1}{16}\|\sqrt{\omega_{\lambda}}\mathbf{U}_{r}\|_{L_{\varphi_{r}}^{2}(\Omega)}^{2}
     +C_{*}\|\sqrt{\y}\mathbf{U}_{r}\|_{L_{\varphi_{r}}^{2}(\Omega)}^{2}.
\end{align*}
Noticing ${\bf supp}~\partial_{y}\mathcal{P}\subset(\check{y},\hat{y})$ and $\mathcal{P}^{-1}\partial_{y}\mathcal{P}=-\frac{1}{2}\rho$ on $[\check{y}_{*},\hat{y}_{*}]$, one has
\begin{align*}
J_{22}
=&~\int_{\mathbb{T}}\int_{\check{y}_{*}}^{\hat{y}_{*}}\varphi_{r}^{2}\rho\mathbf{U}_{r}\partial_{y}\mathbf{U}_{r}dydx
  -2\int_{\mathbb{T}}\int_{[\check{y},\hat{y}]\setminus[\check{y}_{*},\hat{y}_{*}]}\varphi_{r}^{2}\frac{\partial_{y}\mathcal{P}}{\mathcal{P}}\mathbf{U}_{r}\partial_{y}\mathbf{U}_{r}dydx\\
\leq&~C_{*}\|\mathbf{U}_{r}\|_{L_{\varphi_{r}}^{2}(\Omega)}^{2}
     +C\|\rho\mathbf{U}_{r}\|_{L_{\varphi_{r}}^{2}(\Omega)}^{2}
     +\frac{1}{16}\|\partial_{y}\mathbf{U}_{r}\|_{L_{\varphi_{r}}^{2}(\Omega)}^{2}\\
\leq&~C_{*}\|\mathbf{U}_{r}\|_{L_{\varphi_{r}}^{2}(\Omega)}^{2}
     +C\|\mathcal{P}\mathbf{U}_{r}\|_{L_{\varphi_{r}}^{2}(\Omega)}^{2}
     +\frac{1}{16}\|\partial_{y}\mathbf{U}_{r}\|_{L_{\varphi_{r}}^{2}(\Omega)}^{2}.
\end{align*}
To sum up, one has
\begin{align*}
J_{2}
\leq&~C_{*}\|\mathbf{U}_{r}\|_{L_{\varphi_{r}}^{2}(\Omega)}^{2}
     +\frac{1}{16}\|\sqrt{\omega_{\lambda}}\mathbf{U}_{r}\|_{L_{\varphi_{r}}^{2}(\Omega)}^{2}
     +C_{*}\|\sqrt{\y}\mathbf{U}_{r}\|_{L_{\varphi_{r}}^{2}(\Omega)}^{2}\\
    &+C\|\mathcal{P}\mathbf{U}_{r}\|_{L_{\varphi_{r}}^{2}(\Omega)}^{2}
     +\frac{1}{16}\|\partial_{y}\mathbf{U}_{r}\|_{L_{\varphi_{r}}^{2}(\Omega)}^{2}.
\end{align*}

\noindent{\bf \underline{Estimate of~$J_{3}$.}}
As is shown in \eqref{I33}, one has
\begin{align}\label{J31}
&\frac{\partial_{y}^{-1}\left[\partial_{x}\dot{\mathfrak{U}}\right]\varsigma'+\varsigma''-\partial_{y}\eta\partial_{x}^{2}\mathbb{U}+\eta\partial_{x}^{2}\partial_{y}u_{0}}{2(\dot{\mathfrak{U}}+\varsigma)}
 -(1-\chi)\frac{\partial_{y}^{-1}\left[\partial_{x}\mathfrak{U}\right]\varsigma'+\varsigma''-\partial_{y}\eta\partial_{x}^{2}\mathbb{U}+\eta\partial_{x}^{2}\partial_{y}u_{0}}{2\mathcal{P}(\dot{\mathfrak{U}}+\varsigma)}\\
&\frac{-\partial_{y}^{-1}[\partial_{x}\mathfrak{U}]\chi'+\chi'\rho-\chi''}{\mathcal{P}\sqrt{\dot{\mathfrak{U}}+\varsigma}}
+\frac{\partial_{y}^{-1}\left[\partial_{x}\mathfrak{U}\right]\chi'+\chi''}{\mathcal{P}}
\leq C_{*}\iota^{-1}.\nonumber
\end{align}
Since $|\partial_{y}\hat{F}|\leq C_{*}\iota^{-1}(\dot{\mathfrak{U}}+\bar{\varsigma})$ and $|\partial_{y}\rho|\leq C\mathcal{P}^{2}$, one has
\begin{align}\label{J32}
2\partial_{y}\rho-r\partial_{x}\mathfrak{U}
+\frac{\partial_{y}\hat{F}}{2(\dot{\mathfrak{U}}+\varsigma)}
-(1-\chi)\frac{\partial_{y}\hat{F}}{2\mathcal{P}(\dot{\mathfrak{U}}+\varsigma)}\leq C\mathcal{P}^{2}+C_{*}\iota^{-1}.
\end{align}
Combining \eqref{J31} and \eqref{J32} with \eqref{I32}, it follows that
\begin{align*}
-\mathbf{L}_{0}\leq C\mathcal{P}^{2}+C_{*}\iota^{-1},
\end{align*}
which gives
\begin{align*}
J_{3}
\leq C_{*}\iota^{-1}\|\mathbf{U}_{r}\|_{L_{\varphi_{r}}^{2}(\Omega)}^{2}
    +C\|\mathcal{P}\mathbf{U}_{r}\|_{L_{\varphi_{r}}^{2}(\Omega)}^{2}.
\end{align*}

\noindent{\bf \underline{Estimate of~$J_{4}$.}}
Using the trace theorem, we know that
\begin{align*}
J_{4}
=(\mathbf{U}_{r},\varphi_{r}\partial_{y}\varphi_{r}\mathbf{U}_{r})_{L_{x}^{2}}|_{y=0}
\leq C\|\mathbf{U}_{r}\|_{L_{\varphi_{r}}^{2}(\Omega)}^{2}
     +\frac{1}{16}\|\partial_{y}\mathbf{U}_{r}\|_{L_{\varphi_{r}}^{2}(\Omega)}^{2}.
\end{align*}

\noindent{\bf \underline{Estimate of~$J_{5}$.}}
With the help of the boundary condition in \eqref{Un}, by using the trace estimate one can deduce that
\begin{align*}
J_{5}
=-(\varphi_{r}^{2}\partial_{y}\mathbf{U}_{r},\mathbf{U}_{r})_{L_{x}^{2}}|_{y=0}
=& \varphi_{0}^{2}(\rho \mathbf{U}_{r},\mathbf{U}_{r})_{L_{x}^{2}}|_{y=0}
  -\left.\varphi_{0}^{2}(\partial_{x}^{r+1}P,\mathbf{U}_{r})_{L_{x}^{2}}\right|_{y=0}\\
\leq& C_{*}\|\mathbf{U}_{r}\|_{L_{\varphi_{r}}^{2}(\Omega)}^{2}
     +\frac{1}{16}\|\partial_{y}\mathbf{U}_{r}\|_{L_{\varphi_{r}}^{2}(\Omega)}^{2}
     +C\|\partial_{x}^{r+1}P\|_{L_{x}^{2}(\mathbb{T})}^{2},
\end{align*}
where $\varphi_{0}=\varphi(t,0)$.

\noindent{\bf \underline{Estimate of~$J_{6}$.}}
We shall estimate $J_{6}$ in different ways for $r=s$ and $r=s+1$ respectively.

If $r=s$, in view of Lemma \ref{sdotu} and Lemma \ref{rhoP} one has
\begin{align}\label{J6s1}
\left(\mathcal{P}\rho\bar{\mathcal{T}}_{1},\varphi^{2}\mathbf{U}_{s}\right)
=&~\int_{\mathbb{T}}\int_{\check{y}}^{\hat{y}}\varphi^{2}\mathcal{P}\rho\bar{\mathcal{T}}_{1}\mathbf{U}_{s}dydx
  +\int_{\mathbb{T}}\int_{\mathbb{R}_{+}\setminus[\check{y},\hat{y}]}\varphi^{2}\mathcal{P}\rho\bar{\mathcal{T}}_{1}\mathbf{U}_{s}dydx\\
\leq&~C_{*}\|\bar{\mathcal{T}}_{1}\|_{L^{2}_{x}(\mathbb{T})L^{\infty}_{y}[\check{y},\hat{y}]}
      \|\mathcal{P}\mathbf{U}_{s}\|_{L_{\varphi}^{2}(\Omega)}\|\mathcal{P}\varphi\|_{L_{y}^{2}[\check{y},\hat{y}]}
     +C_{*}\|\bar{\mathcal{T}}_{1}\|_{L_{\psi}^{2}(\Omega)}\|\mathbf{U}_{s}\|_{L_{\varphi}^{2}(\Omega)}\nonumber\\
\leq&~C_{*}\|\bar{\mathcal{T}}_{1}\|_{L^{2}_{x}(\mathbb{T})L^{\infty}_{y}[\check{y},\hat{y}]}
      \|\mathcal{P}\mathbf{U}_{s}\|_{L_{\varphi}^{2}(\Omega)}
      \left(\int_{\check{y}}^{\hat{y}}\frac{1}{((y-y_{*})^2+\bar{\varsigma}+t)^{1+\varepsilon_{0}}}dy\right)^{\frac{1}{2}}\nonumber\\
    &+C_{*}\|\bar{\mathcal{T}}_{1}\|_{L_{\psi}^{2}(\Omega)}\|\mathbf{U}_{s}\|_{L_{\varphi}^{2}(\Omega)}\nonumber\\
\leq&~C_{*}\|\mathbf{U}_{s}\|_{L_{\varphi}^{2}(\Omega)}^{2}
     +\frac{1}{2}\|\mathcal{P}\mathbf{U}_{s}\|_{L_{\varphi}^{2}(\Omega)}^{2}
     +\frac{C_{*}}{(\bar{\varsigma}+t)^{\frac{1}{2}+\varepsilon_{0}}}\|\bar{\mathcal{T}}_{1}\|_{L_{x}^{2}(\mathbb{T})L_{y}^{\infty}[\check{y},\hat{y}]}^{2}
     +C_{*}\|\bar{\mathcal{T}}_{1}\|_{L_{\psi}^{2}(\Omega)}^{2}.\nonumber
\end{align}
and
\begin{align}\label{J6s2}
\left(-\mathcal{P}\partial_{y}\bar{\mathcal{T}}_{1},\varphi^{2}\mathbf{U}_{s}\right)
\leq \frac{C_{*}}{(\bar{\varsigma}+t)^{\varepsilon_{0}}}\|\partial_{y}\bar{\mathcal{T}}_{1}\|_{L_{\psi}^{2}(\Omega)}^{2}
    +\frac{1}{2}\|\mathcal{P}\mathbf{U}_{s}\|_{L_{\varphi}^{2}(\Omega)}^{2}.
\end{align}
Recall the definition of $\bar{\mathcal{T}}_{1}$, using Sobolev imbedding inequality one has
\begin{align}\label{J6s3}
&\|\bar{\mathcal{T}}_{1}\|_{L_{x}^{2}(\mathbb{T})L_{y}^{\infty}[\check{y},\hat{y}]}^{2}
+\|\bar{\mathcal{T}}_{1}\|_{L_{\psi}^{2}(\Omega)}^{2}
+\|\partial_{y}\bar{\mathcal{T}}_{1}\|_{L_{\psi}^{2}(\Omega)}^{2}\\
\leq&~C\sum_{i=2}^{s-1}\|\partial_{x}^{s-i+1}\mathbb{U}\partial_{x}^{i}\mathbb{U}\|_{L_{x}^{2}(\mathbb{T})L_{y}^{\infty}[\check{y},\hat{y}]}^{2}
     +C\sum_{i=1}^{s-2}\|\partial_{y}^{-1}[\partial_{x}^{i+1}\mathbb{U}]\partial_{x}^{s-i}\partial_{y}\mathbb{U}\|_{L_{x}^{2}(\mathbb{T})L_{y}^{\infty}[\check{y},\hat{y}]}^{2}\nonumber\\
    &+C\sum_{i=2}^{s-1}\|\partial_{x}^{s-i+1}\mathbb{U}\partial_{x}^{i}\mathbb{U}\|_{L_{\psi}^{2}(\Omega)}^{2}
     +C\sum_{i=1}^{s-2}\|\partial_{y}^{-1}[\partial_{x}^{i+1}\mathbb{U}]\partial_{x}^{s-i}\partial_{y}\mathbb{U}\|_{L_{\psi}^{2}(\Omega)}^{2}\nonumber\\
    &+C\sum_{i=2}^{s-1}\|\partial_{y}\partial_{x}^{s-i+1}\mathbb{U}\partial_{x}^{i}\mathbb{U}\|_{L_{\psi}^{2}(\Omega)}^{2}
     +C\sum_{i=1}^{s-2}\|\partial_{y}^{-1}[\partial_{x}^{i+1}\mathbb{U}]\partial_{x}^{s-i}\partial_{y}^{2}\mathbb{U}\|_{L_{\psi}^{2}(\Omega)}^{2}\nonumber\\
\leq&~C_{*}\|\mathbb{U}\|_{\hat{H}_{\psi}^{s}(\Omega)}^{4}
\leq C_{*}.\nonumber
\end{align}
Combining \eqref{J6s1} with \eqref{J6s2}, from \eqref{J6s3} we have the following estimate of $J_{6}$ when $r=s$,
\begin{align}\label{J6s}
J_{6}
\leq  C_{*}\|\mathbf{U}_{s}\|_{L_{\varphi}^{2}(\Omega)}^{2}
     +\|\mathcal{P}\mathbf{U}_{s}\|_{L_{\varphi}^{2}(\Omega)}^{2}
     +C_{*}\left(1+\frac{1}{(\bar{\varsigma}+t)^{\frac{1}{2}+\varepsilon_{0}}}\right).
\end{align}

If $r=s+1$, the terms involving partial derivativesof order less than $s$ with respect to $x$ in $J_{6}$ can be estimated in a similar way as in the case $r=s$ above, while other terms must be estimated in different manner since they cannot be controlled by $\|\mathbb{U}\|_{\hat{H}_{\psi}^{s}(\Omega)}$. For this purpose, we divide $\mathcal{T}_{1}$ into five parts as follows,
\begin{align*}
\mathcal{T}_{1}
=&-\frac{s(s+1)}{2}\partial_{x}^{2}\mathbb{U}\mathbf{U}_{s}
  +(s+1)\mathcal{P}\left(\rho\partial_{x}^{2}\mathbb{U}-\partial_{y}\partial_{x}^{2}\mathbb{U}\right)\partial_{x}^{s}\mathbb{U}
  +\frac{s(s+1)}{2}\mathcal{P}\partial_{y}^{-1}[\partial_{x}^{s}\mathbb{U}](\partial_{x}^{2}\partial_{y}^{2}\mathbb{U}-\rho\partial_{x}^{2}\partial_{y}\mathbb{U})\\
 &+(s+1)\mathcal{P}\partial_{y}^{-1}[\partial_{x}^{2}\mathbb{U}](\partial_{x}^{s}\partial_{y}^{2}\mathbb{U}-\rho\partial_{x}^{s}\partial_{y}\mathbb{U})
  +\mathcal{P}(\rho\bar{\mathcal{T}}_{1R}-\partial_{y}\bar{\mathcal{T}}_{1R}),
\end{align*}
with
\begin{align*}
\bar{\mathcal{T}}_{1R}(t,x,y)=
 \sum_{i=3}^{s-1}\binom{s+1}{i}\partial_{x}^{s-i+2}\mathbb{U}\partial_{x}^{i}\mathbb{U}
-\sum_{i=2}^{s-2}\binom{s+1}{i}\partial_{y}^{-1}[\partial_{x}^{i+1}\mathbb{U}]\partial_{x}^{s-i+1}\partial_{y}\mathbb{U}.
\end{align*}
Accordingly,
\begin{align}\label{J6}
J_{6}
=&-\frac{s(s+1)}{2}\left(\partial_{x}^{2}\mathbb{U}\mathbf{U}_{s},\hat{\varphi}^{2}\mathbf{U}_{s+1}\right)
  +(s+1)\left(\mathcal{P}\left(\rho\partial_{x}^{2}\mathbb{U}-\partial_{y}\partial_{x}^{2}\mathbb{U}\right)\partial_{x}^{s}\mathbb{U},\hat{\varphi}^{2}\mathbf{U}_{s+1}\right)\\
 &+\frac{s(s+1)}{2}\left(\mathcal{P}\partial_{y}^{-1}[\partial_{x}^{s}\mathbb{U}](\partial_{x}^{2}\partial_{y}^{2}\mathbb{U}-\rho\partial_{x}^{2}\partial_{y}\mathbb{U}),\hat{\varphi}^{2}\mathbf{U}_{s+1}\right)\nonumber\\
 &+(s+1)\left(\mathcal{P}\partial_{y}^{-1}[\partial_{x}^{2}\mathbb{U}](\rho\partial_{x}^{s}\partial_{y}\mathbb{U}-\partial_{x}^{s}\partial_{y}^{2}\mathbb{U}),\hat{\varphi}^{2}\mathbf{U}_{s+1}\right)\nonumber\\
 &+\left(\mathcal{P}(\rho\bar{\mathcal{T}}_{1R}-\partial_{y}\bar{\mathcal{T}}_{1R}),\hat{\varphi}^{2}\mathbf{U}_{s+1}\right).\nonumber
\end{align}
It is clearly that
\begin{align}\label{J61}
-\frac{s(s+1)}{2}\left(\partial_{x}^{2}\mathbb{U}\mathbf{U}_{s},\hat{\varphi}^{2}\mathbf{U}_{s+1}\right)
\leq& C_{*}\|\mathbf{U}_{s}\|_{L_{\varphi}^{2}(\Omega)}\|\mathbf{U}_{s+1}\|_{L_{\hat{\varphi}}^{2}(\Omega)}\\
\leq& C_{*}\|\mathbf{U}_{s}\|_{L_{\varphi}^{2}(\Omega)}^{2}
     +C_{*}\|\mathbf{U}_{s+1}\|_{L_{\hat{\varphi}}^{2}(\Omega)}^{2}.\nonumber
\end{align}
Using Corollary \ref{ncC2}, one can deduce
\begin{align}\label{J62}
&(s+1)\left(\mathcal{P}\left(\rho\partial_{x}^{2}\mathbb{U}-\partial_{y}\partial_{x}^{2}\mathbb{U}\right)\partial_{x}^{s}\mathbb{U},\hat{\varphi}^{2}\mathbf{U}_{s+1}\right)\\
\leq&~C(\|\partial_{x}^{2}\mathbb{U}\|_{L^{\infty}(\Omega)}+\|\partial_{y}\partial_{x}^{2}\mathbb{U}\|_{L^{\infty}(\Omega)})
      \|\mathcal{P}^{2}\partial_{x}^{s}\mathbb{U}\|_{L_{\varphi}^{2}(\Omega)}
      \|\mathbf{U}_{s+1}\|_{L_{\hat{\varphi}}^{2}(\Omega)}\nonumber\\
\leq&~C_{*}\|\mathcal{P}\mathbf{U}_{s}\|_{L_{\varphi}^{2}(\Omega)}
      \|\mathbf{U}_{s+1}\|_{L_{\hat{\varphi}}^{2}(\Omega)}\nonumber\\
\leq&~C_{*}\|\mathbf{U}_{s+1}\|_{L_{\hat{\varphi}}^{2}(\Omega)}^{2}
     +\|\mathcal{P}\mathbf{U}_{s}\|_{L_{\varphi}^{2}(\Omega)}^{2}.\nonumber
\end{align}
Using Lemma \ref{rhoP} and recall the definition of $u_{0}$, one has
$\left(\|\mathcal{P}\partial_{x}^{2}\partial_{y}^{2}\mathbb{U}\|_{L_{x}^{\infty}L_{y}^{2}(\Omega)}
 +\|\mathcal{P}^{2}\partial_{x}^{2}\partial_{y}\mathbb{U}\|_{L_{x}^{\infty}L_{y}^{2}(\Omega)}\right)\leq C$.
Also by Corollary \ref{ncC2} one can obtain
\begin{align}\label{J63}
&\frac{s(s+1)}{2}\left(\mathcal{P}\partial_{y}^{-1}[\partial_{x}^{s}\mathbb{U}](\partial_{x}^{2}\partial_{y}^{2}\mathbb{U}-\rho\partial_{x}^{2}\partial_{y}\mathbb{U}),\hat{\varphi}^{2}\mathbf{U}_{s+1}\right)\\
\leq&~C_{*}\|\partial_{x}^{s}\mathbb{U}\|_{L_{\varphi}^{2}(\Omega)}
       \left(\|\mathcal{P}\partial_{x}^{2}\partial_{y}^{2}\mathbb{U}\|_{L_{x}^{\infty}L_{y}^{2}(\Omega)}+\|\mathcal{P}^{2}\partial_{x}^{2}\partial_{y}\mathbb{U}\|_{L_{x}^{\infty}L_{y}^{2}(\Omega)}\right)
       \|\mathbf{U}_{s+1}\|_{L_{\hat{\varphi}}^{2}(\Omega)}\nonumber\\
\leq&~C_{*}\|\mathcal{P}\mathbf{U}_{s}\|_{L_{\varphi}^{2}(\Omega)}
       \|\mathbf{U}_{s+1}\|_{L_{\hat{\varphi}}^{2}(\Omega)}\nonumber\\
\leq&~C_{*}\|\mathbf{U}_{s+1}\|_{L_{\hat{\varphi}}^{2}(\Omega)}^{2}
     +\|\mathcal{P}\mathbf{U}_{s}\|_{L_{\varphi}^{2}(\Omega)}^{2}.\nonumber
\end{align}
From Corollary \ref{wdotuC} and Lemma \ref{sdotu} one can deduce that for any $t\in[0,t_{3}]$ and $x\in\mathbb{T}$,
\begin{align*}
|\partial_{y}^{-1}[\partial_{x}^{2}\mathbb{U}]|
\leq \int_{0}^{y}\int_{\sigma}^{+\infty}|\partial_{y}\partial_{x}^{2}\mathbb{U}|d\sigma dy
\leq C\int_{0}^{y}\int_{\sigma}^{+\infty}(1+\y t)e^{-\frac{3\delta\y}{4}}d\sigma dy
\leq C.
\end{align*}
As a consequence,  by utilizing Corollary \ref{ncC2} one has
\begin{align}\label{J64'}
&(s+1)\left(\mathcal{P}\partial_{y}^{-1}[\partial_{x}^{2}\mathbb{U}](\rho\partial_{x}^{s}\partial_{y}\mathbb{U}-\partial_{x}^{s}\partial_{y}^{2}\mathbb{U}),\hat{\varphi}^{2}\mathbf{U}_{s+1}\right)\\
\leq&~C\left(\|\mathcal{P}\partial_{x}^{s}\partial_{y}\mathbb{U}\|_{L_{\varphi}^{2}(\Omega)}
             +\|\partial_{x}^{s}\partial_{y}^{2}\mathbb{U}\|_{L_{\varphi}^{2}(\Omega)}\right)
       \|\mathcal{P}\mathbf{U}_{s+1}\|_{L_{\hat{\varphi}}^{2}(\Omega)}\nonumber\\
\leq&~C_{\lambda}\left(\|\mathcal{P}\mathbf{U}_{s}\|_{L_{\varphi}^{2}(\Omega)}
             +\|\partial_{y}\mathbf{U}_{s}\|_{L_{\varphi}^{2}(\Omega)}
             +\|\mathbf{U}_{s}\|_{L_{\varphi}^{2}(\Omega)}\right)
       \|\mathcal{P}\mathbf{U}_{s+1}\|_{L_{\hat{\varphi}}^{2}(\Omega)}\nonumber\\
\leq&~C_{\lambda}\|\mathcal{P}\mathbf{U}_{s}\|_{L_{\varphi}^{2}(\Omega)}^{2}
     +C_{\lambda}\|\partial_{y}\mathbf{U}_{s}\|_{L_{\varphi}^{2}(\Omega)}^{2}
     +C\|\mathcal{P}\mathbf{U}_{s+1}\|_{L_{\hat{\varphi}}^{2}(\Omega)}^{2}
     +C_{\lambda}\|\mathbf{U}_{s}\|_{L_{\varphi}^{2}(\Omega)}^{2}\nonumber.
\end{align}

Finally, one can use a similar estimate as used in the case $r=s$ above to get
\begin{align}\label{J65}
\left(\mathcal{P}(\rho\bar{\mathcal{T}}_{1R}-\partial_{y}\bar{\mathcal{T}}_{1R}),\hat{\varphi}^{2}\mathbf{U}_{s+1}\right)
\leq  C_{*}\|\mathbf{U}_{s+1}\|_{L_{\hat{\varphi}}^{2}(\Omega)}^{2}
     +\|\mathcal{P}\mathbf{U}_{s+1}\|_{L_{\hat{\varphi}}^{2}(\Omega)}^{2}
     +C_{*}\left(1+\frac{1}{(\bar{\varsigma}+t)^{\frac{1}{2}+\varepsilon_{0}}}\right).
\end{align}

Substituting \eqref{J61},\eqref{J62},\eqref{J63},\eqref{J64'} and \eqref{J65} into \eqref{J6}, one can conclude that when $r=s+1$,
\begin{align}\label{J6'}
J_{6}
\leq& C_{\lambda}\|\mathcal{P}\mathbf{U}_{s}\|_{L_{\varphi}^{2}(\Omega)}^{2}
     +C\|\mathcal{P}\mathbf{U}_{s+1}\|_{L_{\hat{\varphi}}^{2}(\Omega)}^{2}
     +C_{*}\|\mathbf{U}_{s}\|_{L_{\varphi}^{2}(\Omega)}^{2}
     +C_{*}\|\mathbf{U}_{s+1}\|_{L_{\hat{\varphi}}^{2}(\Omega)}^{2}\\
    &+C_{\lambda}\|\partial_{y}\mathbf{U}_{s}\|_{L_{\varphi}^{2}(\Omega)}^{2}
     +C_{*}\left(1+\frac{1}{(\bar{\varsigma}+t)^{\frac{1}{2}+\varepsilon_{0}}}\right).\nonumber
\end{align}

\noindent{\bf \underline{Estimate of~$J_{7}$.}}
Recall
\begin{align}\label{J7}
J_{7}
= 2(\eta\mathcal{P}\partial_{x}\rho\partial_{x}^{r+1}\mathbb{U},\varphi_{r}^{2}\mathbf{U}_{r})
 +(\partial_{y}\eta\mathcal{P}\partial_{x}^{r+2}\mathbb{U},\varphi_{r}^{2}\mathbf{U}_{r})
 -(\mathcal{P}\partial_{y}^{-1}[\partial_{x}^{r+1}\tilde{\mathfrak{U}}](\varsigma'-\rho\varsigma),\varphi_{r}^{2}\mathbf{U}_{r}).
\end{align}
By using integration by parts, one has
\begin{align}\label{J71}
2(\eta\mathcal{P}\partial_{x}\rho\partial_{x}^{r+1}\mathbb{U},\varphi_{r}^{2}\mathbf{U}_{r})
=&-2(\eta\partial_{x}(\mathcal{P}\partial_{x}\rho)\partial_{x}^{r}\mathbb{U},\varphi_{r}^{2}\mathbf{U}_{r})
 -2(\eta\mathcal{P}\partial_{x}\rho\partial_{x}^{r}\mathbb{U},\varphi_{r}^{2}\partial_{x}\mathbf{U}_{r})\\
\leq& \|\eta\partial_{x}(\mathcal{P}\partial_{x}\rho)\partial_{x}^{r}\mathbb{U}\|_{L_{\varphi_{r}}^{2}(\Omega)}^{2}
     +C\|\sqrt{\eta}\mathcal{P}\partial_{x}\rho\partial_{x}^{r}\mathbb{U}\|_{L_{\varphi_{r}}^{2}(\Omega)}^{2}\nonumber\\
    &+\frac{1}{16}\|\sqrt{\eta}\partial_{x}\mathbf{U}_{r}\|_{L_{\varphi_{r}}^{2}(\Omega)}^{2}
     +C\|\mathbf{U}_{r}\|_{L_{\varphi_{r}}^{2}(\Omega)}^{2}.\nonumber
\end{align}
In virtue of Lemma \ref{rhoP}, by noticing that $\bar{\varsigma}=\iota\bar{\eta}^{2}$,
from Lemma \ref{nc} one has
\begin{align}\label{J711}
\|\eta\partial_{x}(\mathcal{P}\partial_{x}\rho)\partial_{x}^{r}\mathbb{U}\|_{L_{\varphi_{r}}^{2}(\Omega)}^{2}
\leq& C_{*}\bar{\eta}^{2}\|\partial_{x}^{r}\mathbb{U}\|_{L_{\hat{\Psi}_{r}}^{2}(\Omega)}^{2}\\
\leq& C_{*}\bar{\eta}^{2}\|\mathcal{P}\mathbf{U}_{r}\|_{L_{\varphi_{r}}^{2}(\Omega)}^{2}\nonumber\\
\leq& C_{*}\iota^{-1}\|\mathbf{U}_{r}\|_{L_{\varphi_{r}}^{2}(\Omega)}^{2},\nonumber
\end{align}
and
\begin{align}\label{J712}
\|\sqrt{\eta}\mathcal{P}\partial_{x}\rho\partial_{x}^{r}\mathbb{U}\|_{L_{\varphi_{r}}^{2}(\Omega)}^{2}
\leq& C_{*}\bar{\eta}\|\partial_{x}^{r}\mathbb{U}\|_{L_{\hat{\Psi}_{r}}^{2}(\Omega)}^{2}\\
\leq& C_{*}\iota^{-\frac{1}{2}}\|\mathbf{U}_{r}\|_{L_{\varphi_{r}}^{2}(\Omega)}\|\mathcal{P}\mathbf{U}_{r}\|_{L_{\varphi_{r}}^{2}(\Omega)}\nonumber\\
\leq& C_{*}\iota^{-1}\|\mathbf{U}_{r}\|_{L_{\varphi_{r}}^{2}(\Omega)}^{2}
     +\frac{1}{2}\|\mathcal{P}\mathbf{U}_{r}\|_{L_{\varphi_{r}}^{2}(\Omega)}^{2}.\nonumber
\end{align}
Substitute \eqref{J711} and \eqref{J712} into \eqref{J71}, one has
\begin{align}\label{J71'}
2(\eta\mathcal{P}\partial_{x}\rho\partial_{x}^{r+1}\mathbb{U},\varphi_{r}^{2}\mathbf{U}_{r})
\leq \frac{1}{16}\|\sqrt{\eta}\partial_{x}\mathbf{U}_{r}\|_{L_{\varphi_{r}}^{2}(\Omega)}^{2}
    +C_{*}\iota^{-1}\|\mathbf{U}_{r}\|_{L_{\varphi_{r}}^{2}(\Omega)}^{2}
    +\frac{1}{2}\|\mathcal{P}\mathbf{U}_{r}\|_{L_{\varphi_{r}}^{2}(\Omega)}^{2}.
\end{align}
Also by integration by parts, noticing that \eqref{theta} and Corollary \ref{ncC3}, one has
\begin{align*}
(\partial_{y}\eta\mathcal{P}\partial_{x}^{r+2}\mathbb{U},\varphi_{r}^{2}\mathbf{U}_{r})
=&-(\partial_{y}\eta\partial_{x}^{r+1}\mathbb{U},\varphi_{r}^{2}\partial_{x}\mathbf{U}_{r})\\
\leq&~C\theta t\|\vartheta_{c}\partial_{x}^{r+1}\mathbb{U}\|_{L_{\varphi_{r}}^{2}(\Omega)}
       \|\vartheta_{c}\partial_{x}\mathbf{U}_{r}\|_{L_{\varphi_{r}}^{2}(\Omega)}\\
\leq&~C_{\lambda}\bar{\eta}\theta\left(\sqrt{t}\|\vartheta_{c}\partial_{x}\mathbf{U}_{r}\|_{L_{\varphi_{r}}^{2}(\Omega)}
                                       +\|\vartheta_{c}\mathbf{U}_{r}\|_{L_{\varphi_{r}}^{2}(\Omega)}\right)
      \|\vartheta_{c}\partial_{x}\mathbf{U}_{r}\|_{L_{\varphi_{r}}^{2}(\Omega)}\\
\leq&~\left(\frac{1}{16}+C_{\lambda}\theta\sqrt{t}\right)\|\sqrt{\eta}\partial_{x}\mathbf{U}_{r}\|_{L_{\varphi_{r}}^{2}(\Omega)}^{2}
     +C_{*}C_{\theta}\|\mathbf{U}_{r}\|_{L_{\varphi_{r}}^{2}(\Omega)}^{2}.
\end{align*}
Noticing ${\bf supp}~\varsigma \subset(\check{y}_{*},\hat{y}_{*})$ and $\mathcal{P}\rho=-2\partial_{y}\mathcal{P}$ on $[\check{y}_{*},\hat{y}_{*}]$, by using the integration by parts one has
\begin{align}\label{J72}
-(\mathcal{P}\partial_{y}^{-1}[\partial_{x}^{r+1}\tilde{\mathfrak{U}}](\varsigma'-\rho\varsigma),\varphi_{r}^{2}\mathbf{U}_{r})
=&-(\mathcal{P}\varsigma'\partial_{y}^{-1}[\partial_{x}^{r+1}\tilde{\mathfrak{U}}],\varphi_{r}^{2}\mathbf{U}_{r})
  -2(\partial_{y}\mathcal{P}\varsigma\partial_{y}^{-1}[\partial_{x}^{r+1}\tilde{\mathfrak{U}}],\varphi_{r}^{2}\mathbf{U}_{r})\\
=& (\mathcal{P}\varsigma'\partial_{y}^{-1}[\partial_{x}^{r+1}\tilde{\mathfrak{U}}],\varphi_{r}^{2}\mathbf{U}_{r})
  +2(\mathcal{P}\varsigma\partial_{x}^{r+1}\tilde{\mathfrak{U}},\varphi_{r}^{2}\mathbf{U}_{r})\nonumber\\
 &+2(\mathcal{P}\varsigma\partial_{y}^{-1}[\partial_{x}^{r+1}\tilde{\mathfrak{U}}],\varphi_{r}^{2}\partial_{y}\mathbf{U}_{r})
  +4(\mathcal{P}\varsigma\partial_{y}^{-1}[\partial_{x}^{r+1}\tilde{\mathfrak{U}}],\varphi_{r}\partial_{y}\varphi_{r}\mathbf{U}_{r})\nonumber\\
=& (\mathcal{P}\varsigma'\partial_{y}^{-1}[\partial_{x}^{r+1}\tilde{\mathbb{U}}],\varphi_{r}^{2}\mathbf{U}_{r})
  +2(\mathcal{P}\varsigma\partial_{x}^{r+1}\tilde{\mathbb{U}},\varphi_{r}^{2}\mathbf{U}_{r})\nonumber\\
 &+2(\mathcal{P}\varsigma\partial_{y}^{-1}[\partial_{x}^{r+1}\tilde{\mathbb{U}}],\varphi_{r}^{2}\partial_{y}\mathbf{U}_{r})
  +4(\mathcal{P}\varsigma\partial_{y}^{-1}[\partial_{x}^{r+1}\tilde{\mathbb{U}}],\varphi_{r}\partial_{y}\varphi_{r}\mathbf{U}_{r})\nonumber\\
 &+(\mathcal{P}\varsigma'\partial_{y}^{-1}[\partial_{x}^{r+1}u_{0}],\varphi_{r}^{2}\mathbf{U}_{r})
  +2(\mathcal{P}\varsigma\partial_{x}^{r+1}u_{0},\varphi_{r}^{2}\mathbf{U}_{r})\nonumber\\
 &+2(\mathcal{P}\varsigma\partial_{y}^{-1}[\partial_{x}^{r+1}u_{0}],\varphi_{r}^{2}\partial_{y}\mathbf{U}_{r})
  +4(\mathcal{P}\varsigma\partial_{y}^{-1}[\partial_{x}^{r+1}u_{0}],\varphi_{r}\partial_{y}\varphi_{r}\mathbf{U}_{r})\nonumber.
\end{align}
Since $\sqrt{\bar{\varsigma}}\mathcal{P}\leq C$, $\sqrt{\tilde{\bar{\varsigma}}}\tilde{\mathcal{P}}\leq C$ and $\bar{\varsigma}=\iota\bar{\eta}^{2}$, it follows from Corollary \ref{ncC3} that
\begin{align}\label{J721}
(\mathcal{P}\varsigma'\partial_{y}^{-1}[\partial_{x}^{r+1}\tilde{\mathbb{U}}],\varphi_{r}^{2}\mathbf{U}_{r})
\leq&~C_{*}\bar{\varsigma}\|\vartheta_{c}\partial_{x}^{r+1}\tilde{\mathbb{U}}\|_{L_{\check{\psi}}^{2}(\Omega)}
                 \|\varphi_{r}\|_{L^{2}[\check{y},\hat{y}]}
                 \|\mathcal{P}\mathbf{U}_{r}\|_{L_{\varphi_{r}}^{2}(\Omega)}\\
\leq&~C_{*}\bar{\varsigma}\|\mathcal{P}\omega_{o}^{-1}\|_{L^{\infty}(\Omega)}
      \left(\|\vartheta_{c}\partial_{x}\tilde{\mathbf{U}}_{r}\|_{L_{\tilde{\varphi}_{r}}^{2}(\Omega)}
            +\|\tilde{\mathcal{P}}\tilde{\mathbf{U}}_{r}\|_{L_{\tilde{\varphi}_{r}}^{2}(\Omega)}\right)
      \|\mathcal{P}\mathbf{U}_{r}\|_{L_{\varphi_{r}}^{2}(\Omega)}\nonumber\\
\leq&~C_{*}\bar{\varsigma}\tilde{\bar{\varsigma}}^{-\frac{1}{2}}
      \left(\|\vartheta_{c}\partial_{x}\tilde{\mathbf{U}}_{r}\|_{L_{\tilde{\varphi}_{r}}^{2}(\Omega)}
            +\|\tilde{\mathcal{P}}\tilde{\mathbf{U}}_{r}\|_{L_{\tilde{\varphi}_{r}}^{2}(\Omega)}\right)
      \|\mathcal{P}\mathbf{U}_{r}\|_{L_{\varphi_{r}}^{2}(\Omega)}\nonumber\\
\leq&~\tau\|\sqrt{\tilde{\eta}}\partial_{x}\tilde{\mathbf{U}}_{r}\|_{L_{\tilde{\varphi}_{r}}^{2}(\Omega)}^{2}
     +C_{*}C_{\tau}\bar{\varsigma}^{2}\tilde{\bar{\varsigma}}^{-1}\tilde{\bar{\eta}}^{-1}\|\mathcal{P}\mathbf{U}_{r}\|_{L_{\varphi_{r}}^{2}(\Omega)}^{2}\nonumber\\
    &+C_{*}\bar{\varsigma}\tilde{\bar{\varsigma}}^{-\frac{1}{2}}\|\tilde{\mathcal{P}}\tilde{\mathbf{U}}_{r}\|_{L_{\tilde{\varphi}_{r}}^{2}(\Omega)}
                                              \|\mathcal{P}\mathbf{U}_{r}\|_{L_{\varphi_{r}}^{2}(\Omega)}\nonumber\\
\leq&~\tau\|\sqrt{\tilde{\eta}}\partial_{x}\tilde{\mathbf{U}}_{r}\|_{L_{\tilde{\varphi}_{r}}^{2}(\Omega)}^{2}
     +C_{*}C_{\tau}\|\mathbf{U}_{r}\|_{L_{\varphi_{r}}^{2}(\Omega)}\|\mathcal{P}\mathbf{U}_{r}\|_{L_{\varphi_{r}}^{2}(\Omega)}\nonumber\\
    &+C_{*}\|\tilde{\mathbf{U}}_{r}\|_{L_{\tilde{\varphi}_{r}}^{2}(\Omega)}\|\mathcal{P}\mathbf{U}_{r}\|_{L_{\varphi_{r}}^{2}(\Omega)}\nonumber\\
\leq&~\tau\|\sqrt{\tilde{\eta}}\partial_{x}\tilde{\mathbf{U}}_{r}\|_{L_{\tilde{\varphi}_{r}}^{2}(\Omega)}^{2}
     +C_{*}\|\mathbf{U}_{r}\|_{L_{\varphi_{r}}^{2}(\Omega)}^{2}
     +C_{*}\|\tilde{\mathbf{U}}_{r}\|_{L_{\tilde{\varphi}_{r}}^{2}(\Omega)}^{2}\nonumber\\
    &+\frac{1}{2}\|\mathcal{P}\mathbf{U}_{r}\|_{L_{\varphi_{r}}^{2}(\Omega)}^{2},\nonumber
\end{align}
Similarly, we have
\begin{align}\label{J722}
2\left|(\mathcal{P}\varsigma\partial_{x}^{r+1}\tilde{\mathbb{U}},\varphi_{r}^{2}\mathbf{U}_{r})\right|
\leq& C\bar{\varsigma}\|\tilde{\mathcal{P}}\vartheta_{c}\partial_{x}^{r+1}\tilde{\mathbb{U}}\|_{L_{\tilde{\varphi}_{r}}^{2}(\Omega)}
             \|\mathbf{U}_{r}\|_{L_{\varphi_{r}}^{2}(\Omega)}\\
\leq& C_{*}\tilde{\bar{\varsigma}}\|\vartheta_{c}\partial_{x}^{r+1}\tilde{\mathbb{U}}\|_{L_{\tilde{\hat{\Psi}}_{r}}^{2}(\Omega)}
                 \|\mathbf{U}_{r}\|_{L_{\varphi_{r}}^{2}(\Omega)}\nonumber\\
\leq& C_{*}\left(\|\sqrt{\tilde{\eta}}\partial_{x}\tilde{\mathbf{U}}_{r}\|_{L_{\tilde{\varphi}_{r}}^{2}(\Omega)}
            +\|\tilde{\mathbf{U}}_{r}\|_{L_{\tilde{\varphi}_{r}}^{2}(\Omega)}\right)
      \|\mathbf{U}_{r}\|_{L_{\varphi_{r}}^{2}(\Omega)}\nonumber\\
\leq& \tau\|\sqrt{\tilde{\eta}}\partial_{x}\tilde{\mathbf{U}}_{r}\|_{L_{\tilde{\varphi}_{r}}^{2}(\Omega)}^{2}
     +C_{*}C_{\tau}\|\mathbf{U}_{r}\|_{L_{\varphi_{r}}^{2}(\Omega)}^{2}
     +C_{*}\|\tilde{\mathbf{U}}_{r}\|_{L_{\tilde{\varphi}_{r}}^{2}(\Omega)}^{2}.\nonumber
\end{align}
Moreover, note that $\|\mathcal{P}\omega_{o}^{-1}\|_{L^{\infty}(\Omega)}\leq C\bar{\varsigma}^{-\frac{1}{2}+\frac{\epsilon_{0}}{2}}$ for any $t\in[0,t_{3}]$, by using Corollary \ref{ncC3} we can obtain
\begin{align}\label{J723}
&2(\mathcal{P}\varsigma\partial_{y}^{-1}[\partial_{x}^{r+1}\tilde{\mathbb{U}}],\varphi_{r}^{2}\partial_{y}\mathbf{U}_{r})\\
\leq&~C_{*}\bar{\varsigma}\|\vartheta_{c}\partial_{x}^{r+1}\tilde{\mathbb{U}}\|_{L_{\check{\psi}}^{2}(\Omega)}
       \|\mathcal{P}\varphi_{r}\|_{L_{x}^{\infty}(\mathbb{T})L_{y}^{2}([\check{y},\hat{y}])}\|\partial_{y}\mathbf{U}_{r}\|_{L_{\varphi_{r}}^{2}(\Omega)}\nonumber\\
\leq&~C_{*}\bar{\varsigma}^{\frac{3}{4}-\frac{\varepsilon_{0}}{2}}\|\vartheta_{c}\partial_{x}^{r+1}\tilde{\mathbb{U}}\|_{L_{\check{\psi}}^{2}(\Omega)}
                               \|\partial_{y}\mathbf{U}_{r}\|_{L_{\varphi_{r}}^{2}(\Omega)}\nonumber\\
\leq&~C_{*}\bar{\varsigma}^{\frac{3}{4}-\frac{\varepsilon_{0}}{2}}\|\mathcal{P}\omega_{o}^{-1}\|_{L^{\infty}(\Omega)}
           \left(\|\vartheta_{c}\partial_{x}\tilde{\mathbf{U}}_{r}\|_{L_{\tilde{\varphi}_{r}}^{2}(\Omega)}
            +\|\tilde{\mathcal{P}}\tilde{\mathbf{U}}_{r}\|_{L_{\tilde{\varphi}_{r}}^{2}(\Omega)}\right)
      \|\partial_{y}\mathbf{U}_{r}\|_{L_{\varphi_{r}}^{2}(\Omega)}\nonumber\\
\leq&~C_{*}\bar{\varsigma}^{\frac{3}{4}}\tilde{\bar{\varsigma}}^{-\frac{1}{2}}
           \left(\|\vartheta_{c}\partial_{x}\tilde{\mathbf{U}}_{r}\|_{L_{\tilde{\varphi}_{r}}^{2}(\Omega)}
            +\|\tilde{\mathcal{P}}\tilde{\mathbf{U}}_{r}\|_{L_{\tilde{\varphi}_{r}}^{2}(\Omega)}\right)
      \|\partial_{y}\mathbf{U}_{r}\|_{L_{\varphi_{r}}^{2}(\Omega)}\nonumber\\
\leq&~C_{*}\bar{\varsigma}^{\frac{3}{2}}\tilde{\bar{\varsigma}}^{-1}
           \|\vartheta_{c}\partial_{x}\tilde{\mathbf{U}}_{r}\|_{L_{\tilde{\varphi}_{r}}^{2}(\Omega)}^{2}
     +C_{*}\bar{\varsigma}^{\frac{3}{2}}\tilde{\bar{\varsigma}}^{-1}\|\tilde{\mathcal{P}}\tilde{\mathbf{U}}_{r}\|_{L_{\tilde{\varphi}_{r}}^{2}(\Omega)}^{2}
     +\frac{1}{16}\|\partial_{y}\mathbf{U}_{r}\|_{L_{\varphi_{r}}^{2}(\Omega)}^{2}\nonumber\\
\leq&~C_{*}\sqrt{\iota}\|\sqrt{\tilde{\eta}}\partial_{x}\tilde{\mathbf{U}}_{r}\|_{L_{\tilde{\varphi}_{r}}^{2}(\Omega)}^{2}
     +\frac{1}{16}\|\partial_{y}\mathbf{U}_{r}\|_{L_{\varphi_{r}}^{2}(\Omega)}^{2}
     +C_{*}\|\tilde{\mathbf{U}}_{r}\|_{L_{\tilde{\varphi}_{r}}^{2}(\Omega)}\|\tilde{\mathcal{P}}\tilde{\mathbf{U}}_{r}\|_{L_{\tilde{\varphi}_{r}}^{2}(\Omega)}\nonumber\\
\leq&~C_{*}\sqrt{\iota}\|\sqrt{\tilde{\eta}}\partial_{x}\tilde{\mathbf{U}}_{r}\|_{L_{\tilde{\varphi}_{r}}^{2}(\Omega)}^{2}
     +\frac{1}{16}\|\partial_{y}\mathbf{U}_{r}\|_{L_{\varphi_{r}}^{2}(\Omega)}^{2}
     +C_{*}C_{\tau}\|\tilde{\mathbf{U}}_{r}\|_{L_{\tilde{\varphi}_{r}}^{2}(\Omega)}^{2}
     +\frac{1}{2}\tau\|\tilde{\mathcal{P}}\tilde{\mathbf{U}}_{r}\|_{L_{\tilde{\varphi}_{r}}^{2}(\Omega)}^{2},\nonumber
\end{align}
With the help of \eqref{varphiy}, we can employ a similar estimate as in \eqref{J723} to obtain
\begin{align}\label{J724}
&4(\mathcal{P}\varsigma\partial_{y}^{-1}[\partial_{x}^{r+1}\tilde{\mathbb{U}}],\varphi_{r}\partial_{y}\varphi_{r}\mathbf{U}_{r})\\
\leq&~C_{*}\bar{\varsigma}^{\frac{3}{4}}\tilde{\bar{\varsigma}}^{-\frac{1}{2}}
      \left(\|\vartheta_{c}\partial_{x}\tilde{\mathbf{U}}_{r}\|_{L_{\tilde{\varphi}_{r}}^{2}(\Omega)}
            +\|\tilde{\mathcal{P}}\tilde{\mathbf{U}}_{r}\|_{L_{\tilde{\varphi}_{r}}^{2}(\Omega)}\right)
      \|\varphi_{r}^{-1}\partial_{y}\varphi_{r}^{2}\mathbf{U}_{r}\|_{L^{2}(\Omega)}\nonumber\\
\leq&~C_{*}\bar{\varsigma}^{\frac{3}{2}}\tilde{\bar{\varsigma}}^{-1}
           \|\vartheta_{c}\partial_{x}\tilde{\mathbf{U}}_{r}\|_{L_{\tilde{\varphi}_{r}}^{2}(\Omega)}^{2}
     +C_{*}\bar{\varsigma}^{\frac{3}{2}}\tilde{\bar{\varsigma}}^{-1}\|\tilde{\mathcal{P}}\tilde{\mathbf{U}}_{r}\|_{L_{\tilde{\varphi}_{r}}^{2}(\Omega)}^{2}
     +\frac{1}{16}\|\sqrt{\omega_{\lambda}}\mathbf{U}_{r}\|_{L_{\varphi_{r}}^{2}(\Omega)}^{2}
     +C\|\mathbf{U}_{r}\|_{L_{\varphi_{r}}^{2}(\Omega)}^{2}\nonumber\\
\leq&~C_{*}\sqrt{\iota}\|\sqrt{\tilde{\eta}}\partial_{x}\tilde{\mathbf{U}}_{r}\|_{L_{\tilde{\varphi}_{r}}^{2}(\Omega)}^{2}
     +\frac{1}{16}\|\sqrt{\omega_{\lambda}}\mathbf{U}_{r}\|_{L_{\varphi_{r}}^{2}(\Omega)}^{2}
     +C_{*}\|\mathbf{U}_{r}\|_{L_{\varphi_{r}}^{2}(\Omega)}^{2}\nonumber\\
    &+C_{*}C_{\tau}\|\tilde{\mathbf{U}}_{r}\|_{L_{\tilde{\varphi}_{r}}^{2}(\Omega)}^{2}
     +\frac{1}{2}\tau\|\tilde{\mathcal{P}}\tilde{\mathbf{U}}_{r}\|_{L_{\tilde{\varphi}_{r}}^{2}(\Omega)}^{2}.\nonumber
\end{align}
It remains to estimate the last four terms on the right hand side of \eqref{J72}. Noticing ${\bf supp}\varsigma \subset(\check{y}_{*},\hat{y}_{*})$ and \eqref{varphiy}, we have
\begin{align}\label{J725}
&\quad(\mathcal{P}\varsigma'\partial_{y}^{-1}[\partial_{x}^{r+1}u_{0}],\varphi_{r}^{2}\mathbf{U}_{r})
+2(\mathcal{P}\varsigma\partial_{x}^{r+1}u_{0},\varphi_{r}^{2}\mathbf{U}_{r})\\
&+2(\mathcal{P}\varsigma\partial_{y}^{-1}[\partial_{x}^{r+1}u_{0}],\varphi_{r}^{2}\partial_{y}\mathbf{U}_{r})
+2(\mathcal{P}\varsigma\partial_{y}^{-1}[\partial_{x}^{r+1}u_{0}],\partial_{y}\varphi_{r}^{2}\mathbf{U}_{r})\nonumber\\
\leq& C_{*}\|\mathbf{U}_{r}\|_{L_{\varphi_{r}}^{2}(\Omega)}^{2}
    +\frac{1}{16}\|\partial_{y}\mathbf{U}_{r}\|_{L_{\varphi_{r}}^{2}(\Omega)}^{2}
    +\frac{1}{16}\|\sqrt{\omega_{\lambda}}\mathbf{U}_{r}\|_{L_{\varphi_{r}}^{2}(\Omega)}^{2}
    +C_{*}.\nonumber
\end{align}

Substituting \eqref{J721},\eqref{J722},\eqref{J723},\eqref{J724} and \eqref{J725} into \eqref{J72}, we have
\begin{align}\label{J72'}
 -(\mathcal{P}\partial_{y}^{-1}[\partial_{x}^{r+1}\tilde{\mathfrak{U}}](\varsigma'-\rho\varsigma),\varphi_{r}^{2}\mathbf{U}_{r})
\leq& (2\tau+C_{*}\sqrt{\iota})\|\sqrt{\tilde{\eta}}\partial_{x}\tilde{\mathbf{U}}_{r}\|_{L_{\tilde{\varphi}_{r}}^{2}(\Omega)}^{2}
     +\frac{1}{8}\|\sqrt{\omega_{\lambda}}\mathbf{U}_{r}\|_{L_{\varphi_{r}}^{2}(\Omega)}^{2}\\
    &+\frac{1}{8}\|\partial_{y}\mathbf{U}_{r}\|_{L_{\varphi_{r}}^{2}(\Omega)}^{2}
     +C_{*}C_{\tau}\|\mathbf{U}_{r}\|_{L_{\varphi_{r}}^{2}(\Omega)}^{2}
     +C_{*}C_{\tau}\|\tilde{\mathbf{U}}_{r}\|_{L_{\tilde{\varphi}_{r}}^{2}(\Omega)}^{2}\nonumber\\
    &+\frac{1}{2}\|\mathcal{P}\mathbf{U}_{r}\|_{L_{\varphi_{r}}^{2}(\Omega)}^{2}
     +\tau\|\tilde{\mathcal{P}}\tilde{\mathbf{U}}_{r}\|_{L_{\tilde{\varphi}_{r}}^{2}(\Omega)}^{2}+C_{*}.\nonumber
\end{align}

Now combining \eqref{J71'} and \eqref{J72'} with \eqref{J7}, we can conclude that
\begin{align*}
J_{7}
\leq& \left(\frac{1}{8}+C_{\lambda}\theta\sqrt{t}\right)\|\sqrt{\eta}\partial_{x}\mathbf{U}_{r}\|_{L_{\varphi_{r}}^{2}(\Omega)}^{2}
     +\frac{1}{8}\|\sqrt{\omega_{\lambda}}\mathbf{U}_{r}\|_{L_{\varphi_{r}}^{2}(\Omega)}^{2}
     +\frac{1}{8}\|\partial_{y}\mathbf{U}_{r}\|_{L_{\varphi_{r}}^{2}(\Omega)}^{2}
     +\|\mathcal{P}\mathbf{U}_{r}\|_{L_{\varphi_{r}}^{2}(\Omega)}^{2}\\
    &+\tau\|\tilde{\mathcal{P}}\tilde{\mathbf{U}}_{r}\|_{L_{\tilde{\varphi}_{r}}^{2}(\Omega)}^{2}
     +C_{*}C_{\tau}C_{\theta}\iota^{-1}\|\mathbf{U}_{r}\|_{L_{\varphi_{r}}^{2}(\Omega)}^{2}
     +(2\tau+C_{*}\sqrt{\iota})\|\sqrt{\tilde{\eta}}\partial_{x}\tilde{\mathbf{U}}_{r}\|_{L_{\varphi_{r}}^{2}(\Omega)}^{2}
     +C_{*}C_{\tau}\|\tilde{\mathbf{U}}_{r}\|_{L_{\tilde{\varphi}_{r}}^{2}(\Omega)}^{2}+C_{*}\nonumber.
\end{align*}

\noindent{\bf \underline{Estimate of~$J_{8}$.}}
Taking notice of the fact that $\partial_{x}\partial_{y}^{2}\mathfrak{U}-\rho\partial_{x}\partial_{y}\mathfrak{U}=\partial_{x}\rho(\dot{\mathfrak{U}}+\varsigma)$,
and $|\partial_{x}\rho|\leq C\mathcal{P}$, we have
\begin{align*}
J_{8}
=&~r\left(\mathcal{P}\left(\partial_{x}\partial_{y}^{2}\mathfrak{U}-\rho\partial_{x}\partial_{y}\mathfrak{U}\right)\partial_{y}^{-1}[\partial_{x}^{r}\mathfrak{U}],\varphi_{r}^{2}\mathbf{U}_{r}\right)\\
=&~r\left((\dot{\mathfrak{U}}+\varsigma)\mathcal{P}\partial_{x}\rho\partial_{y}^{-1}[\partial_{x}^{r}\mathbb{U}],\varphi_{r}^{2}\mathbf{U}_{r}\right)
  +r\left((\dot{\mathfrak{U}}+\varsigma)\mathcal{P}\partial_{x}\rho\partial_{y}^{-1}[\partial_{x}^{r}u_{0}],\varphi_{r}^{2}\mathbf{U}_{r}\right)\\
\leq&~C_{*}\|\varphi_{r}(\dot{\mathfrak{U}}+\varsigma)\mathcal{P}^{2}\|_{L_{x}^{\infty}(\mathbb{T})L_{y}^{2}(\mathbb{R_{+}})}
       \|\partial_{x}^{r}\mathbb{U}\|_{L_{\varphi_{r}}^{2}(\Omega)}\|\mathbf{U}_{r}\|_{L_{\varphi_{r}}^{2}(\Omega)}\\
    &+C_{*}\|\varphi_{r}\sqrt{\y}(\dot{\mathfrak{U}}+\varsigma)\mathcal{P}^{2}\|_{L_{x}^{\infty}(\mathbb{T})L_{y}^{2}(\mathbb{R_{+}})}
           \|\mathbf{U}_{r}\|_{L_{\varphi_{r}}^{2}(\Omega)}\\
\leq&~C_{*}\|\partial_{x}^{r}\mathbb{U}\|_{L_{\varphi_{r}}^{2}(\Omega)}\|\mathbf{U}_{r}\|_{L_{\varphi_{r}}^{2}(\Omega)}
     +C_{*}\left(\|\sqrt{\y}\mathbb{U}\|_{\hat{H}_{\psi}^{s}(\Omega)}+1\right)\|\mathbf{U}_{r}\|_{L_{\varphi_{r}}^{2}(\Omega)}\\
\leq&~C_{*}\|\mathcal{P}\mathbf{U}_{r}\|_{L_{\varphi_{r}}^{2}(\Omega)}\|\mathbf{U}_{r}\|_{L_{\varphi_{r}}^{2}(\Omega)}
     +C_{*}\left(\|\sqrt{\y}\mathbb{U}\|_{\hat{H}_{\psi}^{s}(\Omega)}+1\right)\|\mathbf{U}_{r}\|_{L_{\varphi_{r}}^{2}(\Omega)}\\
\leq&~\|\mathcal{P}\mathbf{U}_{r}\|_{L_{\varphi_{r}}^{2}(\Omega)}^{2}
     +C_{*}C_{\tau}\|\mathbf{U}_{r}\|_{L_{\varphi_{r}}^{2}(\Omega)}^{2}
     +\tau\|\sqrt{\y}\mathbb{U}\|_{\hat{H}_{\psi}^{s}(\Omega)}^{2}
     +C_{*}.
\end{align*}

\noindent{\bf \underline{Estimate of~$J_{9}$.}}
We divide $J_{9}$ into four parts,
\begin{align}\label{J9}
J_{9}
=&\left((\partial_{y}^{2}\hat{F}-\rho\partial_{y}\hat{F})\mathcal{P}\frac{\partial_{x}^{r}\mathbb{U}}{\dot{\mathfrak{U}}+\varsigma},\varphi_{r}^{2}\mathbf{U}_{r}\right)
  -2\left(\eta\mathcal{P}\partial_{x}\partial_{y}\mathbb{U}\partial_{x}\rho\frac{\partial_{x}^{r}\mathbb{U}}{\dot{\mathfrak{U}}+\varsigma},\varphi_{r}^{2}\mathbf{U}_{r}\right)\\
 &+\left(\mathcal{P}\left(\varsigma'\partial_{x}\mathfrak{U}+\varsigma'''+\partial_{y}^{-1}[\partial_{x}\mathfrak{U}]\varsigma''-\varsigma\partial_{x}\partial_{y}\mathfrak{U}+\eta\partial_{x}^{2}\partial_{y}^{2}u_{0}\right)\frac{\partial_{x}^{r}\mathbb{U}}{\dot{\mathfrak{U}}+\varsigma},\varphi_{r}^{2}\mathbf{U}_{r}\right)\nonumber\\
 &-\left(\mathcal{P}\rho\left(\partial_{y}^{-1}[\partial_{x}\mathfrak{U}]\varsigma'+\varsigma''+\eta\partial_{x}^{2}\partial_{y}u_{0}\right)\frac{\partial_{x}^{r}\mathbb{U}}{\dot{\mathfrak{U}}+\varsigma},\varphi_{r}^{2}\mathbf{U}_{r}\right).\nonumber
\end{align}

Since $\|\mathcal{P}(\partial_{y}^{2}\hat{F}-\rho\partial_{y}\hat{F})\|_{L^{\infty}(\Omega_{t_{3}})}\leq C_{*}$,
and for any $\gamma\in\Gamma_{2}, |D^{\gamma}\hat{F}|\leq C_{*}e^{-\frac{4}{3}\delta\y}~\text{in}~\Omega_{T_{*}},$
with Corollary \ref{ncC2} one can deduce that
\begin{align}\label{J91}
\left((\partial_{y}^{2}\hat{F}-\rho\partial_{y}\hat{F})\mathcal{P}\frac{\partial_{x}^{r}\mathbb{U}}{\dot{\mathfrak{U}}+\varsigma},\varphi_{r}^{2}\mathbf{U}_{r}\right)
\leq& C_{*}\|\partial_{x}^{r}\mathbb{U}\|_{L_{\hat{\Psi}_{r}}^{2}(\Omega)}\|\mathbf{U}_{r}\|_{L_{\varphi_{r}}^{2}(\Omega)}\\
\leq& C_{*}\|\mathcal{P}\mathbf{U}_{r}\|_{L_{\varphi_{r}}^{2}(\Omega)}\|\mathbf{U}_{r}\|_{L_{\varphi_{r}}^{2}(\Omega)}\nonumber\\
\leq& C_{*}\|\mathbf{U}_{r}\|_{L_{\varphi_{r}}^{2}(\Omega)}^{2}
     +\frac{1}{3}\|\mathcal{P}\mathbf{U}_{r}\|_{L_{\varphi_{r}}^{2}(\Omega)}^{2}.\nonumber
\end{align}
Due to Corollary \ref{wdotuC'} and Lemma \ref{rhoP}, moreover, $\partial_{x}u_{0}=0$ for $y\in[0,\hat{y}]$, one has $|\vartheta_{c}^{2}\partial_{x}\partial_{y}\mathbb{U}|\leq C\dot{\mathfrak{U}}$ in $\Omega_{t_{3}}$, which leads to
\begin{align}\label{J92}
-2\left(\eta\mathcal{P}\partial_{x}\partial_{y}\mathbb{U}\partial_{x}\rho\frac{\partial_{x}^{r}\mathbb{U}}{\dot{\mathfrak{U}}+\varsigma},\varphi_{r}^{2}\mathbf{U}_{r}\right)
\leq& C_{*}\bar{\eta}\|\partial_{x}^{r}\mathbb{U}\|_{L_{\hat{\Psi}_{r}}^{2}(\Omega)}\|\mathbf{U}_{r}\|_{L_{\varphi_{r}}^{2}(\Omega)}\\
\leq& C_{*}\bar{\eta}\|\mathcal{P}\mathbf{U}_{r}\|_{L_{\varphi_{r}}^{2}(\Omega)}\|\mathbf{U}_{r}\|_{L_{\varphi_{r}}^{2}(\Omega)}\nonumber\\
\leq& C_{*}\|\mathbf{U}_{r}\|_{L_{\varphi_{r}}^{2}(\Omega)}^{2}
     +\frac{1}{3}\|\mathcal{P}\mathbf{U}_{r}\|_{L_{\varphi_{r}}^{2}(\Omega)}^{2},\nonumber
\end{align}
Notice that
$$\varsigma'\partial_{x}\mathfrak{U}+\varsigma'''+\partial_{y}^{-1}[\partial_{x}\mathfrak{U}]\varsigma''-\varsigma\partial_{x}\partial_{y}\mathfrak{U}+\eta\partial_{x}^{2}\partial_{y}^{2}u_{0}
\leq C_{*}\iota^{-\frac{1}{2}}\sqrt{\bar{\varsigma}},$$
and
$$\rho\left(\partial_{y}^{-1}[\partial_{x}\mathfrak{U}]\varsigma'+\varsigma''+\eta\partial_{x}^{2}\partial_{y}u_{0}\right)
\leq C_{*}\iota^{-\frac{1}{2}}\sqrt{\bar{\varsigma}},$$
also by Corollary \ref{ncC2} one gets
\begin{align}\label{J93}
&\left(\mathcal{P}\left(\varsigma'\partial_{x}\mathfrak{U}+\varsigma'''+\partial_{y}^{-1}[\partial_{x}\mathfrak{U}]\varsigma''-\varsigma\partial_{x}\partial_{y}\mathfrak{U}+\eta\partial_{x}^{2}\partial_{y}^{2}u_{0}\right)\frac{\partial_{x}^{r}\mathbb{U}}{\dot{\mathfrak{U}}+\varsigma},\varphi_{r}^{2}\mathbf{U}_{r}\right)\\
&-\left(\mathcal{P}\rho\left(\partial_{y}^{-1}[\partial_{x}\mathfrak{U}]\varsigma'+\varsigma''+\eta\partial_{x}^{2}\partial_{y}u_{0}\right)\frac{\partial_{x}^{r}\mathbb{U}}{\dot{\mathfrak{U}}+\varsigma},\varphi_{r}^{2}\mathbf{U}_{r}\right)\nonumber\\
\leq& C_{*}\iota^{-\frac{1}{2}}\sqrt{\bar{\varsigma}}\|\partial_{x}^{r}\mathbb{U}\|_{L_{\hat{\Psi}_{r}}^{2}(\Omega)}
                               \|\mathcal{P}\mathbf{U}_{r}\|_{L_{\varphi_{r}}^{2}(\Omega)}\nonumber\\
\leq& C_{*}\iota^{-\frac{1}{2}}\sqrt{\bar{\varsigma}}\|\mathcal{P}\mathbf{U}_{r}\|_{L_{\varphi_{r}}^{2}(\Omega)}^{2}\nonumber\\
\leq& C_{*}\iota^{-1}\|\mathbf{U}_{r}\|_{L_{\varphi_{r}}^{2}(\Omega)}^{2}
     +\frac{1}{3}\|\mathcal{P}\mathbf{U}_{r}\|_{L_{\varphi_{r}}^{2}(\Omega)}^{2},\nonumber
\end{align}

A combination of \eqref{J91}, \eqref{J92} and \eqref{J93} with \eqref{J9} gives
\begin{align*}
J_{9}
\leq& C_{*}\iota^{-1}\|\mathbf{U}_{r}\|_{L_{\varphi_{r}}^{2}(\Omega)}^{2}
     +\|\mathcal{P}\mathbf{U}_{r}\|_{L_{\varphi_{r}}^{2}(\Omega)}^{2}.
\end{align*}

\noindent{\bf \underline{Estimate of~$J_{10}$.}}
Recall that
\begin{align}\label{J10}
J_{10}
= (\mathcal{P}\varsigma\partial_{x}^{r+1}u,\varphi_{r}^{2}\mathbf{U}_{r})
 +(\mathcal{P}(\rho\partial_{x}^{r}\mathbb{F}-\partial_{y}\partial_{x}^{r}\mathbb{F}),\varphi_{r}^{2}\mathbf{U}_{r})
\end{align}
In view of $u=\mathbb{U}-\tilde{\mathbb{U}}$, one can divide the first term on the right hand side of \eqref{J10} into two parts,
\begin{align}\label{J101}
(\mathcal{P}\varsigma\partial_{x}^{r+1}u,\varphi_{r}^{2}\mathbf{U}_{r})
=& (\mathcal{P}\varsigma\partial_{x}^{r+1}\mathbb{U},\varphi_{r}^{2}\mathbf{U}_{r})
  -(\mathcal{P}\varsigma\partial_{x}^{r+1}\tilde{\mathbb{U}},\varphi_{r}^{2}\mathbf{U}_{r})
\end{align}
As a result of Corollary \ref{ncC3}, it follows
\begin{align}\label{J1011}
(\mathcal{P}\varsigma\partial_{x}^{r+1}\mathbb{U},\varphi_{r}^{2}\mathbf{U}_{r})
\leq& C\bar{\varsigma}\|\vartheta_{c}\mathcal{P}\partial_{x}^{r+1}\mathbb{U}\|_{L_{\varphi_{r}}^{2}(\Omega)}
       \|\mathbf{U}_{r}\|_{L_{\varphi_{r}}^{2}(\Omega)}\\
\leq& C_{*}\bar{\varsigma}\|\vartheta_{c}\partial_{x}^{r+1}\mathbb{U}\|_{L_{\hat{\Psi}_{r}}^{2}(\Omega)}
       \|\mathbf{U}_{r}\|_{L_{\varphi_{r}}^{2}(\Omega)}\nonumber\\
\leq& C_{*}\sqrt{\bar{\varsigma}}
      \left(\|\vartheta_{c}\partial_{x}\mathbf{U}_{r}\|_{L_{\varphi_{r}}^{2}(\Omega)}
            +\|\mathcal{P}\mathbf{U}_{r}\|_{L_{\varphi_{r}}^{2}(\Omega)}\right)
      \|\mathbf{U}_{r}\|_{L_{\varphi_{r}}^{2}(\Omega)}\nonumber\\
\leq& \frac{1}{16}\|\sqrt{\eta}\partial_{x}\mathbf{U}_{r}\|_{L_{\varphi_{r}}^{2}(\Omega)}^{2}
     +C_{*}\|\mathbf{U}_{r}\|_{L_{\varphi_{r}}^{2}(\Omega)}^{2}.\nonumber
\end{align}
Combining \eqref{J1011} and \eqref{J722} with \eqref{J101}, one has
\begin{align}\label{J101'}
(\mathcal{P}\varsigma\partial_{x}^{r+1}\mathbb{U},\varphi_{r}^{2}\mathbf{U}_{r})
\leq& \frac{1}{16}\|\sqrt{\eta}\partial_{x}\mathbf{U}_{r}\|_{L_{\varphi_{r}}^{2}(\Omega)}^{2}
     +C_{*}C_{\tau}\|\mathbf{U}_{r}\|_{L_{\varphi_{r}}^{2}(\Omega)}^{2}\\
    &+\tau\|\sqrt{\tilde{\eta}}\partial_{x}\tilde{\mathbf{U}}_{r}\|_{L_{\tilde{\varphi}_{r}}^{2}(\Omega)}^{2}
     +C_{*}\|\tilde{\mathbf{U}}_{r}\|_{L_{\tilde{\varphi}_{r}}^{2}(\Omega)}^{2}.\nonumber
\end{align}

Next let's turn to the estimate of the second term on the right hand side of \eqref{J10}. By a direct calculation, we decompose $\rho\partial_{x}^{r}\mathbb{F}-\partial_{y}\partial_{x}^{r}\mathbb{F}$ into four parts according to the order of the derivatives of $u$ with respect to $x$ as follows,
\begin{align*}
\rho\partial_{x}^{r}\mathbb{F}-\partial_{y}\partial_{x}^{r}\mathbb{F}
= \partial_{y}\partial_{x}^{r}\left(u\partial_{x}u-(1-h)\partial_{y}^{-1}[\partial_{x}u]\partial_{y}u\right)
  -\rho\partial_{x}^{r}\left(u\partial_{x}u-(1-h)\partial_{y}^{-1}[\partial_{x}u]\partial_{y}u\right)
=: \sum_{i=1}^{4}\mathbb{J}_{i},
\end{align*}
where
\begin{align*}
\mathbb{J}_{1}
=& u\partial_{x}^{r+1}\partial_{y}u
  -\rho u\partial_{x}^{r+1}u
  +(1-h)\partial_{y}^{-1}[\partial_{x}^{r+1}u](\rho\partial_{y}u-\partial_{y}^{2}u)\\
 &+h'\partial_{y}^{-1}[\partial_{x}^{r+1}u]\partial_{y}u
  +h\partial_{x}^{r+1}u\partial_{y}u,
\end{align*}
\begin{align*}
\mathbb{J}_{2}
=& \partial_{x}^{r}u(\partial_{x}\partial_{y}u-\rho\partial_{x}u)
  +r\partial_{x}u(\partial_{x}^{r}\partial_{y}u-\rho\partial_{x}^{r}u)
  -r(1-h)\partial_{y}^{-1}[\partial_{x}^{r}u](\partial_{x}\partial_{y}^{2}u-\rho\partial_{x}\partial_{y}u)\\
 &-(1-h)\partial_{y}^{-1}[\partial_{x}u](\partial_{x}^{r}\partial_{y}^{2}u-\rho\partial_{x}^{r}\partial_{y}u)
  +h'(r\partial_{y}^{-1}[\partial_{x}^{r}u]\partial_{x}\partial_{y}u+\partial_{y}^{-1}[\partial_{x}u]\partial_{x}^{r}\partial_{y}u)\\
 &+h(r\partial_{x}^{r}u\partial_{x}\partial_{y}u+\partial_{x}u\partial_{x}^{r}\partial_{y}u),
\end{align*}
\begin{align*}
\mathbb{J}_{3}
=& r\partial_{x}^{r-1}u(\partial_{x}^{2}\partial_{y}u-\rho\partial_{x}^{2}u)
  +\frac{r(r-1)}{2}\partial_{x}^{2}u(\partial_{x}^{r-1}\partial_{y}u-\rho\partial_{x}^{r-1}u)\\
 &-\frac{r(r-1)}{2}(1-h)\partial_{y}^{-1}[\partial_{x}^{r-1}u](\partial_{x}^{2}\partial_{y}^{2}u-\rho\partial_{x}^{2}\partial_{y}u)
  -r(1-h)\partial_{y}^{-1}[\partial_{x}^{2}u](\partial_{x}^{r-1}\partial_{y}^{2}u-\rho\partial_{x}^{r-1}\partial_{y}u)\\
 &+rh'\partial_{y}^{-1}[\partial_{x}^{2}u]\partial_{x}^{r-1}\partial_{y}u
  +\frac{r(r-1)}{2}h'\partial_{y}^{-1}[\partial_{x}^{r-1}u]\partial_{x}^{2}\partial_{y}u
  +rh\partial_{x}^{2}u\partial_{x}^{r-1}\partial_{y}u
  +\frac{r(r-1)}{2}h\partial_{x}^{r-1}u\partial_{x}^{2}\partial_{y}u,
\end{align*}
and
\begin{align*}
\mathbb{J}_{4}
=& \sum_{i=3}^{r-2}\binom{r}{i}[\partial_{x}^{i}u\partial_{x}^{r-i+1}\partial_{y}u-\rho\partial_{x}^{i}u\partial_{x}^{r-i+1}u]\\
 &-(1-h)\sum_{i=2}^{r-3}\binom{r}{i}\left[\partial_{y}^{-1}[\partial_{x}^{i+1}u]\partial_{x}^{r-i}\partial_{y}^{2}u-\rho\partial_{y}^{-1}[\partial_{x}^{i+1}u]\partial_{x}^{n-i}\partial_{y}u\right]\\
 &+h\sum_{i=2}^{r-3}\binom{r}{i}\partial_{x}^{i+1}u\partial_{x}^{r-i}\partial_{y}u
  +h'\sum_{i=2}^{r-3}\binom{r}{i}\partial_{y}^{-1}[\partial_{x}^{i+1}u]\partial_{x}^{r-i}\partial_{y}u
\end{align*}
Accordingly, one has
\begin{align*}
(\mathcal{P}(\rho\partial_{x}^{r}\mathbb{F}-\partial_{y}\partial_{x}^{r}\mathbb{F}),\varphi_{r}^{2}\mathbf{U}_{r})
= \sum_{i=1}^{4}(\mathcal{P}\mathbb{J}_{i},\varphi_{r}^{2}\mathbf{U}_{r})
=: \sum_{i=1}^{4}\mathfrak{J}_{i}.
\end{align*}
Let's start with the estimate of $\mathfrak{J}_{1}$. We divide $\mathfrak{J}_{1}$ into three parts,
\begin{align*}
\mathfrak{J}_{1}
&= (\mathcal{P}(u\partial_{x}^{r+1}\partial_{y}u-\rho u\partial_{x}^{r+1}u),\varphi_{r}^{2}\mathbf{U}_{r})
  +(\mathcal{P}(1-h)\partial_{y}^{-1}[\partial_{x}^{r+1}u](\rho\partial_{y}u-\partial_{y}^{2}u),\varphi_{r}^{2}\mathbf{U}_{r})\\
&\quad+(\mathcal{P}(h'\partial_{y}^{-1}[\partial_{x}^{r+1}u]\partial_{y}u+h\partial_{x}^{r+1}u\partial_{y}u),\varphi_{r}^{2}\mathbf{U}_{r})\\
&=: \mathfrak{J}_{11}+\mathfrak{J}_{12}+\mathfrak{J}_{13}.
\end{align*}
By a direct calculation, one has
\begin{align}\label{JJ11}
\mathfrak{J}_{11}
=& (\mathcal{P}u(\partial_{x}^{r+1}\partial_{y}\mathbb{U}-\rho \partial_{x}^{r+1}\mathbb{U}),\varphi_{r}^{2}\mathbf{U}_{r})
  -(\mathcal{P}u(\partial_{x}^{r+1}\partial_{y}\tilde{\mathbb{U}}-\rho \partial_{x}^{r+1}\tilde{\mathbb{U}}),\varphi_{r}^{2}\mathbf{U}_{r})\\
=& (\mathcal{P}u\partial_{x}(\mathcal{P}^{-1}\mathbf{U}_{r}),\varphi_{r}^{2}\mathbf{U}_{r})
  +(\mathcal{P}u\partial_{x}\rho\partial_{x}^{r}\mathbb{U},\varphi_{r}^{2}\mathbf{U}_{r})
  -(\mathcal{P}u(\partial_{x}^{r+1}\partial_{y}\tilde{\mathbb{U}}-\tilde{\rho}\partial_{x}^{r+1}\tilde{\mathbb{U}}),\varphi_{r}^{2}\mathbf{U}_{r})\nonumber\\
 &+(\mathcal{P}u(\rho-\tilde{\rho})\partial_{x}^{r+1}\tilde{\mathbb{U}},\varphi_{r}^{2}\mathbf{U}_{r})\nonumber\\
=& (\mathcal{P}u\partial_{x}(\mathcal{P}^{-1}\mathbf{U}_{r}),\varphi_{r}^{2}\mathbf{U}_{r})
  -(\mathcal{P}u\partial_{x}(\tilde{\mathcal{P}}^{-1}\tilde{\mathbf{U}}_{r}),\varphi_{r}^{2}\mathbf{U}_{r})
  +(\mathcal{P}u\partial_{x}\rho\partial_{x}^{r}\mathbb{U},\varphi_{r}^{2}\mathbf{U}_{r})\nonumber\\
 &-(\mathcal{P}u\partial_{x}\tilde{\rho}\partial_{x}^{r}\tilde{\mathbb{U}},\varphi_{r}^{2}\mathbf{U}_{r})
  +(\mathcal{P}u(\rho-\tilde{\rho})\partial_{x}^{r+1}\tilde{\mathbb{U}},\varphi_{r}^{2}\mathbf{U}_{r}).\nonumber
\end{align}
Utilizing \eqref{us-3}, and recalling that $\bar{\eta}=\epsilon_{0}^{2}$ when $\check{c}\neq0$, one has
\begin{align}\label{JJ111}
 &(\mathcal{P}u\partial_{x}(\mathcal{P}^{-1}\mathbf{U}_{r}),\varphi_{r}^{2}\mathbf{U}_{r})
-(\mathcal{P}u\partial_{x}(\tilde{\mathcal{P}}^{-1}\tilde{\mathbf{U}}_{r}),\varphi_{r}^{2}\mathbf{U}_{r})\\
=& (u\partial_{x}\mathbf{U}_{r},\varphi_{r}^{2}\mathbf{U}_{r})
  -(\mathcal{P}\tilde{\mathcal{P}}^{-1}u\partial_{x}\tilde{\mathbf{U}}_{r},\varphi_{r}^{2}\mathbf{U}_{r})
  +(\mathcal{P}u\partial_{x}\mathcal{P}^{-1}\mathbf{U}_{r},\varphi_{r}^{2}\mathbf{U}_{r})
  -(\mathcal{P}u\partial_{x}\tilde{\mathcal{P}}^{-1}\tilde{\mathbf{U}}_{r},\varphi_{r}^{2}\mathbf{U}_{r})\nonumber\\
\leq& C_{*}(\check{c}+\sqrt{\bar{\eta}}t)\|\vartheta_{c}\partial_{x}\mathbf{U}_{r}\|_{L_{\varphi_{r}}^{2}(\Omega)}
                                             \|\mathbf{U}_{r}\|_{L_{\varphi_{r}}^{2}(\Omega)}
     +C_{*}(\check{c}+\sqrt{\bar{\eta}}t)\|\vartheta_{c}\partial_{x}\tilde{\mathbf{U}}_{r}\|_{L_{\tilde{\varphi}_{r}}^{2}(\Omega)}\|\mathbf{U}_{r}\|_{L_{\varphi_{r}}^{2}(\Omega)}\nonumber\\
    &+C_{*}(\check{c}+\sqrt{\bar{\eta}}t)\|\mathcal{P}\|_{L^{\infty}(\Omega)}\|\mathbf{U}_{r}\|_{L_{\varphi_{r}}^{2}(\Omega)}^{2}
     +C_{*}(\check{c}+\sqrt{\bar{\eta}}t)\|\mathcal{P}\|_{L^{\infty}(\Omega)}\|\mathbf{U}_{r}\|_{L_{\varphi_{r}}^{2}(\Omega)}
                                             \|\tilde{\mathbf{U}}_{r}\|_{L_{\tilde{\varphi}_{r}}^{2}(\Omega)}\nonumber\\
\leq& \frac{1}{16}\|\sqrt{\eta}\partial_{x}\mathbf{U}_{r}\|_{L_{\varphi_{r}}^{2}(\Omega)}^{2}
     +\tau\|\sqrt{\tilde{\eta}}\partial_{x}\tilde{\mathbf{U}}_{r}\|_{L_{\varphi_{r}}^{2}(\Omega)}^{2}
     +C_{*}C_{\tau}\iota^{-1}\|\mathbf{U}_{r}\|_{L_{\varphi_{r}}^{2}(\Omega)}^{2}
     +\|\tilde{\mathbf{U}}_{r}\|_{L_{\tilde{\varphi}_{r}}^{2}(\Omega)}^{2}.\nonumber
\end{align}
By Corollary \ref{ncC2}, one can obtain
\begin{align}\label{JJ112}
&(\mathcal{P}u\partial_{x}\rho\partial_{x}^{r}\mathbb{U},\varphi_{r}^{2}\mathbf{U}_{r})
-(\mathcal{P}u\partial_{x}\tilde{\rho}\partial_{x}^{r}\tilde{\mathbb{U}},\varphi_{r}^{2}\mathbf{U}_{r})\\
\leq& C_{*}(\check{c}+\sqrt{\bar{\eta}}t)\|\partial_{x}^{r}\mathbb{U}\|_{L_{\hat{\Psi}_{r}}^{2}(\Omega)}
                                             \|\mathbf{U}_{r}\|_{L_{\varphi_{r}}^{2}(\Omega)}
     +C_{*}(\check{c}+\sqrt{\bar{\eta}}t)\|\partial_{x}^{r}\tilde{\mathbb{U}}\|_{L_{\tilde{\hat{\Psi}}_{r}}^{2}(\Omega)}
                                             \|\mathbf{U}_{r}\|_{L_{\varphi_{r}}^{2}(\Omega)}\nonumber\\
\leq& C_{*}(\check{c}+\sqrt{\bar{\eta}}t)\|\mathcal{P}\mathbf{U}_{r}\|_{L_{\varphi_{r}}^{2}(\Omega)}
                                             \|\mathbf{U}_{r}\|_{L_{\varphi_{r}}^{2}(\Omega)}
     +C_{*}(\check{c}+\sqrt{\bar{\eta}}t)\|\tilde{\mathcal{P}}\tilde{\mathbf{U}}_{r}\|_{L_{\tilde{\varphi}_{r}}^{2}(\Omega)}
                                             \|\mathbf{U}_{r}\|_{L_{\varphi_{r}}^{2}(\Omega)}\nonumber\\
\leq& C_{*}\iota^{-1}\|\mathbf{U}_{r}\|_{L_{\varphi_{r}}^{2}(\Omega)}^{2}
     +\|\tilde{\mathbf{U}}_{r}\|_{L_{\tilde{\varphi}_{r}}^{2}(\Omega)}^{2}.\nonumber
\end{align}
Using Corollary \ref{ncC3} and the fact $\check{c}\leq C\iota^{-1}\bar{\varsigma}$, one can deduce
\begin{align}\label{JJ113}
&(\mathcal{P}u(\rho-\tilde{\rho})\partial_{x}^{r+1}\tilde{\mathbb{U}},\varphi_{r}^{2}\mathbf{U}_{r})\\
\leq&~C_{*}(\check{c}+\sqrt{\bar{\eta}}t)
           \|\vartheta_{c}\partial_{x}^{r+1}\tilde{\mathbb{U}}\|_{L_{\tilde{\hat{\Psi}}_{r}}^{2}(\Omega)}
           \|\mathbf{U}_{r}\|_{L_{\varphi_{r}}^{2}(\Omega)}\nonumber\\
\leq&~C_{*}(\check{c}+\sqrt{\bar{\eta}}t)\|\tilde{\mathcal{P}}\|_{L^{\infty}(\Omega)}
           \left(\|\vartheta_{c}\partial_{x}\tilde{\mathbf{U}}_{r}\|_{L_{\tilde{\varphi}_{r}}^{2}(\Omega)}
                 +\|\tilde{\mathcal{P}}\tilde{\mathbf{U}}_{r}\|_{L_{\tilde{\varphi}_{r}}^{2}(\Omega)}\right)
      \|\mathbf{U}_{r}\|_{L_{\varphi_{r}}^{2}(\Omega)}\nonumber\\
\leq&~C_{*}(\iota^{-1}\sqrt{\bar{\varsigma}}+\sqrt{\bar{\eta}t})
           \|\vartheta_{c}\partial_{x}\tilde{\mathbf{U}}_{r}\|_{L_{\tilde{\varphi}_{r}}^{2}(\Omega)}\|\mathbf{U}_{r}\|_{L_{\varphi_{r}}^{2}(\Omega)}\nonumber\\
    &+C_{*}(\check{c}+\sqrt{\bar{\eta}}t)\|\tilde{\mathcal{P}}\|_{L^{\infty}(\Omega)}
           \|\tilde{\mathcal{P}}\tilde{\mathbf{U}}_{r}\|_{L_{\tilde{\varphi}_{r}}^{2}(\Omega)}\|\mathbf{U}_{r}\|_{L_{\varphi_{r}}^{2}(\Omega)}\nonumber\\
\leq& \tau\|\sqrt{\tilde{\eta}}\partial_{x}\tilde{\mathbf{U}}_{r}\|_{L_{\tilde{\varphi}_{r}}^{2}(\Omega)}^{2}
     +C_{*}C_{\tau}(\iota^{-2}+t)\|\mathbf{U}_{r}\|_{L_{\varphi_{r}}^{2}(\Omega)}^{2}\nonumber\\
    &+C_{*}\iota^{-1}\|\tilde{\mathbf{U}}_{r}\|_{L_{\tilde{\varphi}_{r}}^{2}(\Omega)}\|\mathbf{U}_{r}\|_{L_{\varphi_{r}}^{2}(\Omega)}\nonumber\\
\leq& \tau\|\sqrt{\tilde{\eta}}\partial_{x}\tilde{\mathbf{U}}_{r}\|_{L_{\tilde{\varphi}_{r}}^{2}(\Omega)}^{2}
     +\|\tilde{\mathbf{U}}_{r}\|_{L_{\tilde{\varphi}_{r}}^{2}(\Omega)}^{2}
     +C_{*}C_{\tau}\iota^{-2}\|\mathbf{U}_{r}\|_{L_{\varphi_{r}}^{2}(\Omega)}^{2}.\nonumber
\end{align}
Substituting \eqref{JJ111}, \eqref{JJ112} and \eqref{JJ113} into \eqref{JJ11}, it follows
\begin{align}\label{JJ11'}
\mathfrak{J}_{11}
\leq \frac{1}{16}\|\sqrt{\eta}\partial_{x}\mathbf{U}_{r}\|_{L_{\varphi_{r}}^{2}(\Omega)}^{2}
    +2\tau\|\sqrt{\tilde{\eta}}\partial_{x}\tilde{\mathbf{U}}_{r}\|_{L_{\varphi_{r}}^{2}(\Omega)}^{2}
    +C_{*}C_{\tau}\iota^{-2}\|\mathbf{U}_{r}\|_{L_{\varphi_{r}}^{2}(\Omega)}^{2}
    +C\|\tilde{\mathbf{U}}_{r}\|_{L_{\tilde{\varphi}_{r}}^{2}(\Omega)}^{2}.
\end{align}
The term $\mathfrak{J}_{12}$ should be divided into two parts,
\begin{align}\label{JJ12}
\mathfrak{J}_{12}
= (\mathcal{P}(1-h)\partial_{y}^{-1}[\partial_{x}^{r+1}\mathbb{U}](\rho\partial_{y}u-\partial_{y}^{2}u),\varphi_{r}^{2}\mathbf{U}_{r})
 -(\mathcal{P}(1-h)\partial_{y}^{-1}[\partial_{x}^{r+1}\tilde{\mathbb{U}}](\rho\partial_{y}u-\partial_{y}^{2}u),\varphi_{r}^{2}\mathbf{U}_{r})
\end{align}
For the first part on the right hand side of \eqref{JJ12}, since
$\|\rho\|_{L^{\infty}(\Omega)}\leq C$ in $\mathbb{R}_{+}\setminus[\check{y},\hat{y}]$,
using \eqref{us-3} and Corollary \ref{ncC3}  we divide it into two part to get its estimate as follows,
\begin{align}\label{JJ121}
&(\mathcal{P}(1-h)\partial_{y}^{-1}[\partial_{x}^{r+1}\mathbb{U}](\rho\partial_{y}u-\partial_{y}^{2}u),\varphi_{r}^{2}\mathbf{U}_{r})\\
=&~\int_{\mathbb{T}}\int_{\check{y}}^{\hat{y}}\varphi_{r}^{2}\mathcal{P}(1-h)\partial_{y}^{-1}[\partial_{x}^{r+1}\mathbb{U}](\rho\partial_{y}u-\partial_{y}^{2}u)\mathbf{U}_{r}dydx\nonumber\\
 &+\int_{\mathbb{T}}\int_{\mathbb{R}_{+}\setminus[\check{y},\hat{y}]}\varphi_{r}^{2}\mathcal{P}(1-h)\partial_{y}^{-1}[\partial_{x}^{r+1}\mathbb{U}](\rho\partial_{y}u-\partial_{y}^{2}u)\mathbf{U}_{r}dydx\nonumber\\
\leq&~C_{*}\left(\|\partial_{y}u\|_{L^{\infty}(\mathbb{T}\times[\check{y},\hat{y}])}
                 +\|\partial_{y}^{2}u\|_{L^{\infty}(\mathbb{T}\times[\check{y},\hat{y}])}\right)
          \|\vartheta_{c}\partial_{x}^{r+1}\mathbb{U}\|_{L_{\hat{\Psi}_{c}}^{2}(\Omega)}
          \|\mathcal{P}\varphi_{r}\|_{L_{x}^{\infty}(\mathbb{T})L_{y}^{2}([\check{y},\hat{y}])}
          \|\mathcal{P}\mathbf{U}_{r}\|_{L_{\varphi_{r}}^{2}(\Omega)}\nonumber\\
    &+C_{*}\left(\|\psi\partial_{y}u\|_{L_{x}^{\infty}L_{y}^{2}(\Omega)}+\|\psi\partial_{y}^{2}u\|_{L_{x}^{\infty}L_{y}^{2}(\Omega)}\right)
          \|\vartheta_{c}\partial_{x}^{r+1}\mathbb{U}\|_{L_{\hat{\Psi}_{r}}^{2}(\Omega)}\|\mathbf{U}_{r}\|_{L_{\varphi_{r}}^{2}(\Omega)}\nonumber\\
\leq&~C_{*}(\check{c}+\sqrt{\bar{\eta}}t)
          \left(\int_{\check{y}_{*}}^{\hat{y}_{*}}\frac{1}{((y-y_{*})^2+\bar{\varsigma}+t)^{1+\varepsilon_{0}}}dy\right)^{\frac{1}{2}}
          \|\vartheta_{c}\partial_{x}^{r+1}\mathbb{U}\|_{L_{\hat{\Psi}_{c}}^{2}(\Omega)}
          \|\mathcal{P}\mathbf{U}_{r}\|_{L_{\varphi_{r}}^{2}(\Omega)}\nonumber\\
    &+C_{*}(\check{c}+\sqrt{\bar{\eta}}t)\|\vartheta_{c}\partial_{x}^{r+1}\mathbb{U}\|_{L_{\hat{\Psi}_{r}}^{2}(\Omega)}
          \|\mathbf{U}_{r}\|_{L_{\varphi_{r}}^{2}(\Omega)}\nonumber\\
\leq&~C_{*}(\check{c}\iota^{-\frac{1}{4}-\frac{\varepsilon_{0}}{2}}+\sqrt{\bar{\eta}}t^{\frac{3}{4}-\frac{\varepsilon_{0}}{2}})
           \|\mathcal{P}\omega_{o}^{-1}\|_{L^{\infty}(\Omega)}
           \left(\|\vartheta_{c}\partial_{x}\mathbf{U}_{r}\|_{L_{\varphi_{r}}^{2}(\Omega)}
                 +\|\mathcal{P}\mathbf{U}_{r}\|_{L_{\varphi_{r}}^{2}(\Omega)}\right)
          \|\mathcal{P}\mathbf{U}_{r}\|_{L_{\varphi_{r}}^{2}(\Omega)}\nonumber\\
    &+C_{*}(\check{c}+\sqrt{\bar{\eta}}t)\|\vartheta_{c}\partial_{x}^{r+1}\mathbb{U}\|_{L_{\hat{\Psi}_{r}}^{2}(\Omega)}
          \|\mathbf{U}_{r}\|_{L_{\varphi_{r}}^{2}(\Omega)}\nonumber\\
\leq&~C_{*}(\check{c}\iota^{-\frac{3}{4}}+\sqrt{\bar{\eta}}t^{\frac{1}{4}})
      \left(\|\vartheta_{c}\partial_{x}\mathbf{U}_{r}\|_{L_{\varphi_{r}}^{2}(\Omega)}
            +\|\mathcal{P}\mathbf{U}_{r}\|_{L_{\varphi_{r}}^{2}(\Omega)}\right)
      \|\mathcal{P}\mathbf{U}_{r}\|_{L_{\varphi_{r}}^{2}(\Omega)}\nonumber\\
    &+C_{*}(\check{c}\iota^{-\frac{1}{2}}+\sqrt{\bar{\eta}t})
      \left(\|\vartheta_{c}\partial_{x}\mathbf{U}_{r}\|_{L_{\varphi_{r}}^{2}(\Omega)}
            +\|\mathcal{P}\mathbf{U}_{r}\|_{L_{\varphi_{r}}^{2}(\Omega)}\right)
      \|\mathbf{U}_{r}\|_{L_{\varphi_{r}}^{2}(\Omega)}\nonumber\\
\leq&~\frac{1}{16}\|\sqrt{\eta}\partial_{x}\mathbf{U}_{r}\|_{L_{\varphi_{r}}^{2}(\Omega)}^{2}
     +C_{*}(\sqrt{t}+\sqrt{\bar{\eta}}t^{\frac{1}{4}})\|\mathcal{P}\mathbf{U}_{r}\|_{L_{\varphi_{r}}^{2}(\Omega)}^{2}
     +C_{*}\iota^{-\frac{5}{2}}\|\mathbf{U}_{r}\|_{L_{\varphi_{r}}^{2}(\Omega)}^{2}.\nonumber
\end{align}
Here $\hat{\Psi}_{c}=\y^{-\frac{1}{2}}\Psi_{c}$, with $\Psi_{c}$ given below Corollary \ref{ncC1}. Analogously, it follows from Corollary \ref{ncC3} that
\begin{align}\label{JJ122}
&-(\mathcal{P}\partial_{y}^{-1}[\partial_{x}^{r+1}\tilde{\mathbb{U}}](\rho\partial_{y}u-\partial_{y}^{2}u),\varphi_{r}^{2}\mathbf{U}_{r})\\
\leq& C_{*}(\check{c}\iota^{-\frac{3}{4}}+\sqrt{\bar{\eta}}t^{\frac{1}{4}})
      \left(\|\vartheta_{c}\partial_{x}\tilde{\mathbf{U}}_{r}\|_{L_{\tilde{\varphi}_{r}}^{2}(\Omega)}
            +\|\tilde{\mathcal{P}}\tilde{\mathbf{U}}_{r}\|_{L_{\tilde{\varphi}_{r}}^{2}(\Omega)}\right)
      \|\mathcal{P}\mathbf{U}_{r}\|_{L_{\varphi_{r}}^{2}(\Omega)}\nonumber\\
    &+C_{*}(\check{c}\iota^{-\frac{1}{2}}+\sqrt{\bar{\eta}t})
      \left(\|\vartheta_{c}\partial_{x}\tilde{\mathbf{U}}_{r}\|_{L_{\tilde{\varphi}_{r}}^{2}(\Omega)}
            +\|\tilde{\mathcal{P}}\tilde{\mathbf{U}}_{r}\|_{L_{\tilde{\varphi}_{r}}^{2}(\Omega)}\right)
      \|\mathbf{U}_{r}\|_{L_{\varphi_{r}}^{2}(\Omega)}\nonumber\\
\leq& \tau\|\sqrt{\tilde{\eta}}\partial_{x}\tilde{\mathbf{U}}_{r}\|_{L_{\tilde{\varphi}_{r}}^{2}(\Omega)}
     +C_{*}C_{\tau}\sqrt{t}\|\mathcal{P}\mathbf{U}_{r}\|_{L_{\varphi_{r}}^{2}(\Omega)}^{2}
     +\tau\|\tilde{\mathcal{P}}\tilde{\mathbf{U}}_{r}\|_{L_{\tilde{\varphi}_{r}}^{2}(\Omega)}^{2}
     +C_{*}C_{\tau}\iota^{-\frac{5}{2}}\|\mathbf{U}_{r}\|_{L_{\varphi_{r}}^{2}(\Omega)}^{2}.\nonumber
\end{align}
Combining \eqref{JJ121} and \eqref{JJ122} with \eqref{JJ12}, we have
\begin{align}\label{JJ12'}
\mathfrak{J}_{12}
\leq& \frac{1}{16}\|\sqrt{\eta}\partial_{x}\mathbf{U}_{r}\|_{L_{\varphi_{r}}^{2}(\Omega)}^{2}
     +C_{*}C_{\tau}(\sqrt{t}+\sqrt{\bar{\eta}}t^{\frac{1}{4}})\|\mathcal{P}\mathbf{U}_{r}\|_{L_{\varphi_{r}}^{2}(\Omega)}^{2}\\
    &+C_{*}C_{\tau}\iota^{-\frac{5}{2}}\|\mathbf{U}_{r}\|_{L_{\varphi_{r}}^{2}(\Omega)}^{2}
     +\tau\|\sqrt{\tilde{\eta}}\partial_{x}\tilde{\mathbf{U}}_{r}\|_{L_{\tilde{\varphi}_{r}}^{2}(\Omega)}^{2}
     +\tau\|\tilde{\mathcal{P}}\tilde{\mathbf{U}}_{r}\|_{L_{\tilde{\varphi}_{r}}^{2}(\Omega)}^{2}.\nonumber
\end{align}

Noticing that ${\bf supp}~h(y)\in(\hat{y},+\infty)$, we have
\begin{align}\label{JJ13}
\mathfrak{J}_{13}
=&(\mathcal{P}(h'\partial_{y}^{-1}[\partial_{x}^{r+1}\mathbb{U}]\partial_{y}u+h\partial_{x}^{r+1}\mathbb{U}\partial_{y}u),\varphi_{r}^{2}\mathbf{U}_{r})
 -(\mathcal{P}(h'\partial_{y}^{-1}[\partial_{x}^{r+1}\tilde{\mathbb{U}}]\partial_{y}u+h\partial_{x}^{r+1}\tilde{\mathbb{U}}\partial_{y}u),\varphi_{r}^{2}\mathbf{U}_{r})\\
\leq& C_{*}\left(\|\partial_{y}u\|_{L^{\infty}(\Omega)}+\|\partial_{y}u\|_{L_{\varphi}^{2}(\Omega)}\right)
           \|\vartheta_{c}\partial_{x}^{r+1}\mathbb{U}\|_{L_{\hat{\Psi}_{r}}^{2}(\Omega)}\|\mathbf{U}_{r}\|_{L_{\varphi_{r}}^{2}(\Omega)}\nonumber\\
    &+C_{*}\left(\|\partial_{y}u\|_{L^{\infty}(\Omega)}+\|\partial_{y}u\|_{L_{\varphi}^{2}(\Omega)}\right)
           \|\vartheta_{c}\partial_{x}^{r+1}\tilde{\mathbb{U}}\|_{L_{\tilde{\hat{\Psi}}_{r}}^{2}(\Omega)}\|\mathbf{U}_{r}\|_{L_{\varphi_{r}}^{2}(\Omega)}\nonumber\\
\leq& C_{*}(\check{c}\iota^{-\frac{1}{2}}+\sqrt{\bar{\eta}t})
           \left(\|\vartheta_{c}\partial_{x}\mathbf{U}_{r}\|_{L_{\varphi_{r}}^{2}(\Omega)}+\|\mathcal{P}\mathbf{U}_{r}\|_{L_{\varphi_{r}}^{2}(\Omega)}\right)
           \|\mathbf{U}_{r}\|_{L_{\varphi_{r}}^{2}(\Omega)}\nonumber\\
    &+C_{*}(\check{c}\iota^{-\frac{1}{2}}+\sqrt{\bar{\eta}t})
           \left(\|\vartheta_{c}\partial_{x}\tilde{\mathbf{U}}_{r}\|_{L_{\tilde{\varphi}_{r}}^{2}(\Omega)}+\|\tilde{\mathcal{P}}\tilde{\mathbf{U}}_{r}\|_{L_{\tilde{\varphi}_{r}}^{2}(\Omega)}\right)
           \|\mathbf{U}_{r}\|_{L_{\varphi_{r}}^{2}(\Omega)}\nonumber\\
\leq& \frac{1}{16}\|\sqrt{\eta}\partial_{x}\mathbf{U}_{r}\|_{L_{\varphi_{r}}^{2}(\Omega)}^{2}
     +C_{*}C_{\tau}\iota^{-1}\|\mathbf{U}_{r}\|_{L_{\varphi_{r}}^{2}(\Omega)}^{2}
     +\tau\|\sqrt{\tilde{\eta}}\partial_{x}\tilde{\mathbf{U}}_{r}\|_{L_{\tilde{\varphi}_{r}}^{2}(\Omega)}^{2}
     +\|\tilde{\mathbf{U}}_{r}\|_{L_{\tilde{\varphi}_{r}}^{2}(\Omega)}^{2}.\nonumber
\end{align}

Form \eqref{JJ11'}, \eqref{JJ12'} and \eqref{JJ13} one knows that
\begin{align}\label{JJ1}
\mathfrak{J}_{1}
\leq& \frac{3}{16}\|\sqrt{\eta}\partial_{x}\mathbf{U}_{r}\|_{L_{\varphi_{r}}^{2}(\Omega)}^{2}
     +C_{*}C_{\tau}(\sqrt{t}+\sqrt{\bar{\eta}}t^{\frac{1}{4}})\|\mathcal{P}\mathbf{U}_{r}\|_{L_{\varphi_{r}}^{2}(\Omega)}^{2}
     +C_{*}C_{\tau}\iota^{-\frac{5}{2}}\|\mathbf{U}_{r}\|_{L_{\varphi_{r}}^{2}(\Omega)}^{2}\\
    &+4\tau\|\sqrt{\tilde{\eta}}\partial_{x}\tilde{\mathbf{U}}_{r}\|_{L_{\tilde{\varphi}_{r}}^{2}(\Omega)}^{2}
     +\tau\|\tilde{\mathcal{P}}\tilde{\mathbf{U}}_{r}\|_{L_{\tilde{\varphi}_{r}}^{2}(\Omega)}^{2}
     +C\|\tilde{\mathbf{U}}_{r}\|_{L_{\tilde{\varphi}_{r}}^{2}(\Omega)}^{2}.\nonumber
\end{align}

Next let's turn to the estimate of $\mathfrak{J}_{2}$. We shall estimate $\mathfrak{J}_{2}$ in different ways for $r=s$ and $r=s+1$ respectively.
When $r=s$, by a direct calculation,  with \eqref{us-3} and Corollary \ref{ncC2} one can obtain
\begin{align}\label{JJ2(1)}
\mathfrak{J}_{2}
=& (\mathcal{P}\partial_{x}^{s}u(\partial_{x}\partial_{y}u-\rho\partial_{x}u),\varphi^{2}\mathbf{U}_{s})
  +s(\mathcal{P}\partial_{x}u\partial_{x}^{s}\partial_{y}u,\varphi^{2}\mathbf{U}_{s})
  -s(\mathcal{P}\rho\partial_{x}u\partial_{x}^{s}u,\varphi^{2}\mathbf{U}_{s})\\
 &-s(\mathcal{P}(1-h)\partial_{y}^{-1}[\partial_{x}^{s}u](\partial_{x}\partial_{y}^{2}u-\rho\partial_{x}\partial_{y}u),\varphi^{2}\mathbf{U}_{s})\nonumber\\
 &-(\mathcal{P}(1-h)\partial_{y}^{-1}[\partial_{x}u](\partial_{x}^{s}\partial_{y}^{2}u-\rho\partial_{x}^{s}\partial_{y}u),\varphi^{2}\mathbf{U}_{s})\nonumber\\
 &+s(\mathcal{P}h'\partial_{y}^{-1}[\partial_{x}^{s}u]\partial_{x}\partial_{y}u+\mathcal{P}h\partial_{x}^{s}u\partial_{x}\partial_{y}u,\varphi^{2}\mathbf{U}_{s})\nonumber\\
 &+(\mathcal{P}h'\partial_{y}^{-1}[\partial_{x}u]\partial_{x}^{s}\partial_{y}u+\mathcal{P}h\partial_{x}u\partial_{x}^{s}\partial_{y}u),\varphi^{2}\mathbf{U}_{s})\nonumber\\
\leq&~C_{*}\left(\|\partial_{x}\partial_{y}u\|_{L^{\infty}(\Omega)}+\|\partial_{x}u\|_{L^{\infty}(\Omega)}\right)
           \|\partial_{x}^{s}u\|_{L_{\hat{\Psi}}^{2}(\Omega)}\|\mathbf{U}_{s}\|_{L_{\varphi}^{2}(\Omega)}\nonumber\\
    &+C_{*}\|\partial_{x}u\|_{L^{\infty}(\Omega)}
           \left(\|\partial_{x}^{s}\partial_{y}u\|_{L_{\hat{\Psi}}^{2}(\Omega)}+\|\partial_{x}^{s}\partial_{y}^{2}u\|_{L_{\hat{\Psi}}^{2}(\Omega)}\right)
           \|\mathbf{U}_{s}\|_{L_{\varphi}^{2}(\Omega)}\nonumber\\
    &+C_{*}(\|\partial_{x}\partial_{y}^{2}u\|_{L^{\infty}(\Omega)}+\|\partial_{x}\partial_{y}u\|_{L^{\infty}(\Omega)})
           \|\mathcal{P}\varphi\|_{L_{x}^{\infty}(\mathbb{T})L_{y}^{2}(0,\hat{y})}\|\partial_{x}^{s}u\|_{L_{\hat{\psi}}^{2}(\Omega)}\|\mathcal{P}\mathbf{U}_{s}\|_{L_{\varphi}^{2}(\Omega)}\nonumber\\
\leq&~C_{*}(\check{c}+\sqrt{\bar{\eta}}t)
           \|\tilde{\mathcal{P}}\mathcal{U}_{s}\|_{L_{\tilde{\varphi}}^{2}(\Omega)}\|\mathbf{U}_{s}\|_{L_{\varphi}^{2}(\Omega)}
     +C_{*}(\check{c}\iota^{-1}+\sqrt{\bar{\eta}})
           \|\tilde{\mathcal{P}}\mathcal{U}_{s}\|_{L_{\tilde{\varphi}}^{2}(\Omega)}\|\mathbf{U}_{s}\|_{L_{\varphi}^{2}(\Omega)}\nonumber\\
    &+C_{*}(\check{c}+\sqrt{\bar{\eta}}t)
           \left(\|\mathcal{U}_{s}\|_{L_{\tilde{\varphi}}^{2}(\Omega)}
                 +\|\tilde{\mathcal{P}}\mathcal{U}_{s}\|_{L_{\tilde{\varphi}}^{2}(\Omega)}
                 +\|\partial_{y}\mathcal{U}_{s}\|_{L_{\tilde{\varphi}}^{2}(\Omega)}\right)
           \|\mathbf{U}_{s}\|_{L_{\varphi}^{2}(\Omega)}\nonumber\\
\leq&~\|\mathcal{U}_{s}\|_{L_{\tilde{\varphi}}^{2}(\Omega)}^{2}
     +\tau\|\tilde{\mathcal{P}}\mathcal{U}_{s}\|_{L_{\tilde{\varphi}}^{2}(\Omega)}
     +\tau\|\partial_{y}\mathcal{U}_{s}\|_{L_{\tilde{\varphi}}^{2}(\Omega)}^{2}
     +C_{*}C_{\tau}\iota^{-2}\|\mathbf{U}_{s}\|_{L_{\varphi}^{2}(\Omega)}^{2}.\nonumber
\end{align}
When $r=s+1$, we divide $\mathfrak{J}_{2}$ into four parts as follows,
\begin{align}\label{JJ2(2)}
\mathfrak{J}_{2}
&= (\mathcal{P}\partial_{x}^{s+1}u(\partial_{x}\partial_{y}u-\rho\partial_{x}u),\hat{\varphi}^{2}\mathbf{U}_{s+1})
  +(s+1)(\mathcal{P}\partial_{x}u(\partial_{x}^{s+1}\partial_{y}u-\rho\partial_{x}^{s+1}u),\hat{\varphi}^{2}\mathbf{U}_{s+1})\\
&\quad-(s+1)(\mathcal{P}(1-h)\partial_{y}^{-1}[\partial_{x}^{s+1}u](\partial_{x}\partial_{y}^{2}u-\rho\partial_{x}\partial_{y}u),\hat{\varphi}^{2}\mathbf{U}_{s+1})\nonumber\\
&\quad-(\mathcal{P}(1-h)\partial_{y}^{-1}[\partial_{x}u](\partial_{x}^{s+1}\partial_{y}^{2}u-\rho\partial_{x}^{s+1}\partial_{y}u),\hat{\varphi}^{2}\mathbf{U}_{s+1})\nonumber\\
&\quad+(s+1)(\mathcal{P}\partial_{x}\partial_{y}u(h'\partial_{y}^{-1}[\partial_{x}^{s+1}u]+h\partial_{x}^{s+1}u),\hat{\varphi}^{2}\mathbf{U}_{s+1})\nonumber\\
&\quad+(\mathcal{P}\partial_{x}^{s+1}\partial_{y}u(h'\partial_{y}^{-1}[\partial_{x}u]+h\partial_{x}u),\hat{\varphi}^{2}\mathbf{U}_{s+1})\nonumber\\
&=: \mathfrak{J}_{21}+\mathfrak{J}_{22}+\mathfrak{J}_{23}+\mathfrak{J}_{24}+\mathfrak{J}_{25}+\mathfrak{J}_{26}.\nonumber
\end{align}
By using \eqref{us-3} and Corollary \ref{ncC3} we have
\begin{align}\label{JJ21}
\mathfrak{J}_{21}
=&~(\mathcal{P}\partial_{x}^{s+1}u(\partial_{x}\partial_{y}u-\rho\partial_{x}u),\hat{\varphi}^{2}\mathbf{U}_{s+1})\\
\leq&~C_{*}\left(\|\vartheta_{c}^{-1}\partial_{x}\partial_{y}u\|_{L^{\infty}(\Omega)}+\|\vartheta_{c}^{-1}\partial_{x}u\|_{L^{\infty}(\Omega)}\right)
      \|\vartheta_{c}\tilde{\mathcal{P}}^{2}\partial_{x}^{s+1}u\|_{L_{\tilde{\varphi}}^{2}(\Omega)}\|\mathbf{U}_{s+1}\|_{L_{\hat{\varphi}}^{2}(\Omega)}\nonumber\\
\leq&~C_{*}(\check{c}+\sqrt{\bar{\eta}}t)
      \|\vartheta_{c}\partial_{x}^{s+1}u\|_{L_{\tilde{\hat{\Psi}}}^{2}(\Omega)}\|\mathbf{U}_{s+1}\|_{L_{\hat{\varphi}}^{2}(\Omega)}\nonumber\\
\leq&~C_{*}(\check{c}+\sqrt{\bar{\eta}}t)\|\tilde{\mathcal{P}}\|_{L^{\infty}(\Omega)}
      \left(\|\tilde{\mathcal{P}}\mathcal{U}_{s}\|_{L_{\tilde{\varphi}}^{2}(\Omega)}+\|\vartheta_{c}\partial_{x}\mathcal{U}_{s}\|_{L_{\tilde{\varphi}}^{2}(\Omega)}\right)
      \|\mathbf{U}_{s+1}\|_{L_{\hat{\varphi}}^{2}(\Omega)}\nonumber\\
\leq&~C_{*}(\check{c}\iota^{-\frac{1}{2}}+\sqrt{\bar{\eta}t})
      \left(\|\tilde{\mathcal{P}}\mathcal{U}_{s}\|_{L_{\tilde{\varphi}}^{2}(\Omega)}+\|\vartheta_{c}\partial_{x}\mathcal{U}_{s}\|_{L_{\tilde{\varphi}}^{2}(\Omega)}\right)
      \|\mathbf{U}_{s+1}\|_{L_{\hat{\varphi}}^{2}(\Omega)}\nonumber\\
\leq&~\tau\|\sqrt{\eta}\partial_{x}\mathcal{U}_{s}\|_{L_{\tilde{\varphi}}^{2}(\Omega)}^{2}
     +C_{*}C_{\tau}\iota^{-2}\|\mathbf{U}_{s+1}\|_{L_{\hat{\varphi}}^{2}(\Omega)}^{2}
     +\|\mathcal{U}_{s}\|_{L_{\tilde{\varphi}}^{2}(\Omega)}^{2}.\nonumber
\end{align}
Since $\partial_{y}\partial_{x}^{s+1}u-\tilde{\rho}\partial_{x}^{s+1}u=\partial_{x}(\tilde{\mathcal{P}}^{-1}\mathcal{U}_{s})+\partial_{x}\tilde{\rho}\partial_{x}^{s}u$,
similar to \eqref{JJ21} we have
\begin{align}\label{JJ22}
\mathfrak{J}_{22}
=&~(s+1)(\mathcal{P}\partial_{x}u(\partial_{x}^{s+1}\partial_{y}u-\rho\partial_{x}^{s+1}u),\hat{\varphi}^{2}\mathbf{U}_{s+1})\\
=&~(s+1)(\mathcal{P}\partial_{x}u\partial_{x}(\tilde{\mathcal{P}}^{-1}\mathcal{U}_{s}),\hat{\varphi}^{2}\mathbf{U}_{s+1})
  +(s+1)(\mathcal{P}\partial_{x}u\partial_{x}\tilde{\rho}\partial_{x}^{s}u,\hat{\varphi}^{2}\mathbf{U}_{s+1})\nonumber\\
 &+(s+1)(\mathcal{P}\partial_{x}u(\tilde{\rho}-\rho)\partial_{x}^{s+1}u,\hat{\varphi}^{2}\mathbf{U}_{s+1})\nonumber\\
\leq&~C\|\vartheta_{c}^{-1}\partial_{x}u\|_{L^{\infty}(\Omega)}
       \|\vartheta_{c}\partial_{x}\mathcal{U}_{s}\|_{L_{\varphi}^{2}(\Omega)}
       \|\mathbf{U}_{s+1}\|_{L_{\hat{\varphi}}^{2}(\Omega)}\nonumber\\
    &+C_{*}\|\partial_{x}u\|_{L^{\infty}(\Omega)}
           \|\tilde{\mathcal{P}}\mathcal{U}_{s}\|_{L_{\tilde{\varphi}}^{2}(\Omega)}
           \|\mathbf{U}_{s+1}\|_{L_{\hat{\varphi}}^{2}(\Omega)}\nonumber\\
    &+C_{*}\|\partial_{x}u\|_{L^{\infty}(\Omega)}
           \|\partial_{x}^{s}u\|_{L_{\hat{\Psi}}^{2}(\Omega)}
           \|\mathbf{U}_{s+1}\|_{L_{\hat{\varphi}}^{2}(\Omega)}\nonumber\\
    &+C_{*}\|\vartheta_{c}^{-1}\partial_{x}u\|_{L^{\infty}(\Omega)}
           \|\vartheta_{c}\tilde{\mathcal{P}}^{2}\partial_{x}^{s+1}u\|_{L_{\tilde{\hat{\varphi}}}^{2}(\Omega)}
           \|\mathbf{U}_{s+1}\|_{L_{\hat{\varphi}}^{2}(\Omega)}\nonumber\\
\leq&~C_{*}(\check{c}+\sqrt{\bar{\eta}}t)
      \left(\|\tilde{\mathcal{P}}\mathcal{U}_{s}\|_{L_{\tilde{\varphi}}^{2}(\Omega)}
            +\|\vartheta_{c}\partial_{x}\mathcal{U}_{s}\|_{L_{\tilde{\varphi}}^{2}(\Omega)}\right)
      \|\mathbf{U}_{s+1}\|_{L_{\hat{\varphi}}^{2}(\Omega)}\nonumber\\
    &+C_{*}(\check{c}+\sqrt{\bar{\eta}}t)
      \|\vartheta_{c}\partial_{x}^{s+1}u\|_{L_{\hat{\Psi}}^{2}(\Omega)}\|\mathbf{U}_{s+1}\|_{L_{\tilde{\hat{\varphi}}}^{2}(\Omega)}\nonumber\\
\leq&~\tau\|\sqrt{\eta}\partial_{x}\mathcal{U}_{s}\|_{L_{\tilde{\varphi}}^{2}(\Omega)}^{2}
     +C_{*}C_{\tau}\iota^{-2}\|\mathbf{U}_{s+1}\|_{L_{\hat{\varphi}}^{2}(\Omega)}^{2}
     +\|\mathcal{U}_{s}\|_{L_{\tilde{\varphi}}^{2}(\Omega)}^{2}.\nonumber
\end{align}
Similar to that in \eqref{JJ121}, one has
\begin{align}\label{JJ23}
\mathfrak{J}_{23}
=&-(s+1)(\mathcal{P}(1-h)\partial_{y}^{-1}[\partial_{x}^{s+1}u](\partial_{x}\partial_{y}^{2}u-\rho\partial_{x}\partial_{y}u),\hat{\varphi}^{2}\mathbf{U}_{s+1})\\
=&~(s+1)\int_{\mathbb{T}}\int_{\check{y}}^{\hat{y}}\hat{\varphi}^{2}\mathcal{P}(1-h)\partial_{y}^{-1}[\partial_{x}^{s+1}u](\rho\partial_{x}\partial_{y}u-\partial_{x}\partial_{y}^{2}u)\mathbf{U}_{s+1}dydx\nonumber\\
 &+(s+1)\int_{\mathbb{T}}\int_{\mathbb{R}_{+}\setminus[\check{y},\hat{y}]}\hat{\varphi}^{2}\mathcal{P}(1-h)\partial_{y}^{-1}[\partial_{x}^{s+1}u](\rho\partial_{x}\partial_{y}u-\partial_{x}\partial_{y}^{2}u)\mathbf{U}_{s+1}dydx\nonumber\\
\leq&~C_{*}\left(\|\partial_{x}\partial_{y}u\|_{L^{\infty}(\mathbb{T}\times[\check{y},\hat{y}])}
                 +\|\partial_{x}\partial_{y}^{2}u\|_{L^{\infty}(\mathbb{T}\times[\check{y},\hat{y}])}\right)
          \|\vartheta_{c}\partial_{x}^{s+1}u\|_{L_{\hat{\psi}}^{2}(\Omega)}
          \|\mathcal{P}\hat{\varphi}\|_{L_{x}^{\infty}(\mathbb{T})L_{y}^{2}([\check{y},\hat{y}])}
          \|\mathcal{P}\mathbf{U}_{s+1}\|_{L_{\hat{\varphi}}^{2}(\Omega)}\nonumber\\
    &+C_{*}\left(\|\varphi\partial_{x}\partial_{y}u\|_{L_{x}^{\infty}L_{y}^{2}(\Omega)}+\|\varphi\partial_{x}\partial_{y}^{2}u\|_{L_{x}^{\infty}L_{y}^{2}(\Omega)}\right)
          \|\vartheta_{c}\partial_{x}^{s+1}u\|_{L_{\tilde{\hat{\Psi}}}^{2}(\Omega)}\|\mathbf{U}_{s+1}\|_{L_{\hat{\varphi}}^{2}(\Omega)}\nonumber\\
\leq&~C_{*}(\check{c}+\sqrt{\bar{\eta}}t)
          \left(\int_{\check{y}_{*}}^{\hat{y}_{*}}\frac{1}{((y-y_{*})^2+\bar{\varsigma}+t)^{1+\varepsilon_{0}}}dy\right)^{\frac{1}{2}}
          \|\vartheta_{c}\partial_{x}^{s+1}u\|_{L_{\hat{\psi}}^{2}(\Omega)}
          \|\mathcal{P}\mathbf{U}_{s+1}\|_{L_{\hat{\varphi}}^{2}(\Omega)}\nonumber\\
    &+C_{*}(\check{c}+\sqrt{\bar{\eta}}t)\|\vartheta_{c}\partial_{x}^{s+1}u\|_{L_{\tilde{\hat{\Psi}}}^{2}(\Omega)}
          \|\mathbf{U}_{s+1}\|_{L_{\hat{\varphi}}^{2}(\Omega)}\nonumber\\
\leq&~C_{*}(\check{c}\iota^{-\frac{3}{4}}+\sqrt{\bar{\eta}}t^{\frac{1}{4}})
      \left(\|\vartheta_{c}\partial_{x}\mathcal{U}_{s}\|_{L_{\tilde{\varphi}}^{2}(\Omega)}
            +\|\tilde{\mathcal{P}}\mathcal{U}_{s}\|_{L_{\tilde{\varphi}}^{2}(\Omega)}\right)
      \|\mathcal{P}\mathbf{U}_{s+1}\|_{L_{\hat{\varphi}}^{2}(\Omega)}\nonumber\\
    &+C_{*}(\check{c}\iota^{-\frac{1}{2}}+\sqrt{\bar{\eta}t})
      \left(\|\vartheta_{c}\partial_{x}\mathcal{U}_{s}\|_{L_{\tilde{\varphi}}^{2}(\Omega)}
            +\|\tilde{\mathcal{P}}\mathcal{U}_{s}\|_{L_{\tilde{\varphi}}^{2}(\Omega)}\right)
      \|\mathbf{U}_{s+1}\|_{L_{\hat{\varphi}}^{2}(\Omega)}\nonumber\\
\leq&~\tau\|\sqrt{\eta}\partial_{x}\mathcal{U}_{s}\|_{L_{\tilde{\varphi}}^{2}(\Omega)}^{2}
     +\tau\|\tilde{\mathcal{P}}\mathcal{U}_{s}\|_{L_{\tilde{\varphi}}^{2}(\Omega)}^{2}
     +C_{*}C_{\tau}\sqrt{t}\|\mathcal{P}\mathbf{U}_{s+1}\|_{L_{\hat{\varphi}}^{2}(\Omega)}^{2}
     +C_{*}C_{\tau}\iota^{-\frac{5}{2}}\|\mathbf{U}_{s+1}\|_{L_{\hat{\varphi}}^{2}(\Omega)}^{2}.\nonumber
\end{align}

By using $u=\mathbb{U}-\tilde{\mathbb{U}}$, we decompose $\mathfrak{J}_{24}$ into two parts as follows,
\begin{align}\label{JJ24}
\mathfrak{J}_{24}
=&-(\mathcal{P}(1-h)\partial_{y}^{-1}[\partial_{x}u](\partial_{x}^{s+1}\partial_{y}^{2}u-\rho\partial_{x}^{s+1}\partial_{y}u),\hat{\varphi}^{2}\mathbf{U}_{s+1})\\
=&~(\mathcal{P}(1-h)\partial_{y}^{-1}[\partial_{x}u](\rho\partial_{x}^{s+1}\partial_{y}\mathbb{U}-\partial_{x}^{s+1}\partial_{y}^{2}\mathbb{U}),\hat{\varphi}^{2}\mathbf{U}_{s+1})\nonumber\\
 &-(\mathcal{P}(1-h)\partial_{y}^{-1}[\partial_{x}u](\rho\partial_{x}^{s+1}\partial_{y}\tilde{\mathbb{U}}-\partial_{x}^{s+1}\partial_{y}^{2}\tilde{\mathbb{U}}),\hat{\varphi}^{2}\mathbf{U}_{s+1}).\nonumber
\end{align}
Since $\partial_{x}^{s+1}\partial_{y}^{2}\mathbb{U}-\rho\partial_{x}^{s+1}\partial_{y}\mathbb{U}=\partial_{y}(\mathcal{P}^{-1}\mathbf{U}_{s+1})+\partial_{y}\rho\partial_{x}^{s+1}\mathbb{U}$,
one has
\begin{align}\label{JJ241}
&(\mathcal{P}(1-h)\partial_{y}^{-1}[\partial_{x}u](\rho\partial_{x}^{s+1}\partial_{y}\mathbb{U}-\partial_{x}^{s+1}\partial_{y}^{2}\mathbb{U}),\hat{\varphi}^{2}\mathbf{U}_{s+1})\\
=&-(\mathcal{P}(1-h)\partial_{y}^{-1}[\partial_{x}u]\partial_{y}\mathcal{P}^{-1}\mathbf{U}_{s+1},\hat{\varphi}^{2}\mathbf{U}_{s+1})
  -((1-h)\partial_{y}^{-1}[\partial_{x}u]\partial_{y}\mathbf{U}_{s+1},\hat{\varphi}^{2}\mathbf{U}_{s+1})\nonumber\\
 &-(\mathcal{P}(1-h)\partial_{y}^{-1}[\partial_{x}u]\partial_{y}\rho\partial_{x}^{s+1}\mathbb{U},\hat{\varphi}^{2}\mathbf{U}_{s+1})\nonumber\\
\leq&~C_{*}\|\psi\partial_{x}u\|_{L_{x}^{\infty}L_{y}^{2}(\Omega)}
           \|\mathcal{P}\mathbf{U}_{s+1}\|_{L_{\hat{\varphi}}^{2}(\Omega)}
           \|\mathbf{U}_{s+1}\|_{L_{\hat{\varphi}}^{2}(\Omega)}
     +C\|\psi\partial_{x}u\|_{L_{x}^{\infty}L_{y}^{2}(\Omega)}
       \|\partial_{y}\mathbf{U}_{s+1}\|_{L_{\hat{\varphi}}^{2}(\Omega)}
       \|\mathbf{U}_{s+1}\|_{L_{\hat{\varphi}}^{2}(\Omega)}\nonumber\\
    &+C_{*}\|\psi\partial_{x}u\|_{L_{x}^{\infty}L_{y}^{2}(\Omega)}
           \|\mathcal{P}^{3}\partial_{x}^{s+1}\mathbb{U}\|_{L_{\hat{\varphi}}^{2}(\Omega)}
           \|\mathbf{U}_{s+1}\|_{L_{\hat{\varphi}}^{2}(\Omega)}\nonumber\\
\leq&~C_{*}\iota^{-\frac{1}{2}}\|\mathbf{U}_{s+1}\|_{L_{\hat{\varphi}}^{2}(\Omega)}^{2}
     +\frac{1}{16}\|\partial_{y}\mathbf{U}_{s+1}\|_{L_{\hat{\varphi}}^{2}(\Omega)}^{2}
     +C_{*}(\check{c}\iota^{-\frac{1}{2}}+\sqrt{\bar{\eta}t})
      \|\partial_{x}^{s+1}\mathbb{U}\|_{L_{\hat{\Psi}_{s+1}}^{2}(\Omega)}\|\mathbf{U}_{s+1}\|_{L_{\hat{\varphi}}^{2}(\Omega)}\nonumber\\
\leq&~C_{*}\iota^{-\frac{1}{2}}\|\mathbf{U}_{s+1}\|_{L_{\hat{\varphi}}^{2}(\Omega)}^{2}
     +\frac{1}{16}\|\partial_{y}\mathbf{U}_{s+1}\|_{L_{\hat{\varphi}}^{2}(\Omega)}^{2}
     +C_{*}(\check{c}\iota^{-\frac{1}{2}}+\sqrt{\bar{\eta}t})
      \|\mathcal{P}\mathbf{U}_{s+1}\|_{L_{\hat{\varphi}}^{2}(\Omega)}\|\mathbf{U}_{s+1}\|_{L_{\hat{\varphi}}^{2}(\Omega)}\nonumber\\
\leq&~C_{*}\iota^{-1}\|\mathbf{U}_{s+1}\|_{L_{\hat{\varphi}}^{2}(\Omega)}^{2}
     +\frac{1}{16}\|\partial_{y}\mathbf{U}_{s+1}\|_{L_{\hat{\varphi}}^{2}(\Omega)}^{2}.\nonumber
\end{align}
Parallel to \eqref{JJ241}, one has
\begin{align}\label{JJ242}
&(\mathcal{P}(1-h)\partial_{y}^{-1}[\partial_{x}u](\rho\partial_{x}^{s+1}\partial_{y}\tilde{\mathbb{U}}-\partial_{x}^{s+1}\partial_{y}^{2}\tilde{\mathbb{U}}),\hat{\varphi}^{2}\mathbf{U}_{s+1})\\
=&-(\mathcal{P}(1-h)\partial_{y}^{-1}[\partial_{x}u]\partial_{y}\mathcal{P}^{-1}\tilde{\mathbf{U}}_{s+1},\hat{\varphi}^{2}\mathbf{U}_{s+1})
  -((1-h)\partial_{y}^{-1}[\partial_{x}u]\partial_{y}\tilde{\mathbf{U}}_{s+1},\hat{\varphi}^{2}\mathbf{U}_{s+1})\nonumber\\
 &-(\mathcal{P}(1-h)\partial_{y}^{-1}[\partial_{x}u]\partial_{y}\rho\partial_{x}^{s+1}\tilde{\mathbb{U}},\hat{\varphi}^{2}\mathbf{U}_{s+1})\nonumber\\
\leq&~C_{*}\|\psi\partial_{x}u\|_{L_{x}^{\infty}L_{y}^{2}(\Omega)}
           \|\mathcal{P}\mathbf{U}_{s+1}\|_{L_{\hat{\varphi}}^{2}(\Omega)}
           \|\tilde{\mathbf{U}}_{s+1}\|_{L_{\tilde{\hat{\varphi}}}^{2}(\Omega)}
     +C\|\psi\partial_{x}u\|_{L_{x}^{\infty}L_{y}^{2}(\Omega)}
       \|\partial_{y}\tilde{\mathbf{U}}_{s+1}\|_{L_{{\tilde{\hat\varphi}}}^{2}(\Omega)}
       \|\mathbf{U}_{s+1}\|_{L_{\hat{\varphi}}^{2}(\Omega)}\nonumber\\
    &+C_{*}\|\psi\partial_{x}u\|_{L_{x}^{\infty}L_{y}^{2}(\Omega)}
           \|\tilde{\mathcal{P}}^{2}\partial_{x}^{s+1}\tilde{\mathbb{U}}\|_{L_{\tilde{\hat{\varphi}}}^{2}(\Omega)}
           \|\mathcal{P}\mathbf{U}_{s+1}\|_{L_{\hat{\varphi}}^{2}(\Omega)}\nonumber\\
\leq&~C_{*}(\check{c}\iota^{-\frac{1}{2}}+\sqrt{\bar{\eta}t})
            \|\mathbf{U}_{s+1}\|_{L_{\hat{\varphi}}^{2}(\Omega)}\|\tilde{\mathbf{U}}_{s+1}\|_{L_{\tilde{\hat{\varphi}}}^{2}(\Omega)}
     +\frac{1}{16}\|\partial_{y}\tilde{\mathbf{U}}_{s+1}\|_{L_{\tilde{\hat{\varphi}}}^{2}(\Omega)}^{2}
     +C_{*}\|\mathbf{U}_{s+1}\|_{L_{\hat{\varphi}}^{2}(\Omega)}^{2}\nonumber\\
    &+C_{*}(\check{c}+\sqrt{\bar{\eta}}t)
      \|\partial_{x}^{s+1}\tilde{\mathbb{U}}\|_{L_{\tilde{\hat{\Psi}}_{s+1}}^{2}(\Omega)}\|\mathcal{P}\mathbf{U}_{s+1}\|_{L_{\hat{\varphi}}^{2}(\Omega)}\nonumber\\
\leq&~C_{*}\iota^{-1}\|\mathbf{U}_{s+1}\|_{L_{\hat{\varphi}}^{2}(\Omega)}^{2}
     +\frac{1}{2}\|\tilde{\mathbf{U}}_{s+1}\|_{L_{\tilde{\hat\varphi}}^{2}(\Omega)}^{2}
     +\frac{1}{16}\|\partial_{y}\tilde{\mathbf{U}}_{s+1}\|_{L_{\tilde{\hat{\varphi}}}^{2}(\Omega)}^{2}\nonumber\\
    &+C_{*}(\check{c}+\sqrt{\bar{\eta}}t)
      \|\tilde{\mathcal{P}}\tilde{\mathbf{U}}_{s+1}\|_{L_{\tilde{\hat{\varphi}}}^{2}(\Omega)}\|\mathcal{P}\mathbf{U}_{s+1}\|_{L_{\hat{\varphi}}^{2}(\Omega)}\nonumber\\
\leq&~C_{*}\iota^{-2}\|\mathbf{U}_{s+1}\|_{L_{\hat{\varphi}}^{2}(\Omega)}^{2}
     +\|\tilde{\mathbf{U}}_{s+1}\|_{L_{\tilde{\hat{\varphi}}}^{2}(\Omega)}^{2}
     +\frac{1}{16}\|\partial_{y}\mathbf{U}_{s+1}\|_{L_{\hat{\varphi}}^{2}(\Omega)}^{2}.\nonumber
\end{align}
Substituting \eqref{JJ241} and \eqref{JJ242} into \eqref{JJ24}, one has
\begin{align}\label{JJ24'}
\mathfrak{J}_{24}
\leq C_{*}\iota^{-2}\|\mathbf{U}_{s+1}\|_{L_{\hat{\varphi}}^{2}(\Omega)}^{2}
    +\|\tilde{\mathbf{U}}_{s+1}\|_{L_{\tilde{\hat{\varphi}}}^{2}(\Omega)}^{2}
    +\frac{1}{8}\|\partial_{y}\mathbf{U}_{s+1}\|_{L_{\hat{\varphi}}^{2}(\Omega)}^{2}.
\end{align}

Similar to \eqref{JJ21}, one has
\begin{align}\label{JJ25}
\mathfrak{J}_{25}
=&~(s+1)(\mathcal{P}\partial_{x}\partial_{y}u(h'\partial_{y}^{-1}[\partial_{x}^{s+1}u]+h\partial_{x}^{s+1}u),\hat{\varphi}^{2}\mathbf{U}_{s+1})\\
\leq&~C_{*}\left(\|\partial_{x}\partial_{y}u\|_{L^{\infty}(\Omega)}+\|\partial_{x}\partial_{y}u\|_{L_{x}^{\infty}L_{y}^{2}(\Omega)}\right)
           \|\vartheta_{c}\partial_{x}^{s+1}u\|_{L_{\tilde{\hat{\Psi}}}^{2}(\Omega)}\|\mathbf{U}_{s+1}\|_{L_{\hat{\varphi}}^{2}(\Omega)}\nonumber\\
\leq&~C_{*}(\check{c}+\sqrt{\bar{\eta}}t)
      \|\vartheta_{c}\partial_{x}^{s+1}u\|_{L_{\tilde{\hat{\Psi}}}^{2}(\Omega)}\|\mathbf{U}_{s+1}\|_{L_{\hat{\varphi}}^{2}(\Omega)}\nonumber\\
\leq&~C_{*}(\check{c}+\sqrt{\bar{\eta}}t)\|\tilde{\mathcal{P}}\|_{L^{\infty}(\Omega)}
      \left(\|\tilde{\mathcal{P}}\mathcal{U}_{s}\|_{L_{\tilde{\varphi}}^{2}(\Omega)}+\|\vartheta_{c}\partial_{x}\mathcal{U}_{s}\|_{L_{\tilde{\varphi}}^{2}(\Omega)}\right)
      \|\mathbf{U}_{s+1}\|_{L_{\hat{\varphi}}^{2}(\Omega)}\nonumber\\
\leq&~C_{*}(\check{c}\iota^{-\frac{1}{2}}+\sqrt{\bar{\eta}t})
      \left(\|\tilde{\mathcal{P}}\mathcal{U}_{s}\|_{L_{\tilde{\varphi}}^{2}(\Omega)}+\|\vartheta_{c}\partial_{x}\mathcal{U}_{s}\|_{L_{\tilde{\varphi}}^{2}(\Omega)}\right)
      \|\mathbf{U}_{s+1}\|_{L_{\hat{\varphi}}^{2}(\Omega)}\nonumber\\
\leq&~\tau\|\sqrt{\eta}\partial_{x}\mathcal{U}_{s}\|_{L_{\tilde{\varphi}}^{2}(\Omega)}^{2}
     +C_{*}C_{\tau}\iota^{-2}\|\mathbf{U}_{s+1}\|_{L_{\hat{\varphi}}^{2}(\Omega)}^{2}
     +\|\mathcal{U}_{s}\|_{L_{\tilde{\varphi}}^{2}(\Omega)}^{2}.\nonumber
\end{align}

Since $\partial_{y}\partial_{x}^{s+1}u-\tilde{\rho}\partial_{x}^{s+1}u=\partial_{x}(\tilde{\mathcal{P}}^{-1}\mathcal{U}_{s})+\partial_{x}\tilde{\rho}\partial_{x}^{s}u$, one has
\begin{align*}
\|\vartheta_{c}\partial_{y}\partial_{x}^{s+1}u\|_{L_{\varphi}^{2}(\Omega)}
\leq& C_{*}\|\vartheta_{c}\partial_{x}^{s+1}u\|_{L_{\tilde{\hat{\Psi}}}^{2}(\Omega)}
    +C_{*}\|\vartheta_{c}\partial_{x}^{s}u\|_{L_{\tilde{\hat{\Psi}}}^{2}(\Omega)}
    +\|\vartheta_{c}\partial_{x}(\mathcal{P}^{-1}\mathcal{U}_{s})\|_{L_{\varphi}^{2}(\Omega)}\\
\leq& C_{*}\|\tilde{\mathcal{P}}\|_{L^{\infty}(\Omega)}
           \left(\|\tilde{\mathcal{P}}\mathcal{U}_{s}\|_{L_{\tilde{\varphi}}^{2}(\Omega)}+\|\vartheta_{c}\partial_{x}\mathcal{U}_{s}\|_{L_{\tilde{\varphi}}^{2}(\Omega)}\right)
     +C_{*}\|\tilde{\mathcal{P}}\mathcal{U}_{s}\|_{L_{\tilde{\varphi}}^{2}(\Omega)}\\
    &+C_{*}\|\mathcal{U}_{s}\|_{L_{\tilde{\varphi}}^{2}(\Omega)}
     +C_{*}\|\vartheta_{c}\partial_{x}\mathcal{U}_{s}\|_{L_{\tilde{\varphi}}^{2}(\Omega)},
\end{align*}
moreover, $\vartheta_{c}^{-1}\leq C\psi$ owing to $t\leq \min\{\frac{1}{12\Lambda},\frac{1}{4\theta}\}$, thus
\begin{align}\label{JJ26}
\mathfrak{J}_{26}
\leq&~C_{*}\left(\|\vartheta_{c}^{-1}\partial_{x}u\|_{L^{\infty}(\Omega)}
                 +\|\vartheta_{c}^{-1}\partial_{x}u\|_{L_{x}^{\infty}L_{y}^{2}(\Omega)}\right)
          \|\vartheta_{c}\partial_{y}\partial_{x}^{s+1}u\|_{L_{\varphi}^{2}(\Omega)}\|\mathbf{U}_{s+1}\|_{L_{\hat{\varphi}}^{2}(\Omega)}\\
\leq&~C_{*}(\check{c}\iota^{-1}+\sqrt{\bar{\eta}})
           \left(\|\mathcal{U}_{s}\|_{L_{\tilde{\varphi}}^{2}(\Omega)}+\|\vartheta_{c}\partial_{x}\mathcal{U}_{s}\|_{L_{\tilde{\varphi}}^{2}(\Omega)}\right)
           \|\mathbf{U}_{s+1}\|_{L_{\hat{\varphi}}^{2}(\Omega)}\nonumber\\
\leq&~\tau\|\sqrt{\eta}\partial_{x}\mathcal{U}_{s}\|_{L_{\varphi}^{2}(\Omega)}^{2}
     +C_{*}C_{\tau}\iota^{-2}\|\mathbf{U}_{s+1}\|_{L_{\hat{\varphi}}^{2}(\Omega)}^{2}
     +\|\mathcal{U}_{s}\|_{L_{\tilde{\varphi}}^{2}(\Omega)}^{2}.\nonumber
\end{align}

Combining \eqref{JJ21},\eqref{JJ22},\eqref{JJ23},\eqref{JJ24'},\eqref{JJ25} and \eqref{JJ26} with \eqref{JJ2(2)}, one can conclude that when $r=s+1$,
\begin{align*}
\mathfrak{J}_{2}
\leq&~\frac{1}{8}\|\partial_{y}\mathbf{U}_{s+1}\|_{L_{\hat{\varphi}}^{2}(\Omega)}^{2}
     +C_{*}C_{\tau}\iota^{-\frac{5}{2}}\|\mathbf{U}_{s+1}\|_{L_{\hat{\varphi}}^{2}(\Omega)}^{2}
     +C_{*}C_{\tau}\sqrt{t}\|\mathcal{P}\mathbf{U}_{s+1}\|_{L_{\hat{\varphi}}^{2}(\Omega)}^{2}\\
    &+\|\tilde{\mathbf{U}}_{s+1}\|_{L_{\tilde{\hat{\varphi}}}^{2}(\Omega)}^{2}
     +C\|\mathcal{U}_{s}\|_{L_{\tilde{\varphi}}^{2}(\Omega)}^{2}
     +5\tau\|\sqrt{\eta}\partial_{x}\mathcal{U}_{s}\|_{L_{\tilde{\varphi}}^{2}(\Omega)}^{2}
     +\tau\|\tilde{\mathcal{P}}\mathcal{U}_{s}\|_{L_{\tilde{\varphi}}^{2}(\Omega)}^{2}.
\end{align*}

Now let's focus on the estimate of $\mathfrak{J}_{3}$. When $r=s+1$, the estimate of $\mathfrak{J}_{3}$ is similar to that in \eqref{JJ2(1)},
\begin{align*}
\mathfrak{J}_{3}
=&~(s+1)(\mathcal{P}\partial_{x}^{s}u(\partial_{x}^{2}\partial_{y}u-\rho\partial_{x}^{2}u),\hat{\varphi}^{2}\mathbf{U}_{s+1})
  +\frac{s(s+1)}{2}(\partial_{x}^{2}u\mathcal{U}_{s},\hat{\varphi}^{2}\mathbf{U}_{s+1})\\
 &-\frac{s(s+1)}{2}(\mathcal{P}(1-h)\partial_{y}^{-1}[\partial_{x}^{s}u](\partial_{x}^{2}\partial_{y}^{2}u-\rho\partial_{x}^{2}\partial_{y}u),\hat{\varphi}^{2}\mathbf{U}_{s+1})\nonumber\\
 &+\frac{s(s+1)}{2}(\mathcal{P}\partial_{x}^{2}\partial_{y}u(h'\partial_{y}^{-1}[\partial_{x}^{s}u]+h\partial_{x}^{s}u),\hat{\varphi}^{2}\mathbf{U}_{s+1})\nonumber\\
 &-(s+1)(\mathcal{P}(1-h)\partial_{y}^{-1}[\partial_{x}^{2}u][\partial_{y}(\mathcal{P}^{-1}\mathcal{U}_{s})+\partial_{y}\rho\partial_{x}^{s}u],\hat{\varphi}^{2}\mathbf{U}_{s+1})\nonumber\\
 &+(s+1)(\mathcal{P}\partial_{x}^{s}\partial_{y}u(h'\partial_{y}^{-1}[\partial_{x}^{2}u]+h\partial_{x}^{2}u),\hat{\varphi}^{2}\mathbf{U}_{s+1})\nonumber\\
\leq&~C_{*}(\|\partial_{x}^{2}\partial_{y}u\|_{L^{\infty}(\Omega)}+\|\partial_{x}^{2}u\|_{L^{\infty}(\Omega)})
           \|\partial_{x}^{s}u\|_{L_{\hat{\Psi}}^{2}(\Omega)}\|\mathbf{U}_{s+1}\|_{L_{\hat{\varphi}}^{2}(\Omega)}
     +C\|\partial_{x}^{2}u\|_{L^{\infty}(\Omega)}\|\mathcal{U}_{s}\|_{L_{\tilde{\varphi}}^{2}(\Omega)}\|\mathbf{U}_{s+1}\|_{L_{\hat{\varphi}}^{2}(\Omega)}\nonumber\\
    &+C_{*}(\|\partial_{x}^{2}\partial_{y}^{2}u\|_{L^{\infty}(\Omega)}+\|\partial_{x}^{2}\partial_{y}u\|_{L^{\infty}(\Omega)})
           \|\mathcal{P}\varphi\|_{L_{x}^{\infty}(\mathbb{T})L_{y}^{2}(0,\hat{y})}
           \|\partial_{x}^{s}u\|_{L_{\hat{\psi}}^{2}(\Omega)}\|\mathcal{P}\mathbf{U}_{s+1}\|_{L_{\hat{\varphi}}^{2}(\Omega)}\nonumber\\
    &+C_{*}\|\psi\partial_{x}^{2}u\|_{L_{x}^{\infty}L_{y}^{2}(\Omega)}\|\partial_{y}\mathcal{U}_{s}\|_{L_{\tilde{\varphi}}^{2}(\Omega)}
           \|\mathbf{U}_{s+1}\|_{L_{\hat{\varphi}}^{2}(\Omega)}
     +C_{*}\|\psi\partial_{x}^{2}u\|_{L_{x}^{\infty}L_{y}^{2}(\Omega)}\|\tilde{\mathcal{P}}\mathcal{U}_{s}\|_{L_{\tilde{\varphi}}^{2}(\Omega)}
           \|\mathbf{U}_{s+1}\|_{L_{\hat{\varphi}}^{2}(\Omega)}\nonumber\\
    &+C_{*}\|\psi\partial_{x}^{2}u\|_{L_{x}^{\infty}L_{y}^{2}(\Omega)}\|\partial_{x}^{s}u\|_{L_{\hat{\Psi}}^{2}(\Omega)}
           \|\mathcal{P}\mathbf{U}_{s+1}\|_{L_{\hat{\varphi}}^{2}(\Omega)}\nonumber\\
    &+C_{*}\left(\|\psi\partial_{x}^{2}u\|_{L_{x}^{\infty}L_{y}^{2}(\Omega)}+\|\partial_{x}^{2}u\|_{L^{\infty}(\Omega)}\right)
           \|\partial_{y}\partial_{x}^{s}u\|_{L_{\varphi}^{2}(\Omega)}\|\mathbf{U}_{s+1}\|_{L_{\check{\varphi}}^{2}(\Omega)}\\
\leq&~C_{*}(\check{c}+\sqrt{\bar{\eta}}t)\|\tilde{\mathcal{P}}\mathcal{U}_{s}\|_{L_{\tilde{\varphi}}^{2}(\Omega)}
                                             \|\mathbf{U}_{s+1}\|_{L_{\hat{\varphi}}^{2}(\Omega)}
     +C_{*}\|\mathcal{U}_{s}\|_{L_{\tilde{\varphi}}^{2}(\Omega)}\|\mathbf{U}_{s+1}\|_{L_{\hat{\varphi}}^{2}(\Omega)}\nonumber\\
    &+C_{*}(\check{c}\iota^{-1}+\sqrt{\bar{\eta}})\|\tilde{\mathcal{P}}\mathcal{U}_{s}\|_{L_{\tilde{\varphi}}^{2}(\Omega)}
                                                  \|\mathbf{U}_{s+1}\|_{L_{\hat{\varphi}}^{2}(\Omega)}
     +C_{*}\|\partial_{y}\mathcal{U}_{s}\|_{L_{\tilde{\varphi}}^{2}(\Omega)}\|\mathbf{U}_{s+1}\|_{L_{\hat{\varphi}}^{2}(\Omega)}\nonumber\\
\leq&~\tau\|\tilde{\mathcal{P}}\mathcal{U}_{s}\|_{L_{\tilde{\varphi}}^{2}(\Omega)}
     +\tau\|\partial_{y}\mathcal{U}_{s}\|_{L_{\tilde{\varphi}}^{2}(\Omega)}^{2}
     +C\|\mathcal{U}_{s}\|_{L_{\tilde{\varphi}}^{2}(\Omega)}^{2}
     +C_{*}C_{\tau}\iota^{-2}\|\mathbf{U}_{s+1}\|_{L_{\hat{\varphi}}^{2}(\Omega)}^{2},\nonumber
\end{align*}
while when $r=s$, since $\bar{\varsigma}=\iota\bar{\eta}^{2},\varepsilon_{0}<\frac{1}{4}$ and $\bar{\varsigma}=\iota\epsilon_{0}^{4}$ when $\check{c}\neq 0$, by noticing \eqref{us-3} one gets,
\begin{align*}
\mathfrak{J}_{3}
=&~s(\mathcal{P}\partial_{x}^{s-1}u(\partial_{x}^{2}\partial_{y}u-\rho\partial_{x}^{2}u),\varphi^{2}\mathbf{U}_{s})
  +\frac{s(s-1)}{2}(\mathcal{P}\partial_{x}^{2}u(\partial_{x}^{s-1}\partial_{y}u-\rho\partial_{x}^{s-1}u),\varphi^{2}\mathbf{U}_{s})\\
 &-\frac{s(s-1)}{2}(\mathcal{P}(1-h)\partial_{y}^{-1}[\partial_{x}^{s-1}u](\partial_{x}^{2}\partial_{y}^{2}u-\rho\partial_{x}^{2}\partial_{y}u),\varphi^{2}\mathbf{U}_{s})\\
 &-s(\mathcal{P}(1-h)\partial_{y}^{-1}[\partial_{x}^{2}u](\partial_{x}^{s-1}\partial_{y}^{2}u-\rho\partial_{x}^{s-1}\partial_{y}u),\varphi^{2}\mathbf{U}_{s})\nonumber\\
 &+\frac{s(s-1)}{2}(\mathcal{P}\partial_{x}^{2}\partial_{y}u(h'\partial_{y}^{-1}[\partial_{x}^{s-1}u]+h\partial_{x}^{s-1}u),\varphi^{2}\mathbf{U}_{s})\\
 &+s(\mathcal{P}\partial_{x}^{s-1}\partial_{y}u(h'\partial_{y}^{-1}[\partial_{x}^{2}u]+h\partial_{x}^{2}u),\varphi^{2}\mathbf{U}_{s})\nonumber\\
\leq&~C_{*}(\|\mathcal{P}\varphi\|_{L^{\infty}(\Omega)}\|\partial_{x}^{2}\partial_{y}u\|_{L^{\infty}(\Omega)}
            +\|\mathcal{P}^{2}\varphi\|_{L^{\infty}(\Omega)}\|\partial_{x}^{2}u\|_{L^{\infty}(\Omega)})
           \|u\|_{\hat{H}_{\psi}^{s}(\Omega)}\|\mathbf{U}_{s}\|_{L_{\varphi}^{2}(\Omega)}\nonumber\\
    &+C_{*}\|\mathcal{P}\varphi\|_{L^{\infty}(\Omega)}(1+\|\mathcal{P}\|_{L^{\infty}(\Omega)})
           \left(\|\partial_{x}^{2}u\|_{L^{\infty}(\Omega)}+\|\partial_{x}^{2}u\|_{L_{x}^{\infty}L_{y}^{2}(\Omega)}\right)
           \|u\|_{\hat{H}_{\psi}^{s}(\Omega)}\|\mathbf{U}_{s}\|_{L_{\varphi}^{2}(\Omega)}\nonumber\\
    &+C_{*}(\|\mathcal{P}\varphi\|_{L^{\infty}(\Omega)}\|\partial_{x}^{2}\partial_{y}^{2}u\|_{L_{x}^{\infty}L_{y}^{2}(\Omega)}
            +\|\mathcal{P}^{2}\varphi\|_{L^{\infty}(\Omega)}\|\partial_{x}^{2}\partial_{y}u\|_{L_{x}^{\infty}L_{y}^{2}(\Omega)})
           \|u\|_{\hat{H}_{\psi}^{s}(\Omega)}\|\mathbf{U}_{s}\|_{L_{\varphi}^{2}(\Omega)}\nonumber\\
    &+C_{*}\|\mathcal{P}\varphi\|_{L^{\infty}(\Omega)}\left(\|\partial_{x}^{2}\partial_{y}u\|_{L_{x}^{\infty}L_{y}^{2}(\Omega)}
                                                     +\|\partial_{x}^{2}\partial_{y}u\|_{L^{\infty}(\Omega)}\right)
           \|u\|_{\hat{H}_{\psi}^{s}(\Omega)}\|\mathbf{U}_{s}\|_{L_{\varphi}^{2}(\Omega)}\nonumber\\
\leq&~C_{*}\left(\|\mathcal{P}\varphi\|_{L^{\infty}(\Omega)}+\|\mathcal{P}^{2}\varphi\|_{L^{\infty}(\Omega)}\right)
           \|u\|_{H^{s-3}(\Omega)}\|u\|_{\hat{H}_{\psi}^{s}(\Omega)}\|\mathbf{U}_{s}\|_{L_{\varphi}^{2}(\Omega)}\\
\leq&~\frac{C_{*}(\check{c}+\sqrt{\bar{\eta}}t)}{(\bar{\varsigma}+t)^{1+\frac{\varepsilon_{0}}{2}}}
      \|u\|_{\hat{H}_{\psi}^{s}(\Omega)}\|\mathbf{U}_{s}\|_{L_{\varphi}^{2}(\Omega)}\\
\leq&~\|u\|_{\hat{H}_{\psi}^{s}(\Omega)}^{2}
     +C_{*}\iota^{-\frac{5}{2}}\|\mathbf{U}_{s}\|_{L_{\varphi}^{2}(\Omega)}^{2}.
\end{align*}

Finally, using the method for estimating the term $\mathfrak{J}_{3}$ in the case $r=s$, one can obtain,
\begin{align*}
\mathfrak{J}_{4}
=&~\sum_{i=3}^{r-2}\binom{r}{i}(\mathcal{P}\partial_{x}^{i}u\partial_{x}^{r-i+1}\partial_{y}u,\varphi_{r}^{2}\mathbf{U}_{r})
  -\sum_{i=3}^{r-2}\binom{r}{i}(\mathcal{P}\rho\partial_{x}^{i}u\partial_{x}^{r-i+1}u,\varphi_{r}^{2}\mathbf{U}_{r})\\
 &-\sum_{i=2}^{r-3}\binom{r}{i}(\mathcal{P}(1-h)\partial_{y}^{-1}[\partial_{x}^{i+1}u]\partial_{x}^{r-i}\partial_{y}^{2}u,\varphi_{r}^{2}\mathbf{U}_{r})
  +\sum_{i=2}^{r-3}\binom{r}{i}(\mathcal{P}(1-h)\rho\partial_{y}^{-1}[\partial_{x}^{i+1}u]\partial_{x}^{r-i}\partial_{y}u,\varphi_{r}^{2}\mathbf{U}_{r})\\
 &+\sum_{i=2}^{r-3}\binom{r}{i}(\mathcal{P}h\partial_{x}^{i+1}u\partial_{x}^{r-i}\partial_{y}u,\varphi_{r}^{2}\mathbf{U}_{r})
  +\sum_{i=2}^{r-3}\binom{r}{i}(\mathcal{P}h'\partial_{y}^{-1}[\partial_{x}^{i+1}u]\partial_{x}^{r-i}\partial_{y}u,\varphi_{r}^{2}\mathbf{U}_{r})\\
\leq&~C_{*}\left(\|\mathcal{P}\varphi_{r}\|_{L^{\infty}(\Omega)}+\|\mathcal{P}^{2}\varphi_{r}\|_{L^{\infty}(\Omega)}\right)
           \|u\|_{H^{s-3}(\Omega)}\|u\|_{\hat{H}_{\psi}^{s}(\Omega)}\|\mathbf{U}_{r}\|_{L_{\varphi_{r}}^{2}(\Omega)}\\
\leq&~\frac{C_{*}(\check{c}+\sqrt{\bar{\eta}}t)}{(\bar{\varsigma}+t)^{1+\frac{\varepsilon_{0}}{2}}}
      \|u\|_{\hat{H}_{\psi}^{s}(\Omega)}\|\mathbf{U}_{r}\|_{L_{\varphi_{r}}^{2}(\Omega)}\\
\leq&~\|u\|_{\hat{H}_{\psi}^{s}(\Omega)}^{2}
     +C_{*}\iota^{-\frac{5}{2}}\|\mathbf{U}_{r}\|_{L_{\varphi_{r}}^{2}(\Omega)}^{2}.
\end{align*}

Collecting all the estimate of $\mathfrak{J}_{1}, \mathfrak{J}_{2}, \mathfrak{J}_{3}$ and $\mathfrak{J}_{4}$ above, together with \eqref{J10} and \eqref{J101'} one can draw a conclusion that when $r=s$,
\begin{align*}
J_{10}
\leq& \frac{1}{4}\|\sqrt{\eta}\partial_{x}\mathbf{U}_{s}\|_{L_{\varphi}^{2}(\Omega)}^{2}
     +C_{*}C_{\tau}\iota^{-\frac{5}{2}}\|\mathbf{U}_{s}\|_{L_{\varphi}^{2}(\Omega)}^{2}
     +C_{*}C_{\tau}(\sqrt{t}+\sqrt{\bar{\eta}}t^{\frac{1}{4}})\|\mathcal{P}\mathbf{U}_{s}\|_{L_{\varphi}^{2}(\Omega)}^{2}\\
    &+5\tau\|\sqrt{\tilde{\eta}}\partial_{x}\tilde{\mathbf{U}}_{s}\|_{L_{\tilde{\varphi}}^{2}(\Omega)}^{2}
     +2\tau\|\tilde{\mathcal{P}}\tilde{\mathbf{U}}_{s}\|_{L_{\tilde{\varphi}}^{2}(\Omega)}^{2}
     +\tau\|\tilde{\mathcal{P}}\mathcal{U}_{s}\|_{L_{\tilde{\varphi}}^{2}(\Omega)}
     +\tau\|\partial_{y}\mathcal{U}_{s}\|_{L_{\tilde{\varphi}}^{2}(\Omega)}^{2}\\
    &+C\|\mathcal{U}_{s}\|_{L_{\tilde{\varphi}}^{2}(\Omega)}^{2}
     +C_{*}\|\tilde{\mathbf{U}}_{s}\|_{L_{\tilde{\varphi}}^{2}(\Omega)}^{2}
     +C\|u\|_{\hat{H}_{\psi}^{s}(\Omega)}^{2},
\end{align*}
and when $r=s+1$,
\begin{align*}
J_{10}
\leq&~\frac{1}{4}\|\sqrt{\eta}\partial_{x}\mathbf{U}_{s+1}\|_{L_{\hat{\varphi}}^{2}(\Omega)}^{2}
     +\frac{1}{8}\|\partial_{y}\mathbf{U}_{s+1}\|_{L_{\hat{\varphi}}^{2}(\Omega)}^{2}
     +C_{*}C_{\tau}(\sqrt{t}+\sqrt{\bar{\eta}}t^{\frac{1}{4}})\|\mathcal{P}\mathbf{U}_{s+1}\|_{L_{\hat{\varphi}}^{2}(\Omega)}^{2}\\
    &+C_{*}C_{\tau}\iota^{-\frac{5}{2}}\|\mathbf{U}_{s+1}\|_{L_{\hat{\varphi}}^{2}(\Omega)}^{2}
     +5\tau\|\sqrt{\tilde{\eta}}\partial_{x}\tilde{\mathbf{U}}_{s+1}\|_{L_{\tilde{\hat{\varphi}}}^{2}(\Omega)}^{2}
     +\tau\|\tilde{\mathcal{P}}\tilde{\mathbf{U}}_{s+1}\|_{L_{\tilde{\hat{\varphi}}}^{2}(\Omega)}^{2}\\
    &+C\|\tilde{\mathbf{U}}_{s+1}\|_{L_{\tilde{\hat{\varphi}}}^{2}(\Omega)}^{2}
     +\tau\|\partial_{y}\mathcal{U}_{s}\|_{L_{\tilde{\varphi}}^{2}(\Omega)}^{2}
     +5\tau\|\sqrt{\eta}\partial_{x}\mathcal{U}_{s}\|_{L_{\tilde{\varphi}}^{2}(\Omega)}^{2}
     +C\|\mathcal{U}_{s}\|_{L_{\tilde{\varphi}}^{2}(\Omega)}^{2}\\
    &+2\tau\|\tilde{\mathcal{P}}\mathcal{U}_{s}\|_{L_{\tilde{\varphi}}^{2}(\Omega)}^{2}
     +C\|u\|_{\hat{H}_{\psi}^{s}(\Omega)}^{2}.
\end{align*}

\noindent{\bf \underline{Estimate of~$J_{11}$.}}
We decompose $J_{11}$ into four parts,
\begin{align}\label{J11}
J_{11}
=\sum_{i=1}^{4}(\mathcal{P}\bar{\mathcal{T}}_{8i},\varphi_{r}^{2}\mathbf{U}_{r}),
\end{align}
where
\begin{align*}
\bar{\mathcal{T}}_{81}
=& \rho\sum_{i=2}^{r}\binom{r}{i}\partial_{x}^{r-i+1}\mathbb{U}\partial_{x}^{i}u_{0}
  +\rho\sum_{i=0}^{r-1}\binom{r}{i}\partial_{x}^{r-i+1}u_{0}\partial_{x}^{i}\mathbb{U}\\
 &-\sum_{i=0}^{r-1}\binom{r}{i}\partial_{y}\partial_{x}^{r-i+1}u_{0}\partial_{x}^{i}\mathbb{U}
  -\sum_{i=2}^{r}\binom{r}{i}\partial_{y}\partial_{x}^{r-i+1}\mathbb{U}\partial_{x}^{i}u_{0}\\
 &+\sum_{i=0}^{r-2}\binom{r}{i}\partial_{y}^{-1}[\partial_{x}^{i+1}\mathbb{U}]\partial_{x}^{r-i}\partial_{y}^{2}u_{0}
  -\rho\sum_{i=0}^{r-2}\binom{r}{i}\partial_{y}^{-1}[\partial_{x}^{i+1}\mathbb{U}]\partial_{x}^{r-i}\partial_{y}u_{0},
\end{align*}
\begin{align*}
\bar{\mathcal{T}}_{82}
= \sum_{i=1}^{r-2}\binom{r}{i}\partial_{y}^{-1}[\partial_{x}^{i+1}u_{0}]\partial_{x}^{r-i}\partial_{y}^{2}\mathbb{U}
  -\rho\sum_{i=1}^{r-2}\binom{r}{i}\partial_{y}^{-1}[\partial_{x}^{i+1}u_{0}]\partial_{x}^{r-i}\partial_{y}\mathbb{U},
\end{align*}
\begin{align*}
\bar{\mathcal{T}}_{83}
=\rho\partial_{x}^{r}(\partial_{t}u_{0}+u_{0}\partial_{x}u_{0}+\partial_{x}P),
\end{align*}
and
\begin{align*}
\bar{\mathcal{T}}_{84}
=&-\sum_{i=0}^{r}\binom{r}{i}\partial_{y}\partial_{x}^{r-i+1}u_{0}\partial_{x}^{i}u_{0}
  +\sum_{i=0}^{r-2}\binom{r}{i}\partial_{y}^{-1}[\partial_{x}^{i+1}u_{0}]\partial_{x}^{r-i}\partial_{y}^{2}u_{0}\\
 &-\rho\sum_{i=0}^{r-2}\binom{r}{i}\partial_{y}^{-1}[\partial_{x}^{i+1}u_{0}]\partial_{x}^{r-i}\partial_{y}u_{0}
  +\partial_{y}^{3}\partial_{x}^{r}u_{0}
  -\rho\partial_{x}^{r}\partial_{y}^{2}u_{0}
  -\partial_{t}\partial_{y}\partial_{x}^{r}u_{0}.
\end{align*}
Recall the definition of $u_{0}$, by using Corollary \ref{ncC2} one has
\begin{align}\label{J111}
(\mathcal{P}\bar{\mathcal{T}}_{81},\varphi_{r}^{2}\mathbf{U}_{r})
\leq&~C_{*}\left(\|\mathbb{U}\|_{\hat{H}_{\psi}^{s}(\Omega)}+\|\partial_{x}^{s}\mathbb{U}\|_{L_{\varphi}^{2}(\Omega)}
              +\|\partial_{y}\partial_{x}^{s}\mathbb{U}\|_{L_{\varphi}^{2}(\Omega)}\right)
        \|\mathbf{U}_{r}\|_{L_{\varphi_{r}}^{2}(\Omega)}\\
\leq&~C_{*}\left(\|\mathbb{U}\|_{\hat{H}_{\psi}^{s}(\Omega)}+\|\mathcal{P}\mathbf{U}_{s}\|_{L_{\varphi}^{2}(\Omega)}
              +\|\mathbf{U}_{s}\|_{L_{\varphi}^{2}(\Omega)}\right)
        \|\mathbf{U}_{r}\|_{L_{\varphi_{r}}^{2}(\Omega)}\nonumber\\
\leq&~C_{*}\|\mathbb{U}\|_{\hat{H}_{\psi}^{s}(\Omega)}^{2}
     +\|\mathcal{P}\mathbf{U}_{s}\|_{L_{\varphi}^{2}(\Omega)}^{2}
     +C_{*}\|\mathbf{U}_{r}\|_{L_{\varphi_{r}}^{2}(\Omega)}^{2}
     +C_{*}\|\mathbf{U}_{s}\|_{L_{\varphi}^{2}(\Omega)}.\nonumber
\end{align}

The term $(\mathcal{P}\bar{\mathcal{T}}_{82},\varphi_{r}^{2}\mathbf{U}_{r})$ is the most difficult term to estimate owing to the nonlocal term. Thanks to the factor $\y^{-\frac{1}{2}}$ in the weight function $\varphi$, if $r=s$,
\begin{align}\label{J112-1}
(\mathcal{P}\bar{\mathcal{T}}_{82},\varphi_{r}^{2}\mathbf{U}_{r})
\leq&~C_{*}\left(\|\mathbb{U}\|_{\hat{H}_{\psi}^{s}(\Omega)}+\|\sqrt{\y}\partial_{x}^{s-1}\partial_{y}^{2}\mathbb{U}\|_{L_{\hat{\psi}}^{2}(\Omega)}\right)
        \|\sqrt{\y}\mathbf{U}_{r}\|_{L_{\varphi_{r}}^{2}(\Omega)}\\
\leq&~C_{*}\left(\|\mathbb{U}\|_{\hat{H}_{\psi}^{s}(\Omega)}+\|\mathbb{U}\|_{\hat{H}_{\psi,\y}^{s}(\Omega)}\right)
        \|\sqrt{\y}\mathbf{U}_{r}\|_{L_{\varphi_{r}}^{2}(\Omega)}\nonumber\\
\leq&~C_{*}\|\mathbb{U}\|_{\hat{H}_{\psi}^{s}(\Omega)}^{2}
     +\tau\|\mathbb{U}\|_{\hat{H}_{\psi,\y}^{s}(\Omega)}^{2}
     +C_{*}C_{\tau}\|\sqrt{\y}\mathbf{U}_{r}\|_{L_{\varphi_{r}}^{2}(\Omega)}^{2}.\nonumber
\end{align}
and if $r=s+1$,
\begin{align}\label{J112-2}
(\mathcal{P}\bar{\mathcal{T}}_{82},\varphi_{r}^{2}\mathbf{U}_{r})
\leq&~C_{*}\left(\|\mathbb{U}\|_{\hat{H}_{\psi}^{s}(\Omega)}+\|\partial_{x}^{s-1}\partial_{y}^{2}\mathbb{U}\|_{L_{\hat{\psi}}^{2}(\Omega)}
                 +\|\partial_{x}^{s}\partial_{y}^{2}\mathbb{U}\|_{L_{\varphi}^{2}(\Omega)}\right)
        \|\sqrt{\y}\mathbf{U}_{r}\|_{L_{\varphi_{r}}^{2}(\Omega)}\\
\leq&~C_{*}\left(\|\mathbb{U}\|_{\hat{H}_{\psi}^{s}(\Omega)}+\|\mathbf{U}_{s}\|_{L_{\varphi}^{2}(\Omega)}
                 +\|\mathcal{P}\mathbf{U}_{s}\|_{L_{\varphi}^{2}(\Omega)}+\|\partial_{y}\mathbf{U}_{s}\|_{L_{\varphi}^{2}(\Omega)}\right)
        \|\sqrt{\y}\mathbf{U}_{r}\|_{L_{\varphi_{r}}^{2}(\Omega)}\nonumber\\
\leq&~C_{*}\|\mathbb{U}\|_{\hat{H}_{\psi}^{s}(\Omega)}^{2}
     +C_{*}\|\mathbf{U}_{s}\|_{L_{\varphi}^{2}(\Omega)}^{2}
     +\|\mathcal{P}\mathbf{U}_{s}\|_{L_{\varphi}^{2}(\Omega)}^{2}\nonumber\\
    &+\frac{1}{16}\|\partial_{y}\mathbf{U}_{s}\|_{L_{\varphi}^{2}(\Omega)}^{2}
     +C_{*}\|\sqrt{\y}\mathbf{U}_{r}\|_{L_{\varphi_{r}}^{2}(\Omega)}^{2}.\nonumber
\end{align}
A combination of \eqref{J112-1} with \eqref{J112-2} gives
\begin{align}\label{J112}
(\mathcal{P}\bar{\mathcal{T}}_{82},\varphi_{r}^{2}\mathbf{U}_{r})
\leq&~C_{*}\|\mathbb{U}\|_{\hat{H}_{\psi}^{s}(\Omega)}^{2}
     +C_{*}\|\mathbf{U}_{s}\|_{L_{\varphi}^{2}(\Omega)}^{2}
     +\tau\|\mathbb{U}\|_{\hat{H}_{\psi,\y}^{s}(\Omega)}^{2}\\
    &+\|\mathcal{P}\mathbf{U}_{s}\|_{L_{\varphi}^{2}(\Omega)}^{2}
     +\frac{1}{16}\|\partial_{y}\mathbf{U}_{s}\|_{L_{\varphi}^{2}(\Omega)}^{2}
     +C_{*}C_{\tau}\|\sqrt{\y}\mathbf{U}_{r}\|_{L_{\varphi_{r}}^{2}(\Omega)}^{2}.\nonumber
\end{align}

From \eqref{B} one can deduce
$$\partial_{t}u_{0}+u_{0}\partial_{x}u_{0}+\partial_{x}P=\partial_{t}(u_{0}-U)+(u_{0}-U)\partial_{x}u_{0}+U\partial_{x}(u_{0}-U),$$
moreover, from the definition of $u_{0}$ one can get
$$\partial_{x}^{r}(\partial_{t}u_{0}+u_{0}\partial_{x}u_{0})=0,~\text{in}~\hat{\Omega}_{t_{3}},$$
thus one has
\begin{align}\label{J113}
(\mathcal{P}\bar{\mathcal{T}}_{83},\varphi_{r}^{2}\mathbf{U}_{r})
\leq&~C_{*}\|\mathcal{P}\partial_{x}^{r+1}P\|_{L_{\varphi_{r}}^{2}(\hat{\Omega})}\|\mathcal{P}\mathbf{U}_{r}\|_{L_{\varphi_{r}}^{2}(\Omega)}\\
    &+C_{*}\|\partial_{x}^{r}(\partial_{t}u_{0}+u_{0}\partial_{x}u_{0}+\partial_{x}P)\|_{L_{\psi}^{2}(\Omega)}\|\mathbf{U}_{r}\|_{L_{\varphi_{r}}^{2}(\Omega)}\nonumber\\
\leq&~\frac{C_{*}}{(\bar{\varsigma}+t)^{\frac{1}{2}+\varepsilon_{0}}}\|\partial_{x}^{r+1}P\|_{L^{2}(\mathbb{T})}^{2}
     +\|\mathcal{P}\mathbf{U}_{r}\|_{L_{\varphi_{r}}^{2}(\Omega)}^{2}
     +\|\mathbf{U}_{r}\|_{L_{\varphi_{r}}^{2}(\Omega)}^{2}
     +C_{*}.\nonumber
\end{align}

Finally, it is easy to see that
\begin{align}\label{J114}
(\mathcal{P}\bar{\mathcal{T}}_{84},\varphi_{r}^{2}\mathbf{U}_{r})
\leq \|\mathbf{U}_{r}\|_{L_{\varphi_{r}}^{2}(\Omega)}^{2}
    +C_{*}.
\end{align}

Combining \eqref{J111}, \eqref{J112}, \eqref{J113} and \eqref{J114} with \eqref{J11}, one has
\begin{align*}
J_{11}
\leq&~C_{*}\|\mathbb{U}\|_{\hat{H}_{\psi}^{s}(\Omega)}^{2}
     +C_{*}\|\mathbf{U}_{s}\|_{L_{\varphi}^{2}(\Omega)}^{2}
     +C_{*}\|\mathbf{U}_{r}\|_{L_{\varphi_{r}}^{2}(\Omega)}^{2}
     +2\|\mathcal{P}\mathbf{U}_{s}\|_{L_{\varphi}^{2}(\Omega)}^{2}
     +\|\mathcal{P}\mathbf{U}_{r}\|_{L_{\varphi_{r}}^{2}(\Omega)}^{2}\nonumber\\
    &+\frac{1}{16}\|\partial_{y}\mathbf{U}_{s}\|_{L_{\varphi}^{2}(\Omega)}^{2}
     +C_{*}C_{\tau}\|\sqrt{\y}\mathbf{U}_{r}\|_{L_{\varphi_{r}}^{2}(\Omega)}^{2}
     +\frac{C_{*}}{(\bar{\varsigma}+t)^{\frac{1}{2}+\varepsilon_{0}}}\|\partial_{x}^{r+1}P\|_{L^{2}(\mathbb{T})}^{2}
     +C_{*}.
\end{align*}

Now we are in a position to establish the estimate of $\mathbf{U}_{s}$ and $\mathbf{U}_{s+1}$. Combining the estimates of $J_{0}$,$J_{2}$,...,$J_{11}$ above with \eqref{UnE}, for any $t\in[0,t_{3}]$ we have
\begin{align*}
&\frac{1}{2}\frac{d}{dt}\|\mathbf{U}_{s}\|_{L_{\varphi}^{2}(\Omega)}^{2}
+\frac{3}{4}\Lambda\delta\|\sqrt{\y}\mathbf{U}_{s}\|_{L_{\varphi}^{2}(\Omega)}^{2}
+\frac{\lambda}{64}\|\sqrt{\omega_{\lambda}}\mathbf{U}_{s}\|_{L_{\varphi}^{2}(\Omega)}^{2}
+\|\partial_{y}\mathbf{U}_{s}\|_{L_{\varphi}^{2}(\Omega)}^{2}
+\|\sqrt{\eta}\partial_{x}\mathbf{U}_{s}\|_{L_{\varphi}^{2}(\Omega)}^{2}\\
\leq&~C_{*}C_{\tau}\iota^{-\frac{5}{2}}\|\mathbf{U}_{s}\|_{L_{\varphi}^{2}(\Omega)}^{2}
     +\frac{3}{8}\|\partial_{y}\mathbf{U}_{s}\|_{L_{\varphi}^{2}(\Omega)}^{2}
     +\left[C_{*}C_{\tau}(\sqrt{t}+\sqrt{\bar{\eta}}t^{\frac{1}{4}})+C\right]\|\mathcal{P}\mathbf{U}_{s}\|_{L_{\varphi}^{2}(\Omega)}^{2}\\
    &+\left(\frac{3}{8}+C_{\lambda}\theta\sqrt{t}\right)\|\sqrt{\eta}\partial_{x}\mathbf{U}_{s}\|_{L_{\varphi}^{2}(\Omega)}^{2}
     +\frac{1}{4}\|\sqrt{\omega_{\lambda}}\mathbf{U}_{s}\|_{L_{\varphi}^{2}(\Omega)}^{2}
     +(7\tau+C_{*}\sqrt{\iota})\|\sqrt{\tilde{\eta}}\partial_{x}\tilde{\mathbf{U}}_{s}\|_{L_{\tilde{\varphi}}^{2}(\Omega)}^{2}\\
    &+3\tau\|\tilde{\mathcal{P}}\tilde{\mathbf{U}}_{s}\|_{L_{\tilde{\varphi}}^{2}(\Omega)}^{2}
     +C_{*}C_{\tau}\|\tilde{\mathbf{U}}_{s}\|_{L_{\tilde{\varphi}}^{2}(\Omega)}^{2}
     +\tau\|\partial_{y}\mathcal{U}_{s}\|_{L_{\tilde{\varphi}}^{2}(\Omega)}^{2}
     +\tau\|\tilde{\mathcal{P}}\mathcal{U}_{s}\|_{L_{\tilde{\varphi}}^{2}(\Omega)}^{2}
     +C_{*}C_{\tau}\|\sqrt{\y}\mathbf{U}_{s}\|_{L_{\varphi}^{2}(\Omega)}^{2}\\
    &+\tau\|\sqrt{\y}\mathbb{U}\|_{\hat{H}_{\psi}^{s}(\Omega)}^{2}
     +C\|u\|_{\mathcal{H}_{\psi,\tilde{\varphi}}^{s}(\Omega)}^{2}
     +\frac{C_{*}}{(\bar{\varsigma}+t)^{\frac{1}{2}+\varepsilon_{0}}}\|\partial_{x}^{s+1}P\|_{L^{2}(\mathbb{T})}^{2}
     +C_{*}\left(1+\frac{1}{(\bar{\varsigma}+t)^{\frac{1}{2}+\varepsilon_{0}}}\right),
\end{align*}
and
\begin{align*}
&\frac{1}{2}\frac{d}{dt}\|\mathbf{U}_{s+1}\|_{L_{\hat{\varphi}}^{2}(\Omega)}^{2}
+\frac{3}{4}\Lambda\delta\|\sqrt{\y}\mathbf{U}_{s+1}\|_{L_{\hat{\varphi}}^{2}(\Omega)}^{2}
+\frac{\lambda}{64}\|\sqrt{\omega_{\lambda}}\mathbf{U}_{s+1}\|_{L_{\hat{\varphi}}^{2}(\Omega)}^{2}
+\|\partial_{y}\mathbf{U}_{s+1}\|_{L_{\hat{\varphi}}^{2}(\Omega)}^{2}
+\|\sqrt{\eta}\partial_{x}\mathbf{U}_{s+1}\|_{L_{\hat{\varphi}}^{2}(\Omega)}^{2}\\
\leq&~\frac{7}{16}\|\partial_{y}\mathbf{U}_{s+1}\|_{L_{\hat{\varphi}}^{2}(\Omega)}^{2}
     +\left(\frac{3}{8}+C_{\lambda}\theta\sqrt{t}\right)\|\sqrt{\eta}\partial_{x}\mathbf{U}_{s+1}\|_{L_{\hat{\varphi}}^{2}(\Omega)}^{2}
     +\frac{1}{4}\|\sqrt{\omega_{\lambda}}\mathbf{U}_{s+1}\|_{L_{\hat{\varphi}}^{2}(\Omega)}^{2}
     +C_{*}C_{\tau}\iota^{-\frac{5}{2}}\|\mathbf{U}_{s+1}\|_{L_{\hat{\varphi}}^{2}(\Omega)}^{2}\\
    &+\left[C_{*}C_{\tau}(\sqrt{t}+\sqrt{\bar{\eta}}t^{\frac{1}{4}})+C\right]\|\mathcal{P}\mathbf{U}_{s+1}\|_{L_{\hat{\varphi}}^{2}(\Omega)}^{2}
     +C_{*}C_{\tau}\|\sqrt{\y}\mathbf{U}_{s+1}\|_{L_{\hat{\varphi}}^{2}(\Omega)}^{2}
     +C_{\lambda}\|\mathcal{P}\mathbf{U}_{s}\|_{L_{\varphi}^{2}(\Omega)}^{2}
     +C_{\lambda}\|\partial_{y}\mathbf{U}_{s}\|_{L_{\varphi}^{2}(\Omega)}^{2}\\
    &+C_{*}\|\mathbf{U}_{s}\|_{L_{\varphi}^{2}(\Omega)}^{2}
     +(7\tau+C_{*}\sqrt{\iota})\|\sqrt{\tilde{\eta}}\partial_{x}\tilde{\mathbf{U}}_{s+1}\|_{L_{\tilde{\hat{\varphi}}}^{2}(\Omega)}^{2}
     +3\tau\|\tilde{\mathcal{P}}\tilde{\mathbf{U}}_{s+1}\|_{L_{\tilde{\hat{\varphi}}}^{2}(\Omega)}^{2}
     +C_{*}C_{\tau}\|\tilde{\mathbf{U}}_{s+1}\|_{L_{\tilde{\hat{\varphi}}}^{2}(\Omega)}^{2}\\
    &+\tau\|\partial_{y}\mathcal{U}_{s}\|_{L_{\tilde{\varphi}}^{2}(\Omega)}^{2}
     +5\tau\|\sqrt{\eta}\partial_{x}\mathcal{U}_{s}\|_{L_{\tilde{\varphi}}^{2}(\Omega)}^{2}
     +2\tau\|\tilde{\mathcal{P}}\mathcal{U}_{s}\|_{L_{\tilde{\varphi}}^{2}(\Omega)}^{2}
     +C\|u\|_{\mathcal{H}_{\psi,\tilde{\varphi}}^{s}(\Omega)}^{2}
     +\tau\|\sqrt{\y}\mathbb{U}\|_{\hat{H}_{\psi}^{s}(\Omega)}^{2}\\
    &+\frac{C_{*}}{(\bar{\varsigma}+t)^{\frac{1}{2}+\varepsilon_{0}}}\|\partial_{x}^{s+2}P\|_{L_{x}^{2}(\mathbb{T})}^{2}
     +C_{*}\left(1+\frac{1}{(\bar{\varsigma}+t)^{\frac{1}{2}+\varepsilon_{0}}}\right).
\end{align*}
By taking $\lambda$ sufficiently large in the above inequalities, we obtain the following proposition.
\begin{prop}
Let $\lambda>\lambda_{*}, \Lambda>\Lambda_{*}, \theta>\theta_{*}, \iota\in(0,\iota_{*}), \mathcal{Z}>\mathcal{Z}_{*}$, where $\lambda_{*}, \Lambda_{*},\iota_{*}, \mathcal{Z}_{*}$ are given in Theorem \ref{wpae}, and $\theta_{*}$ is given in Lemma \ref{wdotu}. Under the assumption of Theorem \ref{Phe}, if $\lambda$ is large enough, then the following inequality for any $t\in[0,t_{3}]$ and $\bar{\tau}\in(0,1)$,
\begin{align}\label{phe}
&\frac{1}{2}\frac{d}{dt}\left(\|\mathbf{U}_{s}\|_{L_{\varphi}^{2}(\Omega)}^{2}+\bar{\tau}\|\mathbf{U}_{s+1}\|_{L_{\hat{\varphi}}^{2}(\Omega)}^{2}\right)
 +\frac{3}{4}\Lambda\delta\left(\|\sqrt{\y}\mathbf{U}_{s}\|_{L_{\varphi}^{2}(\Omega)}^{2}+\bar{\tau}\|\sqrt{\y}\mathbf{U}_{s+1}\|_{L_{\hat{\varphi}}^{2}(\Omega)}^{2}\right)\\
&+\frac{\lambda}{128}\|\sqrt{\omega_{\lambda}}\mathbf{U}_{s}\|_{L_{\varphi}^{2}(\Omega)}^{2}
 +\frac{\lambda\bar{\tau}}{128}\|\sqrt{\omega_{\lambda}}\mathbf{U}_{s+1}\|_{L_{\hat{\varphi}}^{2}(\Omega)}^{2}
 +\frac{1}{2}\|\partial_{y}\mathbf{U}_{s}\|_{L_{\varphi}^{2}(\Omega)}^{2}
 +\frac{\bar{\tau}}{2}\|\partial_{y}\mathbf{U}_{s+1}\|_{L_{\hat{\varphi}}^{2}(\Omega)}^{2}\nonumber\\
&+\frac{5}{8}\|\sqrt{\eta}\partial_{x}\mathbf{U}_{s}\|_{L_{\varphi}^{2}(\Omega)}^{2}
 +\frac{5\bar{\tau}}{8}\|\sqrt{\eta}\partial_{x}\mathbf{U}_{s+1}\|_{L_{\hat{\varphi}}^{2}(\Omega)}^{2}\nonumber\\
\leq&~C_{*}C_{\tau}\iota^{-\frac{5}{2}}
      \left(\|\mathbf{U}_{s}\|_{L_{\varphi}^{2}(\Omega)}^{2}+\bar{\tau}\|\mathbf{U}_{s+1}\|_{L_{\hat{\varphi}}^{2}(\Omega)}^{2}\right)
      +C_{*}C_{\tau}\left(\|\sqrt{\y}\mathbf{U}_{s}\|_{L_{\varphi}^{2}(\Omega)}^{2}+\bar{\tau}\|\sqrt{\y}\mathbf{U}_{s+1}\|_{L_{\hat{\varphi}}^{2}(\Omega)}^{2}\right)\nonumber\\
    &+\left[C_{*}C_{\tau}(\sqrt{t}+\sqrt{\bar{\eta}}t^{\frac{1}{4}})+\check{C}\right]
      \left(\|\mathcal{P}\mathbf{U}_{s}\|_{L_{\varphi}^{2}(\Omega)}^{2}+\bar{\tau}\|\mathcal{P}\mathbf{U}_{s+1}\|_{L_{\hat{\varphi}}^{2}(\Omega)}^{2}\right)
     +\check{\mathcal{C}}_{\lambda}\bar{\tau}\|\mathcal{P}\mathbf{U}_{s}\|_{L_{\varphi}^{2}(\Omega)}^{2}\nonumber\\
    &+\check{\mathcal{C}}_{\lambda}\bar{\tau}\|\partial_{y}\mathbf{U}_{s}\|_{L_{\varphi}^{2}(\Omega)}^{2}
     +(7\tau+C_{*}\sqrt{\iota})\left(\|\sqrt{\tilde{\eta}}\partial_{x}\tilde{\mathbf{U}}_{s}\|_{L_{\tilde{\varphi}}^{2}(\Omega)}^{2}
                                     +\bar{\tau}\|\sqrt{\tilde{\eta}}\partial_{x}\tilde{\mathbf{U}}_{s+1}\|_{L_{\tilde{\hat{\varphi}}}^{2}(\Omega)}^{2}\right)\nonumber\\
    &+C_{\lambda}\theta\sqrt{t}\left(\|\sqrt{\eta}\partial_{x}\mathbf{U}_{s}\|_{L_{\varphi}^{2}(\Omega)}^{2}
                             +\bar{\tau}\|\sqrt{\eta}\partial_{x}\mathbf{U}_{s+1}\|_{L_{\hat{\varphi}}^{2}(\Omega)}^{2}\right)\nonumber\\
    &+3\tau\|\tilde{\mathcal{P}}\tilde{\mathbf{U}}_{s}\|_{L_{\tilde{\varphi}}^{2}(\Omega)}^{2}
     +3\tau\bar{\tau}\|\tilde{\mathcal{P}}\tilde{\mathbf{U}}_{s+1}\|_{L_{\tilde{\hat{\varphi}}}^{2}(\Omega)}^{2}
     +C_{*}C_{\tau}\|\tilde{\mathbf{U}}_{s}\|_{L_{\tilde{\varphi}}^{2}(\Omega)}^{2}
     +C_{*}C_{\tau}\bar{\tau}\|\tilde{\mathbf{U}}_{s+1}\|_{L_{\tilde{\hat{\varphi}}}^{2}(\Omega)}^{2}\nonumber\\
    &+2\tau\|\sqrt{\y}\mathbb{U}\|_{\hat{H}_{\psi}^{s}(\Omega)}^{2}
     +C\|u\|_{\mathcal{H}_{\psi,\tilde{\varphi}}^{s}(\Omega)}^{2}
     +\frac{C_{*}}{(\bar{\varsigma}+t)^{\frac{1}{2}+\varepsilon_{0}}}\|\partial_{x}^{s+1}P\|_{L_{x}^{2}(\mathbb{T})}^{2}
     +\frac{C_{*}}{(\bar{\varsigma}+t)^{\frac{1}{2}+\varepsilon_{0}}}\|\partial_{x}^{s+2}P\|_{L_{x}^{2}(\mathbb{T})}^{2}\nonumber\\
    &+2\tau\|\partial_{y}\mathcal{U}_{s}\|_{L_{\tilde{\varphi}}^{2}(\Omega)}^{2}
     +5\tau\|\sqrt{\eta}\partial_{x}\mathcal{U}_{s}\|_{L_{\tilde{\varphi}}^{2}(\Omega)}^{2}
     +3\tau\|\tilde{\mathcal{P}}\mathcal{U}_{s}\|_{L_{\tilde{\varphi}}^{2}(\Omega)}^{2}
     +C_{*}\left(1+\frac{1}{(\bar{\varsigma}+t)^{\frac{1}{2}+\varepsilon_{0}}}\right),\nonumber
\end{align}
where $\check{\mathcal{C}}_{\lambda}$ is a positive constant depending only on $\lambda$, and $\check{C}$ is a positive constant independent of $\lambda$ and $\mathcal{Z}$.
\end{prop}

To prove Theorem \ref{Phe}, let $\theta>\theta^{*}$ be fixed, and write
$\mathbb{S}(t)=\|\mathbf{U}_{s}\|_{L_{\varphi}^{2}(\Omega)}^{2}+\bar{\tau}\|\mathbf{U}_{s+1}\|_{L_{\hat{\varphi}}^{2}(\Omega)}^{2}$,
\begin{align*}
\mathbb{W}(t)
=&\Lambda\delta\|\sqrt{\y}\mathbf{U}_{s}\|_{L_{\varphi}^{2}(\Omega)}^{2}
 +\bar{\tau}\Lambda\delta\|\sqrt{\y}\mathbf{U}_{s+1}\|_{L_{\hat{\varphi}}^{2}(\Omega)}^{2}
 +\frac{\lambda}{128}\|\sqrt{\omega_{\lambda}}\mathbf{U}_{s}\|_{L_{\varphi}^{2}(\Omega)}^{2}
 +\frac{\lambda\bar{\tau}}{128}\|\sqrt{\omega_{\lambda}}\mathbf{U}_{s+1}\|_{L_{\hat{\varphi}}^{2}(\Omega)}^{2}\\
&+\|\partial_{y}\mathbf{U}_{s}\|_{L_{\varphi}^{2}(\Omega)}^{2}
 +\bar{\tau}\|\partial_{y}\mathbf{U}_{s+1}\|_{L_{\hat{\varphi}}^{2}(\Omega)}^{2}
 +\|\sqrt{\eta}\partial_{x}\mathbf{U}_{s}\|_{L_{\varphi}^{2}(\Omega)}^{2}
 +\bar{\tau}\|\sqrt{\eta}\partial_{x}\mathbf{U}_{s+1}\|_{L_{\hat{\varphi}}^{2}(\Omega)}^{2},
\end{align*}
and
\begin{align*}
\mathbb{G}(t)
=&~(7\tau+C_{*}\sqrt{\iota})\left(\|\sqrt{\tilde{\eta}}\partial_{x}\tilde{\mathbf{U}}_{s}\|_{L_{\tilde{\varphi}}^{2}(\Omega)}^{2}
                                     +\bar{\tau}\|\sqrt{\tilde{\eta}}\partial_{x}\tilde{\mathbf{U}}_{s+1}\|_{L_{\tilde{\hat{\varphi}}}^{2}(\Omega)}^{2}\right)\\
 &+3\tau\|\tilde{\mathcal{P}}\tilde{\mathbf{U}}_{s}\|_{L_{\tilde{\varphi}}^{2}(\Omega)}^{2}
  +3\tau\bar{\tau}\|\tilde{\mathcal{P}}\tilde{\mathbf{U}}_{s+1}\|_{L_{\tilde{\hat{\varphi}}}^{2}(\Omega)}^{2}
  +C_{*}C_{\tau}\|\tilde{\mathbf{U}}_{s}\|_{L_{\tilde{\varphi}}^{2}(\Omega)}^{2}
  +C_{*}C_{\tau}\bar{\tau}\|\tilde{\mathbf{U}}_{s+1}\|_{L_{\tilde{\hat{\varphi}}}^{2}(\Omega)}^{2}\nonumber\\
 &+2\tau\|\sqrt{\y}\mathbb{U}\|_{\hat{H}_{\psi}^{s}(\Omega)}^{2}
  +C\|u\|_{\mathcal{H}_{\psi,\tilde{\varphi}}^{s}(\Omega)}^{2}
  +\frac{C_{*}}{(\bar{\varsigma}+t)^{\frac{1}{2}+2\varepsilon_{0}}}\|\partial_{x}^{s+1}P\|_{L_{x}^{2}(\mathbb{T})}^{2}
  +\frac{C_{*}}{(\bar{\varsigma}+t)^{\frac{1}{2}+2\varepsilon_{0}}}\|\partial_{x}^{s+2}P\|_{L_{x}^{2}(\mathbb{T})}^{2}\nonumber\\
 &+2\tau\|\partial_{y}\mathcal{U}_{s}\|_{L_{\tilde{\varphi}}^{2}(\Omega)}^{2}
  +5\tau\|\sqrt{\eta}\partial_{x}\mathcal{U}_{s}\|_{L_{\tilde{\varphi}}^{2}(\Omega)}^{2}
  +3\tau\|\tilde{\mathcal{P}}\mathcal{U}_{s}\|_{L_{\tilde{\varphi}}^{2}(\Omega)}^{2}
  +C_{*}\left(1+\frac{1}{(\bar{\varsigma}+t)^{\frac{1}{2}+\varepsilon_{0}}}\right).\nonumber
\end{align*}

By taking $\lambda$ large enough, one has
\begin{align*}
(\check{C}+2)\|\mathcal{P}\mathbf{U}_{s}\|_{L_{\varphi}^{2}(\Omega)}^{2}
\leq \frac{\lambda}{256}\|\sqrt{\omega_{\lambda}}\mathbf{U}_{s}\|_{L_{\varphi}^{2}(\Omega)}^{2}
    +C\|\mathbf{U}_{s}\|_{L_{\varphi}^{2}(\Omega)}^{2},~\forall~t\in[0,t_{3}],
\end{align*}
and
\begin{align}\label{pus}
\|\mathcal{P}\mathbf{U}_{s}\|_{L_{\varphi}^{2}(\Omega)}^{2}
\leq \frac{\lambda}{256}\|\sqrt{\omega_{\lambda}}\mathbf{U}_{s}\|_{L_{\varphi}^{2}(\Omega)}^{2}
    +\|\mathbf{U}_{s}\|_{L_{\varphi}^{2}(\Omega)}^{2},   ~\forall~t\in[0,t_{3}],
\end{align}
moreover,
\begin{align*}
(\check{C}+2)\|\mathcal{P}\mathbf{U}_{s+1}\|_{L_{\hat{\varphi}}^{2}(\Omega)}^{2}
\leq \frac{\lambda}{256}\|\sqrt{\omega_{\lambda}}\mathbf{U}_{s+1}\|_{L_{\hat{\varphi}}^{2}(\Omega)}^{2}
    +C\|\mathbf{U}_{s+1}\|_{L_{\hat{\varphi}}^{2}(\Omega)}^{2},~\forall~t\in[0,t_{3}].
\end{align*}
and
\begin{align}\label{pus+1}
\|\mathcal{P}\mathbf{U}_{s+1}\|_{L_{\hat{\varphi}}^{2}(\Omega)}^{2}
\leq \frac{\lambda}{256}\|\sqrt{\omega_{\lambda}}\mathbf{U}_{s+1}\|_{L_{\hat{\varphi}}^{2}(\Omega)}^{2}
    +\|\mathbf{U}_{s+1}\|_{L_{\hat{\varphi}}^{2}(\Omega)}^{2},~\forall~t\in[0,t_{3}].
\end{align}
With $\lambda$ fixed, we take $\bar{\tau}$ small enough such that $\check{\mathcal{C}}_{\lambda}\bar{\tau}\leq1$, then under the assumption of Theorem \ref{Phe}, since $|\rho(0,x,y)|\leq C|\rho_{in}(x,y)|\leq C\mathcal{P}_{in}$ in $\Omega$, with such fixed $\lambda$ and $\bar{\tau}$ we can take $\mathcal{Z}$ large enough to obtain
\begin{align}\label{ine}
\mathbb{S}(0)
\leq&~C_{\lambda}\|\mathbf{U}_{s}(0,\cdot)\|_{L_{\hat{\psi}_{in}}^{2}(\Omega)}^{2}
     +C_{\lambda}\bar{\tau}\|\mathbf{U}_{s+1}(0,\cdot)\|_{L_{\check{\psi}_{in}}^{2}(\Omega)}^{2}\\
\leq&~C_{\lambda}\|\mathcal{P}_{in}\partial_{y}\partial_{x}^{s}\tilde{u}_{in}\|_{L_{\hat{\psi}_{in}}^{2}(\Omega)}^{2}
     +C_{\lambda}\|\mathcal{P}_{in}\rho(0,x,y)\partial_{x}^{s}\tilde{u}_{in}\|_{L_{\hat{\psi}_{in}}^{2}(\Omega)}^{2}\nonumber\\
    &+C_{\lambda}\bar{\tau}\|\mathcal{P}_{in}\partial_{y}\partial_{x}^{s+1}\tilde{u}_{in}\|_{L_{\check{\psi}_{in}}^{2}(\Omega)}^{2}
     +C_{\lambda}\bar{\tau}\|\mathcal{P}_{in}\rho(0,x,y)\partial_{x}^{s+1}\tilde{u}_{in}\|_{L_{\check{\psi}_{in}}^{2}(\Omega)}^{2}\nonumber\\
\leq&~C_{\lambda}\|\mathcal{P}_{in}\partial_{y}\partial_{x}^{s}\tilde{u}_{in}\|_{L_{\hat{\psi}_{in}}^{2}(\Omega)}^{2}
     +C_{\lambda}\|\mathcal{P}_{in}^{2}\partial_{x}^{s}\tilde{u}_{in}\|_{L_{\hat{\psi}_{in}}^{2}(\Omega)}^{2}\nonumber\\
    &+C_{\lambda}\bar{\tau}\|\mathcal{P}_{in}\partial_{y}\partial_{x}^{s+1}\tilde{u}_{in}\|_{L_{\check{\psi}_{in}}^{2}(\Omega)}^{2}
     +C_{\lambda}\bar{\tau}\|\mathcal{P}_{in}^{2}\partial_{x}^{s+1}\tilde{u}_{in}\|_{L_{\check{\psi}_{in}}^{2}(\Omega)}^{2}\nonumber\\
\leq&~\frac{\bar{\tau}\mathcal{Z}^{2}}{64}.\nonumber
\end{align}
Moreover, we can take $\tilde{t}\in(0,t_{3}]$ and $\tau, \iota$ small enough such that
\begin{align}\label{G}
\int_{0}^{\tilde{t}}\mathbb{G}(\sigma)d\sigma\leq \frac{\bar{\tau}\mathcal{Z}^{2}}{128}.
\end{align}
For such fixed $\lambda,\tau$ and $\mathcal{Z}$, we take $\bar{t}_{*}\in[0,\tilde{t}]$ small enough to obtain
$$C_{*}C_{\tau}(\sqrt{t}+\sqrt{\bar{\eta}}t^{\frac{1}{4}})\leq 1,~\forall~t\in[0,\bar{t}_{*}],$$
With $\lambda,\tau,\bar{\tau},\theta,\iota$ and $\mathcal{Z}$ fixed, by taking $\bar{t}^{*}\in [0,\bar{t}_{*}]$ small enough such that $C_{\lambda}\theta\sqrt{t}\leq \frac{1}{8}$ for any $t\in[0,\bar{t}^{*}]$ and  $\Lambda$ large enough such that $C_{*}C_{\tau}\leq \frac{\Lambda\delta}{4}$ in \eqref{phe}, one can deduce that for any $t\in[0,\bar{t}^{*}]$,
\begin{align}\label{phe0}
&\frac{1}{2}\frac{d}{dt}\left(\|\mathbf{U}_{s}\|_{L_{\varphi}^{2}(\Omega)}^{2}+\bar{\tau}\|\mathbf{U}_{s+1}\|_{L_{\hat{\varphi}}^{2}(\Omega)}^{2}\right)
 +\frac{\Lambda\delta}{2}\left(\|\sqrt{\y}\mathbf{U}_{s}\|_{L_{\varphi}^{2}(\Omega)}^{2}+\bar{\tau}\|\sqrt{\y}\mathbf{U}_{s+1}\|_{L_{\hat{\varphi}}^{2}(\Omega)}^{2}\right)\\
&+\frac{\lambda}{256}\|\sqrt{\omega_{\lambda}}\mathbf{U}_{s}\|_{L_{\varphi}^{2}(\Omega)}^{2}
 +\frac{\lambda\bar{\tau}}{256}\|\sqrt{\omega_{\lambda}}\mathbf{U}_{s+1}\|_{L_{\hat{\varphi}}^{2}(\Omega)}^{2}
 +\frac{1}{2}\|\partial_{y}\mathbf{U}_{s}\|_{L_{\varphi}^{2}(\Omega)}^{2}
 +\frac{\bar{\tau}}{2}\|\partial_{y}\mathbf{U}_{s+1}\|_{L_{\hat{\varphi}}^{2}(\Omega)}^{2}\nonumber\\
&+\frac{1}{2}\|\sqrt{\eta}\partial_{x}\mathbf{U}_{s}\|_{L_{\varphi}^{2}(\Omega)}^{2}
 +\frac{\bar{\tau}}{2}\|\sqrt{\eta}\partial_{x}\mathbf{U}_{s+1}\|_{L_{\hat{\varphi}}^{2}(\Omega)}^{2}\nonumber\\
\leq&~C_{*}C_{\tau}\iota^{-\frac{5}{2}}
      \left(\|\mathbf{U}_{s}\|_{L_{\varphi}^{2}(\Omega)}^{2}+\bar{\tau}\|\mathbf{U}_{s+1}\|_{L_{\hat{\varphi}}^{2}(\Omega)}^{2}\right)
    +\mathbb{G}(t),\nonumber
\end{align}
which gives
\begin{align}\label{phe1}
\frac{d}{dt}\mathbb{S}(t)+\mathbb{W}(t)
\leq C_{*}C_{\tau}\iota^{-\frac{5}{2}}\mathbb{S}(t)+2\mathbb{G}(t).
\end{align}
Utilizing the Gronwall inequality, from \eqref{phe1} one can deduce that for any $t\in[0,\bar{t}^{*}]$,
\begin{align}\label{phe3}
\mathbb{S}(t)+\int_{0}^{t}\mathbb{W}(\sigma)d\sigma
\leq& C_{*}C_{\tau}\iota^{-\frac{5}{2}}\int_{0}^{t}\mathbb{S}(\sigma)d\sigma
     +\int_{0}^{t}2\mathbb{G}(\sigma)d\sigma
     +\mathbb{S}(0)\\
\leq& C_{*}C_{\tau}\iota^{-\frac{5}{2}}\int_{0}^{t}e^{C_{*}C_{\tau}\iota^{-\frac{5}{2}}\sigma}d\sigma
           \left(\int_{0}^{t}2\mathbb{G}(\sigma)d\sigma+\mathbb{S}(0)\right)\nonumber\\
    &+\int_{0}^{t}2\mathbb{G}(\sigma)d\sigma
     +\mathbb{S}(0).\nonumber
\end{align}
With \eqref{G} and \eqref{ine}, from \eqref{phe3} one can choose $\hat{t}\in (0,\bar{t}^{*})$ small enough such that for any $t\in[0,\hat{t}]$, there holds
\begin{align}\label{phe4}
\mathbb{S}(t)+\int_{0}^{t}\mathbb{W}(\sigma)d\sigma\leq \frac{\bar{\tau}\mathcal{Z}^{2}}{16},
\end{align}
Moreover, from \eqref{pus} and \eqref{pus+1} one can choose $t_{*}\in(0,\hat{t}]$ such that
\begin{align}\label{Pus}
\|\mathcal{P}\mathbf{U}_{s}\|_{L_{t}^{2}L_{\varphi}^{2}(\Omega_{t_{*}})}^{2}
\leq \frac{\lambda}{256}\|\sqrt{\omega_{\lambda}}\mathbf{U}_{s}\|_{L_{t}^{2}L_{\varphi}^{2}(\Omega_{t_{*}})}^{2}
    +\|\mathbf{U}_{s}\|_{L_{t}^{\infty}L_{\varphi}^{2}(\Omega_{t_{*}})}^{2},
\end{align}
and
\begin{align}\label{Pus+1}
\|\mathcal{P}\mathbf{U}_{s+1}\|_{L_{t}^{2}L_{\hat{\varphi}}^{2}(\Omega_{t_{*}})}^{2}
\leq \frac{\lambda}{256}\|\sqrt{\omega_{\lambda}}\mathbf{U}_{s+1}\|_{L_{t}^{2}L_{\hat{\varphi}}^{2}(\Omega_{t_{*}})}^{2}
    +\|\mathbf{U}_{s+1}\|_{L_{t}^{\infty}L_{\hat{\varphi}}^{2}(\Omega_{t_{*}})}^{2}.
\end{align}
Recall the definition of $\mathbb{S}$ and $\mathbb{W}$, from \eqref{Pus},\eqref{Pus+1} and \eqref{phe4} we get \eqref{Hoe1} and \eqref{Hoe2}.
\end{proof}
~~~~~~~~~~

\subsection{Estimates of source term}~

In first two subsections of Section 5 we establish the a posteriori estimates for the approximate solution under some certain conditions of source term $\hat{F}$. In this section, we shall verify that these assumptions are satisfied under the assumption that $u$ is the solution to the problem \eqref{app'} obtained in Theorem \ref{wpae}. Moreover, in order to ensure that the next order corrected increment solution can be obtained by Theorem \ref{wpae}, we need to verify that $\hat{F}$ satisfies all the requirements of Theorem \ref{wpae}, under the a posteriori estimates of the approximate solution.

The main result of this subsection is as follows.
\begin{theorem}
Assume $u$ is the solution to the problem \eqref{app'} on the time interval $[0,t_{*}]$ obtained in Theorem \ref{wpae}, let $\hat{F}$ be given in \eqref{hatF}, then
\begin{align}\label{Fe1}
|D^{\gamma}\hat{F}|\leq C_{*}e^{-\frac{4}{3}\delta\y}~\text{in}~\Omega_{t_{*}}, \forall \gamma\in\Gamma_{9},
\end{align}
\begin{align}\label{Fe2}
|\partial_{y}\hat{F}|\leq C_{*}\iota^{-1}(\dot{\mathfrak{U}}+\bar{\varsigma})~\text{in}~\Omega_{t_{*}},
\end{align}
\begin{align}\label{Fe3}
|\partial_{y}\hat{F}|\leq \bar{\iota}+C|y-y_{*}|+C_{*}t~\text{in}~\Omega_{t_{*}}^{*},
\end{align}
\begin{align}\label{Fe4}
\|\mathcal{P}(\partial_{y}^{2}\hat{F}-\rho\partial_{y}\hat{F})\|_{L^{\infty}(\Omega_{t_{*}})}\leq C_{*}\iota^{-1},
\end{align}
\begin{align}\label{Fe5}
\|\partial_{t}\hat{F}\|_{L_{t}^{\infty}H_{\psi}^{s-4}(\Omega_{t_{*}})}
+\|\hat{F}\|_{L_{t}^{\infty}H_{\psi}^{s-2}(\Omega_{t_{*}})}
\leq C_{*}\bar{\eta}.
\end{align}
Moreover, assume that $\mathfrak{U}$ is the solution to the problem \eqref{App+1} on $[0,t_{*}]$ such that \eqref{Hoe1} and \eqref{Hoe2} holds, then
\begin{align}\label{Fe6}
& \|\hat{\mathbf{F}}_{s}\|_{L_{t}^{2}L_{\hat{\psi}}^{2}(\Omega_{t_{*}})}^{2}
 +\|\hat{F}\|_{L_{t}^{2}\hat{H}_{\psi}^{s}(\Omega_{t_{*}})}^{2}
 +\sum_{i=0}^{\frac{s-1}{2}}\|\partial_{t}^{i}\hat{F}\|_{L_{t}^{\infty}\mathring{H}^{s-1-2i}([0,t_{*}]\times\mathbb{T})}^{2}
 +\|\partial_{x}^{s}\hat{F}\|_{L_{t}^{\infty}\mathring{H}^{0}([0,t_{*}]\times\mathbb{T})}^{2}
\leq C_{*}\bar{\eta},
\end{align}
where $\hat{\mathbf{F}}_{s}=\mathcal{P}\left(\partial_{y}\partial_{x}^{s}\hat{F}-\rho\partial_{x}^{s}\hat{F}\right)$.
\end{theorem}
\begin{proof}
Since $\Lambda t_{*}\leq \frac{1}{12}$, the inequality \eqref{Fe1} is a direct corollary of \eqref{priore}, in virtue of \eqref{psi} and \eqref{h}. In the following we focus on the proof of \eqref{Fe2}--\eqref{Fe6}.

\noindent{\bf \underline{Proof of \eqref{Fe2}.}}
Since
\begin{align*}
\partial_{t}\partial_{y}\hat{F}
=&-\partial_{t}\partial_{x}u\tilde{\varsigma}
  -\partial_{y}^{-1}[\partial_{t}\partial_{x}u]\tilde{\varsigma}'
  -\partial_{t}u\partial_{x}\partial_{y}u
  -u\partial_{t}\partial_{x}\partial_{y}u
  -h'\partial_{y}^{-1}[\partial_{t}\partial_{x}u]\partial_{y}u\\
 &-h'\partial_{y}^{-1}[\partial_{x}u]\partial_{t}\partial_{y}u
  -h\partial_{t}\partial_{x}u\partial_{y}u
  -h\partial_{x}u\partial_{t}\partial_{y}u
  +(1-h)\partial_{y}^{-1}[\partial_{t}\partial_{x}u]\partial_{y}^{2}u\nonumber\\
 &+(1-h)\partial_{y}^{-1}[\partial_{x}u]\partial_{t}\partial_{y}^{2}u,\nonumber
\end{align*}
for any $t\in(0,t_{*}]$, by using \eqref{dte}, \eqref{ue} and Sobolev imbedding theorem, one can deduce
\begin{align}
\|\partial_{t}\partial_{y}\hat{F}\|_{L^{\infty}(\Omega)}
\leq&~C\bar{\tilde{\varsigma}}\|\partial_{t}\partial_{x}u\|_{L^{\infty}(\Omega)}
     +\|\partial_{t}u\|_{L^{\infty}(\Omega)}\|\partial_{x}\partial_{y}u\|_{L^{\infty}(\Omega)}
     +\|u\|_{L^{\infty}(\Omega)}\|\partial_{t}\partial_{x}\partial_{y}u\|_{L^{\infty}(\Omega)}\\
    &+C\|\partial_{t}\partial_{x}u\|_{L^{\infty}(\Omega)}\|\partial_{y}u\|_{L^{\infty}(\Omega)}
     +C\|\partial_{x}u\|_{L^{\infty}(\Omega)}\|\partial_{t}\partial_{y}u\|_{L^{\infty}(\Omega)}\nonumber\\
    &+C\|\psi\partial_{t}\partial_{x}u\|_{L_{x}^{\infty}L_{y}^{2}(\Omega)}\|\partial_{y}^{2}u\|_{L^{\infty}(\Omega)}
     +C\|\psi\partial_{x}u\|_{L_{x}^{\infty}L_{y}^{2}(\Omega)}\|\partial_{t}\partial_{y}^{2}u\|_{L^{\infty}(\Omega)}\nonumber\\
\leq&~C_{*}\bar{\eta},\nonumber
\end{align}
thus,
\begin{align}\label{hatFye}
|\partial_{y}\hat{F}(t,x,y)|
\leq& |\partial_{y}\hat{F}(0,x,y)|+\int_{0}^{t}\|\partial_{t}\partial_{y}\hat{F}(\sigma,x,y)\|_{L^{\infty}(\Omega)}d\sigma\\
\leq& |\partial_{y}\hat{F}(0,x,y)|+C_{*}\bar{\eta}t,~\text{in}~\Omega_{t_{*}}.\nonumber
\end{align}
Notice that
\begin{align*}
\partial_{y}\hat{F}(0,x,y)
=-\partial_{x}\check{u}\tilde{\varsigma}
 -\partial_{y}^{-1}[\partial_{x}\check{u}]\tilde{\varsigma}'
 -\check{u}\partial_{x}\partial_{y}\check{u}
 -h'\partial_{y}^{-1}[\partial_{x}\check{u}]\partial_{y}\check{u}
 -h\partial_{x}\check{u}\partial_{y}\check{u}
 +(1-h)\partial_{y}^{-1}[\partial_{x}\check{u}]\partial_{y}^{2}\check{u},
\end{align*}
one has $|\partial_{y}\hat{F}(0,x,y)|\leq C\check{c}\leq C\iota^{-1}\bar{\varsigma}$ in $\Omega$, which leads to
\begin{align}\label{hatFye'}
|\partial_{y}\hat{F}|
\leq C\iota^{-1}\bar{\varsigma}+C_{*}\bar{\eta}t,~\text{in}~\Omega_{t_{*}}.
\end{align}
By using \eqref{dotub}, from \eqref{hatFye'} one knows
\begin{align}\label{hatFye-1}
|\partial_{y}\hat{F}|
\leq C\iota^{-1}\bar{\varsigma}+C_{*}\bar{\eta}t
\leq C\iota^{-1}\bar{\varsigma}+C_{*}\dot{\mathfrak{U}},~\text{in}~\Omega_{t_{*}}^{*}.
\end{align}
while using \eqref{dotu} and \eqref{Fe1} one has
\begin{align}\label{hatFye-2}
|\partial_{y}\hat{F}|
\leq C_{*}e^{-(1+\theta)\delta\y}
\leq C_{*}\dot{\mathfrak{U}},~\text{in}~\Omega_{t_{*}}\setminus\Omega_{t_{*}}^{*}.
\end{align}
Combining \eqref{hatFye-1} with \eqref{hatFye-2}, one gets \eqref{Fe2}.

\noindent{\bf \underline{Proof of \eqref{Fe3}.}}
Recall that $\check{u}=\tilde{u}_{in}$ or $0$. If $\check{u}=0$, $\partial_{y}\hat{F}(0,x,y)=0$ in $\Omega$, then \eqref{Fe3} follows from \eqref{hatFye}, while
if $\check{u}=\tilde{u}_{in}$, then one can choose $\iota$ small enough to obtain
\begin{align}
|\partial_{y}\hat{F}(0,x,y)|
\leq C\iota
    +C|\partial_{x}\partial_{y}\tilde{u}_{in}|
    +C|\partial_{y}\tilde{u}_{in}|
    +C|\partial_{y}^{2}\tilde{u}_{in}|
\leq \bar{\iota}+C|y-y_{*}|,~\text{in}~\Omega_{t_{*}}^{*},
\end{align}
which gives \eqref{Fe3} together with \eqref{hatFye}.

\noindent{\bf \underline{Proof of \eqref{Fe4}.}}
By a direct calculation, one has
\begin{align}\label{hatFy}
\partial_{y}\hat{F}
=-\partial_{x}u\tilde{\varsigma}
 -\partial_{y}^{-1}[\partial_{x}u]\tilde{\varsigma}'
 -u\partial_{x}\partial_{y}u
 -h'\partial_{y}^{-1}[\partial_{x}u]\partial_{y}u
 -h\partial_{x}u\partial_{y}u
 +(1-h)\partial_{y}^{-1}[\partial_{x}u]\partial_{y}^{2}u,
\end{align}
and
\begin{align}\label{hatFyy}
\partial_{y}^{2}\hat{F}
=&-\tilde{\varsigma}\partial_{x}\partial_{y}u
  -2\partial_{x}u\tilde{\varsigma}'
  -\partial_{y}^{-1}[\partial_{x}u]\tilde{\varsigma}''
  -\partial_{y}u\partial_{x}\partial_{y}u
  -u\partial_{x}\partial_{y}^{2}u
  -h''\partial_{y}^{-1}[\partial_{x}u]\partial_{y}u
  -2h'\partial_{x}u\partial_{y}u\\
 &-2h'\partial_{y}^{-1}[\partial_{x}u]\partial_{y}^{2}u
  -h\partial_{x}\partial_{y}u\partial_{y}u
  +(1-2h)\partial_{x}u\partial_{y}^{2}u
  +(1-h)\partial_{y}^{-1}[\partial_{x}u]\partial_{y}^{3}u.\nonumber
\end{align}
Utilizing \eqref{us-3}, one can deduce from \eqref{hatFy} and \eqref{hatFyy} that for any $t\in[0,t_{*}]$,
\begin{align*}
\|\partial_{y}\hat{F}\|_{L^{\infty}(\Omega)}
\leq&~\bar{\tilde{\varsigma}}\|\partial_{x}u\|_{L^{\infty}(\Omega)}
     +\|u\|_{L^{\infty}(\Omega)}\|\partial_{x}\partial_{y}u\|_{L^{\infty}(\Omega)}
     +\|\partial_{x}u\|_{L^{\infty}(\Omega)}\|\partial_{y}u\|_{L^{\infty}(\Omega)}\\
    &+C\|\psi\partial_{x}u\|_{L_{x}^{\infty}L_{y}^{2}(\Omega)}
       \left(\|\partial_{y}u\|_{L^{\infty}(\Omega)}+\|\partial_{y}^{2}u\|_{L^{\infty}(\Omega)}+\bar{\tilde{\varsigma}}\right)\nonumber\\
\leq&~C\check{c}+C_{*}\bar{\eta}t,\nonumber
\end{align*}
and
\begin{align*}
\|\partial_{y}^{2}\hat{F}\|_{L^{\infty}(\Omega)}
\leq&~\bar{\tilde{\varsigma}}\|\partial_{x}\partial_{y}u\|_{L^{\infty}(\Omega)}
     +C\bar{\tilde{\varsigma}}\|\partial_{x}u\|_{L^{\infty}(\Omega)}
     +2\|\partial_{y}u\|_{L^{\infty}(\Omega)}\|\partial_{x}\partial_{y}u\|_{L^{\infty}(\Omega)}
     +\|u\|_{L^{\infty}(\Omega)}\|\partial_{x}\partial_{y}^{2}u\|_{L^{\infty}(\Omega)}\\
    &+C\|\psi\partial_{x}u\|_{L_{x}^{\infty}L_{y}^{2}(\Omega)}
       \left(\|\partial_{y}u\|_{L^{\infty}(\Omega)}+\|\partial_{y}^{2}u\|_{L^{\infty}(\Omega)}+\|\partial_{y}^{3}u\|_{L^{\infty}(\Omega)}+\bar{\tilde{\varsigma}}\right)\\
    &+C\|\partial_{x}u\|_{L^{\infty}(\Omega)}\left(\|\partial_{y}u\|_{L^{\infty}(\Omega)}+\|\partial_{y}^{2}u\|_{L^{\infty}(\Omega)}\right)\\
\leq&~C\check{c}+C_{*}\bar{\eta}t,
\end{align*}
thus one knows that for any $t\in[0,t_{*}]$,
\begin{align*}
\|\mathcal{P}(\partial_{y}^{2}\hat{F}-\rho\partial_{y}\hat{F})\|_{L^{\infty}(\Omega)}
\leq&~C\|\mathcal{P}\|_{L^{\infty}(\Omega)}^{2}\|\partial_{y}F\|_{L^{\infty}(\Omega)}
     +\|\mathcal{P}\|_{L^{\infty}(\Omega)}\|\partial_{y}^{2}F\|_{L^{\infty}(\Omega)}\\
\leq&~C_{*}(\check{c}+\bar{\eta}t)\|\mathcal{P}\|_{L^{\infty}(\Omega)}^{2}
\leq C_{*}(\iota^{-1}+\bar{\eta})
\leq C_{*}\iota^{-1},
\end{align*}
which gives \eqref{Fe4}.

\noindent{\bf \underline{Proof of \eqref{Fe5}.}}
By using \eqref{priore}, it is easy to get that
\begin{align}\label{Fs-2}
\|\hat{F}\|_{L_{t}^{\infty}H_{\psi}^{s-2}(\Omega_{t_{*}})}
\leq C\bar{\tilde{\varsigma}}\|u\|_{L_{t}^{\infty}H_{\psi}^{s-1}(\Omega_{t_{*}})}
     +C\|u\|_{L_{t}^{\infty}H_{\psi}^{s-1}(\Omega_{t_{*}})}^{2}
\leq C_{*}\bar{\eta}.
\end{align}
By a direct calculation, one has
\begin{align*}
\partial_{t}\hat{F}
=-\partial_{y}^{-1}[\partial_{x}\partial_{t}u]\tilde{\varsigma}
 -\partial_{t}u\partial_{x}u
 -u\partial_{x}\partial_{t}u
 +(1-h)\partial_{y}^{-1}[\partial_{x}\partial_{t}u]\partial_{y}u
 +(1-h)\partial_{y}^{-1}[\partial_{x}u]\partial_{y}\partial_{t}u,~\text{in}~\Omega_{t_{*}},
\end{align*}
thus for any $t\in(0,t_{*}]$, with the help of \eqref{us-3} one has
\begin{align}\label{Fts-4}
\|\partial_{t}\hat{F}\|_{H_{\psi}^{s-4}(\Omega)}
\leq&~C\bar{\tilde{\varsigma}}\|\partial_{t}u\|_{H_{\psi}^{s-3}(\Omega)}
     +C\|u\|_{H_{\psi}^{s-3}(\Omega)}\|\partial_{t}u\|_{H_{\psi}^{s-3}(\Omega)}\\
\leq&~C_{*}\left(\bar{\tilde{\varsigma}}+\|u\|_{H_{\psi}^{s-3}(\Omega)}\right)
           \left(\|u\|_{H_{\psi}^{s-2}(\Omega)}^{2}+\|u\|_{H_{\psi}^{s-1}(\Omega)}+\bar{\kappa}+\|F\|_{H_{\psi}^{s-3}(\Omega)}\right)\nonumber\\
\leq&~C_{*}\bar{\eta}.\nonumber
\end{align}
Combining \eqref{Fs-2} with \eqref{Fts-4}, one has \eqref{Fe5}.

\noindent{\bf \underline{Proof of \eqref{Fe6}.}}
By a direct calculation one has,
\begin{align*}
 \partial_{y}\partial_{x}^{s}\hat{F}-\rho\partial_{x}^{s}\hat{F}
=\mathcal{M}_{1}+\mathcal{M}_{2}+\mathcal{M}_{3},
\end{align*}
where,
\begin{align*}
\mathcal{M}_{1}
=&-u(\partial_{x}^{s+1}\partial_{y}u-\rho\partial_{x}^{s+1}u)
  -\partial_{y}^{-1}[\partial_{x}^{s+1}u](\tilde{\varsigma}'-\rho\tilde{\varsigma})
  +(1-h)\partial_{y}^{-1}[\partial_{x}^{s+1}u](\partial_{y}^{2}u-\rho\partial_{y}u)\\
 &-\partial_{x}^{s+1}u\tilde{\varsigma}
  -h'\partial_{y}^{-1}[\partial_{x}^{s+1}u]\partial_{y}u
  -h\partial_{x}^{s+1}u\partial_{y}u,
\end{align*}
\begin{align*}
\mathcal{M}_{2}
=&-s\partial_{x}u(\partial_{x}^{s}\partial_{y}u-s\rho\partial_{x}^{s}u)
  -\partial_{x}^{s}u(\partial_{x}\partial_{y}u-\rho\partial_{x}u)
  +(1-h)\partial_{y}^{-1}[\partial_{x}u](\partial_{x}^{s}\partial_{y}^{2}u-\rho\partial_{x}^{s}\partial_{y}u)\\
 &+s(1-h)\partial_{y}^{-1}[\partial_{x}^{s}u](\partial_{x}\partial_{y}^{2}u-\rho\partial_{x}\partial_{y}u)
  -sh'\partial_{y}^{-1}[\partial_{x}^{s}u]\partial_{x}\partial_{y}u
  -sh\partial_{x}^{s}u\partial_{x}\partial_{y}u\\
 &-h'\partial_{y}^{-1}[\partial_{x}u]\partial_{x}^{s}\partial_{y}u
  -h\partial_{x}u\partial_{x}^{s}\partial_{y}u,
\end{align*}
and
\begin{align*}
\mathcal{M}_{3}
=&-\sum_{i=2}^{s-1}\binom{s}{i}\partial_{x}^{i}u\partial_{x}^{s-i+1}\partial_{y}u
  -h'\sum_{i=1}^{s-2}\binom{s}{i}\partial_{y}^{-1}[\partial_{x}^{i+1}u]\partial_{x}^{s-i}\partial_{y}u\\
 &-h\sum_{i=1}^{s-2}\binom{s}{i}\partial_{x}^{i+1}u\partial_{x}^{s-i}\partial_{y}u
  +(1-h)\sum_{i=1}^{s-2}\binom{s}{i}\partial_{y}^{-1}[\partial_{x}^{i+1}u]\partial_{x}^{s-i}\partial_{y}^{2}u\\
 &+\rho \sum_{i=2}^{s-1}\binom{s}{i}\partial_{x}^{i}u\partial_{x}^{s-i+1}u
  -(1-h)\rho\sum_{i=1}^{s-2}\partial_{y}^{-1}[\partial_{x}^{i+1}u]\partial_{x}^{s-i}\partial_{y}u.
\end{align*}
Accordingly,
\begin{align}\label{Fs}
\|\hat{\mathbf{F}}_{s}\|_{L_{t}^{2}L_{\hat{\psi}}^{2}(\Omega_{t_{*}})}^{2}
\leq \|\mathcal{P}\mathcal{M}_{1}\|_{L_{t}^{2}L_{\hat{\psi}}^{2}(\Omega_{t_{*}})}^{2}
    +\|\mathcal{P}\mathcal{M}_{2}\|_{L_{t}^{2}L_{\hat{\psi}}^{2}(\Omega_{t_{*}})}^{2}
    +\|\mathcal{P}\mathcal{M}_{3}\|_{L_{t}^{2}L_{\hat{\psi}}^{2}(\Omega_{t_{*}})}^{2}.
\end{align}

Owing to \eqref{rhop1} and \eqref{dotub}, one has $\rho\leq C\mathcal{P}\leq C\tilde{\mathcal{P}}$.
Recalling that $u=\mathbb{U}-\tilde{\mathbb{U}}$ and $\tilde{\mathbb{U}}=0, \bar{\eta}=\epsilon^{2}$ if $\check{c}\neq0$, for any $t\in[0,t_{*}]$, one can obtain by using \eqref{us-3} that
\begin{align}\label{pm11}
&\|\mathcal{P}u(\partial_{x}^{s+1}\partial_{y}u-\rho\partial_{x}^{s+1}u)\|_{L_{\hat{\psi}}^{2}(\Omega)}^{2}\\
\leq&~\|u\mathcal{P}(\partial_{x}^{s+1}\partial_{y}\mathbb{U}-\rho\partial_{x}^{s+1}\mathbb{U})\|_{L_{\hat{\psi}}^{2}(\Omega)}^{2}
     +\|u\mathcal{P}(\partial_{x}^{s+1}\partial_{y}\tilde{\mathbb{U}}-\tilde{\rho}\partial_{x}^{s+1}\tilde{\mathbb{U}})\|_{L_{\hat{\psi}}^{2}(\Omega)}^{2}
     +\|u\mathcal{P}(\rho-\tilde{\rho})\partial_{x}^{s+1}\tilde{\mathbb{U}}\|_{L_{\varphi}^{2}(\Omega)}^{2}\nonumber\\
\leq&~C_{\lambda}\|\sqrt{\y}u\|_{L^{\infty}(\Omega)}^{2}\left(\|\mathbf{U}_{s+1}\|_{L_{\hat{\varphi}}^{2}(\Omega)}^{2}
                                                     +\|\tilde{\mathbf{U}}_{s+1}\|_{L_{\tilde{\hat{\varphi}}}^{2}(\Omega)}^{2}
                                                     +\|\tilde{\mathcal{P}}\tilde{\mathbf{U}}_{s+1}\|_{L_{\tilde{\hat{\varphi}}}^{2}(\Omega)}^{2}\right)\nonumber\\
\leq&~(C\check{c}^{2}+C_{*}\bar{\eta}t^{2})\left(\|\mathbf{U}_{s+1}\|_{L_{\hat{\varphi}}^{2}(\Omega)}^{2}
                                                     +\|\tilde{\mathbf{U}}_{s+1}\|_{L_{\tilde{\hat{\varphi}}}^{2}(\Omega)}^{2}
                                                     +\|\tilde{\mathcal{P}}\tilde{\mathbf{U}}_{s+1}\|_{L_{\tilde{\hat{\varphi}}}^{2}(\Omega)}^{2}\right)\nonumber\\
\leq&~C\check{c}^{2}\|\mathbf{U}_{s+1}\|_{L_{\hat{\varphi}}^{2}(\Omega)}^{2}
     +C_{*}\bar{\eta}t\left(\|\mathbf{U}_{s+1}\|_{L_{\hat{\varphi}}^{2}(\Omega)}^{2}
                                                     +\|\tilde{\mathbf{U}}_{s+1}\|_{L_{\tilde{\hat{\varphi}}}^{2}(\Omega)}^{2}\right)\nonumber\\
\leq& C\bar{\eta}+C_{*}\bar{\eta}t,\nonumber
\end{align}
and
\begin{align}\label{pm12}
&\|\mathcal{P}\partial_{y}^{-1}[\partial_{x}^{s+1}u](\tilde{\varsigma}'-\rho\tilde{\varsigma})\|_{L_{\hat{\psi}}^{2}(\Omega)}^{2}\\
\leq&~\|\mathcal{P}\partial_{y}^{-1}[\partial_{x}^{s+1}\mathbb{U}](\tilde{\varsigma}'-\rho\tilde{\varsigma})\|_{L_{\hat{\psi}}^{2}(\Omega)}^{2}
     +\|\mathcal{P}\partial_{y}^{-1}[\partial_{x}^{s+1}\tilde{\mathbb{U}}](\tilde{\varsigma}'-\rho\tilde{\varsigma})\|_{L_{\hat{\psi}}^{2}(\Omega)}^{2}\nonumber\\
\leq&~C\tilde{\bar{\varsigma}}^{2}\|\mathcal{P}^{2}\|_{L_{x}^{\infty}L_{y}^{2}(\hat{\Omega})}^{2}\|\partial_{x}^{s+1}\mathbb{U}\|_{L_{\check{\psi}}^{2}(\Omega)}^{2}
     +C\tilde{\bar{\varsigma}}^{2}\|\mathcal{P}^{2}\|_{L_{x}^{\infty}L_{y}^{2}(\hat{\Omega})}^{2}\|\partial_{x}^{s+1}\tilde{\mathbb{U}}\|_{L_{\check{\psi}}^{2}(\Omega)}^{2}\nonumber\\
\leq&~C_{*}\sqrt{\iota}\bar{\eta}(\|\mathcal{P}\mathbf{U}_{s+1}\|_{L_{\hat{\varphi}}^{2}(\Omega)}^{2}+\|\tilde{\mathcal{P}}\tilde{\mathbf{U}}_{s+1}\|_{L_{\tilde{\hat{\varphi}}}^{2}(\Omega)}^{2})\nonumber.
\end{align}
Also,
\begin{align}\label{pm13}
&\|\mathcal{P}(1-h)\partial_{y}^{-1}[\partial_{x}^{s+1}u](\partial_{y}^{2}u-\rho\partial_{y}u)\|_{L_{\hat{\psi}}^{2}(\Omega)}^{2}\\
\leq&~\|\mathcal{P}\partial_{y}^{-1}[\partial_{x}^{s+1}\mathbb{U}](\partial_{y}^{2}u-\rho\partial_{y}u)\|_{L_{\check{\psi}}^{2}(\hat{\Omega})}^{2}
     +\|\mathcal{P}\partial_{y}^{-1}[\partial_{x}^{s+1}\tilde{\mathbb{U}}](\partial_{y}^{2}u-\rho\partial_{y}u)\|_{L_{\check{\psi}}^{2}(\hat{\Omega})}^{2}\nonumber\\
\leq&~C(\|\partial_{y}u\|_{L^{\infty}(\Omega)}+\|\partial_{y}^{2}u\|_{L^{\infty}(\Omega)})^{2}\|\mathcal{P}^{2}\|_{L_{x}^{\infty}L_{y}^{2}(\hat{\Omega})}^{2}
       (\|\partial_{x}^{s+1}\mathbb{U}\|_{L_{\check{\psi}}^{2}(\Omega)}^{2}+\|\partial_{x}^{s+1}\tilde{\mathbb{U}}\|_{L_{\check{\psi}}^{2}(\Omega)}^{2})\nonumber\\
\leq&~C_{*}(\check{c}^{2}+\bar{\eta}t^{2})\|\mathcal{P}^{2}\|_{L_{x}^{\infty}L_{y}^{2}(\hat{\Omega})}^{2}
           (\|\mathcal{P}\mathbf{U}_{s+1}\|_{L_{\hat{\varphi}}^{2}(\Omega)}^{2}+\|\tilde{\mathcal{P}}\tilde{\mathbf{U}}_{s+1}\|_{L_{\tilde{\hat{\varphi}}}^{2}(\Omega)}^{2})\nonumber\\
\leq&~C_{*}\bar{\eta}\sqrt{t}(\|\mathcal{P}\mathbf{U}_{s+1}\|_{L_{\hat{\varphi}}^{2}(\Omega)}^{2}+\|\tilde{\mathcal{P}}\tilde{\mathbf{U}}_{s+1}\|_{L_{\tilde{\hat{\varphi}}}^{2}(\Omega)}^{2})
     +C_{*}\iota^{-\frac{5}{2}}\check{c}^{2}\|\mathbf{U}_{s+1}\|_{L_{\hat{\varphi}}^{2}(\Omega)}^{2}\nonumber\\
\leq&~C_{*}\bar{\eta}\sqrt{t}(\|\mathcal{P}\mathbf{U}_{s+1}\|_{L_{\hat{\varphi}}^{2}(\Omega)}^{2}+\|\tilde{\mathcal{P}}\tilde{\mathbf{U}}_{s+1}\|_{L_{\tilde{\hat{\varphi}}}^{2}(\Omega)}^{2})
     +C_{*}\iota^{-\frac{5}{2}}\epsilon^{-2}\bar{\eta}.\nonumber
\end{align}
By using the Corollary \ref{ncC3}, one can deduce
\begin{align}\label{pm14}
\|\mathcal{P}\partial_{x}^{s+1}u\tilde{\varsigma}\|_{L_{\hat{\psi}}^{2}(\Omega)}^{2}
\leq&~C_{\lambda}\tilde{\bar{\varsigma}}^{2}\left(\|\partial_{x}^{s+1}\mathbb{U}\|_{L_{\check{\psi}}^{2}(\Omega)}^{2}+\|\partial_{x}^{s+1}\tilde{\mathbb{U}}\|_{L_{\check{\psi}}^{2}(\Omega)}^{2}\right)\\
\leq&~C_{\lambda}\tilde{\bar{\varsigma}}^{2}\left(\|\mathcal{P}\mathbf{U}_{s+1}\|_{L_{\hat{\varphi}}^{2}(\Omega)}^{2}+\|\tilde{\mathcal{P}}\tilde{\mathbf{U}}_{s+1}\|_{L_{\tilde{\hat{\varphi}}}^{2}(\Omega)}^{2}\right)\nonumber\\
\leq&~C_{*}\tilde{\bar{\varsigma}}\left(\|\mathbf{U}_{s+1}\|_{L_{\hat{\varphi}}^{2}(\Omega)}^{2}+\|\tilde{\mathbf{U}}_{s+1}\|_{L_{\tilde{\hat{\varphi}}}^{2}(\Omega)}^{2}\right)\nonumber\\
\leq&~C_{*}\bar{\eta}.\nonumber
\end{align}
Using \eqref{us-3} and noticing $\tilde{\mathbb{U}}=0$ if $\check{c}\neq0$ again, we have
\begin{align*}
\|\mathcal{P}h'\partial_{y}^{-1}[\partial_{x}^{s+1}u]\partial_{y}u\|_{L_{\hat{\psi}}^{2}(\Omega)}^{2}
\leq&~\|\partial_{y}u\|_{L_{x}^{\infty}L_{y}^{2}(\Omega)}^{2}
      \left(\|\partial_{x}^{s+1}\mathbb{U}\|_{L_{\check{\psi}}^{2}(\Omega)}^{2}+\|\partial_{x}^{s+1}\tilde{\mathbb{U}}\|_{L_{\check{\psi}}^{2}(\Omega)}^{2}\right)\\
\leq&~C_{*}(\check{c}^{2}+\bar{\eta}t^{2})\left(\|\mathcal{P}\mathbf{U}_{s+1}\|_{L_{\hat{\varphi}}^{2}(\Omega)}^{2}+\|\tilde{\mathcal{P}}\tilde{\mathbf{U}}_{s+1}\|_{L_{\tilde{\hat{\varphi}}}^{2}(\Omega)}^{2}\right)\\
\leq&~C_{*}\bar{\eta}t\left(\|\mathbf{U}_{s+1}\|_{L_{\hat{\varphi}}^{2}(\Omega)}^{2}+\|\tilde{\mathbf{U}}_{s+1}\|_{L_{\tilde{\hat{\varphi}}}^{2}(\Omega)}^{2}\right)
     +C_{*}\iota^{-1}\check{c}^{2}\|\mathbf{U}_{s+1}\|_{L_{\hat{\varphi}}^{2}(\Omega)}^{2}\nonumber\\
\leq&~C_{*}\bar{\eta}t
     +C_{*}\iota^{-1}\bar{\eta}.
\end{align*}
Similarly,
\begin{align}\label{pm15}
\|\mathcal{P}h\partial_{x}^{s+1}u\partial_{y}u\|_{L_{\hat{\psi}}^{2}(\Omega)}^{2}
\leq&~C_{\lambda}\|\y\partial_{y}u\|_{L^{\infty}(\Omega)}^{2}
      \left(\|\partial_{x}^{s+1}\mathbb{U}\|_{L_{\check{\psi}}^{2}(\Omega)}^{2}+\|\partial_{x}^{s+1}\tilde{\mathbb{U}}\|_{L_{\check{\psi}}^{2}(\Omega)}^{2}\right)\\
\leq&~C_{*}(\check{c}^{2}+\bar{\eta}t^{2})\left(\|\mathcal{P}\mathbf{U}_{s+1}\|_{L_{\hat{\varphi}}^{2}(\Omega)}^{2}+\|\tilde{\mathcal{P}}\tilde{\mathbf{U}}_{s+1}\|_{L_{\tilde{\hat{\varphi}}}^{2}(\Omega)}^{2}\right)\nonumber\\
\leq&~C_{*}\bar{\eta}t\left(\|\mathbf{U}_{s+1}\|_{L_{\hat{\varphi}}^{2}(\Omega)}^{2}+\|\tilde{\mathbf{U}}_{s+1}\|_{L_{\tilde{\hat{\varphi}}}^{2}(\Omega)}^{2}\right)
     +C_{*}\iota^{-1}\check{c}^{2}\|\mathbf{U}_{s+1}\|_{L_{\hat{\varphi}}^{2}(\Omega)}^{2}\nonumber\\
\leq&~C_{*}\bar{\eta}t+C\iota^{-1}\bar{\eta}.\nonumber
\end{align}
Combining \eqref{pm11}--\eqref{pm15}, one can deduce that for any $t\in[0,t_{*}]$,
\begin{align*}
\|\mathcal{P}\mathcal{M}_{1}\|_{L_{\hat{\psi}}^{2}(\Omega)}^{2}
\leq  C\iota^{-\frac{5}{2}}\bar{\eta}+C_{*}\bar{\eta}t
     +C_{*}(\sqrt{\iota}+\sqrt{t})\bar{\eta}(\|\mathcal{P}\mathbf{U}_{s+1}\|_{L_{\hat{\varphi}}^{2}(\Omega)}^{2}
            +\|\tilde{\mathcal{P}}\tilde{\mathbf{U}}_{s+1}\|_{L_{\tilde{\hat{\varphi}}}^{2}(\Omega)}^{2}),
\end{align*}
which gives
\begin{align}\label{pm1}
\|\mathcal{P}\mathcal{M}_{1}\|_{L_{t}^{2}L_{\hat{\psi}}^{2}(\Omega_{t_{*}})}^{2}\leq C_{*}\bar{\eta}.
\end{align}

For any $t\in[0,t_{*}]$, notice that $\mathcal{U}_{s}=\tilde{\mathcal{P}}\partial_{x}u(\partial_{x}^{s}\partial_{y}u-\tilde{\rho}\partial_{x}^{s}u)$,
by using \eqref{priore} and \eqref{us-3} we have
\begin{align}\label{pm21}
&s\|\mathcal{P}\partial_{x}u(\partial_{x}^{s}\partial_{y}u-\rho\partial_{x}^{s}u)\|_{L_{\hat{\psi}}^{2}(\Omega)}^{2}\\
\leq&~C_{\lambda}\|\partial_{x}u\|_{L^{\infty}(\Omega)}^{2}\|\tilde{\mathcal{P}}(\partial_{x}^{s}\partial_{y}u-\tilde{\rho}\partial_{x}^{s}u)\|_{L_{\varphi}^{2}(\Omega)}^{2}
     +C_{\lambda}\|\partial_{x}u\|_{L^{\infty}(\Omega)}^{2}\|\tilde{\mathcal{P}}(\rho-\tilde{\rho})\partial_{x}^{s}u\|_{L_{\varphi}^{2}(\Omega)}^{2}\nonumber\\
\leq&~C_{\lambda}\|\partial_{x}u\|_{L^{\infty}(\Omega)}^{2}\|\mathcal{U}_{s}\|_{L_{\tilde{\varphi}}^{2}(\Omega)}^{2}
     +C_{\lambda}\|\partial_{x}u\|_{L^{\infty}(\Omega)}^{2}\|\partial_{x}^{s}u\|_{L_{\tilde{\hat{\Psi}}}^{2}(\Omega)}^{2}\nonumber\\
\leq&~C_{*}(\check{c}^{2}+\bar{\eta}t^{2})(\|\mathcal{U}_{s}\|_{L_{\tilde{\varphi}}^{2}(\Omega)}^{2}
                                           +\|\tilde{\mathcal{P}}\mathcal{U}_{s}\|_{L_{\tilde{\varphi}}^{2}(\Omega)}^{2})\nonumber\\
\leq&~C_{*}\bar{\eta}t\|\mathcal{U}_{s}\|_{L_{\tilde{\varphi}}^{2}(\Omega)}^{2}
     +C_{*}\iota^{-1}\check{c}^{2}\|\mathcal{U}_{s}\|_{L_{\tilde{\varphi}}^{2}(\Omega)}^{2}\nonumber\\
\leq&~C_{*}\iota^{-1}\bar{\eta}.\nonumber
\end{align}
In virtue of \eqref{us-3} and Corollary \ref{ncC2}, one has
\begin{align}\label{pm22}
s\|\mathcal{P}\partial_{x}^{s}u(\partial_{x}\partial_{y}u-\rho\partial_{x}u)\|_{L_{\hat{\psi}}^{2}(\Omega)}^{2}
\leq&~C_{*}\left(\|\partial_{x}\partial_{y}u\|_{L^{\infty}(\Omega)}^{2}+\|\partial_{x}u\|_{L^{\infty}(\Omega)}^{2}\right)
          \|\partial_{x}^{s}u\|_{L_{\tilde{\hat{\Psi}}}^{2}(\Omega)}^{2}\\
\leq&~C_{*}(\check{c}^{2}+\bar{\eta}t^{2})\|\tilde{\mathcal{P}}\mathcal{U}_{s}\|_{L_{\tilde{\varphi}}^{2}(\Omega)}^{2}\nonumber\\
\leq&~C_{*}\bar{\eta}t\|\mathcal{U}_{s}\|_{L_{\tilde{\varphi}}^{2}(\Omega)}^{2}
     +C_{*}\iota^{-1}\check{c}^{2}\|\mathcal{U}_{s}\|_{L_{\tilde{\varphi}}^{2}(\Omega)}^{2}\nonumber\\
\leq&~C_{*}\iota^{-1}\bar{\eta}.\nonumber
\end{align}
Since
$\tilde{\mathcal{P}}(\partial_{y}^{2}\partial_{x}^{s}u-\tilde{\rho}\partial_{x}^{s}\partial_{y}u)
=\partial_{y}\mathcal{U}_{s}
+\tilde{\mathcal{P}}\partial_{y}\tilde{\mathcal{P}}^{-1}\mathcal{U}_{s}
+\tilde{\mathcal{P}}\partial_{y}\tilde{\rho}\partial_{x}^{s}u
+\tilde{\mathcal{P}}(\tilde{\rho}-\rho)\mathcal{U}_{s}
+\tilde{\mathcal{P}}(\tilde{\rho}-\rho)\tilde{\rho}\partial_{x}^{s}u$,
one has
\begin{align}\label{pm23}
&\|\mathcal{P}(1-h)\partial_{y}^{-1}[\partial_{x}u](\partial_{x}^{s}\partial_{y}^{2}u-\rho\partial_{x}^{s}\partial_{y}u)\|_{L_{\hat{\psi}}^{2}(\Omega)}^{2}\\
\leq&~C\|\partial_{x}u\|_{L_{x}^{\infty}L_{y}^{2}(\Omega)}^{2}\|\tilde{\mathcal{P}}(\partial_{x}^{s}\partial_{y}^{2}u-\rho\partial_{x}^{s}\partial_{y}u)\|_{L_{\hat{\psi}}^{2}(\Omega)}^{2}\nonumber\\
\leq&~C_{*}\|\partial_{x}u\|_{L_{x}^{\infty}L_{y}^{2}(\Omega)}^{2}\|\tilde{\mathcal{P}}\|_{L^{\infty}(\Omega)}^{2}
       (\|\partial_{y}\mathcal{U}_{s}\|_{L_{\tilde{\varphi}}^{2}(\Omega)}^{2}
        +\|\tilde{\mathcal{P}}\mathcal{U}_{s}\|_{L_{\tilde{\varphi}}^{2}(\Omega)}^{2}
        +\|\partial_{x}^{s}u\|_{L_{\tilde{\hat{\Psi}}}^{2}(\Omega)}^{2})\nonumber\\
\leq&~C_{*}(\check{c}^{2}+\bar{\eta}t^{2})\|\tilde{\mathcal{P}}\|_{L^{\infty}(\Omega)}^{2}
       (\|\partial_{y}\mathcal{U}_{s}\|_{L_{\tilde{\varphi}}^{2}(\Omega)}^{2}+\|\tilde{\mathcal{P}}\mathcal{U}_{s}\|_{L_{\tilde{\varphi}}^{2}(\Omega)}^{2})\nonumber\\
\leq&~C_{*}\iota^{-2}\bar{\eta}
     +C_{*}(\check{c}^{2}\iota^{-1}+\bar{\eta}t)\|\partial_{y}\mathcal{U}_{s}\|_{L_{\tilde{\varphi}}^{2}(\Omega)}^{2}.\nonumber
\end{align}
Utilizing \eqref{us-3} and Corollary \ref{ncC2}, one can deduce
\begin{align}\label{pm24}
&s\|\mathcal{P}(1-h)\partial_{y}^{-1}[\partial_{x}^{s}u](\partial_{x}\partial_{y}^{2}u-\rho\partial_{x}\partial_{y}u)\|_{L_{\hat{\psi}}^{2}(\Omega)}^{2}\\
\leq&~C_{*}\left(\|\partial_{x}\partial_{y}^{2}u\|_{L_{x}^{\infty}L_{y}^{2}(\Omega)}^{2}+\|\partial_{x}\partial_{y}u\|_{L_{x}^{\infty}L_{y}^{2}(\Omega)}^{2}\right)
        \|\mathcal{P}\|_{L^{\infty}(\Omega)}^{4}\|\partial_{x}^{s}u\|_{L_{\tilde{\hat{\Psi}}}^{2}(\Omega)}^{2}\nonumber\\
\leq&~C_{*}(\check{c}^{2}+\bar{\eta}t^{2})\|\mathcal{P}\|_{L^{\infty}(\Omega)}^{4}\|\tilde{\mathcal{P}}\mathcal{U}_{s}\|_{L_{\tilde{\varphi}}^{2}(\Omega)}^{2}\nonumber\\
\leq&~C_{*}(\check{c}^{2}\iota^{-2}+\bar{\eta})\|\tilde{\mathcal{P}}\mathcal{U}_{s}\|_{L_{\tilde{\varphi}}^{2}(\Omega)}^{2}.\nonumber
\end{align}
Similarly, one has
\begin{align}\label{pm25}
s\|\mathcal{P}\partial_{x}\partial_{y}u(h'\partial_{y}^{-1}[\partial_{x}^{s}u]+h\partial_{x}^{s}u)\|_{L_{\hat{\psi}}^{2}(\Omega)}^{2}
\leq&~C_{*}\|\partial_{x}\partial_{y}u\|_{L^{\infty}(\Omega)}^{2}
        \|\mathcal{P}\|_{L^{\infty}(\Omega)}^{2}\|\partial_{x}^{s}u\|_{L_{\tilde{\hat{\Psi}}}^{2}(\Omega)}^{2}\\
\leq&~C_{*}(\check{c}^{2}+\bar{\eta}t^{2})\|\mathcal{P}\|_{L^{\infty}(\Omega)}^{2}\|\tilde{\mathcal{P}}\mathcal{U}_{s}\|_{L_{\tilde{\varphi}}^{2}(\Omega)}^{2}\nonumber\\
\leq&~C_{*}(\check{c}^{2}\iota^{-2}+\bar{\eta})\|\mathcal{U}_{s}\|_{L_{\tilde{\varphi}}^{2}(\Omega)}^{2}\nonumber\\
\leq&~C_{*}(\check{c}^{2}\iota^{-2}+\bar{\eta})\bar{\eta},\nonumber
\end{align}
and
\begin{align}\label{pm26}
&\|\mathcal{P}(h'\partial_{y}^{-1}[\partial_{x}u]+h\partial_{x}u)\partial_{x}^{s}\partial_{y}u\|_{L_{\hat{\psi}}^{2}(\Omega)}^{2}\\
\leq&~C_{\lambda}\left(\|\partial_{x}u\|_{L_{x}^{\infty}L_{y}^{2}(\Omega)}^{2}+\|\partial_{x}u\|_{L^{\infty}(\Omega)}^{2}\right)
       \left(\|\mathcal{U}_{s}\|_{L_{\tilde{\varphi}}^{2}(\Omega)}^{2}+\|\tilde{\mathcal{P}}\tilde{\rho}\partial_{x}^{s}u\|_{L_{\tilde{\varphi}}^{2}(\Omega)}^{2}\right)\nonumber\\
\leq&~C_{\lambda}(\check{c}^{2}+\bar{\eta}t^{2})
       \left(\|\mathcal{U}_{s}\|_{L_{\tilde{\varphi}}^{2}(\Omega)}^{2}+\|\tilde{\mathcal{P}}\mathcal{U}_{s}\|_{L_{\tilde{\varphi}}^{2}(\Omega)}^{2}\right)\nonumber\\
\leq&~C_{\lambda}(\check{c}^{2}\iota^{-1}+\bar{\eta}t)\|\mathcal{U}_{s}\|_{L_{\tilde{\varphi}}^{2}(\Omega)}^{2}\nonumber\\
\leq&~C_{*}\bar{\eta}.\nonumber
\end{align}
Collecting the estimates in \eqref{pm21}--\eqref{pm26}, for any $t\in[0,t_{*}]$ one has
\begin{align*}
\|\mathcal{P}\mathcal{M}_{2}\|_{L_{\hat{\psi}}^{2}(\Omega)}^{2}
\leq C_{*}\iota^{-2}\bar{\eta}
    +C_{*}(\check{c}^{2}\iota^{-1}+\bar{\eta}t)\|\partial_{y}\mathcal{U}_{s}\|_{L_{\tilde{\varphi}}^{2}(\Omega)}^{2}
    +C_{*}(\check{c}^{2}\iota^{-2}+\bar{\eta})\|\tilde{\mathcal{P}}\mathcal{U}_{s}\|_{L_{\tilde{\varphi}}^{2}(\Omega)}^{2}
\end{align*}
which gives
\begin{align}\label{pm2}
\|\mathcal{P}\mathcal{M}_{2}\|_{L_{t}^{2}L_{\hat{\psi}}^{2}(\Omega_{t_{*}})}^{2}
\leq& C_{*}\iota^{-2}\bar{\eta}
     +C_{*}(\check{c}^{2}\iota^{-1}+\bar{\eta}t)\|\partial_{y}\mathcal{U}_{s}\|_{L_{t}^{2}L_{\tilde{\varphi}}^{2}(\Omega_{t_{*}})}^{2}
     +C_{*}(\check{c}^{2}\iota^{-2}+\bar{\eta})\|\tilde{\mathcal{P}}\mathcal{U}_{s}\|_{L_{t}^{2}L_{\tilde{\varphi}}^{2}(\Omega_{t_{*}})}^{2}\\
\leq& C_{*}\iota^{-2}\bar{\eta}.\nonumber
\end{align}

The estimate of $\|\mathcal{P}\mathcal{M}_{3}\|_{L_{t}^{2}L_{\hat{\psi}}^{2}(\Omega_{t_{*}})}^{2}$ is straightforward. For any $t\in[0,t_{*}]$, by using \eqref{priore} and Sobolev embedding theorem one has
\begin{align*}
\|\mathcal{P}\mathcal{M}_{3}\|_{L_{\hat{\psi}}^{2}(\Omega)}^{2}
\leq C\|\mathcal{P}\|_{L^{\infty}(\Omega)}^{4}
      \|u\|_{H_{\psi}^{s-3}(\Omega)}^{2}
      \|u\|_{\hat{H}_{\psi}^{s}(\Omega)}^{2}
\leq C_{*}(\check{c}^{2}\iota^{-2}+\bar{\eta})\|u\|_{\hat{H}_{\psi}^{s}(\Omega)}^{2}
\leq  C_{*}\iota^{-2}\bar{\eta},
\end{align*}
thus
\begin{align}\label{pm3}
\|\mathcal{P}\mathcal{M}_{3}\|_{L_{t}^{2}L_{\hat{\psi}}^{2}(\Omega_{t_{*}})}^{2}
\leq C_{*}\iota^{-2}\bar{\eta}.
\end{align}

Substituting \eqref{pm1},\eqref{pm2} and \eqref{pm3} into \eqref{Fs}, it follows
\begin{align}\label{Fse}
\|\hat{\mathbf{F}}_{s}\|_{L_{t}^{2}L_{\hat{\psi}}^{2}(\Omega_{t_{*}})}^{2}\leq C_{*}\iota^{-2}\bar{\eta}.
\end{align}

Making use of Corollary \ref{ncC3}, one can deduce
\begin{align}\label{FL2Hs}
\|\hat{F}\|_{L_{t}^{2}\hat{H}_{\psi}^{s}(\Omega_{t_{*}})}
\leq&~C\|\partial_{y}u\|_{L_{t}^{2}\hat{H}_{\psi}^{s}(\Omega_{t_{*}})}\|u\|_{L_{t}^{\infty}\hat{H}_{\psi}^{s}(\Omega_{t_{*}})}\\
    &+\left(\bar{\tilde{\varsigma}}+\|u\|_{L_{t}^{\infty}\hat{H}_{\psi}^{s}(\Omega_{t_{*}})}\right)
      \left(\|\partial_{y}\partial_{x}^{s}u\|_{L_{t}^{2}L_{\hat{\psi}}^{2}(\Omega_{t_{*}})}+\|\partial_{x}^{s}u\|_{L_{t}^{2}L_{\hat{\psi}}^{2}(\Omega_{t_{*}})}\right)\nonumber\\
    &+C\|u\|_{L_{t}^{\infty}\hat{H}_{\psi}^{s}(\Omega_{t_{*}})}\left(\|\partial_{y}^{2}\partial_{x}^{s}u\|_{L_{t}^{2}L_{\hat{\psi}}^{2}(\Omega_{t_{*}})}
                            +\|\partial_{y}\partial_{x}^{s}u\|_{L_{t}^{2}L_{\hat{\psi}}^{2}(\Omega_{t_{*}})}
                            +\|\partial_{x}^{s}u\|_{L_{t}^{2}L_{\hat{\psi}}^{2}(\Omega_{t_{*}})}\right)\nonumber\\
    &+C\|u\|_{L_{t}^{\infty}\hat{H}_{\psi}^{s}(\Omega_{t_{*}})}\left(\|\partial_{y}u\|_{L_{t}^{2}\hat{H}_{\psi}^{s}(\Omega_{t_{*}})}
                            +\|u\|_{L_{t}^{\infty}\hat{H}_{\psi}^{s}(\Omega_{t_{*}})}\right)\nonumber\\
\leq&~C_{*}\left(\bar{\tilde{\varsigma}}+\|u\|_{L_{t}^{\infty}\hat{H}_{\psi}^{s}(\Omega_{t_{*}})}\right)
           \left(\|\partial_{y}u\|_{L_{t}^{2}\hat{H}_{\psi}^{s}(\Omega_{t_{*}})}
                 +\|u\|_{L_{t}^{\infty}\mathcal{H}_{\psi,\tilde{\varphi}}^{s}(\Omega_{t_{*}})}\right)\nonumber\\
    &+C_{*}\left(\bar{\tilde{\varsigma}}+\|u\|_{L_{t}^{\infty}\hat{H}_{\psi}^{s}(\Omega_{t_{*}})}\right)
           \left(\|\partial_{y}\mathcal{U}_{s}\|_{L_{t}^{2}L_{\tilde{\varphi}}^{2}(\Omega)}
                 +\|\mathcal{P}\mathcal{U}_{s}\|_{L_{t}^{2}L_{\tilde{\varphi}}^{2}(\Omega)}\right)\nonumber\\
\leq&~C_{*}\bar{\eta}.\nonumber
\end{align}

Finally, let's focus on the estimate of $\hat{F}$ on the boundary $y=0$.
Noticing that ${\bf supp}\tilde{\varsigma}\subset(\check{y},\hat{y})$ and $u|_{y=0}=0$, for any integer $k\in[0,\frac{s-1}{2}], t\in(0,t_{*}]$ and $\alpha\in\Gamma_{s-1}$,
\begin{align*}
\|\partial_{t}^{k}[D^{\alpha}\hat{F}]\|_{L_{x}^{2}(\mathbb{T})}|_{y=0}
\leq&~\|\partial_{t}^{k}[D^{\alpha}(\partial_{y}^{-1}[\partial_{x}u]\partial_{y}u)]\|_{L_{x}^{2}(\mathbb{T})}|_{y=0}
     +\|\partial_{t}^{k}[D^{\alpha}(u\partial_{x}u)]\|_{L_{x}^{2}(\mathbb{T})}|_{y=0}\\
\leq&~C\sum_{i=0}^{k}\sum_{2\leq|\beta|\leq|\alpha|}\|\partial_{t}^{i}D^{\beta}[\partial_{y}^{-1}[\partial_{x}u]]\partial_{t}^{k-i}D^{\alpha-\beta}[\partial_{y}u]\|_{L_{x}^{2}(\mathbb{T})}|_{y=0}\\
    &+C\sum_{i=0}^{k}\sum_{1\leq|\beta|\leq|\alpha|}\|\partial_{t}^{i}D^{\beta}u\partial_{t}^{k-i}D^{\alpha-\beta}[\partial_{x}u]\|_{L_{x}^{2}(\mathbb{T})}|_{y=0}\\
\leq&~C\sum_{i=0}^{k}\|\partial_{t}^{i}[D^{\alpha}u]\|_{L_{x}^{2}(\mathbb{T})}^{2}|_{y=0}.
\end{align*}
It follows from the Sobolev embedding theorem that
\begin{align*}
\|D^{\alpha}\hat{F}\|_{L_{x}^{2}(\mathbb{T})}|_{y=0}
\leq C\|D^{\alpha}u\|_{L_{x}^{2}(\mathbb{T})}^{2}|_{y=0}
\leq C\|u\|_{\ddot{H}^{|\alpha|+1}(\hat{\Omega})}^{2},
\end{align*}
and for $k\geq1$, by using \eqref{ube} one has
\begin{align*}
\|\partial_{t}^{k}[D^{\alpha}F]\|_{L_{x}^{2}(\mathbb{T})}|_{y=0}
\leq&~C_{*}\|u\|_{\ddot{H}^{|\alpha|+2k+1}(\hat{\Omega})}^{2}
     +C_{*}\sum_{i=1}^{k-1}\|\partial_{t}^{i}F\|_{\mathring{H}^{|\alpha|+2k-2i-2}(\mathbb{T})}^{2}
     +C_{*}\bar{\kappa}^{2},
\end{align*}
which gives
\begin{align*}
\|\hat{F}\|_{\mathring{H}^{s-1}(\mathbb{T})}^{2}\leq C\|u\|_{\ddot{H}^{s}(\hat{\Omega})}^{4}\leq C\bar{\eta}^{2},
\end{align*}
and,
\begin{align*}
\sum_{i=1}^{\frac{s-1}{2}}\|\partial_{t}^{i}\hat{F}\|_{\mathring{H}^{s-1-2i}(\mathbb{T})}^{2}
\leq&~C_{*}\|u\|_{\ddot{H}^{s}(\hat{\Omega})}^{4}
     +C_{*}\sum_{i=1}^{\frac{s-1}{2}}\sum_{j=0}^{i-1}\|\partial_{t}^{j}F\|_{\mathring{H}^{s-2j-3}(\mathbb{T})}^{4}
     +C_{*}\bar{\kappa}^{4}\\
\leq&~C_{*}\|u\|_{\ddot{H}^{s}(\hat{\Omega})}^{4}
     +C_{*}\sum_{i=0}^{\frac{s-3}{2}}\|\partial_{t}^{i}F\|_{\mathring{H}^{s-2i-3}(\mathbb{T})}^{4}
     +C_{*}\bar{\kappa}^{4}\\
\leq&~C_{*}\bar{\eta}^{2}.
\end{align*}
Thus
\begin{align}\label{Fbe}
\sum_{i=0}^{\frac{s-1}{2}}\|\partial_{t}^{i}\hat{F}\|_{L_{t}^{\infty}\mathring{H}^{s-1-2i}([0,t_{*}]\times\mathbb{T})}^{2}
\leq C_{*}\bar{\eta}^{2}.
\end{align}
Since $u|_{y=0}=0$, one has $\partial_{x}^{s}\hat{F}|_{y=0}=0$, thus
\begin{align}\label{Fbs+1e}
\|\partial_{x}^{s}\hat{F}\|_{L_{t}^{\infty}\mathring{H}^{0}([0,t_{*}]\times\mathbb{T})}^{2}=0.
\end{align}
A combination of \eqref{Fse}, \eqref{FL2Hs} \eqref{Fbe} and \eqref{Fbs+1e} gives \eqref{Fe6}.

\end{proof}

\section{Proof of the main result}~

In this section, with the help of the well-posedness result of the corrected increment problems obtained in Section 4 and the a posteriori estimates of the approximate solution obtained in Section 5, we shall construct the solution to the problem \eqref{Prandtl} based on the iterative scheme introduced in Section 3.

Let $\eta_{0}=1$ and $\bar{\varsigma}_{0}=\iota$. For any $n\in\mathbb{N}_{+}$, in the problem (AE$_{n}$) and \eqref{An} we set
$\bar{\eta}_{n}=\epsilon_{0}^{2n}, \bar{\varsigma}_{n}=\iota\bar{\eta}_{n}^{2}$,
and correspondingly, $\eta_{n}=\bar{\eta}_{n}\vartheta, \varsigma_{n}=\bar{\varsigma}_{n}\chi_{c}$, moreover, $\kappa_{n}=\eta_{n}-\eta_{n-1}$.
For any $(t,y)\in [0,T]\times[0,+\infty)$ and $n\in\mathbb{N}_{+}$, let $\varphi_{n}=\y^{-\frac{1}{2}}\psi\sqrt{\zeta_{n}}$, with
\begin{align*}
\zeta_{n}(t,y)
=&\chi(y)\left[\frac{e^{-\frac{16\lambda^{4}(\bar{\varsigma}_{n}+t)}{(y-y_{*})^{2}+16\lambda^{2}(\bar{\varsigma}_{n}+t)}}}{[\frac{\lambda^{-2}}{16}(y-y_{*})^{2}+\bar{\varsigma}_{n}+t]^{\varepsilon_{0}}}
  +\left(\frac{\lambda}{\varepsilon_{0}}\right)^{\lambda}e^{-\lambda^{2}}H_{\lambda}\left(\frac{y-y_{*}}{\sqrt{2(\bar{\varsigma}_{n}+t)}}+\frac{\sqrt{2}\varepsilon_{1}}{(\bar{\varsigma}_{n}+t)^{\varepsilon_{1}}}\right)\chi_{\lambda}\left(\frac{y-y_{*}}{\sqrt{\bar{\varsigma}_{n}+t}}\right)\right]\\
 &+1-\chi(y),
\end{align*}
and let
$$\rho_{n}=\frac{\partial_{y}\dot{\mathfrak{U}}_{n-1}+\varsigma_{n}'}{\dot{\mathfrak{U}}_{n-1}+\varsigma_{n}}
~\text{and}~
\mathcal{P}_{n}=\frac{\chi(y)}{\sqrt{\dot{\mathfrak{U}}_{n-1}+\varsigma_{n}}}+1-\chi(y).$$

For any $n\in\mathbb{N}$, let $F_{n}$ and $\mathfrak{U}_{n}$ be defined as in Section 3, and $\mathbb{U}_{n}=\mathfrak{U}_{n}-u_{0}$. In order for the iteration to proceed, we need that the following iteration assumptions hold for any $k\in\mathbb{N}_{+}$ and some $\hat{t}_{*}\in(0,T]$ independent of $k$.

\begin{assumption}[Iterative Assumption]\label{IA}
For some $\hat{t}_{*}\in(0,T]$ and a constant $k\in\mathbb{N}_{+}$, assume the following inequalities hold.

\noindent{\underline{Assumptions on $F_{k-1}$.}}
\begin{align}\label{ITA0}
|D^{\gamma}F_{k-1}|\leq C_{*}e^{-\frac{4}{3}\delta\y}~\text{in}~\Omega_{\hat{t}_{*}}, \forall \gamma\in\Gamma_{9},
\end{align}
\begin{align}\label{ITA1}
|\partial_{y}F_{k-1}|\leq C_{*}\iota^{-1}(\dot{\mathfrak{U}}_{k-1}+\bar{\varsigma}_{k})~\text{in}~\Omega_{\hat{t}_{*}}~\text{and}~
\|\mathcal{P}_{k}(\partial_{y}^{2}F_{k-1}-\rho_{k}\partial_{y}F_{k-1})\|_{L^{\infty}(\Omega_{\hat{t}_{*}})}\leq C_{*}\iota^{-1},
\end{align}
\begin{align}\label{ITA2}
\|\partial_{t}F_{k-1}\|_{L_{t}^{\infty}H_{\psi}^{s-4}(\Omega_{\hat{t}_{*}})}+\|F_{k-1}\|_{L_{t}^{\infty}H_{\psi}^{s-2}(\Omega_{\hat{t}_{*}})}
\leq C_{*}\bar{\eta}_{k-1},
\end{align}
\begin{align}\label{ITA3}
&\|\mathbf{F}_{k-1,s}\|_{L_{t}^{2}L_{\hat{\psi}}^{2}(\Omega_{\hat{t}_{*}})}^{2}
+\|F_{k-1}\|_{L_{t}^{2}\hat{H}_{\psi}^{s}(\Omega_{\hat{t}_{*}})}^{2}
+\sum_{i=0}^{\frac{s-1}{2}}\|\partial_{t}^{i}F_{k-1}\|_{L_{t}^{\infty}\mathring{H}^{s-1-2i}([0,\hat{t}_{*}]\times\mathbb{T})}^{2}\\
&+\|\partial_{x}^{s}F_{k-1}\|_{L_{t}^{\infty}\mathring{H}^{0}([0,\hat{t}_{*}]\times\mathbb{T})}^{2}
\leq C_{*}\bar{\eta}_{k},\nonumber
\end{align}
where $\mathbf{F}_{k-1,s}=\mathcal{P}_{k}\left(\partial_{y}\partial_{x}^{s}F_{k-1}-\rho_{k}\partial_{x}^{s}F_{k-1}\right)$.

\noindent{\underline{Assumptions on $\mathfrak{U}_{k-1}$.}}
\begin{align}\label{ITA4}
\dot{\mathfrak{U}}_{k-1}\geq C\varpi(U-\mathfrak{U}_{k-1})\geq C\varpi e^{-(1+\theta t)\delta\y}~\text{in}~\Omega_{\hat{t}_{*}},
\end{align}
\begin{align}\label{ITA5}
|\partial_{x}\dot{\mathfrak{U}}_{k-1}|\leq C(1+t\y)\dot{\mathfrak{U}}_{k-1}, |\partial_{x}^{2}\dot{\mathfrak{U}}_{k-1}|\leq C(1+t\y)\dot{\mathfrak{U}}_{k-1}~\text{in}~\Omega_{\hat{t}_{*}},
\end{align}
\begin{align}\label{ITA6}
\dot{\mathfrak{U}}_{k-1}\geq C(\varpi+t)~\text{in}~\Omega_{\hat{t}_{*}}^{*}.
\end{align}

\noindent{\underline{Assumptions on $\mathbb{U}_{k-1}$.}}
\begin{align}\label{ITA7}
|\vartheta_{c}^{2}D^{\gamma}\mathbb{U}_{k-1}|\leq C_{*}\y^{-1}\dot{\mathfrak{U}}_{k-1}~\text{in}~\check{\Omega}_{\hat{t}_{*}}^{\hat{y}}, \forall \gamma\in\Gamma_{6}.
\end{align}
\begin{align}\label{ITA8}
& \|\mathbb{U}_{k-1}\|_{L_{t}^{\infty}\mathcal{H}_{\psi,\varphi_{k}}^{s}(\Omega_{\hat{t}_{*}})}
 +\|\partial_{y}\mathbb{U}_{k-1}\|_{L_{t}^{2}\hat{H}_{\psi}^{s}(\Omega_{\hat{t}_{*}})}
 +\|\partial_{y}\mathbf{U}_{k-1,s}\|_{L_{t}^{2}L_{\varphi_{k}}^{2}(\Omega_{\hat{t}_{*}})}
 +\|\mathcal{P}_{k}\mathbf{U}_{k-1,s}\|_{L_{t}^{2}L_{\varphi_{k}}^{2}(\Omega_{\hat{t}_{*}})}\\
&+\|\mathbf{U}_{k-1,s+1}\|_{L_{t}^{2}L_{\hat{\varphi}_{k}}^{2}(\Omega_{\hat{t}_{*}})}
 +\|\mathcal{P}_{k}\mathbf{U}_{k-1,s+1}\|_{L_{t}^{2}L_{\hat{\varphi}_{k}}^{2}(\Omega_{\hat{t}_{*}})}
 +\|\sqrt{\eta_{k-1}}\partial_{x}\mathbf{U}_{k-1,s+1}\|_{L_{t}^{2}L_{\hat{\varphi}_{k}}^{2}(\Omega_{\hat{t}_{*}})}
\leq \mathcal{Z}.\nonumber
\end{align}

\noindent{\underline{Assumptions on $\rho_{k}$ and $\mathcal{P}_{k}$.}}
\begin{align}\label{ITA9}
\rho_{k}\leq -\frac{15}{16}\delta,~\text{in}~\check{\Omega}_{\hat{t}_{*}}^{Y},
\end{align}
\begin{align}\label{ITA10}
|\rho_{k}|\leq C\mathcal{P}_{k},~|\partial_{x}\rho_{k}|\leq C\mathcal{P}_{k},~
|\partial_{y}\rho_{k}|\leq C\mathcal{P}_{k}^{2},~\\
\left|\partial_{x}(\mathcal{P}_{k}\partial_{x}\rho_{k})\right|\leq C(1+t\y)^{2}\mathcal{P}_{k}^{2},~\text{in}~\Omega_{\hat{t}_{*}}.\nonumber
\end{align}
Here the constants $\mathcal{Z}, C_{*}$ and $C$ are independent of $k$.
\end{assumption}

From Theorem \ref{wpae} and Theorem \ref{Phe}, we can prove the following theorem.
\begin{theorem}\label{MRp}
Assume that the iterative assumption Assumption \ref{IA} holds for $n=1$ and $\hat{t}_{*}=T$, and
\begin{align}\label{ITA11}
 \|\tilde{u}_{in}\|_{\hat{H}_{\psi_{in}}^{s}(\Omega)}
+\|\mathcal{P}_{1}(0,x,y)[\partial_{y}\partial_{x}^{s}\tilde{u}_{in}-\rho_{1}(0,x,y)\partial_{x}^{s}\tilde{u}_{in}]\|_{L_{\varphi_{1,in}}^{2}(\Omega)}
\leq \frac{(1-\epsilon_{0})\sqrt{\bar{\eta}_{1}}\mathcal{Z}}{8\epsilon_{0}},
\end{align}
hold for some positive constant $\mathcal{Z}$, where $\varphi_{1,in}=\varphi_{1}(0,y)$, then under the assumption of Theorem \ref{MR}, there exists constants $t_{*}\in(0,T]$ and $Z>0$ such that problem \eqref{Prandtl} admits an unique classical solution $u$ such that $w:=u-u_{0}\in L_{t}^{\infty}\hat{H}_{\psi}^{s}(\Omega_{t_{*}})$ satisfies the following estimate,
\begin{align}\label{ME}
& \|w\|_{L_{t}^{\infty}\hat{H}_{\psi}^{s}(\Omega_{t_{*}})}
 +\|\partial_{y}w\|_{L_{t}^{2}\hat{H}_{\psi}^{s}(\Omega_{t_{*}})}
 +\|\partial_{x}^{s}w\|_{L_{t}^{2}L_{\hat{\psi}}^{2}(\Omega_{t_{*}})}
 +\|\partial_{y}\partial_{x}^{s}w\|_{L_{t}^{2}L_{\hat{\psi}}^{2}(\Omega_{t_{*}})}
 +\|\partial_{y}^{2}\partial_{x}^{s}w\|_{L_{t}^{2}L_{\hat{\psi}}^{2}(\Omega_{t_{*}})}\\
&+\|\partial_{x}^{s+1}w\|_{L_{t}^{2}L_{\check{\psi}}^{2}(\Omega_{t_{*}})}
 +\|\partial_{y}\partial_{x}^{s+1}w\|_{L_{t}^{2}L_{\check{\psi}}^{2}(\Omega_{t_{*}})}
 +\|\partial_{y}^{2}\partial_{x}^{s+1}w\|_{L_{t}^{2}L_{\check{\psi}}^{2}(\Omega_{t_{*}})}
\leq Z.\nonumber
\end{align}
\end{theorem}
\begin{proof}[Proof of Theorem \ref{MRp}]
In fact, by applying Theorem \ref{wpae} and Theorem \ref{Phe} to the problem \eqref{AEn+1} and (A$_{n+1}$) respectively, one can obtain the following lemma.
\begin{lemma}\label{IT}
Under the assumption of Theorem \ref{MR}, there exists uniform positive constants $\lambda, \Lambda, \theta, \iota, \mathcal{Z}$ and $t_{*}$ such that for any $n\in\mathbb{N}_{+}$, \\
(i) Assume the iterative assumption Assumption \ref{IA} holds for $k=n$, and \eqref{ITA11} holds in addition, then the problem \textsc{(AE$_{n}$)} admits an unique solution
$u_{n}\in L_{t}^{\infty}\mathcal{H}_{\psi,\varphi_{n}}^{s}(\Omega_{t_{*}})$ satisfying
\begin{align}\label{ITC1}
&~\|u_{n}\|_{L_{t}^{\infty}\mathcal{H}_{\psi,\varphi_{n}}^{s}(\Omega_{t_{*}})}
 +\|u_{n}\|_{L_{t}^{2}\hat{H}_{\psi,\y}^{s}(\Omega_{t_{*}})}
 +\|\partial_{y}u_{n}\|_{L_{t}^{2}\hat{H}_{\psi}^{s}(\Omega_{t_{*}})}\\
&+\|\partial_{x}u_{n}\|_{L_{t}^{2}\hat{H}_{\psi,\eta}^{s}(\Omega_{t_{*}})}
 +\|\mathcal{P}_{n}\mathcal{U}_{n,s}\|_{L_{t}^{2}L_{\varphi_{n}}^{2}(\Omega_{t_{*}})}
 +\|\sqrt{\y}\mathcal{U}_{n,s}\|_{L_{t}^{2}L_{\varphi_{n}}^{2}(\Omega_{t_{*}})}\nonumber\\
&+\|\partial_{y}\mathcal{U}_{n,s}\|_{L_{t}^{2}L_{\varphi_{n}}^{2}(\Omega_{t_{*}})}
 +\|\sqrt{\eta_{n}}\partial_{x}\mathcal{U}_{n,s}\|_{L_{t}^{2}L_{\varphi_{n}}^{2}(\Omega_{t_{*}})}
\leq \frac{1-\epsilon_{0}}{\epsilon_{0}}\sqrt{\bar{\eta}_{n}}\mathcal{Z}.\nonumber
\end{align}
(ii) $\mathfrak{U}_{n}=\mathfrak{U}_{n-1}+u_{n}=\mathbb{U}_{n}+u_{0}$ is the solution to the problem \eqref{An} on the time interval $[0,t_{*}]$. If
\begin{align}\label{ITA12}
~&\|\mathbb{U}_{n}\|_{L_{t}^{\infty}\hat{H}_{\psi}^{s}(\Omega_{t_{*}})}
+\|\mathbb{U}_{n}\|_{L_{t}^{2}\hat{H}_{\psi,\y}^{s}(\Omega_{t_{*}})}
+\|\mathbf{U}_{n-1,s}\|_{L_{t}^{\infty}L_{\varphi_{n+1}}^{2}(\Omega_{t_{*}})}\\
&+\|\mathbf{U}_{n-1,s+1}\|_{L_{t}^{\infty}L_{\hat{\varphi}_{n}}^{2}(\Omega_{t_{*}})}
+\|\sqrt{\eta_{n}}\partial_{x}\mathbf{U}_{n-1,s}\|_{L_{t}^{2}L_{\varphi_{n}}^{2}(\Omega_{t_{*}})}
+\|\mathcal{P}_{n}\mathbf{U}_{n-1,s}\|_{L_{t}^{2}L_{\varphi_{n}}^{2}(\Omega_{t_{*}})}\nonumber\\
&+\|\sqrt{\eta_{n-1}}\partial_{x}\mathbf{U}_{n-1,s+1}\|_{L_{t}^{2}L_{\hat{\varphi}_{n}}^{2}(\Omega_{t_{*}})}
+\|\mathcal{P}_{n}\mathbf{U}_{n-1,s+1}\|_{L_{t}^{2}L_{\hat{\varphi}_{n}}^{2}(\Omega_{t_{*}})}
\leq \mathcal{Z},\nonumber
\end{align}
then,
\begin{align}\label{ITC2}
&\|\mathbf{U}_{n,s}\|_{L_{t}^{\infty}L_{\varphi_{n+1}}^{2}(\Omega_{t_{*}})}
 +\|\sqrt{\y}\mathbf{U}_{n,s}\|_{L_{t}^{2}L_{\varphi_{n+1}}^{2}(\Omega_{t_{*}})}
 +\|\mathcal{P}_{n+1}\mathbf{U}_{n,s}\|_{L_{t}^{2}L_{\varphi_{n+1}}^{2}(\Omega_{t_{*}})}\\
&+\|\partial_{y}\mathbf{U}_{n,s}\|_{L_{t}^{2}L_{\varphi_{n+1}}^{2}(\Omega_{t_{*}})}
 +\|\sqrt{\eta_{n}}\partial_{x}\mathbf{U}_{n,s}\|_{L_{t}^{2}L_{\varphi_{n+1}}^{2}(\Omega_{t_{*}})}
\leq \mathcal{Z},\nonumber
\end{align}
and
\begin{align}\label{ITC3}
&\|\mathbf{U}_{n,s+1}\|_{L_{t}^{\infty}L_{\hat{\varphi}_{n+1}}^{2}(\Omega_{t_{*}})}
 +\|\sqrt{\y}\mathbf{U}_{n,s+1}\|_{L_{t}^{2}L_{\hat{\varphi}_{n+1}}^{2}(\Omega_{t_{*}})}
 +\|\mathcal{P}_{n+1}\mathbf{U}_{n,s+1}\|_{L_{t}^{2}L_{\hat{\varphi}_{n+1}}^{2}(\Omega_{t_{*}})}\\
&+\|\partial_{y}\mathbf{U}_{n,s+1}\|_{L_{t}^{2}L_{\hat{\varphi}_{n+1}}^{2}(\Omega_{t_{*}})}
 +\|\sqrt{\eta_{n}}\partial_{x}\mathbf{U}_{n,s+1}\|_{L_{t}^{2}L_{\hat{\varphi}_{n+1}}^{2}(\Omega_{t_{*}})}
\leq \mathcal{Z}.\nonumber
\end{align}
moreover, the iterative assumption Assumption \ref{IA} holds for $k=n+1$ on the time interval $[0,t_{*}]$.
\end{lemma}

It is worth noting that all the constants $C, C_{*}, \lambda,\Lambda,\theta,\iota,t_{*}$ and $\mathcal{Z}$ appearing in Lemma \ref{IT} are independent of $n$, thus we can use Lemma \ref{IT} to construct the solution to the problem \eqref{Prandtl}.

Under the assumptions of Theorem \ref{MRp}, it is easy to check that the assumption of Lemma \ref{IT}(i) are satisfied for $n=1$.
On the other hand, since $\mathbb{U}_{1}=u_{1}$, from \eqref{ITC1} one has
\begin{align}
 \|\mathbb{U}_{1}\|_{L_{t}^{\infty}\mathcal{H}_{\psi,\varphi_{1}}^{s}(\Omega_{t_{*}})}
+\|\mathbb{U}_{1}\|_{L_{t}^{2}\hat{H}_{\psi,\y}^{s}(\Omega_{t_{*}})}
\leq \frac{1-\epsilon_{0}}{\epsilon_{0}}\sqrt{\bar{\eta}_{1}}\mathcal{Z}
\leq \mathcal{Z},
\end{align}
which gives \eqref{ITA12} for $n=1$, due to $\mathbb{U}_{0}=0$. Thus the assumption of Lemma \ref{IT}(ii) are satisfied for $n=1$.

Now that the assumptions of Lemma \ref{IT} are satisfied for $n=1$, we can use Lemma \ref{IT} repeatedly to gain a sequence of corrected increment solutions $\{u_{n}\}_{n=1}^{\infty}$, where $u_{n}$ is the unique solution to the corrected increment problem (AE$_{n}$) satisfying \eqref{ITC1}. Accordingly, we obtain a sequence of approximate solutions $\{\mathfrak{U}_{n}\}_{n=1}^{\infty}$ to the problem \eqref{Prandtl}, with $\mathfrak{U}_{n}=\sum\limits_{i=0}^{n}u_{i}$ being the solution to the approximate problem \eqref{An}.
The estimate \eqref{ITC1} tells that
\begin{align}\label{ITC-1}
\sum_{n=1}^{\infty}\|u_{n}\|_{L_{t}^{\infty}\hat{H}_{\psi}^{s}(\Omega_{t_{*}})}
+\sum_{n=1}^{\infty}\|\partial_{y}u_{n}\|_{L_{t}^{2}\hat{H}_{\psi}^{s}(\Omega_{t_{*}})}
\leq\sum_{n=1}^{\infty}\frac{1-\epsilon_{0}}{\epsilon_{0}}\epsilon_{0}^{n}\mathcal{Z}=\mathcal{Z}.
\end{align}
Thus $\sum\limits_{n=1}^{\infty}u_{n}$ converges in the space $L_{t}^{\infty}\hat{H}_{\psi}^{s}(\Omega_{t_{*}})$. Let $u=u_{0}+\sum\limits_{n=1}^{\infty}u_{n}$,
that is $u=\lim\limits_{n\to\infty}\mathfrak{U}_{n}$ in $L_{t}^{\infty}\hat{H}_{\psi}^{s}(\Omega_{t_{*}})$. By taking $k=n$ in \eqref{ITA2} one knows that $\lim\limits_{n\to\infty}\|F_{n}\|_{L_{t}^{\infty}H_{\psi}^{s-2}(\Omega_{t_{*}})}=0$, thus $u$ is a solution to the problem \eqref{Prandtl}.

Notice that $\mathbb{U}_{n}=\mathfrak{U}_{n}-u_{0}$, one knows that $w:=u-u_{0}=\lim\limits_{n\to\infty}\mathbb{U}_{n}$.
The estimates \eqref{ITC2} and \eqref{ITC3} gives us that for any $n\in\mathbb{N}_{+}$,
\begin{align}\label{ITC-2}
& \|\mathbf{U}_{n,s}\|_{L_{t}^{\infty}L_{\hat{\psi}}^{2}(\Omega_{t_{*}})}
 +\|\mathcal{P}_{n+1}\mathbf{U}_{n,s}\|_{L_{t}^{2}L_{\hat{\psi}}^{2}(\Omega_{t_{*}})}
 +\|\partial_{y}\mathbf{U}_{n,s}\|_{L_{t}^{2}L_{\hat{\psi}}^{2}(\Omega_{t_{*}})}\\
&+\|\mathbf{U}_{n,s+1}\|_{L_{t}^{\infty}L_{\check{\psi}}^{2}(\Omega_{t_{*}})}
 +\|\mathcal{P}_{n+1}\mathbf{U}_{n,s+1}\|_{L_{t}^{2}L_{\check{\psi}}^{2}(\Omega_{t_{*}})}
 +\|\partial_{y}\mathbf{U}_{n,s+1}\|_{L_{t}^{2}L_{\check{\psi}}^{2}(\Omega_{t_{*}})}
\leq C_{*},\nonumber
\end{align}
With the help of Corollary \ref{ncC2}, from \eqref{ITC-2} one can deduce that for any $n\in\mathbb{N}_{+}$,
\begin{align}\label{ITC-3}
& \|\partial_{x}^{s}\mathbb{U}_{n}\|_{L_{t}^{2}L_{\hat{\psi}}^{2}(\Omega_{t_{*}})}
 +\|\partial_{y}\partial_{x}^{s}\mathbb{U}_{n}\|_{L_{t}^{2}L_{\hat{\psi}}^{2}(\Omega_{t_{*}})}
 +\|\partial_{y}^{2}\partial_{x}^{s}\mathbb{U}_{n}\|_{L_{t}^{2}L_{\hat{\psi}}^{2}(\Omega_{t_{*}})}
 +\|\partial_{x}^{s+1}\mathbb{U}_{n}\|_{L_{t}^{2}L_{\check{\psi}}^{2}(\Omega_{t_{*}})}\\
&+\|\partial_{y}\partial_{x}^{s+1}\mathbb{U}_{n}\|_{L_{t}^{2}L_{\check{\psi}}^{2}(\Omega_{t_{*}})}
 +\|\partial_{y}^{2}\partial_{x}^{s+1}\mathbb{U}_{n}\|_{L_{t}^{2}L_{\check{\psi}}^{2}(\Omega_{t_{*}})}
\leq C_{*}.\nonumber
\end{align}
Since $\mathbb{U}_{n}$ converges to $w$ in the space $\hat{H}_{\psi}^{s}(\Omega)$, \eqref{ITC-3} together with \eqref{ITC-1} gives \eqref{ME}.

\end{proof}

Now we are in position to prove our main result, by using Theorem \ref{MRp}.
\begin{proof}[Proof of Theorem \ref{MR}]
To prove Theorem \ref{MR}, we need to check that $u_{0}, F_{0}$ and $\tilde{u}_{in}$ satisfy the assumptions in the Theorem \ref{IT}, that is, \eqref{ITA0}--\eqref{ITA11} hold for $n=1$.
In fact, noticing that $\mathfrak{U}_{0}=u_{0}$ and $\mathbb{U}_{0}=0$, by the definition of $u_{0}$ one can check that \eqref{ITA4}--\eqref{ITA10} hold for $n=1$.
Recall the definition of $F_{0}$ and $u_{0}$, \eqref{ITA0} holds clearly. Under the assumption of Theorem \ref{MR} it is easy to obtain
\begin{align}\label{ITA-1}
& \|F_{0}\|_{L_{t}^{2}\hat{H}_{\psi}^{s}(\Omega_{t_{*}})}^{2}
 +\|\partial_{t}F_{0}\|_{L_{t}^{\infty}H_{\psi}^{s-4}(\Omega_{t_{*}})}+\|F_{0}\|_{L_{t}^{\infty}H_{\psi}^{s-2}(\Omega_{t_{*}})}\\
&+\sum_{i=0}^{\frac{s-1}{2}}\|\partial_{t}^{i}F_{0}\|_{L_{t}^{\infty}\mathring{H}^{s-1-2i}(\mathbb{T}\times[0,t_{*}])}^{2}
 +\|\partial_{x}^{s}F_{0}\|_{L_{t}^{\infty}\mathring{H}^{0}(\mathbb{T}\times[0,t_{*}])}^{2}
\leq C_{*}=C_{*}\bar{\eta}_{0},\nonumber
\end{align}
moreover, since $\bar{\varsigma}_{0}=\iota$, one has
\begin{align}\label{ITA-2}
|\partial_{y}F_{0}|\leq C_{*}\leq C_{*}\iota^{-1}(\dot{\mathfrak{U}}_{0}+\bar{\varsigma}_{0})~\text{in}~\Omega_{t_{*}}.
\end{align}
In view of $\|\mathcal{P}_{1}\|_{L^{\infty}(\Omega_{t_{*}})}=\iota^{-\frac{1}{2}}\epsilon_{0}^{-2}$, one has
\begin{align}\label{ITA-3}
\|\mathcal{P}_{1}(\partial_{y}^{2}F_{0}-\rho_{1}\partial_{y}F_{0})\|_{L^{\infty}(\Omega_{t_{*}})}
\leq \|\mathcal{P}_{1}\|_{L^{\infty}(\Omega_{t_{*}})}^{2}
     \left(\|\partial_{y}^{2}F_{0}\|_{L^{\infty}(\Omega_{t_{*}})}+\|\partial_{y}F_{0}\|_{L^{\infty}(\Omega_{t_{*}})}\right)
\leq C_{*}\iota^{-1}.
\end{align}
Notice that $\partial_{y}\partial_{x}^{s}F_{0}=\partial_{x}^{s}F_{0}=0$ in $\hat{\Omega}_{t_{*}}$, one can get
\begin{align}\label{ITA-4}
\|\mathbf{F}_{0,s}\|_{L_{t}^{2}L_{\hat{\psi}}^{2}(\Omega_{t_{*}})}^{2}
=\|\mathcal{P}_{1}\left(\partial_{y}\partial_{x}^{s}F_{0}-\rho_{1}\partial_{x}^{s}F_{0}\right)\|_{L_{t}^{2}L_{\hat{\psi}}^{2}(\check{\Omega}_{t_{*}}^{\hat{y}})}^{2}
\leq C_{*}\bar{\eta}_{0}.
\end{align}
A combination of \eqref{ITA-1}--\eqref{ITA-4} gives \eqref{ITA1}--\eqref{ITA3} for $n=1$. By the definition of $u_{0}$, one can check that $|\rho_{1}(0,x,y)|\leq C\mathcal{P}_{1}(0,x,y)\leq C\mathcal{P}_{in}$ in $\Omega$, thus from \eqref{inA4} one can take $\mathcal{Z}$ large enough to obtain
\begin{align}
&\|\tilde{u}_{in}\|_{\hat{H}_{\psi_{in}}^{s}(\Omega)}
+\|\mathcal{P}_{1}(0,x,y)[\partial_{y}\partial_{x}^{s}\tilde{u}_{in}-\rho_{1}(0,x,y)\partial_{x}^{s}\tilde{u}_{in}]\|_{L_{\varphi_{in}}^{2}(\Omega)}\\
\leq&~\|\tilde{u}_{in}\|_{\hat{H}_{\psi_{in}}^{s}(\Omega)}
     +C_{*}\|\omega_{\varepsilon_{0}}^{in}\mathcal{P}_{in}\partial_{y}\partial_{x}^{s}\tilde{u}_{in}\|_{L_{\hat{\psi}_{in}}^{2}(\Omega)}
     +C_{*}\|\omega_{\varepsilon_{0}}^{in}\mathcal{P}_{in}^{2}\partial_{x}^{s}\tilde{u}_{in}\|_{L_{\hat{\psi}_{in}}^{2}(\Omega)}\nonumber\\
\leq&~\frac{(1-\epsilon_{0})\sqrt{\bar{\eta}_{1}}\mathcal{Z}}{8\epsilon_{0}}.\nonumber
\end{align}

The uniqueness of the solution to the problem \eqref{Prandtl} follows by the $L^{2}$ comparison, we omit the detail.
\end{proof}

\appendix

\section{Proof of the trace estimates}

This section is devoted to the  proof of Lemma \ref{Bc} and Lemma \ref{bc}.
\begin{proof}[Proof of Lemma \ref{Bc}]
The inequality \eqref{Bc1} will be proved by induction.

Let $t\in [0,T_{c}]$. Let $\alpha=(\alpha_{1},\alpha_{2})$ be a multiindex such that $\alpha\in\Gamma_{s-2k+1}$.
Set $D^{\alpha}=\partial_{x}^{\alpha_{1}}\partial_{y}^{\alpha_{2}}$.
Set
$$\Delta_{\tilde{\eta}}:=\partial_{y}^{2}+\tilde{\eta}\partial_{x}^{2}.$$
By the equation \eqref{App}$_{1}$ and noticing that $\tilde{\eta}$ take a constant value in $\hat{\Omega}$, one knows that
\begin{align*}
\partial_{t}[D^{\alpha}\mathfrak{U}]
= \Delta_{\tilde{\eta}}[D^{\alpha}\mathfrak{U}]
 -D^{\alpha}[\mathfrak{U}\partial_{x}\mathfrak{U}-\partial_{y}^{-1}[\partial_{x}\mathfrak{U}]\partial_{y}\mathfrak{U}]
 -D^{\alpha}[\partial_{x}P]
 -D^{\alpha}\mathfrak{F},~\text{in}~\hat{\Omega},
\end{align*}
where $\mathfrak{F}=\tilde{\eta}\partial_{x}^{2}u_{0}+F$.
If $\alpha_{2}=0$ we have $\partial_{t}^{k}[D^{\alpha}\mathfrak{U}]|_{y=0}=0$ since $\mathfrak{U}|_{y=0}=0$. In the following we assume $\alpha_{2}>0$, then
\begin{align}\label{Dalphae}
\partial_{t}[D^{\alpha}\mathfrak{U}]|_{y=0}
= \Delta_{\tilde{\eta}}[D^{\alpha}\mathfrak{U}]|_{y=0}
 -D^{\alpha}[\mathfrak{U}\partial_{x}\mathfrak{U}-\partial_{y}^{-1}[\partial_{x}\mathfrak{U}]\partial_{y}\mathfrak{U}]|_{y=0}
 -D^{\alpha}\mathfrak{F}|_{y=0}.
\end{align}
Noticing $\mathfrak{U}|_{y=0}=\partial_{x}\mathfrak{U}|_{y=0}=0$, by trace theorem we know that for any $\alpha\in\Gamma_{s-2k+1}$, it follows
\begin{align*}
\|\partial_{t}[D^{\alpha}\mathfrak{U}]\|_{L_{x}^{2}(\mathbb{T})}|_{y=0}
\leq \|\mathfrak{U}\|_{\ddot{H}^{|\alpha|+3}(\hat{\Omega})}
    +\|\mathfrak{U}\|_{\ddot{H}^{|\alpha|+1}(\hat{\Omega})}^{2}
    +\|\mathfrak{F}\|_{\mathring{H}^{|\alpha|}(\mathbb{T})}.
\end{align*}
Assume that for any integer $m\in [1,k-1]$ we have
\begin{align}\label{indA}
\|\partial_{t}^{m}[D^{\alpha}\mathfrak{U}]\|_{L_{x}^{2}(\mathbb{T})}|_{y=0}
\leq& C\sum_{i=0}^{m}\|\mathfrak{U}\|_{\ddot{H}^{|\alpha|+2i+1}(\hat{\Omega})}^{m-i+1}
     +C\sum_{\substack{p_{ij}>0,\\ i+j\leq m-1,\\0<i+j+\sum\limits_{i,j} p_{ij}+\sum\limits_{j} q_{j}\leq m}}\prod_{i=0}^{m-1}
      \|\partial_{t}^{i}\mathfrak{F}\|_{\mathring{H}^{|\alpha|+2j}(\mathbb{T})}^{p_{ij}}
      \|\mathfrak{U}\|_{\ddot{H}^{|\alpha|+2j+1}(\hat{\Omega})}^{q_{j}}.
\end{align}
It follows from \eqref{Dalphae} that
\begin{align}\label{Dkalphae}
\partial_{t}^{k}[D^{\alpha}\mathfrak{U}]|_{y=0}
=&\partial_{t}^{k-1}\partial_{t}[D^{\alpha}\mathfrak{U}]|_{y=0}\\
=& \partial_{t}^{k-1}\Delta_{\tilde{\eta}}[D^{\alpha}\mathfrak{U}]|_{y=0}
  -\partial_{t}^{k-1}\left[D^{\alpha}[\mathfrak{U}\partial_{x}\mathfrak{U}-\partial_{y}^{-1}[\partial_{x}\mathfrak{U}]\partial_{y}\mathfrak{U}]\right]|_{y=0}\nonumber\\
 &-\partial_{t}^{k-1}[D^{\alpha}\mathfrak{F}]|_{y=0}.\nonumber
\end{align}
By the inductive assumption \eqref{indA}, we have
\begin{align}\label{Dkalpha1}
&\|\partial_{t}^{k-1}\Delta_{\tilde{\eta}}[D^{\alpha}\mathfrak{U}]\|_{L_{x}^{2}(\mathbb{T})}|_{y=0}\\
\leq& C\sum_{i=0}^{k-1}\|\mathfrak{U}\|_{\ddot{H}^{|\alpha|+2i+3}(\hat{\Omega})}^{k-i}
     +C\sum_{\substack{p_{ij}>0,\\ i+j\leq k-2,\\0<i+j+p_{ij}+q_{j}\leq k-1}}\prod_{i=0}^{k-2}
      \|\partial_{t}^{i}\mathfrak{F}\|_{\mathring{H}^{|\alpha|+2j+2}(\mathbb{T})}^{p_{ij}}
      \|\mathfrak{U}\|_{\ddot{H}^{|\alpha|+2j+3}(\hat{\Omega})}^{q_{j}}\nonumber\\
\leq& C\sum_{i=1}^{k}\|\mathfrak{U}\|_{\ddot{H}^{|\alpha|+2i+1}(\hat{\Omega})}^{k-i+1}
     +C\sum_{\substack{p_{ij}>0,\\ i+j\leq k-1,\\0<i+j+p_{ij}+q_{j}\leq k}}\prod_{i=0}^{k-2}
      \|\partial_{t}^{i}\mathfrak{F}\|_{\mathring{H}^{|\alpha|+2j}(\mathbb{T})}^{p_{ij}}
      \|\mathfrak{U}\|_{\ddot{H}^{|\alpha|+2j+1}(\hat{\Omega})}^{q_{j}}\nonumber\\
\leq& C\sum_{i=0}^{k}\|\mathfrak{U}\|_{\ddot{H}^{|\alpha|+2i+1}(\hat{\Omega})}^{k-i+1}
     +C\sum_{\substack{p_{ij}>0,\\ i+j\leq k-1,\\0<i+j+p_{ij}+q_{j}\leq k}}\prod_{i=0}^{k-1}
      \|\partial_{t}^{i}\mathfrak{F}\|_{\mathring{H}^{|\alpha|+2j}(\mathbb{T})}^{p_{ij}}
      \|\mathfrak{U}\|_{\ddot{H}^{|\alpha|+2j+1}(\hat{\Omega})}^{q_{j}}.\nonumber
\end{align}
A direct calculation gives
\begin{align}\label{Dkalpha20}
&\partial_{t}^{k-1}\left[D^{\alpha}[\mathfrak{U}\partial_{x}\mathfrak{U}-\partial_{y}^{-1}[\partial_{x}\mathfrak{U}]\partial_{y}\mathfrak{U}]\right]|_{y=0}\\
=&\sum_{r=0}^{k-1}\sum_{0<\beta<\alpha}\binom{k-1}{r}\binom{\gamma}{\beta}\partial_{t}^{r}[D^{\beta}\mathfrak{U}]|_{y=0}\partial_{t}^{k-1-r}[D^{\alpha-\beta}[\partial_{x}\mathfrak{U}]]|_{y=0}\nonumber\\
 &-\sum_{r=0}^{k-1}\sum_{\substack{\beta\leq\alpha\\|\beta|>1}}\binom{k-1}{r}\binom{\gamma}{\beta}\partial_{t}^{r}[D^{\beta}\partial_{y}^{-1}[\partial_{x}\mathfrak{U}]]|_{y=0}\partial_{t}^{k-1-r}[D^{\alpha-\beta}[\partial_{y}\mathfrak{U}]]|_{y=0}.\nonumber
\end{align}
By utilizing inductive assumption \eqref{indA}, one has
\begin{align}\label{Dkalpha21}
\|\partial_{t}^{r}[D^{\beta}\mathfrak{U}]\|_{L_{x}^{2}(\mathbb{T})}|_{y=0}
\leq C\sum_{i=0}^{r}\|\mathfrak{U}\|_{\ddot{H}^{|\beta|+2i+1}(\hat{\Omega})}^{r-i+1}
     +C\sum_{\substack{p_{ij}>0,\\ i+j\leq r-1,\\0<i+j+\sum\limits_{i,j} p_{ij}+\sum\limits_{j} q_{j}\leq r}}\prod_{i=0}^{r-1}
      \|\partial_{t}^{i}\mathfrak{F}\|_{\mathring{H}^{|\beta|+2j}(\mathbb{T})}^{p_{ij}}
      \|\mathfrak{U}\|_{\ddot{H}^{|\beta|+2j+1}(\hat{\Omega})}^{q_{j}},
\end{align}
and
\begin{align}\label{Dkalpha22}
\partial_{t}^{k-1-r}[D^{\alpha-\beta}[\partial_{x}\mathfrak{U}]]|_{y=0}
\leq& C\sum_{i=0}^{k-1-r}\|\mathfrak{U}\|_{\ddot{H}^{|\alpha-\beta|+2i+2}(\hat{\Omega})}^{k-r-i}\\
    &+C\sum_{\substack{p_{ij}>0,\\ i+j\leq k-r-2,\\0<i+j+\sum\limits_{i,j} p_{ij}+\sum\limits_{j} q_{j}\leq k-1-r}}\prod_{i=0}^{k-r-2}
      \|\partial_{t}^{i}\mathfrak{F}\|_{\mathring{H}^{|\alpha-\beta|+2j+1}(\mathbb{T})}^{p_{ij}}
      \|\mathfrak{U}\|_{\ddot{H}^{|\alpha-\beta|+2j+2}(\hat{\Omega})}^{q_{j}}.\nonumber
\end{align}
Similarly,
\begin{align}\label{Dkalpha23}
\partial_{t}^{r}[D^{\beta}\partial_{y}^{-1}[\partial_{x}\mathfrak{U}]]|_{y=0}
\leq C\sum_{i=0}^{r}\|\mathfrak{U}\|_{\ddot{H}^{|\beta|+2i+1}(\hat{\Omega})}^{r-i+1}
     +C\sum_{\substack{p_{ij}>0,\\ i+j\leq r-1,\\0<i+j+\sum\limits_{i,j} p_{ij}+\sum\limits_{j} q_{j}\leq r}}\prod_{i=0}^{r-1}
      \|\partial_{t}^{i}\mathfrak{F}\|_{\mathring{H}^{|\beta|+2j}(\mathbb{T})}^{p_{ij}}
      \|\mathfrak{U}\|_{\ddot{H}^{|\beta|+2j+1}(\hat{\Omega})}^{q_{j}},
\end{align}
and
\begin{align}\label{Dkalpha24}
\partial_{t}^{k-1-r}[D^{\alpha-\beta}[\partial_{y}\mathfrak{U}]]|_{y=0}
\leq& C\sum_{i=0}^{k-1-r}\|\mathfrak{U}\|_{\ddot{H}^{|\alpha-\beta|+2i+2}(\hat{\Omega})}^{k-r-i}\\
    &+C\sum_{\substack{p_{ij}>0,\\ i+j\leq k-r-2,\\0<i+j+\sum\limits_{i,j} p_{ij}+\sum\limits_{j} q_{j}\leq k-1-r}}\prod_{i=0}^{k-r-2}
      \|\partial_{t}^{i}\mathfrak{F}\|_{\mathring{H}^{|\alpha-\beta|+2j+1}(\mathbb{T})}^{p_{ij}}
      \|\mathfrak{U}\|_{\ddot{H}^{|\alpha-\beta|+2j+2}(\hat{\Omega})}^{q_{j}}.\nonumber
\end{align}
Substituting the estimates \eqref{Dkalpha21}-\eqref{Dkalpha24} above into \eqref{Dkalpha20}, it follows
\begin{align}\label{Dkalpha2}
&\partial_{t}^{k-1}\left[D^{\alpha}[\mathfrak{U}\partial_{x}\mathfrak{U}-\partial_{y}^{-1}[\partial_{x}\mathfrak{U}]\partial_{y}\mathfrak{U}]\right]|_{y=0}\\
\leq& C\sum_{i=0}^{k}\|\mathfrak{U}\|_{\ddot{H}^{|\alpha|+2i+1}(\hat{\Omega})}^{k-i+1}
     +C\sum_{\substack{p_{ij}>0,\\ i+j\leq k-1,\\0<i+j+p_{ij}+q_{j}\leq k}}\prod_{i=0}^{k-1}
      \|\partial_{t}^{i}\mathfrak{F}\|_{\mathring{H}^{|\alpha|+2j}(\mathbb{T})}^{p_{ij}}
      \|\mathfrak{U}\|_{\ddot{H}^{|\alpha|+2j+1}(\hat{\Omega})}^{q_{j}}.\nonumber
\end{align}
It is clear to see that
\begin{align}\label{Dkalpha3}
\|\partial_{t}^{k-1}[D^{\alpha}\mathfrak{F}]\|_{L_{x}^{2}(\mathbb{T})}|_{y=0}
\leq \|\partial_{t}^{k-1}\mathfrak{F}\|_{\mathring{H}^{|\alpha|}(\mathbb{T})}
\leq \sum_{\substack{p_{ij}>0,\\ i+j\leq k-1,\\0<i+j+p_{ij}+q_{j}\leq k}}\prod_{i=0}^{k-1}
      \|\partial_{t}^{i}\mathfrak{F}\|_{\mathring{H}^{|\alpha|+2j}(\mathbb{T})}^{p_{ij}}
      \|\mathfrak{U}\|_{\ddot{H}^{|\alpha|+2j+1}(\hat{\Omega})}^{q_{j}}.
\end{align}
Combining \eqref{Dkalpha1}, \eqref{Dkalpha2} and \eqref{Dkalpha3} with \eqref{Dkalphae}, we have
\begin{align}\label{DkalphaE}
\|\partial_{t}^{k}[D^{\alpha}\mathfrak{U}]\|_{L_{x}^{2}(\mathbb{T})}|_{y=0}
\leq C\sum_{i=0}^{k}\|\mathfrak{U}\|_{\ddot{H}^{|\alpha|+2i+1}(\hat{\Omega})}^{k-i+1}
    +C\sum_{\substack{p_{ij}>0,\\ i+j\leq k-1,\\0<i+j+p_{ij}+q_{j}\leq k}}\prod_{i=0}^{k-1}
     \|\partial_{t}^{i}\mathfrak{F}\|_{\mathring{H}^{|\alpha|+2j}(\mathbb{T})}^{p_{ij}}
     \|\mathfrak{U}\|_{\ddot{H}^{|\alpha|+2j+1}(\hat{\Omega})}^{q_{j}}.
\end{align}
Since \eqref{DkalphaE} hold for any $\alpha\in\Gamma_{s-2k+1}$, we get \eqref{Bc1}.

Now let us turn to the prove of \eqref{Bc2}. For convenience, let
$$\mathcal{A}:=\partial_{t}-\tilde{\eta}\partial_{x}^{2},$$
and
$$\mathcal{B}:=\mathfrak{U}\partial_{x}+\partial_{y}^{-1}[\partial_{x}\mathfrak{U}]\partial_{y}.$$
By direct calculation, one can obtain from \eqref{App}$_{1}$ that for any $n\geq 2$,
\begin{align*}
\partial_{y}^{2n}\mathfrak{U}|_{y=0}
=& \partial_{y}^{2n-2}(\mathcal{B}\mathfrak{U}+\partial_{x}P+\mathfrak{F})|_{y=0}
  +\mathcal{A}\partial_{y}^{2n-2}\mathfrak{U}|_{y=0}\\
=& \sum_{i=0}^{n-2}\mathcal{A}^{i}\partial_{y}^{2n-2-2i}(\mathcal{B}\mathfrak{U}+\partial_{x}P+\mathfrak{F})|_{y=0}
  +\mathcal{A}^{n-1}\partial_{y}^{2}\mathfrak{U}|_{y=0}\\
=& \sum_{i=0}^{n-2}\mathcal{A}^{i}\partial_{y}^{2n-2-2i}(\mathcal{B}\mathfrak{U}+\mathfrak{F})|_{y=0}
  +\mathcal{A}^{n-1}\partial_{x}P.
\end{align*}
Noticing that
$\partial_{y}(\mathcal{B}\mathfrak{U})|_{y=0}=0$ and for $r\geq 2$,
\begin{align*}
\partial_{y}^{r}(\mathcal{B}\mathfrak{U})|_{y=0}
=& \partial_{y}^{r}(\mathfrak{U}\partial_{x}\mathfrak{U})|_{y=0}
  -\partial_{y}^{r}(\partial_{y}^{-1}[\partial_{x}\mathfrak{U}]\partial_{y}\mathfrak{U})|_{y=0}\\
=& \sum_{i=1}^{r-1}\binom{r}{i}\partial_{y}^{r-i}\mathfrak{U}\partial_{y}^{i}\partial_{x}\mathfrak{U}|_{y=0}
  -\sum_{i=2}^{r}\binom{r}{i}\partial_{y}^{i-1}\partial_{x}\mathfrak{U}\partial_{y}^{r-i+1}\mathfrak{U}|_{y=0}\\
=&\sum_{i=1}^{r-1}\left(\binom{r}{i}-\binom{r}{i+1}\right)\partial_{y}^{r-i}\mathfrak{U}\partial_{y}^{i}\partial_{x}\mathfrak{U}|_{y=0},
\end{align*}
one has
\begin{align}\label{Dknm}
\partial_{t}^{k}\partial_{y}^{2n}\partial_{x}^{m}\mathfrak{U}|_{y=0}
=& \sum_{i=0}^{n-2}\partial_{t}^{k}\mathcal{A}^{i}\partial_{y}^{2n-2-2i}\partial_{x}^{m}(\mathcal{B}\mathfrak{U})
  +\sum_{i=0}^{n-2}\partial_{t}^{k}\mathcal{A}^{i}\partial_{y}^{2n-2-2i}\partial_{x}^{m}\mathfrak{F}|_{y=0}
  +\partial_{t}^{k}\mathcal{A}^{n-1}\partial_{x}^{m+1}P\\
=& \sum_{i=0}^{n-2}\sum_{j=1}^{2n-3-2i}\left(\binom{2n-2-2i}{j}-\binom{2n-2-2i}{j+1}\right)
    \partial_{t}^{k}\mathcal{A}^{i}\partial_{x}^{m}(\partial_{y}^{2n-2-2i-j}\mathfrak{U}\partial_{y}^{j}\partial_{x}\mathfrak{U})|_{y=0}\nonumber\\
 &+\sum_{i=0}^{n-2}\partial_{t}^{k}\mathcal{A}^{i}\partial_{y}^{2n-2-2i}\partial_{x}^{m}\mathfrak{F}|_{y=0}
  +\partial_{t}^{k}\mathcal{A}^{n-1}\partial_{x}^{m+1}P\nonumber\\
=:&L_{1}+L_{2}+L_{3}.\nonumber
\end{align}
Notice that
\begin{align*}
\partial_{t}^{k}\mathcal{A}^{i}\partial_{x}^{m}(\partial_{y}^{2n-2-2i-j}\mathfrak{U}\partial_{y}^{j}\partial_{x}\mathfrak{U})|_{y=0}
=&\sum_{r=0}^{i}\binom{i}{r}\partial_{t}^{r+k}(-\tilde{\eta}\partial_{x}^{2})^{i-r}\partial_{x}^{m}(\partial_{y}^{2n-2-2i-j}\mathfrak{U}\partial_{y}^{j}\partial_{x}\mathfrak{U})|_{y=0}\\
= \sum_{r=0}^{i}\sum_{p=0}^{r+k}\sum_{q=0}^{i-r}\sum_{h=0}^{m}\binom{i}{r}\binom{r}{p}
 &\binom{i-r}{q}\binom{m}{h}\partial_{t}^{p}(-\tilde{\eta}\partial_{x}^{2})^{q}\partial_{x}^{h}\partial_{y}^{2n-2-2i-j}\mathfrak{U}\partial_{t}^{r+k-p}(-\tilde{\eta}\partial_{x}^{2})^{i-r-q}\partial_{x}^{m-h+1}\partial_{y}^{j}\mathfrak{U}|_{y=0},
\end{align*}
by using Sobolev imbedding inequality and \eqref{Bc1}, one has
\begin{align}\label{Dknm1}
\|L_{1}\|_{L_{x}^{2}(\mathbb{T})}
\leq& C\sum_{i=0}^{n-2}\sum_{j=1}^{2n-3-2i}\sum_{r=0}^{i}\sum_{p=0}^{r+k}\sum_{q=0}^{i-r}\sum_{h=0}^{m}
      \|\partial_{t}^{p}\partial_{x}^{2q+h}\partial_{y}^{2n-2-2i-j}\mathfrak{U}\|_{\mathring{L}^{\infty}(\mathbb{T})}
      \|\partial_{t}^{r+k-p}\partial_{x}^{2(i-r-q)+m-h+1}\partial_{y}^{j}\mathfrak{U}\|_{\mathring{H}^{0}(\mathbb{T})}\\
\leq&\mathcal{Q}_{0}\left(\|\mathfrak{U}\|_{\ddot{H}^{2n+2k+m-1}(\hat{\Omega})},\sum_{i=0}^{n-3}\|\partial_{t}^{i}F\|_{\mathring{H}^{2n+2k+m-4-2i}(\mathbb{T})}\right),\nonumber
\end{align}
where $\mathcal{Q}_{0}(\cdot,\cdot)$ is a binary polynomial function. Since
\begin{align*}
\sum_{i=0}^{n-2}\partial_{t}^{k}\mathcal{A}^{i}\partial_{y}^{2n-2-2i}\partial_{x}^{m}\mathfrak{F}|_{y=0}
=& \sum_{i=0}^{n-2}\sum_{r=0}^{i}\binom{i}{r}
   \partial_{t}^{r+k}(-\tilde{\eta}\partial_{x}^{2})^{i-r}\partial_{y}^{2n-2-2i}\partial_{x}^{m}\mathfrak{F}|_{y=0}\\
=& \sum_{i=0}^{n-2}\sum_{r=0}^{i}\binom{i}{r}
   \partial_{t}^{r+k}(-\tilde{\eta}\partial_{x}^{2})^{i-r}\partial_{y}^{2n-2-2i}\partial_{x}^{m}\mathfrak{F}|_{y=0}\\
=& \sum_{i=0}^{n-2}\sum_{r=0}^{i}C\partial_{t}^{r+k}\partial_{y}^{2n-2-2i}\partial_{x}^{m+2i-2r}\mathfrak{F}|_{y=0},
\end{align*}
it follows
\begin{align}\label{Dknm2}
\|L_{2}\|_{L_{x}^{2}(\mathbb{T})}
\leq C\sum_{i=0}^{n-2}\|\partial_{t}^{i+k}\mathfrak{F}\|_{\mathring{H}^{2n+m-2-2i}(\mathbb{T})}.
\end{align}
Similarly, due to
\begin{align*}
\partial_{t}^{k}\mathcal{A}^{n-1}\partial_{x}^{m+1}P
=\sum_{i=0}^{n-1}\binom{n-1}{i}\partial_{t}^{i+k}(-\tilde{\eta}\partial_{x}^{2})^{n-1-i}\partial_{x}^{m+1}P
\end{align*}
one has
\begin{align}\label{Dknm3}
\|L_{3}\|_{L_{x}^{2}(\mathbb{T})}
\leq C\sum_{i=0}^{n-1}\|\partial_{t}^{i+k}\partial_{x}P\|_{H^{2n+m-2-2i}(\mathbb{T})}.
\end{align}
Combining \eqref{Dknm1}, \eqref{Dknm2} and \eqref{Dknm3} with \eqref{Dknm}, we get \eqref{Bc2}.

\end{proof}

\begin{proof}[Proof of Lemma \ref{bc}]
Let $t\in[0,T_{c}]$. Noticing that ${\bf supp}~h(y)\in(\hat{y},+\infty)$, the equation \eqref{app}$_{1}$ reads in $\hat{\Omega}$ as
\begin{align}\label{Papp}
\partial_{t}u
+\mathfrak{U}\partial_{x}u+u\partial_{x}\mathfrak{U}
-(\dot{\mathfrak{U}}+\varsigma)\partial_{y}^{-1}[\partial_{x}u]
-\partial_{y}^{-1}[\partial_{x}\mathfrak{U}]\partial_{y}u
=\partial_{y}^{2}u+\eta\partial_{x}^{2}u+\kappa\partial_{x}^{2}\mathbb{U}+F.
\end{align}
Set
$$\mathcal{A}:=\partial_{t}-\eta\partial_{x}^{2},~\Delta_{\eta}:=\partial_{y}^{2}+\eta\partial_{x}^{2},$$
and
$$\mathcal{B}:=-\mathfrak{U}\partial_{x}-\partial_{x}\mathfrak{U}+(\dot{\mathfrak{U}}+\varsigma)\partial_{y}^{-1}\partial_{x}+\partial_{y}^{-1}[\partial_{x}\mathfrak{U}]\partial_{y}.$$
One can deduce from \eqref{Papp} that for any $k\in\mathbf{N}_{+}$ the following equality holds in $\hat{\Omega}$,
\begin{align*}
\partial_{t}^{k}u
=\sum_{\substack{p_{i}>0,\\~\sum\limits_{i}(i+p_{i})q_{i}=k}}C\left(\prod_{i=0}^{k-1}[\partial_{t}^{i}(\mathcal{B}+\Delta_{\eta})^{p_{i}}]^{q_{i}}\right)u
 +\sum_{\substack{p_{i}>0,\\j+\sum\limits_{i}(i+p_{i})q_{i}=k-1}}C\left(\prod_{i=0}^{k-2}[\partial_{t}^{i}(\mathcal{B}+\Delta_{\eta})^{p_{i}}]^{q_{i}}\partial_{t}^{j}\right)(\kappa\partial_{x}^{2}\mathbb{U}+F),
\end{align*}
which gives for any $\alpha\in\Gamma_{s-2k-1}$ that,
\begin{align}\label{ube}
\|\partial_{t}^{k}[D^{\alpha}u]\|_{L_{x}^{2}(\mathbb{T})}|_{y=0}
\leq&~\tilde{c}_{1}\|u\|_{\ddot{H}^{|\alpha|+2k+1}(\hat{\Omega})}
     +\tilde{c}_{2}\sum_{i=0}^{k-1}\left(\bar{\kappa}\|\partial_{t}^{i}\mathbb{U}\|_{\mathring{H}^{|\alpha|+2k-2i}(\mathbb{T})}+\|\partial_{t}^{i}F\|_{\mathring{H}^{|\alpha|+2k-2i-2}(\mathbb{T})}\right)\\
\leq&~\tilde{c}_{1}\|u\|_{\ddot{H}^{|\alpha|+2k+1}(\hat{\Omega})}
     +\tilde{c}_{2}\sum_{i=0}^{k-1}\|\partial_{t}^{i}F\|_{\mathring{H}^{|\alpha|+2k-2i-2}(\mathbb{T})}
     +\tilde{c}_{3}\bar{\kappa},\nonumber
\end{align}
where $\tilde{c}_{1},\tilde{c}_{3}$ are constants depending on
$$\|\mathfrak{U}\|_{L_{t}^{\infty}\ddot{H}^{|\alpha|+2k+1}(\hat{\Omega}_{T_{c}})},~
 \sum\limits_{i=0}^{k-2}\|\partial_{t}^{i}u_{0}\|_{L_{t}^{\infty}\mathring{H}^{|\alpha|+2k-2i}([0,T_{c}]\times\mathbb{T})}~\text{and}~
 \sum\limits_{i=0}^{k-2}\|\partial_{t}^{i}F\|_{L_{t}^{\infty}\mathring{H}^{|\alpha|+2k-2i-2}([0,T_{c}]\times\mathbb{T})},$$
and
$\tilde{c}_{2}$ is a constant depending on
$$\|\mathfrak{U}\|_{L_{t}^{\infty}\ddot{H}^{|\alpha|+2k-1}(\hat{\Omega}_{T_{c}})},~
\sum\limits_{i=0}^{k-2}\|\partial_{t}^{i}u_{0}\|_{L_{t}^{\infty}\mathring{H}^{|\alpha|+2k-2i-2}([0,T_{c}]\times\mathbb{T})},~\text{and}~ \sum\limits_{i=0}^{k-2}\|\partial_{t}^{i}F\|_{L_{t}^{\infty}\mathring{H}^{|\alpha|+2k-2i-4}([0,T_{c}]\times\mathbb{T})}.$$
By direct calculation, it follows from \eqref{Papp} that
\begin{align*}
\partial_{y}^{2n}u|_{y=0}
=&-\partial_{y}^{2n-2}(\mathcal{B}u+\kappa\partial_{x}^{2}\mathbb{U}+F)|_{y=0}
  +\mathcal{A}\partial_{y}^{2n-2}u|_{y=0}\\
=&-\sum_{i=0}^{n-2}\mathcal{A}^{i}\partial_{y}^{2n-2-2i}(\mathcal{B}u+\kappa\partial_{x}^{2}\mathbb{U}+F)|_{y=0}
  +\mathcal{A}^{n-1}\partial_{y}^{2}u|_{y=0}\\
=&-\sum_{i=0}^{n-2}\mathcal{A}^{i}\partial_{y}^{2n-2-2i}(\mathcal{B}u+\kappa\partial_{x}^{2}\mathbb{U}+F)|_{y=0}
  -\mathcal{A}^{n-1}F|_{y=0}.
\end{align*}
Noticing that
$\partial_{y}\mathcal{B}u|_{y=0}=0$ and for $k\geq 2$,
\begin{align*}
\partial_{y}^{k}\mathcal{B}u|_{y=0}
=&-\partial_{y}^{k}(\mathfrak{U}\partial_{x}u)|_{y=0}
  -\partial_{y}^{k}(u\partial_{x}\mathfrak{U})|_{y=0}
  +\partial_{y}^{k}(\partial_{y}\mathfrak{U}\partial_{y}^{-1}[\partial_{x}u])|_{y=0}
  +\partial_{y}^{k}(\partial_{y}^{-1}[\partial_{x}\mathfrak{U}]\partial_{y}u)|_{y=0}\\
=&-\sum_{i=1}^{k-1}\binom{k}{i}\partial_{y}^{k-i}\mathfrak{U}\partial_{y}^{i}\partial_{x}u|_{y=0}
  -\sum_{i=1}^{k-1}\binom{k}{i}\partial_{y}^{k-i}u\partial_{y}^{i}\partial_{x}\mathfrak{U}|_{y=0}
  +\sum_{i=2}^{k}\binom{k}{i}\partial_{y}^{k-i+1}\mathfrak{U}\partial_{y}^{i-1}\partial_{x}u|_{y=0}\\
 &+\sum_{i=2}^{k}\binom{k}{i}\partial_{y}^{i-1}\partial_{x}\mathfrak{U}\partial_{y}^{k-i+1}u|_{y=0}\\
=& \sum_{i=1}^{k-1}\left(\binom{k}{i+1}-\binom{k}{i}\right)
   \left(\partial_{y}^{k-i}\mathfrak{U}\partial_{y}^{i}\partial_{x}u|_{y=0}+\partial_{y}^{k-i}u\partial_{y}^{i}\partial_{x}\mathfrak{U}|_{y=0}\right),
\end{align*}
one has
\begin{align}\label{Dnm}
\partial_{y}^{2n}\partial_{x}^{m}u|_{y=0}
=&-\sum_{i=0}^{n-2}\mathcal{A}^{i}\partial_{y}^{2n-2-2i}\partial_{x}^{m}\mathcal{B}u
  -\kappa\sum_{i=0}^{n-2}\mathcal{A}^{i}\partial_{y}^{2n-2-2i}\partial_{x}^{m+2}\mathbb{U}|_{y=0}
  -\sum_{i=0}^{n-1}\mathcal{A}^{i}\partial_{x}^{m}F|_{y=0}\\
=& \sum_{i=0}^{n-2}\sum_{j=1}^{2n-3-2i}\left(\binom{2n-2-2i}{j+1}-\binom{2n-2-2i}{j}\right)
    \mathcal{A}^{i}\partial_{x}^{m}(\partial_{y}^{2n-2-2i-j}\mathfrak{U}\partial_{y}^{j}\partial_{x}u)|_{y=0}\nonumber\\
 &+\sum_{i=0}^{n-2}\sum_{j=1}^{2n-3-2i}\left(\binom{2n-2-2i}{j+1}-\binom{2n-2-2i}{j}\right)
    \mathcal{A}^{i}\partial_{x}^{m}(\partial_{y}^{2n-2-2i-j}u\partial_{y}^{j}\partial_{x}\mathfrak{U})|_{y=0}\nonumber\\
 &-\kappa\sum_{i=0}^{n-2}\mathcal{A}^{i}\partial_{x}^{m+2}\partial_{y}^{2n-2-2i}\mathbb{U}|_{y=0}
  -\sum_{i=0}^{n-1}\mathcal{A}^{i}\partial_{y}^{2n-2-2i}\partial_{x}^{m}F|_{y=0}\nonumber\\
:=&\sum_{i=1}^{4}\mathcal{J}_{i}.\nonumber
\end{align}
A direct calculation shows that
\begin{align*}
\mathcal{A}^{i}\partial_{x}^{m}(\partial_{y}^{2n-2-2i-j}\mathfrak{U}\partial_{y}^{j}\partial_{x}u)|_{y=0}
=&\sum_{r=0}^{i}\binom{i}{r}\partial_{t}^{r}(-\eta\partial_{x}^{2})^{i-r}\partial_{x}^{m}(\partial_{y}^{2n-2-2i-j}\mathfrak{U}\partial_{y}^{j}\partial_{x}u)|_{y=0}\\
= \sum_{r=0}^{i}\sum_{p=0}^{r}\sum_{q=0}^{i-r}\sum_{h=0}^{m}\binom{i}{r}\binom{r}{p}
 &\binom{i-r}{q}\binom{m}{h}\partial_{t}^{p}(-\eta\partial_{x}^{2})^{q}\partial_{x}^{h}\partial_{y}^{2n-2-2i-j}\mathfrak{U}\partial_{t}^{r-p}(-\eta\partial_{x}^{2})^{i-r-q}\partial_{x}^{m-h}\partial_{y}^{j}\partial_{x}u|_{y=0},
\end{align*}
which gives
\begin{align*}
\|\mathcal{J}_{1}\|_{L_{x}^{2}(\mathbb{T})}
\leq& C\sum_{i=0}^{n-2}\sum_{j=1}^{2n-3-2i}\sum_{r=0}^{i}\sum_{p=0}^{r}\sum_{q=0}^{i-r}\sum_{h=0}^{m}
      \|\partial_{t}^{p}\partial_{x}^{2q+h}\partial_{y}^{2n-2-2i-j}\mathfrak{U}\|_{\mathring{L}^{\infty}(\mathbb{T})}
      \|\partial_{t}^{r-p}\partial_{x}^{2(i-r-q)+m-h}\partial_{y}^{j}\partial_{x}u\|_{\mathring{H}^{0}(\mathbb{T})}\\
\leq& c_{*}\left(\|u\|_{\ddot{H}^{2n+m-1}(\hat{\Omega})}
                 +\sum_{i=0}^{n-3}\|\partial_{t}^{i}F\|_{\mathring{H}^{2n+m-4-2i}(\mathbb{T})}+\bar{\kappa}\right),
\end{align*}
and similarly,
\begin{align*}
\|\mathcal{J}_{2}\|_{L_{x}^{2}(\mathbb{T})}
\leq& C\sum_{i=0}^{n-2}\sum_{j=1}^{2n-3-2i}\sum_{r=0}^{i}\sum_{p=0}^{r}\sum_{q=0}^{i-r}\sum_{h=0}^{m}
      \|\partial_{t}^{p}\partial_{x}^{2q+h}\partial_{y}^{2n-2-2i-j}u\|_{\mathring{L}^{\infty}(\mathbb{T})}
      \|\partial_{t}^{r-p}\partial_{x}^{2(i-r-q)+m-h}\partial_{y}^{j}\partial_{x}\mathfrak{U}\|_{\mathring{H}^{0}(\mathbb{T})}\\
\leq& c_{*}\left(\|u\|_{\ddot{H}^{2n+m-1}(\hat{\Omega})}
                 +\sum_{i=0}^{n-3}\|\partial_{t}^{i}F\|_{\mathring{H}^{2n+m-4-2i}(\mathbb{T})}+\bar{\kappa}\right),
\end{align*}
where
$c_{*}$ is a constant depending on $\|\mathfrak{U}\|_{\ddot{H}^{2n+m-1}(\hat{\Omega})}$,
$\sum\limits_{i=0}^{n-2}\|\partial_{t}^{i}u_{0}\|_{\mathring{H}^{2n+m-2-2i}(\mathbb{T})}$, and
$\sum\limits_{i=0}^{n-2}\|\partial_{t}^{i}F\|_{\mathring{H}^{2n+m-4-2i}(\mathbb{T})}$.
Noticing that
\begin{align*}
\mathcal{J}_{3}
=-\kappa\sum_{i=0}^{n-2}\sum_{j=0}^{i}\binom{i}{j}\partial_{t}^{j}(-\partial_{x}^{2})^{i-j}\partial_{x}^{m+2}\partial_{y}^{2n-2-2i}\mathbb{U}|_{y=0},
\end{align*}
we have
\begin{align*}
\|\mathcal{J}_{3}\|_{L_{x}^{2}(\mathbb{T})}
\leq& C\bar{\kappa}\sum_{i=0}^{n-2}\sum_{j=0}^{i}
       \|\partial_{t}^{j}\partial_{x}^{m+2i-2j+2}\partial_{y}^{2n-2-2i}\mathbb{U}\|_{\mathring{H}^{0}(\mathbb{T})}\\
\leq& c^{*}\bar{\kappa},
\end{align*}
where
$c^{*}$ is a constant depending on
$$\|\mathfrak{U}\|_{L_{t}^{\infty}\ddot{H}^{2n+m-1}(\hat{\Omega}_{T_{c}})},
\sum\limits_{i=0}^{n-3}\|\partial_{t}^{i}u_{0}\|_{L_{t}^{\infty}\mathring{H}^{2n+m-2i}([0,T_{c}]\times\mathbb{T})},$$
$$\sum\limits_{i=0}^{n-3}\|\partial_{t}^{i}F\|_{L_{t}^{\infty}\mathring{H}^{2n+m-2-2i}([0,T_{c}]\times\mathbb{T})},
\sum\limits_{i=0}^{n-2}\|\partial_{t}^{i}\partial_{x}P\|_{H^{2n+m-2-2i}(\mathbb{T})}.$$
The estimate of $\mathcal{I}_{4}$ is straightforward,
\begin{align*}
\|\mathcal{J}_{4}\|_{L_{x}^{2}(\mathbb{T})}
\leq C\sum_{i=0}^{n-1}\|\partial_{t}^{i}F\|_{\mathring{H}^{2n+m-2i-2}(\mathbb{T})}.
\end{align*}
Substituting the estimate of $\mathcal{J}_{1}, \mathcal{J}_{2}, \mathcal{J}_{3}$ and $\mathcal{J}_{4}$ above into \eqref{Dnm}, one can get the conclusion.
\end{proof}

\noindent{\bf Acknowledgments:}
The first author was partially supported by National Natural Science Foundation of China under Grant No. 12101371, and Shandong Provincial Natural Science Foundation under Grant No. ZR2021QA006. The second author was partially supported by National Natural Science Foundation of China under Grant No. 12331008.

\end{document}